\documentclass[12pt,a4paper]{book} 

\usepackage[utf8]{inputenc}
\usepackage[T1]{fontenc} 
\usepackage{lmodern}  
\usepackage[english]{babel}  
\usepackage{amsmath}  
\usepackage{amsthm}   
\usepackage{amsfonts} 
\usepackage{amssymb}  

\usepackage{booktabs} 
	\usepackage{multirow} 

\usepackage[percent]{overpic}

\usepackage{bm}       
\usepackage{lipsum}   
\usepackage[
	bookmarksnumbered = true,
	colorlinks=true,
	linkcolor=blue,
	anchorcolor=black,
	citecolor=green!50!black,
	urlcolor=blue
	]
	{hyperref}       
\usepackage{graphicx} 
\usepackage{caption}  

\usepackage{adjustbox}

\usepackage{algorithm}

\usepackage[breakable,skins,most]{tcolorbox}

\usepackage{enumitem}

\usepackage{color}
\makeatletter
\def\mathcolor#1#{\@mathcolor{#1}}
\def\@mathcolor#1#2#3{%
	\protect\leavevmode
	\begingroup\color#1{#2}#3\endgroup
}
\makeatother

\definecolor{green1}{rgb}{0.3, 0.8, 0.3}

\newenvironment{Side}[1][]{%
\begin{tcolorbox}[enhanced,breakable,sharp corners,
	boxrule=0mm,right=0pt,
	leftrule=0.3cm,grow to left by=0mm, 
	interior style={fill=white},
	frame style={fill=green, fill opacity=0.15},
	left=0mm,right=0mm,
	before upper={\setlength{\parindent}{1em}%
		\everypar{{\setbox0\lastbox}\@minipagefalse\everypar{}}},
		title  = #1, 
	coltitle=white,
title style={white!85!black},
	]
	\ignorespaces
}
{\unskip
\end{tcolorbox}%
\ignorespacesafterend
}

\newenvironment{GraphicalInterpretationbox}[1][]{%
	\begin{tcolorbox}[enhanced,breakable,sharp corners,
		boxrule=0mm,right=0pt,
		leftrule=0.15cm, rightrule=0.15cm, grow to left by=0mm, 
		toprule=0.15cm,bottomrule=0.15cm,
		interior style={fill=white},
		frame style={fill=green, fill opacity=0.15},
		left=0mm,right=0mm,
		]
		\ignorespaces
		\begin{GraphicalInterpretation}}{\end{GraphicalInterpretation}
		\unskip
	\end{tcolorbox}
	\ignorespacesafterend
}

\newenvironment{exbox}[1][]{%
	\begin{tcolorbox}[enhanced,breakable,sharp corners,
		boxrule=0mm,right=0pt,
		leftrule=0.15cm,rightrule=0.15cm,grow to left by=0mm, 
		toprule=0.15cm,bottomrule=0.15cm,
		interior style={fill=white},
		frame style={fill=green, fill opacity=0.15},
		left=0mm,right=0mm,
		]
		\ignorespaces
		\begin{Example}}{\end{Example}
		\unskip
	\end{tcolorbox}
	\ignorespacesafterend
}

\newenvironment{thmbox}[1][]{%
	\begin{tcolorbox}[
		title  = #1, 
		enhanced,
		breakable=true,
		arc=0mm,
		outer arc=0mm,
		boxrule=0pt,
		toprule=1pt,
		leftrule=0pt,
		bottomrule=1pt, 
		rightrule=0pt,
		left=0.2cm,
		right=0.2cm,
		titlerule=0.5em,
		toptitle=0.1cm,
		bottomtitle=0.2cm,
		top=0.2cm,
		colframe=white!85!black,
		colback=white!85!black,
		coltitle=white,
		title style={white!85!black},
		fonttitle=\bfseries,fontupper=\normalsize,
		after=\par]
	\begin{Theorem}}
	{\end{Theorem}
	\end{tcolorbox}%
}

\newenvironment{lembox}[1][]{%
	\begin{tcolorbox}[
		title  = #1, 
		enhanced,
		breakable=true,
		arc=0mm,
		outer arc=0mm,
		boxrule=0pt,
		toprule=1pt,
		leftrule=0pt,
		bottomrule=1pt, 
		rightrule=0pt,
		left=0.2cm,
		right=0.2cm,
		titlerule=0.5em,
		toptitle=0.1cm,
		bottomtitle=0.2cm,
		top=0.2cm,
		colframe=white!85!black,
		colback=white!85!black,
		coltitle=white,
		title style={white!85!black},
		fonttitle=\bfseries,fontupper=\normalsize,
		after=\par]
		\begin{Lemma}}
		{\end{Lemma}
	\end{tcolorbox}%
}

\newenvironment{proofbox}[1][]{%
	\begin{tcolorbox}[enhanced,breakable,sharp corners,
		boxrule=0mm,right=0pt,
		leftrule=0.3cm,grow to left by=0mm, 
		interior style={fill=white},
		frame style={fill=black, fill opacity=0.15},
		left=0mm,right=0mm,
		before upper={\setlength{\parindent}{1em}%
			\everypar{{\setbox0\lastbox}\@minipagefalse\everypar{}}},
		before=\par,
		before skip=0pt,
		]
		\ignorespaces
		\begin{proof}}{\end{proof}
		\unskip
	\end{tcolorbox}
	\ignorespacesafterend
}	

\newtcbtheorem{BoxAlgo}{Alg.}{
	enhanced,
	breakable=true,
	arc=0mm,
	outer arc=0mm,
	boxrule=0pt,
	toprule=1pt,
	leftrule=0pt,
	bottomrule=2pt, 
	rightrule=0pt,
	left=0.2cm,
	right=0.2cm,
	titlerule=0.5em,
	toptitle=0.1cm,
	bottomtitle=-0.1cm,
	top=0.2cm,
	colframe=white!25!black,
	colback=white,
	coltitle=white,
	title style={white!25!black},
	fonttitle=\bfseries,fontupper=\normalsize
}{algo}

\newtcbtheorem{BoxFun}{function}{
	separator sign none,
enhanced,
breakable=true,
arc=0mm,
outer arc=0mm,
boxrule=0pt,
toprule=2pt,
leftrule=0pt,
bottomrule=2pt, 
rightrule=0pt,
left=0.2cm,
right=0.2cm,
titlerule=1.5pt,
toptitle=0.1em,
bottomtitle=0.0cm,
top=0.2cm,
colframe=black,
colback=white,
coltitle=black,
title style={white!100!black},
fonttitle=\bfseries,fontupper=\normalsize
}{fun}

\usepackage{geometry}
\usepackage{fancyhdr}
\fancypagestyle{plain}{%
	\chead{}     
	\lfoot{}     
	\cfoot{}     
	\rfoot{}     
	\fancyhead{}
	\fancyfoot{} 
	\fancyfoot[RO,LE]{\thepage}
	}

\usepackage{tocbasic}
	\DeclareTOCStyleEntry[
		beforeskip=.2em plus 1pt,
		pagenumberformat=\textbf
		]{tocline}{chapter}

\usepackage[pages=some]{background}
\backgroundsetup{
	scale=1,
	color=black,
	opacity=1,
	angle=0,
	contents={\includegraphics[width=\paperwidth]{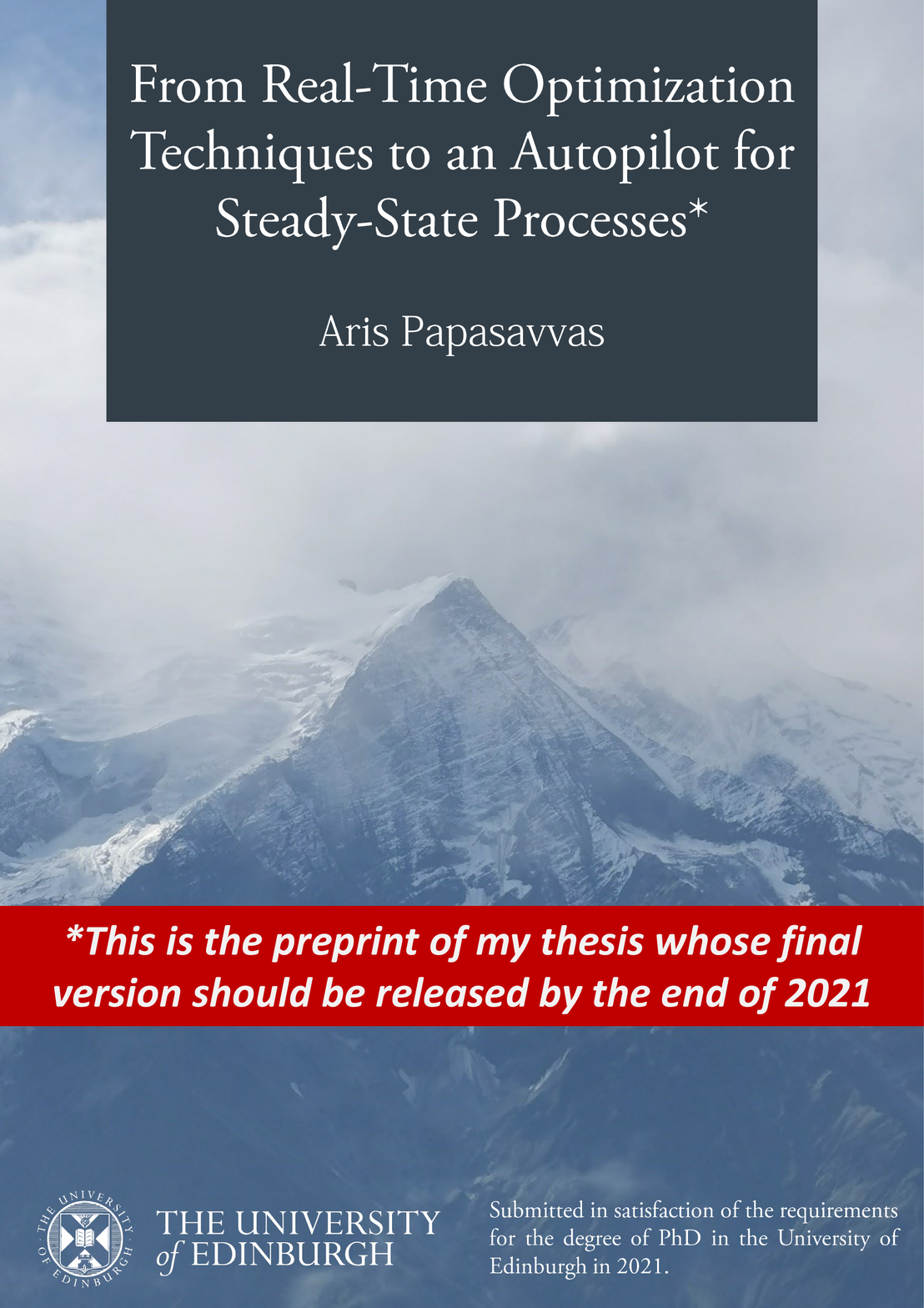}} 
}
\usepackage{tikz}

\usetikzlibrary{arrows.meta}

\theoremstyle{definition} 

\newtheorem{Lemma}{Lemma}
\newtheorem{Theorem}{Theorem}
\newtheorem{Definition}{Definition}
\newtheorem{Corollary}{Corollary}

\newtheorem{Example}{Example}

\newtheorem{Assumption}{Assumption}
\newtheorem{Remark}{Remark}

\newtheorem{GraphicalInterpretation}{Graphical Interpretation}

\usepackage{chngcntr}
\counterwithin{Conjecture}{chapter}
\counterwithin{Theorem}{chapter}
\counterwithin{GraphicalInterpretation}{chapter}
\counterwithin{Definition}{chapter}
\counterwithin{Corollary}{chapter}
\counterwithin{Proposition}{chapter}
\counterwithin{Claim}{chapter}
\counterwithin{Heuristic}{chapter}
\counterwithin{Algorithm}{chapter}
\counterwithin{Example}{chapter}
\counterwithin{Problem}{chapter}
\counterwithin{Assumption}{chapter}
\counterwithin{Remark}{chapter}
\counterwithin{Criterion}{chapter}
\counterwithin{Criteria}{chapter}
\counterwithin{Lemma}{chapter}

\usepackage{mathbbol}
\DeclareSymbolFontAlphabet{\amsmathbb}{AMSb}%

\usetikzlibrary{shapes.misc}

\usepackage[customcolors]{hf-tikz}



\definecolor{magenta_ASP}{rgb}{0.87, 0.81, 0.9}
\definecolor{magenta_2}{rgb}{1, 0.8, 1}
\definecolor{bleu_ASP}{rgb}{0.7, 0.78, 0.9}
\definecolor{vert_ASP}{rgb}{0.66, 0.82, 0.56}

\definecolor{magenta_clair}{rgb}{1,0.8, 1}
\definecolor{magenta_clair2}{rgb}{0.7137, 0.5373, 0.7137}
\definecolor{gold}{rgb}{0.749, 0.5647, 0}
\definecolor{blue_clair2}{rgb}{0.6157, 0.7725, 0.9020}

\definecolor{gris_clair}{rgb}{0.5, 0.5,0.5}
\definecolor{gris_clair2}{rgb}{0.7, 0.7,0.7}
\definecolor{gris_fonce}{rgb}{0.4, 0.4,0.4}
\definecolor{blue_clair}{rgb}{0.4500 ,   0.4500  ,  0.7200}
\definecolor{vert_fonce}{rgb}{0.100 ,   0.4500  ,  0.4500}

\newcommand\Overline[2][1pt]{%
	\begin{tikzpicture}[baseline=(a.base)]
		\node[inner xsep=0pt,inner ysep=1.5pt] (a) {$#2$};
		\draw[line width= #1] (a.north west) -- (a.north east);
	\end{tikzpicture}
}
\newcommand\Underline[2][1pt]{%
	\begin{tikzpicture}[baseline=(a.base)]
		\node[inner xsep=0pt,inner ysep=1.5pt] (a) {$#2$};
		\draw[line width= #1] (a.south west) -- (a.south east);
	\end{tikzpicture}
}

\newcommand\fe{\bm{\mathfrak{f}}}
\newcommand\g{\bm{\mathfrak{g}}}

\usepackage{undertilde}

\let\emptyset\varnothing

\let\b\relax 
\newcommand{\b}{\mathsf{b}}
\newcommand{\erf}{\text{erf}}
\newcommand{\x}{\mathsf{x}}

\newcommand{\cc}{\mathsf{c}}
\newcommand{\CC}{\mathsf{C}}
\newcommand{\n}{\mathsf{n}}
\let\d\relax 
\newcommand{\d}{\mathsf{d}}
\newcommand{\D}{\mathsf{D}}

\newcommand{\f}{\mathsf{f}}

\let\k\relax 
\newcommand{\k}{\mathsf{k}}
\newcommand{\m}{\mathsf{m}}

\let\L\relax 
\newcommand{\L}{\mathsf{L}}

\let\v\relax 
\newcommand{\v}{\mathsf{v}}
\let\V\relax 
\newcommand{\V}{\mathsf{V}}

\usepackage{longtable}
\begin{document}

\begin{titlepage} 
	\BgThispage
\end{titlepage}
\ \ \ 
\newpage
\thispagestyle{empty}
\ \ \ 
\newpage
\thispagestyle{empty}
\null\vspace{\stretch {1}}
\begin{flushright}
\begin{minipage}{6cm}
	\textit{To my sisters, to my parents.}
\end{minipage}
\end{flushright}
\vspace{\stretch{2}}\null

\newpage
\thispagestyle{empty}
\ \ \ 
	\pagenumbering{roman}
	\chapter*{Declaration}
\addcontentsline{toc}{chapter}{Declaration}

I declare that this thesis was composed by myself and that the work contained herein is my own, except where explicitly stated otherwise in the text. 

\begin{flushright}
	\textit{Aris Papasavvas}
\end{flushright}
	\chapter*{Acknowledgment}
\addcontentsline{toc}{chapter}{Acknowledgment}

This thesis would not have been possible without the presence of many people. First of all I would like to thank my supervisor: Dr Gregory Francois, for his constant support, his kindness  and all the freedom he allowed me in my research. Likewise, I would like to express my gratitude to:  Professor Dominique Bonvin, Dr. Alejandro Marchetti, and Dr. Tafarel de Avila Ferreira, for offering me the same ideal working environment, and for introducing me to the exciting field of research that real-time optimization can be.  I am very thankful to my co-supervisor: Dr. Dimitrios Gerogiorgis,  for all the advice he gave me concerning the planning of my PhD and my professional career.

However, a thesis is not only a time and means that is allocated to a research, it is also an environment that is shared and faced with the help of the people around us. These three and a half years would certainly not have been the same without Dr. Pierre Fayon, and the future doctors Robin Delbart, Christophe Floreani, and Alasdair Speakman. 

Finally, I thank my parents Marie-Pierre and Savvas, my sisters  Dora and Skévi, and my almost-brother Nicolas for their support, and in a way their presence, despite the distance and the events that have kept us apart for far too long.
	\chapter*{Lay Summary}
\addcontentsline{toc}{chapter}{Lay Summary}

To explain what the subject of this thesis is, I propose the following metaphor: 

Imagine a bakery whose owner's objective is to produce a certain quantity of a given quality of bread at the lowest price, e.g. a few loaves of very high quality bread to target a small group of exigent consumers or many loaves of low quality bread to supply school canteens.

This bakery is run by a baker, whose job is to ensure that the owner's objectives are met while using the equipment available, which can be summarized as a simple oven (the idea here is not to explain the actual operation of a bakery).  This oven has an interface that allows you to do two things: (i) set the cooking temperature and (ii) set the cooking time. Between this interface and the actual operation of the oven there are so-called regulators that measure the temperature of the oven in real time and according to the deviation from the target temperature, increase or reduce the power consumed by the oven. 

So on one side there is the owner who takes decisions on the long term and at a very low frequency (gives new objectives every month, quarter, or even year). On the other hand there is the oven that makes decisions via its regulators on a very short term and very high frequency (almost continuously). In the middle, there is the baker who must make the link between the owner's objectives and the instructions to be given to the oven. To make this link, the baker can use his expertise (his knowledge of bread baking) as well as his experience in using the tools he has at his disposal, i.e. the oven.  Indeed, the more he uses it, the better he understands its functioning, its imperfections, etc. 

The study of how one can combine this experience with the data collected over time in order to better achieve the imposed objectives (by the owner) given the available tools (the oven with its controllers), is the study of real-time optimization. With this thesis, one tries to discover what the baker should do so that in the shortest possible time and with the least amount of risk, he can identify the recipe that best satisfies the objectives of his boss.  However, this thesis is much more focused on industrial issues where the quality and safety requirements as well as operating costs have much more significant consequences.  

	\chapter*{Abstract}
\addcontentsline{toc}{chapter}{Abstract}

Any industrial system goes along with objectives to be met (e.g. economic performance), disturbances to handle (e.g. market fluctuations, catalyst decay, unexpected variations in uncontrolled flow rates and compositions, …), and uncertainties about its behavior. In response to these, decisions must be taken and instructions be sent to the operators to drive and maintain the plant at satisfactory, yet potentially changing operating conditions.  \\

Over the past thirty years many methods have been created and developed to answer these questions. In particular, the field of Real-Time Optimization (RTO)  has emerged that, among others, encompasses methods that allow the systematic improvement of the performances of the industrial system, using plant measurements and a potentially inaccurate tool to predict its behaviour, generally in the form of a model. Even though the definition of RTO can differ between authors, inside and outside the process systems engineering community, there is currently no RTO method, which is deemed capable of fully automating the aforementioned decision-making process. This thesis consists of a series of contributions in this direction, which brings RTO closer to being capable of a full plant automation. \\

To this end, two main development strategies have been used.  \\ 

With the first strategy, methodological improvements and new ways to implement existing RTO approaches have been proposed. By methodological improvements it is meant both improvement of the current implementation and new ways to alleviate intrinsic weaknesses of RTO techniques. In particular, one recent and popular RTO technique is considered and the first  three scientific contributions are proposed:
\begin{itemize}
	\item The standard implementation of this RTO technique considers a single tuning parameter, i.e. a gain matrix. In this thesis (Chapter 2), a way to tune this free parameter is proposed, and it is also shown, among other things, that this choice will almost always be the best possible.
	\item The same popular RTO technique is considered and is shown to have a non-negligible risk of violating plant constraints. The second main contribution of this thesis (Chapter 3) consists of a structural modification of the way RTO is implemented that significantly reduces this risk.
	\item In Chapter 5, many new ways to implement a RTO method are proposed. More precisely, several ways of correcting a model on the basis of measurements are identified that do not consist of a direct measurement-based correction of model predictions, a choice that is almost systematically made with RTO approaches.
\end{itemize}

The second strategy consists of attempting to build a prototype of what would be an “ideal” RTO method, would it exist. In Chapter 4 the properties that such an ``ideal'' method should have is proposed, together with how these functionalities should interact. As this is an introductory chapter to this new structure, each of these functions is proposed to be as simple as possible. As a result, a new framework labeled Simple Autopilot for Steady-State Processes (S-ASP) is proposed. The same structure is considered in Chapter 6, but this time all its functionalities are developed further and made as effective as possible. The result is a much more complex and integrated RTO method labeled ``Autopilot for Stationary Processes (ASP)'' and which is shown to provide very interesting practical results. \\

In addition to these main contributions, the interested reader will also find in this thesis:
\begin{itemize}
	\item Chapter 1: What is thought to be an exhaustive presentation of the environment with which a RTO approach has to deal, when applied to an industrial process. Additionally, to best of the author’s knowledge, no article to date has provided such a detailed presentation. 
	\item Chapter 2: A gradual and didactic introduction to RTO that does not require any prior knowledge. This chapter has been built to be the best possible starting point for getting familiar with RTO (the way I would have appreciated it at the beginning of my PhD thesis).
\end{itemize}

It should be noted that the several contributions in this PhD can be applied to other research fields such as:
\begin{itemize}
	\item “Pure” optimization: One of the contributions presented in Chapter 2 can be applied to the selection of a ``free'' parameter of a well-known numerical optimization method (sequential quadratic optimization -- SQP), as explained in the conclusion.
	\item Reduced-order model optimization (ROM-Optimization): Replacing the words ``model'' and `plant'' by ``high-fidelity model'' and ``low-fidelity model'', respectively, is sufficient to show that everything that is discussed, proposed and developed in chapter 2 can be applied to ROM-Optimization without needing any modification or adaptation. This is underlined in the conclusion, whereby the conceptual similarities between this field and RTO are discussed, and where it is also shown how these two research fields could be combined. To the best of the author’s knowledge, these similarities have been reported before but never discussed with this level of details.
\end{itemize}

Finally, each contribution is illustrated by means of a case study from chemical engineering, such as the Williams-Otto reactor, the Tennessee Eastman challenge process, and the Williams-Otto process. Moreover, all codes are made available (\href{https://github.com/arpsa/Thesis/}{open source}) and are commented so that they can be used to, e.g. reproduce and/or improve the results presented in this manuscript, or also to apply them to industrial processes. \\

\textbf{Keywords: }Real-time optimization, Decision-making, Optimization, Reduced-order-model optimization, Autopilot for steady-state processes, Operational research.

	\tableofcontents
	\chapter{Introduction}
		\setcounter{page}{1}    
		\pagenumbering{arabic}	

\setcounter{page}{1}    
\pagenumbering{arabic}

\section{Motivation}

\begin{figure}[H]
	\centering
	\includegraphics[width=1\linewidth] {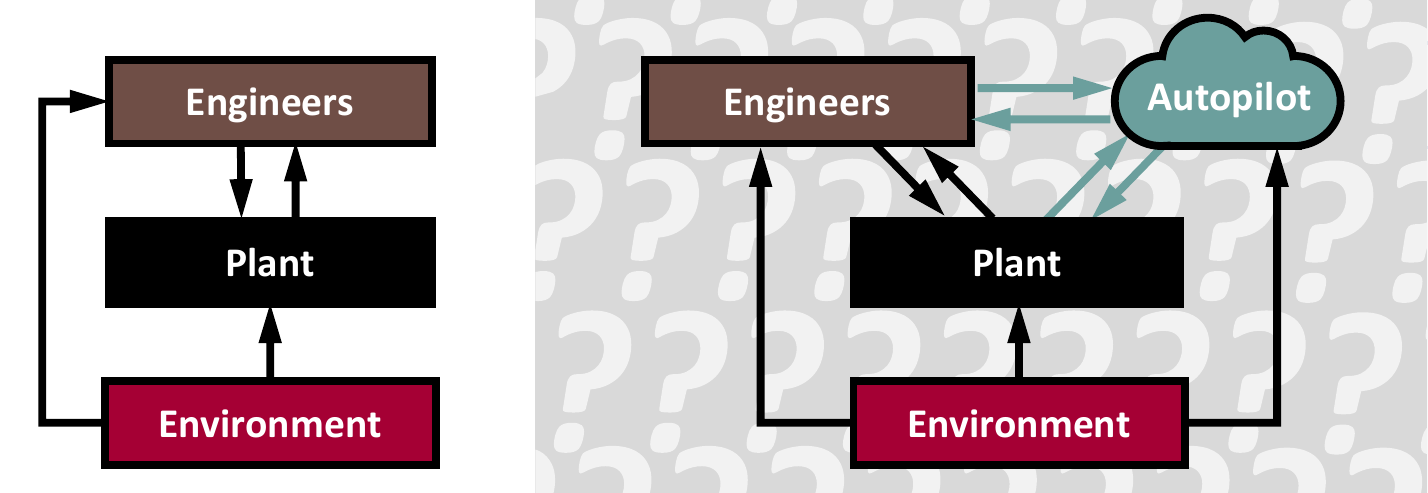}
	\caption{RTO's goal could be to build autopilots to help engineers to drive plants.}
	\label{fig:1___Context}
\end{figure}

Let's place ourselves in the position of an engineer whose mission is to manage a factory and ask ourselves the following question: \textit{What do I have to do to improve the plant's production as much as possible?} 
Which means: 
\begin{itemize}[noitemsep]
	\item To guarantee a predefined quality of production at all times.
	\item To comply at all times with the safety regulations and legislation of the country in which the plant is located
	\item To maximize the profits generated by the plant throughout its life.
\end{itemize}
The set of answers to this question could define a research field called \textit{real-time optimization} (RTO). The ultimate goal of RTO is to build an interface that facilitates the work of this engineer as much as possible while providing a  high level of guarantees on the quality of the plant's monitoring. As in some cases such an interface could manage the piloting of a plant in an autonomous way, one will simply call it ``\textit{Autopilot}'' (see Figure~\ref{fig:1___Context}). 

\section{The autopilot environment}

Before thinking about how the RTO methods work, it is necessary to define what they interact with, i.e. the three other objects in Figure~\ref{fig:1___Context}:

\subsection{The plant}

The plant is a system whose inputs $\bm{x}_p\in\amsmathbb{R}^{n_{x_p}}$ are decision variables $\bm{u}\in\amsmathbb{R}^{n_u}$ and disturbances $\bm{d}_p\in\amsmathbb{R}^{n_{d_p}}$:
\begin{equation}
	\bm{x}_p := \left[\bm{u}^{\rm T}, \bm{d}_p^{\rm T}\right]^{\rm T}.
\end{equation}
More precisely:
\vspace{-\topsep}
\begin{itemize}[noitemsep]	
	\item The decision variables $\bm{u}$ are manipulable, usually the setpoints of the plant's controllers, and   continuous in $\amsmathbb{R}^{n_u}$. \textit{In this work we do not consider cases where decision variables are binary or integer variables.}
	\item The disturbances $\bm{d}_p$ are imposed to the plant by the environment. These disturbances can be, for example, variations in the quality of the raw materials that feed the plant, the effects of the weather on the plant, etc.  
\end{itemize}
The outputs of the plant $\bm{y}_p\in\amsmathbb{R}^{n_y}$ gather all the measured variables. 

The plant is considered to be operating continuously (i.e. processing materials without interruption), and to be designed to operate at steady-state (SS).  The relationship between the inputs $\bm{x}_p$ and the outputs $\bm{y}_p$ at SS is materialized with the function $\bm{f}_p$: 
\begin{equation}
	\bm{y}_p := \bm{f}_p(\bm{x}_p).
\end{equation}

The measures $\widehat{\bm{y}}_p$ of $\bm{y}_p$ are typically polluted by measurement errors $\bm{\epsilon}_y$, and in this work one considers that these errors follow an unbiased normal random distribution: 
\begin{align} \label{eq:1___3_Mesures_Normal_Uncertainty}
	\widehat{\bm{y}}_p := \ & \bm{y}_p +  \bm{\epsilon}_y, &
	\bm{\epsilon}_y \sim \mathcal{N}(\bm{0},\bm{N}_y),
\end{align}
where $\bm{N}_y\in\amsmathbb{R}^{n_y \times n_y}$ is the covariance matrix of the measurement errors $\bm{\epsilon}_y$ (which can itself be a function of the inputs $\bm{x}_p$ of the plant).

In short, the plant receives disturbances $\bm{d}_p$ from the environment, receives instructions $\bm{u}$ either from the engineers or from the autopilot, and returns measurements $\widehat{\bm{y}}_p$ to the engineers and to the autopilot.

\subsection{The environment}

The environment gathers all the elements that are external to the plant and that cannot be manipulated either by the engineers or by the autopilot. Two distinct elements in the environment have been identified:
The first is the disturbances $\bm{d}_p$ that affect the plant, which have already been discussed in the section introducing the plant.   

The second one gathers the operating cost of the plant which can be materialized with a function $\phi(\bm{u},\bm{y}_p)\in\amsmathbb{R}$, and operating constraints that can be materialized with a function $\bm{g}(\bm{u},\bm{y}_p)\in\amsmathbb{R}^{n_g}$. More specifically, operational costs and constraints are functions based on (i) agreements with suppliers and customers, (ii) safety and environmental regulations, or (iii) empirical limitations aimed at extending the life of the plant. 

Essentially, the Environment acts directly on the plant through disturbances  $\bm{d}_p$, and guides the decision making of the Engineers and/or of the Autopilot by providing the functions $\phi$ and $\bm{g}$.  

\subsection{The engineers}

On their side, the Engineers receive the measurements $\widehat{\bm{y}}_p$ from the plant and the functions $\phi$ and $\bm{g}$ from the Environment. They can then use their technical knowledge to build a model $\bm{f}$ of the functions $\bm{f}_p$ to guide their decisions. The inputs $\bm{x}$ of this model gather the decision variables $\bm{u}$ and the subset of the disturbances $\bm{d}_p$ which is part of the measures $\widehat{\bm{y}}_p$.  It is thus necessary to distinguish the measured disturbances $\bm{d}^{\prime}_p$ from those not measured $\bm{d}^{\prime\prime}_p$: 
\begin{equation}
	\bm{d}_p := \left[{\bm{d}^{\prime}_p}^{\rm T},
	{\bm{d}^{\prime\prime}_p}^{\rm T} \right]^{\rm T}. 
\end{equation}
The model of the plant is therefore: 
\begin{align}
	\bm{y} := \ & \bm{f}(\bm{x},\bm{\theta}), &
	\text{where: } \  \bm{x} := \ & \left[{\bm{u}}^{\rm T},
	{\bm{d}^{\prime}_p}^{\rm T} \right]^{\rm T},
\end{align}
where  $\bm{\theta}\in\amsmathbb{R}^{n_{\theta}}$ gathers all the parameters of the model, e.g. the values of the known but unmeasured perturbations, the empirical physical properties of the transformed materials, etc. 

On the basis of this model, the measurements $\widehat{\bm{y}}_p$ and the functions $\phi$ and $\bm{g}$, the Engineers have to choose  the appropriate inputs $\bm{u}$ to be sent to the plant. But they can also give the model to the Autopilot --- if necessary parameterize the Autopilot --- and let it handle $\bm{u}$. In return, the Autopilot should return information to the engineers to facilitate the plant's supervision. 

\section{State of the art}

In order to facilitate the explanation of the state of the art, it has been divided into three areas
\begin{itemize}
	\item \textit{Theoretical RTO}{\color{blue},} which gathers all the contributions that focus on the identification of the optimal operating point of the plant when the environment and the experimental conditions are ideal.
	\item \textit{Practical RTO}{\color{blue},} which gathers all the contributions that focus on maintaining of the properties of theoretical RTO methods when the environment and the experimental conditions are no longer ideal. 
	\item \textit{Efficient RTO}{\color{blue},} which gathers all the contributions that question the fundamental principles of "classical" RTO methods in order to improve their performances. These improvements are generally associated with an increase in the complexity of the RTO algorithm.
\end{itemize}

\subsection{Theoretical RTO}
\label{sec:1_3_1_OTR_Theorique}

	\subsubsection{The theoretical RTO's problem formulation}
	\label{sec:1_3_1_a_Probleme_OTR_Theorique}
	
	Let's consider the ``ideal'' conditions where:
	\begin{itemize}[noitemsep]
		\item The environment does not produce any disturbance $\bm{d}_p = \emptyset$. So the functions of the plant and of the model can be reduced to:
		\begin{align*}
			\bm{y}_p = \ & \bm{f}_p(\bm{u}), &
			\bm{y}   = \ & \bm{f}(\bm{u},\bm{\theta}).
		\end{align*}
		\item There are no measurement errors: $\widehat{\bm{y}}_p=\bm{y}_p.$
		More precisely, one considers that the values of  $\phi_p$, $\bm{g}_p$, and $\bm{f}_p)$ and their gradients $(\nabla_{\bm{u}} \phi_p,\nabla_{\bm{u}}\bm{g}_p,\nabla_{\bm{u}}\bm{f}_p)$ can be estimated without error $\forall \bm{u}\in\amsmathbb{R}^{n_u}$, where:
		\begin{align*}
			\phi_p(\bm{u}) := \ & \phi(\bm{u},\bm{f}_p(\bm{u})),\\
			\bm{g}_p(\bm{u}) := \ & \bm{g}(\bm{u},\bm{f}_p(\bm{u})).
		\end{align*}
	\end{itemize}
	
	Therefore, the theoretical RTO problem consists to (\textit{objective 1}) solve the following nonlinear optimization problem (NLP):
	\begin{align}
		\bm{u}_{p}^{\star}:=  \operatorname{arg}
		\underset{\bm{u}}{\operatorname{min}} \quad &  \phi_p(\bm{u})  \quad 
		\text{s.t.} \quad  \bm{g}_p(\bm{u})  \leq \bm{0},   \nonumber
	\end{align}
	while (\textit{objective 2}) evaluating the plant functions a minimal number of times, and (\textit{objective 3})  avoiding constraints violations during those evaluations. It is important to understand that an evaluation of these functions is an experiment to be conducted in the real world on the actual plant. And usually such an experience is costly in time and money (hence \textit{objective 2}), and presents risks for the plant (hence \textit{objective 3}).  
	
	\subsubsection{Achieving objective 1}
	
	The first RTO method was proposed in 1970 in \cite{Duncanson:1970} and is known as: the two-step approach (TS). As the name suggests, this method consists of repeating two actions until they converge. (i) An update of the parameters $\bm{\theta}$ of the model  based on the most recent measurements of the plant's inputs and outputs. (ii) An updated-model-based optimization to identify the point $\bm{u}$ around which the next experiments must be conducted in order to re-identify the parameters $\bm{\theta}$. Unless the error between the functions  $\bm{f}$ and $\bm{f}_p$ is of parametric nature, i.e. $\exists \bm{\theta}\in\amsmathbb{R}^{n_{\theta}}$ such that $\bm{f}(\bm{u},\bm{\theta})=\bm{f}_p(\bm{u})$, it has been shown that TS does not generally converge on the optimal decisions \cite{Roberts:79}.
	
	The first theoretical RTO method to introduce structural corrections was proposed in 1979 in \cite{Roberts:79} and is known as: integrated system optimization and parameter estimation (ISOPE). This method consists in adding one action to TS. The three actions of ISOPE are: (i) An update of the parameters $\bm{\theta}$ of the model on the basis of the most recent measurements of the inputs and outputs of the plant. (ii - \textit{new}) The identification of an artificial parameter (called ``a modifier'') bringing an affine  correction to the cost function of the model. (iii)  An updated-model-based optimization to identify the point $\bm{u}$ around which the next experiments must be conducted in order to re-identify the parameters $\bm{\theta}$ and the modifier $\nabla_{\bm{u}}\phi_p$. It has been shown that, provided that the optimal decision $\bm{u}_p^{\star}$ does not activate any constraints, ISOPE guarantees that $\bm{u}_p^{\star}$ is identified upon convergence. 
	
	The first theoretical RTO method to guarantee unconditionally that if a convergence is observed then it can only be on the optimal decisions $\bm{u}_p^{\star}$ was proposed in 2005 in \cite{Gao:05}  and is known as: iterative setpoint optimization (ISO). This method drops the idea of updating the model parameters to focus only on updating modifiers that are now applied to both the cost and constraint functions. ISO works in two steps (i) An update of the modifiers to bring an affine correction to the costs and constraints functions of the model. (ii) An updated-model-based optimization to identify the point $\bm{u}$ around which the next experiments must be conducted in view of the re-identification of the modifiers which requires the estimates of $\nabla_{\bm{u}}\phi_p$, $\bm{g}_p$ and $\nabla_{\bm{u}}\bm{g}_p$.
	
	The first theoretical RTO method using a filter to increase the chances converging on the optimal decision $\bm{u}_p^{\star}$ was proposed in 2009 in \cite{Marchetti:09}  and is known as: modifier adaptation (MA). This method is simply ISO to which one adds a filter either on the update of the modifiers or on the movement in the input space $\bm{u}$. 
	
	The first theoretical RTO method to provide convergence properties similar to those of MA \textit{via an indirect correction of the costs and constraints of the model} was proposed in 2009 in \cite{Marchetti:09} and is known as output modifier adaptation (MAy).  This method is mentioned in \cite{Marchetti:09,Marchetti:09b} and is fully analyzed in \cite{Papasavvas:2019a}. MAy works exactly like MA with the difference that the modifiers are used to apply an affine correction to the function $\bm{f}$. And as demonstrated in \cite{Papasavvas:2019a} this is sufficient to guarantee optimality upon convergence. 
	
	The first theoretical RTO method to provide the guarantee that it is always possible to converge on the optimal decision $\bm{u}_p^{\star}$ was proposed in 2013 in \cite{Francois:2013} and can be named as proposed in \cite{Papasavvas:2019}: modifier adaptation with convexified problem (MAc). This method consists of simply applying MA while replacing the model by a convex approximation of itself.

	\subsubsection{Achieving objective 2}
	
All the methods mentioned in the previous section require that at each iteration the values $(\phi_p,\bm{g}_p,\bm{f}_p)$ and the gradients $(\nabla_{\bm{u}}\phi_p,\nabla_{\bm{u}}\bm{g}_p,\nabla_{\bm{u}}\bm{f}_p)$are measured on the plant.  Although it is considered that these values can be measured in a perfect way, measuring them accurately presents an irreducible minimal cost in terms of the number of experiments to be performed. For example, to assess the values  $(\phi_p,\bm{g}_p,\bm{f}_p)$ at a given point requires only one experiment whereas to evaluate $(\nabla_{\bm{u}}\phi_p,\nabla_{\bm{u}}\bm{g}_p,\nabla_{\bm{u}}\bm{f}_p)$ requires at least $n_u+1$ experiments (if one applies the most basic finite difference method). So, if the problem's dimension is high, i.e. $n_u \gg 1$, then each iteration is likely to be costly in terms of number of experiments and it is therefore of interest to find a way to reduce this number. 
	
	This research direction has been studied in \cite{Costello:2016} where it is suggested to evaluate the gradient only in privileged directions. This method is improved in \cite{Singhal:2018} where a strategy that allows  these directions to be adapted in real time to minimize optimality losses (due to incomplete gradient evaluations) at convergence is proposed. To evaluate these directions, these methods require engineers to provide a relevant estimate of the set $\mathbb{\Theta}$ in which the parameters of the model $\bm{\theta}$ must be such that $\forall \bm{u}\in\amsmathbb{R}^{n_u}$,  $\forall \bm{\theta}\in \mathbb{\Theta}$: $\nabla_{\bm{u}}\phi(\bm{u},\bm{f}(\bm{u},\bm{\theta})) \approx \nabla_{\bm{u}}\phi_p(\bm{u})$ and $\nabla_{\bm{u}}\bm{g}(\bm{u},\bm{f}(\bm{u},\bm{\theta})) \approx \nabla_{\bm{u}}\bm{g}_p(\bm{u})$. These methods will therefore be systematically efficient when both the modeling error is purely parametric, and the set  $\mathbb{\Theta}$ is small and contains the ``plant's parameters'', would they exist. Otherwise the convergence performance of these methods is less clear.  
	
	Another way to get around this problem is proposed in \cite{Rodger:2011,Marchetti:2013a}. The idea is not to measure the plant gradient at each iteration but to estimate them with the $n_u+1$ results of the most recent experiments. This reasoning can be pushed further by noting that the measurements of these gradients are used to compute modifiers whose estimation can be assimilated to the resolution of an optimization problem associated with the estimation of the Lagrangian of the plant.  It is on the basis of this observation that \cite{Navia:2015,Blanco:2018} propose to simply replace the repeated plant gradient evaluation by a two-step procedure to be performed at each iteration:
	(Step 1) Choose the new operating point by calculating the minimum of the updated-model. (Step 2) Compute the new modifiers by minimizing the estimated value of the plant Lagrangian that they they imply.

	\subsubsection{Achieving objective 3}
	
	It is clear that being able to converge on the minimum of the plant is an essential property for a theoretical RTO method, but to do so without destroying the plant on the way is also very important.  Ensuring the feasibility of each iteration is a challenge that currently has no satisfactory solution.
	
	The study proposed in \cite{Bunin:2013a,Bunin:2013b} shows that one way to ensure that no experiment violates the constraints of the plant  can be to use  limitations on the size of the iterations based on the plant's Lipschitz constants. Another study \cite{Marchetti:2017} shows that if the constraints of the model are functions that, are at all points, more convex than the functions of the plant, then each iteration produced by MA will be feasible. Of course, in practice neither the plant's Lipschitz constants nor such a model is a priori available.  Naturally, it is possible to replace the Lipschitz constants by overestimates of their values or to use a model so convex that it is undoubtedly more convex than the plant throughout the input space. However, in such cases, applying the methods proposed by \cite{Bunin:2013a,Bunin:2013b} and \cite{Marchetti:2017} would give sequences of iterations that would be feasible but that would converge so slowly on  $\bm{u}_p^{\star}$ that it would not be practical. This means that (in our opinion) these two methods are not really satisfactory in practice. 

	To limit the feasibility issue, it has been proposed to combine the RTO method with a trust region type of constraints \cite{Conn:2000} to keep each consecutive iteration close enough to each other so that they ``secure each other'' \cite{Bunin:2014c,Jonin:2018,Chanona:2019}. Moreover, in practice it is common to add safety margins to the \textit{constraints used} by the RTO method so that small violations do not result in violations of the \textit{real constraints}  \cite{Visser:2000}. Of course, the larger the back-offs are, the less chance there is that the constraints will be violated, the counterpart is a loss of optimality at convergence. \cite{Srinivasan:03b} proposes a method to interactively re-evaluate the back-off used to obtain a better compromise between the risk of violating the constraints and the loss of optimality at convergence. Finally, if the parameters of the model used are associated with an uncertainty of the type $\bm{\theta}\in\mathbb{\Theta}$, where $\mathbb{\Theta}$ is a set of possible parameters. Then, it is possible to consider using robust optimization methods of the type of those discussed in \cite{Srinivasan:02}. 

\subsection{Practical RTO}
	\subsubsection{The practical RTO's problem formulation}
	
	Let's consider the ``un-ideal'' conditions where:
	\begin{itemize}[noitemsep]
		\item The environment can produce disturbances $\bm{d}_p \neq \emptyset$. So the functions of the plant and of the model are: 
		\begin{align*}
			\bm{y}_p = \ & \bm{f}_p(\bm{x}_p), &
			\bm{y}   = \ & \bm{f}(\bm{x},\bm{\theta}).
		\end{align*}
		\item There can be measurement errors: $$\widehat{\bm{y}}_p\neq\bm{y}_p,$$
		therefore, the estimates of $\{\phi_p,\bm{g}_p,\bm{f}_p\}$  and $\{\nabla_{\bm{u}} \phi_p,$ $\nabla_{\bm{u}}\bm{g}_p,\nabla_{\bm{u}}\bm{f}_p\}$ are inaccurate.		
		\item The cost and constraint functions can change over time.
	\end{itemize}
	
	Therefore, the problem of practical RTO is to (\textit{objective 1}) solve the following NLP:
	\begin{align*}
		\bm{u}_{p}^{\star}:=  \operatorname{arg}
		\underset{\bm{u}}{\operatorname{min}} \quad &  \phi_p(\bm{x}_p)   \quad 
		\text{s.t.} \quad  \bm{g}_p(\bm{x})  \leq \bm{0},  
	\end{align*}
	while (\textit{objective 2}) using plant function evaluations only a minimal number of times, and (\textit{objective 3}) avoiding violating the constraints during those evaluations. 
	
	To achieve these three objectives, the idea is generally to apply theoretical RTO methods to which features are added to reduce the effects of measurement errors and disturbances. This combination of theoretical RTO methods with ``additional functionalities'' can be assimilated into the function of the autopilot mentioned above (see Figure~\ref{fig:1___Context}). The following section is dedicated to the ``additional functionalities'' that have been proposed so far.

	\subsubsection{Management of measurement noise and disturbances}
	
	To reduce the sensitivity of the gradient estimation to measurement noise, it has been proposed to base it on a larger number of measurements. For example it is conceivable to use a higher order finite difference method. Based on a similar idea, \cite{Gao:16} proposed to use quadratic approximations to ``combine'' a large number of results from experiments performed in a large enough domain of the input space to filter out the effects of measurement noise. Also, it has been suggested to use non-linear regression methods such as Gaussian process regression \cite{Ferreira:2018} or neural networks \cite{Matias:2019a} to perform such filtering.  
	
	 Disturbances that affect a plant designed to operate at steady-state (SS) can be classified into three categories. Disturbances that are (i) slow, e.g. a degradation, (ii) fast but rare, e.g. a step on the quality of materials used, and (iii) fast and frequent, e.g. a rapid succession of steps.\footnote{If the plant is designed to operate at steady state and is subject to other types of disturbances, it will indeed never reach steady state. Therefore the primary concern should be to improve the design and/or reduce disturbances.} 
	 A disturbance is said to be fast (or slow) if the transition time between two stable values of the disturbance is significantly smaller (or larger) than the largest time constant of the undisturbed plant  $\tau_{max}$.

	\textit{Concerning the disturbances of type (i):} If the plant is only subject to disturbances of type (i), then one way to handle them is to ignore them by considering that any experiment is completed after a fixed duration  $t_{chosen}>\tau_{max}$. After a duration $t_{chosen}$ the measurements made on the plant can then be considered as being steady-state measurements (although this is not the actual case due to the slow disturbances). This approach is proposed and illustrated in \cite{Bunin:2010}.
	
	\textit{Concerning the disturbances of type (ii):} If the plant is subject to (ii) then these are implicitly handled by any iterative RTO method such as ISO, MA, MAy, etc. Indeed, a step on a perturbation during an iteration of one of these methods would only make the information used during this same iteration false. So, at the end of this iteration there is a risk that a wrong decision is taken. However, as these methods do not use memory, the following iterations will be executed correctly and this perturbation will not have any more harmful effect on the decision making.  So, disturbances of type (ii) are to some extent managed. However, if the disturbances are measured and taken into account in the decision making process, then one can avoid these risky iterations and make the decision sequence more relevant as suggested in \cite{Navia:2019}.
	
	\textit{Concerning the disturbances of type  (iii):} Finally, if the plant is subject to  (iii), then it is unlikely to reach a steady state. In this case, all of the methods discussed in this thesis that aim to identify the best possible stable state for the plant may not be suitable. Nevertheless, it is possible to extract principles, ideas and concepts from these methods and apply them to other more appropriate methods. For example, it is possible to combine economic model predictive control (EMPC-\cite{Ellis14}) with MAy as suggested in \cite{Vaccari:2017,Faulwasser:2019,Vaccari:21}.

\subsection{Towards efficient RTO}

	In addition to the theoretical and practical RTO objectives, it is possible to consider performance objectives. One will not attempt to define precisely what the performance of an RTO method is because, as one will see, there are multiple ways of understanding it. 
	
	One can try to reduce the time required to converge on  $\bm{u}_p^{\star}$ by replacing for each iteration the waiting time of the plant stabilization by a shorter fixed waiting time. As a result, it is no longer measurements of the plant's SS that are used to estimate its values and gradients, but transient data. This idea initially applied to TS \cite{Besl:1998,Sequeira:2002,Matias:2018}, has been also applied to MA \cite{Francois:2014,Blanco:2017}. Although it has been shown that using transient measures can accelerate convergence, \cite{Ferreira:2017} shows that such an approach presents risks that can be mitigated by using dynamic models instead of static ones. Indeed, since the measurements are now transient, they can only be pertinently compared to the predictions of a dynamical model.  In addition to these risks, it is clear that such approaches increase the sensitivity of any RTO methods to measurement noise, since the noise filtering is obviously reduced to its minimum. 
	
	Large-scale plant models can generally be represented as networks of submodels. If these submodels are subject to uncertainties, then their effects generally tend to spread to the rest of the network. For example, the inaccurate outputs of an inaccurate submodel can be the inputs of a correct submodel. But since the latter is not evaluated at the correct inputs, its outputs are also incorrect, and so on. This type of error diffusion can be stopped by correcting each sub-model individually as proposed in \cite{Roberts:95} in the context of an ISOPE implementation. In the context of adapting this approach to MA, it has been shown that in addition to reducing the diffusion of errors in the network, the individual correction of each submodel allows a parallelization of the computations that can enable the distribution of the computational load and potentially make some parts of a model private \cite{Milosavljevic:2017,Schneider:17}.
	
	Basic RTO methods use very simple correction functions. Be it ISOPE, ISO, or MA, all use affine correction functions affecting the model in its entire input space. Locally an affine correction based on \textit{local} measures of values and gradient makes sense.  Extending it to the whole input space is questionable. Indeed, if the experimental results give \textit{local} information about the plant, but  the model remains the only \textit{global} data available about it. So, ideally the updated model should only be modified in the vicinity of the experienced points. To achieve this objective, the only proposal that has been made is to use non-linear regression methods such as Gaussian process regression \cite{Ferreira:2018,Chanona:2019} (even though they do not specify the objective they meet which is to use \textit{local} and \textit{global} data for what they are). 
		
	The outputs of the vast majority of RTO methods are set-points for the controllers of the plant. However, there are other options; indeed, \cite{Marlin:97,Engell:07} suggests that the output of an ideal RTO system should be a control structure leading to optimal long-term performance. For the moment, only \cite{Bunin:2013d} provides a RTO method for tuning the parameters of the plant's controllers.   
	
	The majority of RTO methods do not have a stopping mechanism. For example, if MA is applied and the sequence of experiments has converged on $\bm{u}_p^{\star}$, then the gradient of the plant will be constantly re-calculated at $\bm{u}_p^{\star}$. Although these revaluations lead to a form of sub-optimality because one is less often at $\bm{u}_p^{\star}$. They allow any event taking place after the identification $\bm{u}_p^{\star}$ to be managed. However, instead of perpetually perturbing the plant in case something happens, \cite{Ye:2018,Mukkula:2019}  suggest to introduce a stop mechanism when convergence is detected and a restart mechanism when an event occurs. To detect such an event, they suggest using process monitoring theory (e.g. fault detection methods) to detect significant changes in plant behavior and trigger a decision making.
	
\section{Objectives and plan of the thesis}

The objective of this thesis is to build an autopilot allowing the automation of a large number of decisions that the engineers supervising a plant have to make, with an exclusive focus on the search of performances. \\

The remainder of this work is divided into six chapters that can be summarized as follows:  \\

\textbf{Chapter 2:}  In this chapter, a reformulation of a part of the history of the development of RTO methods in which several of our contributions are integrated is proposed. This chapter ends with the proposal of a new theoretical RTO method that offers better capabilities of convergence on optimal decisions than most other methods of the same family.

\textbf{Chapter 3:} Starting from the method proposed in chapter 2, one exposes a weakness common to any RTO method using a filter and one explains how to correct this weakness. The result of this chapter is a new theoretical RTO method combining the efficiency of the one built in chapter 2 with new safety features. (This chapter is an extension of our articles \cite{Papasavvas:2019,Papasavvas:2019b,Papasavvas:2020}.)

\textbf{Chapter 4:} This chapter is dedicated to the introduction to what  the structure of an ideal autopilot could be.  The functions composing this structure as well as the way they interact with each other is discussed. And in order to make the message as clear as possible, one illustrates the functioning of a very simplified version of this structure, which is called \textit{simple} \textit{autopilot for steady-state processes }(S-ASP), on the Williams-Otto reactor, which is a benchmark case study in the RTO community. 

\textbf{Chapter 5:} In this chapter, ones starts from the method proposed in chapter 3 and one tries to push its capabilities further by improving the way the data obtained at the end of each experiment is used. One proposes here to exploit the structure of the model to apply corrections not only to the outputs of the model but directly to variables interconnecting subparts of it. Empirical studies as well as an implementation on the Tennessee Eastman challenge process are made to demonstrate the interest of such approaches. (This chapter is an extension or our articles \cite{Papasavvas:2019a,Papasavvas:2019c})

\textbf{Chapter 6:}  In this chapter, one discusses how it would be possible to exploit the data obtained from the \textit{whole history} of the experiments. Also, one proposes an improved version of the S-ASP that is called \textit{autopilot for steady-state process} (ASP) and one illustrates its functioning on one case study: 
 the Williams-Otto plant. 

\textbf{Chapter 7:} This chapter concludes this thesis and points towards new research directions that may be relevant. 
	\chapter{A RTO algorithm that converges} %
		\label{Chap:2_Vers_Une_meilleure_Convergence}

\section{The problem to solve}
\label{sec:2_1_Definition_PB_ORT_Theorique}

The problem addressed in this chapter is the theoretical RTO problem presented in section~\ref{sec:1_3_1_a_Probleme_OTR_Theorique}. In short, one wants to identify

\begin{align}
	\bm{u}_{p}^{\star}:=  \operatorname{arg}
	\underset{\bm{u}}{\operatorname{min}} \quad &  \phi_p(\bm{u}) := \phi(\bm{u},\bm{f}_p(\bm{u})) \quad 
	\text{s.t.} \quad  \bm{g}_p(\bm{u}) := \bm{g}(\bm{u},\bm{f}_p(\bm{u})) \leq \bm{0},   \label{eq:2_2___1_Plant_PB}
\end{align}
while using only a minimal number of assessments of the values $(\phi_p,\bm{g}_p,\bm{f}_p)$ and gradients $(\nabla_{\bm{u}} \phi_p,\nabla_{\bm{u}}\bm{g}_p,\nabla_{\bm{u}}\bm{f}_p)$. In addition to these measurements, a model is available that provides an estimate, $\bm{f}$, of the plant's mapping, $\bm{f}_p$. 

\section{Using a solver}

\begin{figure}[h]
	\centering
	\includegraphics[width=1\linewidth] {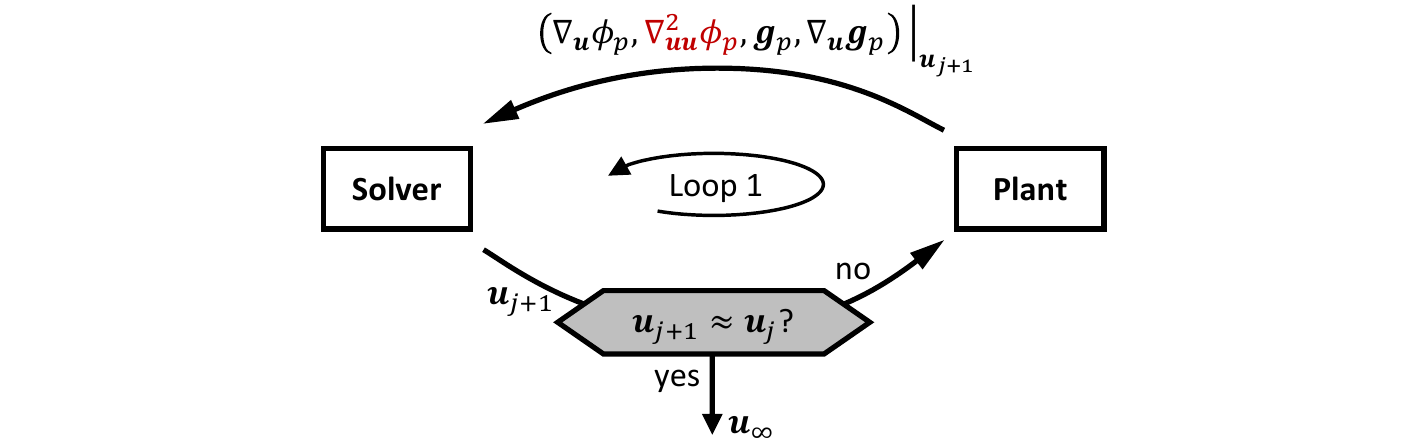}
	\caption{Simplified and schematic operation of a solver}
	\label{fig:2___1_Solver_Plant}
\end{figure}

To solve problem \eqref{eq:2_2___1_Plant_PB}, the simplest option is to connect a standard non-linear optimization problem solver (NLP) directly to the plant. The strategy that these solvers rely on can be summarized as a repetition of evaluations of $\phi$ and $\bm{g}$ in the neighborhood of a point $\bm{u}_j$ allowing the construction of an alternative NLP (approximating \eqref{eq:2_2___1_Plant_PB} around $\bm{u}_j$) and whose solution $\bm{u}_{j+1}$ is ``easy'' to compute. The resulting sequence  $\{\bm{u}_0,...,\bm{u}_{j},\bm{u}_{j+1},...,\bm{u}_{\infty}\}$ is supposed to converge to the desired value, that is $\bm{u}_{\infty} = \bm{u}_p^{\star}$. Let's take for example one of these very popular families of methods: which is sequential quadratic programming (SQP). At each iteration it builds a quadratic approximation of $\phi$ and a linear approximation of $\bm{g}$ which form a quadratic program (QP) whose general writing can be:  
\begin{align}
	\bm{u}_{j+1}:=  \operatorname{arg}
	\underset{\bm{u}}{\operatorname{min}} \quad &  
	\nabla_{\bm{u}}\phi_p|_{\bm{u}_i} (\bm{u}- \bm{u}_i)
	+ \frac{1}{2}(\bm{u}- \bm{u}_i)^{\rm T} 
	\nabla^2_{\bm{uu}}\phi_p|_{\bm{u}_i}
	(\bm{u}- \bm{u}_i)
	\label{eq:2_2___2_QP_PB} \\
	\text{s.t.} \quad & \bm{g}_p|_{\bm{u}_i} + \nabla_{\bm{u}_i}\bm{g}_p|_{\bm{u}_i} (\bm{u}- \bm{u}_i) \leq \bm{0}.   \nonumber
\end{align}
If the active constraints at $\bm{u}_{j+1}$ are known, this problem is relatively simple to solve (see Appendix~\ref{sec:A_2_1_Resoudre_QP_Contraintes_Egalite}).  But since they are usually unknown, they must be identified by following an iterative procedure similar to the one given in Appendix~\ref{sec:A_2_2_Resoudre_QP_Contraintes_Inegalite}. Such procedures may involve constraint violations and usually require a large number of evaluations of the constraint functions (so in our case a large number of experiments to be performed on the plant).

Finally, although such approaches provide very good guarantees on their ability to converge on the plant optimum, they present the following three major problems (if applied directly to an actual system such as an industrial process):
\begin{itemize}
	\item[i.] They use the Hessians of the plant while they are assumed inaccessible.  
	\item[ii.] They almost systematically violate constraints during the identification phase of the active constraints.
	\item[iii.] They require a significantly large number of evaluations of plant values, gradients and Hessians that make the identification of $\bm{u}_p^{\star}$ very slow (when one evaluation means one experiment to be performed on the real process). 
\end{itemize} 

Of course, commercial solvers are generally more complex than the SQP algorithm presented here (and in the appendix). Some can provide better guarantees of non-violation of the constraints during the exploration process they use to identify $\bm{u}_p^{\star}$ (e.g. interior point methods). The counterpart is generally an even slower convergence, i.e. more evaluations of the plant values. 

Without going deeper into the details of these solvers, the main message of this section is that applying a solver directly to a plant is typically not the best idea because of the three major problems stated above, among others. 

\section{Towards the use of models}
\label{sec:2_3_Vers_L_Utilisation_De_Modeles}

\begin{figure}[t]
	\centering
	\includegraphics[width=1\linewidth] {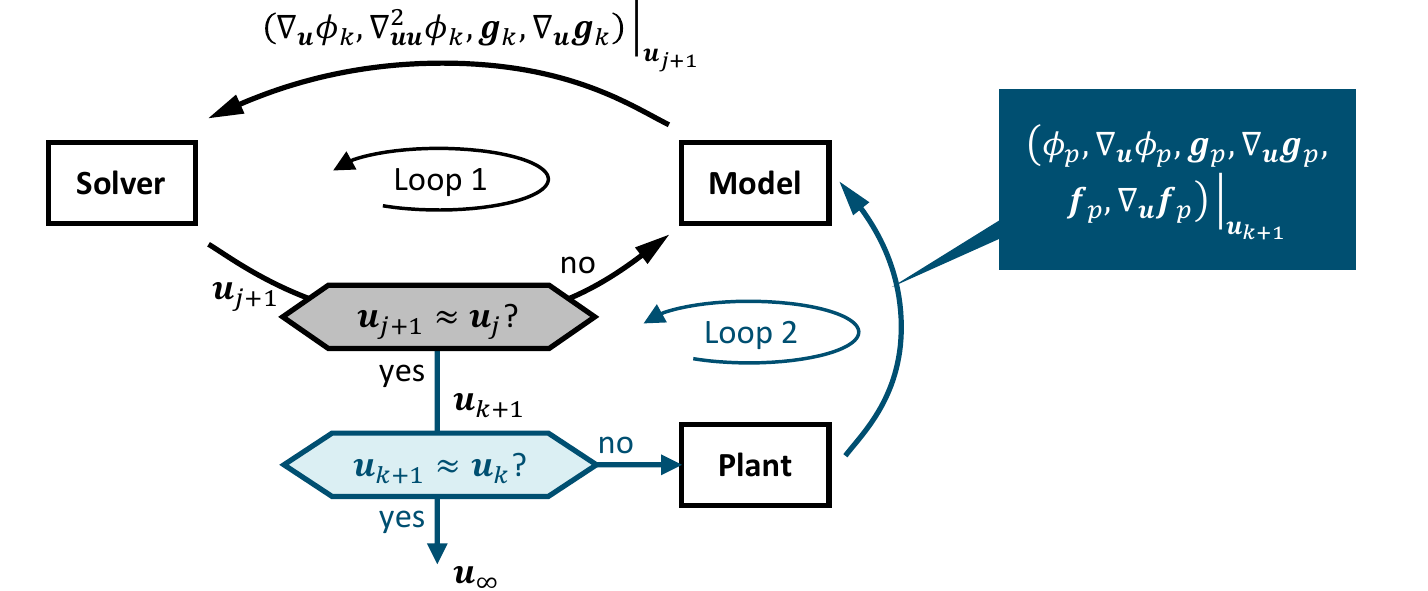}
	\caption{Towards the use of a model to solve the theoretical RTO problem.}
	\label{fig:2___1_Solver_Model_Plant}
\end{figure}

To reduce the effects of these three problems, one can use a model (i) to approximate the plant's Hessians instead of measuring them, (ii) to validate the solver's iterations on the model to reduce the chances of violating the plant constraints, and (iii) to use the global information provided by the model to accelerate the identification of $\bm{u}_p^{\star}$. In this case, the model becomes an intermediate between the solver and the plant and each experiment done on the plant is no longer used by the solver to build simplified representations of the plant, but rather to update the model.  Therefore, solving problem \eqref{eq:2_2___1_Plant_PB} is done through a double iterative mechanism of the type of the one  illustrated in Figure~\ref{fig:2___1_Solver_Model_Plant} where:
\begin{itemize}
	\item Loop 1 corresponds to the optimization of the updated model:
	\begin{align} \bm{u}_{k+1} = \operatorname{arg}
		\underset{\bm{u}}{\operatorname{min}} \quad &  \phi_k(\bm{u},\bm{\theta}_k) 
		 \quad 
		\text{s.t.} \quad \bm{g}_k(\bm{u},\bm{\theta}_k) \leq \bm{0},  \label{eq:2_2___3_Model_PB}
	\end{align}
	where the functions $\phi_k$ and $\bm{g}_k$, and the parameters $\bm{\theta}_k$, are the $k^{\text{th}}$ updates of  $\phi$, $\bm{g}$ and $\bm{\theta}$, respectively.  This mechanism can be an SQP algorithm as presented previously or any other NLP solver. Indeed, as the model is numerical, there is no penalty be it on the number of evaluation of the functions $(\phi_k,\bm{g}_k)$, or on the violation of constraints that are now ``virtual''\footnote{Of course there remains challenges in numerical optimization. But it is not the actually the topic of this thesis. However, if a high-fidelity model that is expensive to evaluate is used, then one recommends to the reader to check section~\ref{sec:7_x_OTR_theotique_et_ROM_Opt}.}. 
	\item Loop 2 corresponds to the update of the model, i.e. the construction of $\phi_{k+1}$, $\bm{g}_{k+1}$ and $\bm{\theta}_{k+1}$, by using the measurements of $\phi_p$, $\nabla_{\bm{u}}\phi_p$, $\bm{g}_p$, and $\nabla_{\bm{u}}\bm{g}_p$ obtained on the plant for each of the experimented operational points $(\bm{u}_0,...,\bm{u}_k)$. 
\end{itemize}
If the mechanism of loop 1 is a ``well known'' standard optimization of a numerical model, then loop 2 is a novelty whose mechanism is to be defined. 

Ideally, the mechanism of loop 2 should be such that $\bm{u}_{\infty} = \bm{u}_{p}^{\star}$ (see Figure~\ref{fig:2___1_Solver_Model_Plant}). It has been shown that: 

\begin{thmbox} \label{thm:2___1_AffineCorrection_Impl_KKTmathcing}
	If at the convergence point $\bm{u}_{\infty}$ of an RTO algorithm the equalities 
	\begin{align}
		\bm{g}_{\infty}|_{\bm{u}_{\infty},\bm{\theta}_{\infty}}  = \ & \bm{g}_p|_{\bm{u}_{\infty}}, &
		\nabla_u \bm{g}_{\infty}|_{\bm{u}_{\infty},\bm{\theta}_{\infty}}  = \ & \nabla_u \bm{g}_p|_{\bm{u}_{\infty}},  & 
		\nabla_u \phi_{\infty}|_{\bm{u}_{\infty},\bm{\theta}_{\infty}}  = \ & \nabla_u \phi_p|_{\bm{u}_{\infty}},
		\label{eq:2___4_NecessaryEqualities}
	\end{align}
	are true, then $\bm{u}_{\infty}$ is a KKT point (Karush–Kuhn–Tucker \cite{Karush:1939,Kuhn:1951}) of the plant, i.e. a minimum (like $\bm{u}_p^{\star}$), a maximum, or a saddle point.
\end{thmbox}
\begin{proofbox} (This is the generalization proposed in \cite{Papasavvas:2019c} of other more specific theorems proposed in \cite{Marchetti:09,Papasavvas:2019a,Milosavljevic:2017}) 
	By definition:
	\begin{align} & \bm{u}_{\infty} := \operatorname{arg}
		\underset{\bm{u}}{\operatorname{min}} \quad   \phi_{\infty}(\bm{u},\bm{\theta}_\infty) 
		\quad 
		\text{s.t.} \quad \bm{g}_\infty(\bm{u},\bm{\theta}_\infty) \leq \bm{0},  \label{eq:2___5_Proof_PB_OPT_Model} \\
		\Rightarrow \ &  \exists \bm{\lambda} \in \amsmathbb{R}^{n_g} \text{ such that } 
		\left\{
		\begin{array}{l}
			\bm{g}_{\infty}|_{\bm{u}_{\infty}, \bm{\theta}_{\infty}} \leq \bm{0},\\
			\bm{\lambda}^{\rm T} \bm{g}_{\infty}|_{\bm{u}_{\infty}, \bm{\theta}_{\infty}} = \bm{0}, \\
			\bm{\lambda} \geq \bm{0}, \\
			\nabla_{\bm{u}}\phi_{\infty}|_{\bm{u}_{\infty},\bm{\theta}_{\infty}} + \bm{\lambda}^{\rm T} \nabla_{\bm{u}}\bm{g}_{\infty}|_{\bm{u}_{\infty},\bm{\theta}_{\infty}} = \bm{0}.
		\end{array}
		\right. \label{eq:2___6_KKT_Model}
	\end{align}
	Indeed, as $\bm{u}_{\infty}$ is a solution of \eqref{eq:2___5_Proof_PB_OPT_Model}, it is a KKT point of the same problem, hence the implication of \eqref{eq:2___6_KKT_Model}. Because of the equalities \eqref{eq:2___4_NecessaryEqualities}, one can replace $(\bm{g}_{\infty}, \nabla_{\bm{u}}\bm{g}_{\infty}, \nabla_{\bm{u}}\phi_{\infty})|_{\bm{u}_{\infty},\bm{\theta}_{\infty}}$ by $(\bm{g}_{p}, \nabla_{\bm{u}}\bm{g}_{p}, \nabla_{\bm{u}}\phi_{p})|_{\bm{u}_{\infty}}$ and observe that $\bm{u}_{\infty}$ is also a KKT point of the plant. In brief:
	\begin{align} 
		\text{\eqref{eq:2___4_NecessaryEqualities}  \& \eqref{eq:2___6_KKT_Model}}  \Rightarrow \ &   
		\left\{
		\begin{array}{l}
			\bm{g}_{p}|_{\bm{u}_{\infty}} \leq \bm{0},\\
			\bm{\lambda}^{\rm T} \bm{g}_{p}|_{\bm{u}_{\infty}} = \bm{0}, \\
			\bm{\lambda} \geq \bm{0}, \\
			\nabla_{\bm{u}}\phi_{p}|_{\bm{u}_{\infty}} + \bm{\lambda}^{\rm T} \nabla_{\bm{u}}\bm{g}_{p}|_{\bm{u}_{\infty}} = \bm{0}.
		\end{array}
		\right. \label{eq:2___7_Proof_End}
	\end{align}
\end{proofbox}
It is then clear that for an RTO method to guarantee the optimality of the plant upon convergence, it is necessary that the equalities \eqref{eq:2___4_NecessaryEqualities} are true. These equalities can be enforced in two different ways which characterize two classes of approaches, direct and indirect approaches. 
	
\textit{Direct approaches} consist of applying correction functions directly to $\phi$ and $\bm{g}$. More precisely, these functions must be corrected as follows:
\begin{subequations}
	\begin{align}
	\phi_{\infty}(\bm{u},\bm{\theta}_{\infty}) := \ & \phi(\bm{u},\bm{f}(\bm{u},\bm{\theta}_{\infty})) + \mu_{\infty}^{\phi}(\bm{u},\bm{\theta}_{\infty})  \\
	\bm{g}_{\infty}(\bm{u},\bm{\theta}_{\infty}) := \ & \bm{g}(\bm{u},\bm{f}(\bm{u},\bm{\theta}_{\infty})) +\bm{\mu}_{\infty}^{\bm{g}}(\bm{u},\bm{\theta}_{\infty}), 
	\end{align}
	where $(\mu_{\infty}^{\phi},\bm{\mu}_{\infty}^{\bm{g}})$ are functions whose Taylor series expansions around $\bm{u}_{\infty}$ are:
	\begin{align}
		\mu_{\infty}^{\phi}(\bm{u},\bm{\theta}_{\infty}) := \ & \nabla_{\bm{u}} (\phi_p - \phi)|_{\bm{u}_{\infty},\bm{\theta}_{\infty}}^{\rm T} (\bm{u}-\bm{u}_{\infty}) + \text{ h.o.t.}, \\
		\bm{\mu}_{\infty}^{\bm{g}}(\bm{u},\bm{\theta}_{\infty}) := \ & 
		( \bm{g}_p- \bm{g})|_{\bm{u}_{\infty},\bm{\theta}_{\infty}} + 
		\nabla_{\bm{u}} (\bm{g}_p -  \bm{g})|_{\bm{u}_{\infty},\bm{\theta}_{\infty}} (\bm{u}-\bm{u}_{\infty})  + \text{ h.o.t.}.
	\end{align}
	so that the equalities \eqref{eq:2___4_NecessaryEqualities} are enforced. 
\end{subequations}

\textit{Indirect approaches} consist of applying correction functions to $\bm{f}$ to indirectly correct the functions $\phi$ and $\bm{g}$. To do this, the functions $\bm{f}$ must be corrected as follows:
\begin{subequations}
	\begin{align}
		\bm{f}_{\infty}(\bm{u},\bm{\theta}_{\infty}) := \ &  \bm{f}(\bm{u},\bm{\theta}_{\infty}) + \bm{\mu}_{\infty}^{\bm{f}}(\bm{u},\bm{\theta}_{\infty}),
	\end{align}
	where 
	\begin{align}
		\bm{\mu}_{\infty}^{\bm{f}}(\bm{u},\bm{\theta}_{\infty}) := \ & 
		( \bm{f}_p- \bm{f})|_{\bm{u}_{\infty},\bm{\theta}_{\infty}} + 
		\nabla_{\bm{u}} (\bm{f}_p -  \bm{f})|_{\bm{u}_{\infty},\bm{\theta}_{\infty}} (\bm{u}-\bm{u}_{\infty})  + \text{ h.o.t.},
		\label{eq:2___15_XXXX}
	\end{align}
	The functions $(\phi,\bm{g})$ are not modified, i.e.
	\begin{align}
		\phi_{\infty}(\bm{u},\bm{\theta}_{\infty}) := \ & \phi(\bm{u},\bm{f}_{\infty}(\bm{u},\bm{\theta}_{\infty})) \\
		\bm{g}_{\infty}(\bm{u},\bm{\theta}_{\infty}) := \ & \bm{g}(\bm{u},\bm{f}_{\infty}(\bm{u},\bm{\theta}_{\infty})). 
	\end{align}
	The fact that \eqref{eq:2___15_XXXX} implies \eqref{eq:2___4_NecessaryEqualities} is not necessarily obvious on the face of it. The mathematical proof of this implication is given in the Appendix~\ref{App:A___CorrectionFonctionComposite} and is inspired of \cite{Papasavvas:2019a}.
\end{subequations}

These two types of approaches are very important because they separate \textit{the whole set of theoretical RTO methods that guarantee the optimal operation of the plant upon convergence} into two distinct families. Now that these two classes of approaches are clearly identified, the simplest possible versions can be built. 

\section{Two very simple methods (ISO-D/I)}
\label{sec:2_4_ISO_ISOy}

To obtain the two ``simplest'' direct and indirect theoretical RTO methods, it is sufficient to make the following four simplifications:
\begin{itemize}
	\item \textbf{(Simplification 1)}  Leave the parameters at predetermined values called nominal values
	\begin{align}
		\bm{\theta}_k := \ & \bm{\theta}_n, & \forall k \in \mathbb{Z},
	\end{align}
	to avoid having to solve parameter identification problems. \textit{Note that from this point on, if it is not specified for which parameter a function is evaluated, it means that the values that are used are the nominal ones, $\bm{\theta}_n$.}
	\begin{align*}
		\bm{f}(\bm{u},\bm{\theta}_n)     \equiv \ & \bm{f}(\bm{u}),  &
		(\cdot)|_{\bm{u}_k,\bm{\theta}_n} \equiv \ & (\cdot)|_{\bm{u}_k}. 
	\end{align*}
	\item \textbf{(Simplification 2)} Enforce the satisfaction of the following equalities at each iteration $k$: 
	\begin{align}
		\bm{g}_{k}|_{\bm{u}_{k}}  = \ & \bm{g}_p|_{\bm{u}_{k}}, &
		\nabla_u \bm{g}_{k}|_{\bm{u}_{k}}  = \ & \nabla_u \bm{g}_p|_{\bm{u}_{k}},  & 
		\nabla_u \phi_{k}|_{\bm{u}_{k}}  = \ & \nabla_u \phi_p|_{\bm{u}_{k}}, & \forall k \in \mathbb{Z}.
		\label{eq:2___11_SimpleApproach_SameCones}
	\end{align}
	to make sure that they are true when $k\rightarrow \infty$ and thus satisfy \eqref{eq:2___4_NecessaryEqualities}.
	\item \textbf{(Simplification 3)} Use the simplest correction functions allowing to obtain the equalities \eqref{eq:2___4_NecessaryEqualities}. In other words, the functions $\mu_k^{\phi}$ and $\bm{\mu}_k^{\bm{g}}$, or $\bm{\mu}_k^{\bm{f}}$, are affine, i.e. their higher order terms $(h.o.t.)$ are fixed to $0$.
	\item \textbf{(Simplification 4)} Use only the measures obtained at $\bm{u}_k$ ($\phi_p,$ $ \nabla_{\bm{u}}\phi_p,\bm{g}_p,\nabla_{\bm{u}}\bm{g}_p$). Those obtained previously at $(\bm{u}_0,...,\bm{u}_{k-1})$ are ignored so that there is no need to manage any database.  
\end{itemize}
This is how Algorithm~\ref{algo:ISO}, which corresponds to the ``direct'' version of iterative setpoint optimization method (ISO-D -- \cite{Alexandrov:98,Gao:05}), can be reconstructed. At the same time the indirect version of this method can be obtained, which is called indirect iterative setpoint optimization (ISO-I) and whose details are given in Algorithm~\ref{algo:ISOy}.\\

\hypersetup{citecolor=white}
\begin{BoxAlgo}{\textbf{Direct Iterative Setpoint optimization (ISO-D -- \cite{Alexandrov:98,Gao:05})}}{ISO} 
	\textbf{Initialization.} Provide $\bm{u}_0$, functions $(\bm{f},\phi,\bm{g})$, and the stoping criteron of step 4).
	\tcblower
		\textbf{for} $k=0 \rightarrow \infty$
		\begin{itemize}[noitemsep]
			\item[1) ]  \textbf{Measure} $(\nabla_{\bm{u}}\phi_p, \bm{g}_p,\nabla_{\bm{u}}\bm{g}_p)|_{\bm{u}_{k}}$ on the plant.
			\item[2) ] \textbf{Construct} the functions $(\phi_k,\bm{g}_k)$: 
			\begin{subequations} \label{eq:3___19_Corrections_ISO}
				\begin{align}
					\phi_{k}(\bm{u}) := \ & \phi(\bm{u},\bm{f}(\bm{u})) + \mu_{k}^{\phi}(\bm{u})  \\
					\bm{g}_{k}(\bm{u}) := \ & \bm{g}(\bm{u},\bm{f}(\bm{u})) +\bm{\mu}_{k}^{\bm{g}}(\bm{u}), 
				\end{align}
				where 
				\begin{align}
					\mu_{k}^{\phi}(\bm{u}) := \ & \nabla_{\bm{u}} (\phi_p - \phi)|_{\bm{u}_{k}}^{\rm T} (\bm{u}-\bm{u}_{k}), \\
					\bm{\mu}_{k}^{\bm{g}}(\bm{u}) := \ & 
					( \bm{g}_p- \bm{g})|_{\bm{u}_k} + 
					\nabla_{\bm{u}} (\bm{g}_p -  \bm{g})|_{\bm{u}_{k}} (\bm{u}-\bm{u}_{k}),
				\end{align}
			\end{subequations}
			\item[3) ]  \textbf{Solve} the model-based optimization problem \eqref{eq:2_2___3_Model_PB} using the definitions \eqref{eq:3___19_Corrections_ISO} to find $\bm{u}_{k+1}$.
			\item[4) ] \textbf{Stop} if $\bm{u}_{k+1}\approx\bm{u}_{k}$ and return $\bm{u}_{\infty} := \bm{u}_{k+1}$.
		\end{itemize}
		\noindent {\bf end}
\end{BoxAlgo}
\hypersetup{citecolor=green!50!black}	

\begin{BoxAlgo}{\textbf{Indirect Iterative Setpoint optimization (ISO-I)}}{ISOy}
	\textbf{Initialization.} Provide $\bm{u}_0$, functions $(\bm{f},\phi,\bm{g})$, and the stoping criteron of step 4).
	\tcblower
	\textbf{for} $k=0 \rightarrow \infty$
	\begin{itemize}[noitemsep]
		\item[1) ]  \textbf{Measure} $(\bm{f}_p,\nabla_{\bm{u}}\bm{f}_p)|_{\bm{u}_{k}}$ on the plant.
		\item[2) ] \textbf{Construct} the functions $(\bm{f}_k,\phi_k,\bm{g}_k)$: 
		\begin{subequations} \label{eq:3___20_Corrections_ISOy}
			\begin{align}
				& \qquad \qquad \qquad \bm{f}_{k}(\bm{u}) :=    \bm{f}(\bm{u}) +\bm{\mu}_{k}^{\bm{f}}(\bm{u}), \\
				& \phi_k(\bm{u}) :=   \phi(\bm{u},\bm{f}_k(\bm{u})), \qquad 
				\bm{g}_k(\bm{u}) :=  \bm{g}(\bm{u},\bm{f}_k(\bm{u})), 
			\end{align}
			where 
			\begin{align}
				\bm{\mu}_{k}^{\bm{f}}(\bm{u}) := \ & 
				( \bm{f}_p- \bm{f})|_{\bm{u}_k} + 
				\nabla_{\bm{u}} (\bm{f}_p -  \bm{f})|_{\bm{u}_{k}} (\bm{u}-\bm{u}_{k}).
			\end{align}
		\end{subequations}
		\item[3) ]  \textbf{Solve} the model-based optimization problem \eqref{eq:2_2___3_Model_PB} using the definition \eqref{eq:3___20_Corrections_ISOy} to find $\bm{u}_{k+1}$.
		\item[4) ] \textbf{Stop} if $\bm{u}_{k+1}\approx\bm{u}_{k}$ and return $\bm{u}_{\infty} := \bm{u}_{k+1}$.
	\end{itemize}
	\noindent {\bf end}
\end{BoxAlgo}

These two methods are particularly interesting because they are the simplest to understand and to implement, and they are at the basis of many other more complex RTO methods. To analyze the capacities of these methods to converge on  $\bm{u}_p^{\star}$, it is necessary to define clearly what is meant by ``capacity of convergence''.  
\section{Three levels of convergence}

Starting off by introducing the concept of asymptotic convergence, a sequence of points $\{\bm{u}_{\ell},\bm{u}_{\ell+1}, ... \bm{u}_{\infty}\}$ converges asymptotically on $\bm{u}_{p}^{\star}$ if, and only if, each of these points progressively gets closer to $\bm{u}_{p}^{\star}$. In more mathematical terms, this sequence can be said to converge asymptotically on $\bm{u}_{p}^{\star}$  starting from an iteration $\ell$ if: 
\begin{align} \label{eq:2___15_DefConvergenceAsymptotique}
	& \forall k\geq \ell: \quad \underset{k\rightarrow\infty}{\operatorname{lim}} \bm{u}_{k} = \bm{u}_{p}^{\star},
\ \text{ and } \ \left\{\begin{array}{ll}
		\| \bm{u}_{k}  - \bm{u}_p^{\star}  \| < \| \bm{u}_{k+1}  - \bm{u}_p^{\star}  \|, & \text{if } \bm{u}_{k} \neq \bm{u}_p^{\star}, \\
		\| \bm{u}_{k}  - \bm{u}_p^{\star}  \| = \| \bm{u}_{k+1}  - \bm{u}_p^{\star}  \| = 0, & \text{if } \bm{u}_{k} = \bm{u}_p^{\star}.
	\end{array}\right. 
\end{align}

\begin{figure}[htp]
	\centering
	\includegraphics[width=1\linewidth] {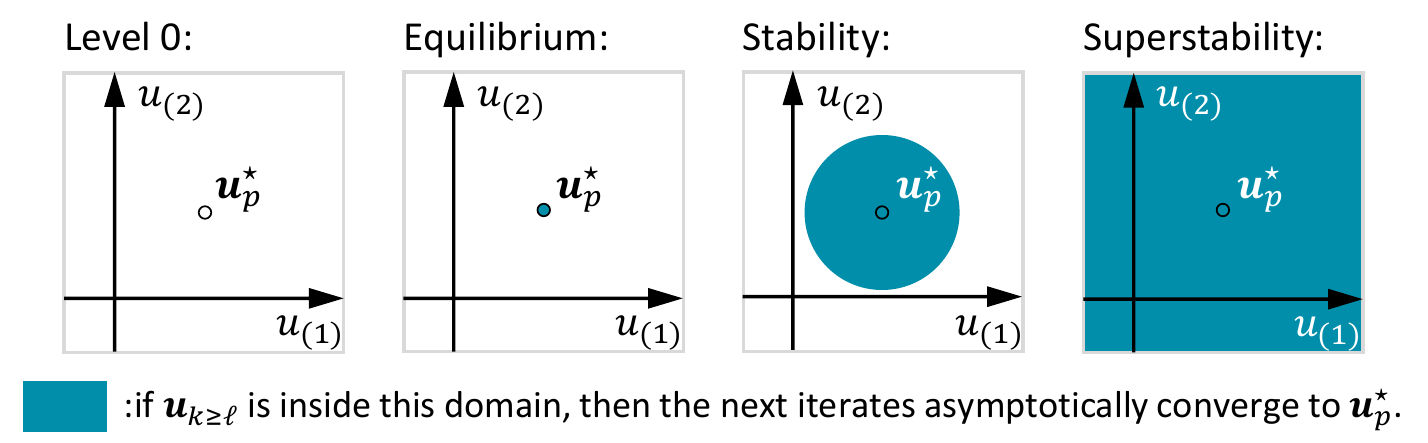}
	\caption{Schematic definition of the 4 levels of convergence.}
	\label{fig:2___2_Domaines}
\end{figure}

Ideally, independently of their initialization $\bm{u}_0$ and after a limited number of iterates, ISO-D/I would enter into such a convergent sequence, i.e.
\begin{equation} \label{eq:2___16_CeQuOnVeut}
	\forall \bm{u}_0\in\amsmathbb{R}^{n_u}, \ \exists \ell\in\mathbb{Z} \quad \text{ such that \eqref{eq:2___15_DefConvergenceAsymptotique}.}
\end{equation} 
Instead of directly analyzing the ability of ISO-D/I to provide \eqref{eq:2___16_CeQuOnVeut}, this condition can be divided into four levels of convergence whose definitions are:
\begin{itemize}[noitemsep]
	\item \textbf{Level 0:} A theoretical RTO method provides level 0 of convergence if no matter what $\bm{u}_0$ is, it is \textit{unable} to guarantee asymptotic convergence of the iterates on $\bm{u}_p^{\star}$.
	\item \textbf{Level 1 -- Equilibrium condition:} A theoretical RTO method satisfies the equilibrium condition if  $\bm{u}_{\ell}=\bm{u}_p^{\star}$ $\Rightarrow$ $\bm{u}_{k\geq\ell}=\bm{u}_p^{\star}$. 
	\item \textbf{Level 2 -- Stability Condition:} A theoretical RTO method satisfies the stability condition if $\exists r\in\amsmathbb{R}>0$ such that if $\bm{u}_{\ell}\in\mathcal{B}(\bm{u}_p^{\star},r)$ then $\bm{u}_{k\geq\ell}$ asymptotically converges on $\bm{u}_p^{\star}$, where $\mathcal{B}(\bm{u}_p^{\star},r)$ is a ball of center $\bm{u}_p^{\star}$ and radius $r>0$. 
	\item \textbf{Level 3 -- Superstability Condition:} A theoretical RTO method satisfies the superstability condition if $\bm{u}_{\ell}\in \amsmathbb{R}^{n_u}$ $\Rightarrow$ $\bm{u}_{k\geq\ell}$ asymptotically converge on $\bm{u}_p^{\star}$.
\end{itemize}
Figure~\ref{fig:2___2_Domaines} provides illustrations of those definitions. Based on these representations, the condition  \eqref{eq:2___16_CeQuOnVeut} for a given $\bm{u}_0$  can be reduced to the expectation that one of the iterations of the sequence $\{\bm{u}_0, \bm{u}_1,...\}$ produced by a theoretical RTO method falls in the blue area (see Figure~\ref{fig:2___2_Domaines}) associated with this method.
So, the larger this domain, i.e. the higher the level of convergence of this method is, the higher the chances to fall in the blue domain and thus to converge on $\bm{u}_{p}^\star$.

With these levels of convergence defined, the question arises of to what level does an ISO-D/I converge and what can be done to raise the level of convergence.

\section{Enforcing the equilibrium condition}
Starting off, determining whether an ISO-D/I provides level 1 convergence; if it does not, determining how it should be modified to provide it. Firstly, analyzing under which conditions ISO-D/I guarantee the satisfaction of the equilibrium condition, the following observation can be made:

\begin{thmbox} \label{thm:2___2_ConditionNiveau1}
	If the RTO method used guarantees the equalities \eqref{eq:2___11_SimpleApproach_SameCones} and if during an iteration $k^{\star}\in\amsmathbb{Z}$, $\bm{u}_{k^{\star}}= \bm{u}_p^{\star}$. Then:
	\begin{itemize}
		\item[i)] The fact that $\bm{u}_p^{\star}$ is a KKT point of the plant makes that it also a KKT point of the model updated at $\bm{u}_p^{\star}$. 
		\item[ii)] The condition for $\bm{u}_{k\geq k^{\star}}= \bm{u}_p^{\star}$ is that:
		\begin{equation}
			\label{eq:2___18_ModelAdequacy_ISO_ISOy}
			\nabla^2_{\bm{u}\bm{u}}\phi_{k^\star}|_{\bm{u}_{p}^{\star}}  + 
			\sum_{i=1}^{n_g}\left[
			\lambda_{p(i)}  \nabla^2_{\bm{u}\bm{u}}g_{(i)k^\star}|_{\bm{u}_{p}^{\star}} \right]  >  0,
		\end{equation}
		where the $\bm{\lambda}_p$ are the KKT-multipliers of the plant.
	\end{itemize}
\end{thmbox}
\begin{proofbox}
	By using the definition of $\bm{u}_p^{\star}$ and  \eqref{eq:2___11_SimpleApproach_SameCones}, one can easily show that $\bm{u}_p^{\star}$ is a KKT point of the model updated at $\bm{u}_p^{\star}$: 
	\begin{align} 
		& \bm{u}_p^{\star} := \operatorname{arg}
		\underset{\bm{u}}{\operatorname{min}} \quad   \phi_{p}(\bm{u}) 
		\quad 
		\text{s.t.} \quad \bm{g}_p(\bm{u}) \leq \bm{0},  \label{eq:2___wegrw} \\
		\text{\eqref{eq:2___wegrw} }  \Rightarrow \ &  \exists \bm{\lambda}_p \in \amsmathbb{R}^{n_g} \text{ such that } 
		\left\{
		\begin{array}{l}
			\bm{g}_{p}|_{\bm{u}_p^{\star}} \leq \bm{0},\\
			\bm{\lambda}_p^{\rm T} \bm{g}_{p}|_{\bm{u}_p^{\star}} = \bm{0}, \\
			\bm{\lambda}_p \geq \bm{0}, \\
			\nabla_{\bm{u}}\phi_{p}|_{\bm{u}_p^{\star}} + \bm{\lambda}_p^{\rm T} \nabla_{\bm{u}}\bm{g}_{p}|_{\bm{u}_p^{\star}} = \bm{0}.
		\end{array}
		\right. \label{eq:2___16_Proof_KKT1_Plant} \\
		\text{\eqref{eq:2___11_SimpleApproach_SameCones} \& \eqref{eq:2___16_Proof_KKT1_Plant} }  \Rightarrow \ & 
		\left\{
		\begin{array}{l}
			\bm{g}_{k^\star}|_{\bm{u}_p^{\star}} \leq \bm{0},\\
			\bm{\lambda}_p^{\rm T} \bm{g}_{k^\star}|_{\bm{u}_p^{\star}} = \bm{0}, \\
			\bm{\lambda}_p \geq \bm{0}, \\
			\nabla_{\bm{u}}\phi_{k^\star}|_{\bm{u}_p^{\star}} + \bm{\lambda}_p^{\rm T} \nabla_{\bm{u}}\bm{g}_{k^\star}|_{\bm{u}_p^{\star}} = \bm{0}.
		\end{array}
		\right. \label{eq:2___16_Proof_KKT1_Modle}
	\end{align}
	So, $\bm{u}_p^{\star}$ is a KKT point of the model updated at this point, and $\bm{\lambda}_p$ are the KKT-multipliers associated to $\bm{u}_p^{\star}$. So, $\bm{u}_p^{\star}$ can be either a minimum, a maximum or a saddle point of the model. In order for it to be a minimum, i.e. $\bm{u}_{k^\star} = \bm{u}_p^{\star} \ \Rightarrow \ \bm{u}_{k^\star+1} = \bm{u}_p^{\star}$, the second order KKT condition must be satisfied, which is the case when \eqref{eq:2___18_ModelAdequacy_ISO_ISOy} is true.
\end{proofbox}

\begin{Remark}
	The model is said to be adequate for ISO-D/I if \eqref{eq:2___18_ModelAdequacy_ISO_ISOy} holds without any intervention. The concept of model adequacy is proposed in \cite{Forbes:94a} where it is applied to TS to demonstrate its weaknesses. It is reused in \cite{Marchetti:09} where it is applied to an improved version of ISO to demonstrate the conceptual advantages of ISO methods over TS. 
\end{Remark}

It is clear that it is a priori impossible to verify whether   \eqref{eq:2___18_ModelAdequacy_ISO_ISOy} is true, because  neither $\bm{u}_p^{\star}$ nor $\bm{\lambda}_p$ are known. However, what one can do is to make this inequality true at all points $\bm{u}\in\amsmathbb{R}^{n_u}$ (instead of only at $\bm{u}_p^{\star}$), and for all $\bm{\lambda}^{\prime}\in\amsmathbb{R}^{n_g+}$ (instead of only for $\bm{\lambda}_p$). To do this, the idea is to replace the functions $\phi_k$ and $\bm{g}_k$ by  approximations that are convex at $\bm{u}_k$  and that are noted $\phi^c_k$ and $\bm{g}^c_k$. These new functions must, $\forall k \in \amsmathbb{Z}$  and $i=1,...,n_g$, be such that:
\begin{align}
	\label{eq:3___21}
	\nabla_{\bm{u}} \phi^c_{k}|_{\bm{u}_k} = \ &  \nabla_{\bm{u}} \phi_{k}|_{\bm{u}_k}, & 
	\bm{g}^c_{k}|_{\bm{u}_k} = \ &  \bm{g}_{k}|_{\bm{u}_k}, & 
	\nabla_{\bm{u}} \bm{g}^c_{k}|_{\bm{u}_k} = \ &  \nabla_{\bm{u}} \bm{g}_{k}|_{\bm{u}_k}, 
\end{align}
\begin{align}
	\label{eq:3___22}
	\nabla^2_{\bm{u}\bm{u}} \phi^c_{k}|_{\bm{u}_k} > \ & 0,  & 
	\nabla^2_{\bm{u}\bm{u}} g^c_{(i)k}|_{\bm{u}_k} > \ & 0,
\end{align}
On one hand, equalities \eqref{eq:3___21} guarantee that Theorems~\ref{thm:2___1_AffineCorrection_Impl_KKTmathcing} and \ref{thm:2___2_ConditionNiveau1} are applicable. 
On the other hand, inequalities  \eqref{eq:3___22} imply that:
\begin{align} \label{eq:3___23}
	\nabla^2_{\bm{u}\bm{u}}\phi^c_{k}|_{\bm{u}_{k}}  + 
	\sum_{i=1}^{n_g}\left[
	\lambda^{\prime}_{(i)}  \nabla^2_{\bm{u}\bm{u}}g^c_{(i)k}|_{\bm{u}_{k}} \right]  > \ & 0, & 
	\forall \bm{\lambda}^{\prime}\in\amsmathbb{R}^{n_g+}, \ \forall \bm{u}_k\in\mathbb{R}^{n_u},
\end{align}
and therefore, by applying Theorem~\ref{thm:2___2_ConditionNiveau1}, one can conclude that  the equilibrium condition is provided.

Now the question that remains is: How to build $(\phi^c_{k},$ $\bm{g}^c_{k})$? It is clear that an infinite number of approaches can be considered. Below, two fairly simple options are proposed: 

\subsection{Option 1: With model preprocessing}

The functions $\phi^c_{k},$ and $\bm{g}^c_{k}$ can be obtained by preprocessing of the model, whereby it is made convex at any point. For ISO-I, this option is not easily applicable as it induces high order corrections (see Appendix~\ref{App:A___CorrectionFonctionComposite} for the details). For ISO-D, this option consists in choosing $\phi^c$ and $g_{(i)}^c$, $\forall i = 1,...,n_g$, in the following way:
\begin{subequations} \label{eq:2___24_ConvexPB}
	\begin{align}
		\label{eq:2___24a_Coutc}
		\phi^c(\bm{u}) :\approx \
		& \phi(\bm{u},\bm{f}(\bm{u},\bm{\theta}_n)), & & 
		\text{such that $\forall \bm{u}\in\amsmathbb{R}^{n_u}$: } \nabla^2_{\bm{u}\bm{u}} \phi^c(\bm{u}) > 0,  \\
		\label{eq:2___24b_Contraintesc}
		g_{(i)}^c(\bm{u}) :\approx \ & g_{(i)}(\bm{u},\bm{f}(\bm{u},\bm{\theta}_n)),  & &
		\text{such that $\forall \bm{u}\in\amsmathbb{R}^{n_u}$: } \nabla^2_{\bm{u}\bm{u}} g_{(i)}^c(\bm{u}) \geq 0.
	\end{align}
\end{subequations}
Since ISO-D applies affine corrections directly to the cost and constraint functions, their initial curvatures are preserved, so  $\forall i = 1,...,n_g$:
\begin{align*}
	\nabla^2_{\bm{u}\bm{u}} \phi^c(\bm{u}) > \ & 0, \forall \bm{u}\in\amsmathbb{R}^{n_u}& &
	\Rightarrow \quad \nabla^2_{\bm{u}\bm{u}} \phi^c_{k}(\bm{u}) >   0, \ \forall (\bm{u},\bm{u}_k)\in\amsmathbb{R}^{n_u} \\
	\nabla^2_{\bm{u}\bm{u}} g_{(i)}^c(\bm{u}) \geq \ & 0, \forall \bm{u}\in\amsmathbb{R}^{n_u}& &
	\Rightarrow \quad \nabla^2_{\bm{u}\bm{u}} g^c_{(i)k}(\bm{u}) \geq  0, \ \forall (\bm{u},\bm{u}_k)\in\amsmathbb{R}^{n_u}.
\end{align*}
Therefore, convex approximations of type \eqref{eq:2___24_ConvexPB} ensure that the equilibrium condition  \eqref{eq:2___18_ModelAdequacy_ISO_ISOy} is satisfied. 

A special case of this strategy is proposed in \cite{Francois:2013} where the cost function is approximated with \eqref{eq:2___24a_Coutc} and where the constraints are simply linearized (which is a special case of \eqref{eq:2___24b_Contraintesc}).

\subsection{Option 2: With successive convexifications}

The functions $\phi^c_{k}$ and $\bm{g}^c_{k}$ can be obtained by applying, at each iteration, a local convexification of the updated functions $\phi_{k}$ and $\bm{g}_{k}$. The idea is to proceed as follows: 
\begin{subequations}\label{eq:3___30_DefinitionMatrices_P}
	\begin{align}
		\phi^c_{k}(\bm{u}) := \ &  \phi_{k}(\bm{u},\bm{f}_{k}(\bm{u})) + \frac{1}{2}
		(\bm{u}-\bm{u}_k)^{\rm T}\bm{Q}_k^{\phi} (\bm{u}-\bm{u}_k), \\
		g^c_{(i)k}(\bm{u}) := \ &  g_{(i)k}(\bm{u},\bm{f}_{k}(\bm{u})) + \frac{1}{2}
		(\bm{u}-\bm{u}_k)^{\rm T}\bm{Q}_k^{g_{(i)}}(\bm{u}-\bm{u}_k),
	\end{align}
	where  	
	\begin{align} 
		\bm{Q}_k^{\phi} :\approx \ &\nabla^2_{uu}\phi_{k}(\bm{u},\bm{f}_{k}(\bm{u})) &&\text{ such that }  \bm{Q}_k^{\phi} > 0, \\
		\bm{Q}_k^{g_{(i)}} :\approx  \ & \nabla^2_{uu}g_{(i)k}(\bm{u},\bm{f}_{k}(\bm{u}))&& \text{ such that }  \bm{Q}_k^{g_{(i)}}\geq 0, \quad \forall i=1,...,n_g.
	\end{align}
\end{subequations}
This strategy taken from \cite{Papasavvas:2020} (section 4.2 of \cite{Papasavvas:2020}) works for both ISO-D and -I. \\

It has been shown that level 1 can be enforced by using the options 1 or 2. Next the stability condition is analyzed.

\section{Enforcing the stability condition}

\subsection{A property of the iterative RTO methods}

The functions $\phi_k$ and $\bm{g}_k$ of the model updated at the iteration $k$ are entirely defined by the nominal functions $\phi$ and $\bm{g}$ and all the results of the experiments performed around the points $\{\bm{u}_0,\bm{u}_1,...,\bm{u}_k\}$. So, in a way, one can consider that the functions   $\phi_k$ and $\bm{g}_k$ are simplified versions of more general functions $\widetilde{g}_{(0)}$ and $\widetilde{\bm{g}}_k$ which can be defined as:
\begin{align} 
	\widetilde{g}_{(0)}(\bm{u},\bm{v}_k,...,\bm{v}_0)|_{\bm{v}_k= \bm{u}_k, ..., \bm{v}_0 = \bm{u}_0} := \ & \phi_k(\bm{u}), \nonumber \\ 
	\widetilde{\bm{g}}(\bm{u},\bm{v}_k,...,\bm{v}_0)|_{\bm{v}_k= \bm{u}_k, ..., \bm{v}_0 = \bm{u}_0} := \ & \bm{g}_k(\bm{u}). \label{eq:2_Generalized_version_function}
\end{align} 
Consider the particular case of purely iterative methods:
\begin{Definition} \textbf{(Purely iterative methods)} \textit{
		A RTO method is said to be \emph{purely iterative}  if it exploits only the results obtained during the last series of experiments, i.e. if: 
		\begin{align} \label{eq:2___25_DefMethodePurementIteratives}
			\phi_k(\bm{u}) = \ &  \widetilde{g}_{(0)}(\bm{u},\bm{v}_k), &
			\bm{g}_k(\bm{u}) = \ &  \widetilde{\bm{g}}(\bm{u},\bm{v}_k).
	\end{align} }
\end{Definition}

\begin{Remark}
	The Simplification 4 used to build ISO-D/I implies that they are purely iterative.   
\end{Remark}
In the case of purely iterative methods, the solution of updated problem \eqref{eq:2_2___3_Model_PB} is a function $\bm{sol}$ of the point $\bm{v}$ around which the latest set of experiments have been conducted. In more mathematical terms, \eqref{eq:2_2___3_Model_PB} becomes: 
\begin{align} 
	\bm{sol}(\bm{v}) := \  & \operatorname{arg}
	\underset{\bm{u}}{\operatorname{min}} \quad   \widetilde{\phi}(\bm{u},\bm{v}),
	\quad \text{s.t.} \quad \widetilde{\bm{g}}(\bm{u},\bm{v}) \leq \bm{0},  \label{eq:2___26_sol_v_def} 
\end{align}	
and the  Lagrange multipliers  associated to  $\bm{sol}(\bm{v})$ are:  
\begin{align} 
	\bm{\lambda}(\bm{v}) := \  & \
	\left\{ 
	\bm{a}\in\amsmathbb{R}^{n_g}\geq \bm{0} \ \left| \
	\widetilde{\phi}(\bm{sol}(\bm{v}),\bm{v}) + \sum_{i=1}^{n_g}
	a_{(i)}
	\widetilde{g}_{(i)}(\bm{sol}(\bm{v}),\bm{v})
	\right.\right\}. \label{eq:2___26_lambda_v_def} 
\end{align}	

These two functions have properties which are given by the following theorem and which are essential for the analysis of the stability of a purely iterative RTO method. 
\begin{thmbox} \textbf{(Properties of purely iterative RTO methods)}
	\label{thm:2___3_Proprietes_Sol}
	\\ 
	If:
	\vspace{-\topsep} 
	\begin{itemize}[]	
		\item A purely iterative RTO method is used;
		\item For any correction point $\bm{v}\in\amsmathbb{R}^{n_u}$ the following relations are true:
		\begin{align}
			\widetilde{\bm{g}}|_{\bm{v},\bm{v}}  = \ & \bm{g}_p|_{\bm{v}}, &
			\nabla_u \widetilde{\bm{g}}|_{\bm{v},\bm{v}}   = \ & \nabla_u \bm{g}_p|_{\bm{v}},  & 
			\nabla_u \widetilde{\phi}|_{\bm{v},\bm{v}}  = \ & \nabla_u \phi_p|_{\bm{v}},
		\end{align}
		\begin{align} \label{eq:2_29_uywgkfyurwkc}
			\nabla^2_{\bm{u}\bm{u}}\widetilde{\phi}|_{\bm{v},\bm{v}}  + 
			\sum_{i=0}^{n_g}\left[
			\lambda^{\prime}_{(i)}  \nabla^2_{\bm{u}\bm{u}}\widetilde{g}_{(i)}|_{\bm{v},\bm{v}} \right]  > \ & 0, & 
			\forall \bm{\lambda}^{\prime}\in\amsmathbb{R}^{n_g+},
		\end{align}
		where the functions $(\widetilde{\phi}, \widetilde{\bm{g}})$ are defined by \eqref{eq:2___25_DefMethodePurementIteratives}; 
		\item For any correction point $\bm{v}\in\amsmathbb{R}^{n_u}$ the problem \eqref{eq:2___26_sol_v_def} is such that  $\bm{sol}(\bm{v})$ is unique and the active constraints at  $\bm{sol}(\bm{v})$ satisfy the Linear Independence Constraint Qualification (LICQ) conditions.
	\end{itemize}
	\vspace{-\topsep} 
	\medskip
	Then the following statements are true:
	\vspace{-\topsep} 
	\medskip
	\begin{itemize}[noitemsep]
		\item[\textbf{A.}] $\bm{sol}(\bm{v})$ and $\bm{\lambda}(\bm{v})$ are globally $\mathcal{C}^0$ and piecewise-$\mathcal{C}^1$ ;
		\item[\textbf{B.}]  If  $\bm{v}^{\bullet}$ is a KKT point of the plant and $\delta \bm{v}$ a vector in $\amsmathbb{R}^{n_u}\backslash \bm{0}$,  then $\bm{sol}(\bm{v}^{\bullet}) = \bm{v}^{\bullet}$ and the value of the directional derivative of $\bm{sol}$
		\begin{align} \label{eq:2___29_DeriveeDirectionnelle_Forward}
			\nabla^S(\bm{v},\delta\bm{v})  := \ &
			\left( \frac{\delta \bm{v} }{\|\delta \bm{v} \|}
			\right)^{\rm T}
			\lim_{\begin{smallmatrix} \alpha \to 0 & \\ \alpha>0 \end{smallmatrix}} \frac{
				\bm{sol}(\bm{v}+\alpha\delta \bm{v})
				- \bm{sol}(\bm{v}) }{\alpha}, \\
			\Big( = \ &
			\left( \frac{\delta \bm{v}}{\|\delta \bm{v}\|}
			\right)^{\rm T} \nabla_{v}\bm{sol}|_{\bm{v}} \frac{\delta \bm{v}}{\|\delta \bm{v}\|}, \quad \text{ si $\bm{sol}$ est $\mathcal{C}^1$ en $\bm{v}$}
			\Big). \nonumber
		\end{align}
		\begin{equation}\label{eq:2___29_NablaS_Signe}
			\text{is: }
			\left\{
			\begin{array}{ll}
				\nabla^S(\bm{v}^{\bullet},\delta\bm{v}) < 1, & \text{if $\bm{v}^{\bullet}$ is a minimum of the plant}, \\
				\nabla^S(\bm{v}^{\bullet},\delta\bm{v}) > 1, & \text{if $\bm{v}^{\bullet}$ is a maximum of the plant}, \\
				\nabla^S(\bm{v}^{\bullet},\delta\bm{v}) \gtrless 1, & \text{if $\bm{v}^{\bullet}$ is a saddle point of the plant.}
			\end{array}
			\right. \qquad \quad 
		\end{equation}
		\item[\textbf{C.}] 
		If $\bm{v}^{\bullet}$ is a minimum of the plant and $\delta\bm{v}$ is a vector orthogonal to the set of active constraints of the plant at $\bm{v}^{\bullet}$, then in order to have $\nabla^{S}(\bm{v},\delta\bm{v})<-1$, the Lagrangian of the plant must be at least ``two times more convex'' than the one of the model in the direction $\delta\bm{v}$, i.e.:  
		\begin{align} \label{eq:2___32_Condition_NablaS_pluspetitque_m1} 
			\frac{\delta\bm{v}^{\rm T}  
				(\nabla_{\bm{uu}}\mathcal{L}|_{\bm{v}^{\bullet},\bm{v}^{\bullet}})^{-1}
				\nabla_{\bm{uu}}\mathcal{L}_p|_{\bm{v}^{\bullet}} \delta\bm{v}}{\|\delta\bm{v}\|^2} > \ & 2, & 
			\Leftrightarrow \quad 
			\nabla^{S}(\bm{v},\delta\bm{v})< \ & -1.
		\end{align}
	\end{itemize}	
\end{thmbox}
\begin{proofbox}
	The proof of this theorem is given in Appendix~\ref{App:A_3_Preuve_THM2_3}. 
\end{proofbox}

When this theorem is applicable, the two successive iterations $\bm{u}_k = \bm{u}_p^{\star} + \delta\bm{u}_k$  and  $\bm{u}_{k+1}$ can be linked  with the following equation: 

\begin{align} 
	& &\bm{u}_{k+1} = \ & \bm{sol}(\bm{u}_p^{\star} + \delta\bm{u}_k), \nonumber \\
	& &             = \ & \bm{sol}(\bm{u}_p^{\star}) + 
	\nabla^S(\bm{u}_p^{\star},\delta\bm{u}_k) \delta\bm{u}_k + \mathcal{O}(\| \delta \bm{u}_k\|^2), \nonumber \\
	& &			 = \ & \bm{u}_p^{\star} + 
	\nabla^S(\bm{u}_p^{\star},\delta\bm{u}_k) \delta\bm{u}_k + \mathcal{O}(\| \delta \bm{u}_k\|^2), \label{eq:2___28_C1_Uasteri} \\
	\Leftrightarrow & & 
	\bm{u}_{k+1} - \bm{u}_p^{\star} = \ &  \nabla^S(\bm{u}_p^{\star},\delta\bm{u}_k) (\bm{u}_{k} - \bm{u}_p^{\star}) + \mathcal{O}(\| \bm{u}_{k} - \bm{u}_p^{\star}\|^2). \label{eq:2___33_rcugbk}
\end{align}
And when  $\delta\bm{u}_k\rightarrow\bm{0}$  the term $\mathcal{O}(\| \bm{u}_{k} - \bm{u}_p^{\star}\|^2)$ is negligible. 
\begin{align} 
	\text{\eqref{eq:2___33_rcugbk} } \Leftrightarrow \ 
	\|\bm{u}_{k+1} - \bm{u}_p^{\star} \|=  \nabla^S(\bm{u}_p^{\star},\delta\bm{u}_k)^2 \|\bm{u}_{k} - \bm{u}_p^{\star}\|. \nonumber
\end{align}
It is possible to derive from this equation the condition for $\bm{u}_{k+1}$ to be closer to $\bm{u}_p^{\star}$ than $\bm{u}_{k}$ and this, for all $\bm{u}_{k}\in\mathcal{B}(\bm{u}_p^{\star},r\rightarrow 0)$, is: 
\begin{align} \label{eq:3___36_kcunrg}
	-1<\nabla^S(\bm{u}_p^{\star},\delta\bm{u}_k)< \ & 1, & 
	\forall \delta\bm{u}_k \in\amsmathbb{R}^{n_u}.
\end{align}
If the model is convex at the correction point, then thanks to the equation \eqref{eq:2___29_NablaS_Signe} of  Theorem~\ref{thm:2___3_Proprietes_Sol}, $\forall \delta\bm{u}_k \in\amsmathbb{R}^{n_u}$,  $\nabla^S(\bm{u}_p^{\star},\delta\bm{u}_k)<1$. However, there is no guarantee that $\forall \delta\bm{u}_k\in\amsmathbb{R}^{n_u}, \nabla^S(\bm{u}_p^{\star},\delta\bm{u}_k) > -1$; thus, the condition \eqref{eq:3___36_kcunrg} can be reduced to:
\begin{align} \label{eq:2___36_Condition_Niveau2_SansFiltre}
	\Underline{\nabla}^S(\bm{u}_p^{\star}) >\ & -1, & \text{where} \qquad 
	\Underline{\nabla}^S(\bm{u}_p^{\star}) := \ & \underset{\delta\bm{u}_k\in\amsmathbb{R}^{n_u}}{\operatorname{min}} \  \nabla^S(\delta\bm{u}_k).
\end{align}
which is therefore the stability condition of all RTO methods satisfying the applicability conditions of Theorem~\ref{thm:2___3_Proprietes_Sol} (e.g. ISO-D/I using convexified models). \\

It is to eliminate the condition \eqref{eq:2___36_Condition_Niveau2_SansFiltre} that \cite{Marchetti:09b,Marchetti:09} propose to apply a filter on the iterations of the decision variables\footnote{It is also shown that an equivalent result can be obtained by applying a filter on the model updates. But this option is not considered in this thesis.}. 

\subsection{Existence of a filter enforcing stability}
A filter on the iterations is implemented as follows: 
\begin{subequations} \label{eq:2___25_ImplementationClassiqueFiltre}
	\begin{enumerate}[noitemsep]	
		\item The minimum of the updated model is computed: 
		\begin{equation} \label{eq:2___27a_ModelBasedPB}
			\bm{u}^{\star}_{k+1} = \operatorname{arg}
			\underset{\bm{u}}{\operatorname{min}} \quad   \phi^c_{k}(\bm{u})  \quad 
			\text{s.t.} \quad  \bm{g}^c_{k}(\bm{u}) \leq \bm{0}.
		\end{equation}
		\item Then one selects $\bm{u}_{k+1}$ between $\bm{u}_k$ and $\bm{u}^{\star}_{k+1}$ with  a filtering gain\footnote{In the literature a filtering matrix  $\bm{K}\in\amsmathbb{R}^{n_u \times n_u}$ whose eigenvalues are in $]0,1]$ is generally used. Here, one indirectly chooses the matrix $\bm{K}=K\bm{I}$, and this choice makes the following developments far more readable.} $K\in]0,1]$:
		\begin{equation} \label{eq:2___27b_Filter}
			\bm{u}_{k+1} = \bm{u}_k + K(\bm{u}^{\star}_{k+1} - \bm{u}_k).
		\end{equation}
	\end{enumerate}
\end{subequations}
In less mathematical terms, implementing a filter on an iteration means choosing the next operating point $\bm{u}_{k+1}$ between the current one $\bm{u}_{k}$ and the minimum of the updated model $\bm{u}_{k+1}^\star$.

It is shown  in \cite{Marchetti:09b,Marchetti:09} that the use of a filter $K$ which is ``sufficiently small'' enforces the stability conditions for any RTO method satisfying the applicability conditions of Theorem~\ref{thm:2___3_Proprietes_Sol}.  To achieve this result, it is sufficient to start from  \eqref{eq:2___27b_Filter}:
\vspace{-\topsep}
\begin{itemize}[noitemsep]
	\item replace $\bm{u}_{k+1}^*$ by \eqref{eq:2___28_C1_Uasteri};
	\item neglect the term $\mathcal{O}(\|\delta \bm{u}_k\|^2)$ to place the study in the neighborhood of $\bm{u}_p^{\star}$ (otherwise it is not the stability but the superstability that would be analyzed -- as discussed later);
	\item add ``$-\bm{u}_p^{\star}$'' on each side of the equality;
\end{itemize}
\vspace{-\topsep}
to get
\vspace{-\topsep}
\begin{align}
	& &	\bm{u}_{k+1} - 	\bm{u}^{\star}_{p} = \ & 
	\bm{u}_k + K (\bm{u}^{\star}_p + 
	\nabla^S(\bm{u}_p^{\star},\delta\bm{u}_k) 
	(\bm{u}_k-\bm{u}_p^{\star}) - \bm{u}_k) - \bm{u}^{\star}_{p}, \nonumber \\
	& &	 = \ & 
	\Big(1 - K \left( 1 - 
	\nabla^S(\bm{u}_p^{\star},\delta\bm{u}_k)
	\right)\Big)
	(\bm{u}_{k} - 	\bm{u}^{\star}_{p} ), \nonumber \\
	\Rightarrow & & 
	\|\bm{u}_{k+1} - 	\bm{u}^{\star}_{p}\| = \ & 
	\Big(1 - K \left( 1 - 
	\nabla^S(\bm{u}_p^{\star},\delta\bm{u}_k)
	\right)\Big)^2 \|\bm{u}_{k} - 	\bm{u}^{\star}_{p} \|.
	\label{eq:2___30_IterFiltreEffet}
\end{align}
Therefore, in order for $\bm{u}_{k+1}$ to be closer to $\bm{u}_p^{\star}$ than $\bm{u}_{k}$ and this, for all $\bm{u}_{k}\in\mathcal{B}(\bm{u}_p^{\star},r\rightarrow 0)$, it is necessary that: 
\vspace{-\topsep}
\begin{align}
	& & \Big(1 - K \left( 1 - 
	\nabla^S(\bm{u}_p^{\star},\delta\bm{u}_k)
	\right)\Big)^2 < \ & 1, \qquad \forall \delta\bm{u}_k\in\amsmathbb{R}^{n_u}, \nonumber \\
	\Leftrightarrow & &  -1 < 
	1 - K \left( 1 - 
	\nabla^S(\bm{u}_p^{\star},\delta\bm{u}_k)
	\right) < \ & 1, \qquad \forall \delta\bm{u}_k\in\amsmathbb{R}^{n_u}, \nonumber  \\
	\Leftrightarrow & &  0 < 
	 K \left( 1 - \nabla^S(\bm{u}_p^{\star},\delta\bm{u}_k) 
	\right) < \ & 2, \qquad \forall \delta\bm{u}_k\in\amsmathbb{R}^{n_u}, \nonumber  \\ 
	\Leftrightarrow & &  0 < 
	K \left( 1 - \Underline{\nabla}^S(\bm{u}_p^{\star}) 
	\right) < \ & 2, \quad \Underline{\nabla}^S(\bm{u}_p^{\star}) :=  \underset{\delta\bm{u}_k\in\amsmathbb{R}^{n_u}}{\operatorname{min}} \  \nabla^S(\delta\bm{u}_k). \label{eq:2___30_ConditionNiveau2}
\end{align}
Since $\bm{u}_p^{\star}$ is a minimum of the plant,
\begin{equation*}
	\text{ \eqref{eq:2___29_NablaS_Signe} }   \Rightarrow \quad
	1 - \nabla^S(\bm{u}_p^{\star},\delta\bm{u}_k) > 0, \ \forall \delta\bm{u}_k\in\amsmathbb{R}^{n_u} \quad 
	 \Rightarrow \quad 1- \Underline{\nabla}^S(\bm{u}_p^{\star}) > 0.
\end{equation*}
Therefore, the stability condition is a condition on $K$:
\begin{equation} \label{eq:2___36_ConditionNiveau2_K}
	0 < K < \frac{2}{1-\Underline{\nabla}^S(\bm{u}_p^{\star}))}, 
\end{equation}
and this set is never empty since $1-\nabla^S(\bm{u}_p^{\star},\delta\bm{u}_k)>0$.

Finally, to choose a ``small enough'' filter one just needs to choose $K$ such that the inequality \eqref{eq:2___36_ConditionNiveau2_K} is satisfied.  However, it is clear that it is a priori impossible to evaluate this inequality because to do so one would need to know both $\bm{u}_p^{\star}$ and the directional derivatives $\nabla^S$ of the function $\bm{sol}$, and neither of them are a priori known.

\subsection{Towards an adaptive filter}

To get around this problem, one proposes (i) to compute at each iteration an approximation $\utilde{\text{\hskip 0.1ex $\nabla$}}_k^S$  of $\Underline{\nabla}^S(\bm{u}_p^{\star})$ and (ii) to choose at each iteration a filter $K_k$ such that: 
\begin{align} \label{eq:2___41_strat1_adaptationK}
		K_k < \ &  \frac{2}{1-\utilde{\text{\hskip 0.1ex $\nabla$}}^S_{k}}, & 
		\utilde{\text{\hskip 0.1ex $\nabla$}}_k^S  := \ &
		\left( \frac{\bm{u}_k - \bm{u}_{k-1}}{\|\bm{u}_k - \bm{u}_{k-1}\|}
		\right)^{\rm T}
		\frac{\bm{u}^{*}_{k+1} - \bm{u}^{*}_{k}}{\| \bm{u}_{k} - \bm{u}_{k-1} \| }, &
		\text{si } \utilde{\text{\hskip 0.1ex $\nabla$}}_k^S < 1.
\end{align}	
Any strategy of adapting the filter $K$ that follows \eqref{eq:2___41_strat1_adaptationK} provides interesting properties in particular:
\begin{thmbox} 		
	\label{thm:2___4_K_Niveau2CasParticulier}
	If the conditions of applicability of the Theorem~\ref{thm:2___3_Proprietes_Sol} are satisfied and if:
	\vspace{-\topsep}
	\begin{itemize}[noitemsep]
		\item The problem is unidimensional, i.e. $n_u=1$;
		\item The function $sol$ is $\mathcal{C}^1$ at $\bm{u}_p^{\star}$;
		\item A filter $K>0$ is applied by following the procedure \eqref{eq:2___25_ImplementationClassiqueFiltre}.
	\end{itemize}
	\vspace{-\topsep}
	Then any filter adaptation strategy that satisfies \eqref{eq:2___41_strat1_adaptationK}  guarantees that the stability condition \eqref{eq:2___36_ConditionNiveau2_K} is always satisfied.
\end{thmbox}
\begin{proofbox} 
	Given that  $u_p^{\star}$ is a minimum of the plant, that the function $sol$ is $\mathcal{C}^1$ at $u_p^{\star}$, and that Theorem~\ref{thm:2___3_Proprietes_Sol} is applicable, the function $sol$ can be approximated by:
	\begin{align} \label{eq:2___42_ApproximationLocale}
		sol(u) = \ & u_p^{\star} + \nabla_u sol|_{u_p^{\star}} (u-u_p^{\star}), &
		\forall \bm{u} \in \amsmathbb{D}_{r} =  [u_p^{\star}-r, \ u_p^{\star}+r].
	\end{align}
	In domain $\amsmathbb{D}_{r}$, there is always a sub-domain $\amsmathbb{D}_{\epsilon} = [u_p^{\star}-\epsilon, \ u_p^{\star}+\epsilon]$ such that whatever is the filter $K>0$ used, $u_k \in \amsmathbb{D}_{\epsilon} \ \Rightarrow \ u_{k+1} \in \amsmathbb{D}_{r}$. Indeed,
	\begin{equation*}
		\begin{array}{rrcl}
			                  & u_p^{\star} - r  < &  u_{k+1}                   & < u_p^{\star} + r  \\ 
			\Leftrightarrow   & u_p^{\star} - r  < &  u_{k} + K(sol(u_k)-u_{k}) & < u_p^{\star} + r  \\
			\Leftrightarrow   & u_p^{\star} - r  < &  u_{k} + K(u_p^{\star} + \nabla_usol|_{u_p^{\star}}(u_k-u_p^{\star})-u_{k}) & < u_p^{\star} + r  \\
			\Leftrightarrow   & u_p^{\star} - \epsilon  < &  u_{k}  & < u_p^{\star} + \epsilon  \\
		\end{array}
	\end{equation*}
	where
	\begin{equation*}
		\epsilon := \frac{r}{|1-K+K\nabla_u sol|_{u_p^{\star}}|}.
	\end{equation*}
	So, if an iterate  $u_k\in\amsmathbb{D}_{\epsilon}$, then $\forall K>0$ the next iterate will satisfy: $u_{k+1}\in\amsmathbb{D}_{\delta}$. Given these two iterations and because (i) locally \eqref{eq:2___42_ApproximationLocale} holds, (ii) $sol$ is $\mathcal{C}^1$, (iii) $n_u=1$, the estimate $\utilde{\text{\hskip 0.1ex $\nabla$}}_k^S$ is equal to the actual value of $\nabla_u sol|_{u_p^{\star}}$.
	Therefore, it is sure that the filter chosen for the iteration  $u_{k+2}$  satisfies \eqref{eq:2___36_ConditionNiveau2_K}. As a result, $u_{k+2}\in\amsmathbb{D}_{r}$ and for all the following iterations \eqref{eq:2___36_ConditionNiveau2_K} will be systematically respected. It follows from this observation that the sequence $\{u_{k+1},...,u_{\infty}\}$ converges asymptotically on $u_p^{\star}$. 
	
	Finally, it is enough that one iterate falls in the domain $\amsmathbb{R}_\epsilon$ so that the filter chosen with \eqref{eq:2___41_strat1_adaptationK} guarantees that the next iterates converge asymptotically on $\bm{u}_p^{\star}$.
	 Noting that $\amsmathbb{D}_r$ is a ball $\mathcal{B}(\bm{u}_p^{\star},r \rightarrow 0)$ in $\amsmathbb{R}$, it can be concluded that any strategy that satisfies \eqref{eq:2___41_strat1_adaptationK} guarantees the satisfaction of the stability condition \eqref{eq:2___36_ConditionNiveau2_K}.
\end{proofbox}

Hereafter, a graphical interpretation of Theorem~\ref{thm:2___4_K_Niveau2CasParticulier} is proposed:
\begin{Side}
\begin{GraphicalInterpretation}   \label{gi:2___1_Niveau2}
	If the conditions of applicability of Theorem~\ref{thm:2___4_K_Niveau2CasParticulier} are satisfied by a RTO method and if no filter is used. 
	 Then, $\forall u_o$ close to a stationary point of the plant $u_p^\bullet$ the  iterations $u_{k>0}$ will follow a behavior similar to one of the 4 cases shown in Figure~\ref{fig:2___3_QuatreScenarios}. The four cases of Figure~\ref{fig:2___3_QuatreScenarios} correspond to situations where:
	\begin{align*}
		\text{case A: } & \ \nabla_u sol|_{u_p^{\star}} \in ]1,+\infty[, & 
		\text{case B: } & \ \nabla_u sol|_{u_p^{\star}} \in ]0,1[, \\
		\text{case C: } & \ \nabla_u sol|_{u_p^{\star}} \in  ]-1,0[, & 
		\text{case D: } & \ \nabla_u sol|_{u_p^{\star}} \in ]-\infty,-1].
	\end{align*}
	
	\begin{minipage}[h]{\linewidth} 	
		\smallskip                   
		\centering                   
		\includegraphics[width=13.5cm] {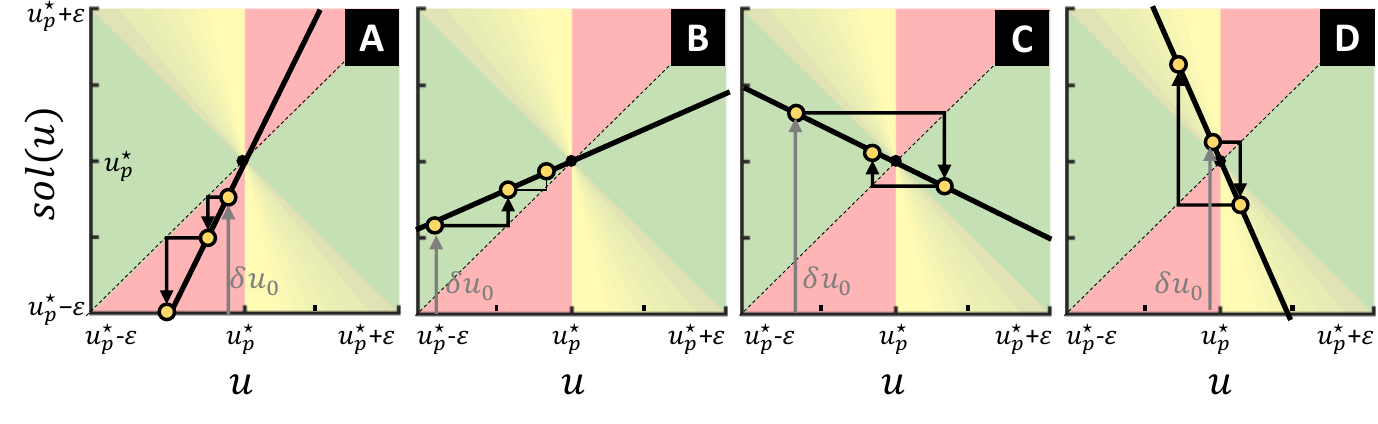}
		\captionof{figure}{Four possible cases when $n_u=1$ and $sol$ is $\mathcal{C}^{1}$ at $u_p^{\bullet}$.}
		\label{fig:2___3_QuatreScenarios} 	
		\smallskip
	\end{minipage}\\

	The black line represents the function $sol(u)$. This line crosses the point $[u_p^\bullet,u_p^\bullet]$ because according to Theorem~\ref{thm:2___3_Proprietes_Sol} $sol(u_p^\bullet)=u_p^\bullet$. The evolution of the iterations is represented by yellow points on this curve. \textit{(One starts from $\delta u_k = \delta u_0 $ to generate the first point $[\delta u_0, sol(u_p^\bullet + \delta u_0)]^{\rm T}$. As no filter is applied, the ordinate of one point is the abscissa of the next one. So, to construct a new point (i.e. move forward one iteration) one must horizontally project a point on the line $y= \delta u_k + u_p^\bullet$ (black dotted line) to then project it vertically on the curve $u_{k+1} = sol(u_p^\bullet + \delta u_k)$.)} Finally, a color code is used to clearly show in which case (A,B,C,D) one is located. If the function $sol$ is in the green domain, then it is the B or C and a filter is not necessary to make the RTO method stable, i.e. $\nabla_u sol|_{u_p^\bullet}\in[-1,1]$ which can be linked to \eqref{eq:3___36_kcunrg}. On the other hand, if it is in the yellow or red areas, then one can observe that the iterations do not converge on $u_p^{\bullet}$. \\
	
	Now, let's build scenarios A' and D' similar to scenarios A and D but using a filter \textit{(i.e. ``where each horizontal projection is reduced'')}.  Figure~\ref{fig:2___4_AllCones} illustrates what is obtained and two things can be observed:

	(Case D':) If the function $sol$ is in the yellow area (i.e. if $\nabla_u sol|_{u_p^\bullet}< -1$) then a filter can provide stability. The figure on the very right of  Figure~\ref{fig:2___4_AllCones} gives a graphical interpretation of the appropriate choices of $K$. 
	
	(Case A':) If the function  $sol$ is in the red area (i.e. if  $\nabla_u sol|_{u_p^\bullet}> 1$), then there is no $K>0$ providing stability. However, according to the Theorem~\ref{thm:2___3_Proprietes_Sol}, $\nabla_u sol|_{u_p^\bullet}> 1$ implies that $u_p^\bullet$ s not a minimum of the plant. \textit{So this particular case, is a case for which the convergence on $u_p^\bullet$ is not desirable. Hence the interest in setting a lower bound on the choice of $K$ to $0$.}  
	\\
	\begin{minipage}[h]{\linewidth} 	
		\smallskip   
		\vspace*{0pt}
		\centering
		\includegraphics[width=13.5cm] {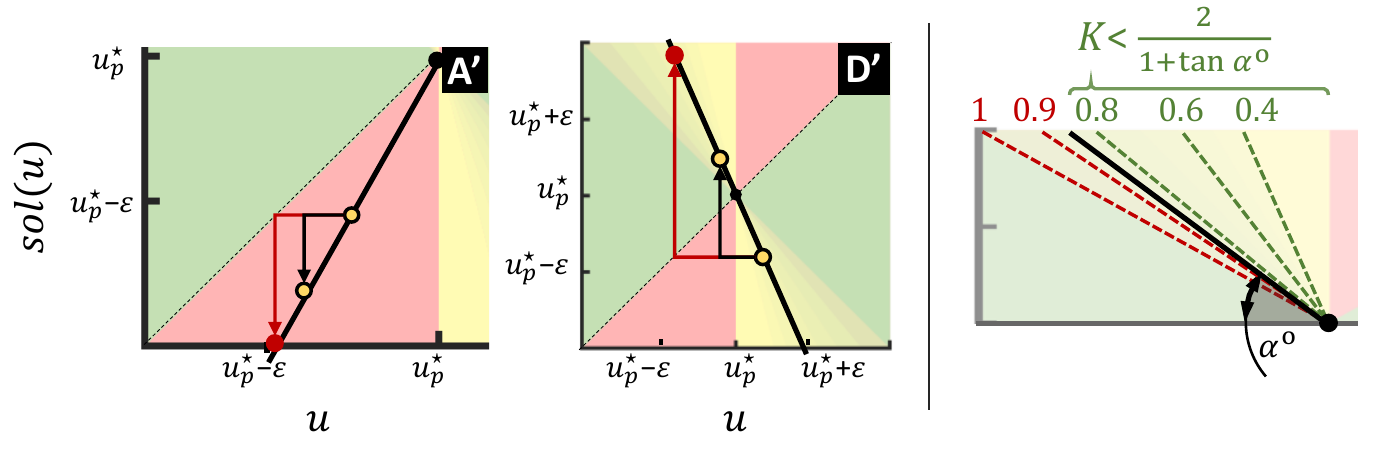}
		\captionof{figure}{Illustration of the effects of the filter.}
		\label{fig:2___4_AllCones}
	\end{minipage}\\
	
	Finally, the green, yellow and red cones whose common vertex is $u_p^{\bullet}$ provide the following information: ``If for a point $u$ the point $[u,sol(u)]$ is within the cone of color
	\vspace{-\topsep}
	\begin{itemize}[noitemsep]	
		\item green, then $sol(u)$ will be \textit{closer} to $u_p^{\bullet}$ than  $u$.
		\item yellow, then $sol(u)$ will be \textit{closer} to $u_p^{\bullet}$ than $u$ \textit{if a small enough filter is used }.
		\item red, then $sol(u)$ will be \textit{farther} from  $u_p^{\bullet}$ than $u$.
	\end{itemize}
	\vspace{-\topsep}
	These color cones are reused later in this chapter. 	
\end{GraphicalInterpretation}
\end{Side}

\subsection{Towards an optimal adaptive filter}

To go further, one can take the equation \eqref{eq:2___30_IterFiltreEffet} and choose the filter $K$ which minimizes the term ``$1 - K \big( 1 - 
\nabla^S(\bm{u}_p^{\star},\delta\bm{u}_k)
\big)$'' so that $\bm{u}_{k+1}$ is as close as possible to $\bm{u}_p^{\star}$. Of course, $\nabla^S(\bm{u}_p^{\star},\delta\bm{u}_k)$ is unknown but can be approximated with \eqref{eq:2___41_strat1_adaptationK}. So, the idea is to define the filter to be used as the solution of the following optimization problem:
\begin{equation} \label{eq:2___43__filtre_Optimal}
	K^\star_k = 
	\operatorname{arg} \ \underset{K}{\operatorname{min}} \quad  \left| 
	1 - K \left[ 
	1 -  \utilde{\text{\hskip 0.1ex $\nabla$}}^S_{k}	\right]  \right| ,
	\quad  \text{s.t.} \quad  K > K_o, \quad  \text{si } \utilde{\text{\hskip 0.1ex $\nabla$}}_k^S < 1.   
\end{equation}	
where $K_o>0$ is the minimum value that the filter can take, which must be chosen by the engineers supervising the plant. The default choice is set to $K_o:=0.1$.   This optimization problem has a simple explicit solution which is:
\begin{align} \label{eq:2___36_OptimizationPB_ChoiceK}
	K^{\star}_k = \ & \max\left\{K_o,   \frac{1}{1- \utilde{\text{\hskip 0.1ex $\nabla$}}^S_{k}} \right\}, &
	\text{if } \utilde{\text{\hskip 0.1ex $\nabla$}}_k^S < 1. 
\end{align}
Notice that by definition   $K_k^\star$ always satisfies the condition \eqref{eq:2___41_strat1_adaptationK} because it is obvious that it is always true that:
\begin{align*}
	\frac{1}{1- \utilde{\text{\hskip 0.1ex $\nabla$}}^S_{k}} < \ & \frac{2}{1- \utilde{\text{\hskip 0.1ex $\nabla$}}^S_{k}}, &
	 \text{if: } \qquad \utilde{\text{\hskip 0.1ex $\nabla$}}^S_{k} <1.
\end{align*}.

Let's illustrate the effects of using this adaptive filter on two purely mathematical examples.
\begin{itemize}
	\item Example~\ref{ex:2___1_Exemple1} is a detailed analysis of the solving of a 1D problem ($n_u=1$) satisfying all the conditions of  Theorem~\ref{thm:2___4_K_Niveau2CasParticulier}. Some practical points related to the implementation, such as convergence management, are discussed. 
	\item Example~\ref{ex:2___2_Kk_nu_geq2} proposes an empirical analyses of the use of the adaptive filter  \eqref{eq:2___36_OptimizationPB_ChoiceK}. It is observed that very interesting behaviors emerges from the use of this adaptive filter. 
\end{itemize}

\begin{exbox}
	\label{ex:2___1_Exemple1}
	\emph{(This is an example used in \cite{Marchetti:09b} (the example 4.3 page 112).)} One considers a problem where the functions of the plant and the model are:
	\begin{align*}   
		\phi(u,y) := \ & y, \ & 
		y_p = f_p(u) := \ & (u-1)^2,  \    & 
		y = f(u)  :=  \ & \frac{u^2}{4},
	\end{align*}
	and where the constraints are simply on the inputs: $-5 \leq u \leq 5$. In this case, it is clear that ISO-D/I are identical and that the equilibrium condition is always satisfied since the cost function is convex at all points and the constraints are known. It is therefore a ``simple'' case study allowing for the effects of the filter to be analyzed. The rest of this example is divided into two parts.  In part 1, a theoretical analysis of the convergence properties of ISO-D/I w.r.t. the filter is proposed.  In part 2, simulation results are shown and some "practical" aspects of the implementation of the adaptive filter are discussed.

	\begin{center}
		\textbf{Part 1: Theoretical analysis}
	\end{center}

	The minimum of the plant is $u_p^{\star} = 1$ and no constraints are active at this point. So the function $sol$ given in Figure~\ref{fig:2___5_Exemple_SolUk} can be approximated around $u_p^{\star}$ in the following way:  
	\begin{align*}
		 & & \widetilde{\phi}(u,u_p^{\star}) := \ & \phi_p|_{u_p^{\star}}  + \nabla_u \phi_p|_{u_p^{\star}} (u - u_p^{\star}) + \frac{1}{2} \nabla^2_{uu} \phi|_{u_p^{\star}} (u - u_p^{\star})^2, \nonumber  \\
		\Rightarrow & & \nabla_u \widetilde{\phi}(u,u_p^{\star}) = \ & \nabla_u \phi_p|_{u_p^{\star}} + \nabla^2_{uu} \phi|_{u_p^{\star}} (u-u_p^{\star}), \\
		\Rightarrow  & & 0 = \ &  \nabla_u \phi_p|_{u_p^{\star}} + \nabla^2_{uu} \phi|_{u_p^{\star}} (sol(u_p^{\star})-u_p^{\star}), \nonumber \\
		\Rightarrow  & & sol(u_p^{\star}) = \ & u_p^{\star} - \frac{\nabla_u \phi_p|_{u_p^{\star}}}{\nabla^2_{uu} \phi|_{u_p^{\star}}}. \label{eq:A___12_Solution}
	\end{align*}
One can observe on Figure~\ref{fig:2___5_Exemple_SolUk} (i) that the function $sol(u)$ is piecewise-$\mathcal{C}^1$  and globally $\mathcal{C}^0$; (ii) that the function $sol(u)$ is $\mathcal{C}^1$ at $u_p^{\star}$; and (iii)  that a filter will be necessary to bring stability since the curve is in the yellow area around $u_p^{\star}$.
\begin{align*}
	\nabla^S(u_p^{\star},1) = \nabla^S(u_p^{\star},-1) = \ & \nabla_{u_k} sol|_{u_p^{\star}} 
	, \\
	= \ &  1- \frac{\nabla^2_{uu} \phi_p|_{u_p^{\star}} \nabla^2_{uu} \phi|_{u_p^{\star}}   -  \nabla^3_{uuu}\phi|_{u_p^{\star}} \nabla_{u} \phi_p|_{u_p^{\star}}  }{\nabla^2_{uu}\phi|_{u_p^{\star}} ^2}, \\
	= \ & 1- \frac{\nabla^2_{uu} \phi_p|_{u_p^{\star}} }{\nabla^2_{uu}\phi|_{u_p^{\star}} } = 1- \frac{2}{1/2} = -3.
\end{align*}
Therefore, in order to satisfy the stability condition one should choose the filter in the interval: $K \in]0,\ \frac{2}{1-(-3)}=\frac{1}{2}[$. And according to \eqref{eq:2___36_OptimizationPB_ChoiceK}, the optimal filter is:
\begin{equation} \label{eq:2___44_Exemple2_1_filtreoptimal}
	K^{\star} = \frac{1}{1-(-3)} = \frac{1}{4}.
\end{equation}
\begin{minipage}[h]{\linewidth}
	\smallskip
	\centering
	\includegraphics[width=4.45cm]{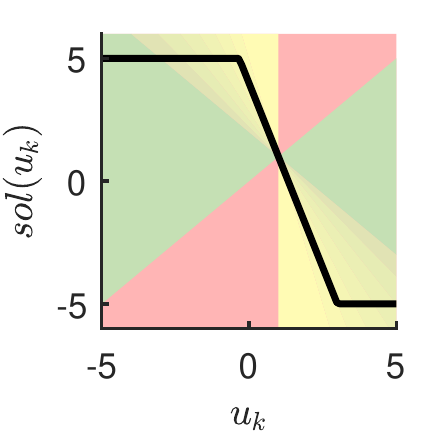} 
	\captionof{figure}{Example~\ref{ex:2___1_Exemple1}: The function $sol(u)$.}
	\label{fig:2___5_Exemple_SolUk} 	
	\smallskip
\end{minipage} \\
	
	\begin{center}
		\textbf{Part 2: Simulations \& Implementation.}
	\end{center}
	
	One runs several simulations starting from the minimum of the model: $u_0=0$, and one considers the three following approaches:  \\

	\textbf{The classical approach (MA and MAy)} which consists in letting the engineers supervising the plant predefine the filter to be used. Figure~\ref{fig:2___3_Approche1Resultats} shows the results obtained for different filter choices:  $K=\{0.1, 0.2, 0.25,0.3, 0.4,$ $0.5\}$. It can be noticed that $K=0.25$ is indeed the best filter since it allows an immediate convergence. On the other hand, all filter $K>0.5$ do not enable convergence. So, depending on the choice that is made this approach may, or may not, converge to the plant optimum. \\
	\begin{minipage}[h]{\linewidth}
		\smallskip
		\centering
		\includegraphics[width=4.45cm]{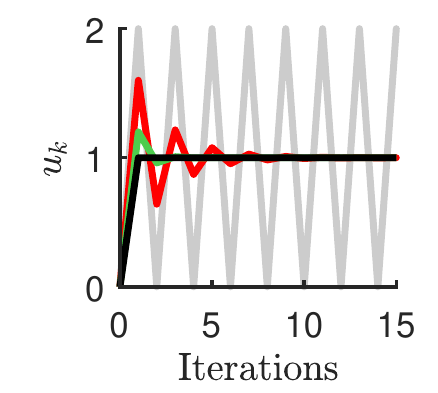}\hskip -0ex
		\includegraphics[width=4.45cm]{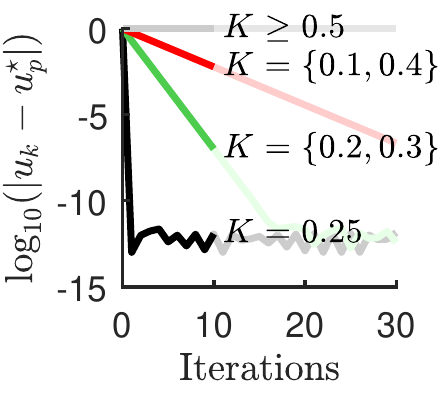}
		\captionof{figure}{The classical approach: Results for multiple filter choices .}
		\label{fig:2___3_Approche1Resultats} 	
		\smallskip
	\end{minipage} \\
	
	\textbf{The maximal approach} which consists in choosing, at each iteration, the largest filter which guarantees stability, i.e. at each iteration the filter that is applied is:
	\begin{align}
		& \text{If } k =0 : & 	K_k = \ &  K_o, \\
		 & \text{If } \utilde{\text{\hskip 0.1ex $\nabla$}}^S_{k} < 1 \text{ et } k > 0  : &
		K_k = \ & \max\Big\{K_o := 0.1, \ 2/(1-\utilde{\text{\hskip 0.1ex $\nabla$}}^S_{k}) - 0.1\Big\}, \label{eq:2___gcerbgfj}\\
		& \text{Otherwise:} & 	K_k = \ &  1. \qquad \text{(No filter)}
	\end{align}
	The	 idea is to estimate the critical value of the filter at the iteration $k$, then select a slightly smaller one (hence ``$-0.1$'') guaranteeing the satisfaction of the stability condition.  \\
	 
	\textbf{The optimal approach} which consists in choosing at each iteration the filter by following the following rules:
	\begin{align*}	
		& \text{If } k =0 : & 	K_k = \ &  K_o, \\
		& \text{If } \utilde{\text{\hskip 0.1ex $\nabla$}}^S_{k} < 1 \text{ and } k > 0  :  &
		K_k = \ & \text{\eqref{eq:2___36_OptimizationPB_ChoiceK}}, \\
		& \text{otherwise:} & 	K_k = \ &  1. \qquad \text{(No filter)}
	\end{align*}
	
	\textbf{Simulation results:} The maximal and optimal approaches require that a past iteration exists to allows the estimation of $\utilde{\text{\hskip 0.1ex $\nabla$}}^S_{k}$. However, at the first iteration it does not exist. So, one must  predefine a filter to use for the first iteration. In this example, and in the rest of this thesis one uses $K_o=0.1$. The simulation results of these two approaches are given on Figure~\ref{fig:2___4_Maximaliste_vs_Optimaliste}. \\
	
	One can see that the optimal approach converges to the optimum of the plant in 2 iterations while finding the optimal value of the filter. The maximal approach is much slower as it requires 10 iterations (or 45 if a very precise estimate of $u_p^{\star}$ is desired). \\
	\begin{minipage}[h]{\linewidth}
		\vspace*{0pt}
		\centering
		\includegraphics[width=4.45cm]{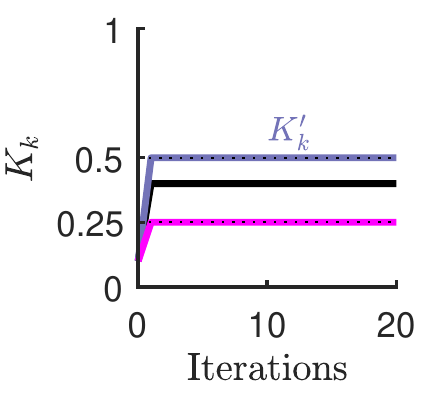}\hskip -0ex
		\includegraphics[width=4.45cm]{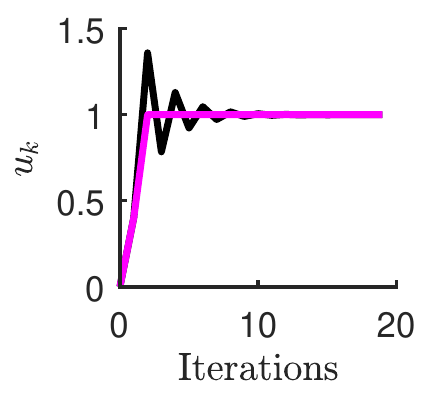}\hskip -0ex
		\includegraphics[width=4.45cm]{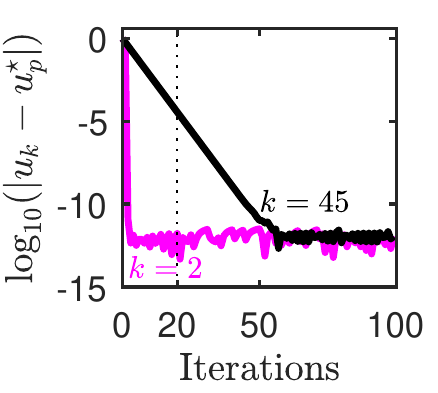}
		\captionof{figure}{Comparison of maximal (in \textbf{black} and \textcolor{blue_clair}{\textbf{blue}}) and  optimal (in \textcolor{magenta}{\textbf{magenta}}) approaches. The blue curve shows what the maximal approach identifies as the limit value of the filter.}
		\label{fig:2___4_Maximaliste_vs_Optimaliste}
	\end{minipage}\\

	\textbf{Remark on the management of the convergence:} One has noticed that when the iterates $u_k$ are very close to the solution $u_p^{\star}$, the maximal and optimal approaches end up choosing filters  $K_k$ in a quasi-random way. In fact, it is the computation of $\utilde{\text{\hskip 0.1ex $\nabla$}}^S_{k}$ with \eqref{eq:2___41_strat1_adaptationK}  which becomes erratic because the fraction  $(u^*_{k+1}-u^*_{k})/(||u_{k}-u_{k-1}||)$ becomes a division of numbers so small that their values is more or less a digital noise related to the precision of the computer. Not managing this can lead to unexpected behaviors. Therefore, one suggests to handle this practical issue with the following improve filter selector: 
	\begin{align}
		& \text{If } k =0 : & & K_k =  K_o, \nonumber \\
		& \text{If } \utilde{\text{\hskip 0.1ex $\nabla$}}^S_{k} < 1 \ \text{ and } \  k > 0  \ \text{ et } \  |u_{k}-u_{k-1}| \leq a: & &  K_k = K_{k-1}, \nonumber \\
		& \text{If } \utilde{\text{\hskip 0.1ex $\nabla$}}^S_{k} < 1 \ \text{ and } \  k > 0  \  \text{ et } \  |u_{k}-u_{k-1}|>a: & &
		K_k = \eqref{eq:2___gcerbgfj} \text{ or }
		\text{\eqref{eq:2___36_OptimizationPB_ChoiceK}} \nonumber \\
		& \text{Otherwise:}  & & K_k =  1, \ \text{(No filter)}.
		\label{eq:2___39_GestionConvergence}
	\end{align}
where 
\begin{equation*}
	a:=  \text{solver's precision}\times100 = 10^{-12}\times 100 = 10^{-10}.
\end{equation*}
\textit{(To obtain the accuracy of the solver one just needs to launch a simulation where instead of optimizing the plant, one optimizes a modified version of the model.  One will then see around which values the of iteration size  $\| u_{k} - u_{k-1} \|$ converge. In this simulation they converge to $\approx10^{-12}$). }
	\begin{center}
		\textbf{Conclusion}
	\end{center}
	
	This example shows that contrary to the insight one may have, filtering alone does not necessarily imply a decrease in the speed of convergence.  Indeed, it has been shown in a mathematical study that a well-chosen filter $K=0.25$ can enable much faster converge than a larger filter $K=0.4$ (and, based on several other tests not presented here, this observation is independent of the starting point).
\end{exbox}

	\begin{exbox} \label{ex:2___2_Kk_nu_geq2}
		\textbf{(Is $\bm{K_k^\star}$ still an appropriate choice when $n_u$>$1$?)} One considers the following generic RTO problem: 
		\begin{align*}   
			\phi(u,y) := \ & y, \ & 
			y_p = f_p(u) := \ & \frac{(\bm{u}-\bm{1}_{n_u})^{\rm T} \bm{A} (\bm{u}-\bm{1}_{n_u})}{2},  \    & 
			y = f(u)  :=  \ & \frac{\bm{u}^{\rm T} \bm{u}}{4},
		\end{align*}
		where $\bm{1}_{n_u} := [1,...,1]^{\rm T}\in \amsmathbb{R}^{n_u}$, and where clearly $\bm{u}_p^{\star} = \bm{1}_{n_u}$. Based on simulation results, the behavior emerging from the use of the filter update strategy \eqref{eq:2___36_OptimizationPB_ChoiceK}-\eqref{eq:2___39_GestionConvergence} is analyzed empirically. In the first part, one considers the case $n_u = 2$, and in second part  one considers the case $n_u = 3$. 
		
		\begin{center}
			\textbf{Part 1: ``$\bm{n_u = 2}$''.}
		\end{center}
		
		One sets:
		\begin{equation*}
			\bm{A} = \left(\begin{array}{cc}
				2 & 0 \\
				0 & 1
			\end{array}
			\right),
		\end{equation*}
		and run a simulation using ISO with the filter update strategy  \eqref{eq:2___36_OptimizationPB_ChoiceK}-\eqref{eq:2___39_GestionConvergence}. Figure~\ref{fig:2___9_ResultsSimEtude2D} shows the simulation results.\\
		\begin{minipage}[h]{\linewidth}
			\vspace*{0pt}
			\centering
			\includegraphics[width=4.45cm]{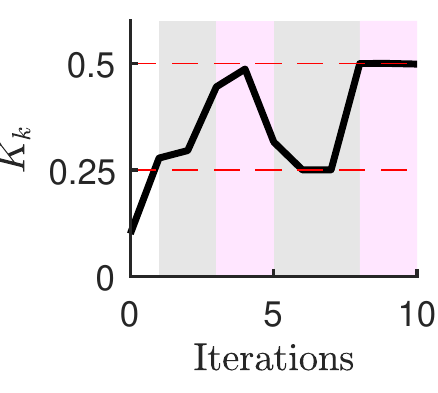}\hskip -0ex
			\includegraphics[width=4.45cm]{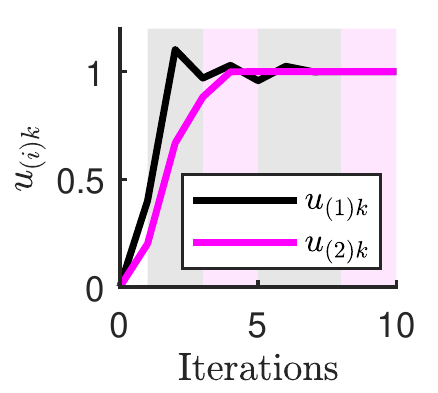}\hskip -0ex
			\includegraphics[width=4.45cm]{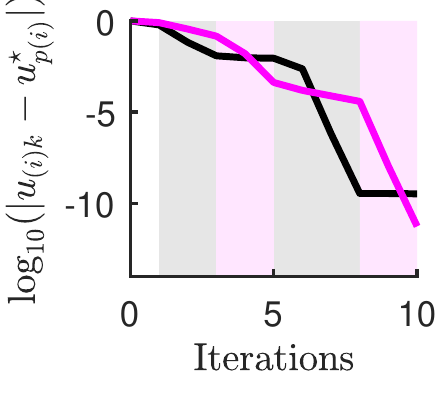}
			\captionof{figure}{Simulation results with the optimal filter adaptation.}
			\label{fig:2___9_ResultsSimEtude2D}
		\end{minipage}\\
	
		One can see that when the value of the filter is close to $1/2$ (highlighted with a magenta background), the distance $\|u_{(2)}-u_{p(2)}^\star\|$ decreases significantly.  On the contrary, when it is close to $1/4$ (highlighted with a gray background), it is the distance $\|u_{(1)}-u_{p(1)}^\star\|$ that decreases significantly. Moreover, one can observe that the filter is systematically close to either $1/2$ or $1/4$. 
		These two values $\{1/4,\ 1/2\}$ are not meaningless as they correspond to the optimal filter values when working only in the subspaces $\bm{u} = \bm{u}_p^{\star} + \alpha[1, 0]^{\rm T}$ and $\bm{u} = \bm{u}_p^{\star} + \alpha[0, 1]^{\rm T}$ with $\alpha\in\amsmathbb{R}$.
		\begin{itemize}
			\item In the subspace $\bm{u} = \bm{u}_p^{\star} + \alpha[1, 0]^{\rm T}$ with $\alpha\in\amsmathbb{R}\backslash 0$:
			\begin{align*}
				&  \alpha[1, 0]\nabla_{u} \bm{sol}|_{\bm{u}_p^{\star}}  = 1- \frac{ 
					\alpha[1, 0] \nabla^2_{\bm{uu}}\phi_p|_{\bm{u}_p^{\star}} \alpha[1, 0]^{\rm T} 
				}{
					\alpha[1, 0] \nabla^2_{\bm{uu}}\phi|_{\bm{u}_p^{\star}} \alpha[1, 0]^{\rm T} 
				} = 1- \frac{2}{1/2} = -3, \\
				&  \Rightarrow \ K^*_k  = \min\left( 1,\  \frac{1}{1- (-3) }\right) = \frac{1}{4}. 
			\end{align*}
			\item In the subspace $\bm{u} = \bm{u}_p^{\star} + \alpha[0, 1]^{\rm T}$ with $\alpha\in\amsmathbb{R}\backslash 0$:
			\begin{align*}
				& \alpha[0, 1]\nabla_{u} \bm{sol}|_{\bm{u}_p^{\star}}  = 1- \frac{ 
					\alpha[0, 1] \nabla^2_{\bm{uu}}\phi_p|_{\bm{u}_p^{\star}} \alpha[0, 1]^{\rm T} 
				}{
					\alpha[0, 1] \nabla^2_{\bm{uu}}\phi|_{\bm{u}_p^{\star}} \alpha[0, 1]^{\rm T} 
				} = 1- \frac{1}{1/2} = -1, \\
				& \Rightarrow \ K^*_k  = \min\left( 1,\  \frac{1}{1- (-1) }\right) = \frac{1}{2}.
			\end{align*}
		\end{itemize}
	
		Therefore, it seems  that the use of the strategy \eqref{eq:2___36_OptimizationPB_ChoiceK}-\eqref{eq:2___39_GestionConvergence} implicitly implies that the solving of the 2-D problem ($n_u=2$) is done through a process that transforms this problem into a succession of 1-D problems in subspaces comparable to $\bm{u} = \bm{u}_p^{\star} + \alpha[1, 0]^{\rm T}$ et $\bm{u} = \bm{u}_p^{\star} + \alpha[0, 1]^{\rm T}$ with $\alpha\in\amsmathbb{R}$. \\
		
		On Figure~\ref{fig:2___1-_ResultatsSim} those results are compared with the ones obtained by applying several fixed filters $K=\{0.25, \ 0.33, \ 0.5\}$. The facts are clear, \eqref{eq:2___36_OptimizationPB_ChoiceK}-\eqref{eq:2___39_GestionConvergence} outperforms any fixed filter strategy. Note that the fixed filter $K=0.33$ has been selected by hand and is close to the best you can get with a fixed filter for this problem.  \\
		\begin{minipage}[h]{\linewidth}
			\vspace*{0pt}
			\centering
			\includegraphics[width=4.45cm]{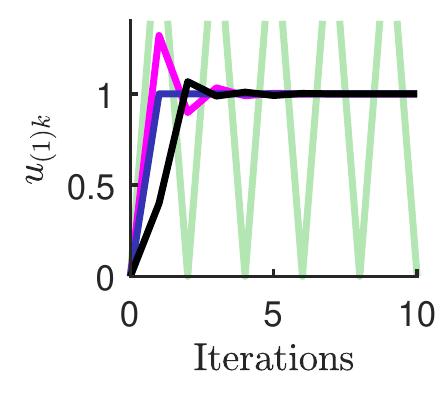}\hskip -0ex
			\includegraphics[width=4.45cm]{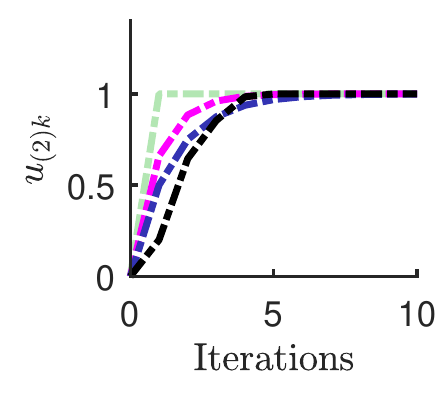}\hskip -0ex
			\includegraphics[width=4.45cm]{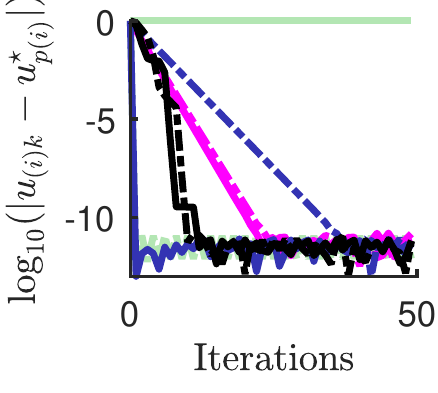}
			\captionof{figure}{The optimal approach (in \textbf{black}) is better than any fixed filter approach (in \textcolor{vert_fonce}{\textbf{green}} for $K=0.5$, \textcolor{blue}{\textbf{blue}} for $K=0.25$, and \textcolor{magenta}{\textbf{magenta}} for $K=0.33$).}
			\label{fig:2___1-_ResultatsSim}
		\end{minipage}\\

		\textit{(To go further: By varying the values of the elements of the diagonal of $\bm{A}$ comparable results are obtained. However, if the matrix $\bm{A}$ is not diagonal, i.e. if a rotation is applied to it, then the results will remain  similar but the value of $K_k^{\star}$ will oscillate between the values of the optimal filters in subspaces  $\bm{u} = \bm{u}_p^{\star} +\alpha \bm{v}_1$ et $\bm{u} = \bm{u}_p^{\star} +\alpha \bm{v}_2$, where $\alpha\in\amsmathbb{R}$, and where $(\bm{v}_1, \bm{v}_2)$ are the eigenvectors of $\bm{A}$.)}
		
		\begin{center}
			\textbf{Part 2: ``$\bm{n_u = 3}$''.}
		\end{center}
		
		One now considers a 3-D case, and one wants to check whether the empirical observations made on the 2-D problems remain applicable. One sets:
		\begin{equation*}
			\bm{A} = \left(\begin{array}{ccc}
				2 & 0 & 0\\
				0 & 1 & 0\\
				0 & 0 & 2/3
			\end{array}
			\right),
		\end{equation*}
		and run a simulation using ISO with the filter update strategy  \eqref{eq:2___36_OptimizationPB_ChoiceK}-\eqref{eq:2___39_GestionConvergence}.  Figure~\ref{fig:2___11_ResulatsSim} shows the simulation results. 
		\\
		\begin{minipage}[h]{\linewidth}
			\vspace*{0pt}
			\centering
			\includegraphics[width=4.45cm]{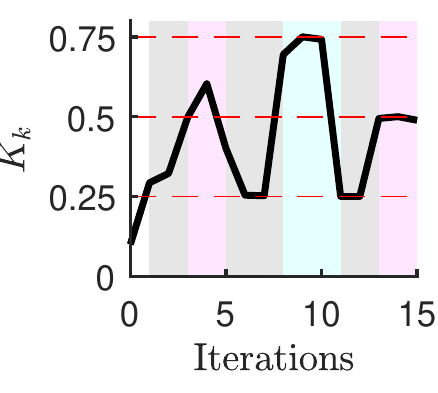}\hskip -0ex
			\includegraphics[width=4.45cm]{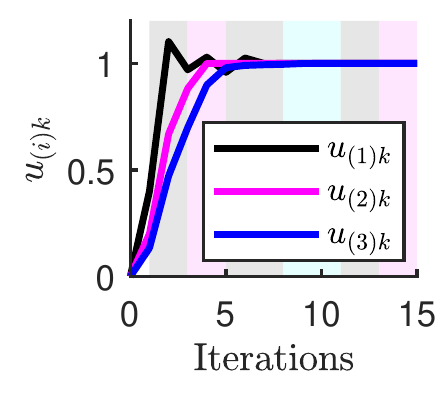}\hskip -0ex
			\includegraphics[width=4.45cm]{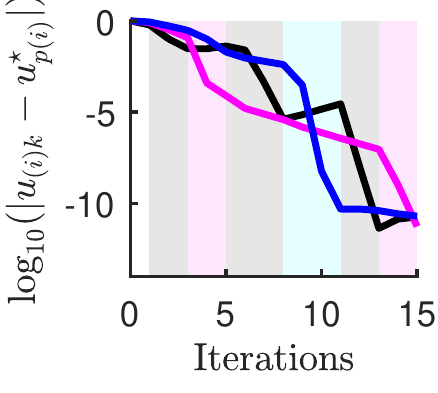}
			\captionof{figure}{Simulation results using the optimal filter adaptation.}
			\label{fig:2___11_ResulatsSim}
		\end{minipage}\\
	
		The analysis of these results is similar to the analysis done for the 2-D case. Therefore, one can skip the technical details and go directly to the calculation of the optimal filters:
		\begin{itemize}
			\item In the subspace $\bm{u} = \bm{u}_p^{\star} + \alpha[1, 0, 0]^{\rm T}$ with $\alpha\in\amsmathbb{R}\backslash 0$:  $K^*_k  = 1/4$.
			\item In the subspace $\bm{u} = \bm{u}_p^{\star} + \alpha[0, 1, 0]^{\rm T}$ with $\alpha\in\amsmathbb{R}\backslash 0$: $K^*_k  = 1/2$.	
			\item In the subspace $\bm{u} = \bm{u}_p^{\star} + \alpha[0, 0, 1]^{\rm T}$ with
			 $\alpha\in\amsmathbb{R}\backslash 0$: $K^*_k  = 3/4$.
		\end{itemize}
		One finds the 3 levels on Figure~\ref{fig:2___11_ResulatsSim}. This validates the extension of the empirical analysis from 2-D cases to higher dimensional cases.  Moreover, as for part 1, the strategy  \eqref{eq:2___36_OptimizationPB_ChoiceK}-\eqref{eq:2___39_GestionConvergence} is superior to any classical strategy which consists in fixing the filter  $K$  from the start (see Figure~\ref{fig:2___12_ResultatsSim}). The results obtained with the adaptive filter are in black and those obtained with the fixed filter $K=0.34$ (manually identified as the ``best'' fixed filter) are in magenta. 
		
		\begin{center}
			\textbf{Conclusion}
		\end{center}
	
		The filter update strategy \eqref{eq:2___36_OptimizationPB_ChoiceK}-\eqref{eq:2___39_GestionConvergence} leads to the emergence of interesting behaviors when the problem dimension is greater than 1. In such cases, the value of the filter appears to alternate between its optimal values in sub-dimensions of the $n_u$-dimensional problem. 
		
		\textit{(Larger dimensions with matrices $\bm{A}$ which are non-diagonal lead to the same observations on testing.)}
		
		\begin{minipage}[h]{\linewidth}
			\vspace*{0pt}
			\centering
			\includegraphics[width=4.45cm]{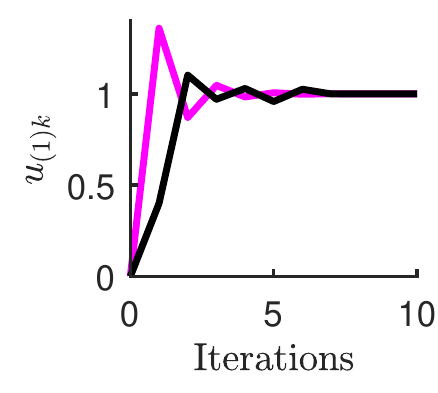}\hskip -0ex
			\includegraphics[width=4.45cm]{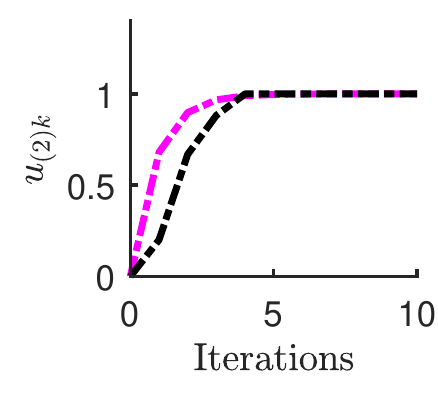}\hskip -0ex
			\includegraphics[width=4.45cm]{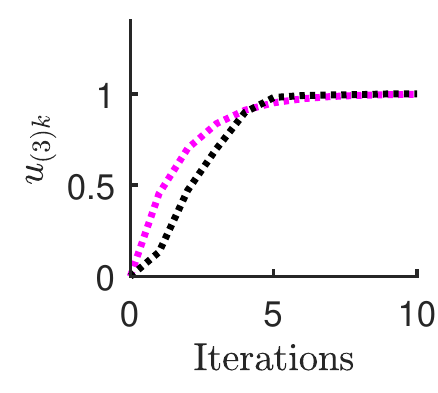} 
			
			\includegraphics[width=4.45cm]{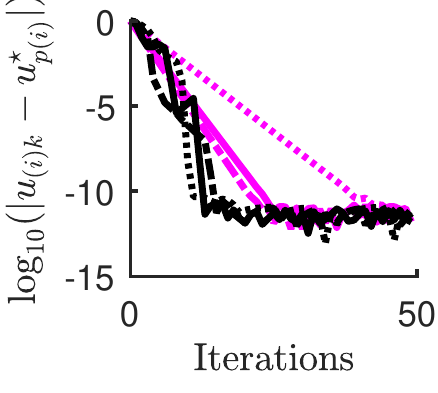}
			\captionof{figure}{The adaptive approach (in \textbf{black}) is better than any fixed filter approach (whose best results are plotted in \textcolor{magenta}{\textbf{magenta}}).}
			\label{fig:2___12_ResultatsSim}
		\end{minipage}
	\end{exbox}

\begin{Remark} \textbf{(Should one use an upper bound for the filter: $K<1$?)} \label{rem:2___1_UtiliteBorneSupK}
	When one uses the strategy \eqref{eq:2___36_OptimizationPB_ChoiceK}-\eqref{eq:2___39_GestionConvergence} it may happen that: 
		\begin{equation}\label{eq:2___36_OptimizationPB_ChoiceK_2}
		 \max\left\{0.1, \frac{1}{1- \utilde{\text{\hskip 0.1ex $\nabla$}}^S_{k}} \right\} > 1. 
	\end{equation}
	This fact calls into question the classical approach (of MA and MAy) of choosing a filter between $0$ and $1$. This limit can easily be imposed by replacing \eqref{eq:2___36_OptimizationPB_ChoiceK} by 
	\begin{align} \label{eq:2___50_JCBJSRCFHYUDFCR}
		K^{\star}_k = \ & \min\left\{1,\  \max\left\{K_o,   \frac{1}{1- \utilde{\text{\hskip 0.1ex $\nabla$}}^S_{k}} \right\}\right\}, &
		\text{si } \utilde{\text{\hskip 0.1ex $\nabla$}}_k^S < 1. 
	\end{align}
	In the following example, the strategy \eqref{eq:2___36_OptimizationPB_ChoiceK}-\eqref{eq:2___39_GestionConvergence} is compared to the one using an upper bound on the filter \eqref{eq:2___39_GestionConvergence}-\eqref{eq:2___50_JCBJSRCFHYUDFCR} on a case where \eqref{eq:2___36_OptimizationPB_ChoiceK_2} holds. The results clearly show the superiority of the unbounded approach.  
\end{Remark}

\begin{exbox} \textbf{(Comparison of \eqref{eq:2___36_OptimizationPB_ChoiceK}-\eqref{eq:2___39_GestionConvergence} versus  \eqref{eq:2___39_GestionConvergence}-\eqref{eq:2___50_JCBJSRCFHYUDFCR})}
		One takes the Example~\ref{ex:2___1_Exemple1} and one changes the model $f$:
		\begin{align*}   
			\phi(u,y) := \ & y, \ & 
			y_p = f_p(u) := \ & (u-1)^2,  \    & 
			y = f(u)  :=  \ & 5u^2.
		\end{align*}
	This change of model has the effect of changing the functions $sol$ associated to this RTO problem which becomes the one represented on the Figure~\ref{fig:2___13_exemple_Ksup1_FonctionSOL}.
	\begin{minipage}[h]{\linewidth}
		\smallskip
		\centering
		\includegraphics[width=4.45cm]{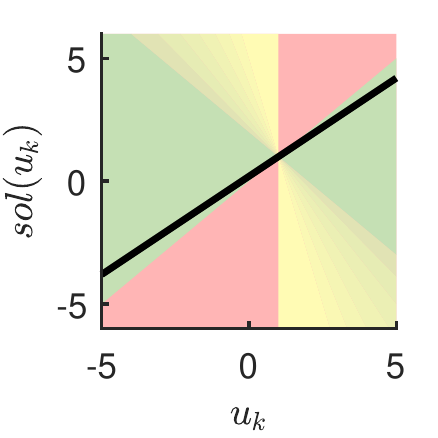} 
		\captionof{figure}{The function $sol(u_k)$.}
		\label{fig:2___13_exemple_Ksup1_FonctionSOL} 	
		\smallskip
	\end{minipage} \\

Now, let's implement ISO with the strategies \eqref{eq:2___36_OptimizationPB_ChoiceK}-\eqref{eq:2___39_GestionConvergence} and \eqref{eq:2___39_GestionConvergence}-\eqref{eq:2___50_JCBJSRCFHYUDFCR} and let's see the results: see Figure~\ref{fig:2___14_ResulatsSim}.\\
	
	These results show very clearly that the strategy  \eqref{eq:2___36_OptimizationPB_ChoiceK}-\eqref{eq:2___39_GestionConvergence} leads to a much faster convergence than \eqref{eq:2___39_GestionConvergence}-\eqref{eq:2___50_JCBJSRCFHYUDFCR}, in 2 iterations against 111 for a precise convergence and 2 iterations against 25 for an "approximate" convergence. Thus, this simulation supports the idea of dropping the upper bound on the filter $K$. \\
	\begin{minipage}[h]{\linewidth}
		\vspace*{0pt}
		\centering
		\includegraphics[width=4.45cm]{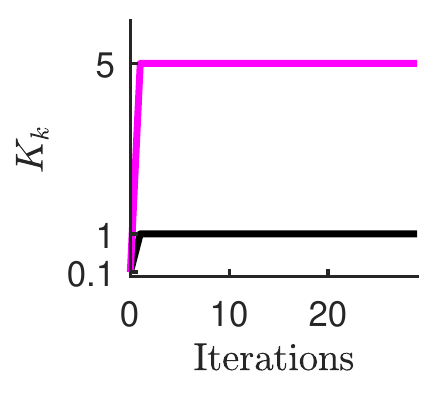}\hskip -0ex
		\includegraphics[width=4.45cm]{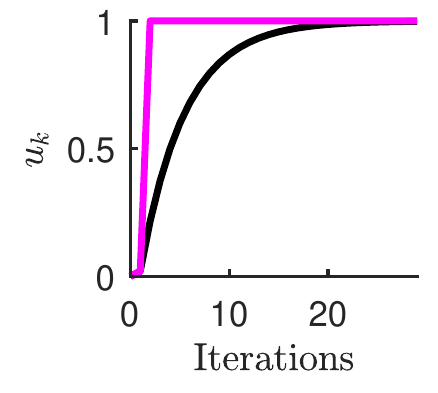}\hskip -0ex
		\includegraphics[width=4.45cm]{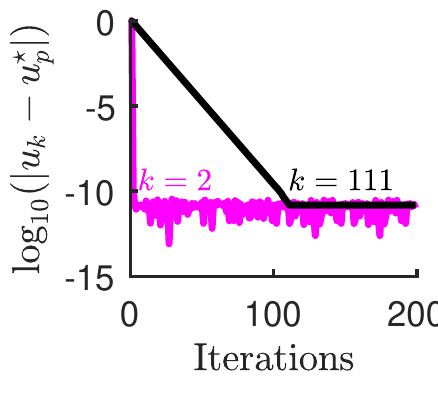}
		\captionof{figure}{Results of the strategy \eqref{eq:2___39_GestionConvergence}-\eqref{eq:2___50_JCBJSRCFHYUDFCR} (in \textbf{black}) and of the strategy \eqref{eq:2___36_OptimizationPB_ChoiceK}-\eqref{eq:2___39_GestionConvergence} (in \textcolor{magenta}{\textbf{magenta})}. }
		\label{fig:2___14_ResulatsSim}
	\end{minipage}
	\end{exbox}

At this stage one has shown that the equilibrium condition of ISO-D/I can be enforced with the use of locally convex models. One has also shown that the stability can be enforced with the use of a filter on the inputs \eqref{eq:2___27b_Filter} whose value can be automatically tuned with the procedure \eqref{eq:2___39_GestionConvergence}-\eqref{eq:2___50_JCBJSRCFHYUDFCR} or \eqref{eq:2___36_OptimizationPB_ChoiceK}-\eqref{eq:2___39_GestionConvergence}. Let us now move on to the analysis of the superstability.

\section{Enforcing the superstability condition}

To determine the superstability condition, one must go back to the beginning of the stability analysis. In fact, the only thing one has to do is to start from  \eqref{eq:2___27b_Filter} and: 
\vspace{-\topsep}
\begin{itemize}[noitemsep]
	\item replace $\bm{u}_{k+1}^*$ by \eqref{eq:2___28_C1_Uasteri};
	\item \textit{not neglect} term $\mathcal{O}(\|\delta \bm{u}_k\|^2)$ \textit{to not place the study only in the neighborhood} of $\bm{u}_p^{\star}$;
	\item add ``$-\bm{u}_p^{\star}$'' on both sides of the equality;
\end{itemize}
\vspace{-\topsep}
to find:
\begin{align}
	\bm{u}_{k+1} - 	\bm{u}^{\star}_{p} = \ & 
	\bm{u}_k + K (\bm{u}^{\star}_p + 
	\nabla^S(\delta\bm{u}_k) 
	(\bm{u}_k-\bm{u}_p^{\star}) + \mathcal{O}(\|\delta \bm{u}_k\|^2) - \bm{u}_k) - \bm{u}^{\star}_{p}, \nonumber \\
	= \ & 
	\left(1 - K \left( 1 - 
	\nabla^S(\delta\bm{u}_k) 
	- \mathcal{O}(\|\delta \bm{u}_k\|)
	\right)\right)
	(\bm{u}_{k} - 	\bm{u}^{\star}_{p} ).
	\label{eq:2___41_IterFiltreEffet_Global}
\end{align}
This equation is similar to \eqref{eq:2___30_IterFiltreEffet} but the difference is that it takes into account the higher order terms  $\mathcal{O}(\|\delta \bm{u}_k\|)$ of \eqref{eq:2___28_C1_Uasteri}. The superstability condition is then:
\begin{equation} \label{eq:2___42_ConditionNiveau3}
	-1 < 
	1 - K \left( 1 - 
\left(\nabla^S(\delta\bm{u}_k)
	 + \mathcal{O}(\|\delta \bm{u}_k\|) \right) \right) < 1, \qquad \forall \delta\bm{u}_k\in\amsmathbb{R}^{n_u}.
\end{equation}
At first glance, the term ``$\nabla^S(\delta\bm{u}_k)
+ \mathcal{O}(\|\delta \bm{u}_k\|)$'' may seem complicated to understand, but in fact it has a simple geometric meaning.  It is the ``slope \footnote{The tangent of the angle between the line passing through $[\bm{u}_p^{\star}{}^{\rm T}, \bm{sol}(\bm{u}_p^{\star})^{\rm T}]^{\rm T}$ and 
	$[(\bm{u}_p^{\star}+\delta\bm{u}_k)^{\rm T}, \bm{sol}(\bm{u}_p^{\star}+\delta\bm{u}_k)^{\rm T}]^{\rm T}$ and the line passing through $[\bm{u}_p^{\star}{}^{\rm T}, \bm{0}^{\rm T}]^{\rm T}$ and 
	$[(\bm{u}_p^{\star}+\delta\bm{u}_k)^{\rm T}, \bm{0}^{\rm T}]^{\rm T}$.
}'' of the line passing through the points $\bm{u}_p^{\star}$ and  $\bm{sol}(\bm{u}_p^{\star} + \delta \bm{u}_k)$. It is therefore proposed to replace it by a term $\Delta^S(\delta\bm{u}_k)$ which is less intimidating and which is defined as follows: 
\begin{equation*}
	\Delta^S(\delta\bm{u}_k) :=  \nabla^S(\delta\bm{u}_k) + \mathcal{O}(||\delta \bm{u}_k||)
	= 
	\left(\frac{\delta \bm{u}_k}{\|\delta \bm{u}_k\|}\right)^{\rm T} 
	\frac{\bm{sol}(\bm{u}_p^{\star} + \delta \bm{u}_k) - \bm{u}_p^{\star}}{\|\delta \bm{u}_k\|}
\end{equation*}
The condition \eqref{eq:2___42_ConditionNiveau3} can then be rewritten more simply as:
\begin{subequations}
	\label{eq:2___42_Niveau3}
	\begin{align}
		0 < K<\frac{2}{1-\Underline{\Delta}^S}, &&
		\text{ si: }  \Overline{\Delta}^S < 1,
	\end{align}
	where 
	\begin{align}
		\Underline{\Delta}^S := \ & \underset{\delta\bm{u}_k}{\operatorname{min}}\{\Delta^S(\delta\bm{u}_k)\}, & 
		\Overline{\Delta}^S := \ & \underset{\delta\bm{u}_k}{\operatorname{max}}\{\Delta^S(\delta\bm{u}_k)\}.
	\end{align}
\end{subequations}
As these calculations and reasoning are neither very intuitive, nor easy to visualize, one proposes the following graphic interpretation: 

\begin{GraphicalInterpretationbox}	\label{gi:2___2_Niveau3}
	Let's take the simplest possible case where $n_u=1$ and where $sol(u)$ is globally $\mathcal{C}^1$.
	Figure~\ref{fig:2___12_InterpretationGraphiqueNiveau3} shows what the function $sol$  could look like (the \textbf{black} curve). The green, yellow and red cones shown are the same as the ones used in the graphical interpretation~\ref{gi:2___1_Niveau2}. the values of  $\Underline{\Delta}^S$ and $\Overline{\Delta}^S<1$ are given graphically as the minimum and maximum slopes between $[u_p^{\star},u_p^{\star}]$ and all the points $[u,sol(u)]$ on the curve $sol$.  \\
	
	It should be clear that if all the points $[u,sol(u)]$ are in the green cone then superstability is guaranteed without requiring a filter. If all the points $[u,sol(u)]$ are in the green and yellow cones then superstability is guaranteed only if the filter is small enough. To understand this at a geometrical level, one must associate Figures~\ref{fig:2___4_AllCones} (image on the right) and \ref{fig:2___12_InterpretationGraphiqueNiveau3}. 
	
	\begin{minipage}[h]{\linewidth}
		\vspace*{0pt}
		\centering 
		\includegraphics[width=13.35cm]{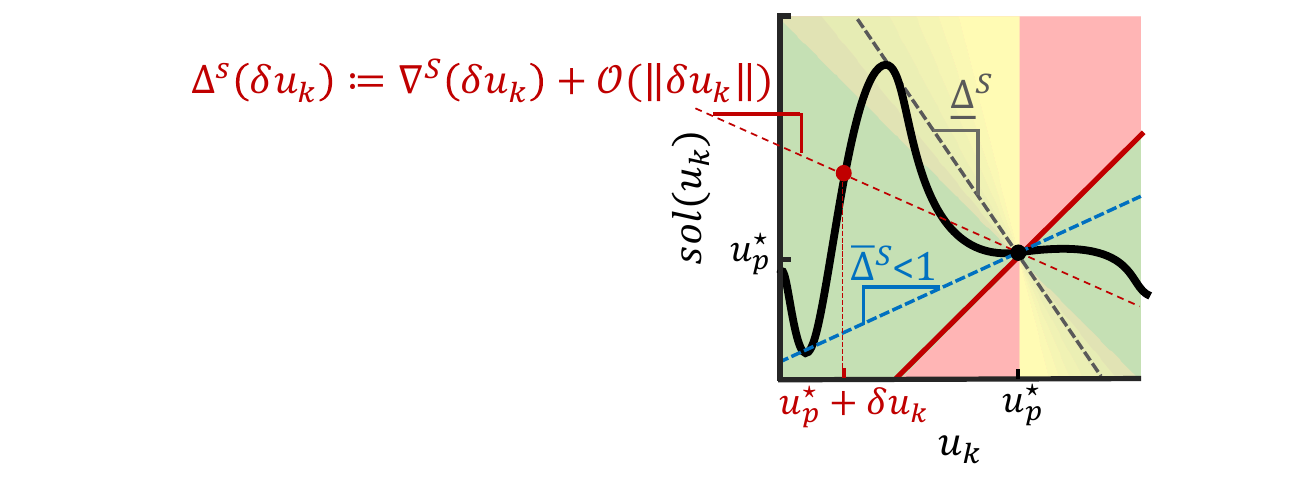}\
		\captionof{figure}{Graphical interpretation of the superstability condition}
		\label{fig:2___12_InterpretationGraphiqueNiveau3}
	\end{minipage}
\end{GraphicalInterpretationbox}	

Unfortunately, it is impossible to guarantee \eqref{eq:2___42_Niveau3} without questioning what was proposed to guarantee the equilibrium (and thus indirectly the stability). Indeed, forcing the updated model to be convex at any correction point makes the plant's minimums, maximums and saddle points fixed points of ISO-D/I. In other words, if $\bm{u}_p^{\bullet}$ is a stationary point -- a minimum, maximum or saddle point -- of the plant, then $\bm{sol}(\bm{u}_p^{\bullet}) = \bm{u}_p^{\bullet}$.  So, if the plant has several stationary points, then it is certain that $\Overline{\Delta}^S\geq 1$ and therefore the condition \eqref{eq:2___42_Niveau3} cannot be satisfied. Nevertheless, even though the plant maximums and saddle points are fixed points of the improved versions of ISO-D and ISO-I, they are ``\textit{unstable}'' fixed points. Indeed, the Theorem~\ref{thm:2___3_Proprietes_Sol} through the equations~\eqref{eq:2___29_NablaS_Signe} shows that the function $\bm{sol}$  around such points has directional derivatives $\geq 1$, and one knows that stability is only provided if they are $<1$. So, such points only satisfy the equilibrium conditions, and any iteration close to them, but not equal to them, is likely to lead to a sequence of iterations that diverge, despite the use of local convexifications and the filter. This instability is illustrated in the following example:

	\begin{exbox}\label{ex:2__4_PtsStable_Et_Instables}
		Let be an unconstrained RTO problem where $n_u=1$, where the cost function is 
		\begin{align*}   
			\phi(u,y) := \ & y,
		\end{align*}
		and where the plant and the model are: 
		\begin{align*} 
			y_p = f_p(u) := \ &  \frac{1}{4}-\frac{1}{4}u^2+\frac{1}{4}u^3+\frac{1}{16}u^4 - \frac{3}{20}u^5 + \frac{1}{24}u^6, & 
			\ \ \ y = f(u)  :=  \ & \frac{u^2}{4}.
		\end{align*}
		Figure~\ref{fig:2___13_Exemple2_Fonctions} shows what these curves look like and one can observe that this plant has two minimums ($u_p^{\bullet}=\{-1,\ 2\}$), one saddle-point ($u_p^{\bullet}=1$) and one maximum ($u_p^{\bullet}=0$). Moreover, one can see that the model is convex and therefore satisfies the equilibrium condition.   
		
		\begin{minipage}[h]{\linewidth}
			\vspace*{0pt}
			\centering 
			\includegraphics[width=8cm]{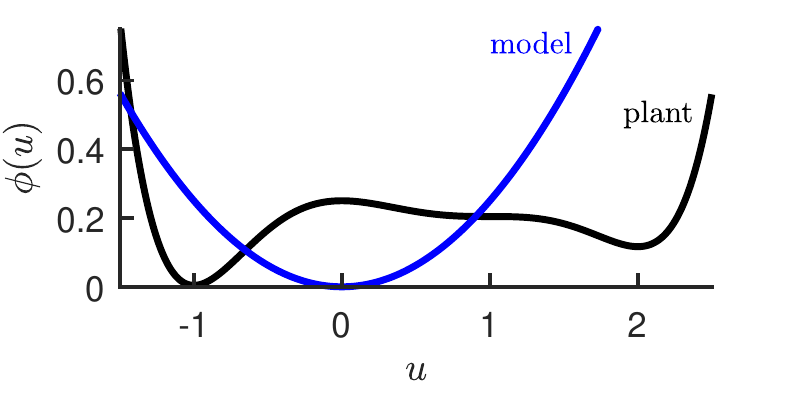}
			\captionof{figure}{The cost functions.}
			\label{fig:2___13_Exemple2_Fonctions}
		\end{minipage}\\
	
		The model updated at a point $u_k$ is:
		\begin{equation*}
			\phi_k(u) = \phi_p(u_k) + \nabla_u\phi_p|_{u_k}(u-u_k) + \frac{1}{4}(u-u_k)^2.
		\end{equation*} 
		Figure~\ref{fig:2___14_Exemple2_LocalCone} gives the solution $sol(u_k)$ of the optimization problem based on this model updated at $u_k$ for all $u_k$. In addition to this curve, the green, yellow, and red cones of the Graphical Interpretation~\ref{gi:2___1_Niveau2} have been added at the four fixed points (where $sol(u_k)=u_k$) so that it is possible to visually judge whether their are stable.  Clearly, only the minimums of the plant $u_p^{\star}=\{-1,\ 2\}$ meet the necessary conditions of stability ($sol(u_k)$ is locally in the green and yellow areas, so there is a filter that stabilizes these points). On the other hand, any point $u_k$ in the neighborhood of the maximum of the plant $u_p^{\bullet}=0$ leads to a $u_{k+1}$ which goes away from it (because $sol(u_k)$ is locally in the red cone). Concerning the saddle point $u_p^{\bullet}=1$  one can see that the function $sol$ is tangential to the boundary between the red and green cones at this point.  On the left side it is tangent but in the green domain while on the right side it is tangent but in the red domain. This means that the points on the left of $u_p^{\bullet}=1$ would lead to a step towards $u_p^{\bullet}=1$, and the points on the right of $u_p^{\bullet}=1$ would lead to a step away from $u_p^{\bullet}=1$. \\
		
		\begin{minipage}[h]{\linewidth}
			\vspace*{0pt}
			\centering 
			\includegraphics[width=9cm]{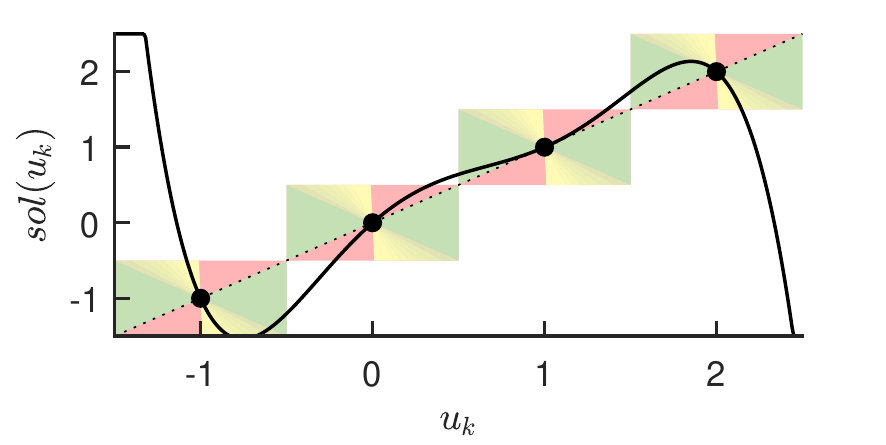}
			\captionof{figure}{Graphical analysis of the convergence capabilities on various types of fixed points of the improved versions of ISO-D/I.}
			\label{fig:2___14_Exemple2_LocalCone}
		\end{minipage}\\
	\end{exbox}

So, as one can see with the Example~\ref{ex:2__4_PtsStable_Et_Instables}, the validity of the condition  $\Overline{\Delta}^S\geq 1$ depends on the number of stationary points in the plant and this number cannot be manipulated. If one decides to ignore the condition $\Overline{\Delta}^S\geq 1$, there remains only the condition:  $0 < K< 2/(1-\Underline{\Delta}^S)$. 

The value of $\Underline{\Delta}^S$ is in practice not accessible because a priori neither $\bm{u}_p^{\star}$ nor $\bm{sol}(\bm{u})$ are known, $\forall \bm{u}\in\amsmathbb{R}^{n_u}$. However, an idea to estimate $\Underline{\Delta}^S$ could be the following:
\begin{equation} \label{eq:2___43_Une_Borne_Superieur_Pour_Niveau3}
	\Underline{\Delta}^S \approx \utilde{\text{\hskip 0.1ex $\Delta$}}^S_k
	:= \underset{\ell=1,...,k}{\operatorname{min}} \left\{
	\left(\frac{\bm{u}_{\ell} - \bm{u}_{\ell-1}}{\|\bm{u}_{\ell} - \bm{u}_{\ell-1}\|}\right)^{\rm T} 
	\frac{\bm{sol}(\bm{u}_{\ell}) - \bm{sol}(\bm{u}_{\ell-1})}{\|\bm{u}_{\ell} - \bm{u}_{\ell-1}\|} \right\}.
\end{equation}
Basically, $\utilde{\text{\hskip 0.1ex $\Delta$}}^S_k$ is the smallest slope observed between two consecutive points until iteration $k$. Of course, there is no guarantee that  $\utilde{\text{\hskip 0.1ex $\Delta$}}^S_k \leq \Underline{\Delta}^S$, but selecting the filter in this way can only increase the chances of convergence. So one can choose:
\begin{align}
	K_k = \ & 
		\max\left\{0.1, \min\left\{  \frac{1}{1- \utilde{\text{\hskip 0.1ex $\nabla$}}^S_{k}}, \ \frac{2}{1- \utilde{\text{\hskip 0.1ex $\Delta$}}^S_k } \right\}\right\}  & & 
		\text{if: } 
		\utilde{\text{\hskip 0.1ex $\Delta$}}^S_k < 1 \text{ and } 
		\utilde{\text{\hskip 0.1ex $\nabla$}}^S_{k} <1, \nonumber \\
		= \ &  \max\left\{0.1,  \  \frac{2}{1- \utilde{\text{\hskip 0.1ex $\Delta$}}^S_k } \right\} & &
		\text{if: }
		\utilde{\text{\hskip 0.1ex $\Delta$}}^S_k < 1. \label{eq:2___44_OptimizationPB_ChoiceK_AvecAm}
\end{align}
or
\begin{align}
	K_k = \ & 
	\max\left\{0.1,\  \min\left\{ 1, \ \frac{1}{1- \utilde{\text{\hskip 0.1ex $\nabla$}}^S_{k}}, \ \frac{2}{1- \utilde{\text{\hskip 0.1ex $\Delta$}}^S_k } \right\}\right\}  & & 
	\text{if: } 
	\utilde{\text{\hskip 0.1ex $\Delta$}}^S_k < 1 \text{ and } 
	\utilde{\text{\hskip 0.1ex $\nabla$}}^S_{k} <1, \nonumber \\
	= \ &  \max\left\{0.1,  \ \min\left\{ 1, \ \frac{2}{1- \utilde{\text{\hskip 0.1ex $\Delta$}}^S_k } \right\} \right\} & &
	\text{if: }
	\utilde{\text{\hskip 0.1ex $\Delta$}}^S_k < 1. \label{eq:2___44_OptimizationPB_ChoiceK_AvecAm_2}
\end{align}
However, doing so is not really necessary because the counterpart associated with this attempt to enforce superstability is a decrease in the positive effects of the adaptive filter illustrated in Example~\ref{ex:2___2_Kk_nu_geq2}. The following example illustrates this last point:
\begin{exbox} 
		If one restarts from the case $n_u=3$ of  Example~\ref{ex:2___2_Kk_nu_geq2}, and  instead of applying the strategy  \eqref{eq:2___36_OptimizationPB_ChoiceK}-\eqref{eq:2___39_GestionConvergence}, \eqref{eq:2___39_GestionConvergence}-\eqref{eq:2___44_OptimizationPB_ChoiceK_AvecAm} is applied. Then the results of Figure~\ref{fig:2____15_Kk_AvecAm} are obtained. On the left the magenta curve  (\textcolor{magenta}{\rule{0.5cm}{0.1cm}}) gives the  $K_k$  obtained with  \eqref{eq:2___39_GestionConvergence}-\eqref{eq:2___44_OptimizationPB_ChoiceK_AvecAm} and those obtained with \eqref{eq:2___36_OptimizationPB_ChoiceK}-\eqref{eq:2___39_GestionConvergence} are plotted in black (\rule{0.5cm}{0.1cm}). The gray area (\hskip 0.4ex\textcolor{gray}{\rule{0.5cm}{0.3cm}}) gives the upper bound for the choice of $K_k$, i.e. $2/(1-\utilde{\text{\hskip 0.1ex $\Delta$}}^S_k)$. The curve on the right gives the evolution of the distance between the iterates $(u_{(1)k}, u_{(2)k}, u_{(3)k})$ and the optimum $\bm{u}_p^{\star}$. Unlike in the previous example, the different $u_{(i)k}$ are not discriminated since this would not have much interest here.  
		 Clearly, one can see that the addition of the boundary  $K_k\leq2/(1-\utilde{\text{\hskip 0.1ex $\Delta$}}^S_k)$ doubles the convergence speed from 15 to 30 iterations.  
		\\
		\begin{minipage}[h]{\linewidth}
			\vspace*{0pt}
			\centering
			\includegraphics[width=4.45cm]{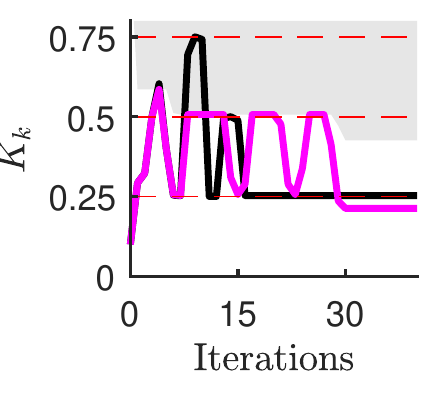}\hskip -0ex
			\includegraphics[width=4.45cm]{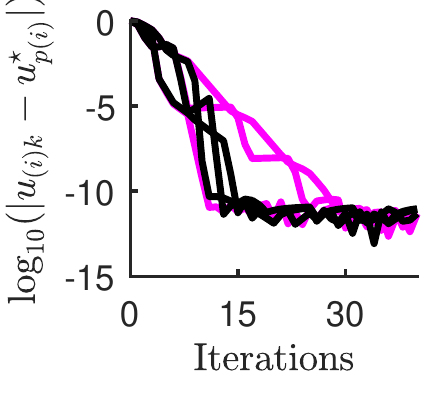}
			\captionof{figure}{Strategy \eqref{eq:2___36_OptimizationPB_ChoiceK}-\eqref{eq:2___39_GestionConvergence} vs. Strategy \eqref{eq:2___39_GestionConvergence}-\eqref{eq:2___44_OptimizationPB_ChoiceK_AvecAm}.}
			\label{fig:2____15_Kk_AvecAm}
		\end{minipage}\\
\end{exbox}

\section{Modifier-Filter-Curvature Adaptation} 

If one applies the set of improvements discussed in the previous three subsections to ISO-D/I, one obtains two new methods that guarantee that the equilibrium and stability conditions are always satisfied. These two methods are called \textit{modifier-filter-curvature adaptation direct (MFCA-D)} and  \textit{modifier-filter-curvature Adaptation indirect (MFCA-I)} which are described hereafter. 

\begin{BoxAlgo}{\textbf{Modifier-Filter-Curvature Adaptation Direct (MFCA-D)}}{MFCA}
	\textbf{Initialization.} Provide $\bm{u}_0$, $\bm{K}_0(=0.1)$, $\bm{u}^*_0(=\bm{u}_0)$, functions $(\bm{f},\phi,\bm{g})$, and the stopping criterion of step 6).
	\tcblower
	\textbf{for} $k=0 \rightarrow \infty$
	\begin{itemize}[noitemsep]
		\item[1) ]  \textbf{Measure} $(\nabla_{\bm{u}}\phi_p, \bm{g}_p,\nabla_{\bm{u}}\bm{g}_p)|_{\bm{u}_{k}}$ on the plant.
		\item[2) ] \textbf{Update} the functions $(\phi,\bm{g})$  with \eqref{eq:3___19_Corrections_ISO}.
		\item[3) ]  \textbf{Make} these cost and constraint functions convex at $\bm{u}_k$: \eqref{eq:3___30_DefinitionMatrices_P}.
		\item[4) ]  \textbf{Compute}  $\bm{u}^\star_{k+1}$ with \eqref{eq:2___27a_ModelBasedPB}.
		\item[5) ] \textbf{Compute} The filter with:
	\end{itemize}	 
		\begin{align} 
			K_k :=  \ &  \left\{
			\begin{array}{ll}
				K_o,     & \text{if } k =0, \\
				K_{k-1}  & \text{if } \utilde{\text{\hskip 0.1ex $\nabla$}}^S_{k} < 1 \ \text{ and } \  k > 0  \ \text{ and }  |\bm{u}_{k}-\bm{u}_{k-1}|\leq \bm{a}, \\
				\max\left\{K_o,\frac{1}{1-\utilde{\text{\hskip 0.1ex $\nabla$}}^S_{k}}\right\}
				&  \text{if } \utilde{\text{\hskip 0.1ex $\nabla$}}^S_{k} < 1 \ \text{ and } \  k > 0  \ \text{ and }    |\bm{u}_{k}-\bm{u}_{k-1}|> \bm{a}, \\
				1, & \text{otherwise.}
			\end{array}
			\right. 
			\label{eq:2___45_LE_FILTRE}\\
			\utilde{\text{\hskip 0.1ex $\nabla$}}_k^S  := \ &
			\left( \frac{\bm{u}_k - \bm{u}_{k-1}}{\|\bm{u}_k - \bm{u}_{k-1}\|}
			\right)^{\rm T}
			\frac{\bm{u}^{\star}_{k+1} - \bm{u}^{\star}_{k}}{\| \bm{u}_{k} - \bm{u}_{k-1} \| }. \label{eq:2___59_Nabla_S}
		\end{align}
	\vspace{-\topsep}
	\begin{itemize}[noitemsep]
		\item[6) ] \textbf{Compute}  $\bm{u}_{k+1} = \bm{u}_{k} + K_k(\bm{u}^*_{k+1}-\bm{u}_{k})$.
		\item[7) ] \textbf{Stop} if $\bm{u}_{k+1}\approx\bm{u}_{k}$ and return $\bm{u}_{\infty} := \bm{u}_{k+1}$.
	\end{itemize}
	\noindent {\bf end}
\end{BoxAlgo}

\begin{BoxAlgo}{\textbf{Modifier-Filter-Curvature Adaptation Indirect (MFCA-I)}}{MFCAy}
	\textbf{Initialization.}  Provide $\bm{u}_0$, $\bm{K}_0(=0.1)$, $\bm{u}^*_0(=\bm{u}_0)$, functions $(\bm{f},\phi,\bm{g})$, and the stopping criterion of step 6).
	\tcblower
	\textbf{for} $k=0 \rightarrow \infty$
	\begin{itemize}[noitemsep]
		\item[1) ]  \textbf{Measure} $(\bm{f}_p,\nabla_{\bm{u}}\bm{f}_p)|_{\bm{u}_{k}}$ on the plant.
		\item[2) ]
		\textbf{Update} the functions $(\bm{f},\phi,\bm{g})$ with \eqref{eq:3___20_Corrections_ISOy}.  
		\item[3) ]  \textbf{Make} the cost and constraint functions convex at $\bm{u}_k$: \eqref{eq:3___30_DefinitionMatrices_P}.
		\item[4) ]  \textbf{Compute}  $\bm{u}^\star_{k+1}$ with \eqref{eq:2___27a_ModelBasedPB}.
		\item[5) ] \textbf{Compute} the filter $	K_k = \text{\eqref{eq:2___45_LE_FILTRE}}$.
		\item[6) ] \textbf{Compute}  $\bm{u}_{k+1} = \bm{u}_{k} + K_k(\bm{u}^\star_{k+1}-\bm{u}_{k})$.
		\item[7) ] \textbf{Stop} if $\bm{u}_{k+1}\approx\bm{u}_{k}$ and return $\bm{u}_{\infty} := \bm{u}_{k+1}$.
	\end{itemize}
	\noindent {\bf end}
\end{BoxAlgo}

\noindent
Where one can summarize the filter update strategy as follows: 
\begin{itemize}  
	\item If ``$k=0$'', i.e. there is no previous iteration that can be used to estimate the gradient of the function $sol$. Then, the selected filter is the $K_o$.
	\item If ``$|\bm{u}_k-\bm{u}_{k-1}|\leq a$'', i.e. the size of the last step is less than or equal to the precision of the computer so one cannot estimate the gradient of the function $sol$. Then, the selected filter is the $K_{k-1}$.
	\item If ``$\utilde{\text{\hskip 0.1ex $\nabla$}}_k^S\geq 1$'', i.e. the estimated gradient of the function $sol$ is greater than $1$, so one is near a maximum or a saddle point of the plant, so applying a filter is useless. So, no filter needs to applied $K_k = 1$.
\end{itemize}  

\section{Conclusion}

This chapter starts from the formulation of the theoretical RTO problem. Then, the majority of the fundamental principles and methods of some of the most important contributions to theoretical RTO are rediscovered and improved in all their aspects, through a completely new progressive and didactic re-explanation. Indeed, thanks to a deep an intuitive explanation of how a filter works in the context of RTO, a new very efficient and meaningful filtering protocol has been introduced.  However, while this chapter is mainly focused on the ability to converge on plant's minimums and the speed of this convergence, the next chapter provides answers on how to make this convergence more secure.

	\chapter{A more secure RTO algorithm}
		\label{Chap:3_KMA}	

\section{An observation}

\begin{figure}[htp]
	\centering
	\includegraphics[width=13.35cm] {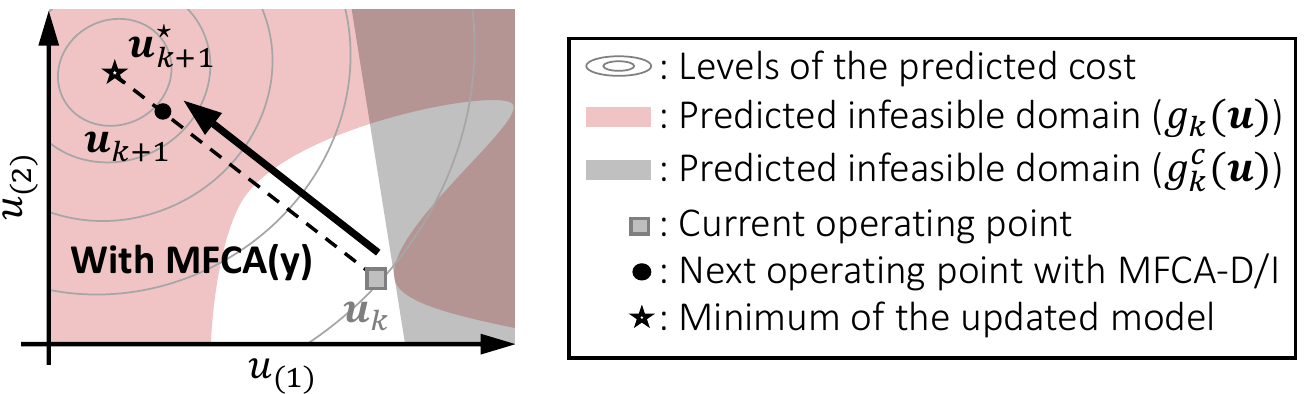}
	\caption{Illustration of the potential risk of convexifying locally concave constraints}
	\label{fig:3___1_Concept_idea2}
\end{figure}

The MFCA-D/I methods proposed at the end of the previous chapter are considered and the following thought experiments are performed:

Considering the case illustrated on Figure~\ref{fig:3___1_Concept_idea2}. This is a schematic representation of a 2-D optimization problem whose cost function  $\phi_{k}$ is convex and whose  constraint function $g_k$ is concave at $\bm{u}_k$. Therefore, if MFCA-D/I is applied, then the constraint $g_k$ would be linearized at $\bm{u}_k$ to obtain $g^c_k$. If $g^c_k$ is a relevant approximation of $g_k$ in the region of $\bm{u}_k$, it is clearly not true globally. In fact, a large area of the input space that was predicted to be infeasible becomes feasible. This presents the non-negligible risk that MFCA-D/I selects   $\bm{u}_{k+1}^{\star}$ and possibly $\bm{u}_{k+1}$  in this unfeasible area, as illustrated in Figure~\ref{fig:3___1_Concept_idea2}.  This makes it more likely that $\bm{u}_{k+1}$ violates the constraints of the plant. 

\section{An idea}

To counter this problem, an initial option would be to combine the optimization step \eqref{eq:2___27a_ModelBasedPB} with the filtering step  \eqref{eq:2___27b_Filter} in order to enforce $\bm{u}_{k+1}$ to be feasible according to the \textit{original} updated model  (i.e. the unconvexified one). This idea leads to two new RTO algorithms:
\vspace{-\topsep}
\begin{itemize}[noitemsep]
	\item  Alg~\ref{algo:K_MFCAv01}:  MFCA-D with filter-based constraints (KMFCA-D),
	\item  Alg~\ref{algo:K_MFCAyv01}: MFCA-I with filter-based constraints (KMFCA-I).
\end{itemize}
\begin{BoxAlgo}{\textbf{MFCA-D with filter-based constraints (KMFCA-D)}}{K_MFCAv01}
	\textbf{Initialization.} Provide $\bm{u}_0$, $\bm{K}_0(=0.1)$, $\bm{u}^*_0(=\bm{u}_0)$, $\bm{a}$, functions $(\bm{f},\phi,\bm{g})$, and the stopping criterion of step 5).
	\tcblower
	\textbf{for} $k=0 \rightarrow \infty$
	\begin{itemize}[noitemsep]
		\item[1) ]  \textbf{Measure} $(\nabla_{\bm{u}}\phi_p, \bm{g}_p,\nabla_{\bm{u}}\bm{g}_p)|_{\bm{u}_{k}}$ on the plant.
		\item[2) ] \textbf{Update} $\phi$ and $\bm{g}$ with \eqref{eq:3___19_Corrections_ISO}. 
		\item[3) ]  \textbf{Compute}  $\phi^c_k$ and $\bm{g}^c_k$ with \eqref{eq:3___30_DefinitionMatrices_P}.
		\item[4) ]  \textbf{Compute} $(\bm{u}^{\star}_{k+1},K_k,\bm{u}_{k+1})$ as follows:
	\end{itemize}
	\vspace{-\topsep}
			\begin{flalign} 
				\bm{u}^{\star}_{k+1} := \ &  \operatorname{arg}
				\underset{\bm{u}}{\operatorname{min}}  \quad    \phi^c_{k}(\bm{u}), \nonumber\\ 
				& \qquad   \text{s.t.} \quad  \bm{g}^c_{k}(\bm{u}) \leq \bm{0}, \nonumber \\
               	& \phantom{\qquad   \text{s.t.} \quad} \bm{g}_{k}\big(\bm{u}_k +  filt_k(\bm{u})(\bm{u} - \bm{u}_k ) \big) \leq \bm{0}, \label{eq:3___1_New_PB_OPT} \\
               	K_k := \ & filt_k(\bm{u}^\star_{k+1}), \label{eq:3___2_New_filter}\\
               	\bm{u}_{k+1} := \ & \bm{u}_k + K_k (\bm{u}^\star_{k+1} - \bm{u}_k), \label{eq:3___3_New_point} \\
               	\text{where}: \qquad\quad  
               	\nonumber \\
               \hspace{-3mm} filt_k(\bm{u}) := \ &  \left\{
               	\begin{array}{l@{}l}
               		\hspace{-2mm} K_o,     & \text{if } k =0, \\
               		\hspace{-2mm}	K_{k-1},  & \text{if } nab_k(\bm{u}) < 1, \   k > 0   \text{ and }   |\bm{u}_{k}-\bm{u}_{k-1} |\leq \bm{a}  \\
               		\hspace{-2mm} \max\left\{K_o,\frac{1}{1-nab_k(\bm{u})}\right\}\hspace{-1mm}, \
               		&  \text{if } nab_k(\bm{u}) < 1 , \    k > 0   \text{ and }  |\bm{u}_{k}-\bm{u}_{k-1}|>\bm{a},   \\
               		\hspace{-2mm} 1, & \text{otherwise.}
               	\end{array}
               	\right. \label{eq:3___4_fonction_filt} \\
               \hspace{-3mm}	nab_k(\bm{u})  := \ &
               	 \frac{\left(\bm{u}_k - \bm{u}_{k-1}\right)^{\rm T} 
               	 \left(\bm{u} - \bm{u}^{*}_{k} \right)
                }{\|\bm{u}_k - \bm{u}_{k-1}\|^2}.
			\end{flalign}
		\begin{itemize}[noitemsep]
			\item[5) ] \textbf{Stop} if $\bm{u}_{k+1}\approx\bm{u}_{k}$ and return $\bm{u}_{\infty} := \bm{u}_{k+1}$.
	\end{itemize}
	\noindent {\bf end}
\end{BoxAlgo}

\begin{BoxAlgo}{\textbf{MFCA-I with filter-based constraints  (KMFCA-I)}}{K_MFCAyv01}
	\textbf{Initialization.}  Provide $\bm{u}_0$, $\bm{K}_0(=0.1)$, $\bm{u}^*_0(=\bm{u}_0)$, $\bm{a}$, functions $(\bm{f},\phi,\bm{g})$, and the stoping criteron of step 5).
	\tcblower
	\textbf{for} $k=0 \rightarrow \infty$
	\begin{itemize}[noitemsep]
		\item[1) ]  \textbf{Measure} $(\bm{f}_p,\nabla_{\bm{u}}\bm{f}_p)|_{\bm{u}_{k}}$ on the plant.
		\item[2) ]
		\textbf{Update}  $\bm{f}$, $\phi$, $\bm{g}$ with \eqref{eq:3___20_Corrections_ISOy}.  
		\item[3) ]  \textbf{Compute}  $\phi^c_k$, $\bm{g}^c_k$ with \eqref{eq:3___30_DefinitionMatrices_P}.
		\item[4) ]  \textbf{Compute}  $(\bm{u}^{\star}_{k+1},K_k,\bm{u}_{k+1})$ with \eqref{eq:3___1_New_PB_OPT}, \eqref{eq:3___2_New_filter}, and  \eqref{eq:3___3_New_point}.
		\item[5) ] \textbf{Stop} if $\bm{u}_{k+1}\approx\bm{u}_{k}$ and return $\bm{u}_{\infty} := \bm{u}_{k+1}$.
	\end{itemize}
	\noindent {\bf end}
\end{BoxAlgo}
The functions $filt_k(\bm{u})$ and $nab_k(\bm{u})$ are the functional versions of  \eqref{eq:2___45_LE_FILTRE} and \eqref{eq:2___59_Nabla_S}. Their argument, $\bm{u}$, is the argument of optimization problem \eqref{eq:3___1_New_PB_OPT} and their outputs are what would be \eqref{eq:2___45_LE_FILTRE} and \eqref{eq:2___59_Nabla_S} if they were evaluated at   $\bm{u}^{\star}_{k+1} = \bm{u}$. 

By construction, the KMFCA-D/I algorithms avoid the problem discussed earlier and illustrated in~\ref{fig:3___1_Concept_idea2} as one can see in the following example:

\begin{exbox}\label{ex:3___1_Risque_de_lineariser_contraintes_concave}
	\textbf{(Illustration of the risks associated with the linearization of locally concave constraints)}
	One considers the following theoretical RTO problem: 
	\begin{align*}
		\phi(u,\bm{y}) := \ & y_{(1)},  &  y_{p(1)} := \ &  5 -14u +8u^2, \\
		& & y_{(1)} := \ & 0.2 - 2u +5u^2, \\
		g(u,\bm{y})    := \ & y_{(2)},  
		&  y_{p(2)} := \ & -2 + 7 u - 29 u^2 + 32 u^3,\\
		& & y_{(2)} := \ & -2 + 9 u - 35 u^2 + 40 u^3,
	\end{align*}
	where $u\in[0,1]$. Figure~\ref{fig:Exemple_3_1_Fonctions} illustrates what these functions are.\\
	\begin{minipage}[h]{\linewidth}
		\vspace*{0pt}
		\centering 
		\includegraphics[width=8cm]{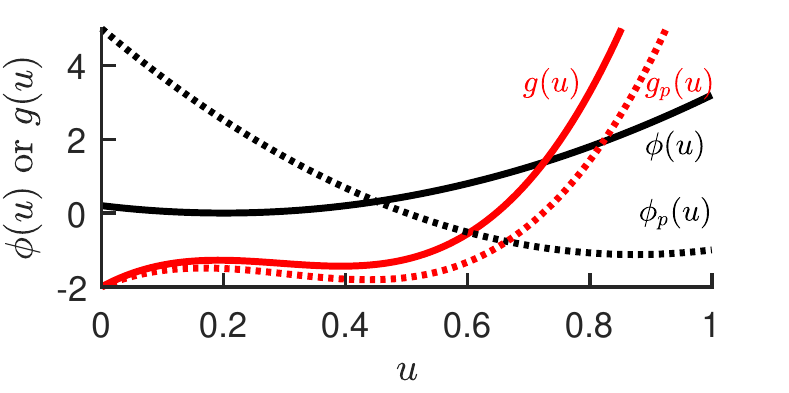}
		\vspace{-4mm}
		\captionof{figure}{ Costs and constraints functions  of the plant and the model.}
		\label{fig:Exemple_3_1_Fonctions}
	\end{minipage}\\
	
	The minimum of the nominal model is $u_0=0.2$. This is therefore the point at which one initializes the MFCA and KMFCA. The simulation results are shown on Figure~\ref{fig:Exemple_3_1_Results1}. \\
	\begin{minipage}[h]{\linewidth}
		\vspace*{0pt}
		\centering
		\includegraphics[width=4.45cm]{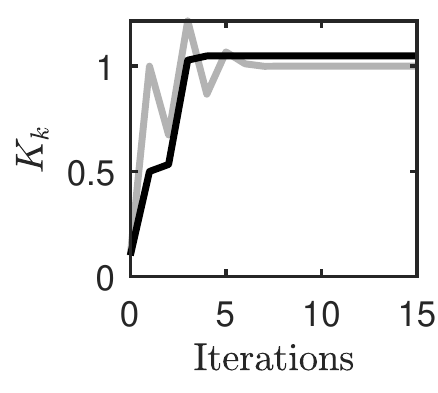}\hskip -0ex
		\includegraphics[width=4.45cm]{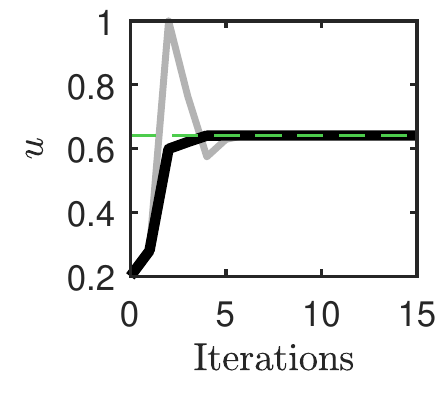}\hskip -0ex
		\includegraphics[width=4.45cm]{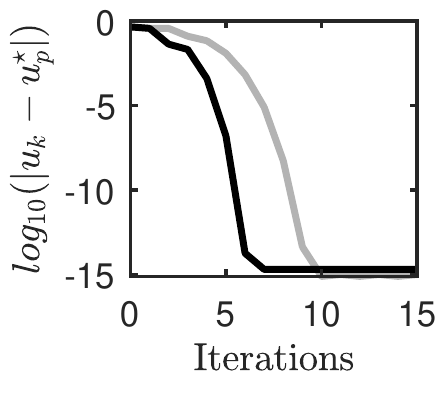} \\
		\begin{minipage}[b]{8.9cm}
			\includegraphics[width=4.45cm]{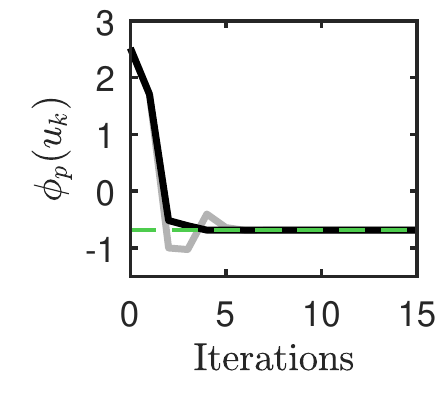}\hskip -0ex
			\includegraphics[width=4.45cm]{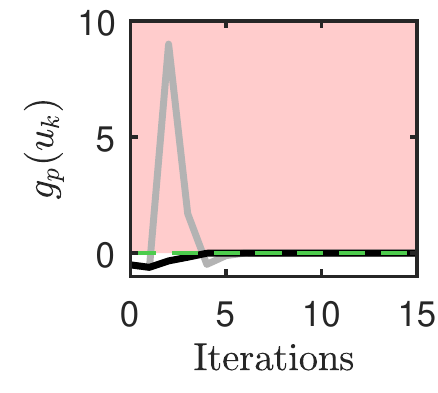}
		\end{minipage} 
		\begin{minipage}[b][3.4cm][t]{4.45cm}
			\noindent
			\textcolor{gray}{\raisebox{0.5mm}{\rule{0.5cm}{0.1cm}}}         : MFCA, \\
			\textcolor{black}{\raisebox{0.5mm}{\rule{0.5cm}{0.1cm}}}       : KMFCA, \\
			\textcolor{red!30!white}{\raisebox{-0.5mm}{\rule{0.5cm}{0.3cm}}} : Unfeasible area, \\
			\textcolor{green1}{\raisebox{0.5mm}{\rule{0.2cm}{0.05cm}\hspace{0.1cm}\rule{0.2cm}{0.05cm}}} : Optimum.
		\end{minipage} 
		\vspace{-4mm}
		\captionof{figure}{Simulation results}
		\label{fig:Exemple_3_1_Results1}
	\end{minipage} \\
	
	One can observe that KMFCA provides much better performances than MFCA. Indeed, MFCA significantly violates the constraints, and if one observes its curve $log_{10}(|u_k-u_p^{\star}|)$, one can see that it is roughly shifted by five iterations to the right. To explain these observations, a detailed analysis is presented in Figures~\ref{fig:Exemple_3_1_Results2} and \ref{fig:Exemple_3_1_Results3}.  \\
	\begin{minipage}[h]{\linewidth}
		\vspace*{0pt}
		\centering
		\includegraphics[trim={0 0 0.2cm 0},clip,height=3.9cm]{
			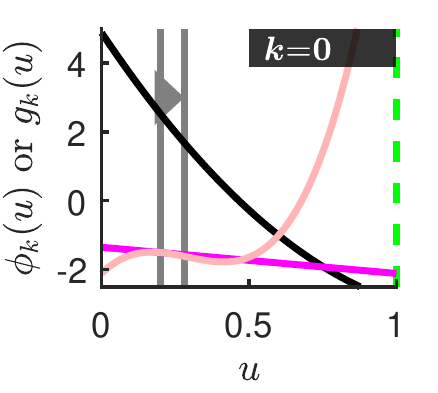}\hskip -0ex
		\includegraphics[trim={0.9cm 0 0.2cm 0},clip,height=3.9cm]{
			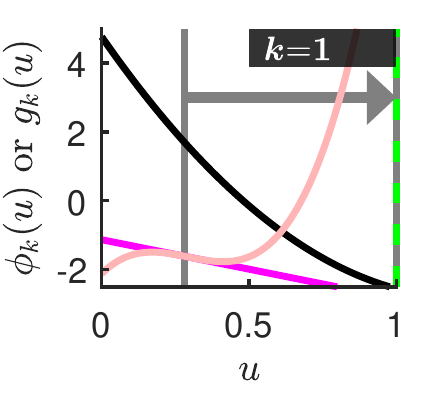}\hskip -0ex
		\includegraphics[trim={0.9cm 0 0.2cm 0},clip,height=3.9cm]{
			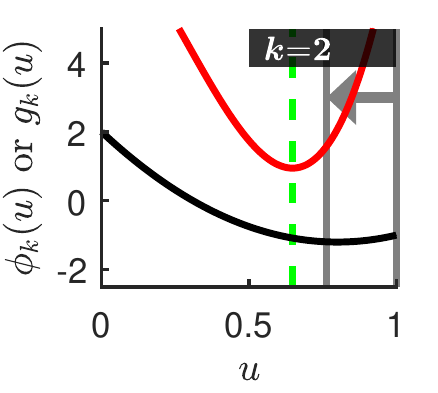}\hskip -0ex
		\includegraphics[trim={0.9cm 0 0.2cm 0},clip,height=3.9cm]{
			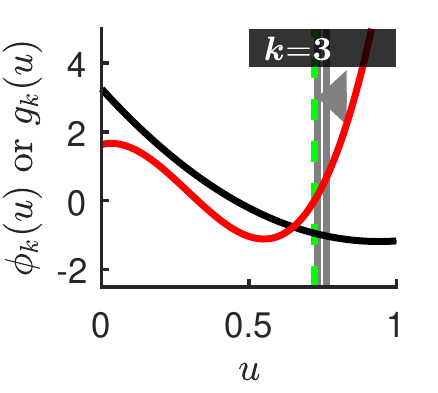} \\
		\begin{minipage}[b]{11cm}
			\includegraphics[trim={0 0 0.2cm 0},clip,height=3.9cm]{
				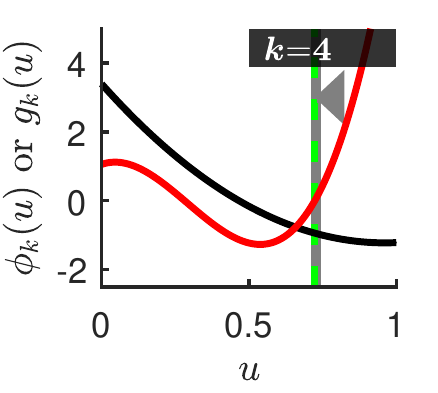}\hskip -0ex
			\includegraphics[trim={0.9cm 0 0.2cm 0},clip,height=3.9cm]{
				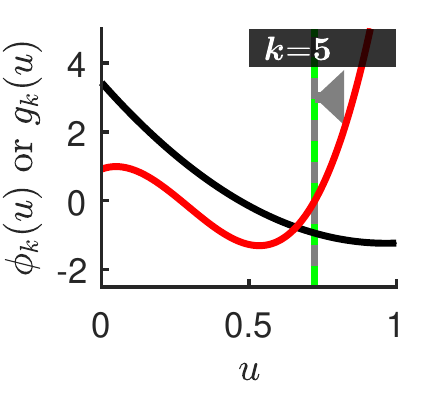}\hskip -0ex
			\includegraphics[trim={0.9cm 0 0.2cm 0},clip,height=3.9cm]{
				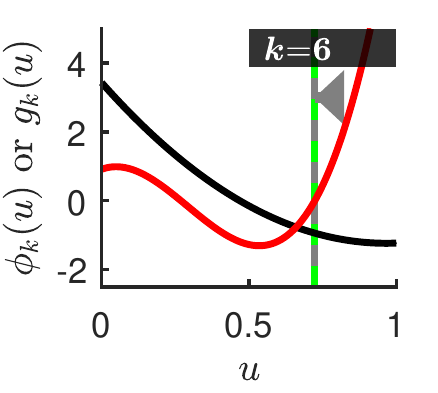}\hskip -0ex
		\end{minipage} 
		\hspace{-0.6cm}
		\begin{minipage}[b][3.7cm][t]{3.3cm}
			\noindent
			\textcolor{black}{\raisebox{0.5mm}{\rule{0.5cm}{0.1cm}}} : $\phi_k(u)$, \\
			\textcolor{magenta}{\raisebox{0.5mm}{\rule{0.5cm}{0.1cm}}} : $g^c_k(u)$, \\
			\textcolor{red}{\raisebox{0.5mm}{\rule{0.5cm}{0.1cm}}} : $g_k(u)$, \\
			\textcolor{red!50!white}{\raisebox{0.5mm}{\rule{0.5cm}{0.1cm}}} : $g_k(u)$ unused, \\
			\resizebox{0.5cm}{0.4cm}{\textcolor{gray}{\raisebox{-1mm}{\rule{0.1cm}{0.45cm}}\hspace{-1mm}\begin{tikzpicture}
						\draw[-{Triangle[width=0.2cm,length=0.25cm]}, line width=0.1cm](0,0) -- (0.7, 0);
					\end{tikzpicture}\hspace{-1mm}\raisebox{-1mm}{\rule{0.1cm}{0.45cm}}}}
			: $u_k \rightarrow u_{k+1}$, \\
			\phantom{}\hspace{0.2cm}\raisebox{-1mm}{\textcolor{green}{\rule{0.1cm}{0.45cm}}}\hspace{0.2cm}\hspace{-0.5cm}\textcolor{white}{\raisebox{0.75mm}{\rule{0.5cm}{0.1cm}}} 
			: $u^{\star}_{k+1}$, 
		\end{minipage} 
		\vspace{-4mm}
		\captionof{figure}{Details of MFCA iterations.}
		\label{fig:Exemple_3_1_Results2}
	\end{minipage}\\
	
	Figures~\ref{fig:Exemple_3_1_Results2} clearly illustrates the problem that MFCA faces. At iterations $k=0$ and $k=1$ the function $g_k$ is concave at $u_k$ and is therefore linearized, hence the magenta and light red curves for these two iterations. Since $g_k$ is ignored in favor of $g^c_k$, MFCA can aim at the points $u_1^{\star}$ and $u_2^{\star}$ that are not feasible according to $g_k$ and thus generate risky iterations. Hence the large constraint violation.  \\
	
	On the other hand, when KMFCA is used, although the constraint $g_k$ is linearized for the choice of $u_1^{\star}$ and $u_2^{\star}$, its satisfaction is imposed at the points  $u_{1}$ and $u_2$; hence, better management of constraints and the absence of violation \\

	\begin{minipage}[h]{\linewidth}
		\vspace*{0pt}
		\includegraphics[trim={0 0 0.2cm 0},clip,width=4.25cm]{
			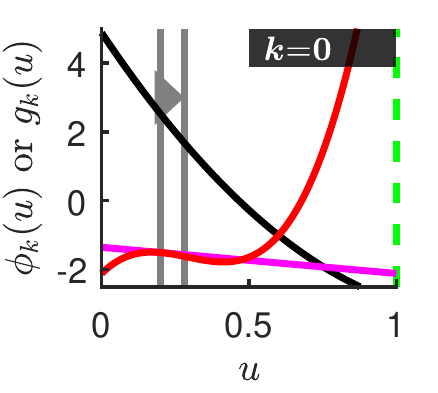}\hskip -0ex
		\includegraphics[trim={1.1cm 0 0.2cm 0},clip,width=3.15cm]{
			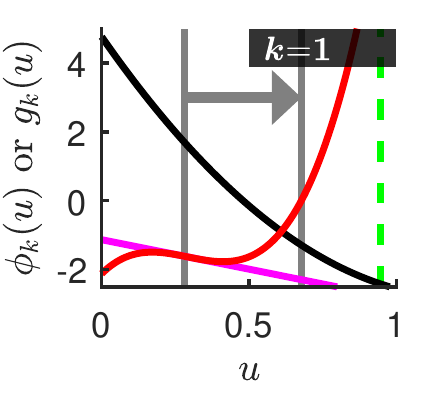}\hskip -0ex
		\includegraphics[trim={1.1cm 0 0.2cm 0},clip,width=3.15cm]{
			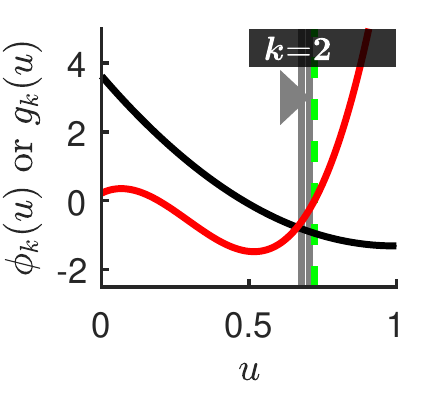}\hskip -0ex
		\includegraphics[trim={1.1cm 0 0.2cm 0},clip,width=3.15cm]{
			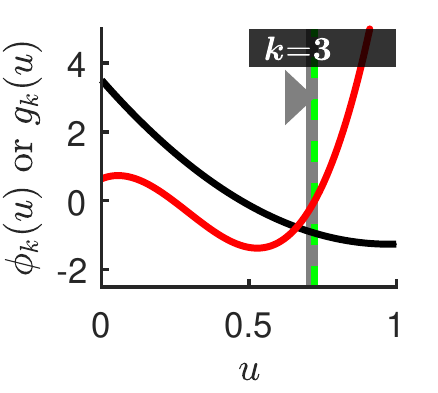} 
		\vspace{-4mm}
		\captionof{figure}{Details of KMFCA iterations (Legend: see Figure~\ref{fig:Exemple_3_1_Results2}).}
		\label{fig:Exemple_3_1_Results3}
	\end{minipage}
\end{exbox}

However, to validate these improvements, one needs evidence that the satisfaction of the equilibrium and stability conditions is not lost.   

\subsection{Effects on the equilibrium condition}
\label{sec:3_2_2_Effet_sur_la_convergence_de_Niveau_1}

To show that the satisfaction of the equilibrium condition is not lost, an analysis of the geometrical properties of Problems~\eqref{eq:2_2___1_Plant_PB} and \eqref{eq:3___1_New_PB_OPT} updated at $\bm{u}_k$ in the vicinity of  $\bm{u}_k$ is required. Indeed, it is clear that the fact that \eqref{eq:3___1_New_PB_OPT} has twice as many constraints as the plant \eqref{eq:2_2___1_Plant_PB} implies that, globally, the feasible domains of these two optimization problems are different, i.e.:
\begin{align*}
	\bm{g}_p(\bm{u}) \leq \ & \bm{0}, & 
	& \not\Leftrightarrow & 
	\left( \begin{array}{c}
		\bm{g}^c_k(\bm{u}) \\
		\bm{g}_k(\bm{u}_k + filt_k(\bm{u})(\bm{u}-\bm{u}_k))
	\end{array} \right) \leq \ & \bm{0}, \quad \forall \bm{u}\in\amsmathbb{R}^{n_u}.
\end{align*}
However, this does not mean that around $\bm{u}_k$ they are not similar. In the previous chapter is was shown that, by forcing these two problems to have similarities at $\bm{u}_k$ (see \eqref{eq:2___11_SimpleApproach_SameCones} and \eqref{eq:3___21}), the validity of the conditions of equilibrium and stability of MFCA-D/I can be enforced. So, one chooses to analyze these two optimization problems at $\bm{u}_k$ through their cones of feasible directions (CDF) at that point. These two cones, called $\text{CFD}_{\text{p,k}}$ for the plant and $\text{CFD}_{\text{k}}$ for the updated model, can be defined mathematically as follows:
\begingroup
\allowdisplaybreaks
\begin{align}
	\text{CFD}_{\text{p,k}} := \ &  \left\{ 
	\bm{d}\in\amsmathbb{R}^{n_u} \ | \ \bm{g}_p(\bm{u}_k)\leq \bm{0}, \ \nabla_{\bm{u}}g_{p(i)}|_{\bm{u}_k}^{\rm T}\bm{d}<0 \text{ if }  g_{p(i)}|_{\bm{u}_k} = 0, \forall i
	\right\}, \label{eq:3___2_CFDk}\\
	\text{CFD}_{\text{k}} := \ &  \big\{ 
	\bm{d}\in\amsmathbb{R}^{n_u} \ | \ \bm{g}_k(\bm{u}_k)\leq \bm{0}, \ \nabla_{\bm{u}}g^c_{k(i)}|_{\bm{u}_k}^{\rm T}\bm{d}<0  \text{ if }  g^c_{k(i)}|_{\bm{u}_k} = 0, \text{ and } ...   \nonumber \\ 
	& 
	\ \  \nabla_{\bm{u}}g_{k(i)}(\bm{u}_k + filt_k(\bm{u})  (\bm{u}-\bm{u}_k))|_{\bm{u}_k}^{\rm T}\bm{d}<0 \text{ if }  g_{k(i)}(\bm{u}_k + filt_k(\bm{u}) ... \nonumber \\ 
	& \ \ (\bm{u}-\bm{u}_k)|_{\bm{u}_k} = 0, \forall i
	\big\},\label{eq:3___2_CFDk_prime}
\end{align}\endgroup
and they can be geometrically interpreted as the set of directions along which one can move starting from a feasible point $\bm{u}_k$ while remaining in the feasible domain. If one can prove that $\text{CFD}_{\text{p,k}} = \text{CFD}_{\text{k}}$, then it would mean that locally the feasible and non-feasible areas of  \eqref{eq:2_2___1_Plant_PB} and  \eqref{eq:3___1_New_PB_OPT} are the same, i.e.:
\begin{align} \label{eq:3___8_MemesCones_equivalence}
	\bm{g}_p(\bm{u}) \leq \ & \bm{0}, & 
	& \Leftrightarrow & 
	\left( \begin{array}{c}
		\bm{g}^c_k(\bm{u}) \\
		\bm{g}_k(\bm{u}_k + filt_k(\bm{u})(\bm{u}-\bm{u}_k))
	\end{array} \right) \leq \ & \bm{0}, \quad \forall \bm{u}\in\mathcal{B}(\bm{u}_k,r\rightarrow 0).
\end{align}

With the following Lemma, it is shown that this is always the case.
\begin{lembox} \label{lem:3___1_CFD_MFCQ}
	Consider the optimization problem of the plant given by
	\eqref{eq:2_2___1_Plant_PB},
	and the one of the model updated at $\bm{u}_k$ given by \eqref{eq:3___1_New_PB_OPT}, such that equalities \eqref{eq:2___11_SimpleApproach_SameCones} and \eqref{eq:3___21} are true. Then:   
	\begin{itemize}[noitemsep]	 
		\item The cone of feasible directions of \eqref{eq:2_2___1_Plant_PB} at $\bm{u}_k$: $\text{CFD}_\text{p,k} := \text{\eqref{eq:3___2_CFDk}}$, and the cone of feasible directions of \eqref{eq:3___1_New_PB_OPT} at $\bm{u}_k$:  $\text{CFD}_\text{k} := \text{\eqref{eq:3___2_CFDk_prime}}$, are equal: 
		\begin{equation}
			\text{CFD}_\text{p,k} = \text{CFD}_\text{k}.
		\end{equation}
		\item If Problem~\eqref{eq:2_2___1_Plant_PB} is LICQ at $\bm{u}_k$, then Problem~\eqref{eq:3___1_New_PB_OPT} is MFCQ (Mangasarian-Fromovitz Constaints Qualification) at $\bm{u}_k$. .
	\end{itemize}
\end{lembox}
\begin{proofbox}
	\begin{subequations}
	Starting by proving the first statement. To do so, the following three observations are made:
	\begin{align} \label{eq:3___5_a}
		g_{k(i)}(\bm{u}_k + \underbrace{filt_k(\bm{u})(\bm{u}-\bm{u}_k)|_{\bm{u}_k}}_{=\bm{0}} = g_{k(i)}|_{\bm{u}_k}. 
	\end{align}
	\begin{align}
		&  \nabla_{\bm{u}}g_{k(i)}(\bm{u}_k + filt_k(\bm{u})(\bm{u}-\bm{u}_k))|_{\bm{u}_k} \nonumber \\ 
		& \qquad \qquad \qquad \qquad = 
		\left(
		\nabla_{\bm{u}}filt_k|_{\bm{u}_k}(\bm{u}_k-\bm{u}_k)	
		\hspace{-1mm} + \hspace{-1mm}
		filt_k|_{\bm{u}_k}	\right) \nabla_{\bm{u}} g_{k(i)}|_{\bm{u}_k}, \nonumber \\
		 & \qquad\qquad \qquad  \qquad = 
		filt_k|_{\bm{u}_k}	 \nabla_{\bm{u}} g_{k(i)}|_{\bm{u}_k}. \label{eq:3___5_b}
	\end{align}
	\begin{align}
		&\text{\eqref{eq:3___21}} \Rightarrow &
		\bm{g}_k|_{\bm{u}_{k}} = \ & \bm{g}^c_k|_{\bm{u}_{k}}, & 
		\nabla_{\bm{u}} \bm{g}_k|_{\bm{u}_{k}} = \ & \nabla_{\bm{u}}\bm{g}^c_k|_{\bm{u}_{k}}, \label{eq:3___5_c}
	\end{align}
	By combining \eqref{eq:3___5_a}, \eqref{eq:3___5_b}, \eqref{eq:3___5_c}, and by construction $filt_k|_{\bm{u}_k}>0$, the following is stated:
	\begin{align}
		& 
		\begin{array}{r@{}lcr@{}l}
			g^c_{k(i)}|_{\bm{u}_k} = \ & 0, & \Rightarrow &
			g^c_{k(i)}(\bm{u}_k + filt_k(\bm{u})(\bm{u}-\bm{u}_k)|_{\bm{u}_k}= \ & 0 , \\
			\nabla_{\bm{u}} g_{k(i)}|_{\bm{u}_k}^{\rm T}\bm{d} >\ & 0, &  \Rightarrow & 
			\nabla_{\bm{u}}g_{k(i)}(\bm{u}_k + filt_k(\bm{u})(\bm{u}-\bm{u}_k))|_{\bm{u}_k}^{\rm T}\bm{d} >\ & 0.
		\end{array} \nonumber
	\end{align}
	So the definition \eqref{eq:3___2_CFDk_prime} of 	$\text{CFD}_{\text{k}}$ can be reduced to:
	\begin{align}
	  \text{CFD}_{\text{k}} = 
	  \big\{ 
	  \bm{d}\in\amsmathbb{R}^{n_u} \ | \ \nabla_{\bm{u}}g^c_{k(i)}|_{\bm{u}_k}^{\rm T}\bm{d}<0  \text{ if }  g^c_{k(i)}|_{\bm{u}_k} = 0, \forall i \big\}.
	  \label{eq:3___6_c}
	\end{align}
	As equalities \eqref{eq:2___11_SimpleApproach_SameCones} are assumed to be true, i.e. $g^c_{k(i)}|_{\bm{u}_k} = \  g_{p(i)}|_{\bm{u}_k}$ and  $\nabla_{\bm{u}}g^c_{k(i)}|_{\bm{u}_k} =  \nabla_{\bm{u}}g_{p(i)}|_{\bm{u}_k}$,
	the definition of $\text{CFD}_{\text{k}}$, given by \eqref{eq:3___6_c}, can be reformulated as: 
	\begin{align}
		\text{CFD}_{\text{k}} = 
		\big\{ 
		\bm{d}\in\amsmathbb{R}^{n_u} \ | \ \nabla_{\bm{u}}g_{p(i)}|_{\bm{u}_k}^{\rm T}\bm{d}<0  \text{ if }  g_{p(i)}|_{\bm{u}_k} = 0, \forall i \big\} = \text{CFD}_{\text{p,k}}.
		\label{eq:3___6_d}
	\end{align}
	Which proves the first statement of the lemma.
	 
	Now one can easily prove the second statement. To prove that Problem~\eqref{eq:3___1_New_PB_OPT} is MFCQ at $\bm{u}_k$ for all $\bm{u}_k\in\amsmathbb{R}^{n_u}$, one must show that $\forall\bm{u}_k\in\amsmathbb{R}^{n_u}:$  $\text{CFD}_{\text{k}} \neq \emptyset$. By definition, the CFD of a LICQ problem  is never empty \cite{Chong:2001}. So, if \eqref{eq:2_2___1_Plant_PB} is LICQ at $\bm{u}_k$, then $\text{CFD}_{\text{p,k}}\neq \emptyset$. As according to \eqref{eq:3___6_d},  $\text{CFD}_{\text{p,k}}= \text{CFD}_{\text{k}}$, one can conclude that: $\text{CFD}_{\text{k}}\neq \emptyset$.
	So the Problem~\eqref{eq:3___1_New_PB_OPT} is always MFCQ at $\bm{u}_k$.
\end{subequations}  
\end{proofbox}

On the basis of this local geometric similarity between Problems~\eqref{eq:3___1_New_PB_OPT} updated at $\bm{u}_k$ and \eqref{eq:2_2___1_Plant_PB}, it can be shown that KMFCA-D/I converges only on a KKT point of the plant, see Theorem~\ref{thm:3___1_KKT_matching_at_conv}. Also, the equilibrium condition of KMFCA-D/I can be identified, see Theorem~\ref{thm:3___2_ConditionNiveau1}.
\begin{thmbox} \label{thm:3___1_KKT_matching_at_conv}
	If KMFCA-D/I converges on the limit value $\bm{u}_\infty := 	\underset{k\rightarrow \infty}{\lim} \bm{u}_k$, then 
	$\exists \bm{\lambda} \in\amsmathbb{R}^{2n_g}$ such that $(\bm{u}_\infty, \bm{\lambda})$ is a KKT point of the model updated at $\bm{u}_\infty$, and such that $(\bm{u}_\infty, \bm{\lambda}_p)$ is a KKT point of the plant with $\bm{\lambda}_p = \bm{\lambda}^{\prime} + filt_{\infty}|_{\bm{u}_{\infty}}\bm{\lambda}^{\prime\prime}$, where $(\bm{\lambda}^{\prime},\bm{\lambda}^{\prime\prime})\in\amsmathbb{R}^{n_g}$ and $	\bm{\lambda} := [
	\bm{\lambda}^{\prime \rm T}, \
	\bm{\lambda}^{\prime\prime \rm T}]^{\rm T}$.
\end{thmbox}
\begin{proofbox}
	Thanks to Lemma~\ref{lem:3___1_CFD_MFCQ}, one knows that Problem~\eqref{eq:3___1_New_PB_OPT} is MFCQ at $\bm{u}_k$ $\forall k$, and therefore at $\bm{u}_\infty$. So, according to \cite{Chong:2001}: 
	\begingroup
	 \allowdisplaybreaks
	\begin{align} 
		& \bm{u}_{\infty} := \operatorname{arg}
		\underset{\bm{u}}{\operatorname{min}} \quad   \phi_{\infty}(\bm{u}) %
		\quad 
		\text{s.t.} \quad \bm{g}^c_\infty(\bm{u}) \leq \bm{0}, \nonumber \\
		& \phantom{ \bm{u}_{\infty} := \operatorname{arg}
			\underset{\bm{u}}{\operatorname{min}} \quad   \phi_{\infty}(\bm{u}) %
			\quad 
			\text{s.t.} \quad}  \bm{g}_\infty\big(\bm{u}_\infty + filt_\infty(\bm{u})(\bm{u}-\bm{u}_\infty)\big) \leq \bm{0},  \label{eq:3___7}\\
		\nonumber \\
		\hspace{-0.1mm} \Rightarrow  &  \exists \bm{\lambda} \in \amsmathbb{R}^{2n_g}:
		\left\{
		\begin{array}{l}
			\left( \begin{array}{c}
				\bm{g}^c_\infty|_{\bm{u}_{\infty}} \\
				\bm{g}_\infty(\bm{u}_\infty + filt_\infty(\bm{u})(\bm{u}-\bm{u}_\infty))|_{\bm{u}_{\infty}} 
			\end{array} \right) \leq \bm{0},  \\
			\bm{\lambda}^{\rm T} 
			\left( \begin{array}{c}
				\bm{g}^c_\infty|_{\bm{u}_{\infty}} \\
				\bm{g}_\infty(\bm{u}_\infty + filt_\infty(\bm{u})(\bm{u}-\bm{u}_\infty))|_{\bm{u}_{\infty}}
			\end{array} \right) = 0, \\
			\bm{\lambda} \geq\bm{0}, \\
			\nabla_{\bm{u}} \phi_{\infty}|_{\bm{u}_{\infty}} + \bm{\lambda}^{\rm T} 
			\left(
				\hspace{-1mm}
				 \begin{array}{c}
				\nabla_{\bm{u}} \bm{g}^c_{\infty}|_{\bm{u}_{\infty}} \\
				\nabla_{\bm{u}} \bm{g}_\infty(\bm{u}_\infty + filt_\infty(\bm{u})(\bm{u}-\bm{u}_\infty))|_{\bm{u}_{\infty}}
			\end{array} 
			\hspace{-1mm} 
			\right)
			= \bm{0}, 
		\end{array}
		\right. \label{eq:3___8}
	\end{align}
 	\endgroup
	where  $\bm{\lambda} := (\bm{\lambda}^{\prime \rm T},  \bm{\lambda}^{\prime\prime\rm T})^{\rm T} \in \amsmathbb{R}^{2n_g}$,   $\bm{\lambda}^{\prime} \in \amsmathbb{R}^{n_g}$ are the Lagrange multipliers associated to the constraints  $\bm{g}^c_{\infty}|_{\bm{u}_{\infty}}$, and $\bm{\lambda}^{\prime\prime} \in \amsmathbb{R}^{n_g}$ are the Lagrange multipliers associated with the constraints  $\bm{g}_{\infty}(\bm{u}_\infty + filt_\infty(\bm{u})(\bm{u}-\bm{u}_\infty))|_{\bm{u}_{\infty}}$. Noting that:
	\begingroup
	\allowdisplaybreaks
	\begin{align}
		\bm{g}_\infty\big(\bm{u}_\infty + filt_\infty(\bm{u})(\bm{u}-\bm{u}_\infty)\big)\big|_{\bm{u}_{\infty}} = \ & \bm{g}_\infty|_{\bm{u}_{\infty}}, \label{eq:3___9_kbrgcs}\\
		\nabla_{\bm{u}} \bm{g}_\infty\big(\bm{u}_\infty + filt_\infty(\bm{u})(\bm{u}-\bm{u}_\infty)\big)\big|_{\bm{u}_{\infty}} = \ & \big(\nabla_{\bm{u}} filt_\infty|_{\bm{u}_{\infty}} (\bm{u}_{\infty} - \bm{u}_{\infty}) +  ... \nonumber \\ 
		 & filt_\infty|_{\bm{u}_{\infty}}\big) \nabla_{\bm{u}} \bm{g}_\infty|_{\bm{u}_{\infty}},  \nonumber\\
		  = \ &filt_\infty|_{\bm{u}_{\infty}}\nabla_{\bm{u}} \bm{g}_\infty|_{\bm{u}_{\infty}}. \label{eq:3___10_iuachngb}
	\end{align}
	\endgroup
	and that 
	\begin{align} \label{eq:3___15_sgljxbhv}
		 &\text{\eqref{eq:3___21}} \Rightarrow &
		\bm{g}_\infty|_{\bm{u}_{\infty}} = \ & \bm{g}^c_\infty|_{\bm{u}_{\infty}}, & 
		\nabla_{\bm{u}} \bm{g}_\infty|_{\bm{u}_{\infty}} = \ & \nabla_{\bm{u}}\bm{g}^c_\infty|_{\bm{u}_{\infty}},
	\end{align}
	system \eqref{eq:3___8} can be re-written as follows:
	\begin{equation} \label{eq:3___9}
		 \left. \begin{array}{c}
		 	\text{\eqref{eq:3___8}} \\
		 	\text{\eqref{eq:3___9_kbrgcs}} \\
		 	\text{\eqref{eq:3___10_iuachngb}}\\
		 	\text{\eqref{eq:3___15_sgljxbhv}}
		 \end{array} \right\}
		 \Rightarrow  \exists \bm{\lambda} \in \amsmathbb{R}^{2n_g}:
		\left\{
		\begin{array}{l}
			\bm{g}_\infty|_{\bm{u}_{\infty}}  \leq \bm{0},  \\
			\big(\bm{\lambda}^{\prime} +
			\bm{\lambda}^{\prime\prime}\big)^{\rm T} 	\bm{g}_\infty|_{\bm{u}_{\infty}} = 0, \\
			\bm{\lambda} \geq\bm{0}, \\
			\nabla_{\bm{u}} \phi_{\infty}|_{\bm{u}_{\infty}} + \big(\bm{\lambda}^{\prime} +
			filt_\infty|_{\bm{u}_{\infty}}
			\bm{\lambda}^{\prime\prime}\big)^{\rm T} 
			\nabla_{\bm{u}} \bm{g}_{\infty}|_{\bm{u}_{\infty}}
			= \bm{0}. 
		\end{array}
		\right.
	\end{equation}
	In addition, as $\bm{g}_\infty|_{\bm{u}_{\infty}}  \leq \bm{0}$,  $\bm{\lambda} \geq\bm{0}$ and $filt_\infty|_{\bm{u}_{\infty}}>0$:
	\begin{align}
		\big(\bm{\lambda}^{\prime} +
		\bm{\lambda}^{\prime\prime}\big)^{\rm T} 	\bm{g}_\infty|_{\bm{u}_{\infty}} = \ &  0, & 
		\Leftrightarrow & &
			\sum_{i=1}^{n_g} \lambda^{\prime}_{(i)} g_{\infty(i)}|_{\bm{u}_{\infty}} +
		\lambda^{\prime\prime}_{(i)} g_{\infty(i)}|_{\bm{u}_{\infty}} = \ & 0, \nonumber \\ 
		 & & 
		\Leftrightarrow & &
		\lambda^{\prime}_{(i)} g_{\infty(i)}|_{\bm{u}_{\infty}} = 
		\lambda^{\prime\prime}_{(i)} g_{\infty(i)}|_{\bm{u}_{\infty}}
		= \ & 0, \ \forall i, \nonumber  \\
		& & 
		\Leftrightarrow & &
		\lambda^{\prime}_{(i)} g_{\infty(i)}|_{\bm{u}_{\infty}} = 
		filt_\infty|_{\bm{u}_{\infty}} \lambda^{\prime\prime}_{(i)} g_{\infty(i)}|_{\bm{u}_{\infty}}
		= \ & 0, \ \forall i, \nonumber \\
		& & 
		\Leftrightarrow & &
		\big(\bm{\lambda}^{\prime} +
		filt_\infty|_{\bm{u}_{\infty}} \bm{\lambda}^{\prime\prime}\big)^{\rm T} 	\bm{g}_\infty|_{\bm{u}_{\infty}} = \ &  0. \label{eq:3___10_b}
	\end{align}
	So, \eqref{eq:3___9} can be rewritten: 
	\begin{equation} \label{eq:3___10}
		\left. \begin{array}{c}
			\text{\eqref{eq:3___9}} \\
			\text{\eqref{eq:3___10_b}}
		\end{array} \right\}
		\Rightarrow   \exists \bm{\lambda} \in \amsmathbb{R}^{2n_g}:
		\left\{
		\begin{array}{l}
			\bm{g}_\infty|_{\bm{u}_{\infty}}  \leq \bm{0},  \\
			\big(\bm{\lambda}^{\prime} +
			filt_\infty|_{\bm{u}_{\infty}} \bm{\lambda}^{\prime\prime}\big)^{\rm T} 	\bm{g}_\infty|_{\bm{u}_{\infty}} = 0, \\
			\bm{\lambda} \geq\bm{0}, \\
			\nabla_{\bm{u}} \phi_{\infty}|_{\bm{u}_{\infty}} + \big(\bm{\lambda}^{\prime} +
			filt_\infty|_{\bm{u}_{\infty}}
			\bm{\lambda}^{\prime\prime}\big)^{\rm T} 
			\nabla_{\bm{u}} \bm{g}_{\infty}|_{\bm{u}_{\infty}}
			= \bm{0}. 
		\end{array}
		\right.
	\end{equation}
	As KMFCA-D/I enforces the equalities \eqref{eq:2___11_SimpleApproach_SameCones}, one can replace $(\bm{g}_{\infty}, \nabla_{\bm{u}}\bm{g}_{\infty},$\\ $ \nabla_{\bm{u}}\phi_{\infty})|_{\bm{u}_{\infty}}$ with $(\bm{g}_{p}, \nabla_{\bm{u}}\bm{g}_{p}, \nabla_{\bm{u}}\phi_{p})|_{\bm{u}_{\infty}}$,  and as
	\begin{align} \label{eq:3___11}
		\bm{\lambda} \geq \ & \bm{0}, & 
		\Rightarrow & & 
		\big(\bm{\lambda}^{\prime} +
		filt_\infty|_{\bm{u}_{\infty}} \bm{\lambda}^{\prime\prime}\big) \geq \ & \bm{0},
	\end{align} 
	it is observed that $\bm{u}_{\infty}$ is a KKT point of the plant; indeed:
	\begin{align} 	\label{eq:3___19_sdkcuggs}
		\left. \begin{array}{c}
			\text{\eqref{eq:3___10}} \\
			\text{\eqref{eq:2___11_SimpleApproach_SameCones}} \\
			\text{\eqref{eq:3___11}}
		\end{array} \right\}
		\Rightarrow  \exists \bm{\lambda} \in \amsmathbb{R}^{2n_g}: &   
		\left\{
		\begin{array}{l}
			\bm{g}_{p}|_{\bm{u}_{\infty}} \leq \bm{0},\\
			\big(\bm{\lambda}^{\prime} +
			filt_\infty|_{\bm{u}_{\infty}} \bm{\lambda}^{\prime\prime}\big)^{\rm T} \bm{g}_{p}|_{\bm{u}_{\infty}} = \bm{0}, \\
			\big(\bm{\lambda}^{\prime} +
			filt_\infty|_{\bm{u}_{\infty}} \bm{\lambda}^{\prime\prime}\big) \geq \bm{0}, \\
			\nabla_{\bm{u}}\phi_{p}|_{\bm{u}_{\infty}} + \big(\bm{\lambda}^{\prime} +
			filt_\infty|_{\bm{u}_{\infty}} \bm{\lambda}^{\prime\prime}\big)^{\rm T} \nabla_{\bm{u}}\bm{g}_{p}|_{\bm{u}_{\infty}} = \bm{0}.
		\end{array}
		\right. 
	\end{align}
	Therefore, if KMFCA-D/I converges on $\bm{u}_{\infty}$, then  $(\bm{u}_{\infty},\bm{\lambda})$ is a KKT point of the updated model, and  $\big(\bm{u}_{\infty},\bm{\lambda}_p=(\bm{\lambda}^{\prime} +
	filt_\infty|_{\bm{u}_{\infty}} \bm{\lambda}^{\prime\prime})\big)$ is a KKT point of the plant. 
\end{proofbox}

\begin{thmbox}
	\textbf{(Equilibrium condition of KMFCA-D/I)} \label{thm:3___2_ConditionNiveau1}
	Consider the plant optimization problem given by \eqref{eq:2_2___1_Plant_PB}, its solution ($\bm{u}_p^{\star}$), the RTO algorithms: KMFCA-D/I, and the following notations:
	\begin{itemize}[leftmargin=*]
		\item The functions $\phi_{\star}$ and $\bm{g}_{\star}$ are the functions of the model updated at $\bm{u}_p^{\star}$, and Problem~\eqref{eq:3___1_New_PB_OPT}${}|_{\star}$ refers to Problem~\eqref{eq:3___1_New_PB_OPT} updated at $\bm{u}_p^{\star}$,
		\item $\bm{N}_{\star} \in \amsmathbb{R}^{n_u \times n_b}$  is a matrix whose $n_b$ columns form an orthonormal basis of the null space of the Jacobian of active constraints of Problem~\eqref{eq:3___1_New_PB_OPT}${}|_{\star}$ at $\bm{u}_p^{\star}$, i.e. $\forall i$  such that  $g_{\star(i)}$  is an active constraint of Problem~\eqref{eq:3___1_New_PB_OPT}${}|_{\star}$ at $\bm{u}_p^{\star}$:
		\begin{align}
			\nabla_{\bm{u}} g_{\star(i)}|_{\bm{u}_p^{\star}} \bm{N}_{\star} = \ & \bm{0}.
		\end{align}
		\item $\bm{N}_{p}\in \amsmathbb{R}^{n_u \times n_c}$  is a matrix whose $n_c$ columns form an orthonormal basis of the null space of the Jacobian of active constraints of Problem~\eqref{eq:2_2___1_Plant_PB} at $\bm{u}_p^{\star}$, i.e. $\forall i$  such that  $g_{p(i)}$  is an active constraint of Problem~\eqref{eq:2_2___1_Plant_PB}${}|_{\star}$ at $\bm{u}_p^{\star}$:
		\begin{align}
			\nabla_{\bm{u}} g_{p(i)}|_{\bm{u}_p^{\star}} \bm{N}_{p} = \ & \bm{0}.
		\end{align}
	\end{itemize}
	Then, the equilibrium condition of KMFCA-D/I is that 
	\begin{equation} \label{eq:3___16_iunefxux}
		\bm{N}_p^{\rm T}
		\left( 
		\nabla^2_{\bm{uu}}\phi^c_{\star}|_{\bm{u}_p^{\star}}
		+
		\sum_{i=1}^{n_g}  
		\lambda_{p(i)}
		\nabla^2_{\bm{uu}} g^{c}_{\star(i)}|_{\bm{u}_p^{\star}}
		\right)
		\bm{N}_p > 0.
	\end{equation}
	By default this condition is satisfied because  $\nabla^2_{\bm{uu}}\phi_{k}|_{\bm{u}_k}>0$ and $\nabla^2_{\bm{uu}} g^{c}_{\star(i)}|_{\bm{u}_p^{\star}}\geq 0$ $\forall k$. 
\end{thmbox}
\nointerlineskip
\begin{proofbox}
		For the equilibrium to be possible, the 1st and 2nd order optimality conditions of  Problem~\eqref{eq:3___1_New_PB_OPT}${}|_{\star}$ must be satisfied at $\bm{u}_p^{\star}$. 
		
		\begin{center}
			\textit{Part 1  --  1st order optimality conditions:} 
		\end{center}
		\noindent
		By definition:
		\begin{align} 
			& \bm{u}_p^{\star} := \operatorname{arg}
			\underset{\bm{u}}{\operatorname{min}} \quad   \phi_{p}(\bm{u}) 
			\quad 
			\text{s.t.} \quad \bm{g}_p(\bm{u}) \leq \bm{0},  \label{eq:3___wegrw} 
		\end{align}
		\begin{align} 
			\text{\eqref{eq:3___wegrw} }  \Leftrightarrow \ &  \exists \bm{\lambda}_p \in \amsmathbb{R}^{n_g} \text{ tel que } 
			\left\{
			\begin{array}{l}
				\bm{g}_{p}|_{\bm{u}_p^{\star}} \leq \bm{0},\\
				\bm{\lambda}_p^{\rm T} \bm{g}_{p}|_{\bm{u}_p^{\star}} = \bm{0}, \\
				\bm{\lambda}_p \geq \bm{0}, \\
				\nabla_{\bm{u}}\phi_{p}|_{\bm{u}_p^{\star}} + \bm{\lambda}_p^{\rm T} \nabla_{\bm{u}}\bm{g}_{p}|_{\bm{u}_p^{\star}} = \bm{0}.
			\end{array}
			\right. \label{eq:3___14_Proof_KKT1_Plant} 
		\end{align}
		As KMFCA-D/I implies that the functions of  Problem~\eqref{eq:3___1_New_PB_OPT}${}|_{\star}$ satisfy both  \eqref{eq:2___11_SimpleApproach_SameCones} and 
		\eqref{eq:3___21}: 
		\begin{align} \label{eq:3___15}
			\left.\begin{array}{c}
				\bm{g}_{p}|_{\bm{u}_p^{\star}} \leq \bm{0} \\
				\text{\eqref{eq:2___11_SimpleApproach_SameCones}} \\
				\text{\eqref{eq:3___21}}
			\end{array}\right\}
			& &
			\Rightarrow & &
			\left(\begin{array}{c}
				\bm{g}^c_{\star}|_{\bm{u}_p^{\star}} \\
				\bm{g}_{\star}\big(\bm{u}_p^{\star} + filt_{\star}(\bm{u})(\bm{u}-\bm{u}_p^{\star})\big)\big|_{\bm{u}_p^{\star}}
			\end{array}\right) = 
			\left(\begin{array}{c}
				\bm{g}_{p}|_{\bm{u}_p^{\star}} \\
				\bm{g}_{p}|_{\bm{u}_p^{\star}}
			\end{array}\right) \leq \bm{0}.
		\end{align}			
		$\bm{\lambda}\in\amsmathbb{R}^{2n_g}$ is defined as a variable which one separates into two sub-vectors $(\bm{\lambda}^{\prime},\bm{\lambda}^{\prime\prime})\in\amsmathbb{R}^{n_g}$, as follows: 
		\begin{equation} \label{eq:3___16}
			\bm{\lambda} := 
			\left(\begin{array}{c}
				\bm{\lambda}^{\prime} \\
				\bm{\lambda}^{\prime\prime}
			\end{array}\right) := 
			\left(\begin{array}{c}
				\bm{\lambda}_p \\
				\bm{0} 
			\end{array}\right) \geq \bm{0}
		\end{equation}
		As $\bm{\lambda}_p\geq \bm{0}$, it is clear that $\bm{\lambda}:=$ \eqref{eq:3___16} implies that  $\bm{\lambda} \geq \bm{0}$.	 So: 
		\begin{align} \label{eq:3___27_kugyfcbx}
			\bm{\lambda}^{\rm T}\left(\begin{array}{c}
				\bm{g}^c_{\star}|_{\bm{u}_p^{\star}} \\
				\bm{g}_{\star}\big(\bm{u}_p^{\star} + filt_{\star}(\bm{u})(\bm{u}-\bm{u}_p^{\star})\big)\big|_{\bm{u}_p^{\star}}
			\end{array}\right) = 
			\bm{\lambda}_p^{\rm T} \bm{g}^c_{\star}|_{\bm{u}_p^{\star}} =
			\bm{\lambda}_p^{\rm T} \bm{g}_{p}|_{\bm{u}_p^{\star}} = 0.
		\end{align}
		In addition, as $\nabla_{\bm{u}} \phi_{p}|_{\bm{u}_p^{\star}} = \nabla_{\bm{u}} \phi^c_{\star}|_{\bm{u}_p^{\star}}$, and  $\nabla_{\bm{u}} \bm{g}_{p}|_{\bm{u}_p^{\star}} = \nabla_{\bm{u}} \bm{g}^c_{\star}|_{\bm{u}_p^{\star}}$,
		the following is observed:          
			\begin{align}
			& & \nabla_{\bm{u}} \phi_{p}|_{\bm{u}_p^{\star}} + 
			\bm{\lambda}_p^{\rm T}
			\nabla_{\bm{u}} \bm{g}_{p}|_{\bm{u}_p^{\star}} = \ & \bm{0}, \nonumber \\
			\Leftrightarrow && \nabla_{\bm{u}} \phi^c_{\star}|_{\bm{u}_p^{\star}} + 
			\bm{\lambda}_p^{\rm T}
			\nabla_{\bm{u}} \bm{g}^c_{\star}|_{\bm{u}_p^{\star}} = \ & \bm{0}, \nonumber \\
			\Leftrightarrow && \nabla_{\bm{u}} \phi^c_{\star}|_{\bm{u}_p^{\star}} + 
			\left(
			\begin{array}{c}
				\bm{\lambda}_p \\ 
				\bm{0}
			\end{array}
			\right)^{\rm T}
			\left(
			\begin{array}{c}
				\nabla_{\bm{u}} \bm{g}^c_{\star}|_{\bm{u}_p^{\star}} \\
				\nabla_{\bm{u}} \bm{g}_{\star}\big(\bm{u}_p^{\star} + filt_{\star}(\bm{u})(\bm{u}-\bm{u}_p^{\star})\big)\big|_{\bm{u}_p^{\star}}
			\end{array}
			\right) = \ & \bm{0}. \label{eq:3___18}
		\end{align}
		The union of equations \eqref{eq:3___15}-\eqref{eq:3___16}-\eqref{eq:3___27_kugyfcbx}-\eqref{eq:3___18} shows that if $(\bm{u}_p^{\star},\bm{\lambda}_p)$ is a KKT point of the plant, then  
		$(\bm{u}_p^{\star},\bm{\lambda}= [\bm{\lambda}_p^{\rm T}, \ \bm{0}^{\rm T}]^{\rm T})$ is a KKT point of Problem~\eqref{eq:3___1_New_PB_OPT}${}|_{\star}$.

		\begin{center}
			\textit{Part 2 -- 2nd order optimality condition:}
		\end{center} 
		\noindent
		As Problem~\eqref{eq:3___1_New_PB_OPT}${}|_{\star}$ is MFCQ at $\bm{u}_p^{\star}$, the 2nd order optimality condition at $\bm{u}_p^{\star}$ is, the existence of a set of KKT multipliers $\bm{\lambda}\in\amsmathbb{R}^{2n_g}$ such that the 1st order optimality conditions are satisfied and such that (\cite{Chong:2001}):  
		\begin{equation} \label{eq:3___29_iueyfc}
			\bm{N}_{\star}^{\rm T}
			\nabla^2_{\bm{uu}} \mathcal{L}_{\star}|_{\bm{u}_p^{\star}}
			\bm{N}_{\star} > \bm{0},
		\end{equation}
		where $\mathcal{L}_{\star}$ is the Lagrangian of Problem~\eqref{eq:3___1_New_PB_OPT}${}|_{\star}$, i.e. 
		\begin{equation} \label{eq:3___20}
			\mathcal{L}_{\star}(\bm{u}) =  \phi_{\star}(\bm{u}) + \sum_{i=1}^{n_g} \lambda^{\prime}_{(i)} g^c_{\star(i)}(\bm{u}) + 
			\lambda^{\prime\prime}_{(i)}g_{\star(i)}\big(\bm{u}_p^{\star} + filt_{\star}(\bm{u})(\bm{u}-\bm{u}_p^{\star})\big).
		\end{equation}
		
		With the candidate $\bm{\lambda}:=$ \eqref{eq:3___16}, for which the validity of the 1st order optimality conditions have been proved in part 1, \eqref{eq:3___20} becomes:
		\begin{equation} 
			\mathcal{L}_{\star}(\bm{u}) =  \phi^c_{\star}(\bm{u}) + \sum_{i=1}^{n_g} 
			\lambda_{p(i)}
			g^c_{\star(i)}(\bm{u}).
		\end{equation}
		According to Lemma~\ref{lem:3___1_CFD_MFCQ}, one can replace  $\bm{N}_{\star}$ by $\bm{N}_{p}$. Indeed, the fact that $CFD_{\star} = CFD_{p,\star}$ implies that the all the gradients of the active constraints of Problems~\eqref{eq:3___1_New_PB_OPT}${}|_{\star}$ and \eqref{eq:2_2___1_Plant_PB} 
		at $\bm{u}_p^{\star}$  belong to the same subspace of $\amsmathbb{R}^{n_u}$. Hence: $\bm{N}_{\star} = \bm{N}_{p}$. One can therefore conclude that the condition \eqref{eq:3___29_iueyfc} is:
		\begin{equation}
			\bm{N}_{\star}^{\rm T}
			\nabla^2_{\bm{uu}} \mathcal{L}_{\star}|_{\bm{u}_p^{\star}}
			\bm{N}_{\star} = 
			\bm{N}_{p}^{\rm T}
			\left(
			\nabla^2_{\bm{uu}}\phi^c_{\star}|_{\bm{u}_p^{\star}} + \sum_{i=1}^{n_g} 
			\lambda_{p(i)}
			\nabla^2_{\bm{uu}} g^c_{\star(i)}|_{\bm{u}_p^{\star}}
			\right)
			\bm{N}_{p}.
		\end{equation}
		As $\forall k$ $\nabla^2_{\bm{uu}}\phi^c_{k}|_{\bm{u}_k} > 0$ and $\nabla^2_{\bm{uu}}g^c_{k(i)}|_{\bm{u}_k} \geq  0$, this condition is valid at $\bm{u}_p^{\star}$. 

		\begin{center}
			\textit{Conclusion}
		\end{center}
		
		In summary, it is shown that the couple $(\bm{u}_p^{\star},\bm{\lambda}= [\bm{\lambda}_p^{\rm T}, \ \bm{0}^{\rm T}]^{\rm T})$ satisfies the 1st and 2nd order optimality conditions of Problem~\eqref{eq:3___1_New_PB_OPT}${}|_{\star}$.  The existence of this point is sufficient to show that KMFCA-D/I guaranties the satisfaction of the equilibrium condition in all circumstances.
\end{proofbox}

At this stage, it has been demonstrated that KMFCA-D/I provides the same guarantees as MFCA-D/I with respect to the validity of the equilibrium conditions. Now, moving on to the stability conditions. 

\subsection{Effects on the stability condition}

As KMFCA-D/I are purely iterative methods, the solution of problem \eqref{eq:3___1_New_PB_OPT} is a function $\bm{sol}^{\prime}$ of the point $\bm{v}$ around which the most recent experiments have been conducted. In more mathematical terms, \eqref{eq:3___1_New_PB_OPT} can be written as:
\begin{align} 
	\bm{sol}^{\prime}(\bm{v}) := \ &  \operatorname{arg}
	\underset{\bm{u}}{\operatorname{min}}  \quad    \widetilde{\phi}^c(\bm{u},\bm{v}), \nonumber\\ 
	& \qquad   \text{s.t.} \quad  \widetilde{\bm{g}}^c(\bm{u},\bm{v}) \leq \bm{0}, \nonumber \\
	& \phantom{\qquad   \text{s.t.} \quad} \widetilde{\bm{g}}\big(\bm{v} +  filt_k(\bm{u})(\bm{u} - \bm{v} ), \bm{v} \big) \leq \bm{0},  \label{eq:3___33_PBoptKMFCA}
\end{align}
where it is recalled that
\begin{align} \label{eq:3___34_def_fonction_tilde}
	\phi_k^c(\bm{u}) :=   \ &  \widetilde{\phi}^c(\bm{u},\bm{v}_k), &
	\bm{g}_k^c(\bm{u}) := \ &  \widetilde{\bm{g}}^c(\bm{u},\bm{v}_k), &
	\bm{g}_k(\bm{u}) := \ &  \widetilde{\bm{g}}(\bm{u},\bm{v}_k).
\end{align}
Then, a Theorem similar to Theorem~\ref{thm:2___3_Proprietes_Sol} can be formulated:
\begin{thmbox}
	\label{thm:3___3_ProprietesFoncytion_solplus}
	If: 
	\vspace{-\topsep}
	\begin{itemize}[noitemsep]
		\item A purely iterative RTO method using \eqref{eq:3___33_PBoptKMFCA} is used,
		\item For any correction point $\bm{v}\in\amsmathbb{R}^{n_u}$ the following equalities are true: 
		
		\begin{align*}
			\widetilde{\bm{g}}|_{\bm{v},\bm{v}}  =  \widetilde{\bm{g}}^c|_{\bm{v},\bm{v}}  = \ & \bm{g}_p|_{\bm{v}}, &
			\nabla_u \widetilde{\bm{g}}|_{\bm{v},\bm{v}}   = \nabla_u \widetilde{\bm{g}}^c|_{\bm{v},\bm{v}}   = \ & \nabla_u \bm{g}_p|_{\bm{v}},  & 
			\nabla_u \widetilde{\phi}|_{\bm{v},\bm{v}}  = \ & \nabla_u \phi_p|_{\bm{v}},
		\end{align*}
		\begin{align*} 
			\nabla^2_{\bm{u}\bm{u}}\widetilde{\phi}|_{\bm{v},\bm{v}}  + 
			\sum_{i=0}^{n_g}\left[
			\lambda^{\prime}_{(i)}  \widetilde{g}_{(i)}|_{\bm{v},\bm{v}} \right]  > \ & 0, & 
			\forall \bm{\lambda}^{\prime}\in\amsmathbb{R}^{n_g+},
		\end{align*}
		where the functions  $(\widetilde{\phi},\widetilde{\bm{g}})$  are defined by \eqref{eq:3___34_def_fonction_tilde}.
		\item For any correction point $\bm{v}\in\amsmathbb{R}^{n_u}$, Problem~\ref{eq:3___33_PBoptKMFCA} is such that $\bm{sol}^{\prime}(\bm{v})$ is unique.
	\end{itemize}
	Then, the following statements are true: 
	\begin{itemize}
		\item[\textbf{A. }] $\bm{sol}^{\prime}(\bm{v})$ is $\mathcal{C}^0$ and piecewise-$\mathcal{C}^1$ in the neighborhood of the KKT points of the plant;
		\item[\textbf{B. }] If $\bm{v}_p^{\bullet}$ is a KKT point of the plant and $\delta\bm{v}$ is a vector in $\amsmathbb{R}^{n_u}\backslash\bm{0}$, then $\bm{sol}^{\prime}(\bm{v}_p^{\bullet})= \bm{v}_p^{\bullet}$ and the value of the directional derivative of  $\bm{sol}^{\prime}$: 
		\begin{align}
			\nabla^{S\prime}(\bm{v},\delta\bm{v})  := \ &
			\left( \frac{\delta \bm{v} }{\|\delta \bm{v} \|}
			\right)^{\rm T}
			\lim_{\begin{smallmatrix} \alpha \to 0 & \\ \alpha>0 \end{smallmatrix}} \frac{
				\bm{sol}^{\prime}(\bm{v}+\alpha\delta \bm{v})
				- \bm{sol}^{\prime}(\bm{v}) }{\alpha}, \\
			\Big( = \ &
			\left( \frac{\delta \bm{v}}{\|\delta \bm{v}\|}
			\right)^{\rm T} \nabla_{v}\bm{sol}^{\prime}|_{\bm{v}} \frac{\delta \bm{v}}{\|\delta \bm{v}\|}, \quad \text{ si $\bm{sol}^{\prime}$ est $\mathcal{C}^1$ en $\bm{v}$}
			\Big). \nonumber
		\end{align}
		\begin{equation}
			\text{is: }
			\left\{
			\begin{array}{ll}
				\nabla^{S\prime}(\bm{v}_p^{\bullet},\delta\bm{v}) < 1, & \text{if $\bm{v}_p^{\bullet}$ is a minimum of the plant}, \\
				\nabla^{S\prime}(\bm{v}_p^{\bullet},\delta\bm{v}) > 1, & \text{if $\bm{v}_p^{\bullet}$ is a maximum of the plant}, \\
				\nabla^{S\prime}(\bm{v}_p^{\bullet},\delta\bm{v}) \gtrless 1, & \text{if $\bm{v}_p^{\bullet}$ is a saddle-point of the plant},
			\end{array}
			\right. \qquad \quad 
		\end{equation}
		\item[\textbf{C. }]
		If $\bm{v}_p^{\star}$ is a minimum of the plant and $\delta\bm{v}$ is a vector  orthogonal to the set of active constraints of the plant at $\bm{v}_p^{\star}$ such that $\nabla^{S+}(\bm{v}_p^{\star}\delta\bm{v}) <-1$, then the Lagrangian of the plant must be at least two times more convex than the one of the model in the direction $\delta\bm{v}$, i.e.: 
		\begin{equation} 
			\delta\bm{v}^{\rm T} \nabla^2_{\bm{uu}} (\mathcal{L}_p|_{\bm{v}_p^{\star}} - 2 \mathcal{L}|_{\bm{v}_p^{\star},\bm{v}_p^{\star}}) \delta\bm{v} >0
			\ \Rightarrow 
			\ \nabla^{S\prime}(\bm{v}_p^{\star},\delta\bm{v})<-1.
		\end{equation}
	\end{itemize}
\end{thmbox}
\begin{proofbox}
	The proof of this Theorem is very similar to the one of Theorem~\ref{thm:2___3_Proprietes_Sol}. The only significant differences are that (i) the names and the number of  constraint are not the same, and (ii) the constraints of Problem~\eqref{eq:3___33_PBoptKMFCA} are not LICQ but MFCQ. This only affects the fact that in this case the KKT multipliers are no longer functions of $\bm{v}$. 
\end{proofbox}

When this theorem is applicable one can link two successive iterations $\bm{u}_k = \bm{u}_p^{\star} + \delta\bm{u}$ and $\bm{u}_{k+1}$ with the following equation: 
\begin{align} 
	& &\bm{u}_{k+1} = \ & \bm{sol}^{\prime}(\bm{u}_p^{\star} + \delta\bm{u}_k), \nonumber \\
	& &             = \ & \bm{sol}^{\prime}(\bm{u}_p^{\star}) + 
	\nabla^{S\prime}(\bm{u}_p^{\star},\delta\bm{u}_k) \delta\bm{u}_k + \mathcal{O}(\| \delta \bm{u}_k\|^2), \nonumber \\
	& &			 = \ & \bm{u}_p^{\star} + 
	\nabla^{S\prime}(\bm{u}_p^{\star},\delta\bm{u}_k) \delta\bm{u}_k + \mathcal{O}(\| \delta \bm{u}_k\|^2), \label{eq:3___40_C1_Uasteri} \\
	\Leftrightarrow & & 
	\bm{u}_{k+1} - \bm{u}_p^{\star} = \ &  \nabla^{S\prime}(\bm{u}_p^{\star},\delta\bm{u}_k) (\bm{u}_{k} - \bm{u}_p^{\star}) + \mathcal{O}(\| \bm{u}_{k} - \bm{u}_p^{\star}\|^2). \label{eq:3___41_rcugbk}
\end{align}
Clearly,  one can notice the great similarity between  \eqref{eq:3___40_C1_Uasteri}-\eqref{eq:3___41_rcugbk}
and \eqref{eq:2___28_C1_Uasteri}-\eqref{eq:2___33_rcugbk}. Therefore, an argument similar to the one used in Chapter~\ref{Chap:2_Vers_Une_meilleure_Convergence} starting from  equation \eqref{eq:2___33_rcugbk} can be used  to demonstrate two things. First:

\begin{thmbox}
	If the conditions of applicability of the Theorem~\ref{thm:3___3_ProprietesFoncytion_solplus} are satisfied, and if:
	\vspace{-\topsep}
	\begin{itemize}[noitemsep]
		\item The problem is 1-D ($n_u = 1$); 
		\item The function $sol^{\prime}$ is $\mathcal{C}^1$ at $\bm{u}_p^{\star}$; 
		\item The filter $K>0$ is applied as follows: 
		\begin{equation*}
			u_{k+1} = u_k + K(sol(u_k)-u_k);
		\end{equation*}
	\end{itemize}
	Then, any filter adaptation strategy that satisfies:
	 \begin{align} \label{eq:3___42_nablasplus}
		K_k < \ &  \frac{2}{1-\utilde{\text{\hskip 0.1ex $\nabla$}}^{S}_{k}}, & 
		\utilde{\text{\hskip 0.1ex $\nabla$}}^{S}_{k}  := \ &
		\left( \frac{\bm{u}_k - \bm{u}_{k-1}}{\|\bm{u}_k - \bm{u}_{k-1}\|}
		\right)^{\rm T}
		\frac{\bm{u}^{*}_{k+1} - \bm{u}^{*}_{k}}{\| \bm{u}_{k} - \bm{u}_{k-1} \| }, &
		\text{if } \utilde{\text{\hskip 0.1ex $\nabla$}}_k^S < 1.
	\end{align}	
 	guarantee that the stability condition is satisfied.
\end{thmbox}
\begin{proofbox}
	The proof of this theorem is identical to the proof of Theorem~\ref{thm:2___4_K_Niveau2CasParticulier} where one would replace the function $sol$ by $sol^{\prime}$.
\end{proofbox}
\noindent
Second, the ``optimal'' filter that maximizes the convergence speed on $\bm{u}_p^{\star}$ when one is close to $\bm{u}_p^{\star}$  is:
\begin{align} \label{eq:3___44_OptimizationPB_ChoiceK}
	K^{\star}_k = \ & \max\left\{K_o,   \frac{1}{1- \utilde{\text{\hskip 0.1ex $\nabla$}}^{S}_{k}} \right\}, &
	\text{if } \utilde{\text{\hskip 0.1ex $\nabla$}}_k^S < 1. 
\end{align}

Finally, if \eqref{eq:2___41_strat1_adaptationK} and \eqref{eq:3___42_nablasplus} are compared, then it can be seen that $\utilde{\text{\hskip 0.1ex $\nabla$}}^{S\prime}_{k}$ and $\utilde{\text{\hskip 0.1ex $\nabla$}}^{S}_{k}$ have the same definitions. So, this filter adaptation strategy, based on the functions $filt_k$ and $nab_k$, guarantees that KMFCA-D/I validates the stability conditions. Moreover, as the optimization problem \eqref{eq:3___1_New_PB_OPT} specifies that the filtered point based on the $filt_k$ and $nab_k$ must be feasible according to the updated model, the application of the filter does not present the danger shown in Figure~\ref{fig:3___1_Concept_idea2}. Finally, all the observations and graphical interpretations made in the previous chapter for MFCA-D/I are applicable to KMFCA-D/I.  Now, looking into how the performances of KMFCA-D/I can be improved.   

\section{Management of irrelevant constraints} 

Convexifying the concave constraints enforces the satisfaction of the equilibrium condition, however this action can also have negative effects on the performance of the algorithm. The following example illustrates this point. 

\begin{exbox} \label{ex:3___2_Contournement_Contrainte}
	\textbf{(Bypassing an inactive constraint)}
	Considering the following theoretical RTO problem: 
	\begin{align*}
		\phi_p(\bm{u}) := \ & - 
		\left(\begin{array}{c}
			4 \\
			2
		\end{array}\right)^{\rm T}\bm{u} + \frac{1}{2}\bm{u}^{\rm T} 
		\left(\begin{array}{cc}
			8 & 0 \\
			0 & 4
		\end{array}\right)\bm{u},  \\
		\phi(\bm{u}) := \ &  - 
		\left(\begin{array}{c}
			0.1 \\
			0.1
		\end{array}\right)^{\rm T}\bm{u} + \frac{1}{2}\bm{u}^{\rm T} 
		\left(\begin{array}{cc}
			2 & 0 \\
			0 & 2
		\end{array}\right)\bm{u}, \\
		g(\bm{u}) := g_p(\bm{u}) := \ &  -2 + 
		\left(\begin{array}{c}
			6 \\
			10
		\end{array}\right)^{\rm T}\bm{u} - \frac{1}{2}\bm{u}^{\rm T} 
		\left(\begin{array}{cc}
			20 & 0 \\
			0  & 36
		\end{array}\right)\bm{u}, 
	\end{align*}
	where $u_{(1)}\in[0,1]$ and $u_{(2)}\in[0,1]$. Figure~\ref{fig:Exemple_3_2_maps} illustrates those functions.\\
	\begin{minipage}[h]{\linewidth}
		\vspace*{0pt}
		\centering 
		\includegraphics[width=6cm]{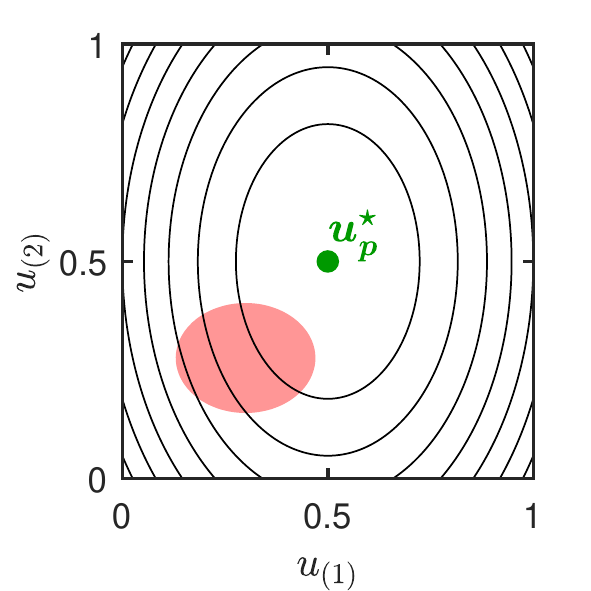}
		\includegraphics[width=6cm]{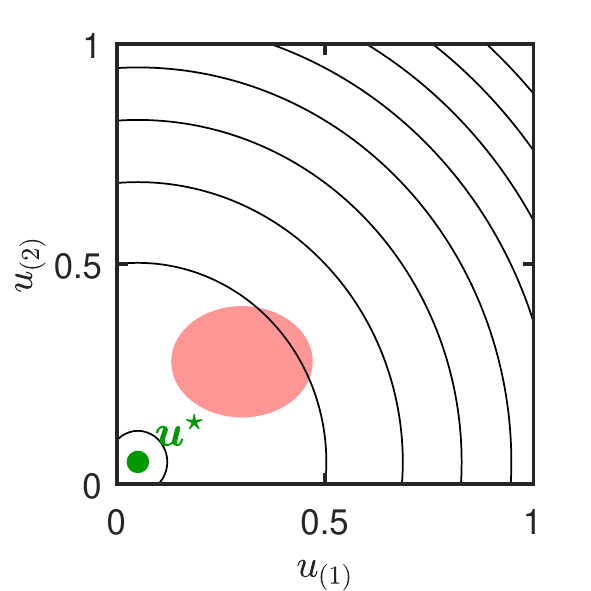}
		\vspace{-3mm}
		\captionof{figure}{Plant (left) and model's (right) cost and constraint functions}
		\label{fig:Exemple_3_2_maps}
	\end{minipage}\\
	
	Initializing MFCA and KMFCA at $\bm{u}^{\star}$ and launching several simulations. The results of these simulations are given in Figure~\ref{fig:3___9_Exemple_3_2_Results1}. \\
	
	\begin{minipage}[h]{\linewidth}
		\vspace*{0pt}
		\centering
		\includegraphics[width=4.45cm]{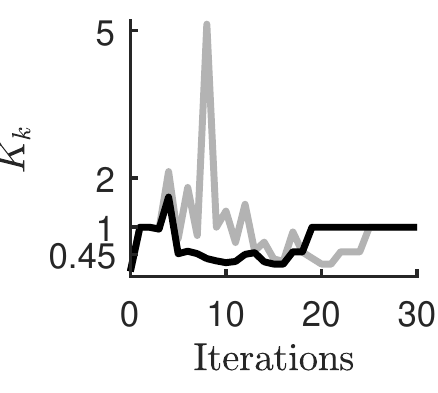}\hskip -0ex
		\includegraphics[width=4.45cm]{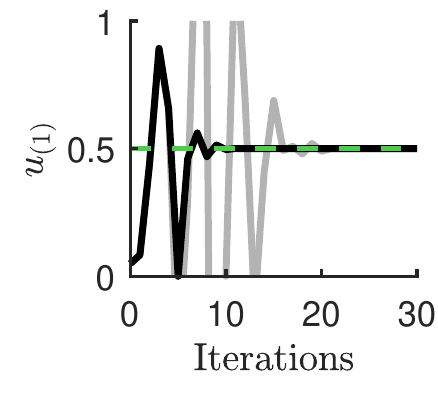}\hskip -0ex
		\includegraphics[width=4.45cm]{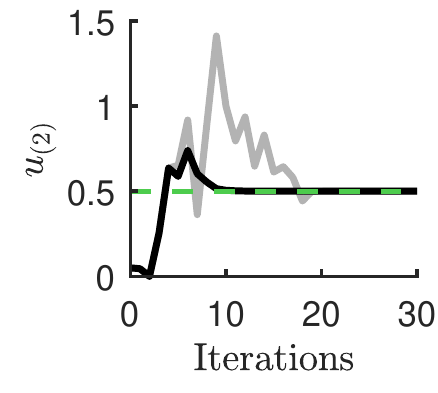} \\
		\includegraphics[width=4.45cm]{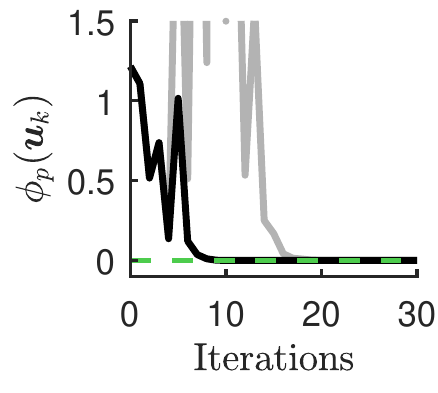}\hskip -0ex
		\includegraphics[width=4.45cm]{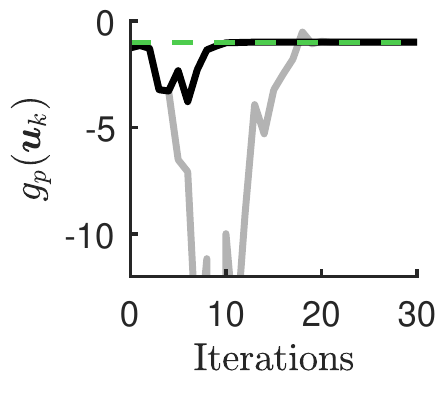}\hskip -0ex
		\includegraphics[width=4.45cm]{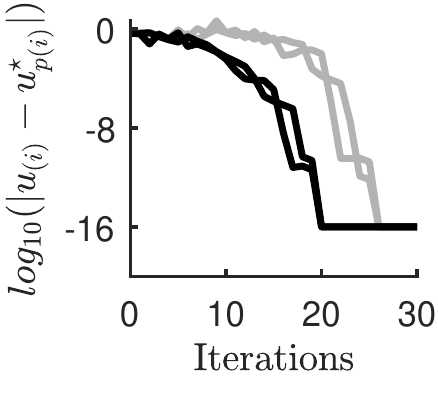}\\
		\textcolor{gray}{\raisebox{0.5mm}{\rule{0.5cm}{0.1cm}}}    : \text{MFCA,} \ 
		\textcolor{black}{\raisebox{0.5mm}{\rule{0.5cm}{0.1cm}}} : \text{KMFCA,} \ 
		\textcolor{green1}{\raisebox{0.5mm}{\rule{0.2cm}{0.1cm}\hspace{0.1cm}\rule{0.2cm}{0.1cm}}} : Solution
		\vspace{-2mm}
		\captionof{figure}{Simulation results.}
		\label{fig:3___9_Exemple_3_2_Results1}
	\end{minipage} \\
	
	Two main observations can be made. First, MFCA does not violate the constraint $g_p$ but the boundaries on the inputs $u_{(1)}$ and $u_{(2)}$ are violated.  The extent of the violation on $u_{(1)}$ is not shown because it would affect the readability of the graph (too large), but the one on $u_{(2)}$ is clearly shown.  As for the previous example, this critical violation of perfectly known constraints is due to the absence of a constraint on the choice of the \textit{applied} point $\bm{u}_{k+1}$.  To better understand what MFCA does at each iteration, Figure~\ref{fig:3___10_Exemple_3_2_Results2} illustrates the optimization problems solved by MFCA for the first 24 iterations.\\
	\textit{(For each sub-figure in Figure~\ref{fig:3___10_Exemple_3_2_Results2} the horizontal and vertical axes correspond to $u_{(1)}$ and $u_{(2)}$, respectively. To increase the readability of these sub-figures, the axes are not marked and not labeled.)}\\
	
	\begin{minipage}[h]{\linewidth}
		\vspace*{0pt}
		{\centering
			\includegraphics[trim={0.35cm 0.25cm 0.2cm  0.15cm },clip,width=2.75cm]{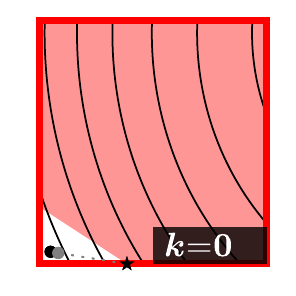}\hskip -0ex
			\includegraphics[trim={0.35cm 0.25cm 0.2cm  0.15cm },clip,width=2.75cm]{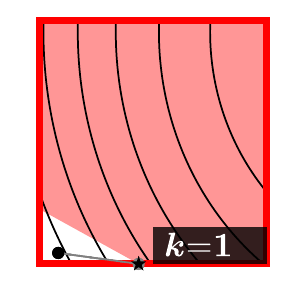}\hskip -0ex
			\includegraphics[trim={0.35cm 0.25cm 0.2cm  0.15cm },clip,width=2.75cm]{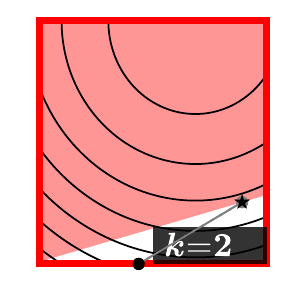}\hskip -0ex
			\includegraphics[trim={0.35cm 0.25cm 0.2cm  0.15cm },clip,width=2.75cm]{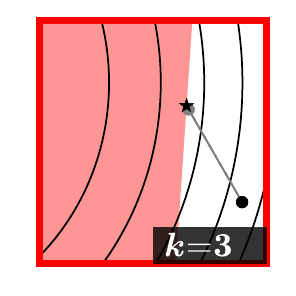}\hskip -0ex
			\includegraphics[trim={0.35cm 0.25cm 0.2cm  0.15cm },clip,width=2.75cm]{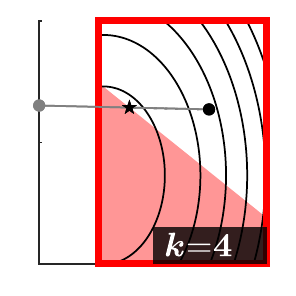} \\
			\includegraphics[trim={0.35cm 0.25cm 0.2cm  0.15cm },clip,width=2.75cm]{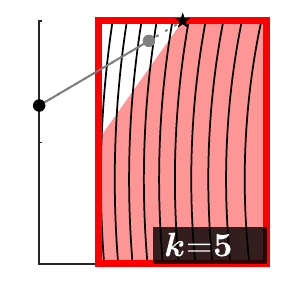}\hskip -0ex
			\includegraphics[trim={0.35cm 0.25cm 0.2cm  0.15cm },clip,width=2.75cm]{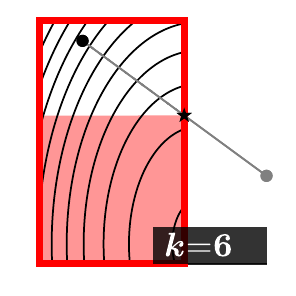}\hskip -0ex
			\includegraphics[trim={0.35cm 0.25cm 0.2cm  0.15cm },clip,width=2.75cm]{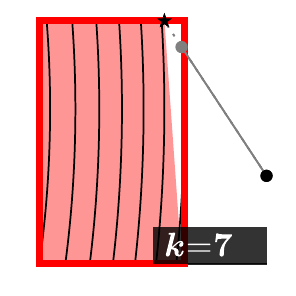}\hskip -0ex
			\includegraphics[trim={0.35cm 0.25cm 0.2cm  0.15cm },clip,width=2.75cm]{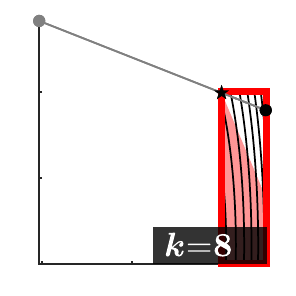}\hskip -0ex
			\includegraphics[trim={0.35cm 0.25cm 0.2cm  0.15cm },clip,width=2.75cm]{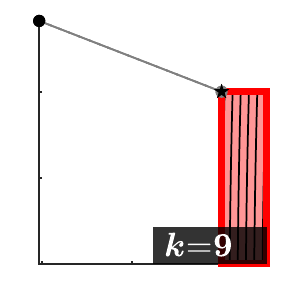} \\
			\includegraphics[trim={0.35cm 0.25cm 0.2cm  0.15cm },clip,width=2.75cm]{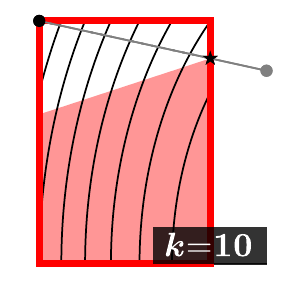}\hskip -0ex
			\includegraphics[trim={0.35cm 0.25cm 0.2cm  0.15cm },clip,width=2.75cm]{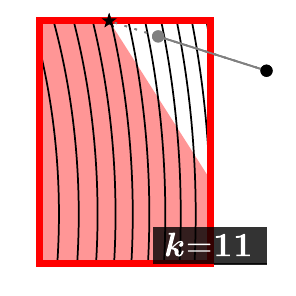}\hskip -0ex
			\includegraphics[trim={0.35cm 0.25cm 0.2cm  0.15cm },clip,width=2.75cm]{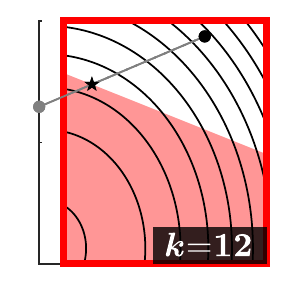}\hskip -0ex
			\includegraphics[trim={0.35cm 0.25cm 0.2cm  0.15cm },clip,width=2.75cm]{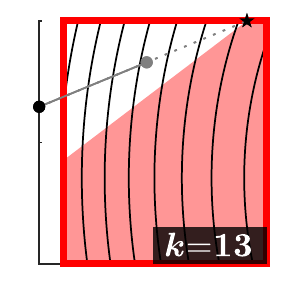}\hskip -0ex
			\includegraphics[trim={0.35cm 0.25cm 0.2cm  0.15cm },clip,width=2.75cm]{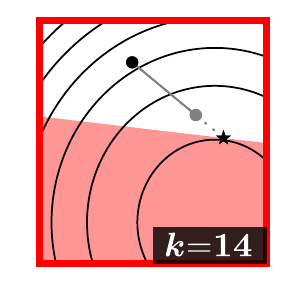} \\
			\includegraphics[trim={0.35cm 0.25cm 0.2cm  0.15cm },clip,width=2.75cm]{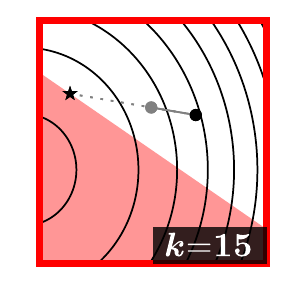}\hskip -0ex
			\includegraphics[trim={0.35cm 0.25cm 0.2cm  0.15cm },clip,width=2.75cm]{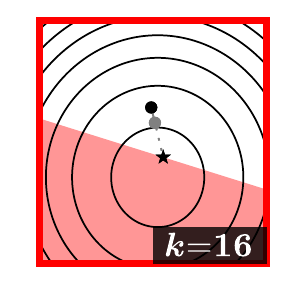}\hskip -0ex
			\includegraphics[trim={0.35cm 0.25cm 0.2cm  0.15cm },clip,width=2.75cm]{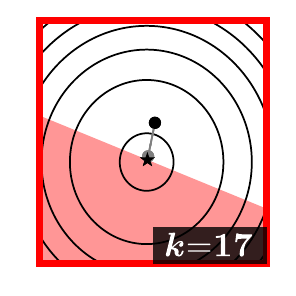}\hskip -0ex
			\includegraphics[trim={0.35cm 0.25cm 0.2cm  0.15cm },clip,width=2.75cm]{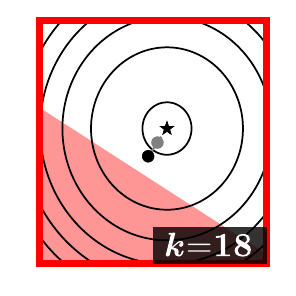}\hskip -0ex
			\includegraphics[trim={0.35cm 0.25cm 0.2cm  0.15cm },clip,width=2.75cm]{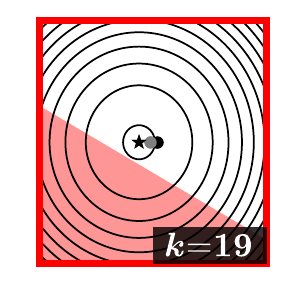} \\
			\includegraphics[trim={0.35cm 0.25cm 0.2cm  0.15cm },clip,width=2.75cm]{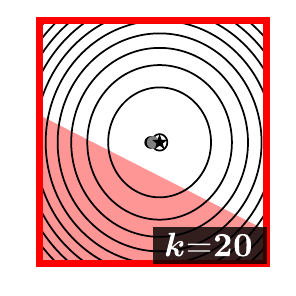}\hskip -0ex
			\includegraphics[trim={0.35cm 0.25cm 0.2cm  0.15cm },clip,width=2.75cm]{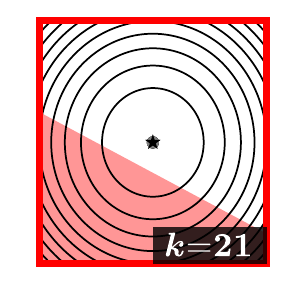}\hskip -0ex
			\includegraphics[trim={0.35cm 0.25cm 0.2cm  0.15cm },clip,width=2.75cm]{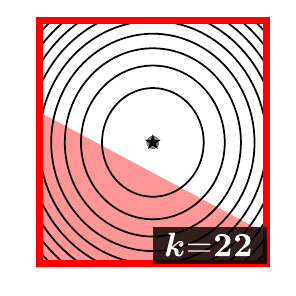}\hskip -0ex
			\includegraphics[trim={0.35cm 0.25cm 0.2cm  0.15cm },clip,width=2.75cm]{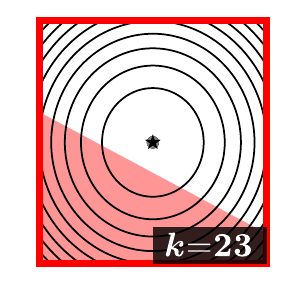}\hskip -0ex
			\includegraphics[trim={0.35cm 0.25cm 0.2cm  0.15cm },clip,width=2.75cm]{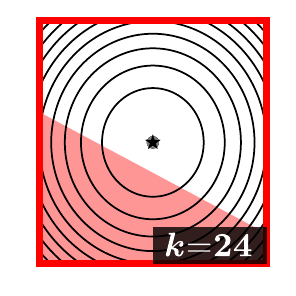} \\
		}
		$\bullet$: $\bm{u}_k$, \
		\textcolor{gray}{$\bullet$}: $\bm{u}_{k+1}$, \ $\bm{\star}$: $\bm{u}^{\star}_{k+1}$, \
		\resizebox{0.5cm}{0.3cm}{\textcolor{red}{$\bm{\square}$}}: Boundaries on $\bm{u}^{\star}_{k+1}$, \ 
		\textcolor{red!50!white}{\raisebox{-0.5mm}{\rule{0.5cm}{0.3cm}}}: Area forbidden to $\bm{u}^{\star}_{k+1}$
		\captionof{figure}{Details of the iterations of \textbf{MFCA}.}
		\label{fig:3___10_Exemple_3_2_Results2}
	\end{minipage} \\

	The first thing one notices on Figure~\ref{fig:3___10_Exemple_3_2_Results2} is that some iterations leave the red frame (the boundaries on $\bm{u}$). This is something that in practice would be totally unacceptable and would immediately disqualify MFCA from any real implementation. Ignoring this problem, a second observation is that convexifying the constraint induces a spreading of the constraint over a large domain of the input space which in turn induces a progressive (i.e. slow) bypass.  \\
	
	Proceeding to the detailed analysis of the KMFCA iterations based on Figure~\ref{fig:3___11_Exemple_3_2_Results3}. It can be observed that each iteration respects the boundaries on the inputs as well as the constraint $g_k$. So, apparently the filter-based constraints solves the main problem of MFCA. However, KMFCA suffers from the same weakness as MFCA in its inability to rapidly bypass the constrained sector. It should be recalled that this weakness is present despite the fact that (i) there is no modeling error on the constraint, and (ii) the constraint is not active at $\bm{u}_p^{\star}$.\\
	\begin{minipage}[h]{\linewidth}
		\vspace*{0pt}
		{\centering
			\includegraphics[trim={0.35cm 0.25cm 0.2cm  0.15cm },clip,width=2.75cm]{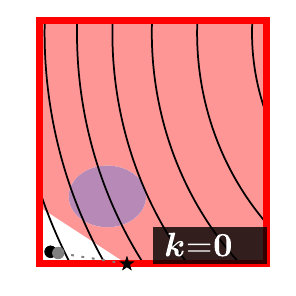}\hskip -0ex
			\includegraphics[trim={0.35cm 0.25cm 0.2cm  0.15cm },clip,width=2.75cm]{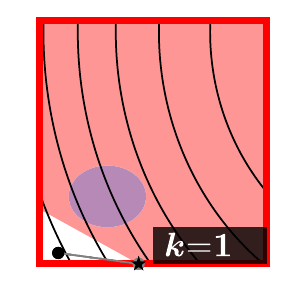}\hskip -0ex
			\includegraphics[trim={0.35cm 0.25cm 0.2cm  0.15cm },clip,width=2.75cm]{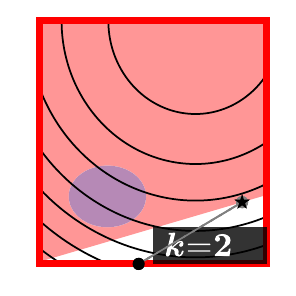}\hskip -0ex
			\includegraphics[trim={0.35cm 0.25cm 0.2cm  0.15cm },clip,width=2.75cm]{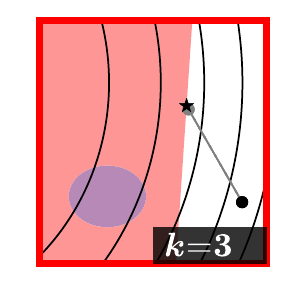}\hskip -0ex
			\includegraphics[trim={0.35cm 0.25cm 0.2cm  0.15cm },clip,width=2.75cm]{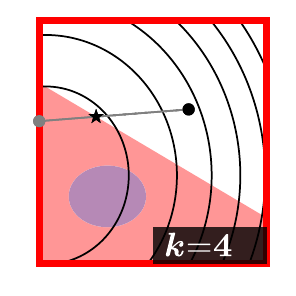} \\
			\includegraphics[trim={0.35cm 0.25cm 0.2cm  0.15cm },clip,width=2.75cm]{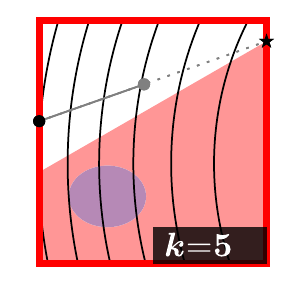}\hskip -0ex
			\includegraphics[trim={0.35cm 0.25cm 0.2cm  0.15cm },clip,width=2.75cm]{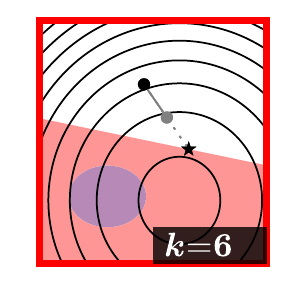}\hskip -0ex
			\includegraphics[trim={0.35cm 0.25cm 0.2cm  0.15cm },clip,width=2.75cm]{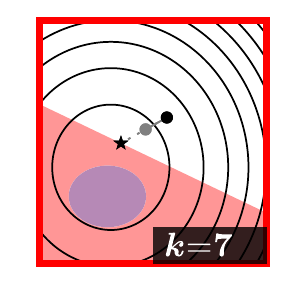}\hskip -0ex
			\includegraphics[trim={0.35cm 0.25cm 0.2cm  0.15cm },clip,width=2.75cm]{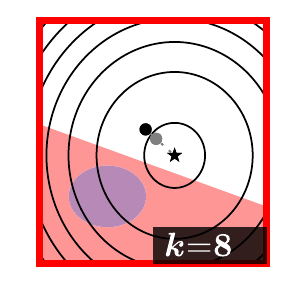}\hskip -0ex
			\includegraphics[trim={0.35cm 0.25cm 0.2cm  0.15cm },clip,width=2.75cm]{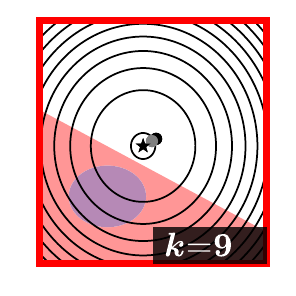}  \\
		}$\bullet$:$\bm{u}_k$, \  
		\textcolor{gray}{$\bullet$}: $\bm{u}_{k+1}$, \ $\bm{\star}$:$\bm{u}^{\star}_{k+1}$, \
		\resizebox{0.5cm}{0.3cm}{\textcolor{red}{$\bm{\square}$}}: 
		Boundaries on $\bm{u}^{\star}_{k+1}$ and $\bm{u}_{k+1}$,
		\
		\textcolor{red!50!white}{\raisebox{-0.5mm}{\rule{0.5cm}{0.3cm}}}:
		Area forbidden to $\bm{u}^{\star}_{k+1}$, 
		\
		\textcolor{magenta_clair2}{\raisebox{-0.5mm}{\rule{0.5cm}{0.3cm}}}:
		Area forbidden to $\bm{u}^{\star}_{k+1}$ and $\bm{u}_{k+1}$.
		\captionof{figure}{Details of the iterations of \textbf{KMFCA}.}
		\label{fig:3___11_Exemple_3_2_Results3}
	\end{minipage} 
\end{exbox}

Example~\ref{ex:3___1_Risque_de_lineariser_contraintes_concave} confirms the critical importance of checking  the feasibility of the updated model at $\bm{u}_{k+1}$ at the iteration  $k$.   Furthermore, it clearly disqualifies MFCA from any practical application by showing that even the boundaries on the plant inputs can be violated.
  However, the systematic convexification of the set of concave constraints at $\bm{u}_k$ can give disproportionate importance to constraints that would otherwise be irrelevant in the decision making. To avoid such performance loss, one idea could be to linearize a concave constraint at  $\bm{u}_k$ if, and only if, it has an effect on decision making, i.e. if it is predicted to be active at $\bm{u}_{k+1}^{\star}$ or at  $\bm{u}_{k+1}$. The following  procedure to identify the constraints to  convexify is proposed: 
\begin{BoxFun*}{f
		\vspace{-5mm}
		\begin{equation} \label{eq:3___43_CaID}
			\hspace{-4mm}\text{$[A] = $ \textsc{CaID}$(\bm{g}_k,n_g, \bm{u}_{k+1},\bm{u}^{\star}_{k+1}, b)$} \qquad \qquad \qquad 
		\end{equation}
	}
	\textbf{Inputs:} $\bm{g}_k$: the constraint functions of the  model updated at $\bm{u}_k$; $n_g$: the number of constraints; $\bm{u}_{k+1}^\star,\bm{u}_{k+1}$: the optimum of the updated model and the filtered point; and $b$: a scalar to be chosen smaller than $0$ and which tends to $0$. \\ 
	\textbf{Outputs:} $A$: the indices of the concave constraints that should influence the decision making (the ones that should be convexified).
	\tcblower
	$A = [\ ]$; \\
	\textbf{\textcolor{blue}{for}} $i = 1:n_g$
	\vspace{-\topsep}\begin{itemize}[noitemsep]
		\item[] \textbf{\textcolor{blue}{if}} $\nabla_{\bm{uu}}^2g_{k(i)}|_{\bm{u}_k} \not\geq 0$ \textbf{\textcolor{blue}{\&\&}} $\big( g_{k(i)}|_{\bm{u}_{k+1}} \geq b$ 
		\textbf{\textcolor{blue}{||}} $g_{k(i)}|_{\bm{u}_{k+1}^{\star}} \geq b \big)$
		\begin{itemize}[noitemsep]
			\item[] $A = [A;i] $
		\end{itemize}
		\item[] \textbf{\textcolor{blue}{end}}
	\end{itemize}\vspace{-\topsep}
	\textbf{\textcolor{blue}{end}}	
\end{BoxFun*}

 If one applies this idea to KMFCA-D/I, then two new methods that avoid performance losses due to constraint bypass are obtained. These two methods are called \textit{modifier-filter-(active)curvature adaptation direct with filter-based constraints (KMFCaA-D)} and  \textit{modifier-filter-(active)curvature adaptation indirect with filter-based constraints (KMFCaA-I)} and are detailed hereafter:

\begin{BoxAlgo}{KMFCaA-D}{KMFCaA}
	\textbf{Initialization.} Provide $\bm{u}_0$, $\bm{K}_0(=0.1)$, $\bm{u}^*_0(=\bm{u}_0)$, $\bm{a}$, functions $(\bm{f},\phi,\bm{g})$, and the stopping criterion of step 5).
	\tcblower
	\textbf{for} $k=0 \rightarrow \infty$
	\begin{itemize}[noitemsep]
		\item[1) ]  \textbf{Measure} $(\nabla_{\bm{u}}\phi_p, \bm{g}_p,\nabla_{\bm{u}}\bm{g}_p)|_{\bm{u}_{k}}$ on the plant.
		\item[2) ] \textbf{Update} $\phi$ and $\bm{g}$ with \eqref{eq:3___19_Corrections_ISO}. 
		\item[3) ]  \textbf{Compute}  $\phi^c_k$ and $\bm{g}^c_k)$ with \eqref{eq:3___30_DefinitionMatrices_P}.
		\item[4) ]  \textbf{Compute} $(\bm{u}^{\star}_{k+1},K_k,\bm{u}_{k+1})$ as follows:
	\end{itemize}
	\vspace{-\topsep}\begingroup
	\allowdisplaybreaks
	\begin{flalign} 
		\bm{u}^{\star}_{k+1} := \ &  \operatorname{arg}
		\underset{\bm{u}}{\operatorname{min}}  \quad    \phi^c_{k}(\bm{u}), \nonumber\\ 
		& \qquad   \text{s.t.} \quad  \bm{g}_{k}(\bm{u}) \leq \bm{0}, \nonumber \\
		& \phantom{\qquad   \text{s.t.} \quad} \bm{g}_{k}\big(\bm{u}_k +  filt_k(\bm{u})(\bm{u} - \bm{u}_k ) \big) \leq \bm{0}, \label{eq:3___43_New_PB_OPT} \\
		K_k := \ & filt_k(\bm{u}^\star_{k+1}), \label{eq:3___44_New_filter}\\
		\bm{u}_{k+1} := \ & \bm{u}_k + K_k (\bm{u}^\star_{k+1} - \bm{u}_k), \label{eq:3___45_New_point} \\
		\text{where}: \qquad\quad  
		\nonumber \\
		\hspace{-6mm} filt_k(\bm{u}) := \ &  \left\{
		\begin{array}{l@{}l}
			\hspace{-2mm} K_o,     & \text{if } k =0, \\
			\hspace{-2mm}	K_{k-1},  & \text{if } nab_k(\bm{u}) < 1 \ \& \ k > 0   \ \& \  |\bm{u}_{k}-\bm{u}_{k-1} |\leq \bm{a}  \\
			\hspace{-2mm} \max\left\{K_o,\frac{1}{1-nab_k(\bm{u})}\right\}\hspace{-1mm}, \
			&  \text{if } nab_k(\bm{u}) < 1   \ \& \  k > 0 \   \& \ |\bm{u}_{k}-\bm{u}_{k-1}|>\bm{a},   \\
			\hspace{-2mm} 1, & \text{otherwise.}
		\end{array}
		\right. \label{eq:3___46_fonction_filt} \\
		\hspace{-3mm}	nab_k(\bm{u})  := \ &
		\frac{\left(\bm{u}_k - \bm{u}_{k-1}\right)^{\rm T} 
			\left(\bm{u} - \bm{u}^{*}_{k} \right)
		}{\|\bm{u}_k - \bm{u}_{k-1}\|^2}. \label{eq:3___48_fonction_nab}
	\end{flalign}\endgroup
	\begin{itemize}[noitemsep]
		\item[5) ]  \textbf{Identify}  with the function \textsc{CaID} the set $A$ of constraints that are concave at $\bm{u}_k$ and active at $\bm{u}_{k+1}$ or $\bm{u}_{k+1}^{\star}$.
		\item[6) ] \textbf{If} $A\neq\emptyset$  \textbf{recompute} $(\bm{u}^{\star}_{k+1},K_k,\bm{u}_{k+1})$ with:
		\begin{flalign} 
			\bm{u}^{\star}_{k+1} := \ &  \operatorname{arg}
			\underset{\bm{u}}{\operatorname{min}}  \quad    \phi^c_{k}(\bm{u}), \nonumber \\ 
			& \qquad   \text{s.t.} \quad  g_{k(i)}(\bm{u}) \leq 0, & & \forall i \in \{1,..,n_g\} \backslash A \nonumber \\
			& \phantom{\qquad   \text{s.t.} \quad} g^c_{k(i)}(\bm{u}) \leq 0, & & \forall i \in \{1,..,n_g\} \cap A, \nonumber \\
			& \phantom{\qquad   \text{s.t.} \quad} \bm{g}_{k}\big(\bm{u}_k +  K_k (\bm{u})(\bm{u} - \bm{u}_k ) \big) \leq \bm{0}, \nonumber \\
			K_k := \ & \text{\eqref{eq:3___44_New_filter}}, \quad  
			\bm{u}_{k+1} := 
			 \text{\eqref{eq:3___45_New_point}}.\label{eq:4___49_Recalculer_u_K_u}
		\end{flalign}
		\item[7) ] \textbf{Stop} if $\bm{u}_{k+1}\approx\bm{u}_{k}$ and return $\bm{u}_{\infty} := \bm{u}_{k+1}$.
	\end{itemize}
	\noindent {\bf end}
\end{BoxAlgo}

\begin{BoxAlgo}{KMFCaA-I}{KMFCaAy}
	\textbf{Initialization.}  Provide $\bm{u}_0$, $\bm{K}_0(=0.1)$, $\bm{u}^*_0(=\bm{u}_0)$, functions $(\bm{f},\phi,\bm{g})$, and the stoping criteron of step 6).
	\tcblower
	\textbf{for} $k=0 \rightarrow \infty$
	\begin{itemize}[noitemsep]
		\item[1) ]  \textbf{Measure} $(\bm{f}_p,\nabla_{\bm{u}}\bm{f}_p)|_{\bm{u}_{k}}$ on the plant.
		\item[2) ]
		\textbf{Update} the functions $\bm{f}$ and $\phi,\bm{g}$ with \eqref{eq:3___20_Corrections_ISOy}.  
		\item[3) ]  \textbf{Compute}  $\phi^c_k$ and $\bm{g}^c_k$ with \eqref{eq:3___30_DefinitionMatrices_P}.
		\item[4) ]  \textbf{Compute}  $(\bm{u}^{\star}_{k+1},K_k,\bm{u}_{k+1})$ with \eqref{eq:3___43_New_PB_OPT}, \eqref{eq:3___44_New_filter}, and  \eqref{eq:3___45_New_point}.
		\item[5) ]  \textbf{Identify} with the function \textsc{CaID} the set $A$ of constraints that are concave at $\bm{u}_k$ and active at $\bm{u}_{k+1}$ or $\bm{u}_{k+1}^{\star}$.
		\item[6) ] \textbf{If} $S\neq\emptyset$  \textbf{recompute} $(\bm{u}^{\star}_{k+1},K_k,\bm{u}_{k+1})$ with \eqref{eq:4___49_Recalculer_u_K_u}.
		\item[7) ] \textbf{Stop} if $\bm{u}_{k+1}\approx\bm{u}_{k}$ and return $\bm{u}_{\infty} := \bm{u}_{k+1}$.
	\end{itemize}
	\noindent {\bf end}
\end{BoxAlgo}

Let's try  this approach on the problem of  example~\ref{ex:3___1_Risque_de_lineariser_contraintes_concave} to validate the ability of KMFCaA-D/I to efficiently bypass constrained domains.

\begin{exbox}	
	Consider the RTO problem of Example~\ref{ex:3___1_Risque_de_lineariser_contraintes_concave} and apply KMFCaA. The simulation results are given in Figure~\ref{fig:3___10_Exemple_3_3_Results1} where they can be compared them with those obtained with MFCA and KMFCA. It is clear that KMFCaA performs much better than the other two methods. Figure~\ref{fig:3___11_Exemple_3_3_Results5}  shows the details of the optimization problems that KMFCaA solves at each iteration and it can be seen that the convexification of the concave constraint is much more parsimonious than with the other methods. Hence, the better performance of KMFCaA.\\
	\begin{minipage}[h]{\linewidth}
		\vspace*{0pt}
		\centering
		\includegraphics[width=4.45cm]{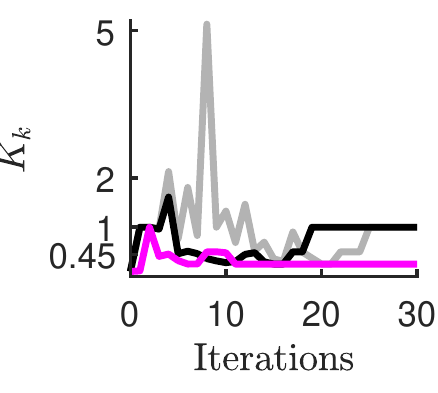}\hskip -0ex
		\includegraphics[width=4.45cm]{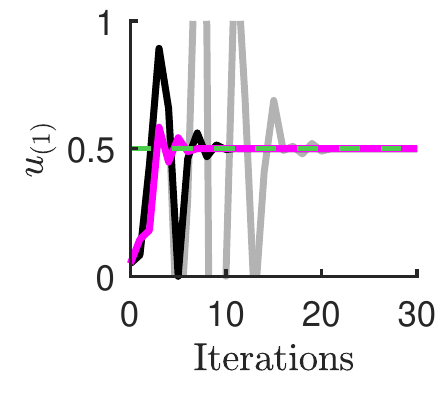}\hskip -0ex
		\includegraphics[width=4.45cm]{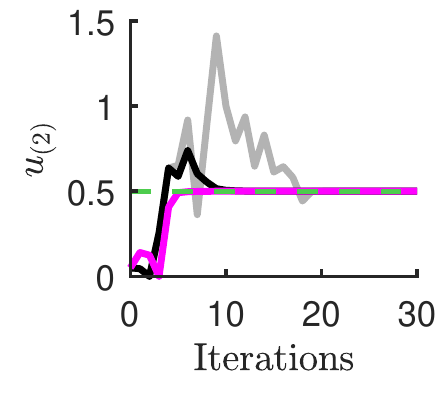} \\
		\includegraphics[width=4.45cm]{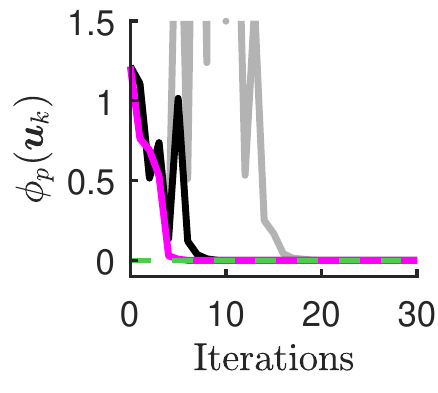}\hskip -0ex
		\includegraphics[width=4.45cm]{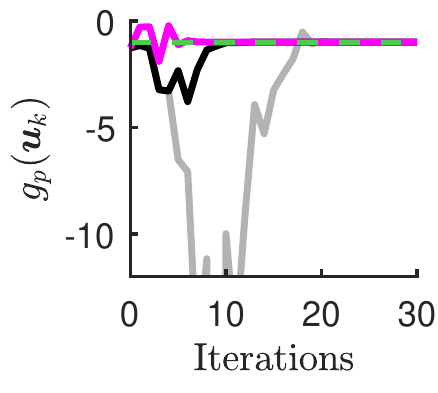}\hskip -0ex
		\includegraphics[width=4.45cm]{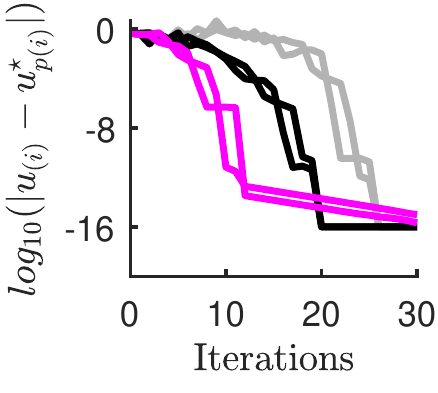}\\
		\textcolor{gray}{\raisebox{0.5mm}{\rule{0.5cm}{0.1cm}}}    : \text{MFCA,} \ 
		\textcolor{black}{\raisebox{0.5mm}{\rule{0.5cm}{0.1cm}}} : \text{KMFCA,} \ 
		\textcolor{magenta}{\raisebox{0.5mm}{\rule{0.5cm}{0.1cm}}}    : \text{KMFCaA}, \ 
		\textcolor{green1}{\raisebox{0.5mm}{\rule{0.2cm}{0.1cm}\hspace{0.1cm}\rule{0.2cm}{0.1cm}}} : Solution
		\vspace{-2mm}
		\captionof{figure}{Simulation results.}
		\label{fig:3___10_Exemple_3_3_Results1}
	\end{minipage} \\

	\begin{minipage}[h]{\linewidth}
		\vspace*{0pt}
		{\centering
			\includegraphics[trim={0.35cm 0.25cm 0.2cm  0.15cm },clip,width=2.75cm]{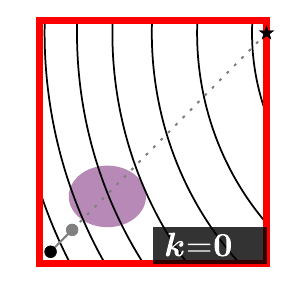}\hskip -0ex
			\includegraphics[trim={0.35cm 0.25cm 0.2cm  0.15cm },clip,width=2.75cm]{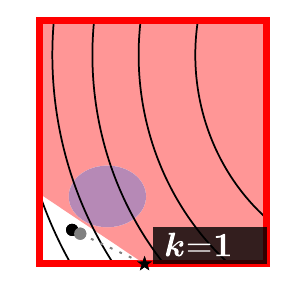}\hskip -0ex
			\includegraphics[trim={0.35cm 0.25cm 0.2cm  0.15cm },clip,width=2.75cm]{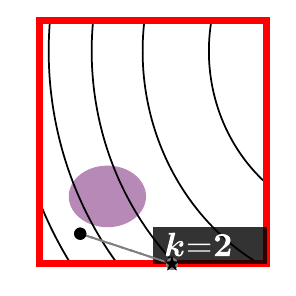}\hskip -0ex
			\includegraphics[trim={0.35cm 0.25cm 0.2cm  0.15cm },clip,width=2.75cm]{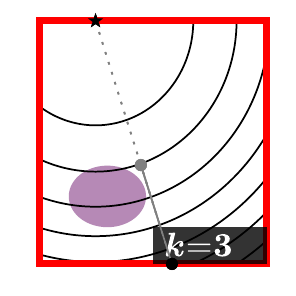}\hskip -0ex
			\includegraphics[trim={0.35cm 0.25cm 0.2cm  0.15cm },clip,width=2.75cm]{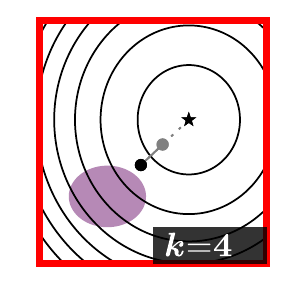} \\
			\includegraphics[trim={0.35cm 0.25cm 0.2cm  0.15cm },clip,width=2.75cm]{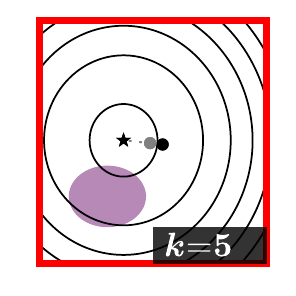}\hskip -0ex
			\includegraphics[trim={0.35cm 0.25cm 0.2cm  0.15cm },clip,width=2.75cm]{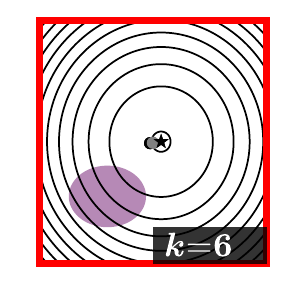}\hskip -0ex
			\includegraphics[trim={0.35cm 0.25cm 0.2cm  0.15cm },clip,width=2.75cm]{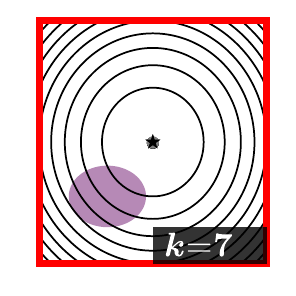}\hskip -0ex
			\includegraphics[trim={0.35cm 0.25cm 0.2cm  0.15cm },clip,width=2.75cm]{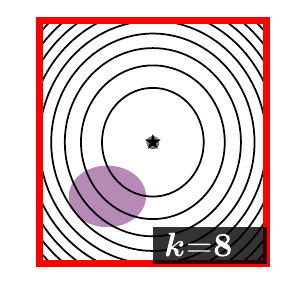}\hskip -0ex
			\includegraphics[trim={0.35cm 0.25cm 0.2cm  0.15cm },clip,width=2.75cm]{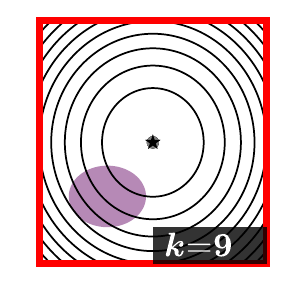} \\
		}$\bullet$:$\bm{u}_k$, \  
		\textcolor{gray}{$\bullet$}: $\bm{u}_{k+1}$, \ $\bm{\star}$:$\bm{u}^{\star}_{k+1}$, \
		\resizebox{0.5cm}{0.3cm}{\textcolor{red}{$\bm{\square}$}}: 
		Boundaries on $\bm{u}^{\star}_{k+1}$ and $\bm{u}_{k+1}$,
		\
		\textcolor{magenta_clair2}{\raisebox{-0.5mm}{\rule{0.5cm}{0.3cm}}}:
		Area forbidden to $(\bm{u}^{\star}_{k+1}, \bm{u}_{k+1})$		.
		\captionof{figure}{Details on the iterations of \textbf{KMFCaA}.}
		\label{fig:3___11_Exemple_3_3_Results5}
	\end{minipage} 
\end{exbox}

\section{The plant's constraints}

Before concluding this chapter, it is important to emphasize that KMFCA and KMFCaA reduce the risks of violating plant constraints by enforcing that the point applied at each iteration to be feasible according to the updated model. However, reducing these risks does not mean that these risks are absent. The following example illustrates this point.

\begin{exbox}\textbf{(Plant constraints can be violated)}
	Consider the following theoretical OTR problem:  
	\begin{align*}
		\phi_p(\bm{u}) := \ & 
		\frac{1}{2}\left(\bm{u} - \left(\begin{array}{c}1 \\ 0.3
		\end{array}\right)\right)^{\rm T} 
		\left(\begin{array}{cc}
			2 & 1 \\
			1 & 2
		\end{array}\right)
		\left(\bm{u} - \left(\begin{array}{c}1 \\ 0.3
		\end{array}\right)\right),  \\
		\phi(\bm{u}) := \ &  
		\frac{1}{2}\left(\bm{u} - \left(\begin{array}{c} 0.5 \\ 0.6
		\end{array}\right)\right)^{\rm T} 
		\left(\begin{array}{cc}
			2    & -1.1 \\
			-1.1 & 2
		\end{array}\right)
		\left(\bm{u} - \left(\begin{array}{c} 0.5 \\ 0.6
		\end{array}\right)\right),  \\
		g_p(\bm{u}) :=  \ &  0.3 
		-\frac{1}{2}\left(\bm{u} - \left(\begin{array}{c} 1 \\ 0
		\end{array}\right)\right)^{\rm T} 
		\left(\begin{array}{cc}
			1 & 0 \\
			0 & 1
		\end{array}\right)
		\left(\bm{u} - \left(\begin{array}{c} 1 \\ 0
		\end{array}\right)\right), \\
		g(\bm{u}) :=  \ &  
		0.75 -\frac{1}{2}\left(\bm{u} - \left(\begin{array}{c} 1.2 \\ 0.1
		\end{array}\right)\right)^{\rm T} 
		\left(\begin{array}{cc}
			1.7 & 0.51 \\
			0.51 & 1.7
		\end{array}\right)
		\left(\bm{u} - \left(\begin{array}{c} 1.2 \\ 0.1
		\end{array}\right)\right),
	\end{align*}
	where $u_{(1)}\in[0,1]$ and $u_{(2)}\in[0,1]$. Figure~\ref{fig:Exemple_3_4_maps} illustrates those functions.\\
	\begin{minipage}[h]{\linewidth}
		\vspace*{0pt}
		\centering 
		\includegraphics[width=5cm]{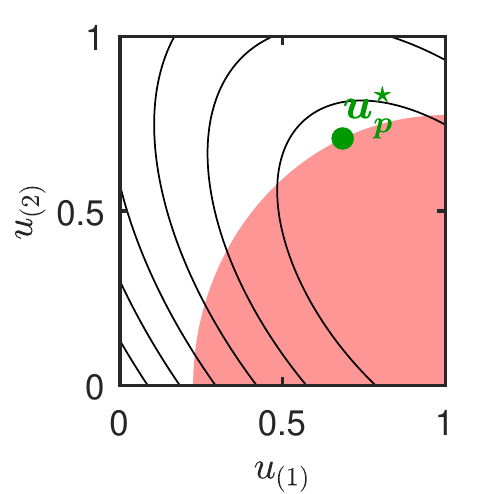}
		\includegraphics[width=5cm]{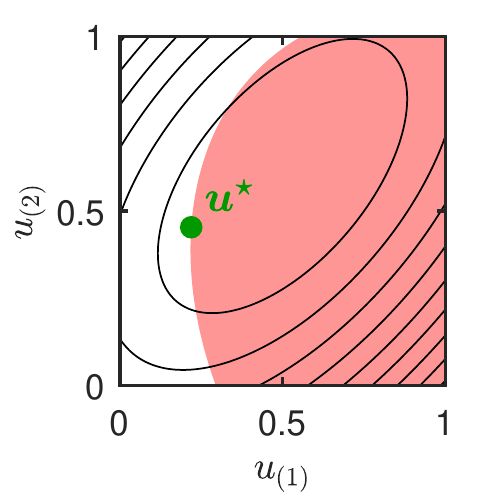}
		\vspace{-3mm}
		\captionof{figure}{Plant (left) and model's (right) cost and constraint functions}
		\label{fig:Exemple_3_4_maps}
	\end{minipage}\\

	MFCA, KMFCA and adKMA are initialized at $\bm{u}^{\star}$ and several several simulations are launched. The results are given on Figure~\ref{fig:3___13_Exemple_3_4_Results1}. \\
	
	It can observed that all the methods converge on the optimum of the plant and that on the path they take to reach this optimum:
	\begin{itemize} 
		\item MFCA violates the boundaries on $u_{(1)}$ and $u_{(2)}$. This kind of problem has already been discussed in Example~\ref{ex:3___2_Contournement_Contrainte}.
		\item KMFCA and KMFCaA violate the plant constraints at one iteration. The reason for this violation is simply that these two methods do not use the actual constraint functions of the plant but only those of the updated model. Figures~\ref{fig:3___14_Exemple_3_4_Results1} and \ref{fig:3___15_Exemple_3_4_Results2} show what are the updated-model-based optimization problems that KMFCA and KMFCaA use at each iteration. What can be observed on these figures is that $\bm{u}_k$ (the point in grey) is systematically feasible according to the updated model. What one can conclude is that the feasibility according to the updated model is not sufficient to ensure the feasibility of the plant. Which is quite logical.
	\end{itemize} 
	\begin{minipage}[h]{\linewidth}
		\vspace*{0pt}
		{\centering
		\includegraphics[width=4.45cm]{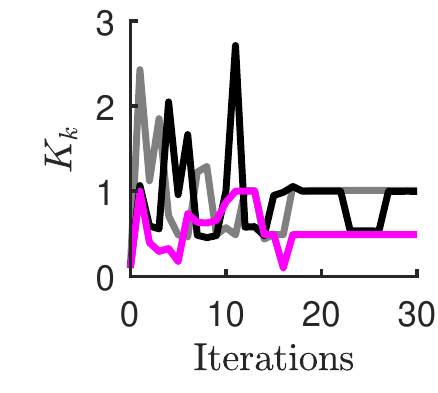}\hskip -0ex
		\includegraphics[width=4.45cm]{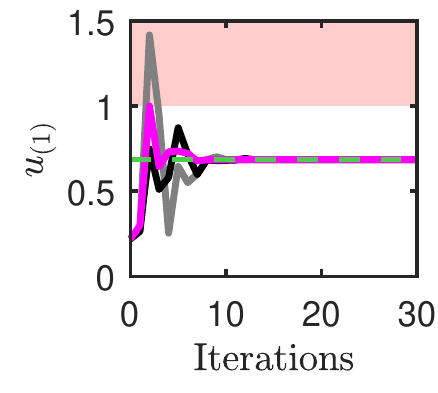}\hskip -0ex
		\includegraphics[width=4.45cm]{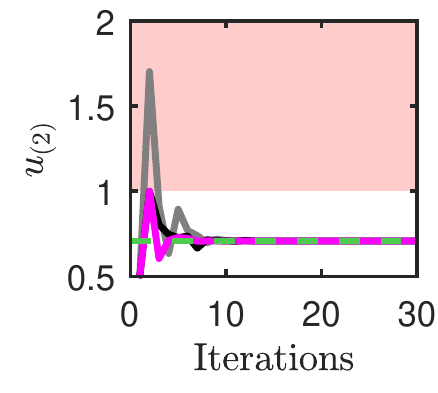} \\
		\includegraphics[width=4.45cm]{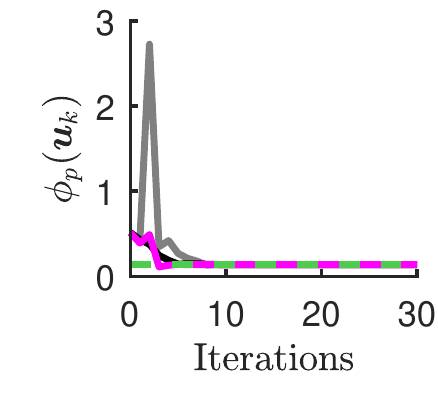}\hskip -0ex
		\includegraphics[width=4.45cm]{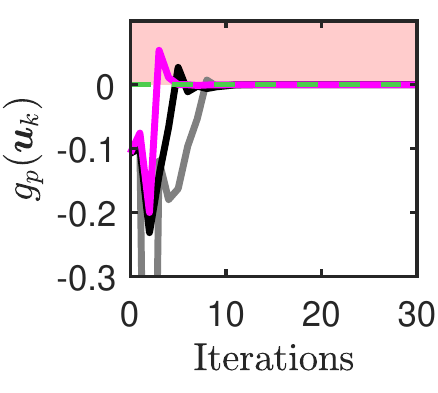}\hskip -0ex
		\includegraphics[width=4.45cm]{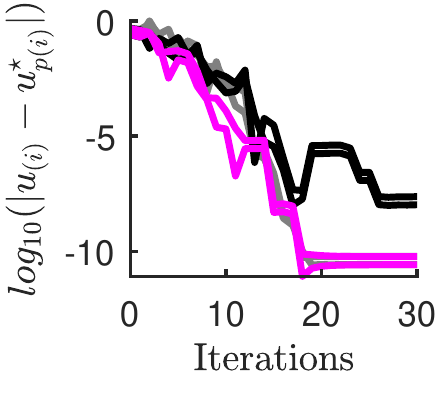} \\
	}
		\textcolor{gray}{\raisebox{0.5mm}{\rule{0.5cm}{0.1cm}}}    : \text{MFCA,} \ 
		\textcolor{black}{\raisebox{0.5mm}{\rule{0.5cm}{0.1cm}}} : \text{KMFCA,} \ 
		\textcolor{magenta}{\raisebox{0.5mm}{\rule{0.5cm}{0.1cm}}}    : \text{KMFCaA}, \ 
		\textcolor{green1}{\raisebox{0.5mm}{\rule{0.2cm}{0.1cm}\hspace{0.1cm}\rule{0.2cm}{0.1cm}}} : Solution
		\vspace{-2mm}
		\captionof{figure}{Simulation results.}
		\label{fig:3___13_Exemple_3_4_Results1}
	\end{minipage} \\
	
	\begin{minipage}[h]{\linewidth}
		\vspace*{0pt}
		{\centering
			\includegraphics[trim={0.35cm 0.25cm 0.2cm  0.15cm },clip,width=2.75cm]{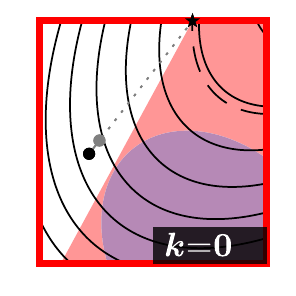}\hskip -0ex
			\includegraphics[trim={0.35cm 0.25cm 0.2cm  0.15cm },clip,width=2.75cm]{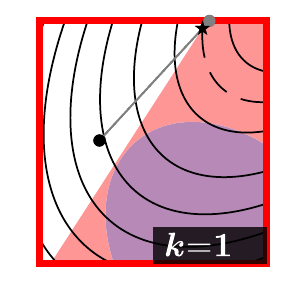}\hskip -0ex
			\includegraphics[trim={0.35cm 0.25cm 0.2cm  0.15cm },clip,width=2.75cm]{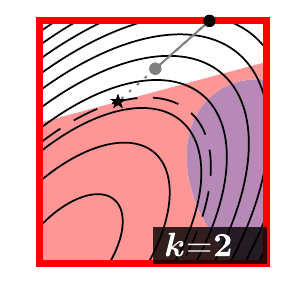}\hskip -0ex
			\includegraphics[trim={0.35cm 0.25cm 0.2cm  0.15cm },clip,width=2.75cm]{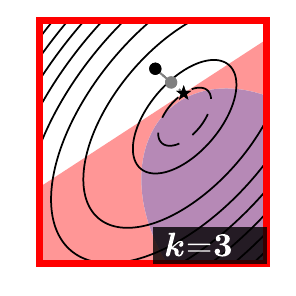}\hskip -0ex
			\includegraphics[trim={0.35cm 0.25cm 0.2cm  0.15cm },clip,width=2.75cm]{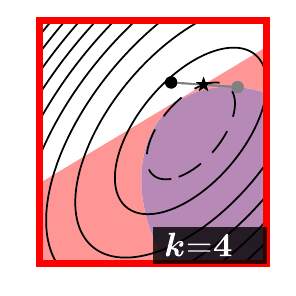} \\
			\includegraphics[trim={0.35cm 0.25cm 0.2cm  0.15cm },clip,width=2.75cm]{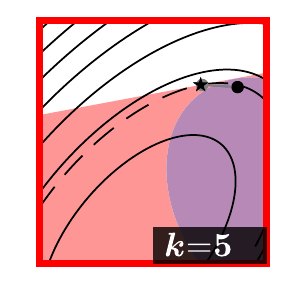}\hskip -0ex
			\includegraphics[trim={0.35cm 0.25cm 0.2cm  0.15cm },clip,width=2.75cm]{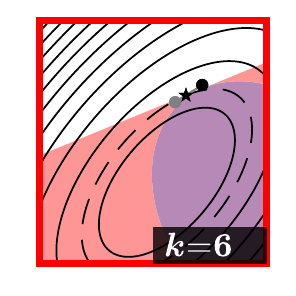}\hskip -0ex
			\includegraphics[trim={0.35cm 0.25cm 0.2cm  0.15cm },clip,width=2.75cm]{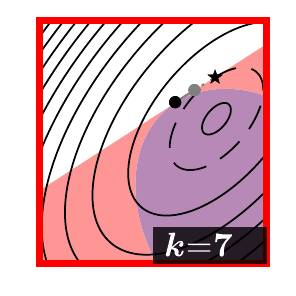}\hskip -0ex
			\includegraphics[trim={0.35cm 0.25cm 0.2cm  0.15cm },clip,width=2.75cm]{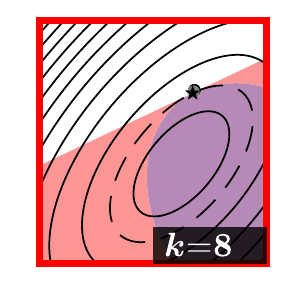}\hskip -0ex
			\includegraphics[trim={0.35cm 0.25cm 0.2cm  0.15cm },clip,width=2.75cm]{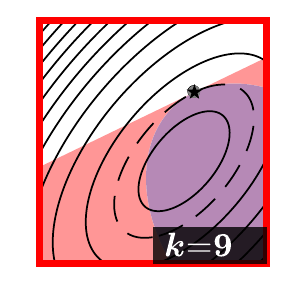} \\
		}$\bullet$:$\bm{u}_k$, \  
		\textcolor{gray}{$\bullet$}:$\bm{u}_{k+1}$, \ $\bm{\star}$:$\bm{u}^{\star}_{k+1}$, \
		\resizebox{0.5cm}{0.3cm}{\textcolor{red}{$\bm{\square}$}}:Boundaries on $\bm{u}$,
		\
		\textcolor{magenta_clair2}{\raisebox{-0.5mm}{\rule{0.5cm}{0.3cm}}}:Area forbidden to $(\bm{u}^{\star}_{k+1}, \bm{u}_{k+1})$
		\captionof{figure}{Details on the iterations of \textbf{KMFCA}.}
		\label{fig:3___14_Exemple_3_4_Results1}
	\end{minipage}

\begin{minipage}[h]{\linewidth}
	\vspace*{0pt}
	{\centering
		\includegraphics[trim={0.35cm 0.25cm 0.2cm  0.15cm },clip,width=2.75cm]{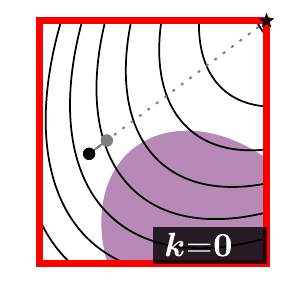}\hskip -0ex
		\includegraphics[trim={0.35cm 0.25cm 0.2cm  0.15cm },clip,width=2.75cm]{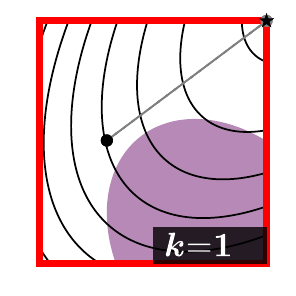}\hskip -0ex
		\includegraphics[trim={0.35cm 0.25cm 0.2cm  0.15cm },clip,width=2.75cm]{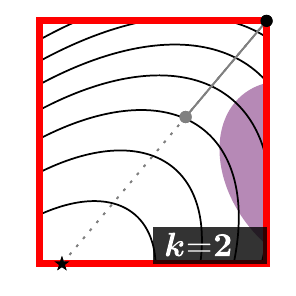}\hskip -0ex
		\includegraphics[trim={0.35cm 0.25cm 0.2cm  0.15cm },clip,width=2.75cm]{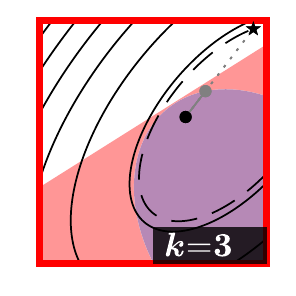}\hskip -0ex
		\includegraphics[trim={0.35cm 0.25cm 0.2cm  0.15cm },clip,width=2.75cm]{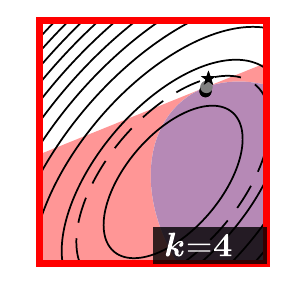} \\
		\includegraphics[trim={0.35cm 0.25cm 0.2cm  0.15cm },clip,width=2.75cm]{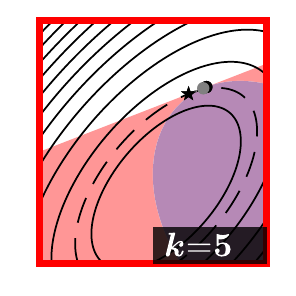}\hskip -0ex
		\includegraphics[trim={0.35cm 0.25cm 0.2cm  0.15cm },clip,width=2.75cm]{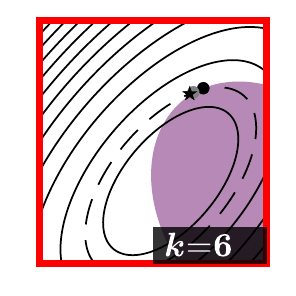}\hskip -0ex
		\includegraphics[trim={0.35cm 0.25cm 0.2cm  0.15cm },clip,width=2.75cm]{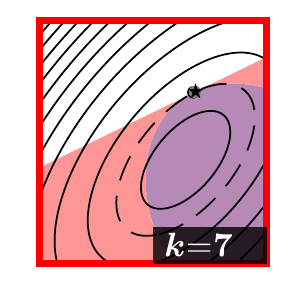}\hskip -0ex
		\includegraphics[trim={0.35cm 0.25cm 0.2cm  0.15cm },clip,width=2.75cm]{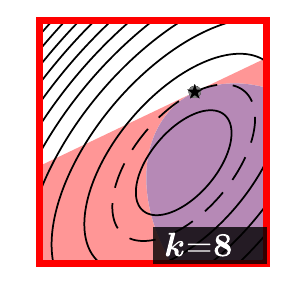}\hskip -0ex
		\includegraphics[trim={0.35cm 0.25cm 0.2cm  0.15cm },clip,width=2.75cm]{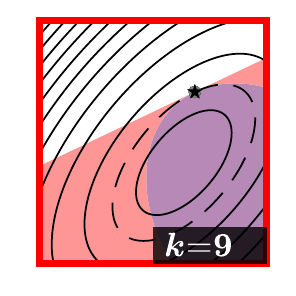} \\
	}$\bullet$:$\bm{u}_k$, \  
	\textcolor{gray}{$\bullet$}:$\bm{u}_{k+1}$, \ $\bm{\star}$:$\bm{u}^{\star}_{k+1}$, \
	\resizebox{0.5cm}{0.3cm}{\textcolor{red}{$\bm{\square}$}}:Boundaries on $\bm{u}$,
	\
	\textcolor{magenta_clair2}{\raisebox{-0.5mm}{\rule{0.5cm}{0.3cm}}}:Area forbidden to $(\bm{u}^{\star}_{k+1}, \bm{u}_{k+1})$
	\captionof{figure}{Details on the iterations of \textbf{KMFCaA}.}
	\label{fig:3___15_Exemple_3_4_Results2}
\end{minipage} 
\end{exbox}

\section{Conclusion}

The conceptual weakness of MFCA that is identified at the beginning of the chapter has been observed on several mathematical examples. From these example one can conclude that if MFCA provides robust guarantees on its ability to converge on the optimum of the plant, the path it takes can be very (too) risky in terms of constraint satisfaction. For example, it has been observed that even the boundaries on the plant's inputs can be violated, although they are a priori  easy to respect and their satisfaction is of primary importance. A simple improvement of MFCA which consists of (to some extent) merging the optimization and filtering steps of MFCA to enforce the points $\bm{u}_{k+1}$ and $\bm{u}_{k+1}^{\star}$ to be simultaneously feasible according to the updated model has been proposed. The resulting method is named KMFCA where the term K refers to the improvement on the use of the filter. KMFCA is shown to (i) have the same convergence properties as MFCA, (ii) provide much better guarantees of feasibility at each iteration $\bm{u}_k$ than MFCA, and (iii) have a minimal increase of the complexity of the algorithm. Then, it is shown that the performance of KMFCA can be slightly increased through a more parsimonious convexification of the concave constraints. Finally, through a last example, it has been reminded that, despite all the properties that MFCA and KMFCA offer, they do not guarantee the perfect and systematic satisfaction of the plant constraints. 
 
Through Chapters 2 and 3, the minimal and essential structure of a theoretical RTO method has been developed.
\textit{Minimal}, because it uses the minimum amount of information necessary to make 
convergence to the optimum of the plant \textit{possible} and \textit{stable}.
\textit{Essential}, because this information is used to make the least risky choices possible.  
In Chapters 5 and 6, it is discussed how to go further an integrate more information (than is strictly necessary) into the decision making process. However, before moving on to these more complex RTO methods, it is relevant to first build a simple autopilot based on the ``simple'' decision making procedure of KMFCaA.

	\chapter{A Simple Autopilot for Steady-State Processes}
		\label{Chap:4_S_ASP}	 

\section{Definition of an autopilot}

\begin{figure}[H]
	\includegraphics[width=14cm] {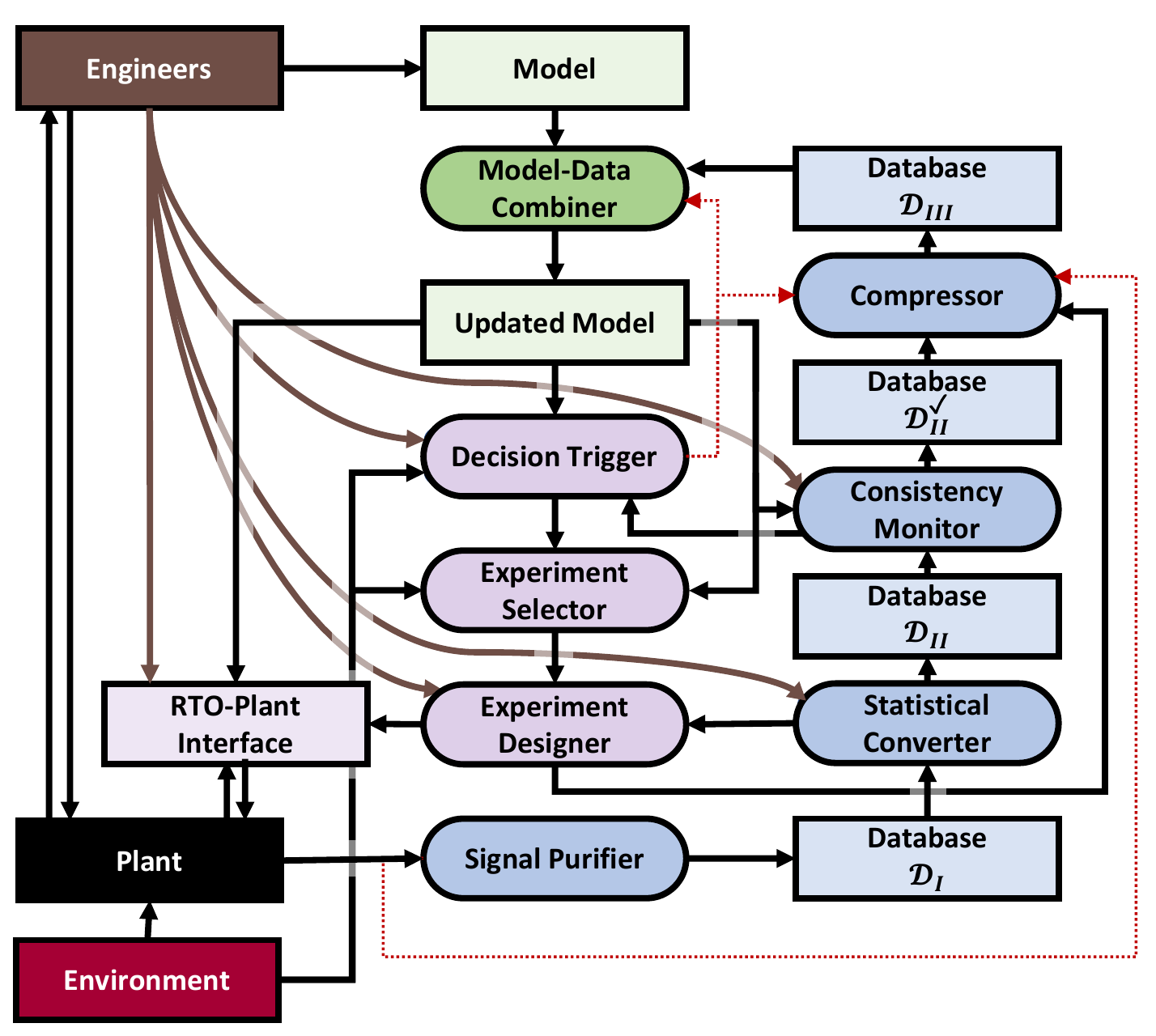}
	\begin{tikzpicture}
		\draw[-Triangle, very thick](0, 0) -- (1, 0);
	\end{tikzpicture}: Information flow,
	\raisebox{-1mm}{\resizebox{1cm}{0.4cm}{$\bm{\square}$}}: objects and  
	\tcbox[on line,boxsep=0pt,left=2pt,right=2pt,top=2pt,bottom=2pt,arc=5pt,colback      = white,colframe=black]{\phantom{aaa}}: functions related to \textcolor{bleu_ASP}{\raisebox{-0.5mm}{\rule{0.5cm}{0.3cm}}}:  data management,
	\textcolor{magenta_ASP}{\raisebox{-0.5mm}{\rule{0.5cm}{0.3cm}}}: decision making,
	\textcolor{vert_ASP}{\raisebox{-0.5mm}{\rule{0.5cm}{0.3cm}}}:  model-data fusion.
	\caption{General structure diagram of an ASP}
	\label{fig:4___1_StructureRTO_GENERALE}
\end{figure}
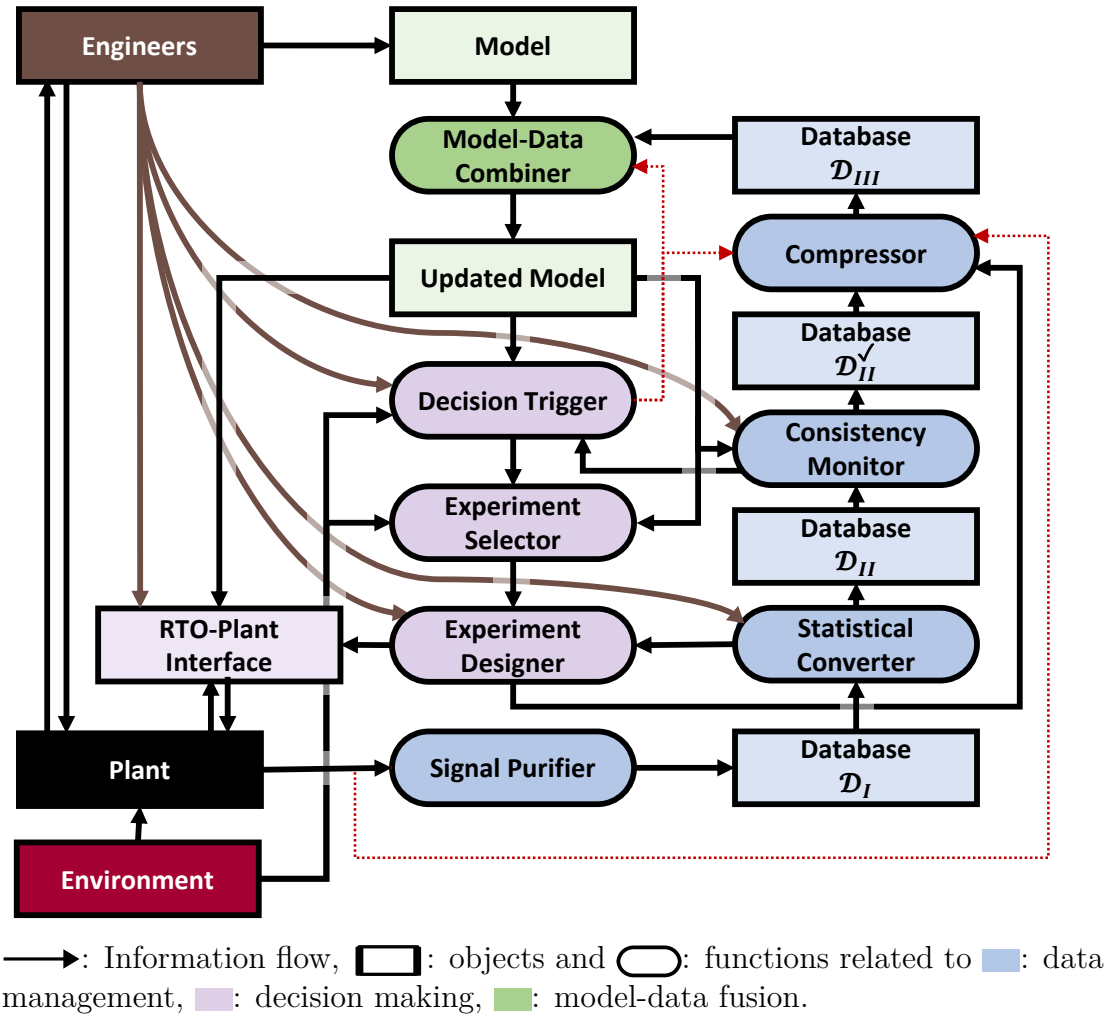

In this chapter, the \textit{the simplest possible} autopilot for stationary processes (ASP) is built, based on the theoretical RTO method introduced in the previous chapter (KMFCaA). To do this, dividing the work an ASP must perform into sub-functions and sub-objects as shown in Figure~\ref{fig:4___1_StructureRTO_GENERALE} is required. In this chapter, each of these sub-functions of an ASP has a section dedicated to it, giving a detailed presentation of its role as well as the simplest way to fulfill it.  Following these presentations, these sub-functions are combined to form what is called a \textit{simple autopilot for steady-state processes} (S-ASP). Then, it is used to optimize a simulation of the William-Otto reactor \cite{Williams:60}, which is a benchmark problem in the RTO community. 

\clearpage 

\section{Data management}

\subsection{The signal purifier}

\begin{figure}[H]
	\centering
	\includegraphics[width=14cm]{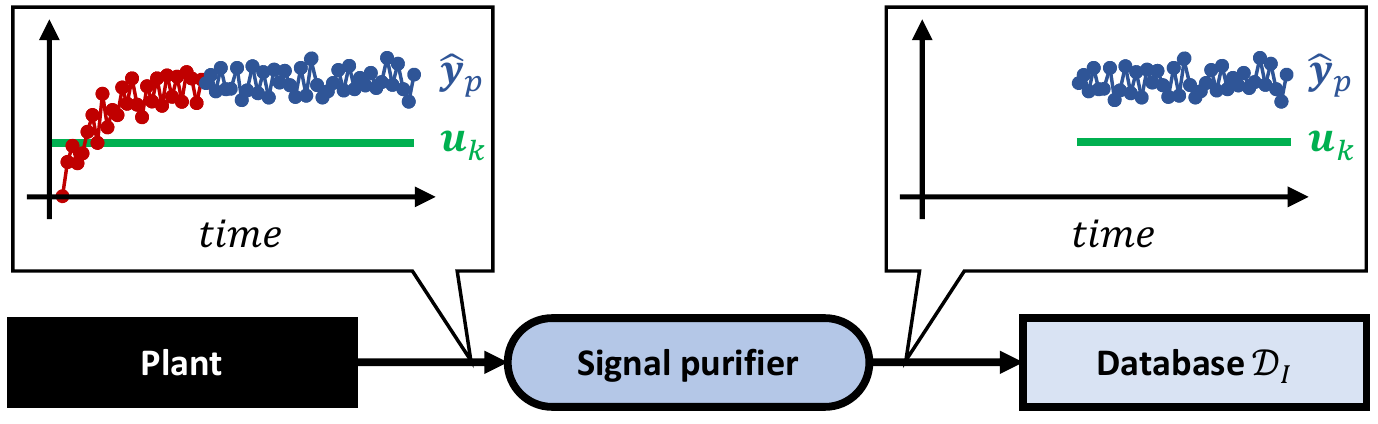} \\
	\textcolor{red!80!black}{$\bullet$}: $\widehat{\bm{y}}_p$ at transient state; \
	\textcolor{blue!70!black}{$\bullet$}: $\widehat{\bm{y}}_p$ at steady-state; \
	\textcolor{green!50!black}{\raisebox{0.5mm}{\rule{0.5cm}{0.1cm}}}: Plant's inputs, \
	\caption{Information flow entering and leaving the signal purifier}
	\label{fig:1_Purificateur de signal}
\end{figure}

As illustrated in Figure~\ref{fig:1_Purificateur de signal}, the role of the signal purifier is to ensure database $\mathcal{D}_{I}$ only receives the measures $\widehat{\bm{y}}_p$ associated with a steady state of the plant.  Clearly, this behavior is analogous to that of a classical  steady-state detectors (SS) to which one would add a switch function opening the data flow when the plant is at SS  and cutting it off the rest of the time.  Several techniques have been proposed to identify a SS, such as: F-like test \cite{Cao:95,Cao:97,Shrowti:10}, Wavelet transform \cite{Jiang:2000,Jiang:03,Korbel:2014}, Polynomial-based \cite{LeRoux:2008}, and ARX-based \cite{Rincon:2015}).  However, the details of their functioning are not given here because SS detectors were not studied as part of this body of work. Instead, and in order to simplify the following developments, it is always assumed that:
\begin{Assumption}
	One considers that the SS detector is ideal and that it always perfectly identifies when the plant is at SS.   
\end{Assumption}

\subsection{The statistical converter}

\begin{figure}[H]
	\centering
	\includegraphics[width=14cm]{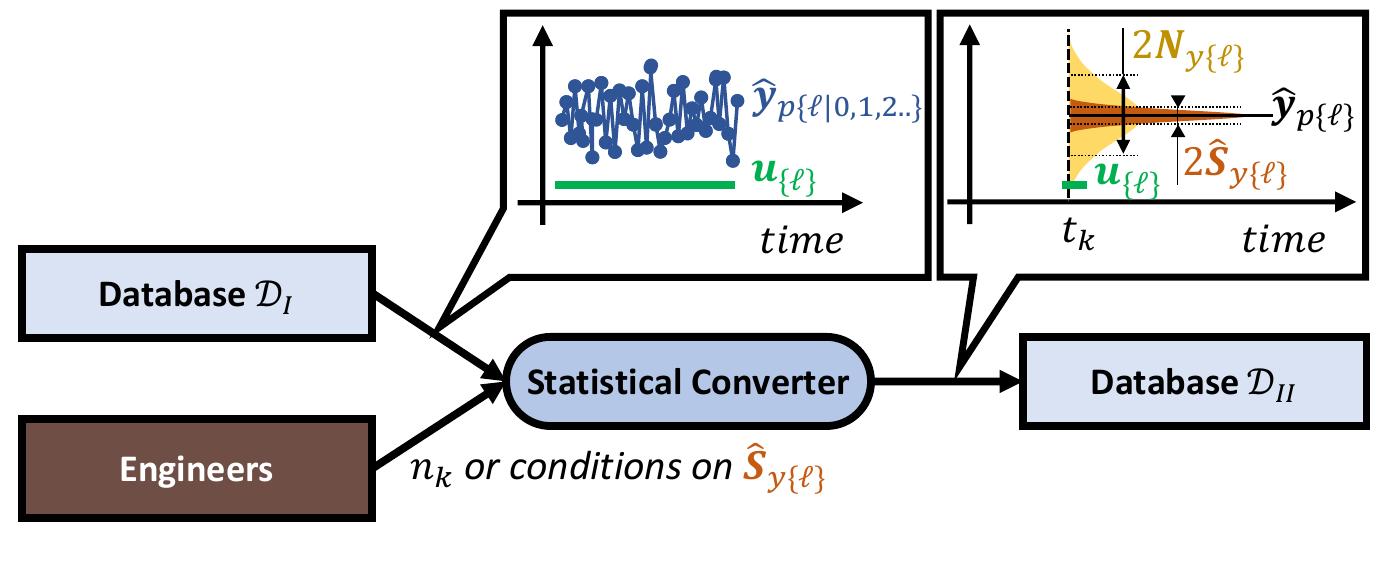}
	\caption{Information flow entering and leaving the statistical converter}
	\label{fig:2_Convertisseur_Statistique}
\end{figure}

As illustrated in Figure~\ref{fig:2_Convertisseur_Statistique}, the role of the statistical converter is to transform the raw measurement sequences of $\mathcal{D}_{I}$ into \textit{usable}  statistical values. Indeed, from the moment the SS is detected, the raw measurements $\widehat{\bm{y}}_p$ of the plant's SS are accumulated to form a sequence of raw data:   $\widehat{\bm{y}}_{p\{\ell|i\}}$, where the index ``$p\{\ell|i\}$'' refers to the $i$-th measurement of plant $p$ during the experiment $\ell$. It has been assumed that all the measures, $\widehat{\bm{y}}_p$, are corrupted by measurement errors whose statistical distribution is normal and unbiased  (see \eqref{eq:1___3_Mesures_Normal_Uncertainty}). So, these measures follow:
\begin{equation}
	\{\widehat{\bm{y}}_{p\{\ell|0\}}, ..., \widehat{\bm{y}}_{p\{\ell|n_{\ell}\}} \} \sim \mathcal{N}(\bm{y}_{p\{\ell\}}, \bm{N}_{y\{\ell\}}),
\end{equation}
where $\bm{y}_{p\{\ell\}}$ is the actual (de-noised) plant outputs, and $\bm{N}_{y\{\ell\}}$ is the covariance matrix characterizing the noise affecting the measurements of  $\bm{y}_{p\{\ell\}}$.
One decides to use a Bayesian approach to estimate $\bm{y}_{p\{\ell\}}$ from these $n_{\ell}+1$ measures. Identifying the mean value of a signal can be interpreted as an affine regression problem where the linear term is fixed to $0$. Therefore, the development presented in Chapter 1 of \cite{Rasmussen:2006} is used to perform this regression. The Bayesian regression problem is then the following: 
\begin{equation}
	\begin{array}{ll}
		\text{Data matrix:}         &  \widehat{\bm{Y}}_{p\{\ell\}}:= [\widehat{\bm{y}}_{p\{\ell|1\}}, ...,\widehat{\bm{y}}_{p\{\ell|n_{\ell}\}}]^{\rm T}, \\
		\text{Design matrix:}       &  \bm{X} := \bm{1}_{1\times n_{\ell}}, \\
		\text{Gaussian prior:}      & p(\bm{y}_{p\{\ell\}}) \sim \mathcal{N}(\widehat{\bm{y}}_{p\{\ell|0\}}, \bm{N}_{y\{\ell\}}),  \\
		\text{Gaussian likelyhood:} & p(\widehat{\bm{y}}_{p\{\ell|i\}}|\bm{y}_{p\{\ell\}}) \sim \mathcal{N}(\bm{y}_{p\{\ell\}},  \bm{N}_{y\{\ell\}}), \ \forall i = 0,..,n_{\ell},
	\end{array}
\end{equation}
where $\bm{1}_{1\times n_k}$ is a matrix of ones of size  $1\times n_{\ell}$.

The initial belief (the Gaussian prior) on $\bm{y}_{p\{\ell\}}$ is chosen as a normal distribution centered on the first measurement $\widehat{\bm{y}}_{p\{\ell|0\}}$ and whose covariance matrix is  $\bm{N}_{y\{\ell\}}$ (if it is unknown, then see the section~\ref{rem:3_1_Si_Npk_Inconnu}). Then, the idea is to update this initial belief with the $n_{\ell}$ measure $\{\widehat{\bm{y}}_{p\{\ell|1\}},...,\widehat{\bm{y}}_{p\{\ell|n_{\ell}\}}\}$. This updated belief is what is commonly called the Gaussian posterior, which is: 
\begin{equation}
	p(\bm{y}_{p\{\ell\}} | \bm{X},\widehat{\bm{Y}}_{p\{\ell\}}) \sim \mathcal{N}(\widehat{\bm{y}}_{p\{\ell\}}, \widehat{\bm{S}}_{p\{\ell\}}),
\end{equation}
where
\begin{equation} \label{eq:4___4_def_yS}
	\begin{array}{rll}
		\widehat{\bm{S}}_{p\{\ell\}} :=  & (\bm{N}_{y\{\ell\}}^{-1} \bm{XX}^{\rm T} + \bm{N}_{y\{\ell\}}^{-1})^{-1} & =  \dfrac{\bm{N}_{y\{\ell\}}}{n_{\ell}+1},	\\
		\widehat{\bm{y}}_{p\{\ell\}} :=  & \bm{N}_{p\{\ell\}}^{-1} \widehat{\bm{S}}_{y\{\ell\}}(\bm{X}\widehat{\bm{Y}}_{p\{\ell\}}+\widehat{\bm{y}}_{p\{\ell|0\}}^{\rm T})^{\rm T} & = \dfrac{1}{n_{\ell}+1} \displaystyle\sum_{a=0}^{n_{\ell}} \widehat{\bm{y}}_{p\{\ell|a\}}.
	\end{array}
\end{equation}
where $\widehat{\bm{y}}_{p\{\ell\}}$ is the estimate $\bm{y}_{p\{\ell\}}$ and $\widehat{\bm{S}}_{y\{\ell\}}$ is the uncertainty associated to this estimate. 

The greater $n_{\ell}$, the closer matrix  $\widehat{\bm{S}}_{y\{\ell\}}$ gets to $\bm{0}$. It is therefore possible to control the uncertainty of the estimates $\widehat{\bm{y}}_{p\{\ell\}}$ by choosing a  $n_{\ell}$ large  enough. This choice is made either directly by the experiment  designer  (presented in section~\ref{sec:4___Designer_d_experience}), or by the statistical converter on the basis of conditions on $\widehat{\bm{S}}_{y\{\ell\}}$, e.g. an upper bound on its eigenvalues. \\ 

It has just been shown that the quality of the estimate $\widehat{\bm{y}}_{p\{\ell\}}$ is a quantity that can be controlled. This is an important point to emphasize. Applying RTO to a system often goes with the temptation of reacting as soon as SS is detected to immediately look for the next set of operating conditions, just like a controller is supposed to react immediately to a control error. Doing so means making the choice, neither necessarily nor trivially, to set $n_{\ell}=0$ and not to control the quality of this the estimate $\widehat{\bm{y}}_{p\{\ell\}}$, which is key to the quest for the optimal plant inputs. On the other hand, setting $n_{\ell}> 0$ allows, to, e.g. mitigate the impact of measurement noise to the estimate $\widehat{\bm{y}}_{p\{\ell\}}$.\\
%

Finally, any group of measures $\{\widehat{\bm{y}}_{p\{\ell\}},$ $\widehat{\bm{S}}_{p\{\ell\}},\bm{N}_{p\{k\}}\}$ is always associated  with a date $t_{\{\ell\}}$ which corresponds to the average of the dates $\{t_{\{\ell|0\}}, ..., t_{\{\ell|n_{\ell}}\}\}$ at which the measures $\{\widehat{\bm{y}}_{p\{\ell|0\}}, ..., \widehat{\bm{y}}_{p\{\ell|n_{\ell}\}}\}$ have been made, i.e. 
\begin{equation}
	t_{\{\ell\}} := \frac{1}{n_{\ell}+1}\sum_{a=0}^{n_{\ell}}t_{\{\ell|a\}}.
\end{equation}

\subsubsection{If the covariance of the measurement noise is unknown}
 \label{rem:3_1_Si_Npk_Inconnu} 
	In the previous analysis it was assumed that the covariance matrix of the noise $\bm{N}_{y\{\ell\}}$ is known.  If this is not the case, then it can be deduced from the data $\widehat{\bm{Y}}_{p\{\ell\}}$. Indeed, identifying $\bm{N}_{y\{\ell\}}$ from data is similar to the model selection problem presented in Chapter 5 of \cite{Rasmussen:2006}. However, since a different (simpler) regression method is used here, the approach must be slightly modified in its application. The idea is to choose the matrix  $\bm{N}_{p\{\ell\}}$  that maximizes the likelihood of the observations, i.e.
	\begin{equation} \label{eq:4___5_Pb_Opt_Id_Noise}
		\bm{N}_{y\{\ell\}} := \operatorname{arg}
		\underset{\bm{N}}{\operatorname{max}} \sum_{a=0}^{n_{\ell}} p(\widehat{\bm{y}}_{p\{\ell|a\}}|\bm{X},\bm{Y}_{p\{\ell\}},\bm{N}) \ \text{s.t.} \  \bm{N}=\bm{N}^{\rm T}, \ \bm{N}\bm{v}>0, \forall \bm{v}\in\amsmathbb{R}^{n_y},
	\end{equation} 
	where the constraints impose that $\bm{N}$ must be symmetrical and positive definite (an uncertainty cannot be negative).  To solve this problem, it is proposed to simplify it by considering only the one-dimensional case ($n_y=1$ so $\bm{N}\in\amsmathbb{R}^{n_y\times n_y}$ is replaced by $\nu^2\in\amsmathbb{R}$), and one considers the 1st order optimality condition:	
	\begin{equation} \label{eq:4___5_Compute_nu}
	\frac{\partial}{\partial \nu^2} \sum_{a=0}^{n_{\ell}} p(\widehat{y}_{p\{\ell|a\}}|\bm{X},\bm{Y}_{p\{\ell\}},\nu^2)\big|_{\nu^2 = \nu_{y\{\ell\}}^2} = 0.
	\end{equation}
	
	\begin{align*}
		p(\widehat{y}_{p\{\ell|a\}}|\bm{X},\bm{Y}_{p\{\ell\}},\nu^2) := \ & 
		\frac{1}{\nu^2\sqrt{2\pi}} \exp\left(
		-\frac{1}{2}\left(
		\frac{\widehat{y}_{p\{\ell|a\}} - \widehat{y}_{p\{\ell\}}}{\nu^2}
		\right)^2
		\right), \\
		\frac{\partial}{\partial \nu^2}	p(\widehat{y}_{p\{\ell|a\}}|\bm{X},\bm{Y}_{p\{\ell\}},\nu^2) := \ &
		\frac{1}{\nu^2\sqrt{2\pi}}
		\left(
		\frac{(\widehat{y}_{p\{\ell|a\}} - \widehat{y}_{p\{\ell\}})^2 }{\nu^2}  -
		1
		\right) 
		\exp\left(
		\frac{\text{-}1}{2}\left(
		\frac{\widehat{y}_{p\{\ell|a\}} - \widehat{y}_{p\{\ell\}}}{\nu^2}
		\right)^2
		\right),
	\end{align*}
	Therefore, the values of  $\nu_{y\{\ell\}}^2$ for which \eqref{eq:4___5_Compute_nu} is true are the solutions of: 
	$$
		\frac{1}{\nu_{y\{\ell\}}^2\sqrt{2\pi}}
		\left(
		\frac {1}{\nu_{y\{\ell\}}2}\sum_{a=0}^{n_{\ell}} \left((\widehat{y}_{p\{\ell|a\}} - \widehat{y}_{p\{\ell\}})^2\right) - (n_{\ell}+1)
		\right) = 0,
	$$
	which are (solution 1:) $\nu_{y\{\ell\}}^2\rightarrow \infty$ and (solution 2:)
	\begin{align}
		& & \nu_{y\{\ell\}}^2 = \ & \dfrac{1}{n_{\ell}+1}\sum_{a=0}^{n_{\ell}} \left((\widehat{y}_{p\{\ell|a\}} - \widehat{y}_{p\{\ell\}})^2\right), \nonumber \\
		\big(\text{\eqref{eq:4___4_def_yS}} \Rightarrow\big) & &      = \ & \dfrac{1}{n_{\ell}+1}\sum_{a=0}^{n_{\ell}} \left(\left(\widehat{y}_{p\{\ell|a\}} - \dfrac{1}{n_{\ell}+1} \displaystyle\sum_{a=0}^{n_{\ell}} \widehat{y}_{p\{\ell|a\}}\right)^2\right). \label{eq:4___7_yufdhs}
	\end{align}
	Clearly, solution 1 is a minimum, and solution 2 is the desired  maximum. One can notice that \eqref{eq:4___7_yufdhs} is the definition of the experimental variance. Therefore, one proposes to generalize this result to multidimensional cases. So the matrix $\bm{N}_{y\{\ell\}}$  to use should be the experimental estimate of the covariance of the points $\{\widehat{\bm{y}}_{p\{\ell|0\}}, ..., \widehat{\bm{y}}_{p,\{\ell|n_{\ell}\}} \}$.

\subsection{The consistency monitor}

\begin{figure}
	\centering
	\includegraphics[width=14cm]{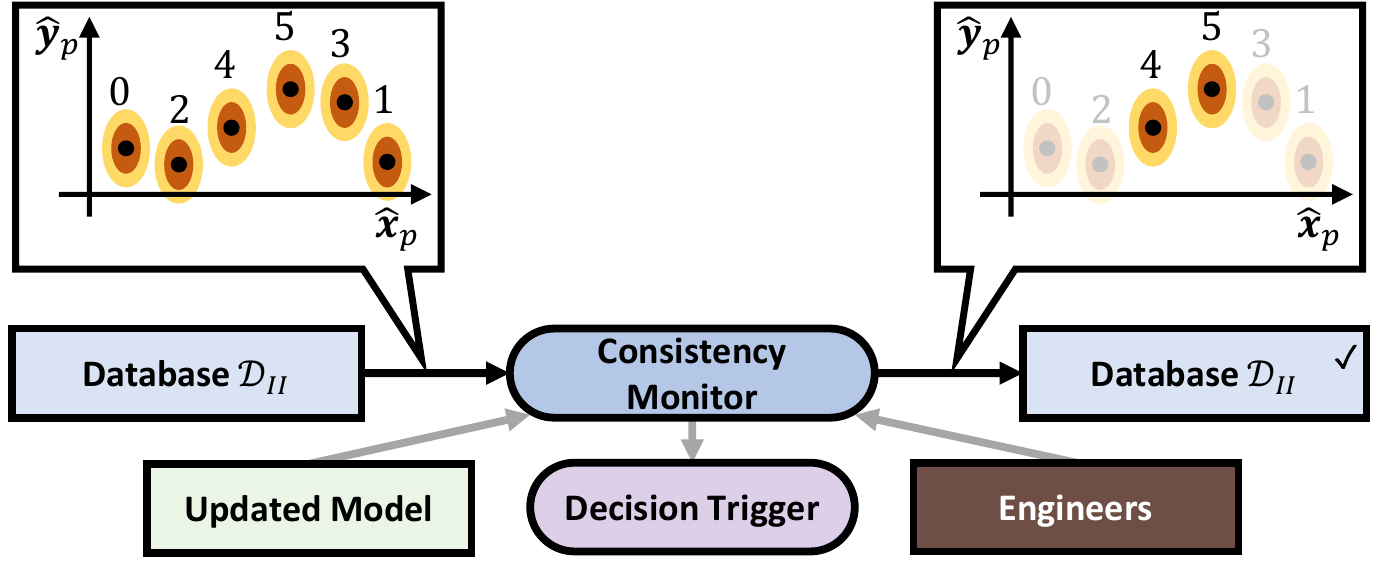}
	\caption{Flow of information entering and leaving the consistency monitor at iteration 5 (hence the numbering of the points from 0 to 5.)}
	\label{fig:3_Controleur_Coherence}
\end{figure}

As illustrated in Figure~\ref{fig:3_Controleur_Coherence}, the role of the consistency monitor is to remove outdated data from $\mathcal{D}_{II}$ to prevent them from being used for the decision making. Indeed, this database at the iteration $k$ contains the entire history of the observations, i.e.
\begin{equation}
	\mathcal{D}_{II\{\ell\}} := \{X_{\{\ell\}}, Y_{\{\ell\}}, S_{\{\ell\}}, N_{\{\ell\}}, T_{\{\ell\}} \},
\end{equation}
where 
\begin{equation*}
	\begin{array}{r@{}c@{}l@{}c@{}l@{}c@{}r}
		X_{\{\ell\}} & := & \{ & \widehat{\bm{x}}_{\{0\}}     & ,... \ , & \widehat{\bm{x}}_{\{\ell\}}      & \},\\
		N_{\{\ell\}} & := & \{ & \bm{N}_{p\{0\}} & ,... \ , & \bm{N}_{p\{\ell\}} & \}, 
	\end{array}
	\quad
	\begin{array}{r@{}c@{}l@{}c@{}l@{}c@{}r}
		Y_{\{\ell\}} & := & \{ & \widehat{\bm{y}}_{p\{0\}} & ,... \ , & \widehat{\bm{y}}_{p\{\ell\}} & \}, \\
		T_{\{\ell\}} & := & \{ & t_{\{0\}} & ,... \ , & t_{\{\ell\}} & \}. 
	\end{array}
	\quad 
	\begin{array}{r@{}c@{}l@{}c@{}l@{}c@{}r}
		S_{\{\ell\}} & := & \{ & \widehat{\bm{S}}_{p\{0\}} & ,... \ , & \widehat{\bm{S}}_{p\{\ell\}} & \},  \\
		\ &
	\end{array}
\end{equation*}
However, it is possible that, over the course of time and/or because of events such as the aging of the installation or the degradation of its components, some of these results are no longer up to date. The results of an experiment $i$ are considered outdated if, for $\widehat{\bm{x}}_{\{\ell+1\}}= \widehat{\bm{x}}_{\{i\}}$, the observations $(\widehat{\bm{y}}_{p\{\ell+1\}},\widehat{\bm{N}}_{y\{\ell+1\}})$ and $(\widehat{\bm{y}}_{p\{i\}},\widehat{\bm{N}}_{y\{i\}})$ are significantly different from one another. The term ``significantly'' is deliberately left vague since a satisfactory definition has never been proposed (at least to the author's knowledge). A definition will be proposed in Chapter~\ref{Chap:6_ASP} (where a more complex autopilot is built). 

The simplest way to perform this task is to consider all the observations that were performed before the experiment  $\ell-n_{\bm{u}}$ as outdated.  This choice is not meaningless as it reduces the amount of data in  $\mathcal{D}_{II\{\ell\}}$ to the minimal number necessary to compute the gradients of the plant. Indeed, to compute a gradient in a $n_{\bm{u}}$-dimensional  space, at least $n_{\bm{u}}+1$ measurements are required. The \textit{validated} database $\mathcal{D}_{II\{\ell\}}^{\bm{\checkmark}} $ is then defined as:
\begin{equation}
	\mathcal{D}_{II\{\ell\}}^{\bm{\checkmark}} := \{X_{\{\ell\}}^{\bm{\checkmark}}, Y_{\{\ell\}}^{\bm{\checkmark}}, S_{\{\ell\}}^{\bm{\checkmark}}, N_{\{\ell\}}^{\bm{\checkmark}}, T_{\{\ell\}}^{\bm{\checkmark}}\},
\end{equation}
where 
\begin{equation*}
	\begin{array}{r@{}c@{}l@{}c@{}l@{}c@{}r}
		X_{\{\ell\}}^{\bm{\checkmark}} & := & \{ & \widehat{\bm{x}}_{\{\ell-n_{\bm{u}}\}}     & ,... \ , & \widehat{\bm{x}}_{\{\ell\}}     & \},\\
		S_{\{\ell\}}^{\bm{\checkmark}} & := & \{ & \widehat{\bm{S}}_{p\{\ell-n_{\bm{u}}\}} & ,... \ , & \widehat{\bm{S}}_{p\{\ell\}} & \},  \\
		T_{\{\ell\}}^{\bm{\checkmark}} & := & \{ & t_{\{\ell-n_{\bm{u}}\}} & ,... \ , & t_{\{\ell\}} & \}. 
	\end{array}
	\qquad
	\begin{array}{r@{}c@{}l@{}c@{}l@{}c@{}r}
		Y_{\{\ell\}}^{\bm{\checkmark}} & := & \{ & \widehat{\bm{y}}_{p\{\ell-n_{\bm{u}}\}} & ,... \ , & \widehat{\bm{y}}_{p\{\ell\}} & \}, \\
		N_{\{\ell\}}^{\bm{\checkmark}} & := & \{ & \bm{N}_{p\{\ell-n_{\bm{u}}\}} & ,... \ , & \bm{N}_{p\{\ell\}} & \}, \\
		\
	\end{array}
\end{equation*}

Of course, this simple approach will not prevent the occasional use of outdated data.  Indeed, if an event changes the behavior of the plant between the experiments $\ell-1$ and  $\ell$, then it will be necessary to carry out $n_u$ additional  experiments so that the measures of before and after this event do not coexist in $\mathcal{D}_{II\{\ell\}}^{\bm{\checkmark}}$.  So, this simple strategy is:
\begin{itemize}
	\item able to reject the effects of changes in the plant's behavior thanks to its limited memory, and
	\item sensitive to changes in the plant's behavior at the time of their occurrence until at least $n_u$ additional experiments are performed.
\end{itemize}

Finally, in Figure~\ref{fig:3_Controleur_Coherence} the consistency monitor is also linked:
\begin{itemize}
	\item to the updated model because it would be possible to use it to identify outdated data as will be proposed in Chapter~\ref{Chap:6_ASP}.  However, as simplicity is sought, this possibility is not considered here.   
	\item to the decision trigger to, for example, notify it that a change in the plant's behavior has been observed or that a measured disturbance has changed significantly, and that therefore a decision must be made in response to this event. However, since the simple version of the consistency monitor does not identify changes in the plant's behavior, this information flow is not considered. 
	\item to the engineers so that they  have the option to manually add or remove data. However, in this simple version of the coherency monitor, this possibility is not included and such an information flow is not considered. 
\end{itemize}

\subsection{The compressor (standard operation)}

\begin{figure}
	\centering
	\includegraphics[width=14cm]{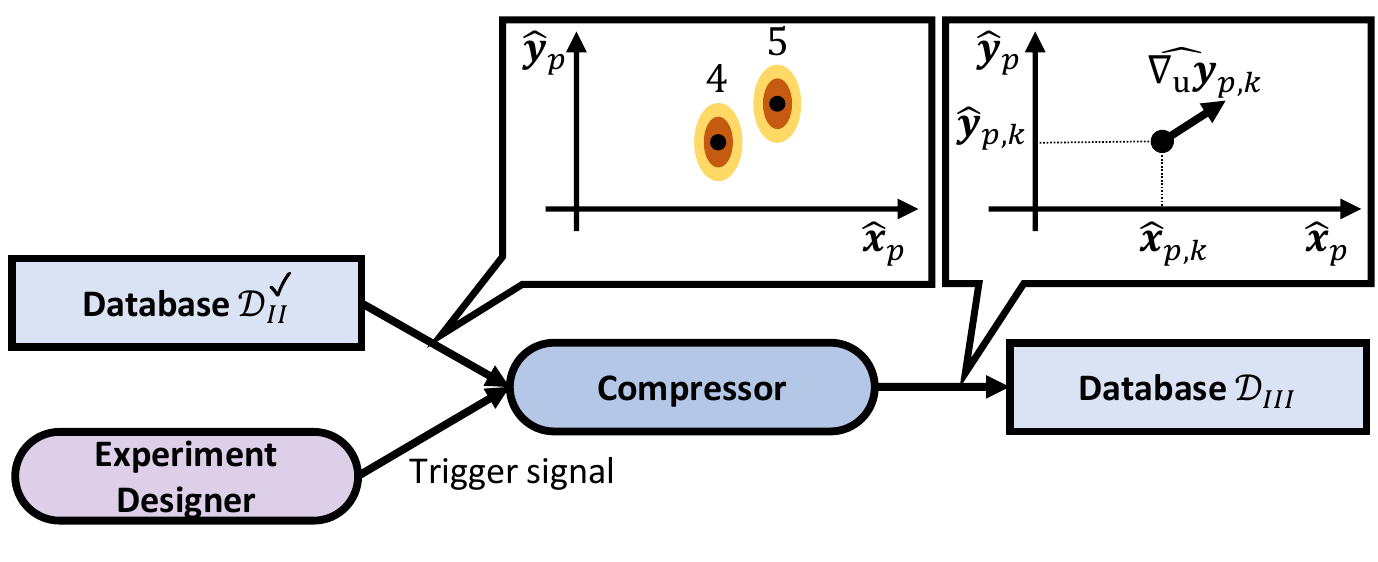}
	\caption{Flow of information entering and leaving the compressor}
	\label{fig:4_Compresseur}
\end{figure}

As illustrated on Figure~\ref{fig:4_Compresseur}, the role of the compressor is to convert the results of the database $\mathcal{D}_{II\{\ell\}}^{\checkmark}$ into a smaller database with a minimal loss of information. 

A simple way to accomplish this task is to replace the similar experiments, i.e. the experiments whose  $\widehat{\bm{x}}_{p\{i\}}$ with $i\in\{0,1,...,\ell\}$ are close, by a value and a gradient. This general idea will be defined and fully explained in Chapter~\ref{Chap:6_ASP}. In the framework of S-ASP, this strategy can be summarized by the use of the finite difference method each time the experiment designer gives the instruction to do so. When this instruction is given, the values and gradients of the plant are computed as follows:
\begin{align}
	\widehat{\bm{x}}_{p,k} := \ & \widehat{\bm{x}}_{p\{\ell-n_u\}}, \nonumber \\
	\widehat{\bm{y}}_{p,k} := \ & \widehat{\bm{y}}_{p\{\ell-n_u\}}, \label{eq:4___3_11_y_bar}\\
	\widehat{\nabla_{\bm{u}} \bm{y}}_{p,k} := \ & 
	\left(
	\begin{array}{ccc}
		\dfrac{\widehat{\bm{y}}_{p\{\ell-n_u+1\}}-\widehat{\bm{y}}_{p\{\ell-n_u\}}}{\|\widehat{\bm{x}}_{p\{\ell-n_u+1\}}-\widehat{\bm{x}}_{p\{\ell-n_u\}}\|},   & \hdots , &   \dfrac{\widehat{\bm{y}}_{p\{\ell\}}-\widehat{\bm{y}}_{p\{\ell-n_u\}}}{\|\widehat{\bm{x}}_{p\{\ell\}}-\widehat{\bm{x}}_{p\{\ell-n_u\}}\|}
	\end{array}
	\right). \label{eq:4___3_12_nabla_y_bar}
\end{align}
The uncertainties on $(\widehat{\bm{x}}_{p,k}, \widehat{\bm{y}}_{p,k}, \widehat{\nabla_{\bm{u}} \bm{y}}_{p,k})$ can be estimated but as they are not going to be used afterwards, these calculations are not made. It is clear that the name compressor for this function that transforms $n_u+1$ points  $\{\widehat{\bm{y}}_{p\{\ell-n_u\}}, ...,$ $ \widehat{\bm{y}}_{p\{\ell\}}\}$$\in\amsmathbb{R}^{n_y}$  into one point $\widehat{\bm{y}}_{p,k}$$\in\amsmathbb{R}^{n_y}$ and a matrix $\widehat{\nabla_{\bm{u}} \bm{y}}_{p,k}$$\in\amsmathbb{R}^{n_y\times n_u}$ may seem inappropriate since  $n_u+1$ data are replaced by $n_u+1$ data, the justification for this name will be clearer in the Chapter~\ref{Chap:6_ASP}.
However, one can already envision that if one were to make S-ASP a bit more complex by using a central finite difference method to evaluate the plant gradients (which would involve minor modifications of the consistency monitor, of the experiments designer, and of equation \eqref{eq:4___3_12_nabla_y_bar}), then this method would indeed compress $2n_u+1$  points in $\amsmathbb{R}^{n_y}$ into one point in $\amsmathbb{R}^{n_y}$ and a matrix in $\amsmathbb{R}^{n_y\times n_u}$. This is indeed a reduction of the number of variables in  $\mathcal{D}_{II\{\ell\}}^{\checkmark}$  from $(2n_u+1)n_y$ to $(n_u+1)n_y$.
 
\section{The model-data combiner}

\begin{figure}
	\centering
	\includegraphics[width=14cm]{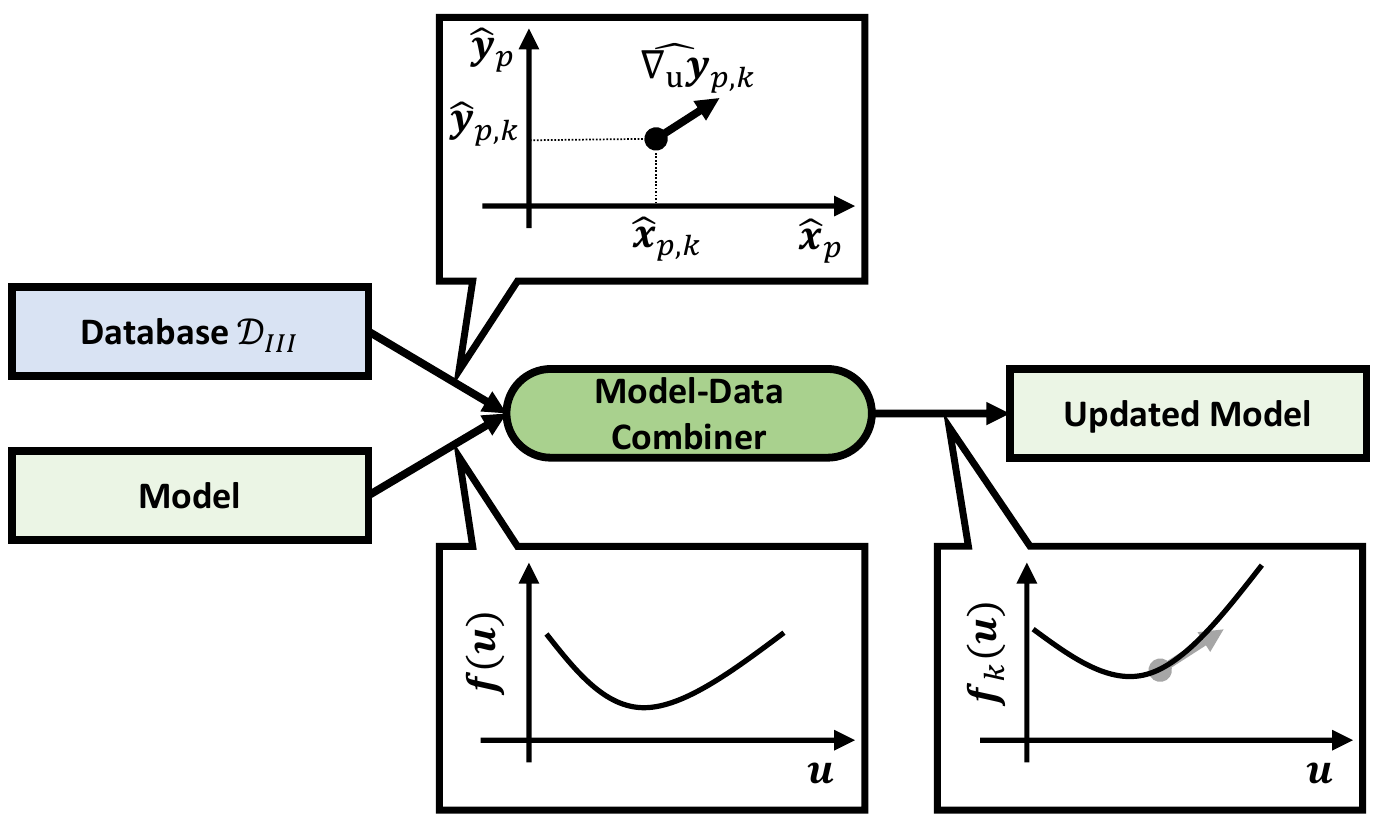}
	\caption{Information flow entering and leaving the Model-Data Combiner (indirect approach)}
	\label{fig:5_Combineur_donnees_Model}
\end{figure}

As illustrated on Figure~\ref{fig:5_Combineur_donnees_Model}, the role of the model-data combiner is to combine the data of $\mathcal{D}_{III}$ and the model built by the engineers to create an updated model. The latter should be able to accurately predict the optimality conditions of the plant, specifically when the S-ASP converges on an operating point.  

The properties that this updated model should have to enable the emergence of this behavior have been identified in Chapter~\ref{Chap:2_Vers_Une_meilleure_Convergence} with Theorem~\ref{thm:2___1_AffineCorrection_Impl_KKTmathcing}, and one has proposed simple ways to enforce these conditions with the \textit{direct} and \textit{indirect} approaches based on the four simplifications listed in section~\ref{sec:2_4_ISO_ISOy}. There are therefore two \textit{simple} approaches to build a model-data combiner:

(Direct approach): By using the measurements $\{\widehat{\bm{x}}_{p,k}, \widehat{\bm{y}}_{p,k}, \widehat{\nabla_{\bm{u}}\bm{y}}_{p,k}\}$ to apply an affine correction to the model's cost and constraint predictions: 
\begin{align}
	\phi_k(\bm{u}) := \ & 	\phi(\bm{u},\bm{f}(\bm{u})) + 
	\mu_{k}^{\phi}(\bm{u}), \\
	\bm{g}_k(\bm{u}) := \ & \bm{g}(\bm{u},\bm{f}(\bm{u})) + 
	\mu_{k}^{\bm{g}}(\bm{u}), 
\end{align}
where:
\begin{align}
	\mu_{k}^{\phi}(\bm{u})	:= \ & 
	\left(
		\left(
			\frac{\partial \phi}{\partial \bm{u}} + 
			\frac{\partial \phi}{\partial \bm{y}}  \widehat{\nabla_{\bm{u}}\bm{y}}_{p,k}
		\right)
		\Big|_{\substack{
			\bm{u}=\widehat{\bm{u}}_k\\ \bm{y}=\widehat{\bm{y}}_{p,k}
		}}
		-
		\left( 
			\frac{\partial \phi}{\partial \bm{u}} +  
			\frac{\partial \phi}{\partial \bm{y}}  \nabla_{\bm{u}}\bm{f}
		\right)
		\Big|_{\substack{
			\bm{u}=\widehat{\bm{u}}_k\\ \bm{y}=\bm{f}(\widehat{\bm{u}}_k)
		}}
	\right) (\bm{u}-\widehat{\bm{u}}_k), \nonumber \\
	\bm{\mu}_{k}^{\bm{g}}(\bm{u})	:= \ & 
	\left(
		\left(
			\frac{\partial \bm{g}}{\partial \bm{u}} + 
			\frac{\partial \bm{g}}{\partial \bm{y}}  \widehat{\nabla_{\bm{u}}\bm{y}}_{p,k}
		\right)
		\Big|_{\substack{
				\bm{u}=\widehat{\bm{u}}_k\\ \bm{y}=\widehat{\bm{y}}_{p,k}
		}}
		-
		\left( 
			\frac{\partial \bm{g}}{\partial \bm{u}} +  
			\frac{\partial \bm{g}}{\partial \bm{y}}  \nabla_{\bm{u}}\bm{f}
		\right)
		\Big|_{\substack{
				\bm{u}=\widehat{\bm{u}}_k\\ \bm{y}=\bm{f}(\widehat{\bm{u}}_k)
		}}
	\right) (\bm{u}-\widehat{\bm{u}}_k) ... \nonumber \\
	 & + \bm{g}\big(\widehat{\bm{u}}_k,\widehat{\bm{y}}_{p,k}\big) 
	 	- \bm{g}\big(\widehat{\bm{u}}_k, \bm{f}(\widehat{\bm{u}}_k)\big).
\end{align}
It is obvious that in the ideal case where: 
\begin{equation} \label{eq:4_16_cas_ideal}
	\widehat{\nabla_{\bm{u}}\bm{y}}_{p,k} 
	= 
	\nabla_{\bm{u}}\bm{f}_{p} \big|_{
			\bm{u}=\widehat{\bm{u}}_k
	} ,
\end{equation}
the following equations are true: 
\begin{align}
	\bm{g}_k(\widehat{\bm{u}}_k)                 = \ & \bm{g}\big(\widehat{\bm{u}}_k, \bm{f}_p(\widehat{\bm{u}}_k)\big), & 
	\nabla_{\bm{u}}\bm{g}_k\big(\widehat{\bm{u}}_k) = \ & \nabla_{\bm{u}}\bm{g}\big(\widehat{\bm{u}}_k, \bm{f}_p(\widehat{\bm{u}}_k)\big) \label{eq:4_17_cas_idea_2}\\
	 & & 
	\nabla_{\bm{u}}\phi_k\big(\widehat{\bm{u}}_k) = \ & \nabla_{\bm{u}}\phi\big(\widehat{\bm{u}}_k, \bm{f}_p(\widehat{\bm{u}}_k)\big).\label{eq:4_18_cas_idea_3}
\end{align}
As this is an ideal case, it is unlikely that in practice equalities \eqref{eq:4_16_cas_ideal}-\eqref{eq:4_18_cas_idea_3} are true. However, if they are not a strict equality, one can hope that they are good approximations.

(Indirect approach):  By using the measurements $\{\widehat{\bm{x}}_{p,k}, \widehat{\bm{y}}_{p,k}, \widehat{\nabla_{\bm{u}}\bm{y}}_{p,k}\}$ to apply an affine correction to the model outputs: 
\begin{align}
	\bm{f}_k(\bm{u}) := \ & \bm{f}(\bm{u}) + 
	\bm{\mu}_{k}^{\bm{y}}(\bm{u}), 
\end{align}
where:
\begin{align}
	\bm{\mu}_{k}^{\bm{y}}(\bm{u})	:= \ & 
	\overline{\bm{y}}_{p,k}
	- \bm{f}(\widehat{\bm{u}}_k) +
	\left(
		\widehat{\nabla_{\bm{u}}\bm{y}}_{p,k}
		- 
		\nabla_{\bm{u}}\bm{f}|_{\widehat{\bm{u}}_k}
	\right) (\bm{u}-\overline{\bm{u}}_k).
\end{align}
Again, it is clear that in the ideal case where \eqref{eq:4_16_cas_ideal} is true, then 
\begin{align}
	\bm{f}_k(\widehat{\bm{u}}_k)                 = \ & \bm{f}_p(\widehat{\bm{u}}_k), & 
	\nabla_{\bm{u}}\bm{f}_k\big(\overline{\bm{u}}_k) = \ &  \bm{f}_p(\widehat{\bm{u}}_k), \label{eq:4_19_cas_idea}
\end{align}
and the functions $\bm{g}_k(\bm{u}) := \bm{g}\big(\bm{u},\bm{f}_k(\bm{u})\big)$, and $\phi_k(\bm{u}) := \phi\big(\bm{u},\bm{f}_k(\bm{u})\big)$ satisfy the equalities \eqref{eq:4_17_cas_idea_2} and \eqref{eq:4_18_cas_idea_3}. But as for the direct approach, it is unlikely that in practice equalities  \eqref{eq:4_16_cas_ideal} and \eqref{eq:4_19_cas_idea} are perfectly true. In general, it will be approximations, and thus the condition of applicability of the Theorem~\ref{thm:2___1_AffineCorrection_Impl_KKTmathcing} are only approximately verified.

\section{The decision making}

The construction of decision-related functions is strongly influenced by the conception  one may have of the RTO-plant interface. In the vast majority of RTO related research articles, this interface is chosen as ``simple'' setpoints for plant controllers. Since this ``simple'' choice is appropriate for S-ASP, one does not discuss possible alternatives in this chapter. Such a discussion is left for the Chapter~\ref{Chap:6_ASP}.

\subsection{The decision trigger}
\label{sec:4_4_1_Declencheur_de_decision}

\begin{figure}[H]
	\centering
	\includegraphics[width=14cm]{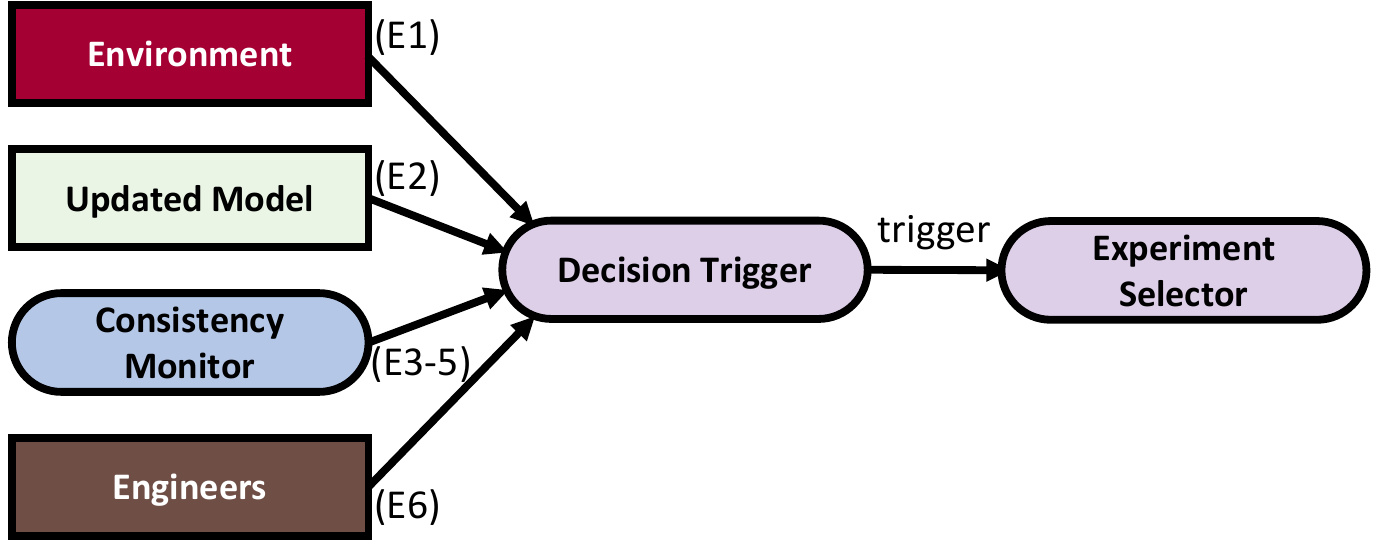}
	\caption{The flow of information entering and leaving the decision trigger}
	\label{fig:6_Declencheur_de_decision}
\end{figure}

As illustrated in Figure~\ref{fig:6_Declencheur_de_decision}, the role of the decision trigger is to trigger a decision-making by the experiment selector whenever one of these 6 events (E) occurs:
\begin{itemize}
	\item (E1) When a cost or a constraint changes.    
	\item (E2) When a new updated model is available. 
	\item (E3) When one is not satisfied with the current experiment, e.g. if it violates constraints. 
	\item (E4) When a ``significant'' variation of a measured disturbance is measured.
	\item (E5) When an anomaly in the plant's behavior is detected. e.g. failure or other unexplained variations.
	\item (E6) When the engineers want it.
\end{itemize} 
In fact, when one of these events occurs, it is generally desired that the autopilot reacts immediately. 

In the framework of the S-ASP, only the events (E1), (E2) and (E6) are considered. Indeed, these events are extremely easy to detect since it is enough to detect a human intervention on a software or the variation of a virtual object. On the other hand, the other events are much more complex to define mathematically and therefore to identify automatically. For example, measured disturbances vary constantly due to measurement noise and their random nature.  It is therefore necessary to define what a ``significant'' variation is and this is anything but obvious.  That's why in this S-ASP, their detection is not integrated into the autopilot and this task is left to the engineers who, when they detect them, can manually trigger a decision-making. 
\begin{itemize}
	\item If the engineers detect  (E3):  Then the last measurement of the plant outputs are considered as measurements of its static state (even if they are measurements of its transient state), and the emergency modes of the compressor and the Model-data combiner are activated (see section \ref{sec:4_7_1_URGENCE}) and the resulting updated model is used to define the new experiments to be implemented.    
	\item If the engineers detect (E4) or (E5): Then the autopilot is reinitialized, i.e. the databases are emptied and the updated model returns to the nominal model. 
\end{itemize}

\begin{figure}
	\centering
	\includegraphics[width=14cm]{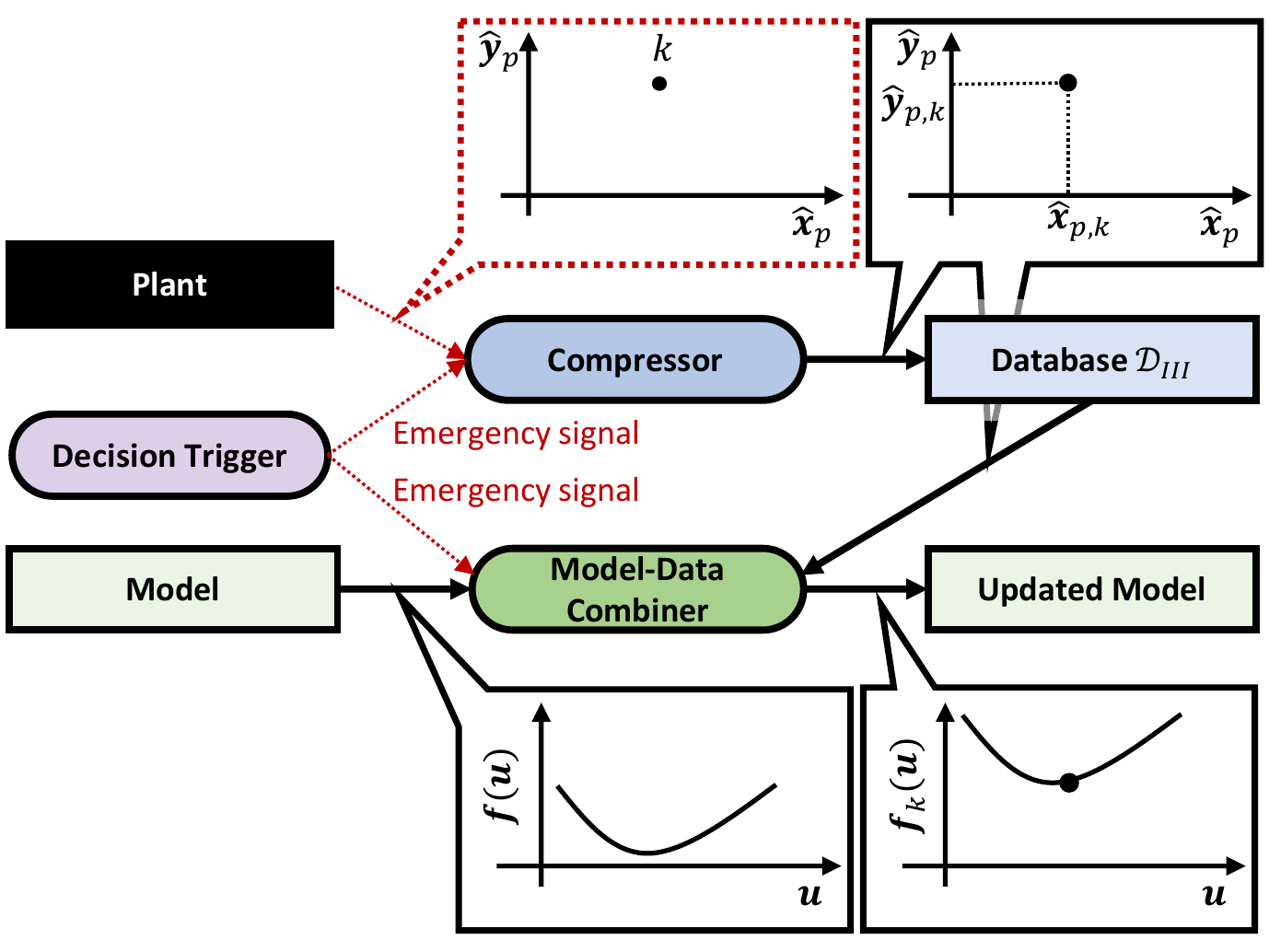}
	\caption{The flows of information when the emergency mode is activated}
	\label{fig:7_Mode_URGENCE}
\end{figure}

\subsubsection{Emergency mode -- (E3) detected}
\label{sec:4_7_1_URGENCE}

When the emergency mode is activated, the database $\mathcal{D}_{III}$ is emptied and the compressor takes the most recent measurement of the plant outputs (whether it is at steady state or not) and add them to the emptied database $\mathcal{D}_{III}$.   Then, the model-data combiner uses this single measurement to apply a 0th order correction to the nominal model.

(Direct approach): By using the measures $\{\widehat{\bm{x}}_{p,k}, \widehat{\bm{y}}_{p,k}\}$ to apply an offset correction to the predictions of the model constraints:  
\begin{align}
	\phi_k(\bm{u}) := \ & 	\phi\big(\bm{u},\bm{f}(\bm{u})\big), \\
	\bm{g}_k(\bm{u}) := \ & \bm{g}\big(\bm{u},\bm{f}(\bm{u})\big) + \bm{g}\big(\widehat{\bm{u}}_k,\widehat{\bm{y}}_{p,k}\big) 
	- \bm{g}\big(\widehat{\bm{u}}_k, \bm{f}(\widehat{\bm{u}}_k)\big).
\end{align}

(Indirect approach): By using the measures $\{\widehat{\bm{x}}_{p,k}, \widehat{\bm{y}}_{p,k}\}$ to apply an offset correction  to the model output predictions:
\begin{align}
	\bm{f}_k(\bm{u}) := \ &  \bm{f}(\bm{u}) +  \widehat{\bm{y}}_{p,k} -\bm{f}(\widehat{\bm{u}}_k), \\
	\phi_k(\bm{u}) := \ & 	\phi\big(\bm{u},\bm{f}_k(\bm{u})\big), \\
	\bm{g}_k(\bm{u}) := \ & \bm{g}\big(\bm{u},\bm{f}_k(\bm{u})\big).
\end{align}

\begin{Remark}
	When the emergency mode is activated, the direct approach's reaction is comparable to what a RTO method called constraints adaptation (CA -\cite{Chachuat:08}) would do. And CA is particularly well suited for maintaining a plant in its feasibility domain. 
\end{Remark}

\subsection{The experiment selector}

\begin{figure}[h!]
	\centering
	\includegraphics[width=14cm]{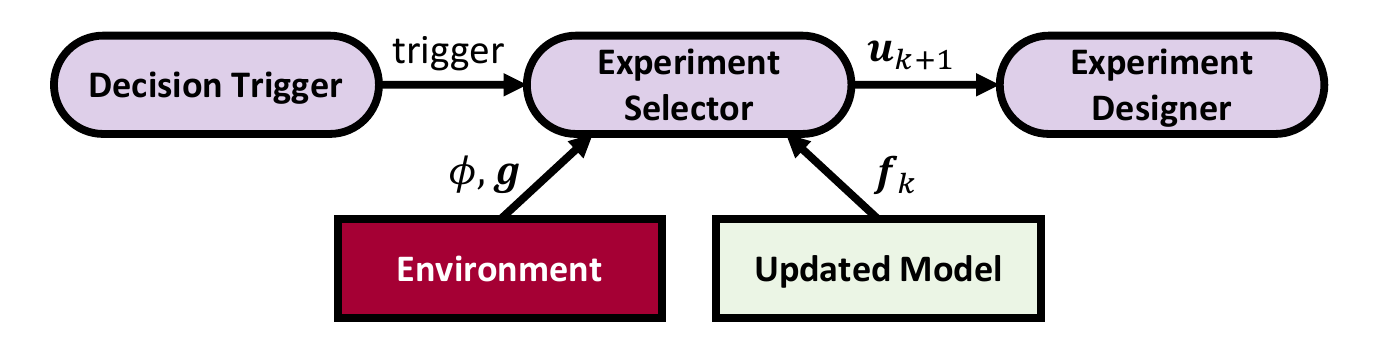}
	\caption{The flows of information entering and leaving the experience selector}
	\label{fig:9_Selecteur_d_experience}
\end{figure}

When the decision trigger decides that a decision should be made, the role of the experiment selector is to select a new operating point around which the future experiments should be done  (see Figure~\ref{fig:9_Selecteur_d_experience}).

A fairly simple and efficient way to make this choice is to apply steps 3) to 6) of KMFCaA (or KMFCA). Indeed, thanks to the Theorems of the Chapters~~\ref{Chap:2_Vers_Une_meilleure_Convergence} and \ref{Chap:3_KMA}, it is known that such a way of choosing the iterations offers good \textit{theoretical} guarantees of convergence provided that the plant values and gradients are available without errors. As in practice these values are not available without approximation errors, one does not expect to reach these theoretical performances, but one hopes to get close to them. 

\subsection{The experiment designer}
\label{sec:4___Designer_d_experience}
\begin{figure}[H]
	\centering
	\includegraphics[width=14cm]{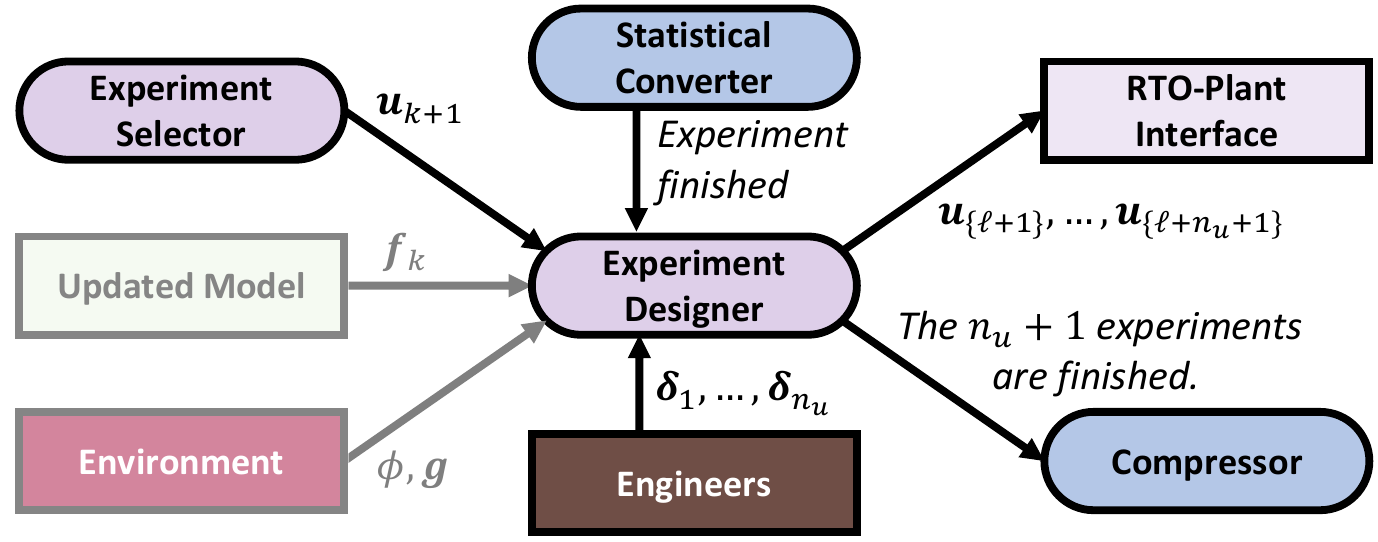}
	\caption{The flows of information entering and leaving the experiment designer}
	\label{fig:10designersexperience}
\end{figure}

When the experiment designer is instructed by the experiment selector to perform experiments around the point $\bm{u}_{k+1}$, its role is to define what these experiments are (i.e. the points to explore) and how they should be conducted (e.g. the time to spend at the explored point to obtain sufficiently accurate estimates). 

The easiest way to do this is to always do the same thing. As S-ASP is partially based on KMFCaA, one knows that at each iteration, the values and gradients of the plant are required given that one wants to apply simplification 2 of Section~\ref{sec:2_4_ISO_ISOy}. So, an idea would be to implement a forward finite difference method which means giving the RTO-plant interface the following successive instructions: 
\begin{align*}
	\bm{u}_{\{\ell+1\}} := \ & \bm{u}_{k+1}, & 
	\bm{u}_{\{\ell+2\}} := \ & \bm{u}_{k+1}+\bm{\delta}_1,&
	& ..., &
	\bm{u}_{\{\ell+n_u+1\}} := \ & \bm{u}_{k+n_u+1}+\bm{\delta}_{n_u}.
\end{align*}
where $\bm{\delta}_i$ is $\forall i \in\{1,n_u\}$ a vector in the direction of the $i$-the axis of the input space, and the norm of this vector must be chosen by the engineers. A strategy for choosing these values is proposed in \cite{Costello:2016}\footnote{A step which is too small could be too affected by the measurement noise, and a step which is too large could be too affected by the curvature of the model functions. So choosing the $\bm{\delta}_i$ is trying to make the best compromise between these two phenomena.}. The transition from one experiment to the next takes place each time the statistical converter informs the experiment designer that an experiment has been completed. Once all the $n_u+1$ experiments are finished, the experiment designer informs the compressor that the database $\mathcal{D}_{II}$ is ``ready to be compressed''.

On the Figure~\ref{fig:4___1_StructureRTO_GENERALE}, the experiment designer is also connected to the updated model and to the environment because one could consider using them to, for example, impose the updated-model-based feasibility of $\{\bm{u}_{\{\ell+2\}}, ...,\bm{u}_{\{\ell+n_u+1\}} \}$. Indeed, if one knows that $\bm{u}_{\{\ell+1\}}=\bm{u}_{k+1}$  is feasible according to the updated model (see Chapter~\ref{Chap:3_KMA}), this is not the case for the other points $\{\bm{u}_{\{\ell+2\}}, ...,\bm{u}_{\{\ell+n_u+1\}} \}$. Since introducing such a feature would make the experiment designer and compressor more complex, it is not included in S-ASP.

Finally, stand-by functionality is added, which is triggered by the following condition: 
\begin{equation} \label{eq:4___27_Critere_arret_S_ASP}
	|\bm{u}_{k+1} - \bm{u}_k| \leq \bm{a},
\end{equation}
where $\bm{a}\in\amsmathbb{R}^{n_u +}$ is a vector that characterizes the stopping criterion of the autopilot's exploration and that must be defined by the engineers. In other words, when inequality \eqref{eq:4___27_Critere_arret_S_ASP} holds, the autopilot stops running experiments until one of the events that triggers the decision trigger occurs. In the case of S-ASP, only events (E1) and (E6) can restart the experimentation phase.

\section{Summary}

\begin{figure}[H]
	\centering
	\includegraphics[width=14cm]{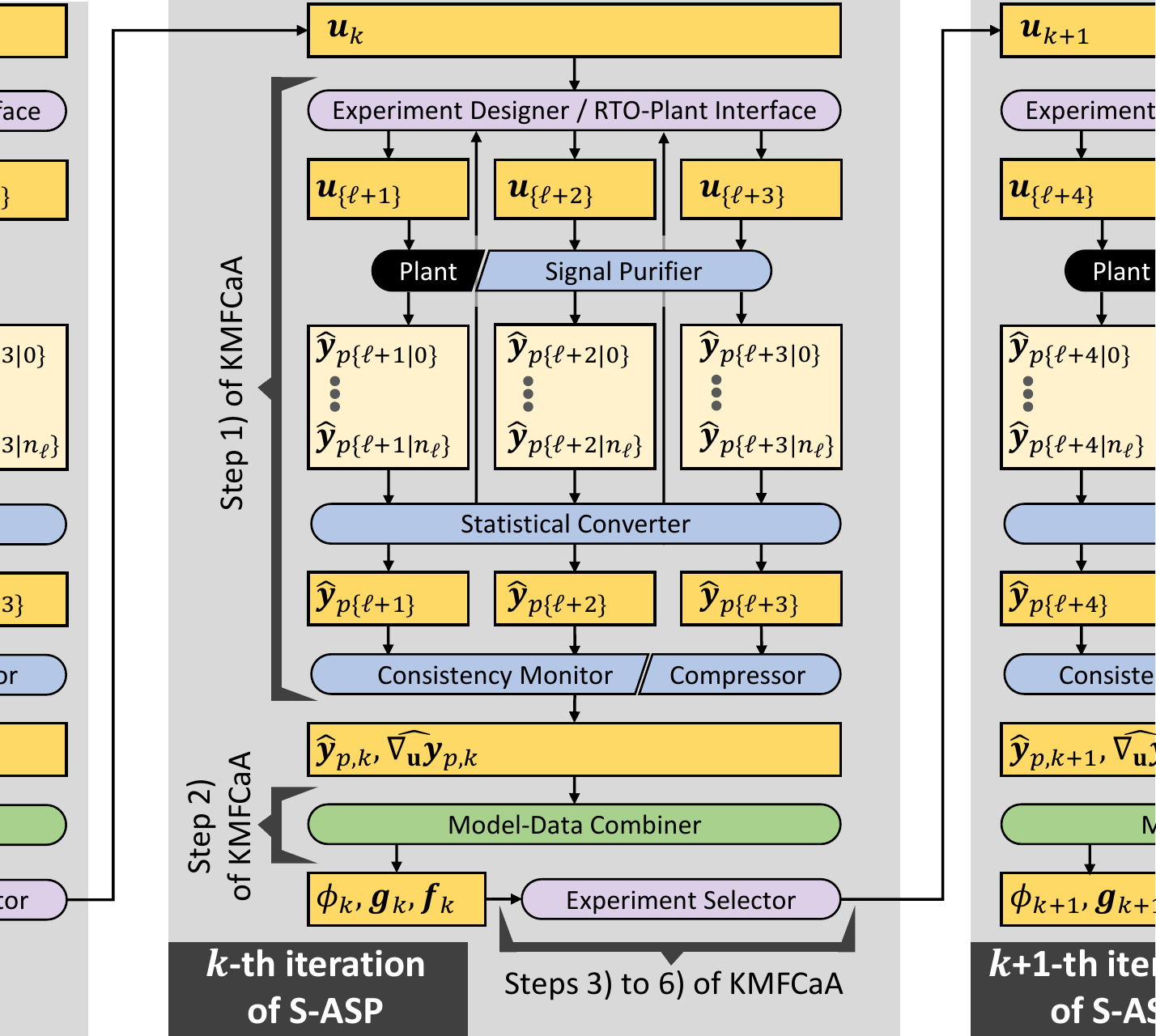}
	\caption{Description of the notations and variables used in one iteration of S-ASP}
	\label{fig:StructureRTO_GENERALE_V2}
\end{figure}

In summary, the different functions of the S-ASP can be combined to, to some extent, find the different steps of the Algorithm~\ref{algo:KMFCaA} (or \ref{algo:KMFCaAy}), as illustrated on Figure~\ref{fig:StructureRTO_GENERALE_V2}:
 \vspace{-\topsep}
\begin{itemize}[noitemsep]
	\item The step 1) of KMFCaA is performed  by the whole data management part plus the experiment designer and the RTO-plant interface, 
	\item The step 2) of KMFCaA is performed by the model-data combiner,
	\item The steps 3) to 6) are performed by the experiment selector,
	\item The step 7) is performed by the decision trigger.  
\end{itemize}
Finally, in addition to Figure~\ref{fig:4___1_StructureRTO_GENERALE} which provides a global view of the interactions between the sub-functions of the S-ASP, Figure~\ref{fig:StructureRTO_GENERALE_V2} offers another global view more focused on the variables, their notations, and which assimilates parts of S-ASP to the different steps of the Algorithm~\ref{algo:KMFCaA} (or \ref{algo:KMFCaAy}).  

\clearpage 

\section{Implementation of S-ASP}

\subsection{The William-Otto reactor}

The standard benchmark case study for RTO considered here is the continuous stirred-tank reactor of \cite{Williams:60}, where the three following reactions take place (for the plant): 
\begin{align}
	\text{A} + \text{B} & \overset{k_{p,1}}{\longrightarrow}   \text{C},           &  k_{p,1} = A_{p,1}e^{-B_{p,1}/(T_R+273.15)} \label{eq:reaction1} \\
	\text{C} + \text{B} & \overset{k_{p,2}}{\longrightarrow}   \text{P} + \text{E}, &  k_{p,2} = A_{p,2}e^{-B_{p,2}/(T_R+273.15)} \\
	\text{P} + \text{C} & \overset{k_{p,3}}{\longrightarrow}   \text{G}.           &  k_{p,3} = A_{p,3}e^{-B_{p,3}/(T_R+273.15)}, \label{eq:reaction3}
\end{align}

The reactants A and B are fed separately to the reactor, with mass flowrates of $F_A$ and $F_B$, respectively. Components P and E are the desired products, C is an intermediate product and G is an undesired by-product. The reactor is operated isothermally at a controlled temperature $T_R$. Steady-state mass balances can be found in \cite{Zhang:00}. The same optimization problem as in \cite{Marchetti:2017} is considered, wherein the input variables are $\bm{u} = [F_{A},F_{B}, T_{R}]^{\rm T}$ and the outputs are $\bm{y}_p=[X_E, X_P, X_G]^{\rm T}$ with $X_i$ denoting the mass fraction of component $i$. 

There is significant structural plant-model mismatch since the model \emph{only considers two reactions} \cite{Forbes:96, Marchetti:2017}. 
The objective is to maximize profit at steady-state, while satisfying an upper bound on $X_G$ and input bounds:
\begin{align}
	\underset{\bm{u}}{\operatorname{max}} \quad &  \phi(\bm{u},\bm{y}_p) = (1143.38 X_P + 25.92 X_E) (F_A+F_B) - 76.23 F_A - 114.34 F_B \label{eq:cost_wotto} \\
	\text{s.t.} \quad & {g}(\bm{y}_p) = X_G - 0.08 \leq {0}, \ \ \ \text{if mode = 1,} \label{eq:constr_wotto} \\
	& F_{A} \in  [3,\ 4.5] \ \ (\text{kg/s}), \quad F_{B} \in  [6,\ 11] \ \ (\text{kg/s}),
	\quad 
	T_R \in [80,\ 105]\ \ (^{\circ}{\rm C}). \nonumber
\end{align} 
The reactor can operate in two modes. \textit{Mode 1} involves a constraint on the concentration of G in the reactor outflow, while \textit{mode 2} does not.  

\subsection{Simulation}

Some information about the implementation: 
\vspace{-\topsep}
\begin{itemize}[noitemsep]
	\item S-ASP is implemented and two simplifying assumptions are made: (i) the SS detector is ideal, and (ii) no measurement uncertainty.
	\item An indirect approach is chosen for the model-data combiner.
	\item Only the responses to events (E1) and (E2) are implemented in the decision trigger. 
	\item The emergency mode is not illustrated. 
	\item The experiment selector corresponds to steps 3) to 6) of KMFCaA-I.
	\item The experiment designer uses the following perturbation values for the estimation of the plant's gradients: $\bm{\delta}_1 = 10^{-7}$, $\bm{\delta}_2 = 10^{-7}$, $\bm{\delta}3 = 10^{-6}$.
	\item No stand-by mode is implemented.
\end{itemize}
The considered scenario consists in simulating 30 iterations in mode 1 and 30 in mode 2. The results are shown in Figure~\ref{fig:4_Casestudy_WO_Sim_results}. 
\clearpage

\noindent
\begin{minipage}[h]{\linewidth}
	\vspace*{0pt}
	\centering
	\includegraphics[width=4.45cm]{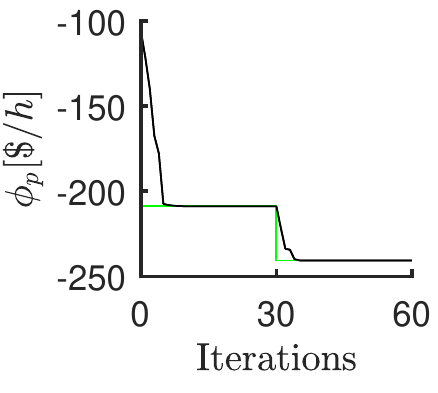}\hskip -0ex
	\includegraphics[width=4.45cm]{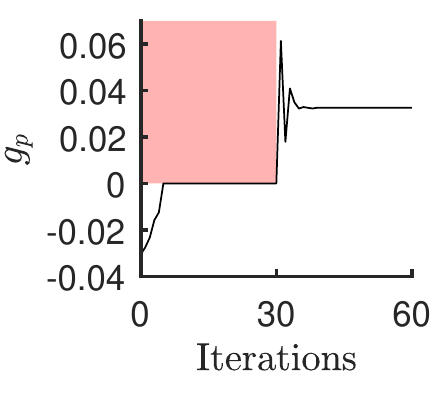}\hskip -0ex
	\includegraphics[width=4.45cm]{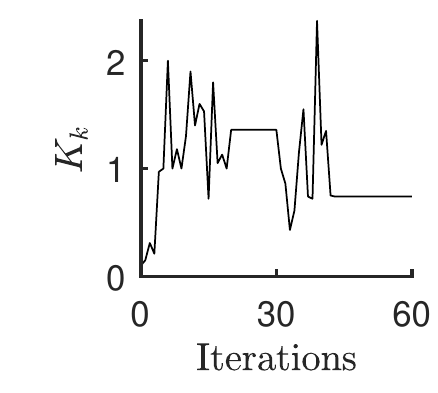} \\
	
	\includegraphics[width=4.45cm]{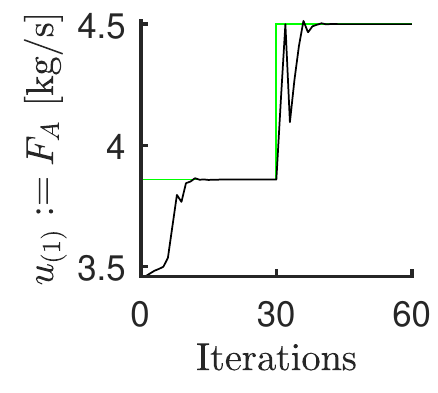}\hskip -0ex
	\includegraphics[width=4.45cm]{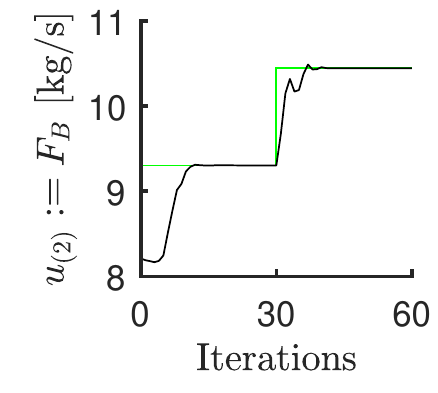}\hskip -0ex
	\includegraphics[width=4.45cm]{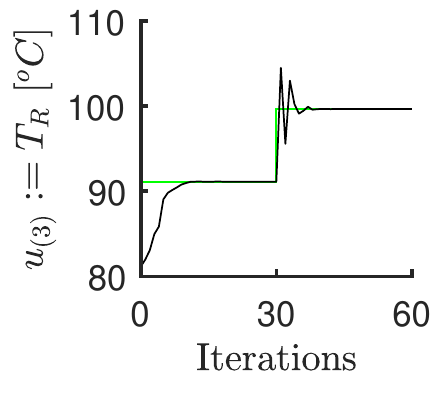} \\
	
	\captionof{figure}{Simulation results}
	\label{fig:4_Casestudy_WO_Sim_results}
\end{minipage} \\

One can see that S-ASP is able to converge on the optimum of the plant in about ten iterations (either for mode 1 or 2). One can also see that the value of the adaptive filter can be greater than 1 and that it systematically ends up converging on a limit value.  Finally, one can see that S-ASP is able to efficiently manage a variation of the environment corresponding to a change of its constraints.

\section{Conclusion}

In this section the structure of an autopilot for a plant designed to run on SS is presented. It is discussed what should be that compose it the functions, as well as what are the interactions that these functions should have in order to efficiently drive a plant. Finally, this structure is partially illustrated on a classic RTO case study.

	\chapter{An algorithm that makes better use of data}
		\label{Chap:5_IMA}

\section{An observation}

\subsection{The classical approach}

\begin{figure}[H]
	\vspace*{0pt}
	{\centering
		\includegraphics[width=13.8cm]{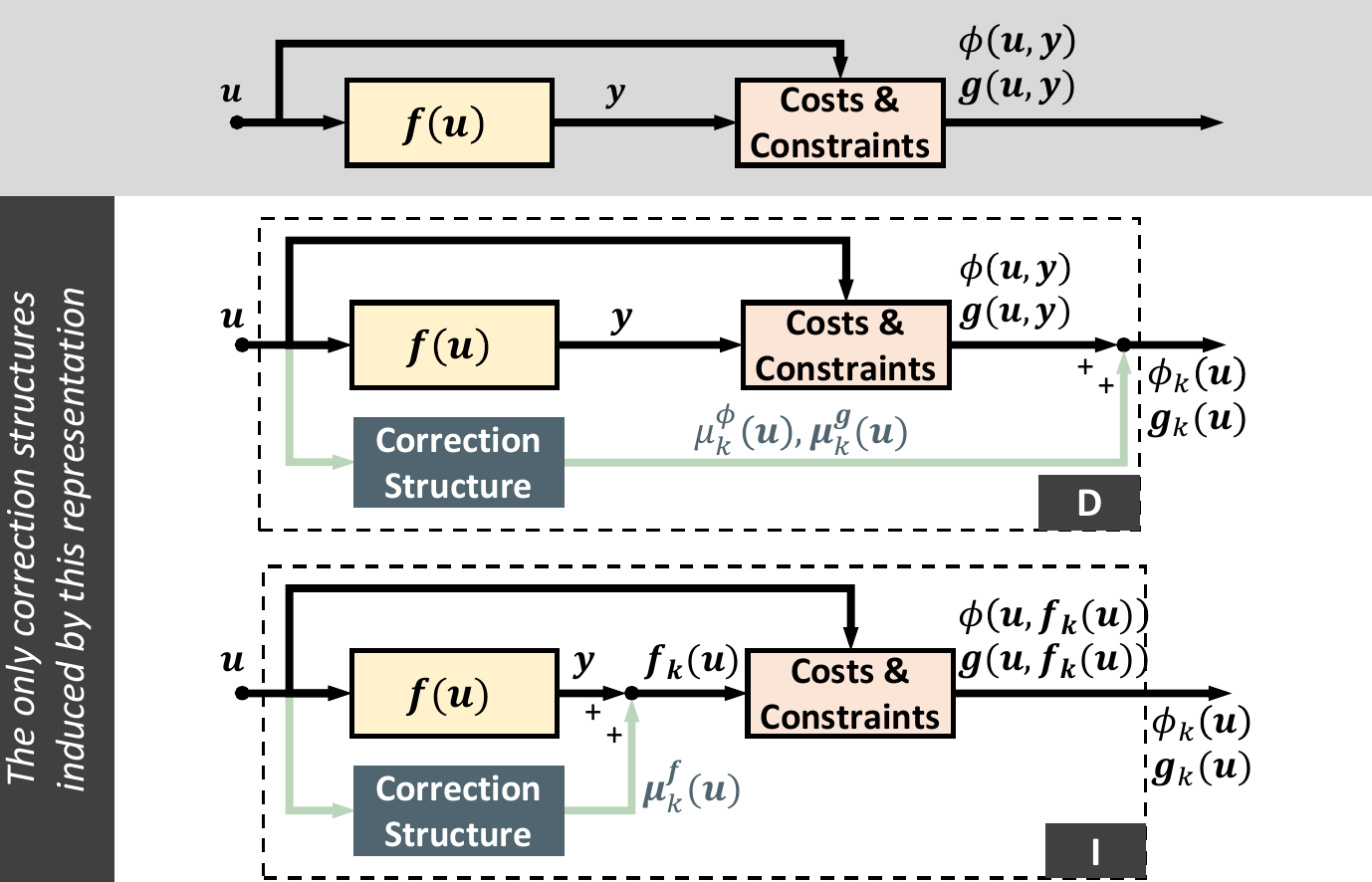} \\
	}
	\caption{Classic direct (D) and indirect (I) correction structures}
	\label{fig:5___1_Observation2}
\end{figure}

The set of ways to correct the predictions of the costs and constraints of a model is entirely predefined by the representation one uses of this model. If, for instance, one \textit{chooses} to interpret it as a ``simple'' mapping between $\bm{u}$ and measures $\bm{y}=\bm{f}(\bm{u})$, then only two types of correction structure are possible. 
They are shown on Figure~\ref{fig:5___1_Observation2} and are called direct (D) and indirect (I) structures in reference to the direct and indirect approaches which are introduced in Chapter~\ref{Chap:2_Vers_Une_meilleure_Convergence}.

\subsection{An alternative approach} 

However, one can also \textit{choose} to represent the model as a network of submodels interconnected by measured or manipulated variables. To obtain such a representation, a model preprocessing method inspired of \cite{Papasavvas:2019c} is proposed, which can be summarized by the four steps illustrated on Figure~\ref{fig:Preprocessing}: 

\begin{figure}[H]
	\centering
	\includegraphics[width=14cm]{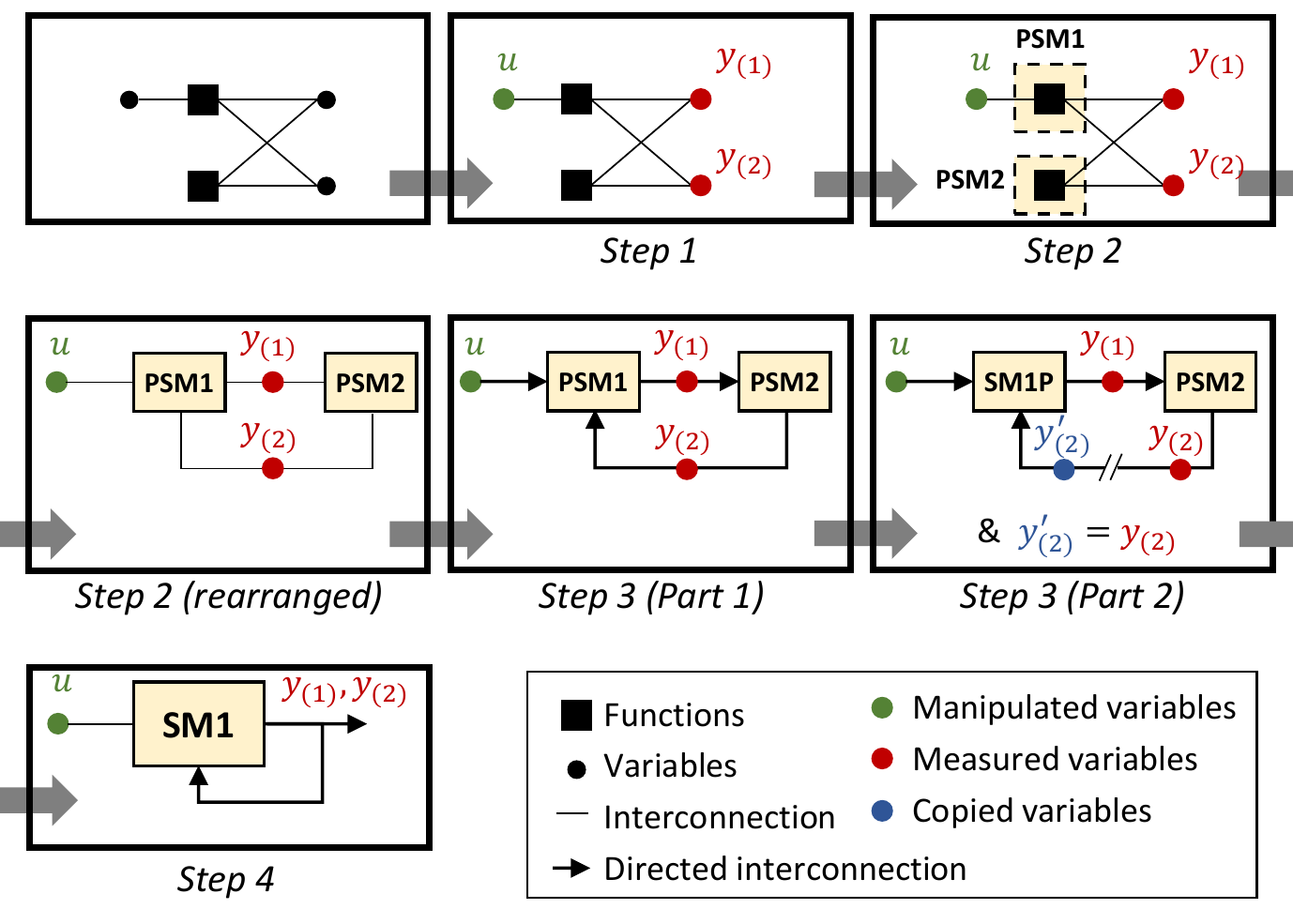}
	\caption{Introduction to the four steps of the model preprocessing  
	}
	\label{fig:Preprocessing}
\end{figure}

\textbf{Step 1 (Variable Clustering):} The variables of the model must be classified as either decision variables $\bm{u}$, or measured variables $\bm{y}$, or states (variables internal to the model). 

\textbf{Step 2 (Network Grouping):} The model can then be represented as a network of equations interconnecting these variables. A subset of this network that satisfies the following conditions is what one calls a primary submodel (PSM): (i) it contains only state variables and functions; (ii) it is connected only to decision or measured variables.

\textbf{Step 3 (Construction of a directed network and selection of the computational sequence):} Performing a simulation consists in finding the set of $\bm{y}$ which, for a given $\bm{u}$, satisfy all the equations of the model. To perform this task, one must first choose the inputs and outputs of each PSM. That is, transform the undirected network of couplings between the variables and functions of the model into a directed network. Then, one must choose how to solve the resulting system of equations  \cite{Barkley:1972,Upadhye:1975,Motard:1981}. 
There are also methods to simultaneously find the optimal directed network and sequence of computations for a given simulation  \cite{Montagna:1988}. However, as these methods are beyond the scope of this thesis, step 3 is systematically performed in an arbitrary way.

\textbf{Step 4 (Construction of the submodels (SMs)):}  Engineers can then aggregate the SMPs as they wish to form submodels (SMs). Indeed, if the model representation with PSMs is the most fragmented representation possible, one leaves the freedom to the engineers to reduce this fragmentation in order to either simplify the resulting RTO algorithm or to force the model consistency condition (which is defined afterwards, see Hypothesis~\ref{ass:5___1_Coherence_Modele}).

\begin{figure}[p]
	\centering
	\includegraphics[width=14cm]{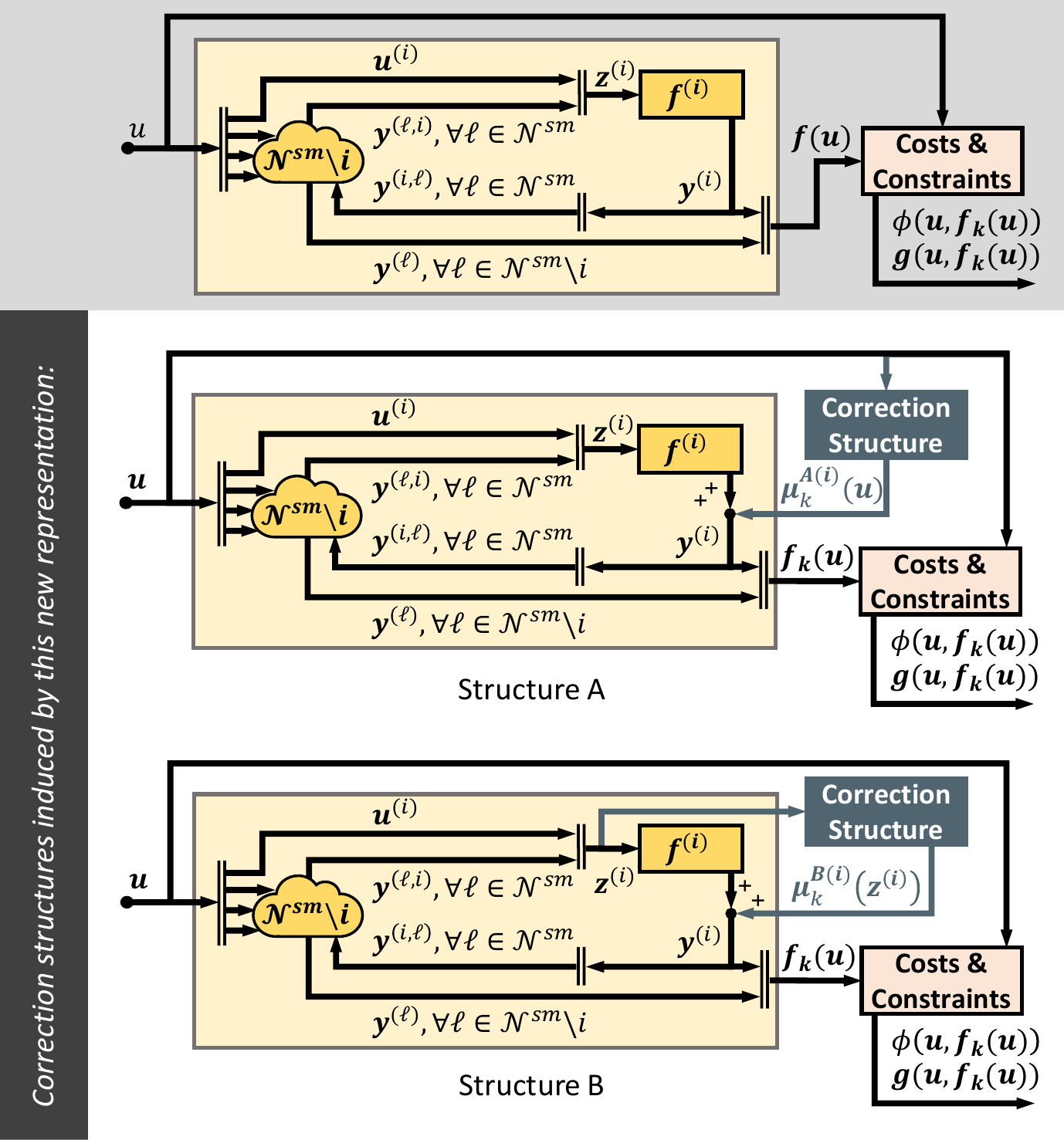}
	\caption{Correction structures A and B}
	\label{fig:5___2_Observation}
\end{figure}

This new system of representation reveals two new (indirect) corrective structures that are named structure A and B. These are illustrated in Figure~\ref{fig:5___2_Observation} where new notations are used whose definitions are the following: 
\begin{itemize}
\item The function: $\bm{f}$ is separated into a network $\bm{\mathcal{N}}^{SM}$ of $i=\{1,...,n_{SM}\}$ functions $\bm{f}^{(i)}$ which are the submodels (SMs). 
\item The inputs of the submodel $i$ (SM$i$) are noted $\bm{z}^{(i)}$ and are composed of\vspace{-\topsep}
\begin{itemize}[noitemsep]
	\item the subset of inputs $\bm{u}$ that are inputs of the SM$i$ which are named $\bm{u}^{(i)}$,
	\item the outputs of the SMs $\ell\in\bm{\mathcal{N}}^{SM}$ that are also inputs of SM$i$ which are named $\bm{y}^{(\ell,i)}$. 
\end{itemize}
In short:
\begin{equation}
	\bm{z}^{(i)} := \left[\bm{u}^{(i)\rm T}, \ \bm{y}^{(1,i)\rm T}, ..., \bm{y}^{(n_{SM},i)\rm T}\right]^{\rm T}.
\end{equation}
\item The outputs of SM$i$ are a subset of the outputs $\bm{y}$  which are named $\bm{y}^{(i)}$. 
\end{itemize}

\begin{figure}[p]
	\centering
	\includegraphics[width=14cm]{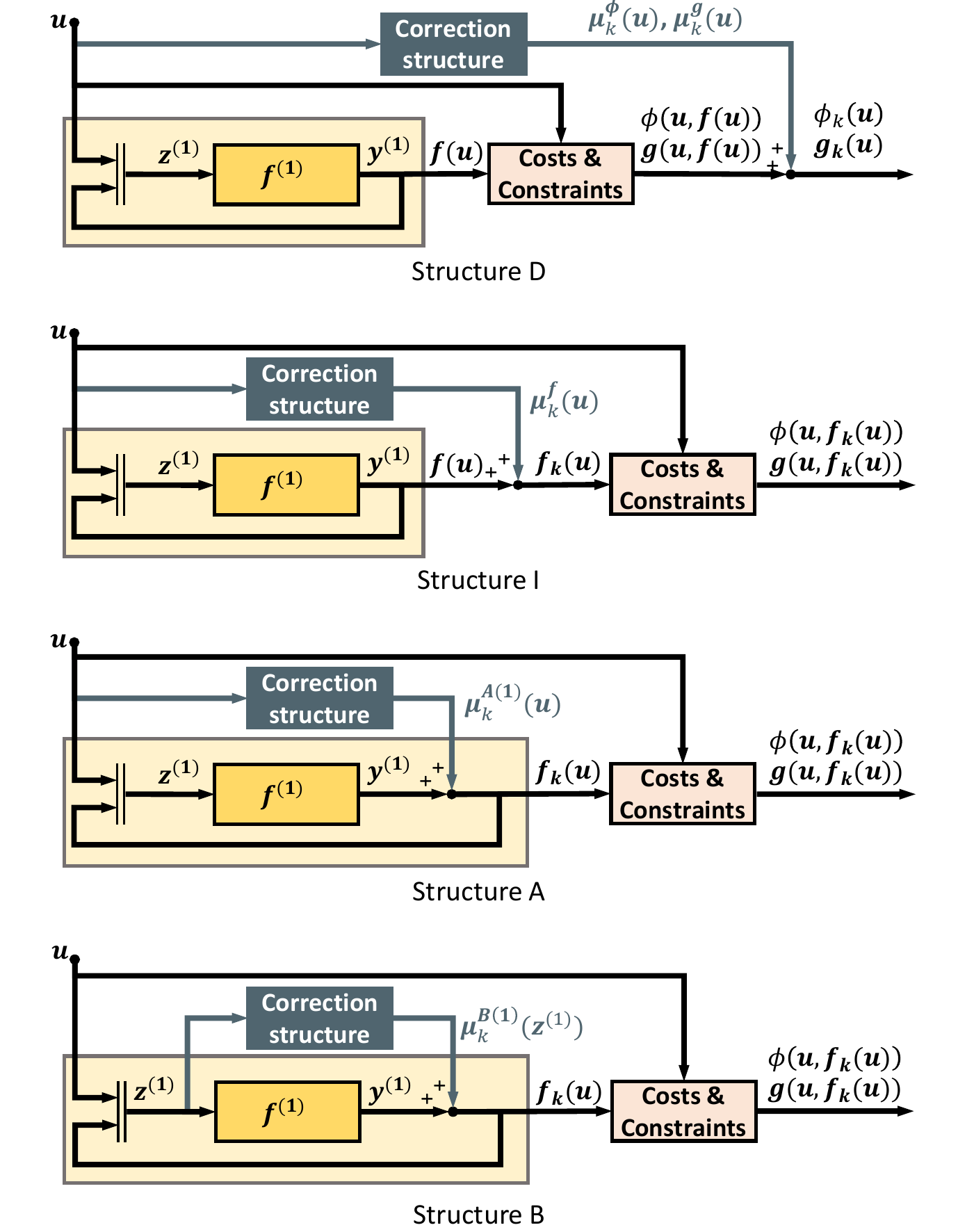}
	\caption{The four correction structures for the case where $n_{SM} = 1$}
	\label{fig:5___3_Observation_Un_SM}
\end{figure}

Finally, with Figure~\ref{fig:5___3_Observation_Un_SM} a graphical comparison of these four correction structures is proposed for the simple case where the number of SMs is only $1$. 

Now that one has identified two new possible correction structures, let us see if it is possible to use them to indirectly implement affine corrections on the costs and constraints predictions. 

\clearpage

\section{The structure A}

\subsection{Splitting the plant's function}

\begin{figure}[H]
	\centering
	\includegraphics[width=14cm]{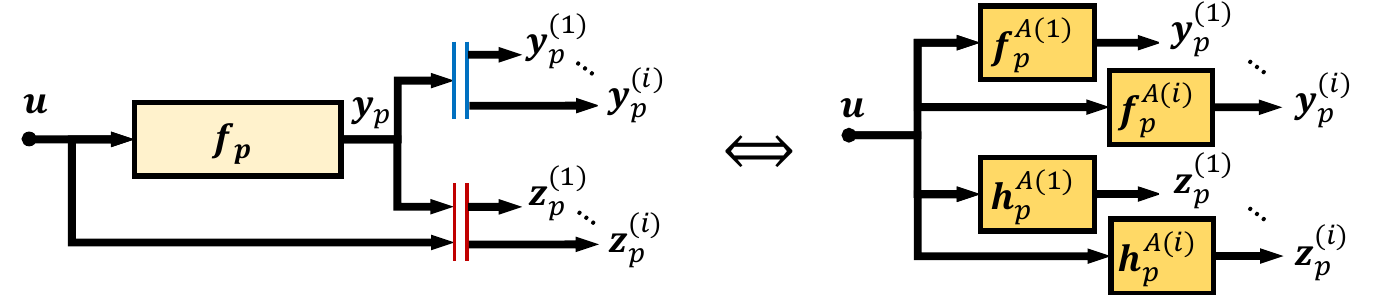}
	\begin{minipage}{6cm}
		a) ``Standard'' representation
	\end{minipage}
	\begin{minipage}{6cm}
		\quad b) Alternative representation 
	\end{minipage}
	\caption{Redefinition of the plant for the structure A}
	\label{fig:5___4_Plant}
\end{figure}

Now, what one wants is to change the representation of the plant so that:
\vspace{-\topsep}
\begin{itemize}[noitemsep]
	\item Part of its inputs are the inputs of the structure A, i.e. $\bm{u}$, and part of its outputs are the outputs of the SMs that the structure A corrects, i.e. $\bm{y}_p^{(i)}$ $\forall i$.
	\item The other part of its inputs are the plant inputs, i.e. $\bm{u}$, and the other part of its outputs are the inputs of the structure A, i.e. $\bm{z}_p^{(i)}$ $\forall i$.
\end{itemize}
It is quite simple to perform this task since it is sufficient to add a demixer at the level of the outputs of $\bm{f}_p$ to get the $\bm{y}_p^{(i)}$, $\forall i$, and a mixer-demixer combining the inputs $\bm{u}$ and the outputs $\bm{y}_p$ to get the $\bm{z}_p^{(i)}$, $\forall i$. Figure~\ref{fig:5___4_Plant} a)  illustrates the actions of this mixer (in blue) and of this mixer-demixer (in red).  One can then define a new structure to represent the plant. This is shown in Figure~\ref{fig:5___4_Plant} b) and is composed of functions $\bm{f}_p^{A(i)}$, $\forall i$, linking $\bm{u}$ to $\bm{y}_p^{(i)}$, and of functions $\bm{h}_p^{A(i)}$, $\forall i$, linking  $\bm{u}$ to $\bm{z}_p^{(i)}$, $\forall i$.  To finish the introduction of this new representation of the plant, one would like to insist on the fact that the existence of the functions $\bm{f}_p^{A(i)}$ and $\bm{h}_p^{A(i)}$, $\forall i$, does not imply any assumption since only mixers and demixers were used to define them.

\subsection{Computation of the correction function of a sub-model}

It is considered that a SM$i$ is ``corrected to the 1st order''  at a point $\bm{u}_k$, if for all inputs $\bm{z}_{p,k}^{(i)} = \bm{h}_p^{A(i)}(\bm{u}_k)$,
the values and gradients of its corrected outputs are equal to those of the plant, i.e:
\begin{align}
	\bm{f}\big|_{\bm{z}_{p,k}^{(i)}} + \bm{\mu}_k^{A(i)}\big|_{\bm{u}_k}  = \ & \bm{f}_p^{A(i)}\big|_{\bm{u}_k}, \label{eq:5___2_SM_cor_val} \\
	\nabla_{\bm{z}^{(i)}} \bm{f}^{(i)}\big|_{\bm{z}_{p,k}^{(i)}}
	\nabla_{\bm{u}} \bm{h}_p^{A(i)}\big|_{\bm{u}_k} 
	+ 
	\nabla_{\bm{u}} \bm{\mu}_k^{A(i)}\big|_{\bm{u}_k}
	 = \  
	& \nabla_{\bm{u}} \bm{f}_p^{A(i)}\big|_{\bm{u}_k}. \label{eq:5___3_SM_cor_grad}
\end{align} 
For these equalities to be systematically true, the correction function $\bm{\mu}_k^{A(i)}$ must be chosen such that it corresponds to the following Taylor series extension: 
\begin{align}
	\bm{\mu}_k^{A(i)}(\bm{u}) :=  \ &
	\bm{f}_p^{A(i)}\big|_{\bm{u}_k} - \bm{f}\big|_{\bm{z}_{p,k}^{(i)}}
	 + 
	\left(
	\nabla_{\bm{u}} \bm{f}_p^{A(i)}\big|_{\bm{u}_k} 
	-
	\nabla_{\bm{z}^{(i)}} \bm{f}^{(i)}\big|_{\bm{z}_{p,k}^{(i)}}
	\nabla_{\bm{u}} \bm{h}_p^{A(i)} \big|_{\bm{u}_k}
	\right)\left(
	\bm{u}-\bm{u}_k
	\right)... \nonumber  \\
	&
	+ 
	\mathcal{O}(\|\bm{u}-\bm{u}_k\|^2).
\end{align} 
The easiest choice is that $\bm{\mu}_k^{A(i)}$ is an affine function by setting $\mathcal{O}(\|\bm{u}-\bm{u}_k\|^2)=\bm{0}$.

\section{The structure B}

\subsection{Splitting the plant's function}
\begin{figure}[h]
	\centering
	\includegraphics[width=14cm]{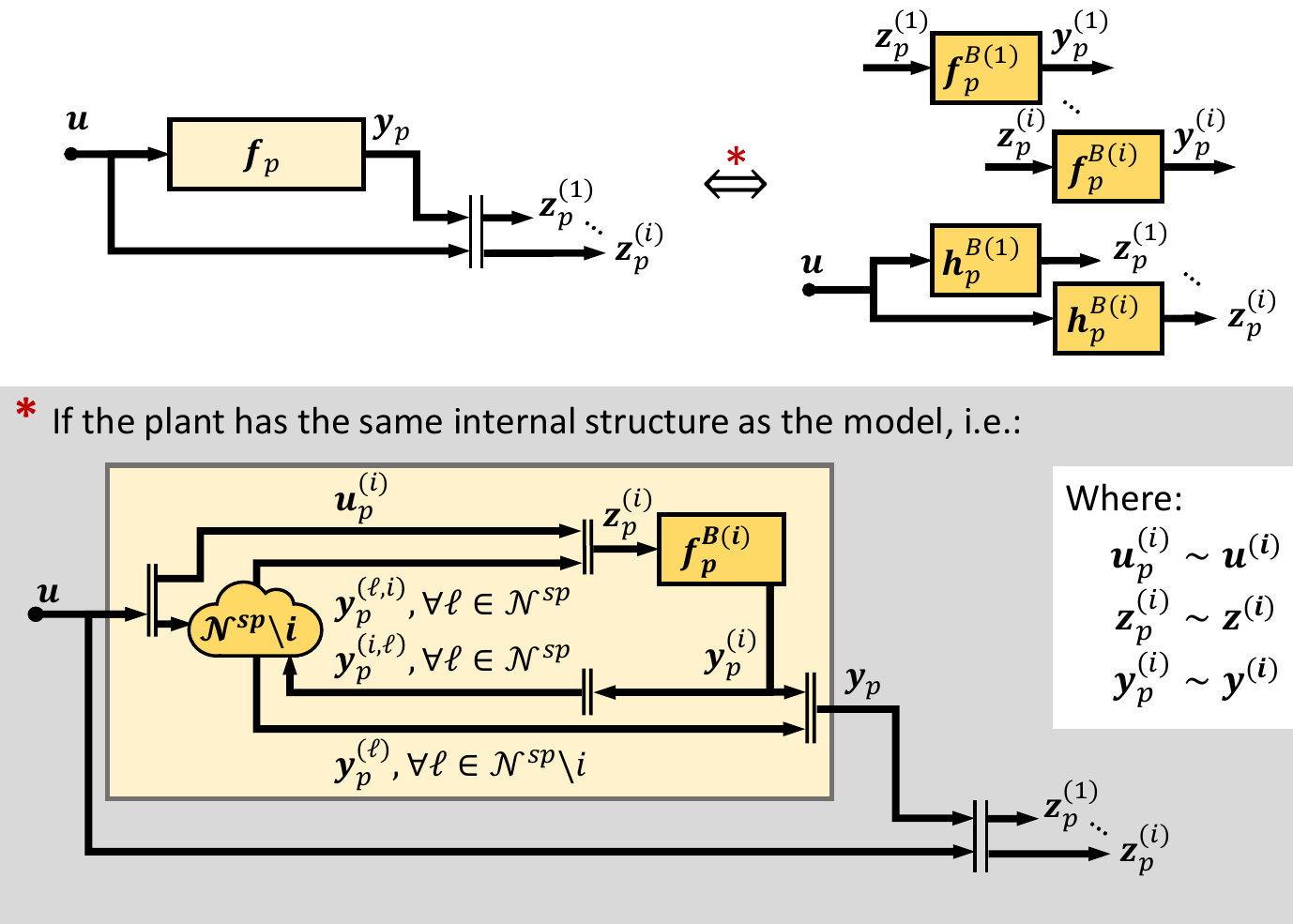}
	\caption{Redefinition of the plant for the structure B}
	\label{fig:5___5_Plant_B}
\end{figure}

Such as for the structure A,  one wants to change the representation of the plant so that:
\vspace{-\topsep}
\begin{itemize}[noitemsep]
	\item Part of its inputs are the inputs of the structures B ($\bm{z}^{(i)}, \forall i$) and part of its outputs are the outputs that the structure B corrects ($\bm{y}^{(i)}, \forall i$).
	\item The other part of its inputs are the plant inputs ($\bm{u}$) and the other part of its outputs are the input of the structures B ($\bm{z}^{(i)}$).
\end{itemize}
To perform this task, it is necessary to assume that the structure of the function $\bm{f}_p$ is the same as the one of the model $\bm{f}$, as shown in Figure~\ref{fig:5___5_Plant_B}.  In other words, one must assume that if an analytical version of the function $\bm{f}_p$ were known, then:
\vspace{-\topsep}
\begin{itemize}[noitemsep]
	\item Applying the three steps of the model preprocessing to this function would yield a network $\bm{\mathcal{N}}^{PSP}$ of $n_{PSP}$ primary subplants (PSP)  that could be aggregated into $n_{SP}=n_{SM}$ subplants (SP). 
	\item Each of these subplants, $i\in\bm{\mathcal{N}}^{SP}$, could then be associated with a submodel of the same index, $i\in\bm{\mathcal{N}}^{SM}$.
	\item The inputs and outputs of a subplant, $\{\bm{z}_p^{(i)},\bm{y}_p^{(i)}\}$, and its associated submodel, $\{\bm{z}^{(i)},\bm{y}^{(i)}\}$, are ``similar'', i.e.    
	\begin{equation}
		\{\bm{z}_p^{(i)},\bm{y}_p^{(i)}\} \sim \{\bm{z}^{(i)},\bm{y}^{(i)}\}.
	\end{equation}
	\textit{(Where two vectors of physical values are said to be similar if they are of the same size and if the values which they compose have the same physical meanings.)}
	\item Functions $\bm{f}_p^{B(i)}$, $ \forall i\in\mathcal{N}^{SP}$, are the equivalent to the functions of the submodels $\bm{f}^{(i)}$, $ \forall i\in\mathcal{N}^{SM}$. 
	\item Functions $\bm{h}_p^{B(i)}$, $ \forall i\in\mathcal{N}^{SP}$ are the result of the action of a mixer-demixer on the inputs $\bm{u}$ and outputs $\bm{y}_p$. It is clear that $\bm{h}_p^{B(i)} = \bm{h}_p^{A(i)}$, $\forall i\in\mathcal{N}^{SM}$. 
\end{itemize}
This assumption is in fact conditional on existence of the functions $\bm{f}_p^{B(i)}$, $ \forall i\in\mathcal{N}^{SP}$, which was introduced in \cite{Papasavvas:2019c} as the model consistency condition: 
\begin{Assumption} \label{ass:5___1_Coherence_Modele}
	\textbf{(Model consistency condition)}
	The internal structure (i.e. the result of the preprocessing) of the model $\bm{f}$ is a network of functions $\bm{f}^{(i)}$, $\forall i \in \mathcal{N}^{SM}$ which correspond to correlations observable on the plant at steady-state and which can be called $\bm{f}_p^{B(i)}$, $\forall i \in \mathcal{N}^{SP}$ with $\mathcal{N}^{SP}\sim \mathcal{N}^{SM}$. Therefore, for each SM$i$ there exist a similar SP$i$: 
	\begin{equation}
	\forall i \in \mathcal{N}^{SM}, \ \ \exists \text{SP}i\sim \text{SM}i.
	\end{equation}
	The similarity of these networks means that they contain the same number of agents (or nodes) and that the interactions that link them are of the same topology (i.e. ``the same nodes are interconnected'') and nature (i.e. ``the same variables are sent through those interconnections'').  
\end{Assumption}
This assumption and its effects are illustrated and discussed in the case studies of sections 
\ref{sec:5_5_3_Model_structurellement_Faux_mais_coherent} and
\ref{sec:5_5_4_Model_structurellement_Faux_et_incoherent}.\textit{(But before moving on to these case studies, let us finish the theoretical analysis of structures A and B.)}

\subsection{Computation of the submodels correction functions}

As for structure A, it is considered that a SM$i$ is ``corrected to the 1st order'' at a point  $\bm{u}_k$ if for the inputs $\bm{z}_{p,k}^{(i)} = \bm{h}_p^{B(i)}(\bm{u}_k)$,
the values and gradients of its corrected outputs are equal to those of the plant, i.e. if:

\begin{align}
	\bm{f}\big|_{\bm{z}_{p,k}^{(i)}} + \bm{\mu}_k^{B(i)}\big|_{\bm{z}_{p,k}^{(i)}}  = \ & \bm{f}_p^{B(i)}\big|_{\bm{z}_{p,k}^{(i)}}, \label{eq:5___4_SM_cor_val_B} \\
	\left(\nabla_{\bm{z}^{(i)}} \bm{f}^{(i)}\big|_{\bm{z}_{p,k}^{(i)}}
	+ 
	\nabla_{\bm{z}^{(i)}} \bm{\mu}_k^{B(i)}\big|_{\bm{z}_{p,k}^{(i)}}
	\right) \nabla_{\bm{u}} \bm{h}_p^{B(i)}\big|_{\bm{u}_k} 
	= \  
	& \nabla_{\bm{z}^{(i)} } \bm{f}_p^{B(i)}\big|_{\bm{z}_{p,k}^{(i)}}\nabla_{\bm{u}} \bm{h}_p^{B(i)}\big|_{\bm{u}_k}. \label{eq:5___5_SM_cor_grad_B}
\end{align} 
For these equalities to be systematically true, the correction function, $\bm{\mu}_k^{B(i)}$, must be chosen such that they correspond to the following Taylor series expansions: 
\begin{align}
	\bm{\mu}_k^{B(i)}(\bm{z}^{(i)}) := \ &\bm{f}_p^{B(i)}\big|_{\bm{z}_{p,k}^{(i)}} - 
	\bm{f}\big|_{\bm{z}_{p,k}^{(i)}} + 
	\left(
	\nabla_{\bm{z}^{(i)}} \bm{f}_p^{B(i)}\big|_{\bm{z}_{p,k}^{(i)}} 
	- \nabla_{\bm{z}^{(i)}} \bm{f}^{(i)}\big|_{\bm{z}_{p,k}^{(i)}}
	\right) \bm{Z}^{(i)}_k
	\big(
	\bm{z}^{(i)} -  ... \nonumber  \\
	& \bm{z}^{(i)}_{p,k}
	\big)+ \mathcal{O}(\|\bm{z}^{(i)} - \bm{z}^{(i)}_{p,k}\|^2).
\end{align}
where the matrix  $\bm{Z}_k^{(i)}$ is a basis of the space in which $\bm{z}^{(i)}$ moves around  $\bm{z}_{p,k}^{(i)}$ when $\bm{u}$ is perturbed in all directions of $\amsmathbb{R}^{n_u}$ around $\bm{u}_k$. In more mathematical terms: 
\begin{equation}
	\bm{Z}_k^{(i)} := 
	\nabla_{\bm{u}} \bm{h}_p^{B(i)}\big|_{\bm{u}_k}  
	\left(
	\nabla_{\bm{u}} \bm{h}_p^{B(i)}\big|_{\bm{u}_k}
	\right)^{+}.
\end{equation}
Again, the simplest choice is that $\bm{\mu}_k^{B(i)}$ is an affine function by setting $\mathcal{O}(\|\bm{z}^{(i)} - \bm{z}^{(i)}_{p,k}\|^2)=\bm{0}$.

\section{Global effects of local corrections}

In order to analyze the global effects of all these local corrections, one proposes to start by analyzing the union of two \textit{systems} $\fe$ and $\g$ that are ``corrected to the 1st order''  at a  point $\bm{u}_k$; specifically where, a \textit{system} is defined as any union of SMs. In the theorem that follows, this union forms a system ``$\fe\cup\g$'' which is correct to the 1st order at  $\bm{u}_k$ is shown.

\begin{thmbox} \label{thm:5___1_Combinaison_Systems_Affienment_Corrects}
	Let a model $\bm{f}$ be of a plant $\bm{f}_p$. Consider two of its systems defined by the functions $[\bm{c}^{\rm T},\bm{v}^{\rm T}]^{\rm T} = \fe(\bm{a},\bm{d})$ and $[\bm{d}^{\rm T},\bm{w}^{\rm T}]^{\rm T} = \g(\bm{c},\bm{b})$, where the variables $\bm{a}$, $\bm{b}$,  $\bm{v}$ and $\bm{w}$ connect them to the rest of the model and the variables $\bm{c}$ and $\bm{d}$ interconnect them. Figure~\ref{fig:5___6_Image_Preuve_THM_5_1} illustrates this notation. The equivalent of these systems in the plant is noted $[\bm{c}_p^{\rm T},\bm{v}_p^{\rm T}]^{\rm T} = \fe_p(\bm{a}_p,\bm{d}_p)$ and $[\bm{d}_p^{\rm T},\bm{w}_p^{\rm T}]^{\rm T} = \g_p(\bm{c}_p,\bm{b}_p)$. 
	
	\begin{minipage}[h]{\linewidth}
		\begin{overpic}[width=14cm]{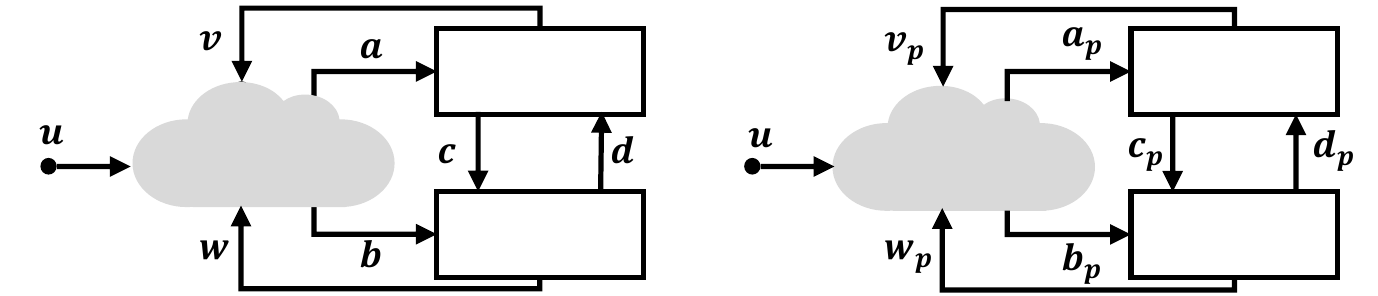}
			\put (38.5,15.5) {$\fe$}
			\put (38.5,3.5) {$\g$}
			\put (88.5,15.5) {$\fe_p$}
			\put (88.5,3.5) {$\g_p$}
			\put (11.5,8.4) {$\bm{\mathcal{N}}^{SM}\backslash\{\fe,\g\}$}
			\put (61.5,8.4) {$\bm{\mathcal{N}}^{SP}\backslash\{\fe_p,\g_p\}$}
		\end{overpic}
		\captionof{figure}{Notations used in this Theorem}
		\label{fig:5___6_Image_Preuve_THM_5_1}
	\end{minipage} \\
	
	If the two systems $\fe$ and $\g$ are correct to the first order at $\bm{u}_k$, i.e: 
	\begin{align}
		\fe(\bm{a},\bm{d}) = \ & \fe_p(\bm{a}_{p,k},\bm{d}_{p,k}) + (\bm{u}-\bm{u}_{k})^{\rm T} \Big(\nabla_a \fe(\bm{a},\bm{d}_{p,k})|_{\bm{a}_{p,k}} \nabla_u \bm{a}|_{\bm{u}_k} +... \nonumber \\
		&  \nabla_d \fe(\bm{a}_{p,k},\bm{d})|_{\bm{d}_{p,k}} \nabla_u \bm{d}|_{\bm{u}_{p,k}} \Big) + 
		\mathcal{O}(\| \bm{u}-\bm{u}_k\|^2),  \label{eq:LocGlo_1} \\
		\g(\bm{c},\bm{b}) = \ & \g_p(\bm{c}_{p,k},\bm{b}_{p,k}) + (\bm{u}-\bm{u}_{k})^{\rm T} \Big(\nabla_c \g(\bm{c},\bm{d}_{b,k})|_{\bm{c}_{p,k}} \nabla_u \bm{c}|_{\bm{u}_k} +... \nonumber \\
		&  \nabla_b \g(\bm{c}_{p,k},\bm{b})|_{\bm{b}_{p,k}} \nabla_u \bm{b}|_{\bm{u}_{p,k}} \Big) +  \mathcal{O}(\| \bm{u}-\bm{u}_k\|^2),  \label{eq:LocGlo_2}
	\end{align}
	where: 
	\begin{align}
		\nabla_a \fe(\bm{a},\bm{d}_{p,k})|_{\bm{a}_{p,k}} 
		\nabla_u \bm{a}_p|_{\bm{u}_k} \hspace{-1mm}
		+ 
		\nabla_d \fe(\bm{a}_{p,k},\bm{d})|_{\bm{d}_{p,k}}
		\nabla_u \bm{d}_p|_{\bm{u}_k}
		= \ &
		\nabla_u \fe_p|_{\bm{u}_{k}} , \label{eq:LocGlo_3}  \\
		\nabla_c \g(\bm{c},\bm{b}_{p,k})|_{\bm{c}_{p,k}} \hspace{-1mm}
		\nabla_u \bm{c}_p|_{\bm{u}_k}+
		\nabla_b \g(\bm{c}_{p,k},\bm{b})|_{\bm{b}_{p,k}}
		\nabla_u \bm{b}_p|_{\bm{u}_k}
		= \  &
		\nabla_u \g_p|_{\bm{u}_{k}}.  \label{eq:LocGlo_4}
	\end{align}
	Then, the union of these two systems $\fe \cup \g$ is also correct to the first order at $\bm{u}_k$.
\end{thmbox}
\begin{proofbox}
	To determine whether the union $\fe \cup \g$ is correct to the first order, one  evaluates them at the following  inputs:
	\begin{align}
		\bm{a} = \ & \bm{a}_{p,k} + (\bm{u}-\bm{u}_k)^{\rm T} \nabla_u \bm{a}_p|_{\bm{u}_k} + \mathcal{O}(\|\bm{u}-\bm{u}_k\|^2), \label{eq:LocGlo_5}\\ 
		\bm{b} = \ & \bm{b}_{p,k} + (\bm{u}-\bm{u}_k)^{\rm T} \nabla_u \bm{b}_p|_{\bm{u}_k} + \mathcal{O}(\|\bm{u}-\bm{u}_k\|^2). \label{eq:LocGlo_6}
	\end{align}
	If for these inputs one finds that $\bm{c}$, $\bm{d}$, $\bm{v}$, and $\bm{w}$ follow
	\begin{align}
		\bm{c} = \ & \bm{c}_{p,k} + (\bm{u}-\bm{u}_k)^{\rm T} \nabla_u \bm{c}_p|_{\bm{u}_k} + \mathcal{O}(\|\bm{u}-\bm{u}_k\|^2), \label{eq:LocGlo_7}\\ 
		\bm{d} = \ & \bm{d}_{p,k} + (\bm{u}-\bm{u}_k)^{\rm T} \nabla_u \bm{d}_p|_{\bm{u}_k} + \mathcal{O}(\|\bm{u}-\bm{u}_k\|^2), \label{eq:LocGlo_8}\\ 
		\bm{v} = \ & \bm{v}_{p,k} + (\bm{u}-\bm{u}_k)^{\rm T} \nabla_u \bm{v}_p|_{\bm{u}_k} + \mathcal{O}(\|\bm{u}-\bm{u}_k\|^2), \label{eq:LocGlo_9}\\ 
		\bm{w} = \ & \bm{w}_{p,k} + (\bm{u}-\bm{u}_k)^{\rm T} \nabla_u \bm{w}_p|_{\bm{u}_k} + \mathcal{O}(\|\bm{u}-\bm{u}_k\|^2), \label{eq:LocGlo_10}
	\end{align}
	then one will conclude that the union $\fe \cup \g$ is correct to the first order.
	
	To analyze how the two systems interact, one introduces the reduced functions $\bm{c} =\fe_c(\bm{a},\bm{d})$ and $\bm{d} =\g_d(\bm{b},\bm{c})$ which are the functions $\fe$ and $\g$ from which the outputs  $\bm{v}$ and $\bm{w}$ have been removed. Interconnecting the systems $\fe_c$ and $\g_d$ gives: 
	\begin{equation} \label{eq:LocGlo_11}
		\bm{d} = \g_d(\fe_c(\bm{a},\bm{d}),\bm{b}). 
	\end{equation}
	The solution of this system for $\bm{u}=\bm{u}_k$ (i.e. $\bm{a}=\bm{a}_{p,k}$ and  $\bm{b}=\bm{b}_{p,k}$ because of \eqref{eq:LocGlo_5} and \eqref{eq:LocGlo_6}) is obvious because: 
	\begin{align}
		\text{\eqref{eq:LocGlo_1}} \Rightarrow	\quad   &\fe_c(\bm{a}_{p,k},\bm{d}_{p,k}) = \fe_{p,c}(\bm{a}_{p,k},\bm{d}_{p,k}) = \bm{c}_{p,k}, \label{eq:LocGlo_12} \\
		\text{\eqref{eq:LocGlo_2}} \Rightarrow	\quad   &\g_d(\bm{c}_{p,k},\bm{b}_{p,k}) = \g_{p,d}(\bm{c}_{p,k},\bm{b}_{p,k}) = \bm{d}_{p,k}, \label{eq:LocGlo_13} \\
		\text{Def.} \Rightarrow	\quad   &  \bm{d}_{p,k} =\g_{p,d}(\fe_{c,p}(\bm{a}_{p,k},\bm{d}_{p,k}),\bm{b}_{p,k}), \label{eq:LocGlo_14}
	\end{align}
	where the last equality is a definition of the plant. Combining \eqref{eq:LocGlo_12}-\eqref{eq:LocGlo_14} gives: 
	\begin{equation} \label{eq:LocGlo_15}
		\bm{d}_{p,k} =\g_d(\fe_c(\bm{a}_{p,k},\bm{d}_{p,k}),\bm{b}_{p,k}).
	\end{equation}
	In other words $\bm{d}_{p,k}$ is a solution of \eqref{eq:LocGlo_11} for $\bm{u} = \bm{u}_k$. Let us now analyze how this solution varies when $\bm{u}$ varies around $\bm{u}_k$, i.e. let's compute $\nabla_u \bm{d}|_{\bm{u}_k}$.
	\begin{align}
		\nabla_u \bm{d}|_{\bm{u}_k} = \ & \nabla_u \g_d(\fe_c(\bm{a},\bm{d}),\bm{b})|_{\bm{u}_k} \nonumber \\
		= \ & 
		\Big(
		\nabla_a \fe_c(\bm{a},\bm{d}_{p,k})|_{\bm{a}_{p,k}} \nabla_u \bm{a}|_{\bm{u}_k} + 
		\nabla_a \fe_c(\bm{a}_{p,k},\bm{d})|_{\bm{d}_{p,k}}  \nabla_u \bm{d}|_{\bm{u}_k}
		\Big) ... \nonumber \\
		& \nabla_c \g_d(\bm{c},\bm{b}_{p,k})|_{\bm{c}_{p,k}} 
		+
		\nabla_b \g_d(\bm{c}_{p,k},\bm{b})|_{\bm{b}_{p,k}} \nabla_u \bm{b}|_{\bm{u}_k}. \label{eq:LocGlo_16}
	\end{align}
	These functions are evaluated at $\{\bm{a}_{p,k},\bm{b}_{p,k},\bm{d}_{p,k}\}$ because \eqref{eq:LocGlo_5}, \eqref{eq:LocGlo_6} and \eqref{eq:LocGlo_15}. One can observe that: 
	\begin{align}
		\text{\eqref{eq:LocGlo_5}} \Rightarrow \  & \nabla_u \bm{a}|_{\bm{u}_k} = \nabla_u \bm{a}_p|_{\bm{u}_k},   \label{eq:LocGlo_17} \\ 
		\text{\eqref{eq:LocGlo_6}} \Rightarrow \  & \nabla_u \bm{b}|_{\bm{u}_k} = \nabla_u \bm{b}_p|_{\bm{u}_k}, \label{eq:LocGlo_18} \\ 
		\left.\begin{array}{c}
		\text{\eqref{eq:LocGlo_3}} \\ \text{\eqref{eq:LocGlo_17}} \\ \text{\eqref{eq:LocGlo_18}}
		\end{array}\right\} \Rightarrow \  &  
		\begin{array}{r}
			\nabla_a \fe_c(\bm{a},\bm{d}_{p,k})|_{\bm{a}_{p,k}} \nabla_u \bm{a}|_{\bm{u}_k} + 
			\nabla_a \fe_c(\bm{a}_{p,k},\bm{d})|_{\bm{d}_{p,k}}  ... \\  
			\nabla_u \bm{d}|_{\bm{u}_k}
			= 
			\nabla_u \fe_{p,c}|_{\bm{u}_{k}},
		\end{array}
	 	\label{eq:LocGlo_19}
		\\
		\left.\begin{array}{c}
			\text{Def.} \\ 
			\text{\eqref{eq:LocGlo_19}}
		\end{array}\right\}
		\Rightarrow \  &\nabla_u \bm{c}|_{\bm{u}_k} = \nabla_u \fe_{p,c}|_{\bm{u}_k}. \label{eq:LocGlo_20} 
	\end{align}
	From these observations, one can reformulate \eqref{eq:LocGlo_16}:
	\begin{align}
		& & \text{\eqref{eq:LocGlo_16}} = \ & 
		\nabla_u \fe_{p,c}|_{\bm{u}_{k}}
		\nabla_c \g_d(\bm{c},\bm{b}_{p,k})|_{\bm{c}_{p,k}}  + 
		\nabla_b \g_d(\bm{c}_{p,k},\bm{b})|_{\bm{b}_{p,k}}
		\nabla_u \bm{b}|_{\bm{u}_k}, \nonumber\\
		\left.\begin{array}{c}
			\text{\eqref{eq:LocGlo_4}} \\
			\text{\eqref{eq:LocGlo_18}} 
		\end{array}\right\}
		\Rightarrow & & = \ &
		\left(
		\nabla_u \fe_{p,c}|_{\bm{u}_{k}} -
		\nabla_u \bm{c}|_{\bm{u}_k}
		\right) 
		\nabla_c \g_d(\bm{c},\bm{b}_{p,k})|_{\bm{c}_{p,k}}|_{\bm{u}_k} + 
		\nabla_u \g_{p,d}|_{\bm{u}_{k}},\nonumber  \\
		\eqref{eq:LocGlo_20} \Rightarrow & & = \ & \nabla_u \g_{p,d}|_{\bm{u}_{k}} = \nabla_u \bm{d}_p|_{\bm{u}_k}, \nonumber \\
		\Leftrightarrow & & \nabla_u\bm{d}|_{\bm{u}_k} = \ & \nabla_u\bm{d}_p|_{\bm{u}_k}. \label{eq:LocGlo_21}
	\end{align}
	Equalities \eqref{eq:LocGlo_15} and \eqref{eq:LocGlo_21} are sufficient to show that the union of systems that are correct to the first order at $\bm{u}_k$ is a system that is also correct to the first order at $\bm{u}_k$.  Indeed, it has been shown that under such conditions, $\bm{d}|_{\bm{u}_k} =\bm{d}_{p,k}$ and $\nabla_u\bm{d}|_{\bm{u}_k} =  \nabla_u\bm{d}_p|_{\bm{u}_k}$. So, the inputs to the system $\fe$ are correct to the first order and so are its outputs  $\bm{c}$ and $\bm{v}$.  As a result, the inputs $\bm{b}$ and $\bm{c}$ of  $\g$  are correct to the first order and so the outputs  $\bm{d}$ and $\bm{w}$, which  concludes this proof.
\end{proofbox}

On the basis of this theorem one can show that:

\begin{thmbox}
	Consider a model that is split into SMs. Assume that the corrections are applied such that all SMs are corrected to the first order at $\bm{u}_k$. Then, the updated model is correct to the first order at $\bm{u}_k$.
\end{thmbox}
\begin{proofbox}
	This theorem proven by induction. 
	
	\textit{Base case:}  On one hand, let's take the system $\fe^{(1)}$, composed of the submodel  $1$, and on the other hand the rest of the model $\bm{N}^{SM}\backslash\fe^{(1)}$, which is denoted $\g$ to keep analogous notations with Theorem~\ref{thm:5___1_Combinaison_Systems_Affienment_Corrects}. Let's apply corrections to these two systems so that they are correct to the first order at $\bm{u}_k$. Then, equalities \eqref{eq:LocGlo_1}-\eqref{eq:LocGlo_4} are true because:
	\begin{itemize}
		\item by definition the two systems are corrected to the first order and equalities  \eqref{eq:LocGlo_1} and \eqref{eq:LocGlo_2}) hold; 
		\item in this particular configuration $\bm{u}=\bm{a}=\bm{b}=\bm{a_p}=\bm{b}_p$, so \eqref{eq:LocGlo_3} and \eqref{eq:LocGlo_4} are necessarily true.  
	\end{itemize}
	One can therefore apply Theorem~\ref{thm:5___1_Combinaison_Systems_Affienment_Corrects} to show that the union of these systems, that are corrected separately,  is correct to the first order, and as this union is the model, it is itself correct to the first order. 
	
	\textit{Induction step:} Given that both $\fe^{(1)}$, ..., $\fe^{(\ell)}$, and $\g:=\bm{N}^{SM}\backslash\{\fe^{(1)}, ..., \fe^{(\ell)}\}$ are all correct to the first order at $\bm{u}_k$, and whose union forms a model that is also correct to the first order at $\bm{u}_k$. There are several ways to ensure that  $\g$ is correct to the first order at $\bm{u}_k$. One way to achieve this goal is to separate  $\g$ into two systems, one of which would be the SM$\ell+1$ denoted $\fe^{(\ell+1)}$  and the other would be $\g^{\prime}:=\bm{\mathcal{N}}^{SM}\backslash\{\fe^{(1)}, ..., \fe^{(\ell+1)}\}$, and to apply corrections to them individually so that they are both correct to the first order at $\bm{u}_k$. Indeed, Theorem ~\ref{thm:5___1_Combinaison_Systems_Affienment_Corrects}  proves that any union of  systems that are correct to the first order leads  a system correct to the first order.  
	
	Finally, by starting from the basic case and by repeating the induction step until the system $\g$  is emptied of all submodels, one proves this Theorem.
\end{proofbox}
With this theorem, one has shown that no matter what correction structure one chooses to correct a SM, as long as all SMs are corrected\footnote{In order to simplify the explanation, it is not always specified that it is a correction \textit{to the first order at the correction point}.}, the model is also corrected. It is worth noting that it is not necessary to use the same structure (A or B) for all the SMs (which makes this theorem more general than the one proposed in \cite{Papasavvas:2019c}). 

Also, one would like to point out that no matter what correction structure is used, all the theorems developed for KMFCaA remain applicable as long as this structure implies that ultimately the model is corrected. Only structure B may not validate this condition when the model is not consistent. 
So, these correction structures only affect KMFCaA at its step 2). So, the way of choosing the filter introduced in Chapter~\ref{Chap:2_Vers_Une_meilleure_Convergence}, as well as the way of implementing the constraints introduced in the Chapter~\ref{Chap:3_KMA} are not affected by these choices.
Afterwards, one says that one uses the KMFCaA-A (or KMFCaA-B) method if one implements KMFCaA with the correction structure A (or B) on \textit{all} the SMs.
Of course, it is perfectly possible to use structure A on a part of the SMs and structure B on the rest. In such a case one would say that one uses a KMFCaA-AB method (but in the framework of this thesis, this option is not analyzed -- despite its potential).  

At this stage, it has been shown that structures other than the commonly used one exist and can be used to obtain convergence properties similar to those of KMFCaA-D/I. Now the question is: \textit{``Which structure should be used?''}

Giving an analytical answer to this question is far from obvious.  If one has compared KMFCaA-D/I through analyses using Taylor series in \cite{Papasavvas:2019a}, such an approach was possible  only because the model was ``modeled'' with a single equation. In the situation of the structures A and B, the model can be any network of more or less non-linear systems, which considerably complicates (or even makes impossible) an analysis of the same nature.  Since no satisfactory analytical approach has been found, it is decided to use empirical studies to give an insight to the answer to this question. 

\section{Empirical analyses of structures D, I, A and B}

In this section, one compares the correction structures D, I, A and B as well as the associated RTO algorithms KMFCaA-D, -I, -A, and -B.

Proceeding in four studies, the impact of the choice of the correction structure is studied on
\begin{itemize}
	\item (Study 1) systems connected in open loop with a consistent model. 
	\item (Study 2) systems connected in a closed loop with a consistent model.
	\item (Study 3) systems connected in open loop with a consistent model \textit{despite a ``significant'' structural error}. 
	\item (Study 4) systems connected in series with an inconsistent model,
\end{itemize}
which represent almost all configurations that can be encountered in practice. Hence the interest in analyzing each of them.  

\subsection{Study 1: Systems in open loop \& Consistent model}

\begin{minipage}[h]{\linewidth}
	\includegraphics[width=14cm]{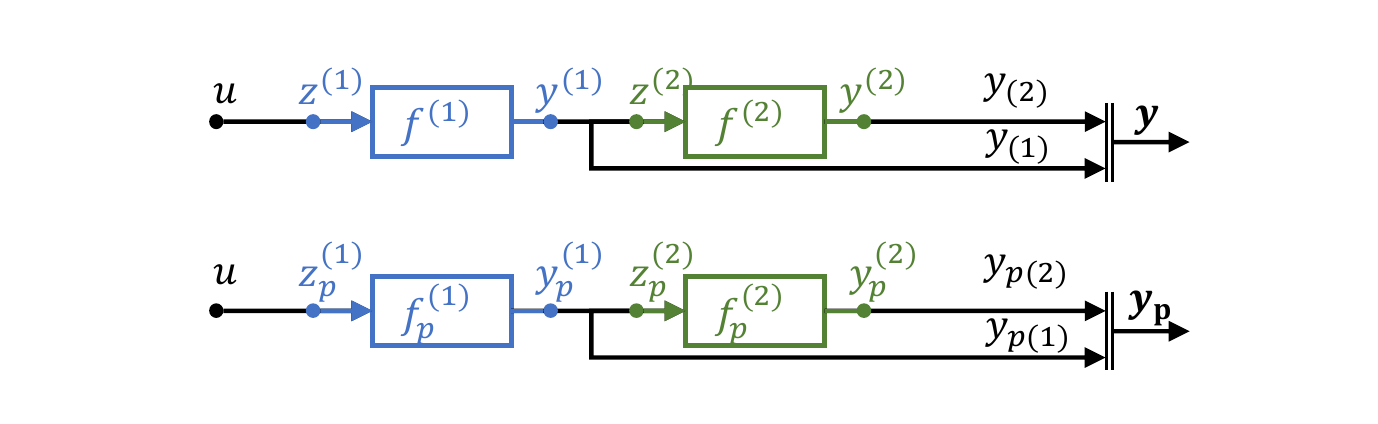}
	\captionof{figure}{Graphical description of the RTO problem of Study 1}
	\label{fig:5___8_Systemes_en_serie}
\end{minipage}\\

One considers the theoretical RTO problem illustrated in Figure~\ref{fig:5___8_Systemes_en_serie} and using the following functions: 
	\begin{align*}
		f^{(1)} := \ &\theta_1 + \theta_2u^{(1)} + \theta_3u^{(1)2}, & f_p^{(1)} := \ & 2 -2u_p^{(1)} + \theta_{p1}u_p^{(1)2},   \\
		 u^{(1)} := \ & u  ,& u_p^{(1)} := \ & u , \\
		f^{(2)} := \ & \theta_4 + \theta_5u^{(2)} + \theta_6u^{(2)2}, & f_p^{(2)} := \ & 2 -2u_p^{(2)} + \theta_{p2}u_p^{(2)2}, \\
		u^{(2)} :=\ & y^{(1)}, & u_p^{(2)} := \ & y_p^{(1)}, \\
		\bm{y}  :=\ & [y^{(1)};y^{(2)}]^{\rm T}
	\end{align*}
	\begin{align*}
		\phi_1(u,\bm{y}) := \ & y_{(2)},   & \phi_2(u,\bm{y}) := \ & uy_{(2)}, \\
		\phi_3(u,\bm{y}) := \ & y_{(2)}^2, & \phi_4(u,\bm{y}) := \ & 1.5u^2 + 0.05uy_{(2)} + 0.005y_{(2)}^2,
	\end{align*}
	where $u\in[-1,1]$. Note that the three cost functions $\{\phi_1,\phi_2,\phi_3\}$ are the three most simple non-linear terms (w.r.t.  $u$) that a cost function (or constraint function) can have; hence, the interest to observe them separately. The function $\phi_4$ is just a simple combination that could be a realistic cost function. 
	The model is parameterized with the parameters $\{\theta_1,...,\theta_6\}$ and the plant with the parameters $\{\theta_{p1},\theta_{p2}\}$. In this example, one considers 3 scenarios:
	\begin{itemize}
		\item \textit{Scenario 1 (A$\rightarrow$NL):}  The functions $\{f^{(1)},f_p^{(1)}\}$ are affine (A),  and the functions $\{f^{(2)},f_p^{(2)}\}$ are non-linear (NL).
		\item \textit{Scenario 2 (NL$\rightarrow$A):} The functions $\{f^{(1)},f_p^{(1)}\}$ are NL,  and the functions $\{f^{(2)},f_p^{(2)}\}$ are A.
		\item \textit{Scenario 3 (NL$\rightarrow$NL):} The functions $\{f^{(1)},f_p^{(1)},f^{(2)},f_p^{(2)}\}$ are NL.
	\end{itemize}
The model and plant parameters considered for each of these scenarios are given in Table~\ref{tab:5_1_Exemple_5_1_Parametres_Modeles_Usine}.  The Figures~\ref{fig:5_9_Exemple_5_1_sc1_Results_1}, \ref{fig:5_12_Exemple_5_1_sc2_Results} and \ref{fig:5_15_Exemple_5_1_Results} show what are the functions of the plant and the models that are used in these studies.   

\noindent
	\begin{minipage}[h]{\linewidth}
		\centering
		\begin{tabular}{cccccccccc}
			\toprule
			Scenario  & $\theta_1$        & $\theta_2$          & $\theta_3$        & $\theta_4$        & $\theta_5$          & $\theta_6$        & $\theta_{p1}$ & $\theta_{p2}$ & Number of models \\
			\midrule
			1         & $\mathbb{\Theta}$ & $\mathbb{\Theta}^-$ & $0$               & $\mathbb{\Theta}$ & $\mathbb{\Theta}^-$ & $\mathbb{\Theta}$ & 0             & 2             & 7776\\
			2         & $\mathbb{\Theta}$ & $\mathbb{\Theta}^-$ & $\mathbb{\Theta}$ & $\mathbb{\Theta}$ & $\mathbb{\Theta}^-$ & $0 $              & 2             & 0             & 7776\\
			3         & $\mathbb{\Theta}$ & $\mathbb{\Theta}^-$ & $\mathbb{\Theta}$ & $\mathbb{\Theta}$ & $\mathbb{\Theta}^-$ & $\mathbb{\Theta}$ & 2             & 2             & 46656 \\
			\bottomrule
		\end{tabular}%
		\captionof{table}{Definitions of the models and of the plant used for each scenario. One defines: $\mathbb{\Theta}:=\{1,1.4,1.8,2.2,2.6,3\}$ and  $\mathbb{\Theta}^-:=\{-3,-2.6,-2.2,-1.8,$ $-1.4,-1\}$.}
		\label{tab:5_1_Exemple_5_1_Parametres_Modeles_Usine}%
	\end{minipage} \\

	These three scenarios cover all the situations that can occur when one has a SM$i$ which is connected to all the other SMs ($\mathcal{N}\backslash i$) \textit{without feedback}. Indeed, for each of these cases one can consider that $f^{(1)}$ (or $f^{(2)}$) is the SM$1$ and $f^{(2)}$ (or $f^{(1)}$) is  $\mathcal{N}\backslash i$. By analyzing the results one gets for each of these scenarios with the correction structures D, I, A and B, empirical rules on the choice of the correction structure are deduced.  
	
	One starts by quantifying the ``local quality'' of the correction associated to each structure. One already knows that, by construction, they all bring affine corrections to the predictions of the costs and constraints.  So, what can be used to differentiate them is the higher order corrections such as the prediction error on the Hessian of the cost function at the correction point. To this end, one has evaluated the distribution of these errors for each model, correction structure, and cost functions at 11 points uniformly distributed in the input space, i.e. $u=\{-1,-0.8,...,0.8,1\}$. Then, by applying a statistical analysis to the obtained results one gets the graphs of Figures~\ref{fig:5_9_Exemple_5_1_sc1_Results_1}, \ref{fig:5_12_Exemple_5_1_sc2_Results}, and \ref{fig:5_15_Exemple_5_1_Results}. In addition to these statistical analyses, Figures~\ref{fig:5_10_Exemple_5_1_sc1_Results_2}, \ref{fig:5_13_Exemple_5_1_sc2_Results_2}, and \ref{fig:5_16_Exemple_5_1_Results_2} provide, for each structure, the percentage of cases for which no other structure brings better predictions of the Hessian.
	
	What can be easily concluded from all these results is that if SM$i$ is in series with the rest of the SMs and the model is consistent then one can rank the correction structures from the most to the least efficient as follows:  
	\begin{align}
		& B>A>I>D, & \text{if the systems are in series  and  $\forall i$ $\exists$ SP$i$ $\sim$ SM$i$.} 
	\end{align}
	In the case of the scenario 1, the structures A and B provide the same results. This is logical since in the cases, the SM1 is perfectly corrected, and therefore correcting the SM2 from  $u^{(1)}$ or $u^{(2)}$ is the same.  
	In the case of the scenario 2, it is the structures A and I  that give the same results for similar reasons.  
	These two observations can easily be verified with a Taylor series development of the corrected models. Finally, in all the other cases the results of structure A are between those of structures B and I, as it can be seen in Figure~\ref{fig:5_15_Exemple_5_1_Results}. This explains why A seems to underperform compared to I or B (as it is usually between the two of them, it is usually outperformed by one of them), and this is why structure A may appear to underperform on Figure~\ref{fig:5_16_Exemple_5_1_Results_2}.
	
	One might ask what is the point of high order corrections in the context of RTO. In an attempt to answer to  this question, one runs simulations starting at 11 points uniformly distributed over the entire input space, i.e. $u_0 =\{-1,-0.8,...,0.8,1\}$, and using all the models of each scenario. The results of these simulations are given in Figure~\ref{fig:5_17_Exemple_5_1_Results_3}.  If one compares these results with the analysis previously made on the quality of correction, it appears quite clearly that the methods which provide the best high-order corrections are those that converge \textit{generally} the fastest on the optimum of the plant and this independent of the model, the cost function, and the starting point $u_0$. \textit{Thus, the interest of high order corrections seems to be to increase the speed of convergence}. If one cannot demonstrate this statement explicitly, one considers that this statistical result is sufficient to support this statement.
	
Finally, one can see that structures I, A, and B almost systematically perform significantly better than structure D, which is, nevertheless, the most used/mentioned structure in the research articles.
	
	\clearpage

\noindent
	\begin{minipage}[h]{\linewidth}
	\begin{center}
		\textbf{Scenario 1: A$\rightarrow$NL}
	\end{center}
	
	\vspace{-3mm}
	\begin{minipage}[h]{\linewidth}
		\vspace*{0pt}
		{\centering
			\begin{minipage}[t]{4.45cm}{\centering%
					\includegraphics[width=4.45cm]{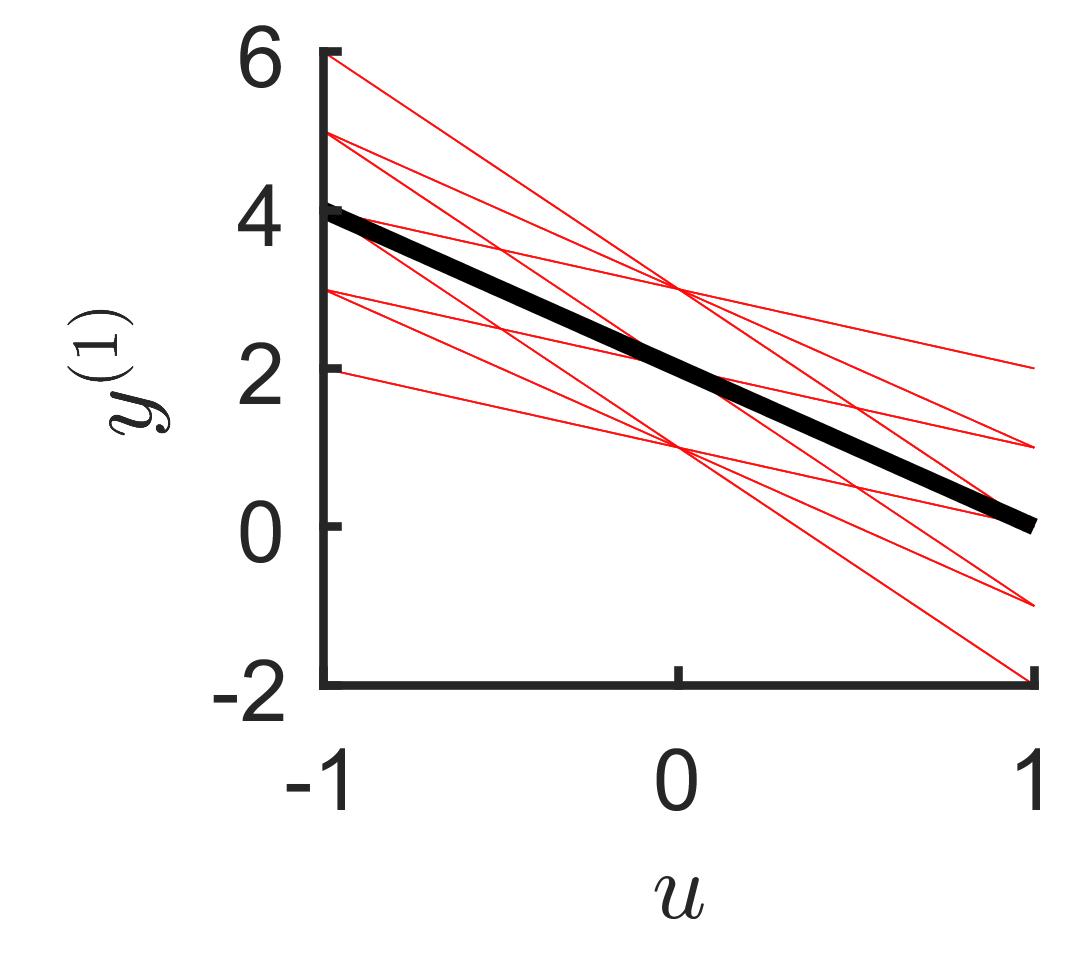}\\
					a) Functions $f^{(1)}$ and $f_p^{(1)}$}
			\end{minipage}\hskip -0ex
			\begin{minipage}[t]{5cm}\centering%
				\includegraphics[width=4.45cm]{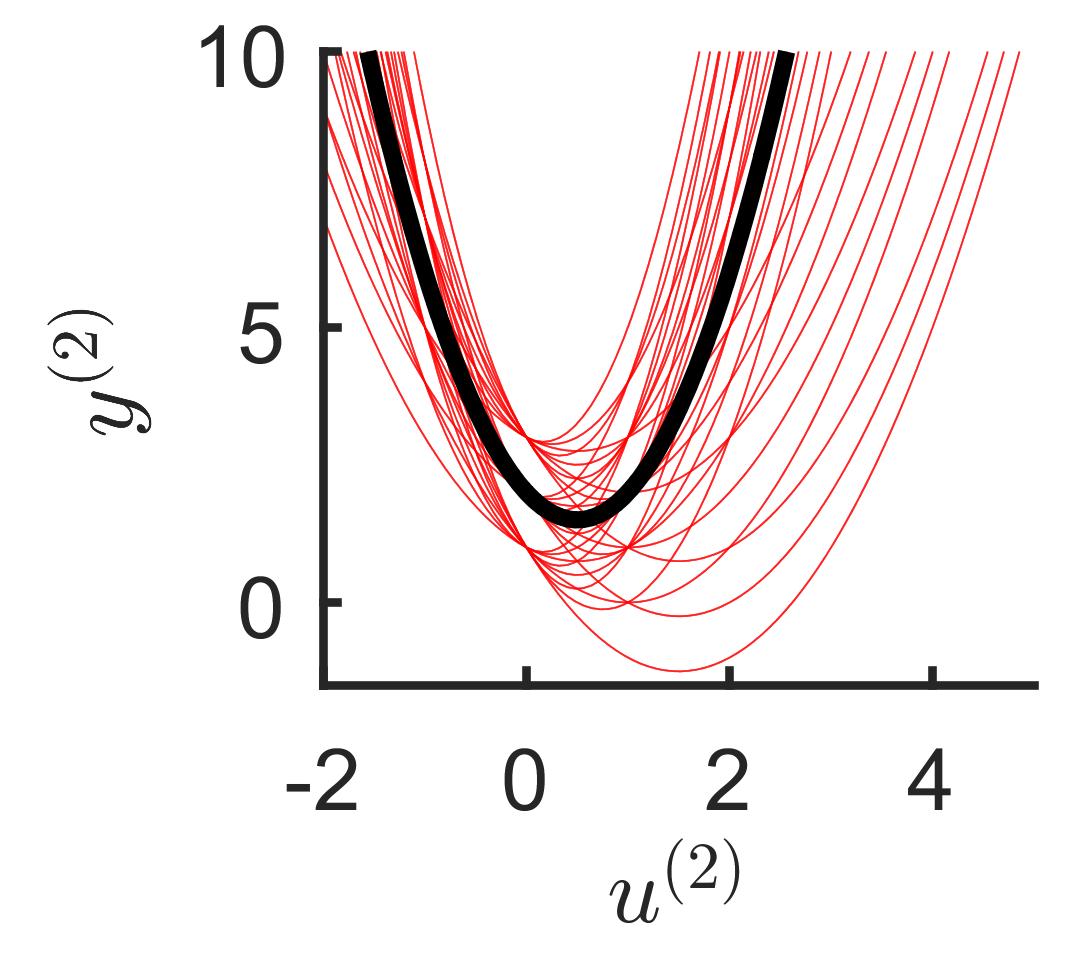}\\
				b) Functions $f^{(2)}$ and $f_p^{(2)}$
			\end{minipage}\hskip -0ex
			\begin{minipage}[t]{4.45cm}\centering%
				\includegraphics[width=4.45cm]{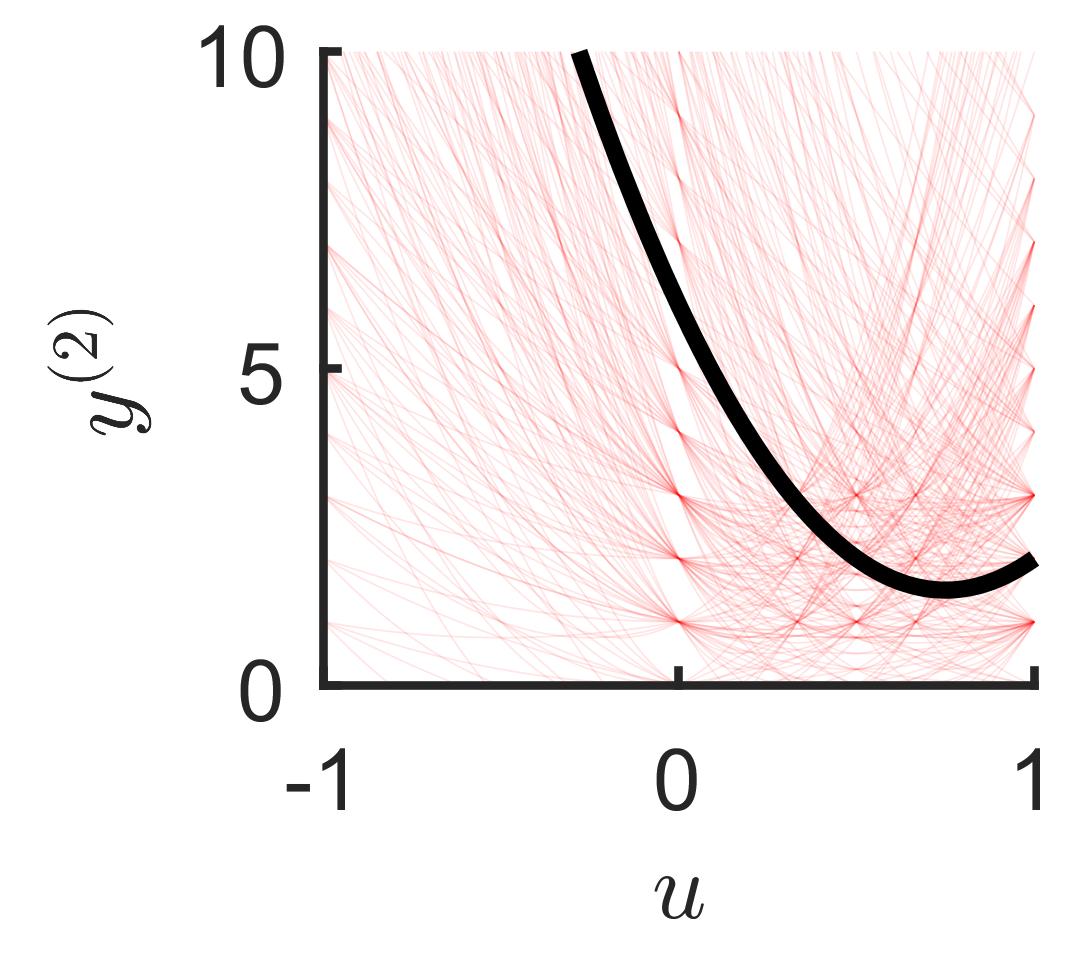}\\
				c) Functions $f$ and $f_p$
			\end{minipage} \\
			
			\medskip
			
			\begin{minipage}[h]{\linewidth} \centering
				\includegraphics[trim={0.1cm 0.2cm 0.2cm  0.1cm },clip,width=3.5cm]{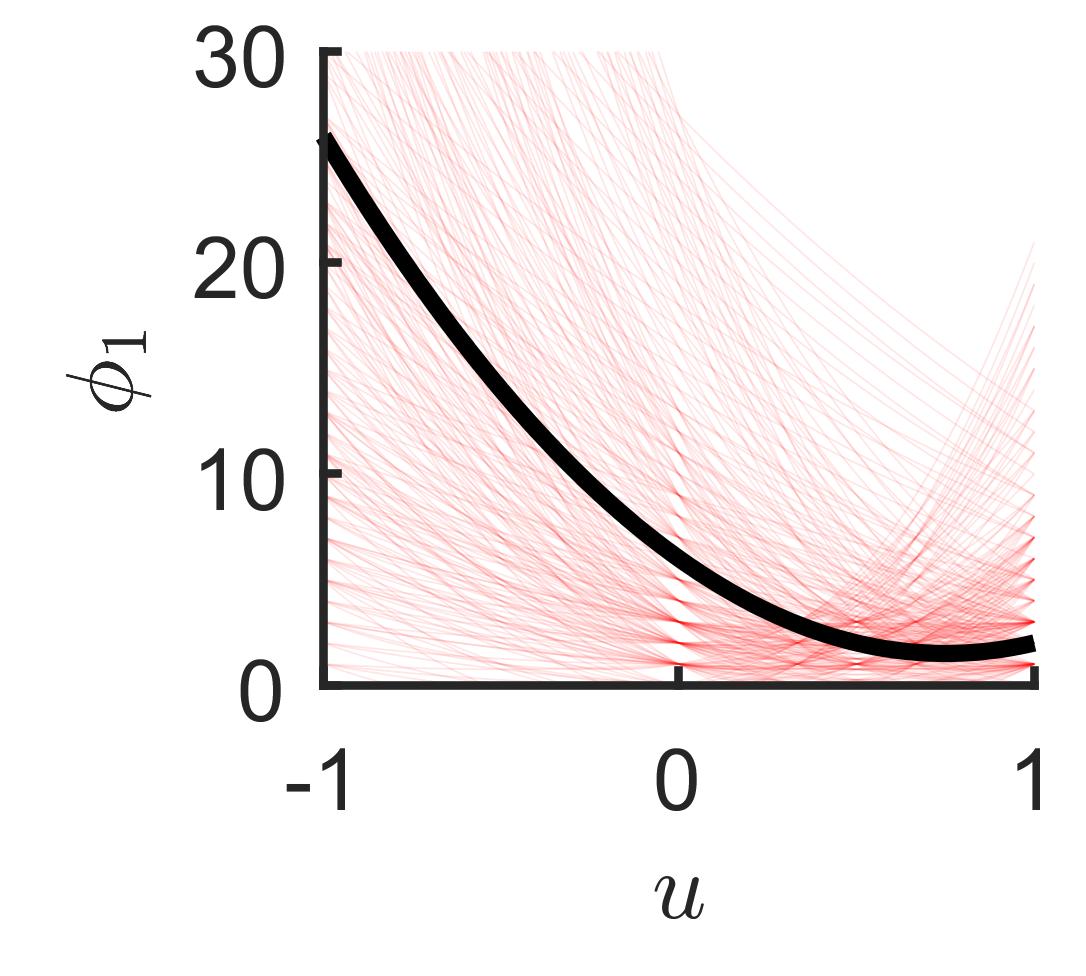}\hskip -0ex
				\includegraphics[trim={0.1cm 0.2cm 0.2cm  0.1cm },clip,width=3.5cm]{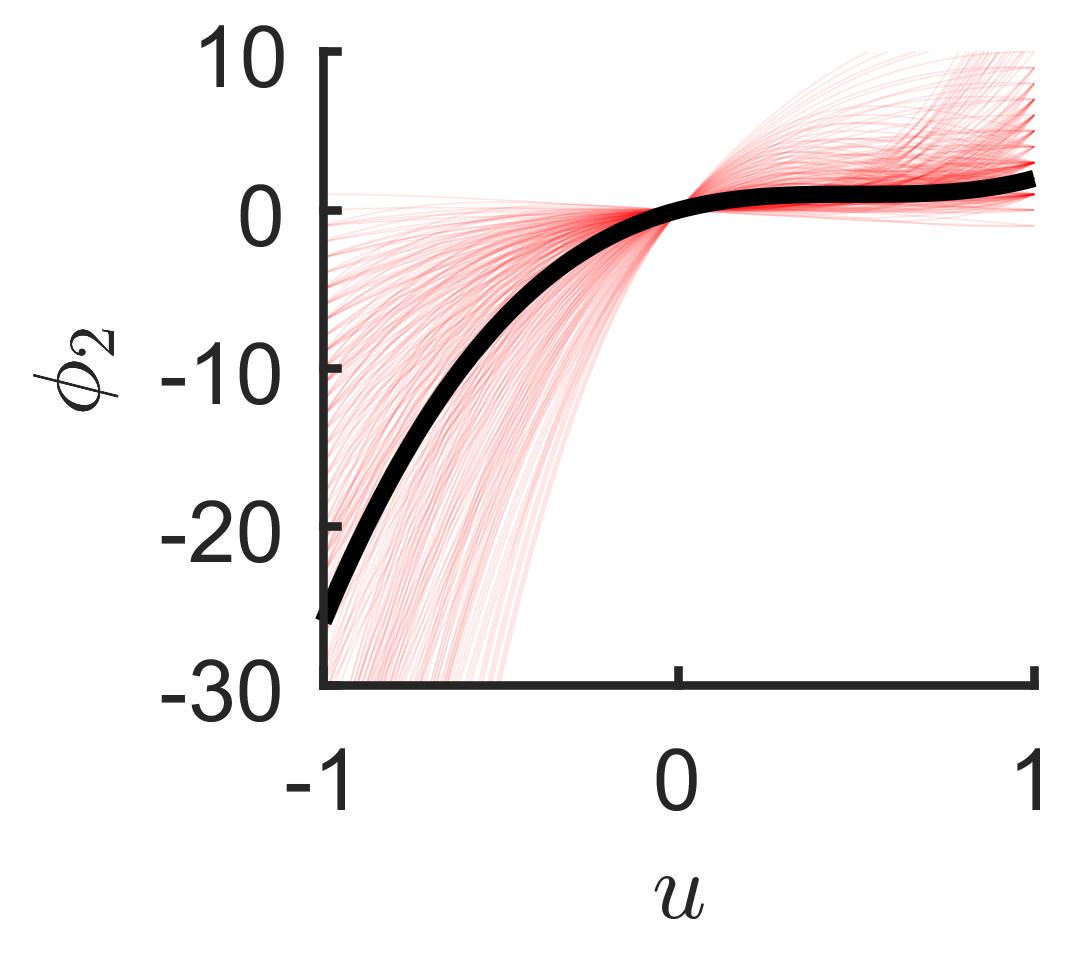}\hskip -0ex
				\includegraphics[trim={0.1cm 0.2cm 0.2cm  0.1cm },clip,width=3.5cm]{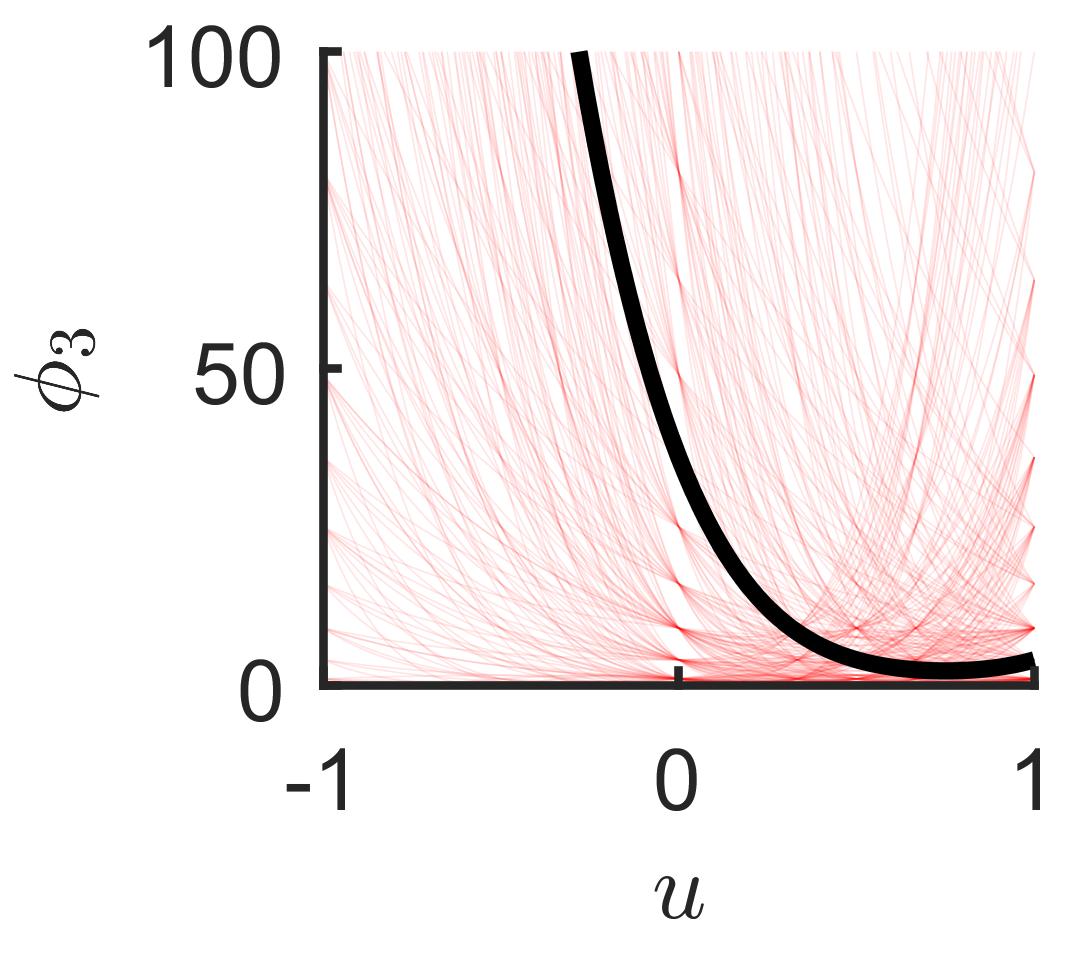}\hskip -0ex
				\includegraphics[trim={0.1cm 0.2cm 0.2cm  0.1cm },clip,width=3.5cm]{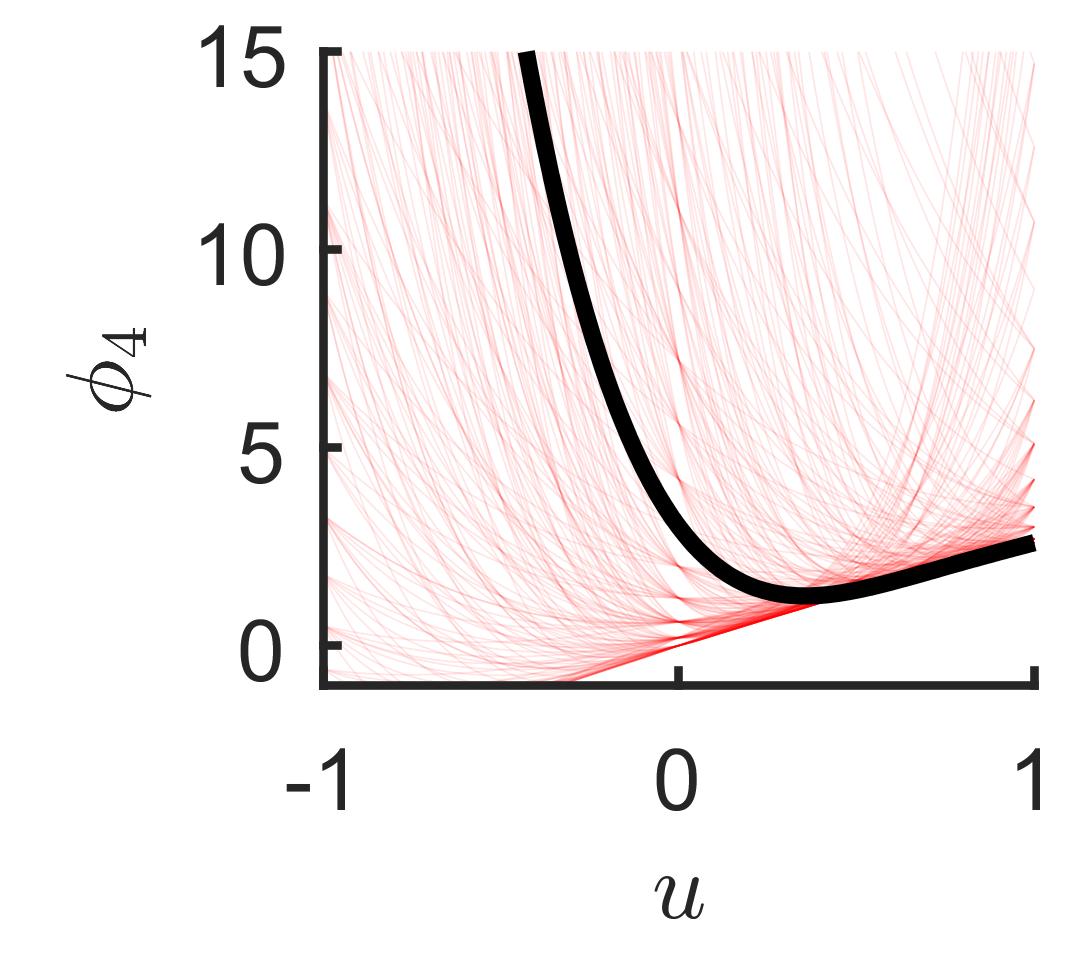}\\
				d) Functions $\phi_1$, $\phi_2$, $\phi_3$, and  $\phi_4$,
			\end{minipage} 
			\medskip
		}	
		\textcolor{red}{\raisebox{0.5mm}{\rule{0.5cm}{0.05cm}}}   : Model, 
		\textcolor{black}{\raisebox{0.5mm}{\rule{0.5cm}{0.1cm}}} : Plant.
		\vspace{-2mm}
		\captionof{figure}{Sc.1: Graphical description of the RTO problems}
		\label{fig:5_8_Exemple_5_1_sc1_Plant_and_Model}
	\end{minipage} \\
	
	\medskip
	
	\begin{minipage}[h]{\linewidth}
		\vspace*{0pt}
		{\centering
			\includegraphics[trim={2.75cm 0.2cm 0.2cm  0.2cm },clip,width=3.45cm]{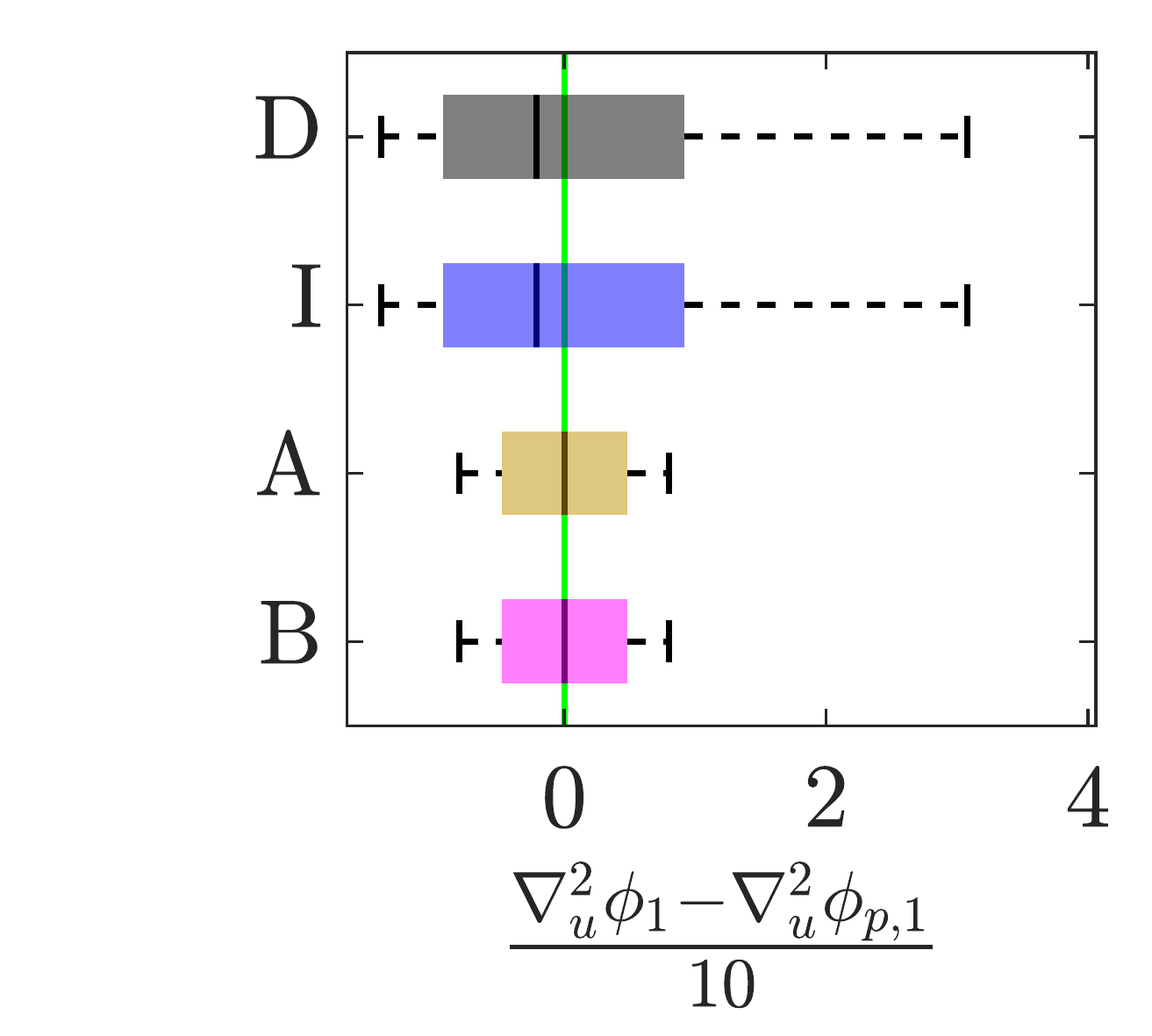}\hskip -0ex
			\includegraphics[trim={2.75cm 0.2cm 0.2cm  0.2cm },clip,width=3.45cm]{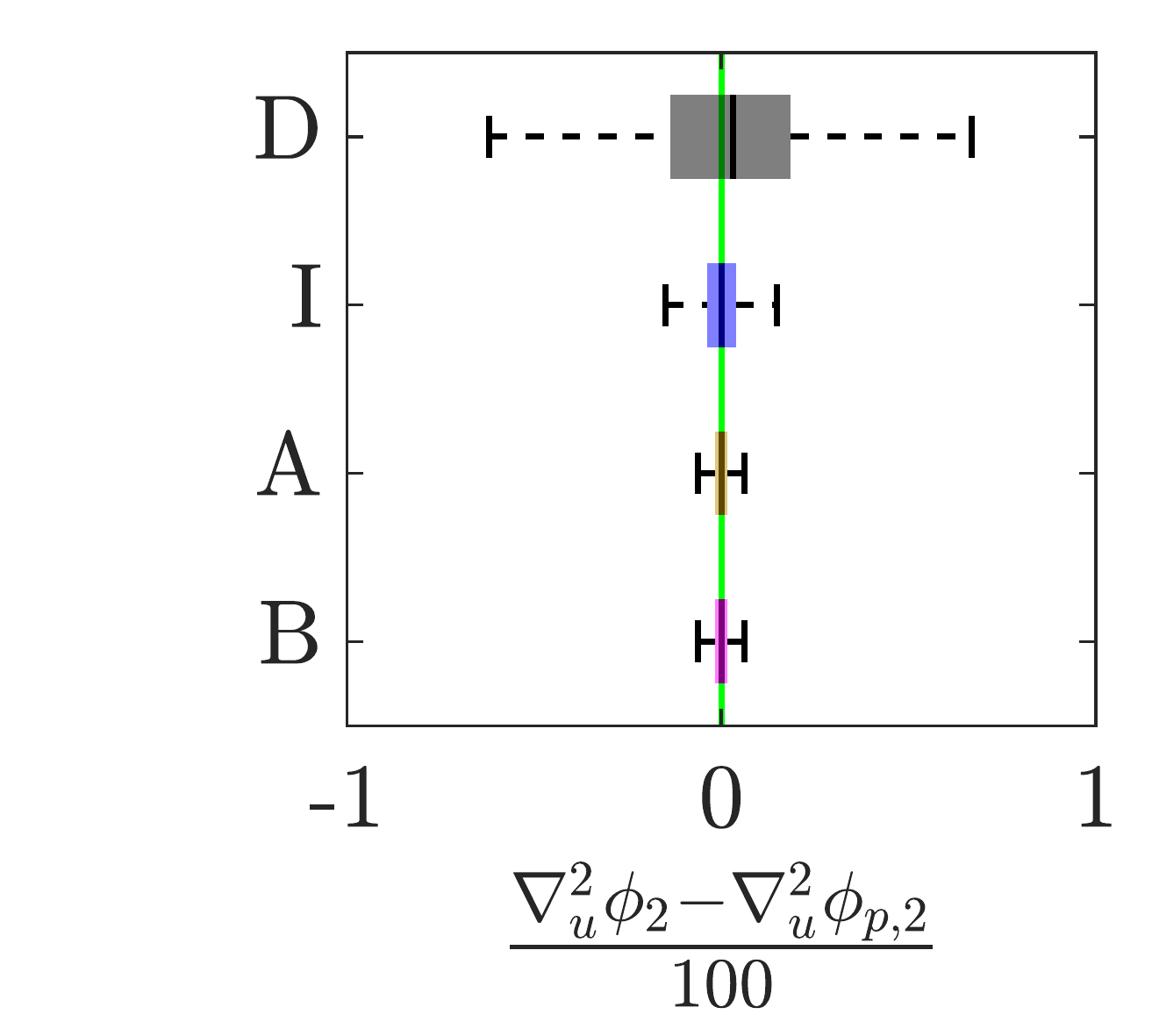}\hskip -0ex
			\includegraphics[trim={2.75cm 0.2cm 0.2cm  0.2cm },clip,width=3.45cm]{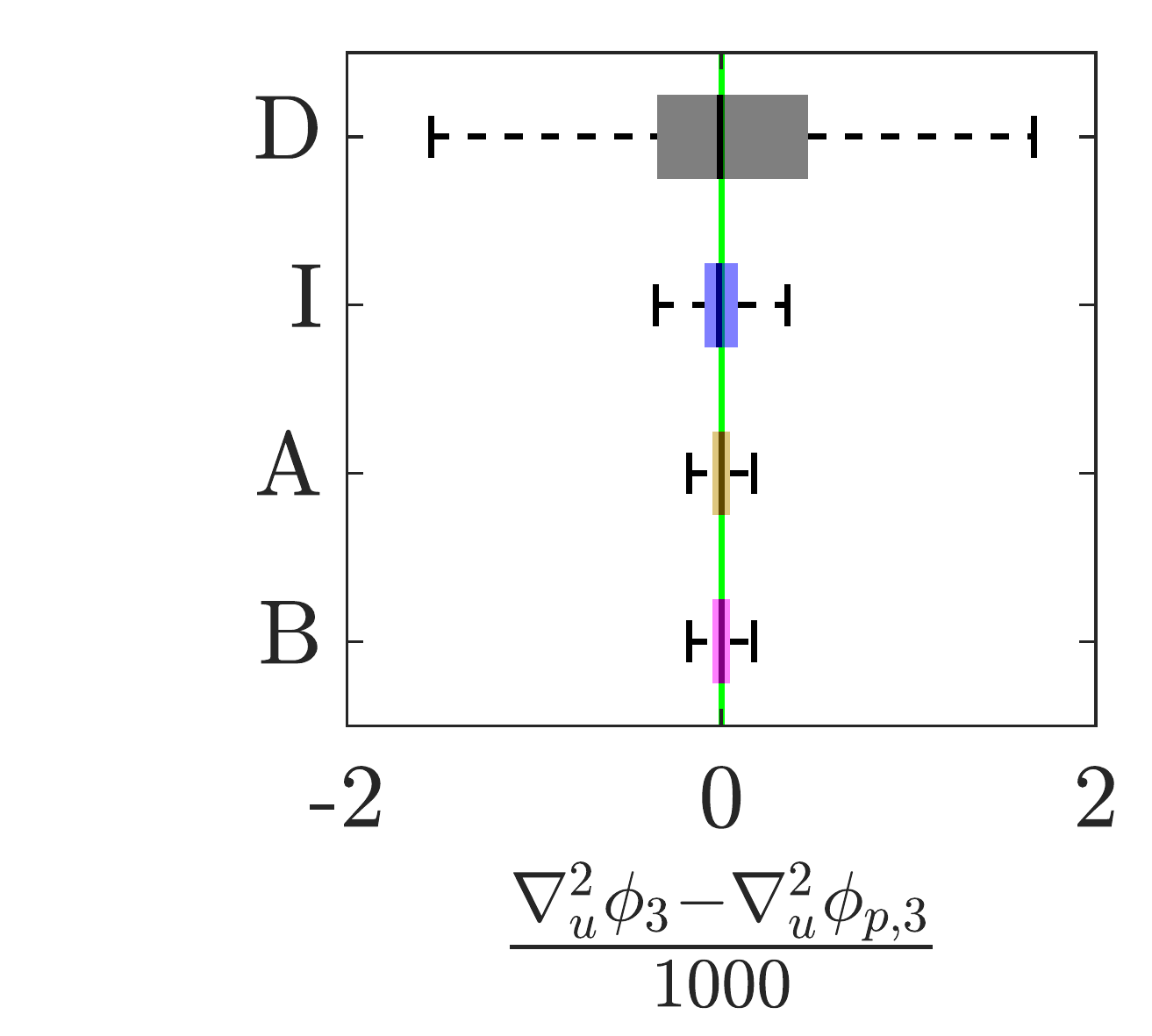}\hskip -0ex
			\includegraphics[trim={2.75cm 0.2cm 0.2cm  0.2cm },clip,width=3.45cm]{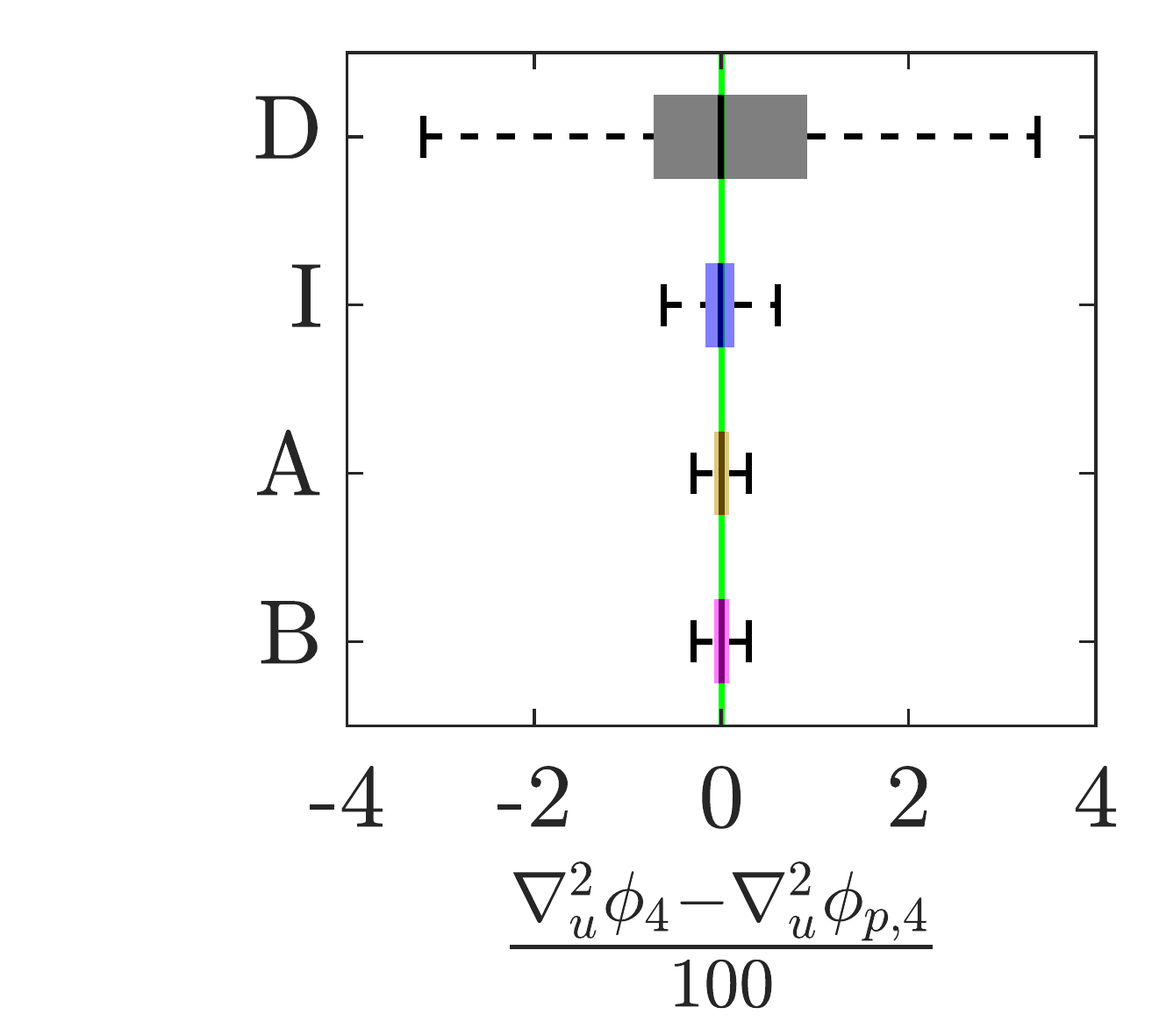}
		}
		\vspace{-2mm}
		\captionof{figure}{Sc.1:  Statistical distributions of the prediction errors on the  Hessian of the plant's cost functions at the correction point for the structures D, I, A, and B.}
		\label{fig:5_9_Exemple_5_1_sc1_Results_1}
	\end{minipage} \\
	\begin{minipage}[h]{\linewidth}
		\vspace*{0pt}
		{\centering
			\includegraphics[trim={0.7cm 0.2cm 0.2cm  0.2cm },clip,width=3.45cm]{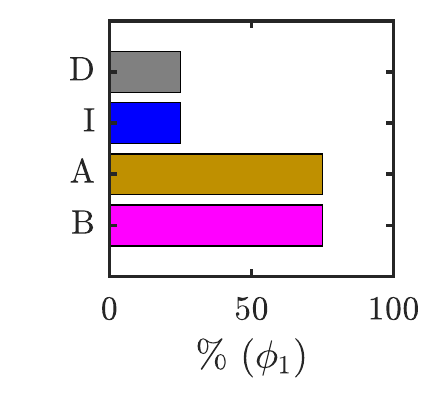}\hskip -0ex
			\includegraphics[trim={0.7cm 0.2cm 0.2cm  0.2cm },clip,width=3.45cm]{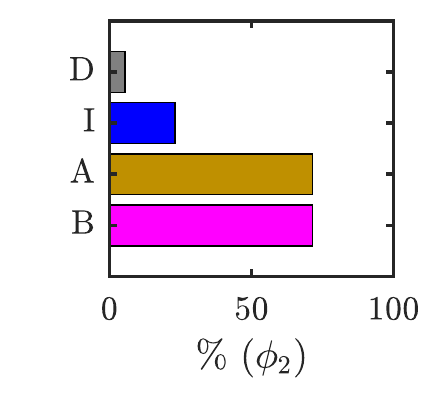}\hskip -0ex
			\includegraphics[trim={0.7cm 0.2cm 0.2cm  0.2cm },clip,width=3.45cm]{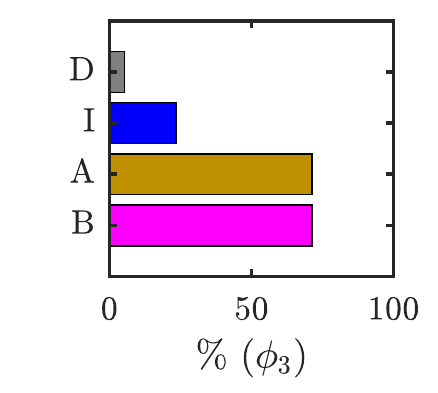}\hskip -0ex
			\includegraphics[trim={0.7cm 0.2cm 0.2cm  0.2cm },clip,width=3.45cm]{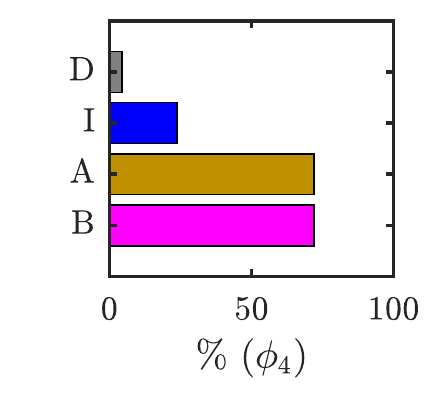} \\
		}
		\vspace{-2mm}
		\captionof{figure}{Sc.1: For each correction structure one gives here the percentage of cases for which no other structure provides better results  \textit{(if two structures provide the same best result then both take the point)}}
		\label{fig:5_10_Exemple_5_1_sc1_Results_2}
	\end{minipage} \\
\end{minipage} 

\noindent
	\begin{minipage}[h]{\linewidth}
	\begin{center}
		\textbf{Scenario 2: NL$\rightarrow$A}
	\end{center}
	
	\vspace{-3mm}
	\begin{minipage}[h]{\linewidth}
		\vspace*{0pt}
		{\centering
			\begin{minipage}[t]{4.45cm}{\centering%
					\includegraphics[width=4.45cm]{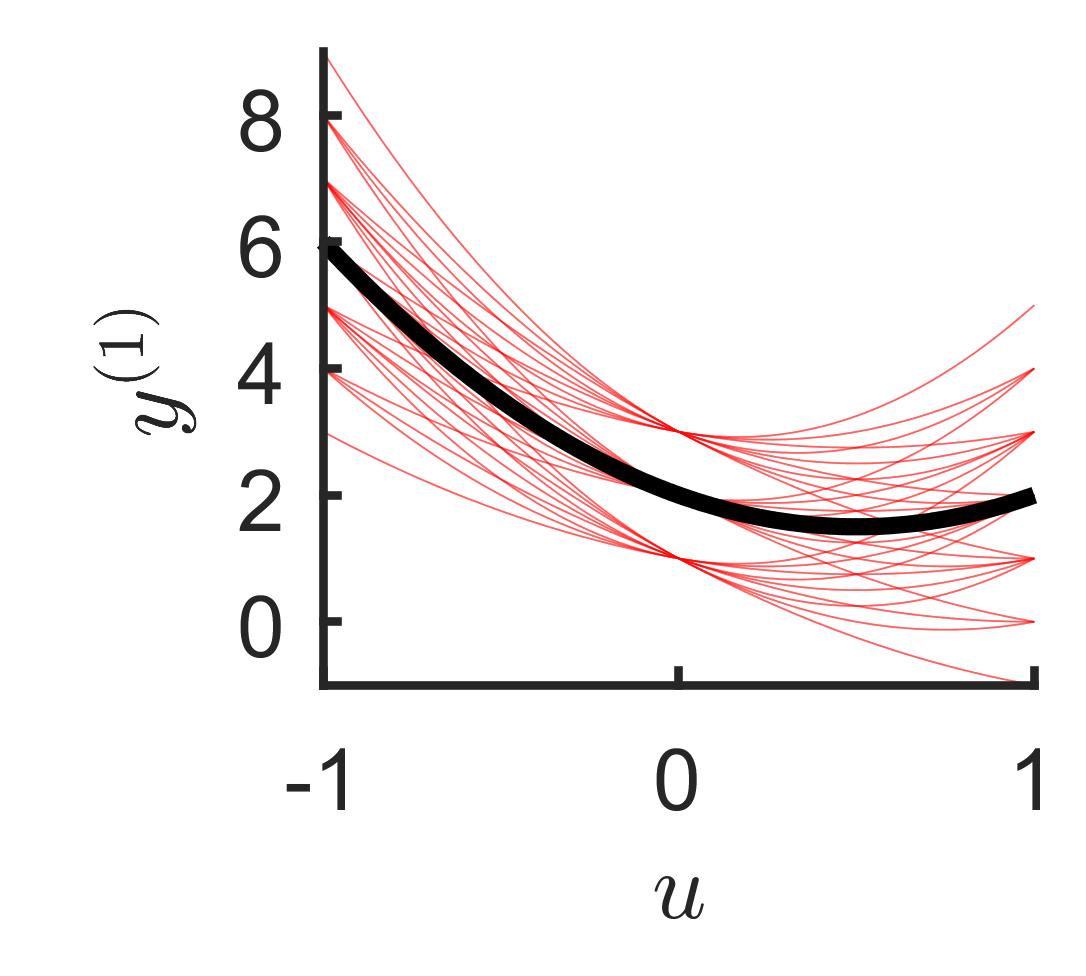}\\
					a) Function $f^{(1)}$ and $f_p^{(1)}$}
			\end{minipage}\hskip -0ex
			\begin{minipage}[t]{4.45cm}\centering%
				\includegraphics[width=4.45cm]{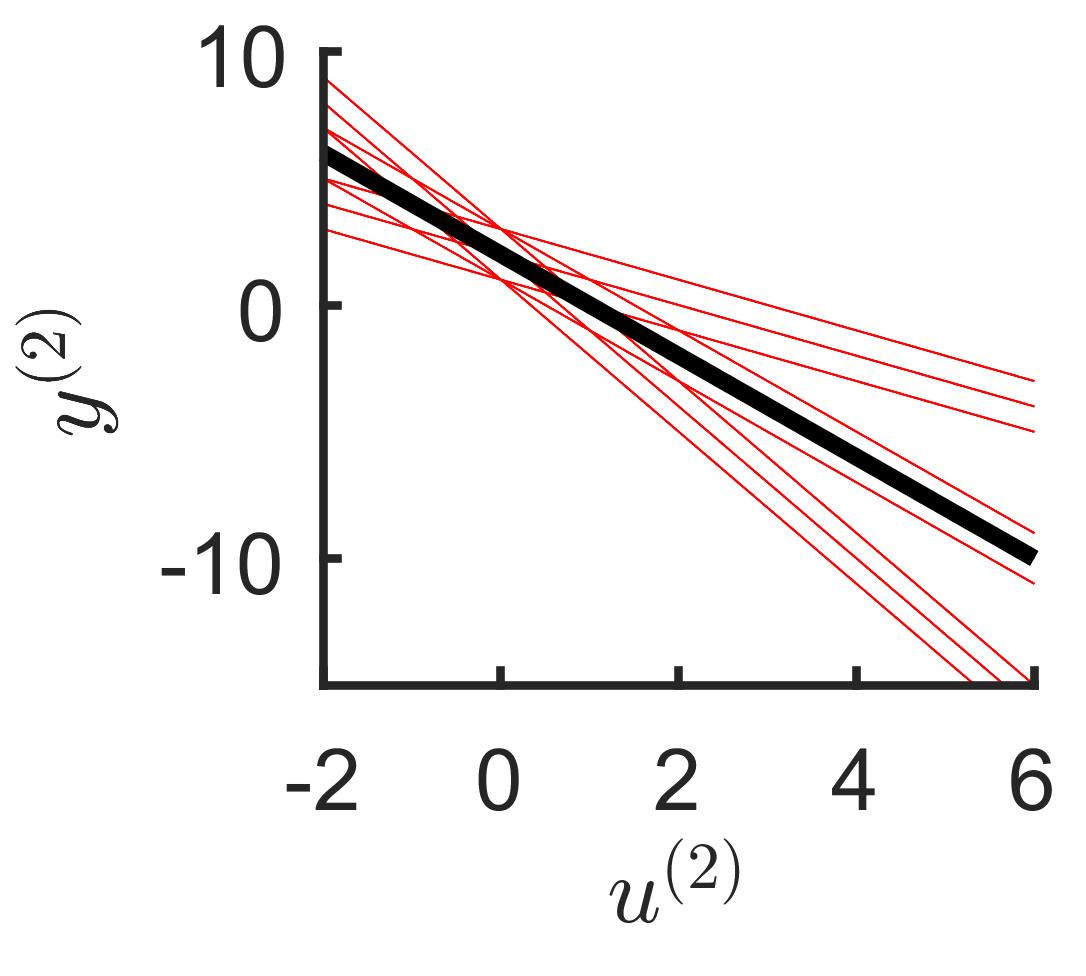}\\
				b) Function $f^{(2)}$ and $f_p^{(2)}$
			\end{minipage}\hskip -0ex
			\begin{minipage}[t]{4.45cm}\centering%
				\includegraphics[width=4.45cm]{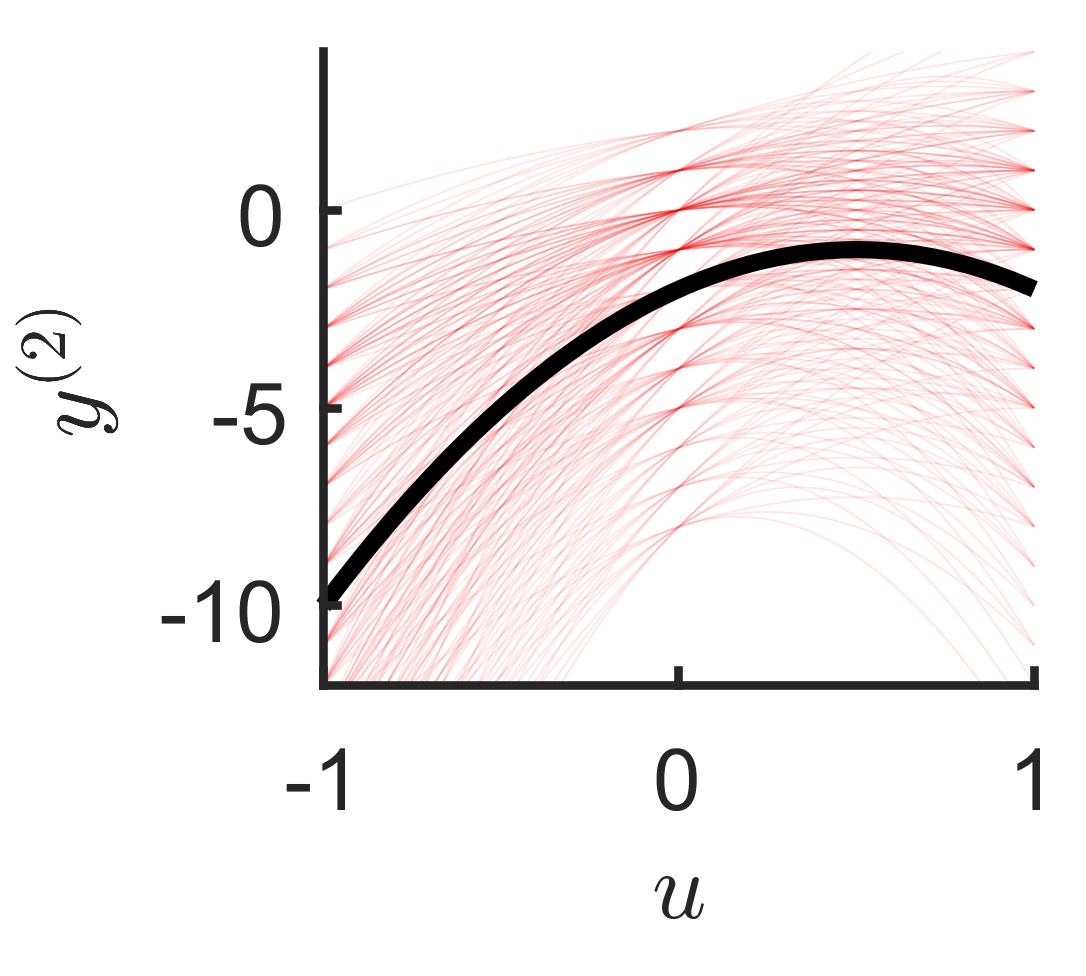}\\
				c) Function $f$ and $f_p$
			\end{minipage} \\
			
			\medskip
			
			\begin{minipage}[h]{\linewidth} \centering
				\includegraphics[trim={0.1cm 0.2cm 0.2cm  0.1cm },clip,width=3.5cm]{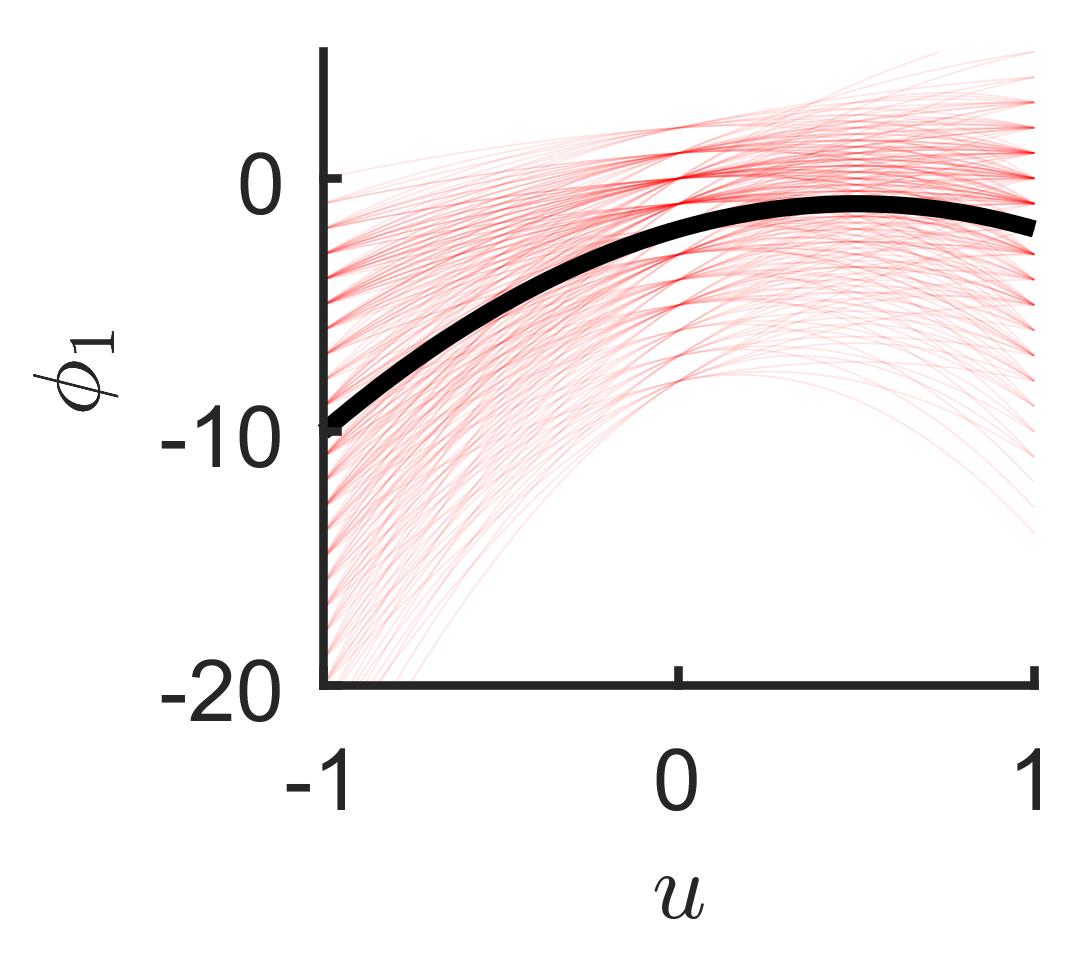}\hskip -0ex
				\includegraphics[trim={0.1cm 0.2cm 0.2cm  0.1cm },clip,width=3.5cm]{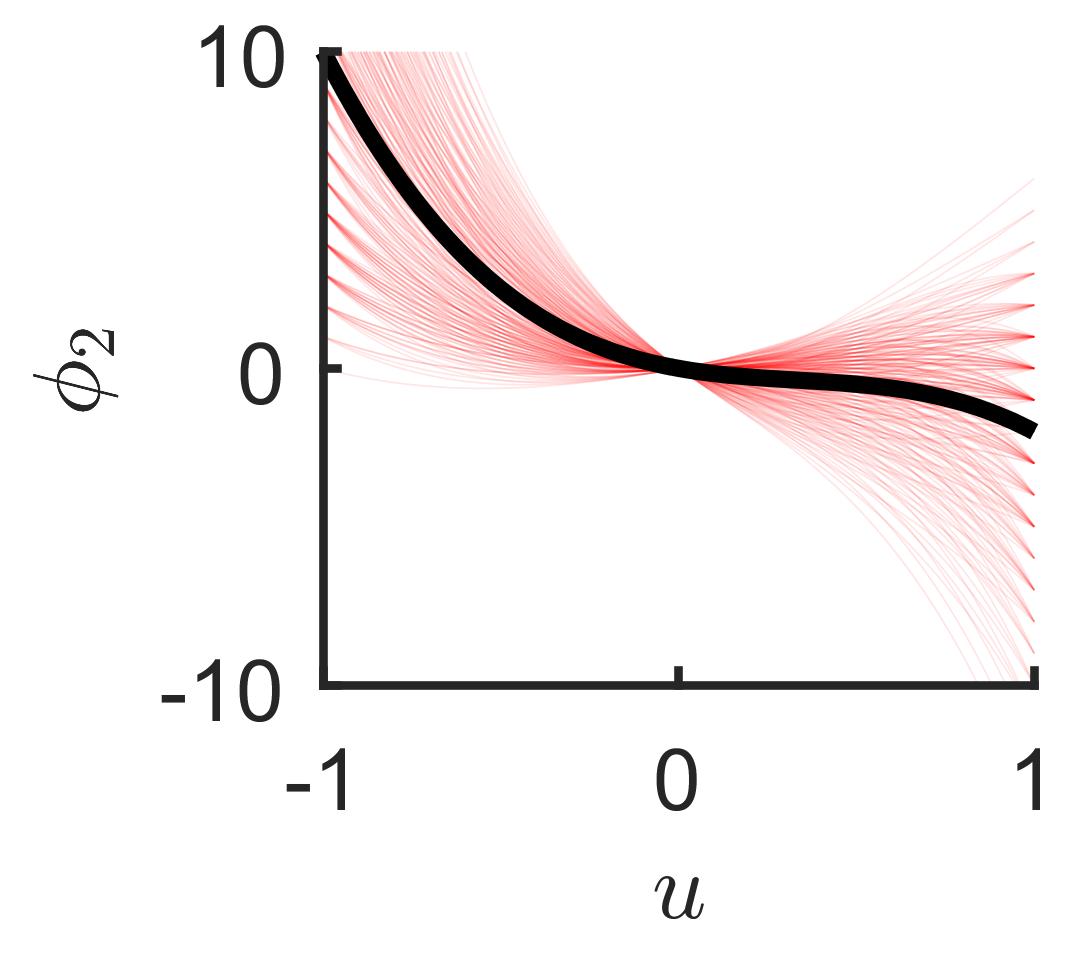}\hskip -0ex
				\includegraphics[trim={0.1cm 0.2cm 0.2cm  0.1cm },clip,width=3.5cm]{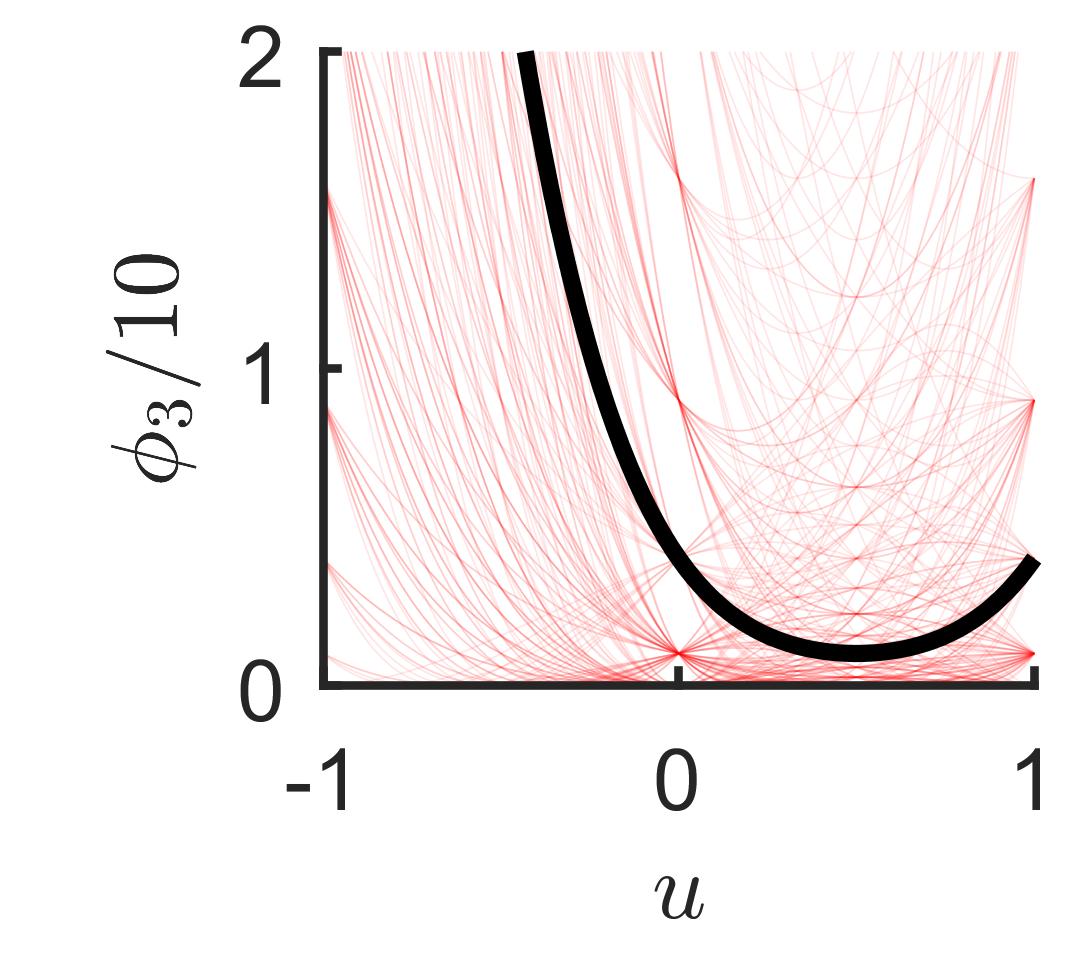}\hskip -0ex
				\includegraphics[trim={0.1cm 0.2cm 0.2cm  0.1cm },clip,width=3.5cm]{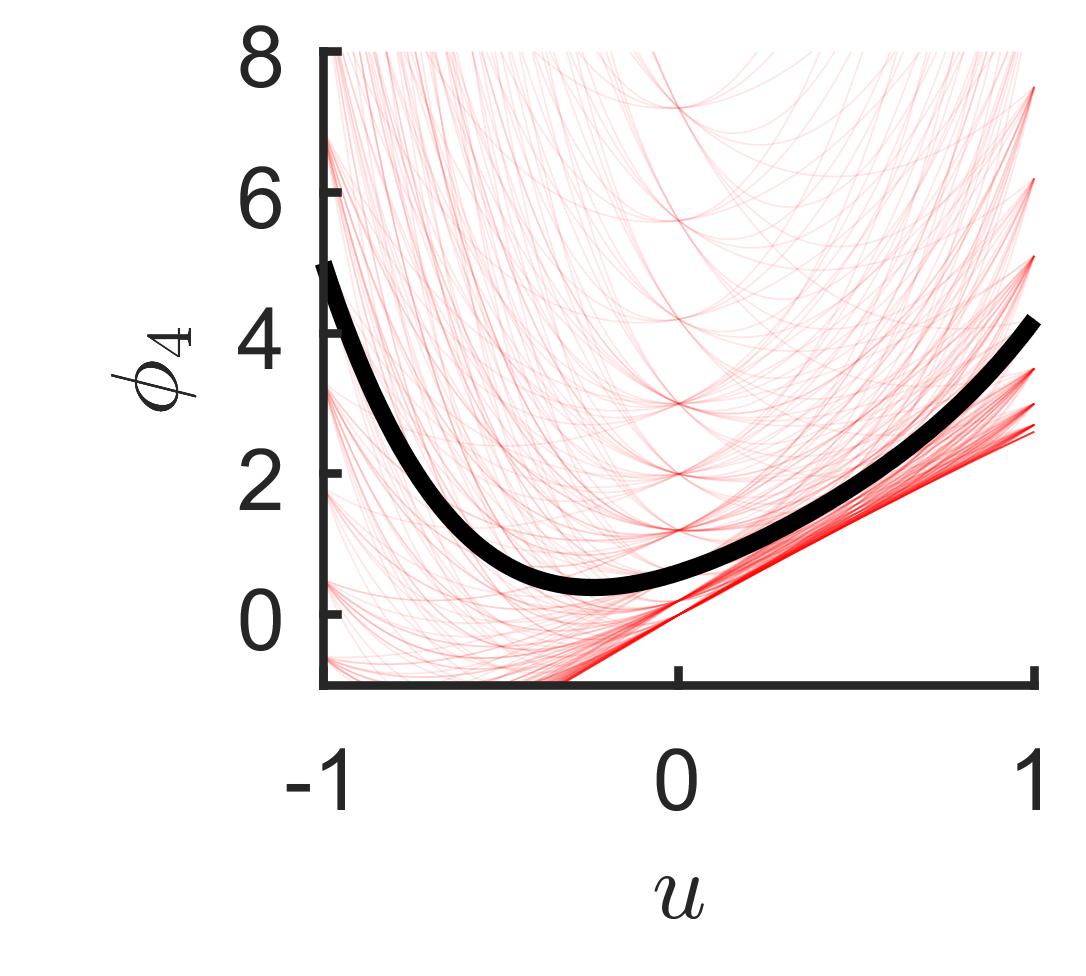}\\
				d) Function $\phi_1$, $\phi_2$, $\phi_3$, et  $\phi_4$,
			\end{minipage} 
			
			\medskip
			
		}	
		\textcolor{red}{\raisebox{0.5mm}{\rule{0.5cm}{0.05cm}}}   : Model, 
		\textcolor{black}{\raisebox{0.5mm}{\rule{0.5cm}{0.1cm}}} : Plant.
		\vspace{-2mm}
		\captionof{figure}{Sc.2:Graphical description of the RTO problems}
		\label{fig:5_11_Exemple_5_1_sc2_Plant_and_Model}
	\end{minipage} \\
	
	\medskip
	
	\begin{minipage}[h]{\linewidth}
		\vspace*{0pt}
		{\centering
			\includegraphics[trim={2.75cm 0.2cm 0.2cm  0.2cm },clip,width=3.45cm]{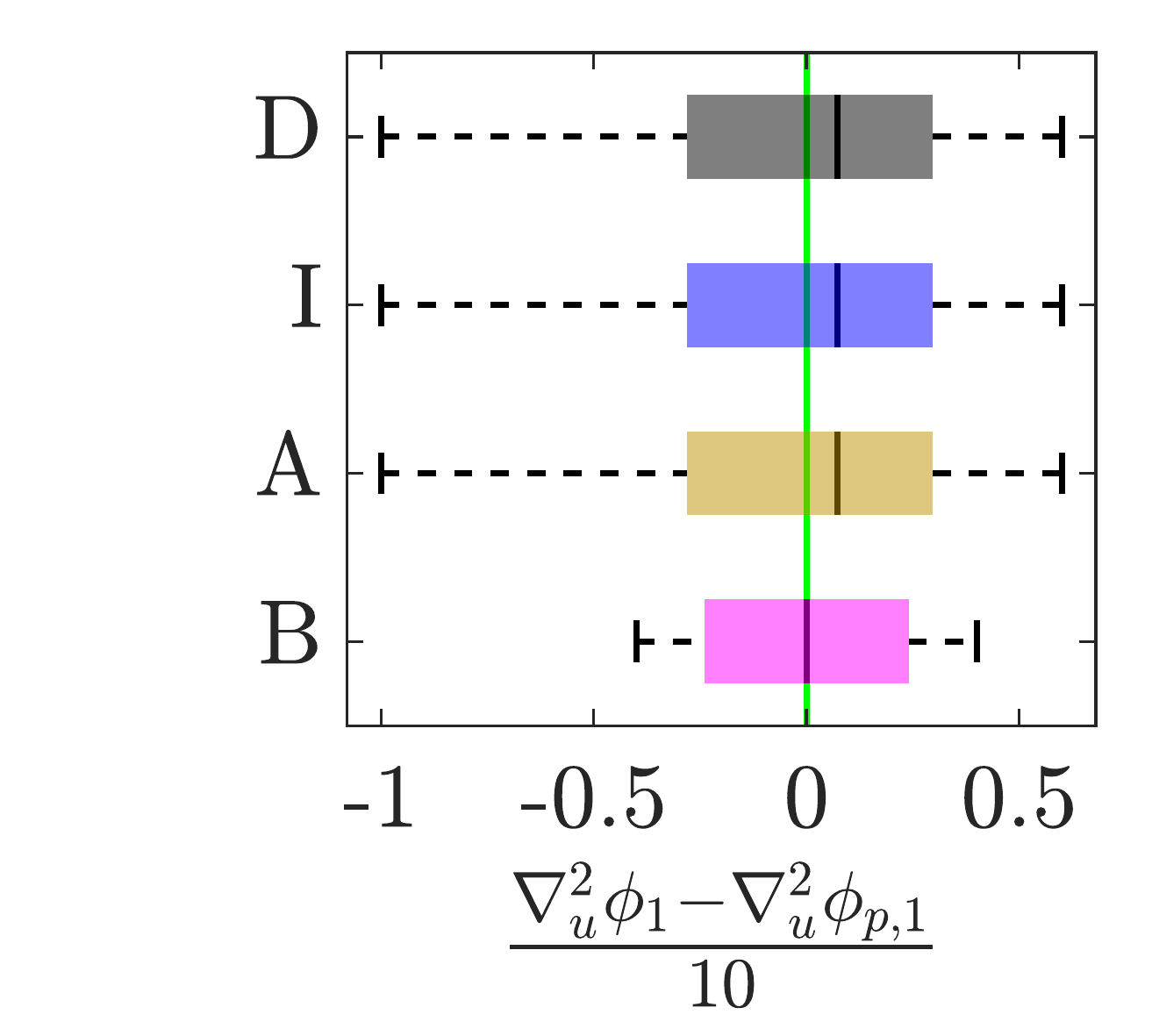}\hskip -0ex
			\includegraphics[trim={2.75cm 0.2cm 0.2cm  0.2cm },clip,width=3.45cm]{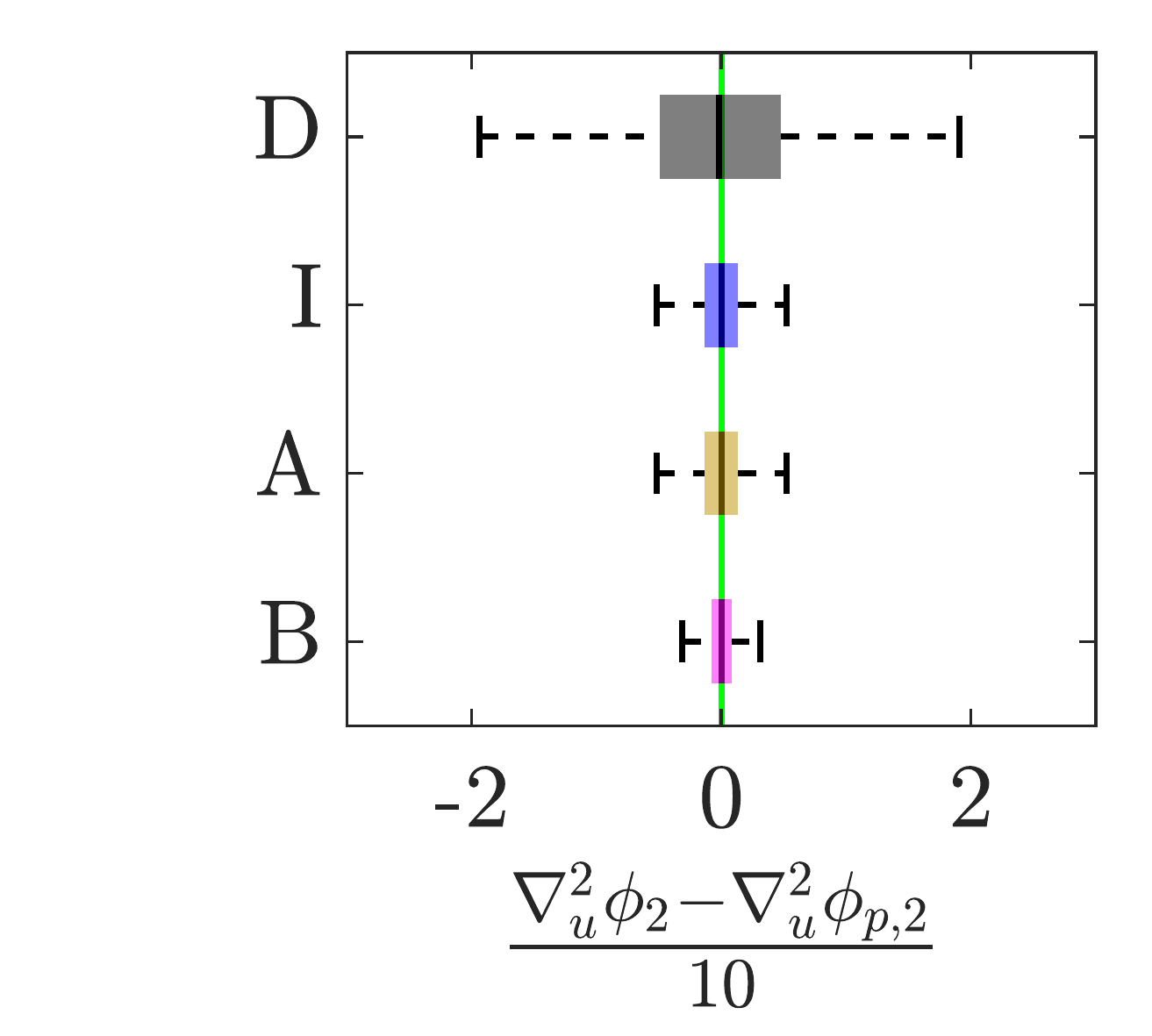}\hskip -0ex
			\includegraphics[trim={2.75cm 0.2cm 0.2cm  0.2cm },clip,width=3.45cm]{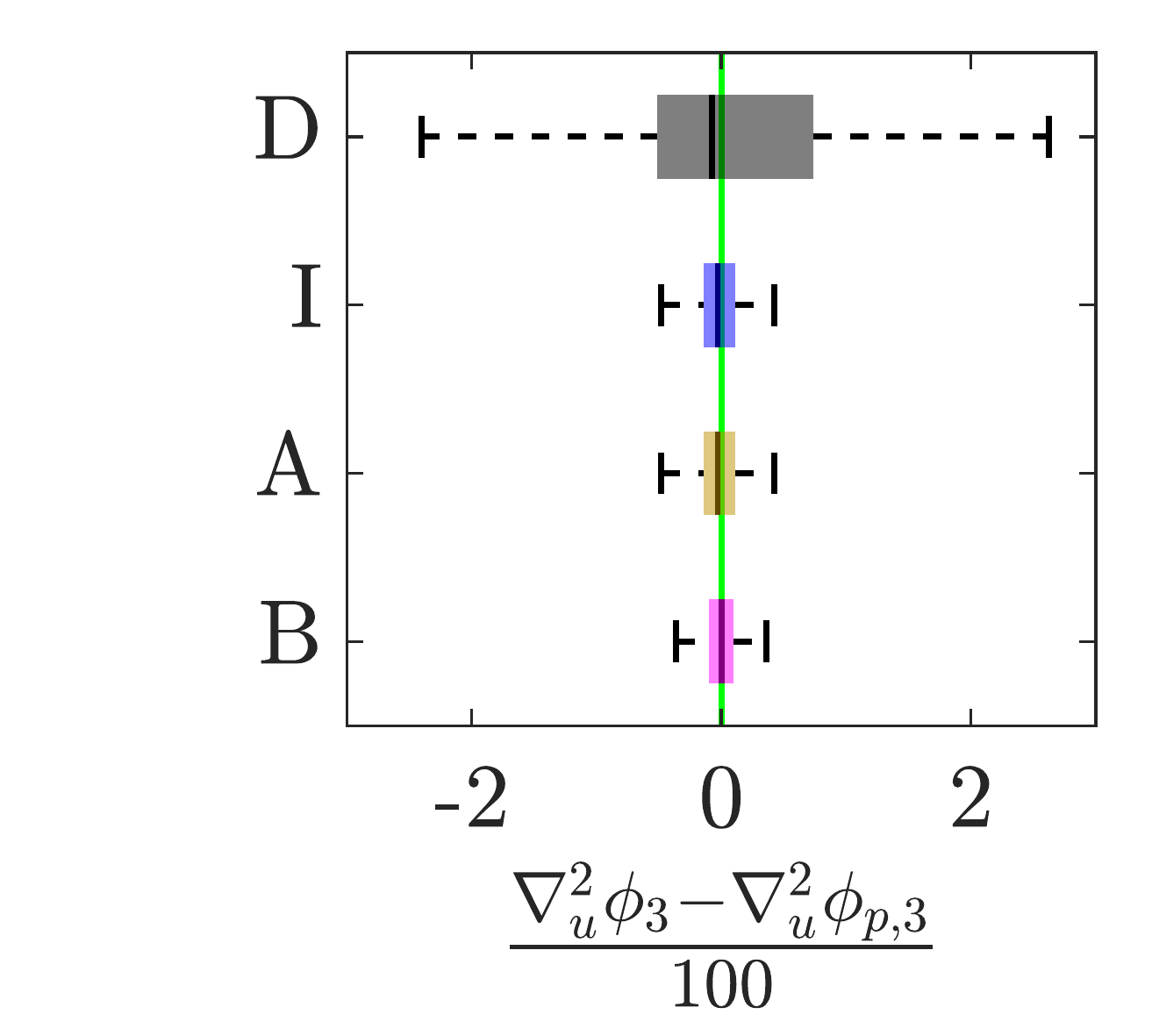}\hskip -0ex
			\includegraphics[trim={2.75cm 0.2cm 0.2cm  0.2cm },clip,width=3.45cm]{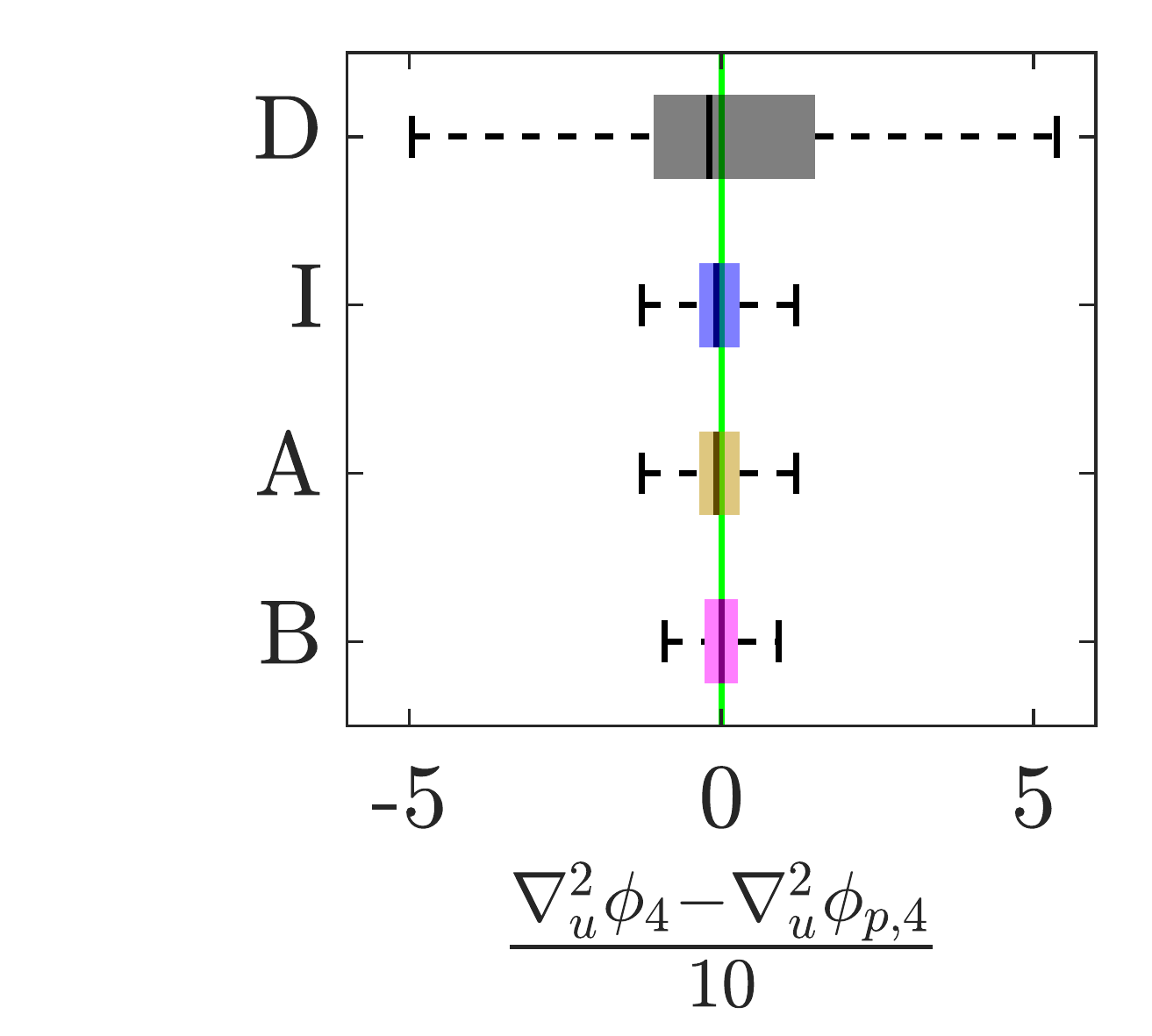}
		}
		\vspace{-2mm}
		\captionof{figure}{Sc.2: Statistical distributions of the prediction errors on the  Hessian of the plant cost's functions at the correction point for the structures D, I, A, and B.}
		\label{fig:5_12_Exemple_5_1_sc2_Results}
	\end{minipage} \\
	\begin{minipage}[h]{\linewidth}
		\vspace*{0pt}
		{\centering
			\includegraphics[trim={0.7cm 0.2cm 0.2cm  0.2cm },clip,width=3.45cm]{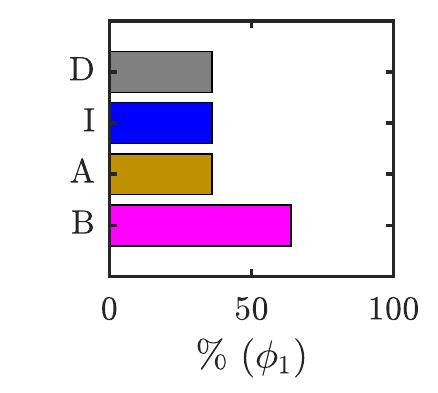}\hskip -0ex
			\includegraphics[trim={0.7cm 0.2cm 0.2cm  0.2cm },clip,width=3.45cm]{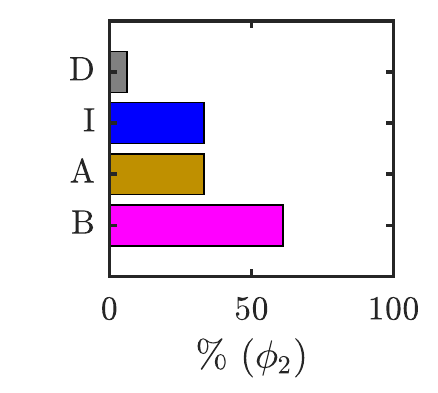}\hskip -0ex
			\includegraphics[trim={0.7cm 0.2cm 0.2cm  0.2cm },clip,width=3.45cm]{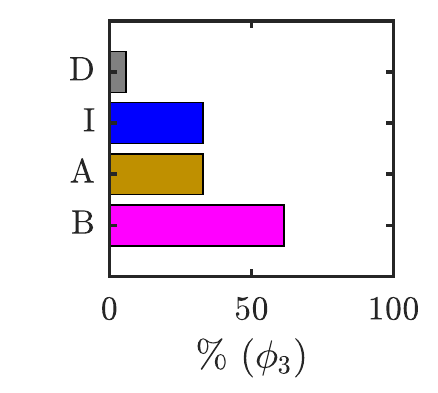}\hskip -0ex
			\includegraphics[trim={0.7cm 0.2cm 0.2cm  0.2cm },clip,width=3.45cm]{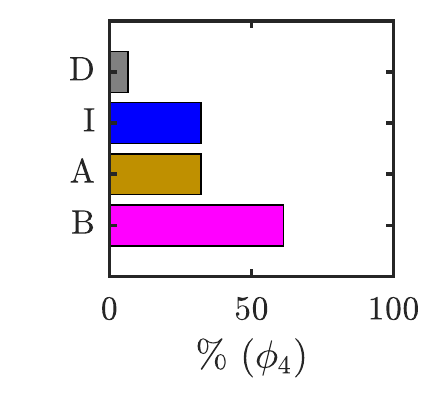} \\
		}
		\vspace{-2mm}
		\captionof{figure}{Sc.2: For each correction structure one gives here the percentage of cases for which no other structure provides better results  \textit{(if two structures provide the same best result then both take the point)}}
		\label{fig:5_13_Exemple_5_1_sc2_Results_2}
	\end{minipage} \\
\end{minipage} 

	\noindent
	\begin{minipage}[h]{\linewidth}
		\begin{center}
			\textbf{Scenario 3: NL$\rightarrow$NL}
		\end{center}
	
		\vspace{-3mm}
		\begin{minipage}[h]{\linewidth}
			\vspace*{0pt}
			{\centering
				\begin{minipage}[t]{4.45cm}{\centering%
					\includegraphics[width=4.45cm]{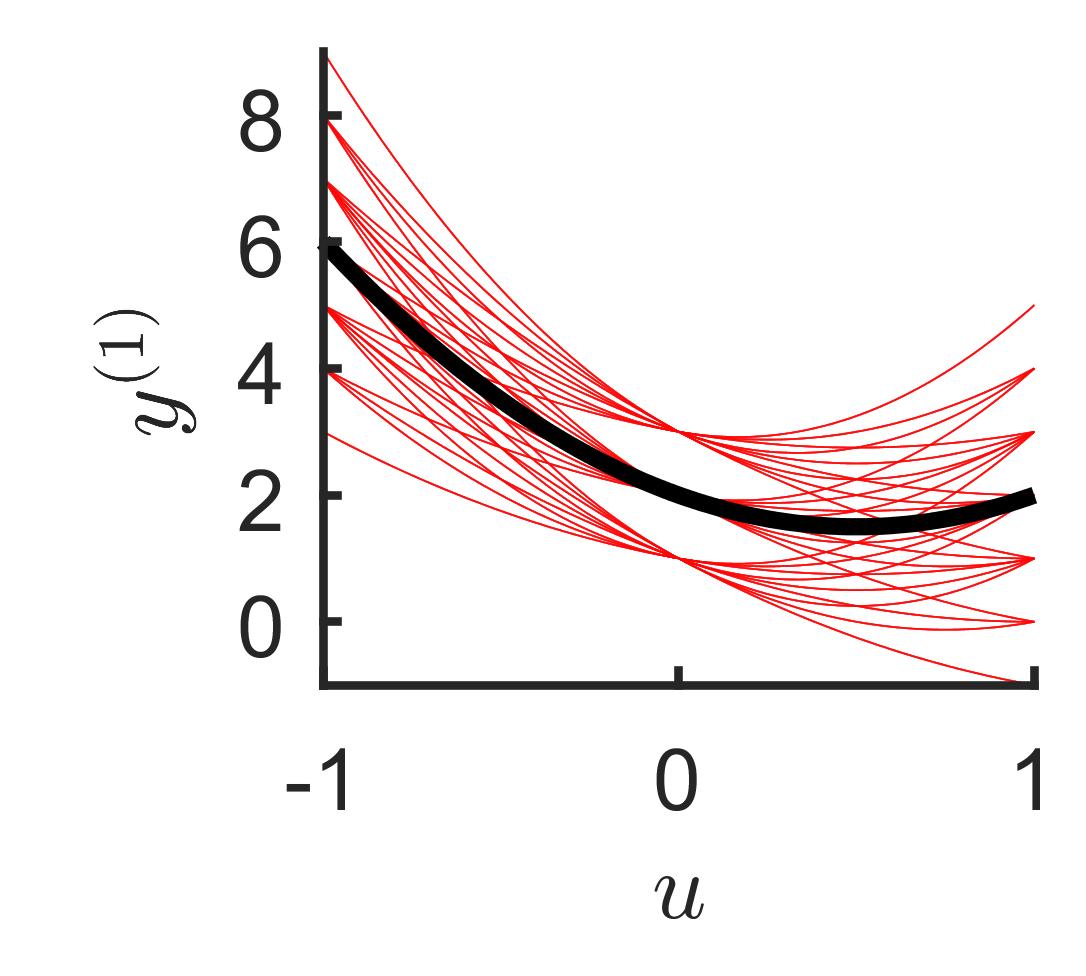}\\
					a) Functions $f^{(1)}$ and $f_p^{(1)}$}
				\end{minipage}\hskip -0ex
				\begin{minipage}[t]{5cm}\centering%
					\includegraphics[width=4.45cm]{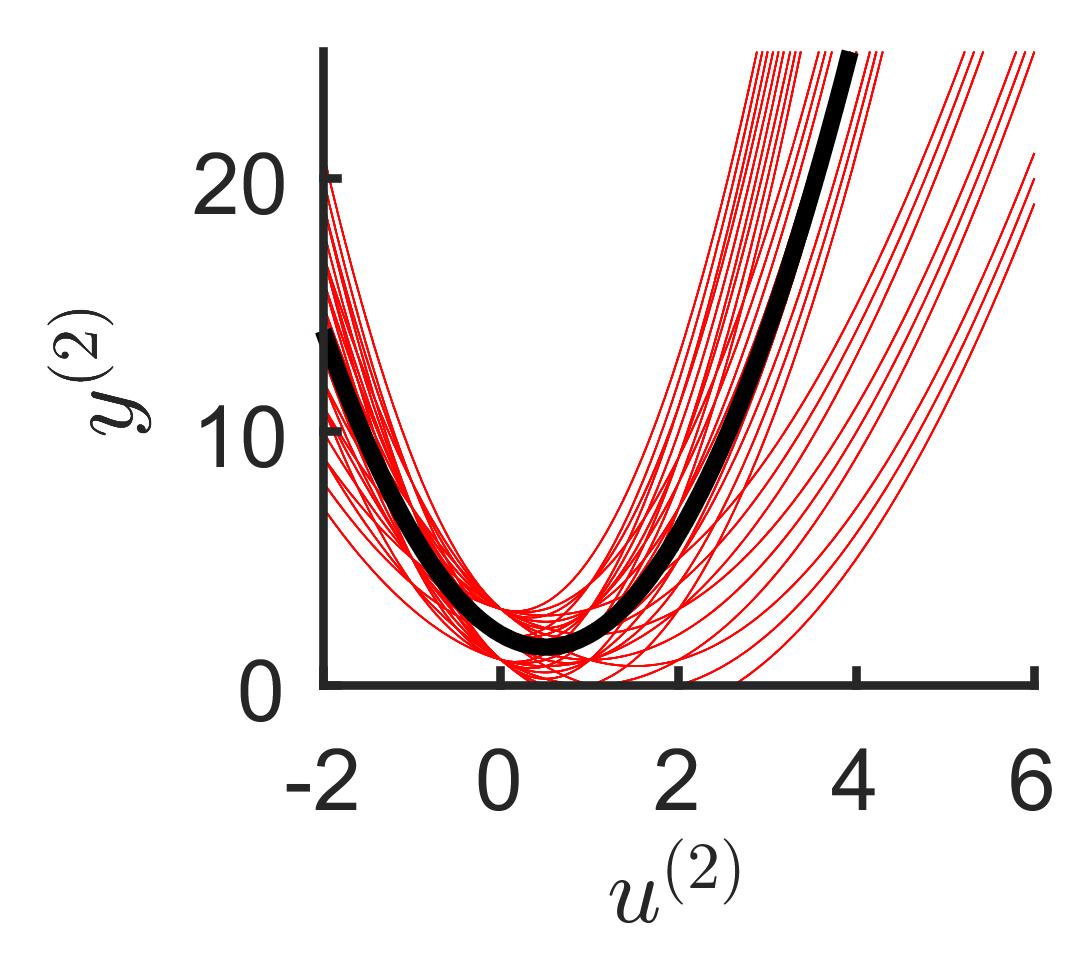}\\
					b) Functions $f^{(2)}$ and $f_p^{(2)}$
				\end{minipage}\hskip -0ex
				\begin{minipage}[t]{4.45cm}\centering%
					\includegraphics[width=4.45cm]{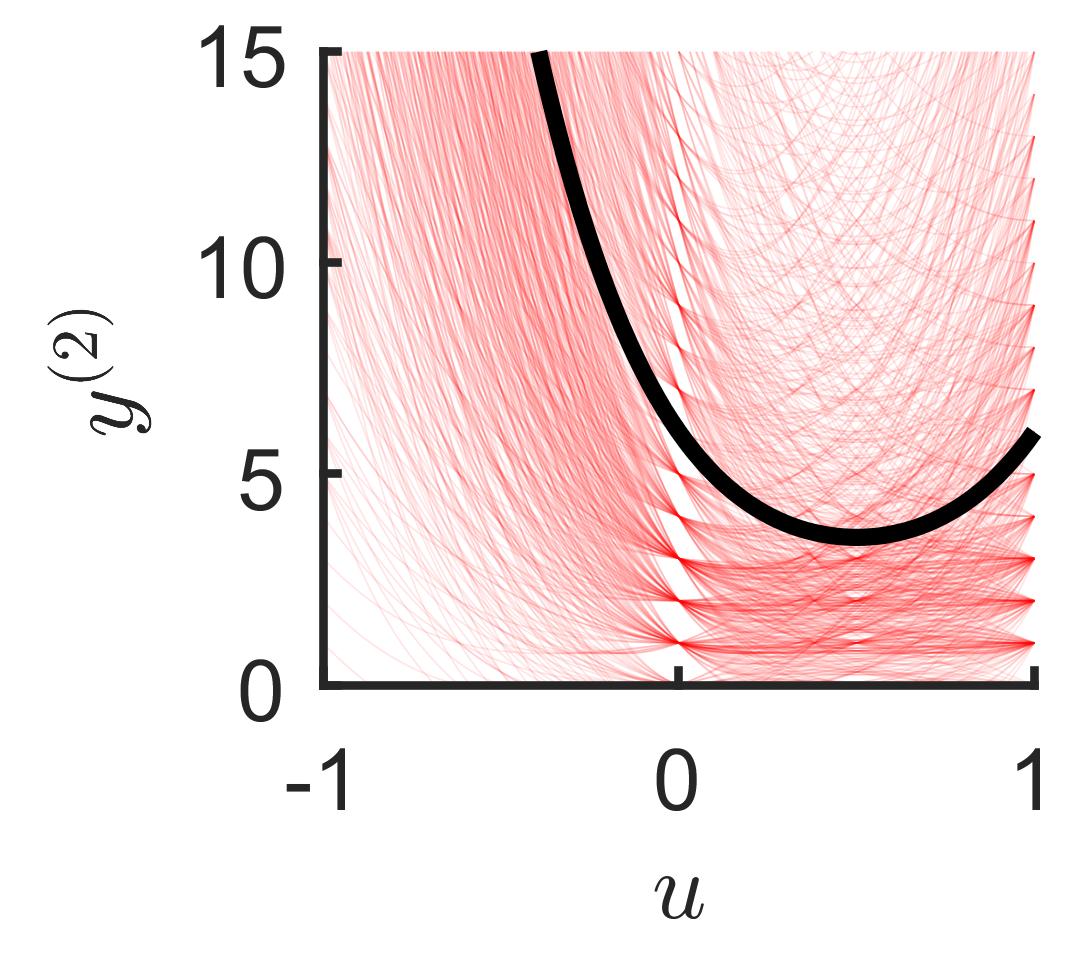}\\
					c) Functions $f$ and $f_p$
				\end{minipage} \\
			
				\medskip
			
				\begin{minipage}[h]{\linewidth} \centering
					\includegraphics[trim={0.1cm 0.2cm 0.2cm  0.1cm },clip,width=3.5cm]{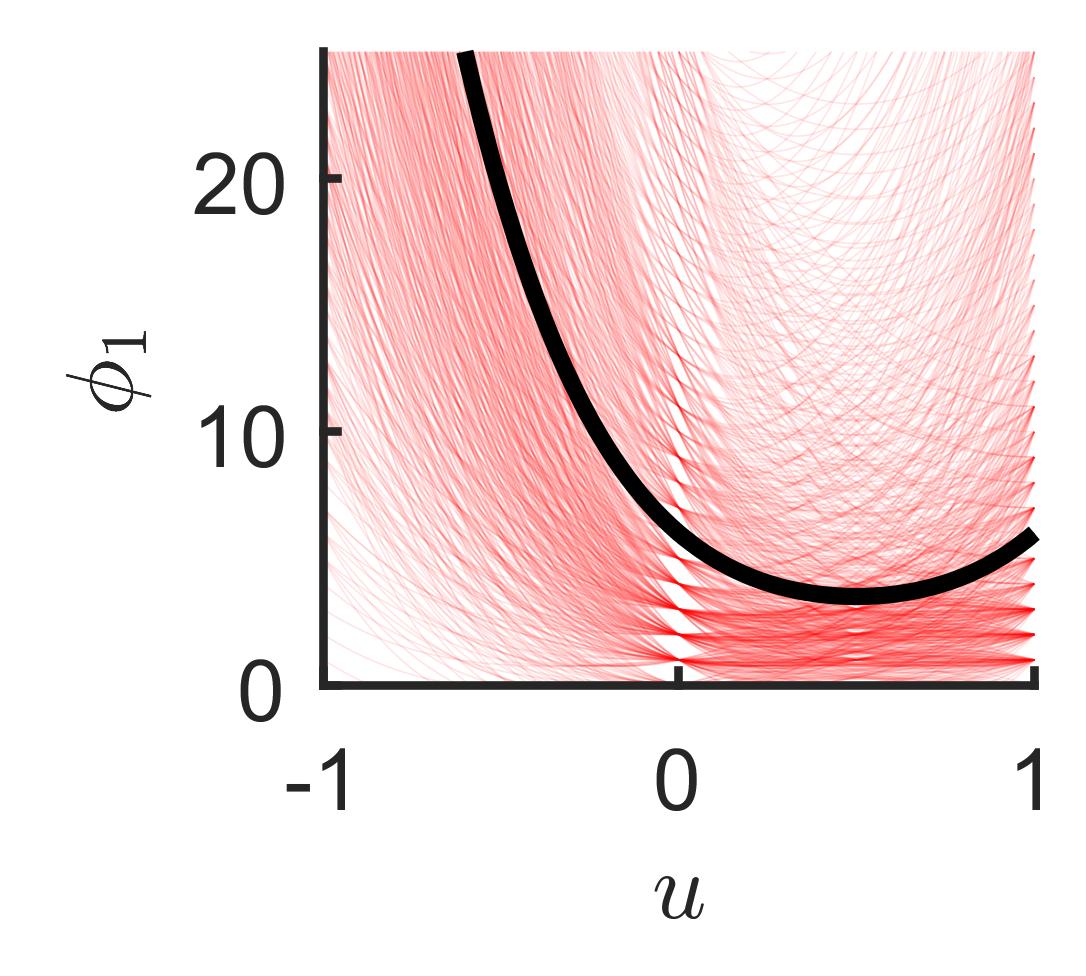}\hskip -0ex
					\includegraphics[trim={0.1cm 0.2cm 0.2cm  0.1cm },clip,width=3.5cm]{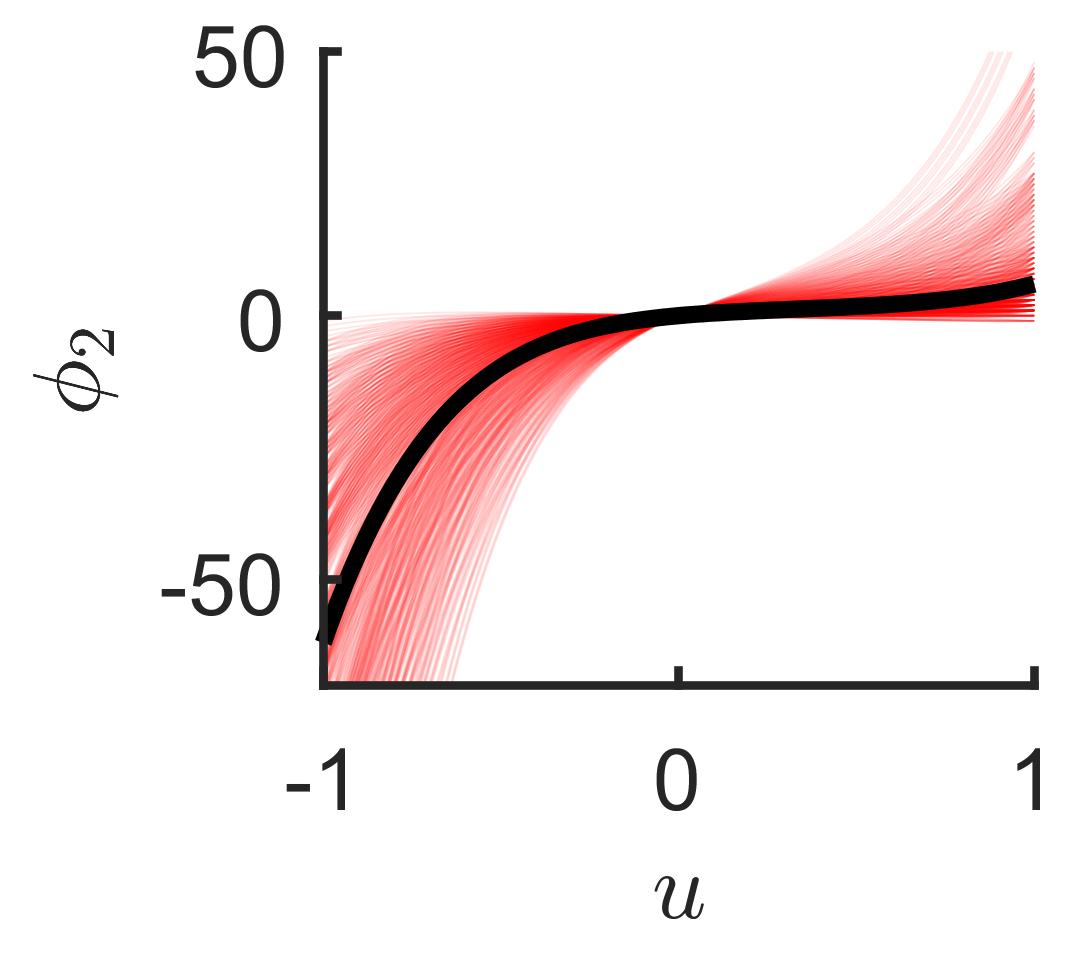}\hskip -0ex
					\includegraphics[trim={0.1cm 0.2cm 0.2cm  0.1cm },clip,width=3.5cm]{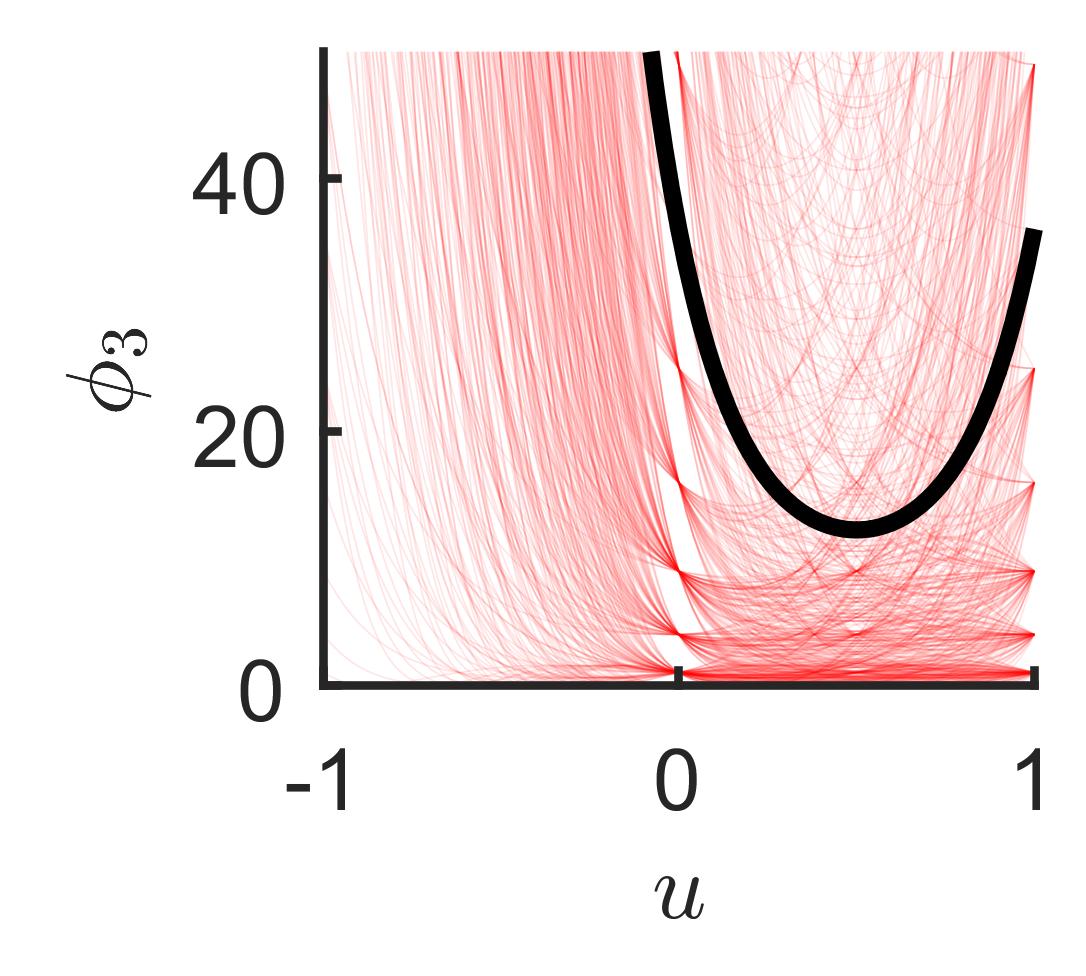}\hskip -0ex
					\includegraphics[trim={0.1cm 0.2cm 0.2cm  0.1cm },clip,width=3.5cm]{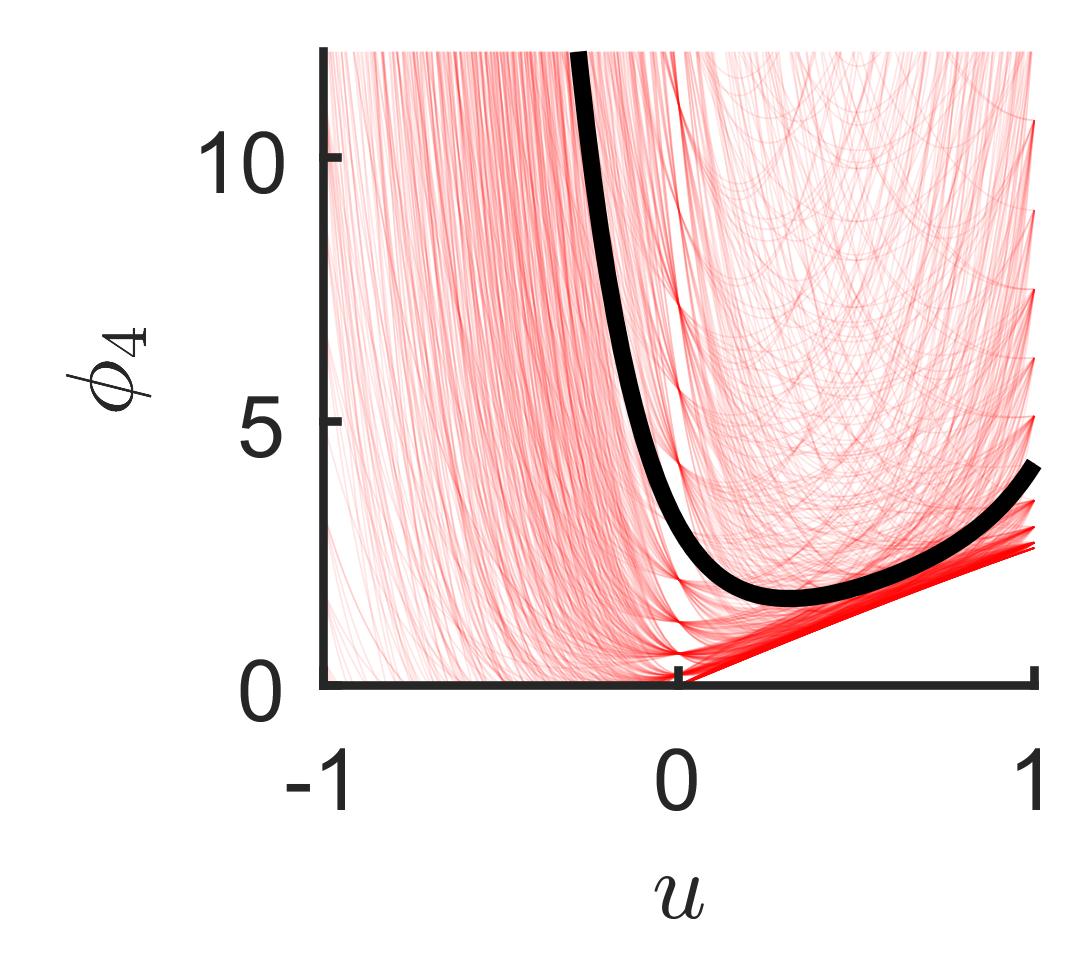}\\
					d) Functions $\phi_1$, $\phi_2$, $\phi_3$, and  $\phi_4$,
				\end{minipage} 
			
				\medskip
	
				}	
				\textcolor{red}{\raisebox{0.5mm}{\rule{0.5cm}{0.05cm}}}   : Model, 
				\textcolor{black}{\raisebox{0.5mm}{\rule{0.5cm}{0.1cm}}} : Plant.
				\vspace{-2mm}
				\captionof{figure}{Sc.3: Graphical description of the RTO problems.}
				\label{fig:5_14_Exemple_5_1_Plant_and_Model}
		\end{minipage} \\
		
		\medskip
		
		\begin{minipage}[h]{\linewidth}
			\vspace*{0pt}
			{\centering
				\includegraphics[trim={2.75cm 0.2cm 0.2cm  0.2cm },clip,width=3.45cm]{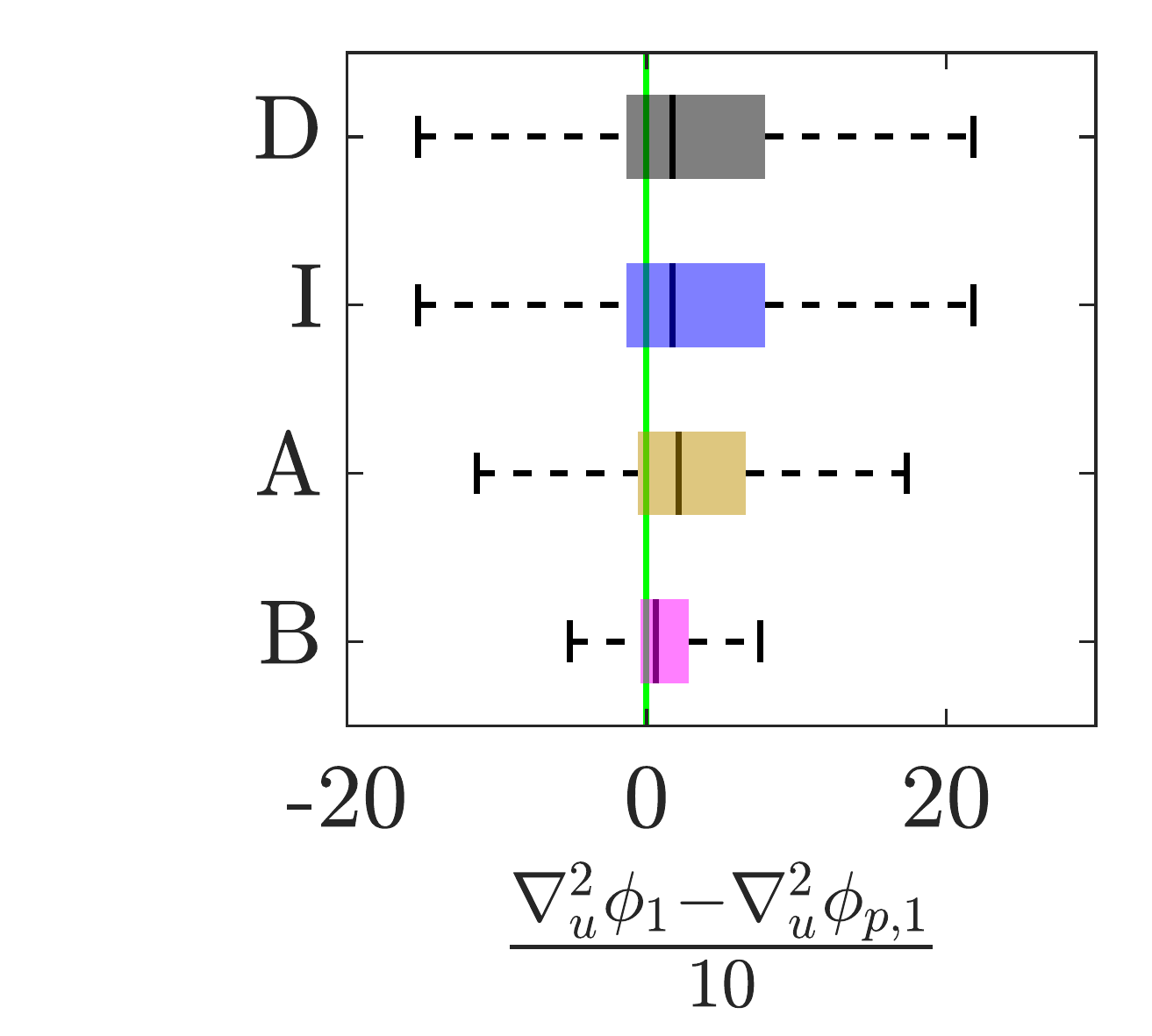}\hskip -0ex
				\includegraphics[trim={2.75cm 0.2cm 0.2cm  0.2cm },clip,width=3.45cm]{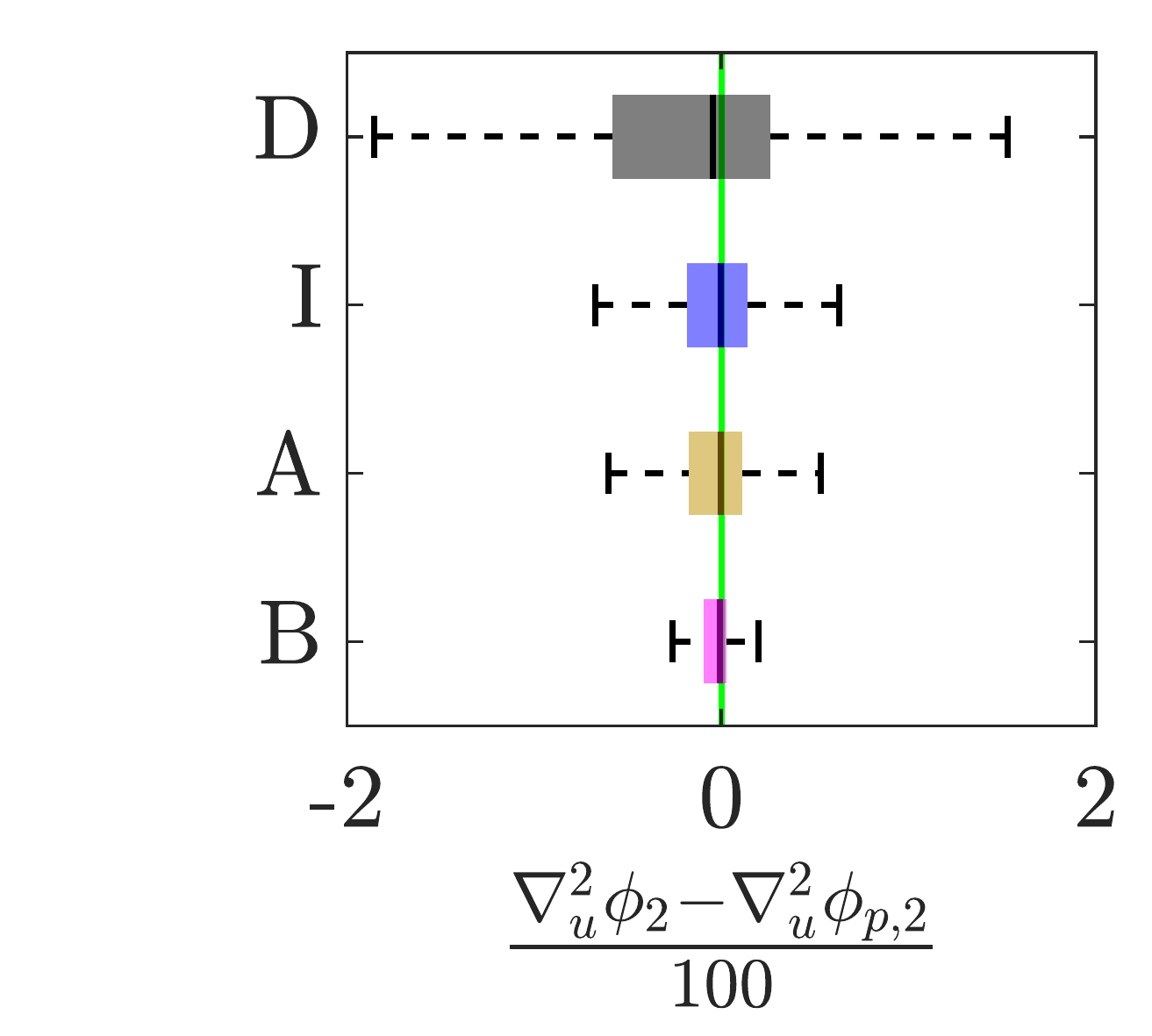}\hskip -0ex
				\includegraphics[trim={2.75cm 0.2cm 0.2cm  0.2cm },clip,width=3.45cm]{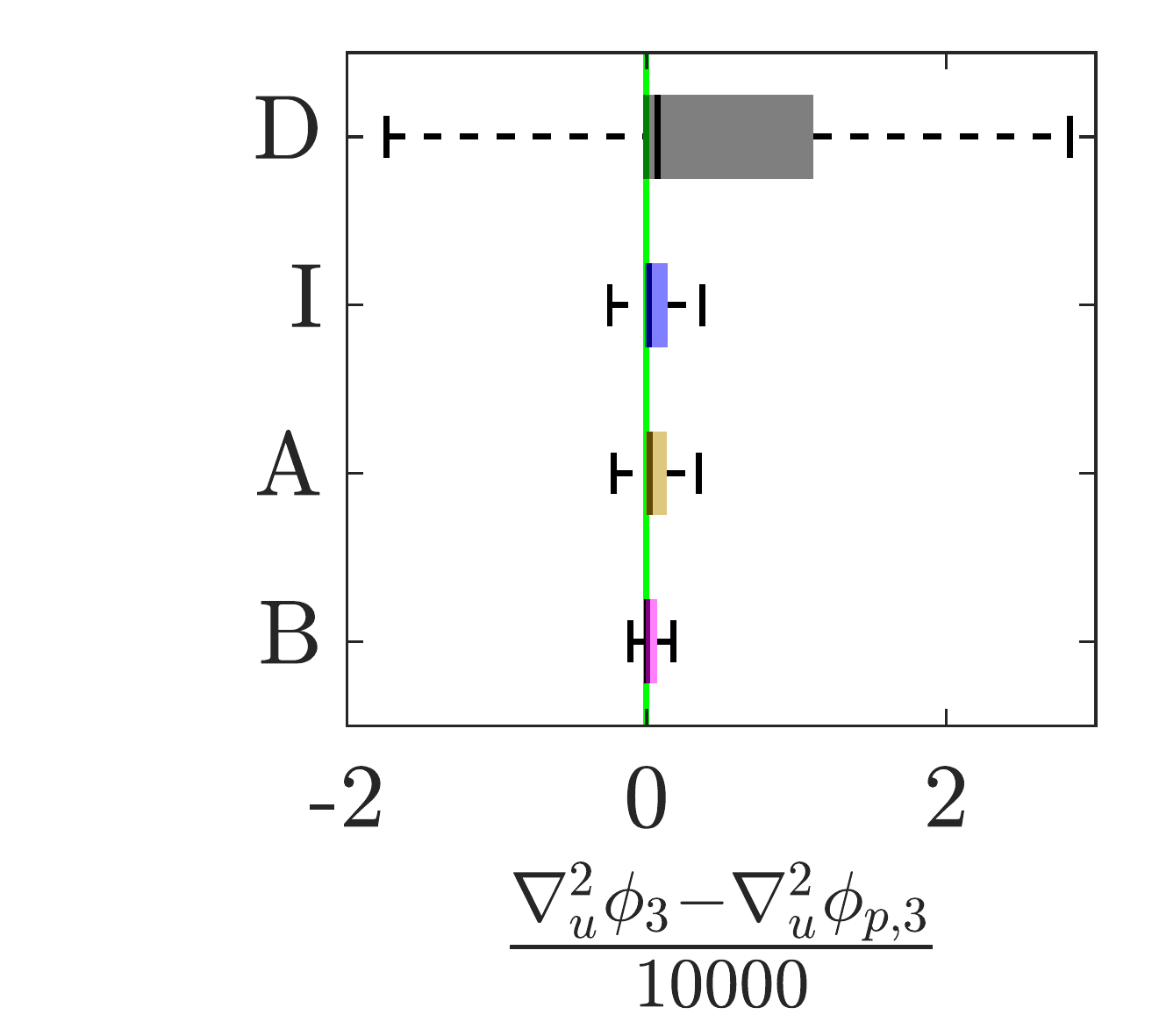}\hskip -0ex
				\includegraphics[trim={2.75cm 0.2cm 0.2cm  0.2cm },clip,width=3.45cm]{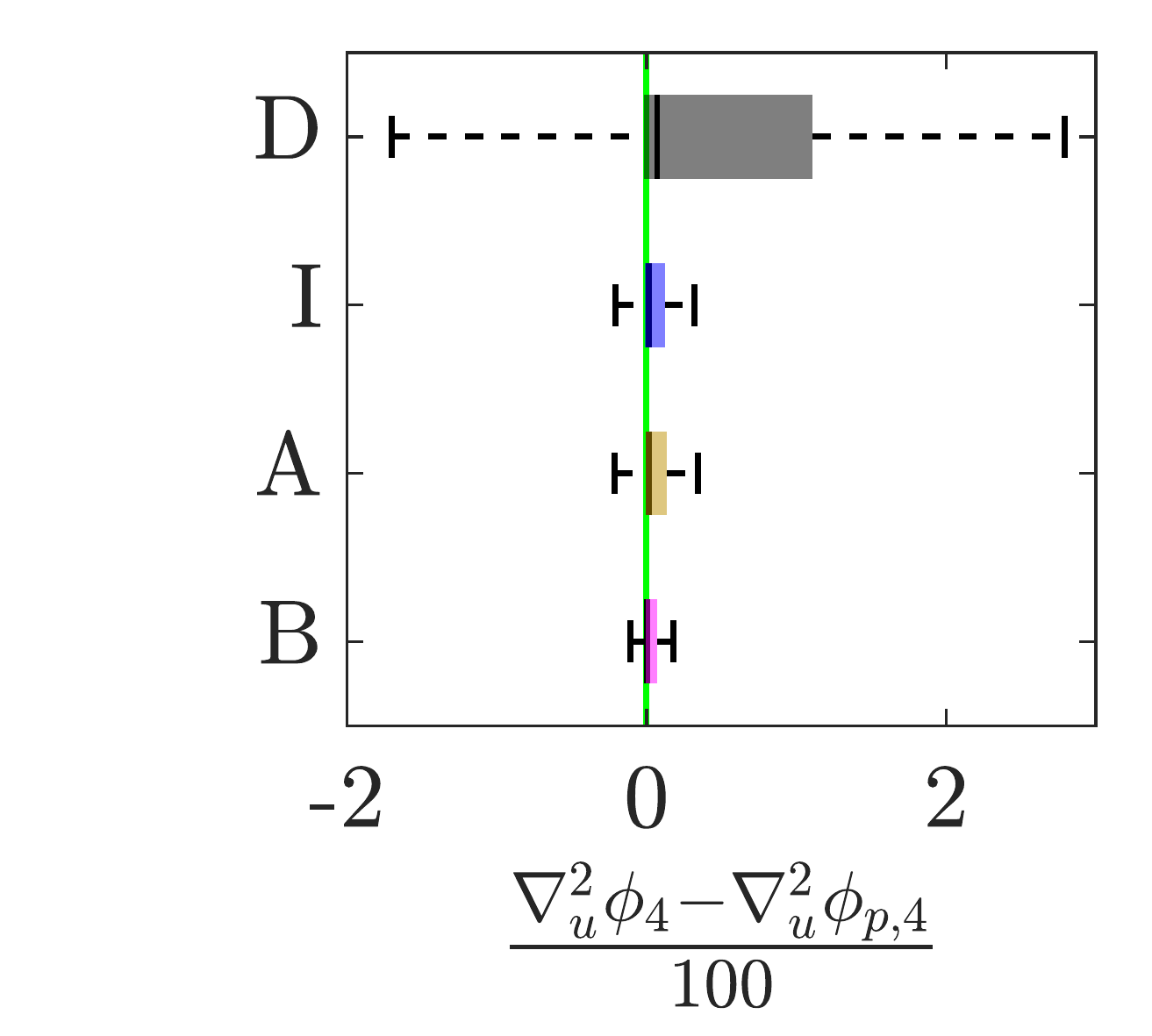}
			}
		\vspace{-2mm}
		\captionof{figure}{Sc.3: Statistical distributions of the prediction errors on the  Hessian of the plant's cost functions at the correction point for the structures D, I, A, and B.}
		\label{fig:5_15_Exemple_5_1_Results}
		\end{minipage} \\
		\begin{minipage}[h]{\linewidth}
			\vspace*{0pt}
			{\centering
				\includegraphics[trim={0.7cm 0.2cm 0.2cm  0.2cm },clip,width=3.45cm]{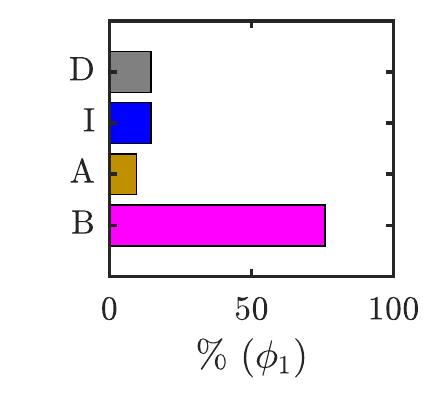}\hskip -0ex
				\includegraphics[trim={0.7cm 0.2cm 0.2cm  0.2cm },clip,width=3.45cm]{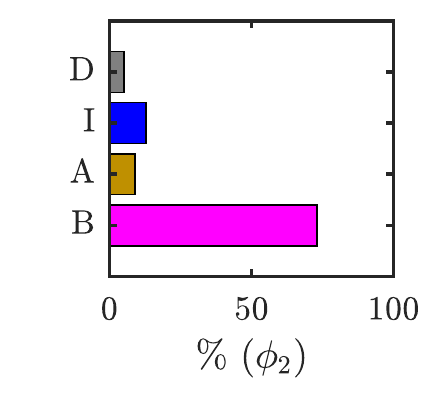}\hskip -0ex
				\includegraphics[trim={0.7cm 0.2cm 0.2cm  0.2cm },clip,width=3.45cm]{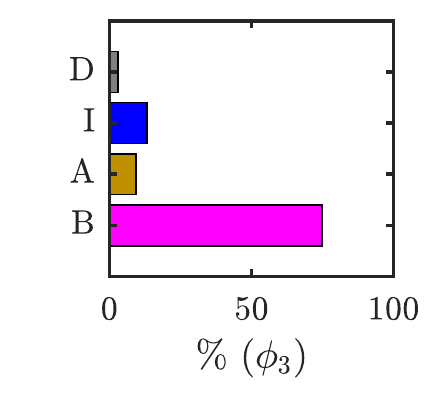}\hskip -0ex
				\includegraphics[trim={0.7cm 0.2cm 0.2cm  0.2cm },clip,width=3.45cm]{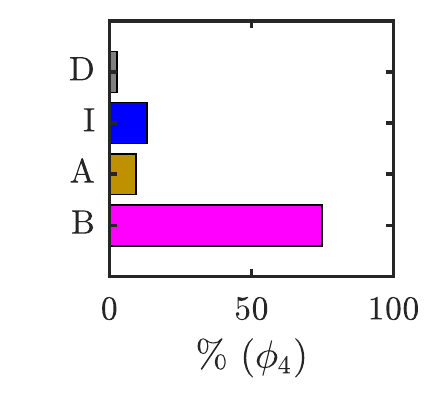} \\
			}
			\vspace{-2mm}
			\captionof{figure}{Sc.3: For each correction structure one gives here the percentage of cases for which no other structure provides better results  \textit{(if two structures provide the same best result then both take the point)}}
			\label{fig:5_16_Exemple_5_1_Results_2}
		\end{minipage} \\
	\end{minipage} 
\noindent
\begin{minipage}[h]{\linewidth}
	\vspace*{0pt}
	{\centering
		\vspace{-5mm}
		\includegraphics[width=3.45cm]{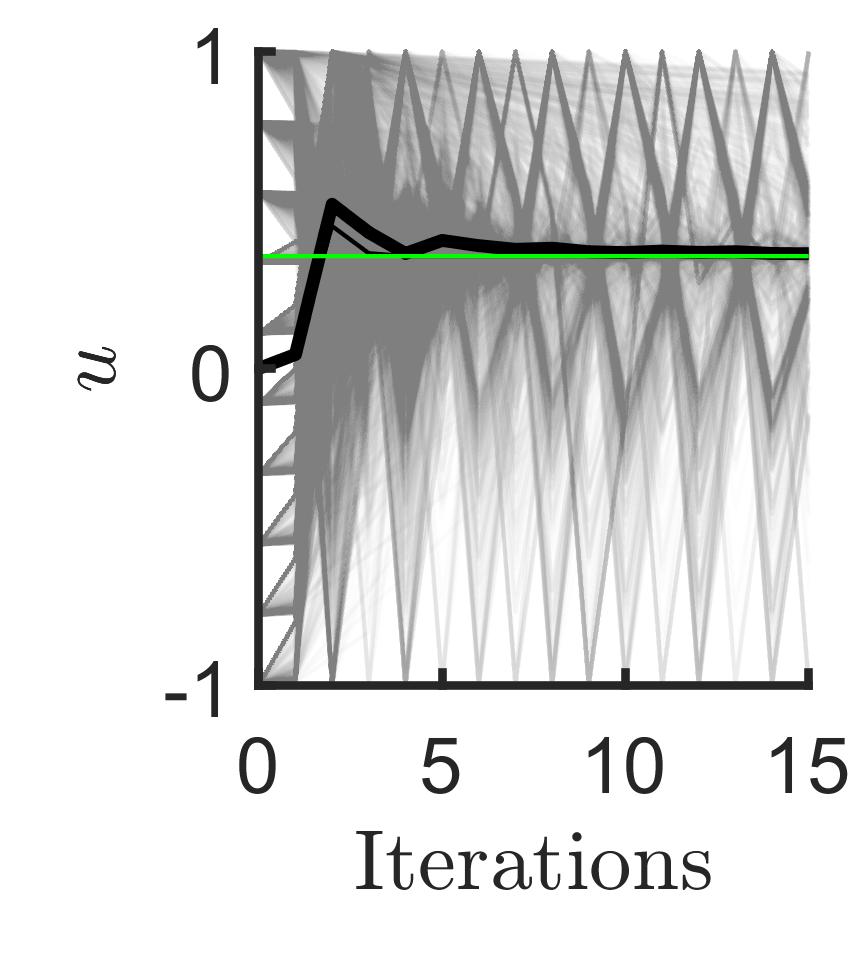}\hskip -0ex
		\includegraphics[width=3.45cm]{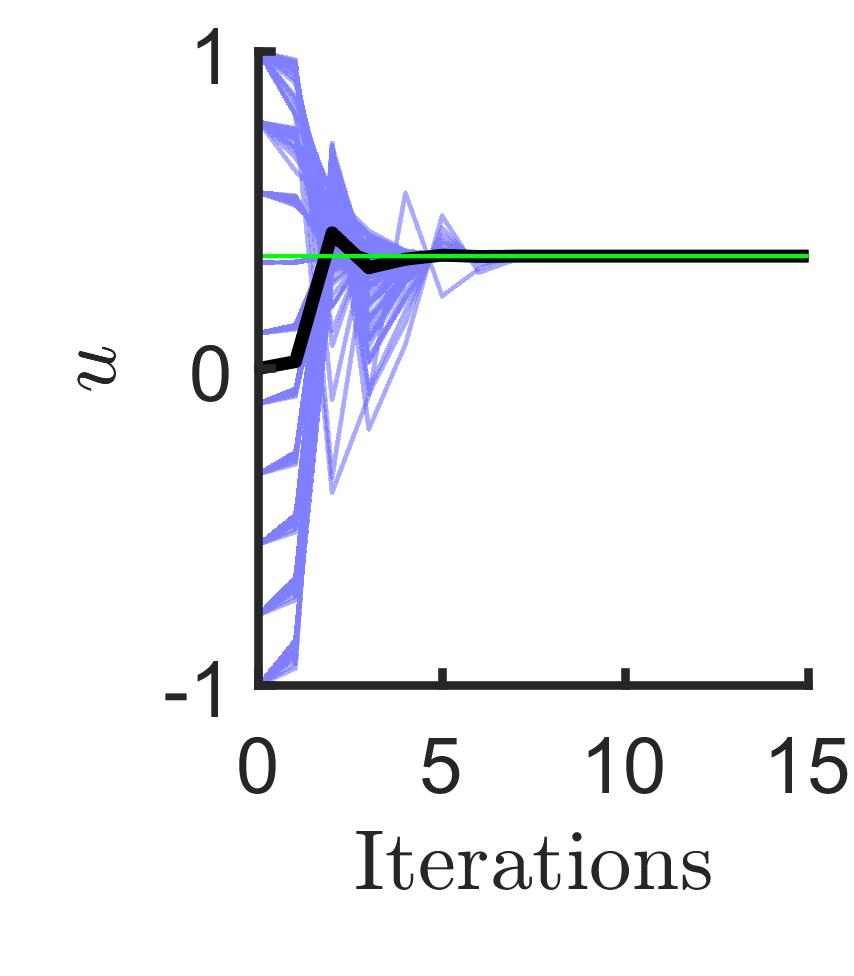}\hskip -0ex
		\includegraphics[width=3.45cm]{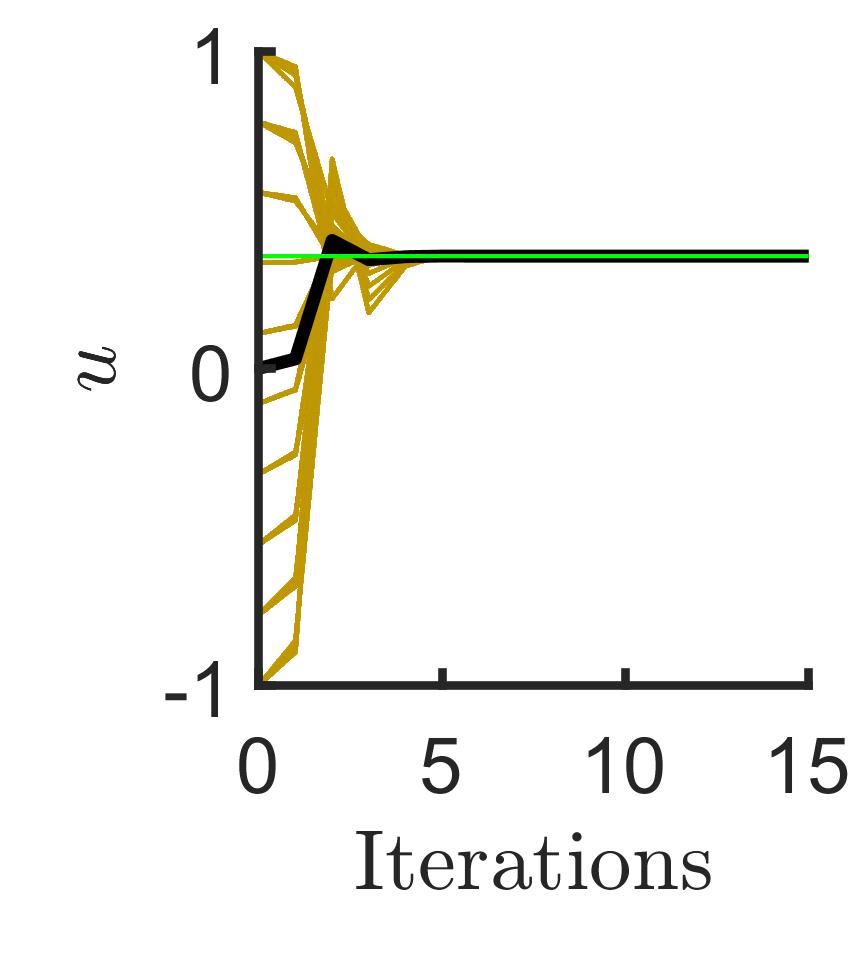}\hskip -0ex
		\includegraphics[width=3.45cm]{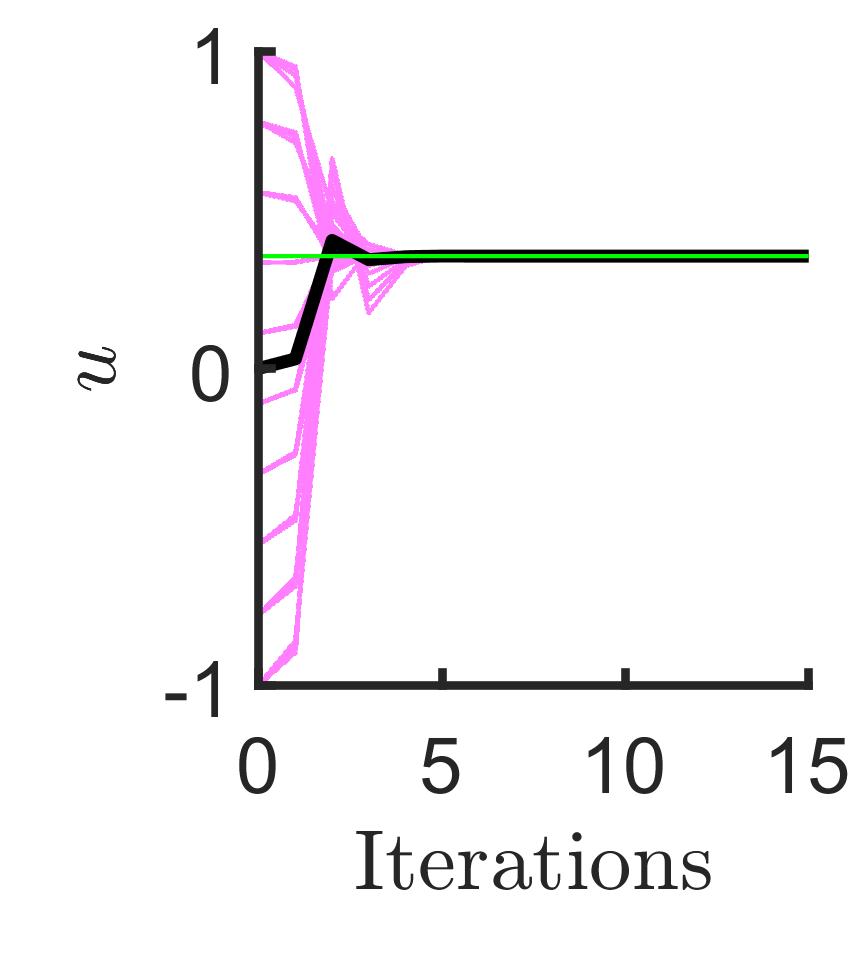} \\
		\vspace{-1.1cm}
		\includegraphics[width=3.45cm]{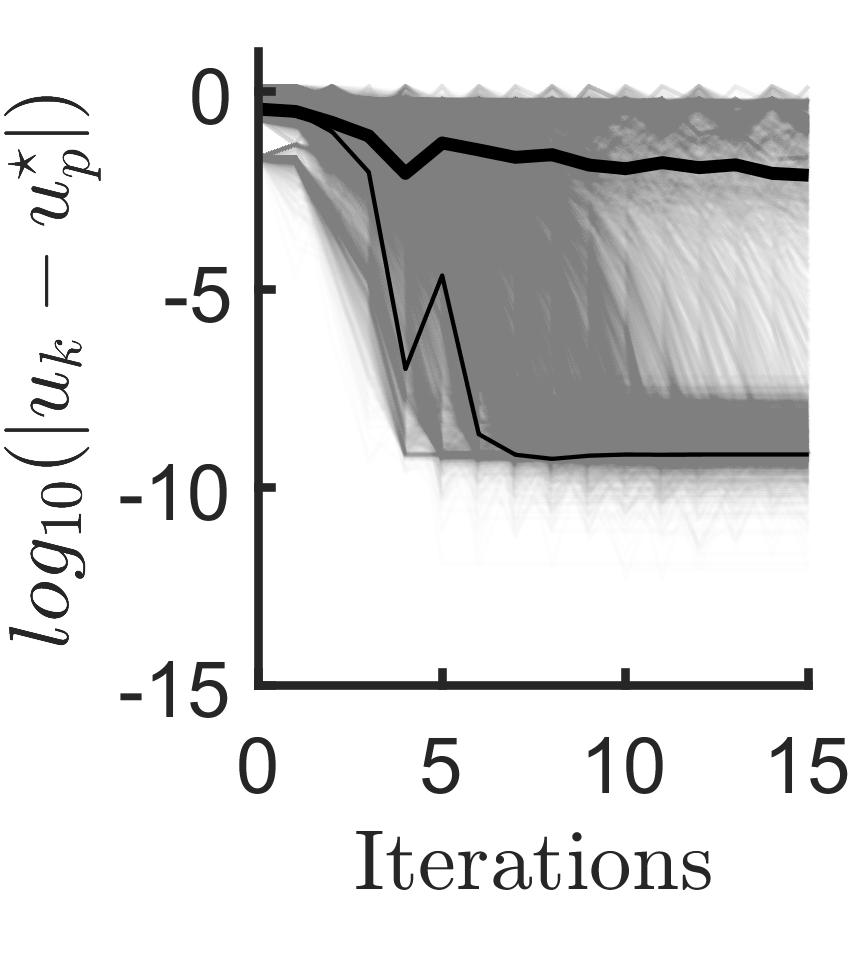}\hskip -0ex
		\includegraphics[width=3.45cm]{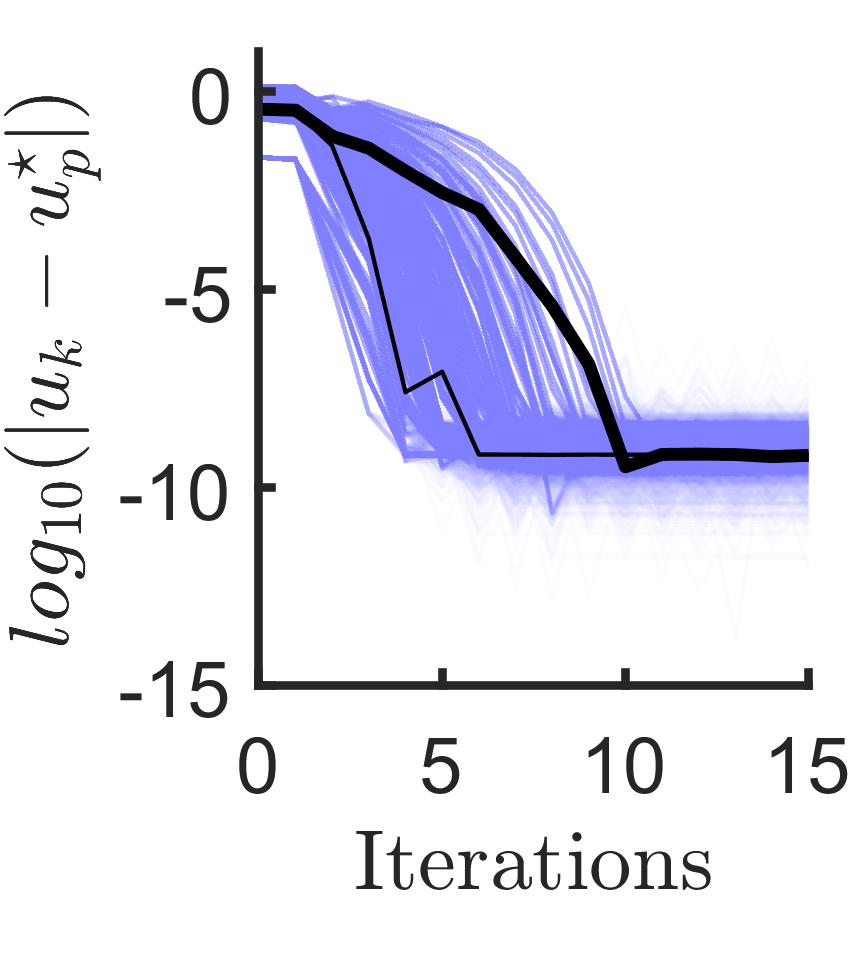}\hskip -0ex
		\includegraphics[width=3.45cm]{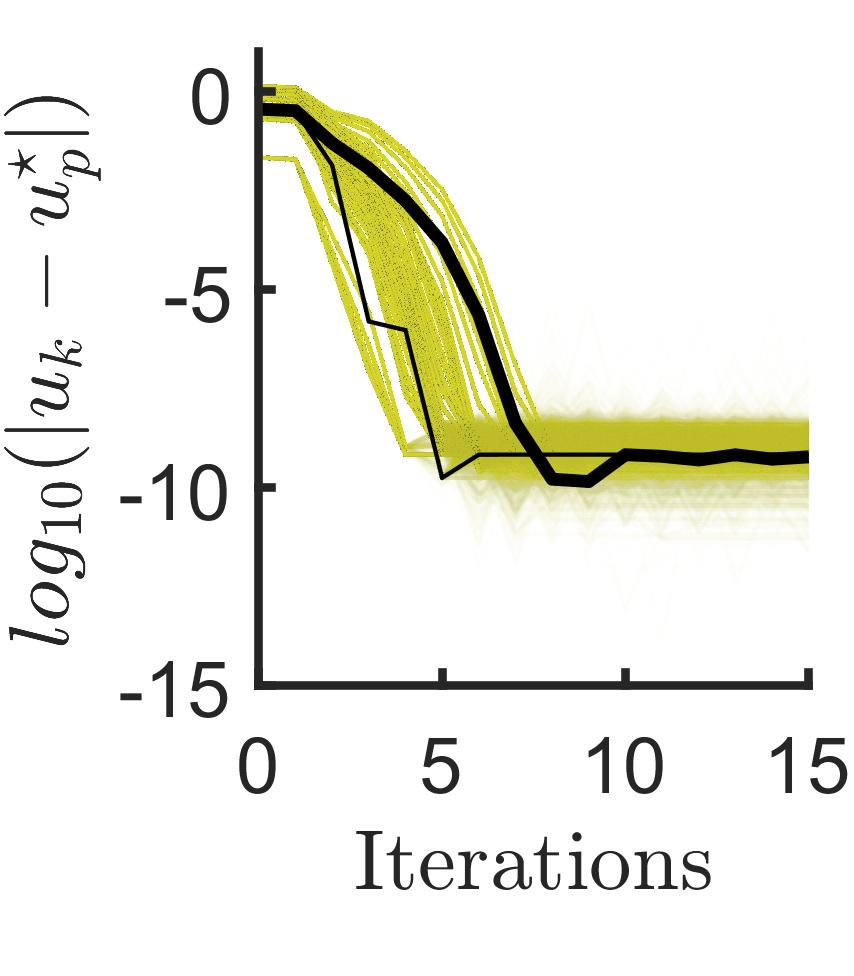}\hskip -0ex
		\includegraphics[width=3.45cm]{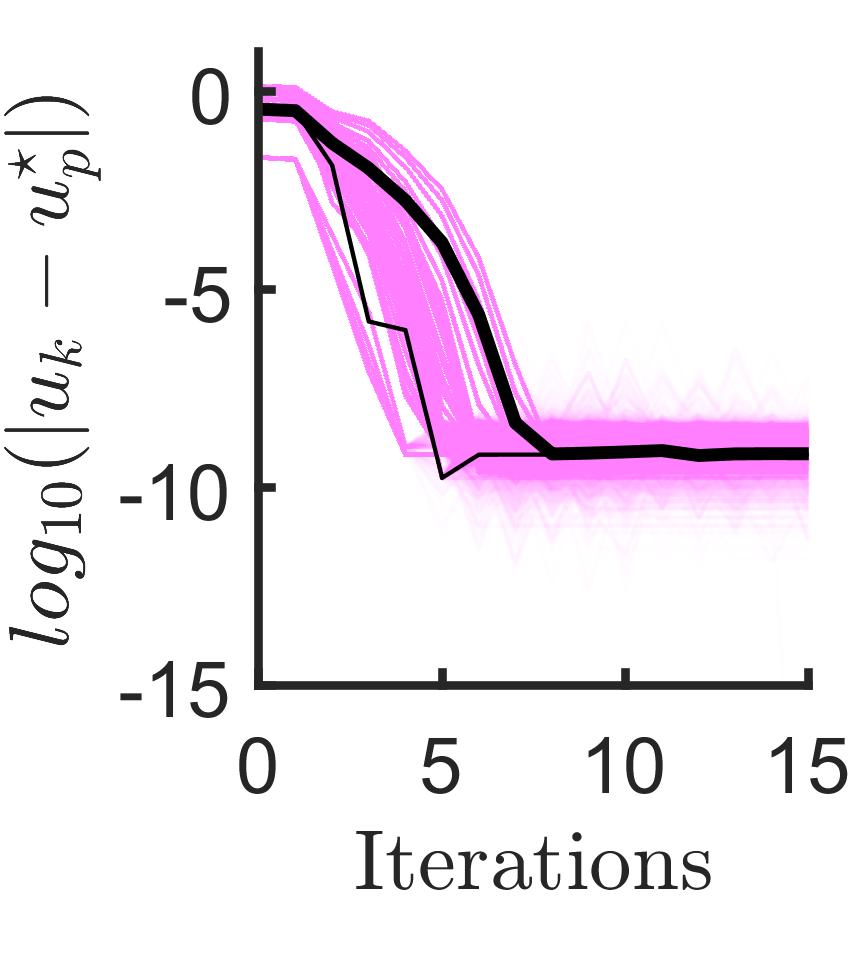}\\
	}
	\vspace{-4mm}
	a) \textbf{Sc.1} \textit{(Structures  D (\textcolor{gris_clair}{\raisebox{0.5mm}{\rule{0.3cm}{0.05cm}}}), I (\textcolor{blue}{\raisebox{0.5mm}{\rule{0.3cm}{0.05cm}}}), 
		A (\textcolor{gold}{\raisebox{0.5mm}{\rule{0.3cm}{0.05cm}}}),
		B (\textcolor{magenta}{\raisebox{0.5mm}{\rule{0.3cm}{0.05cm}}}), 
		mean  (\textcolor{black}{\raisebox{0.5mm}{\rule{0.3cm}{0.1cm}}}),
		median (\textcolor{black}{\raisebox{0.5mm}{\rule{0.3cm}{0.05cm}}}))}\\
	{\centering
		\includegraphics[width=3.45cm]{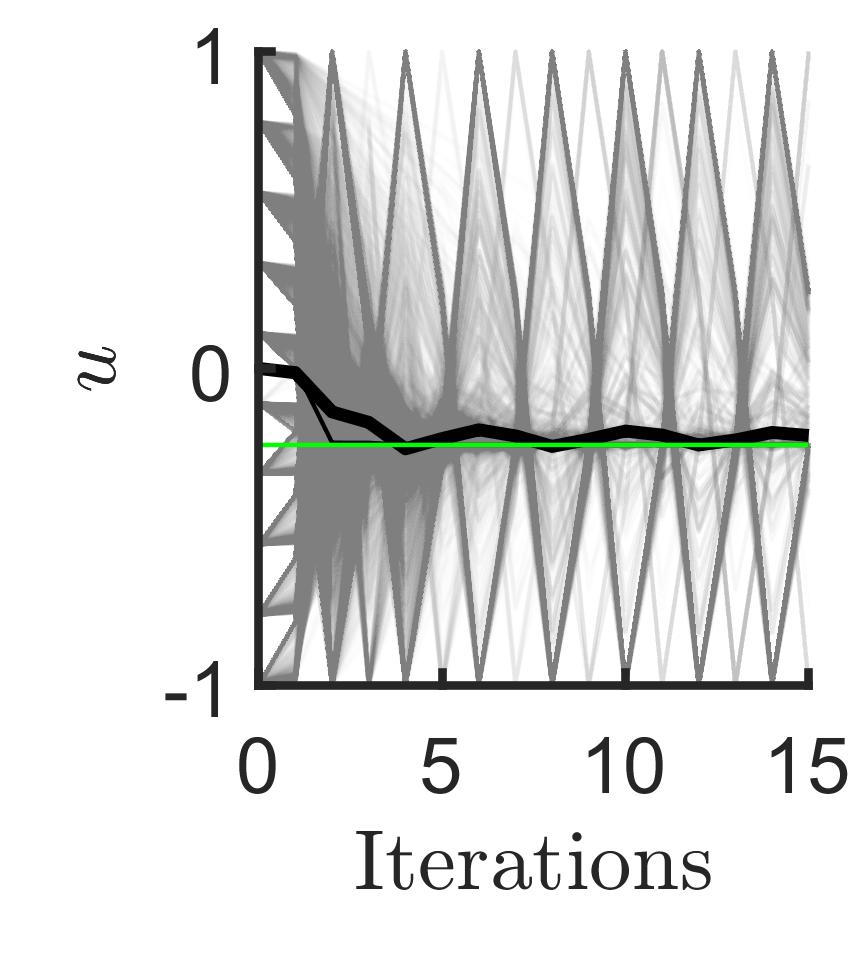}\hskip -0ex
		\includegraphics[width=3.45cm]{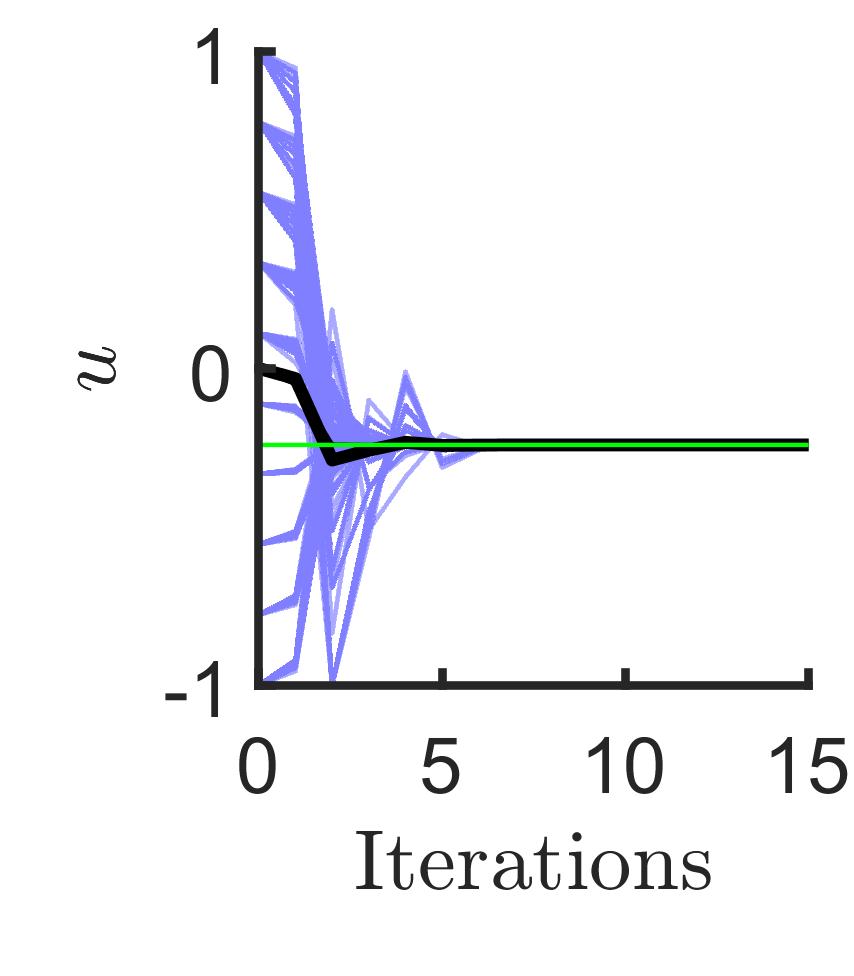}\hskip -0ex
		\includegraphics[width=3.45cm]{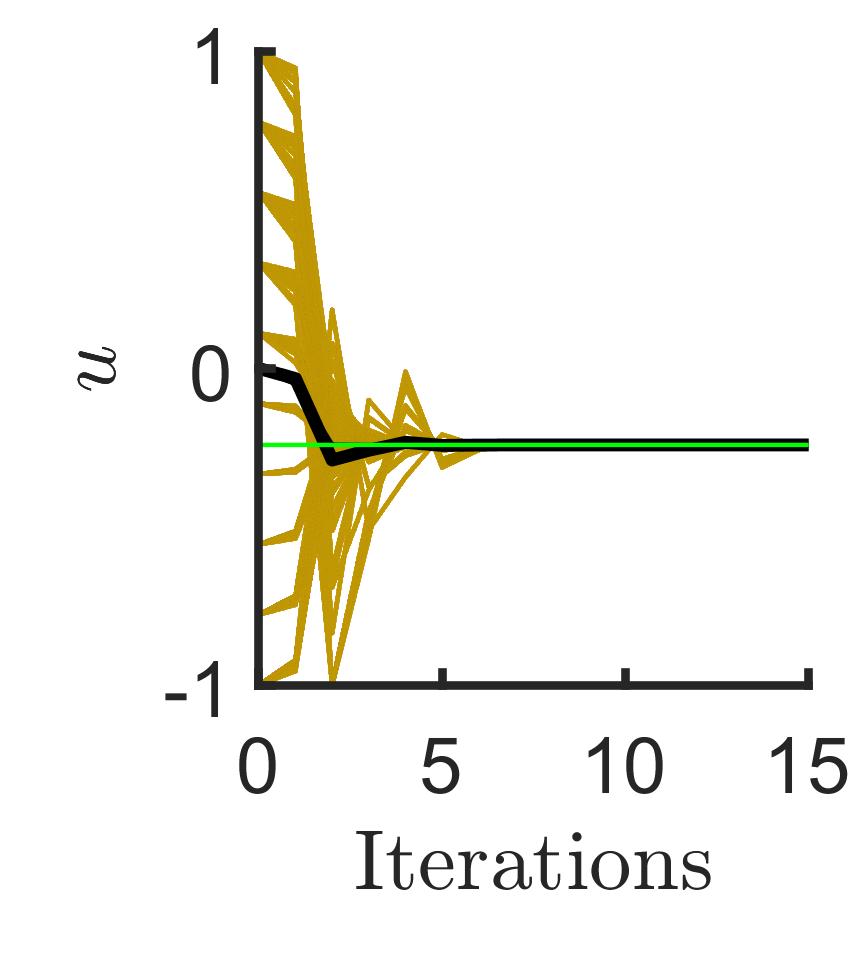}\hskip -0ex
		\includegraphics[width=3.45cm]{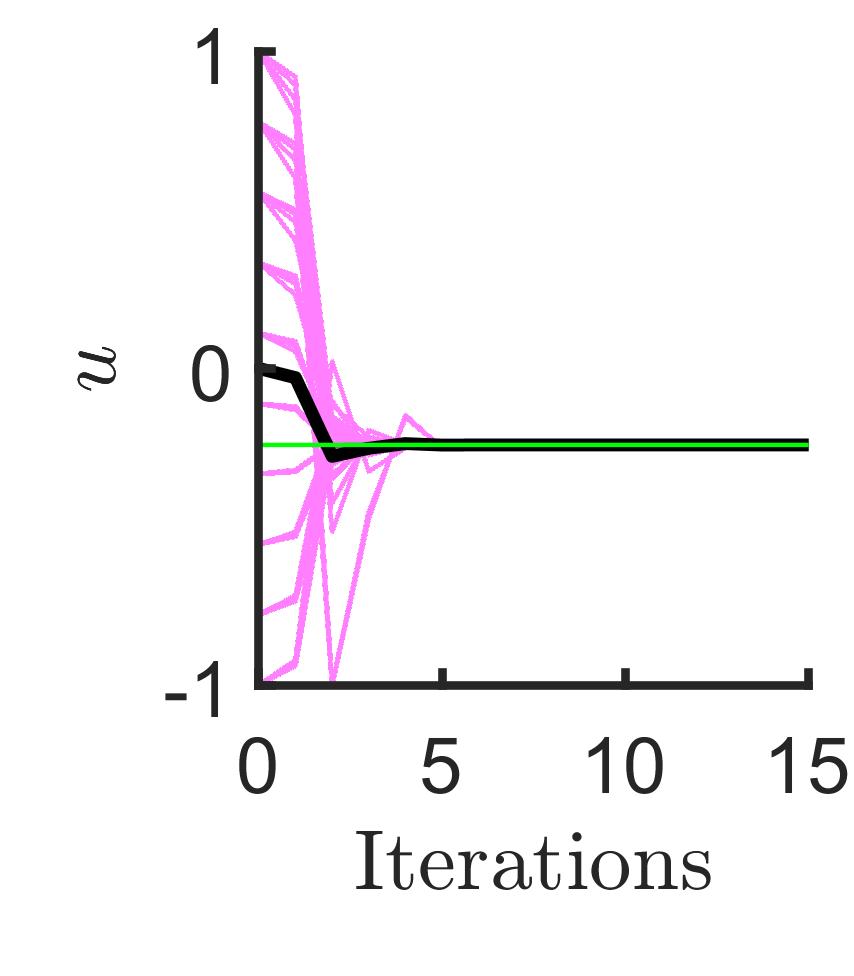} \\
		\vspace{-1.1cm}
		\includegraphics[width=3.45cm]{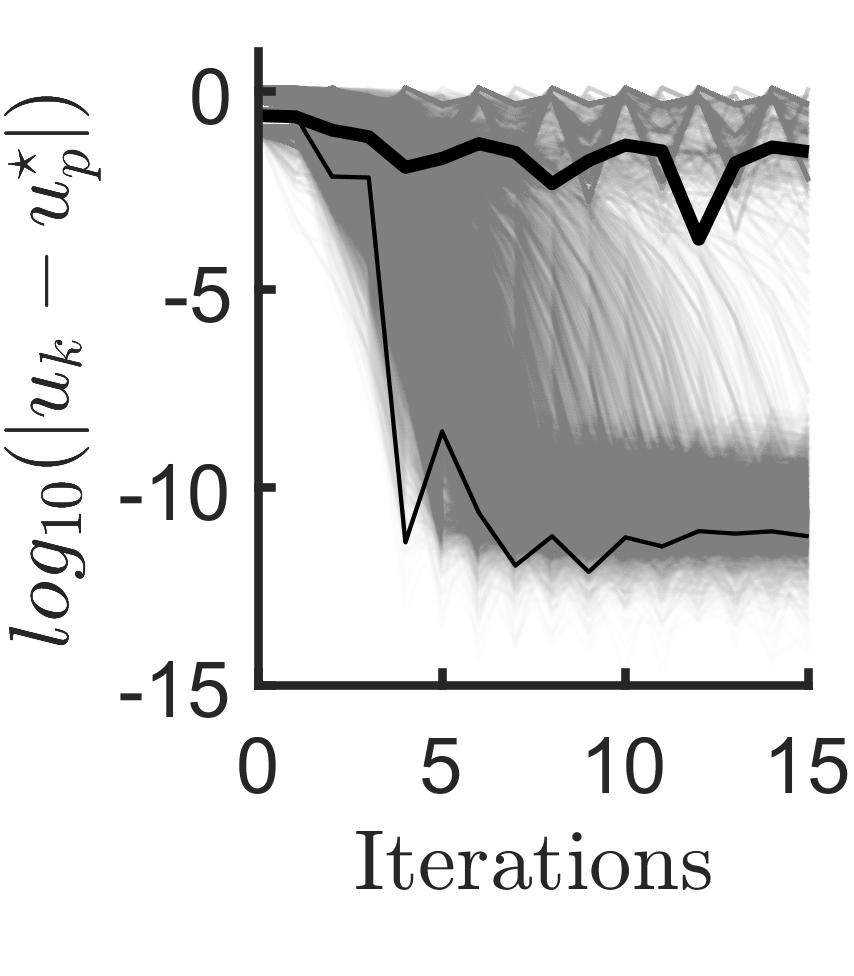}\hskip -0ex
		\includegraphics[width=3.45cm]{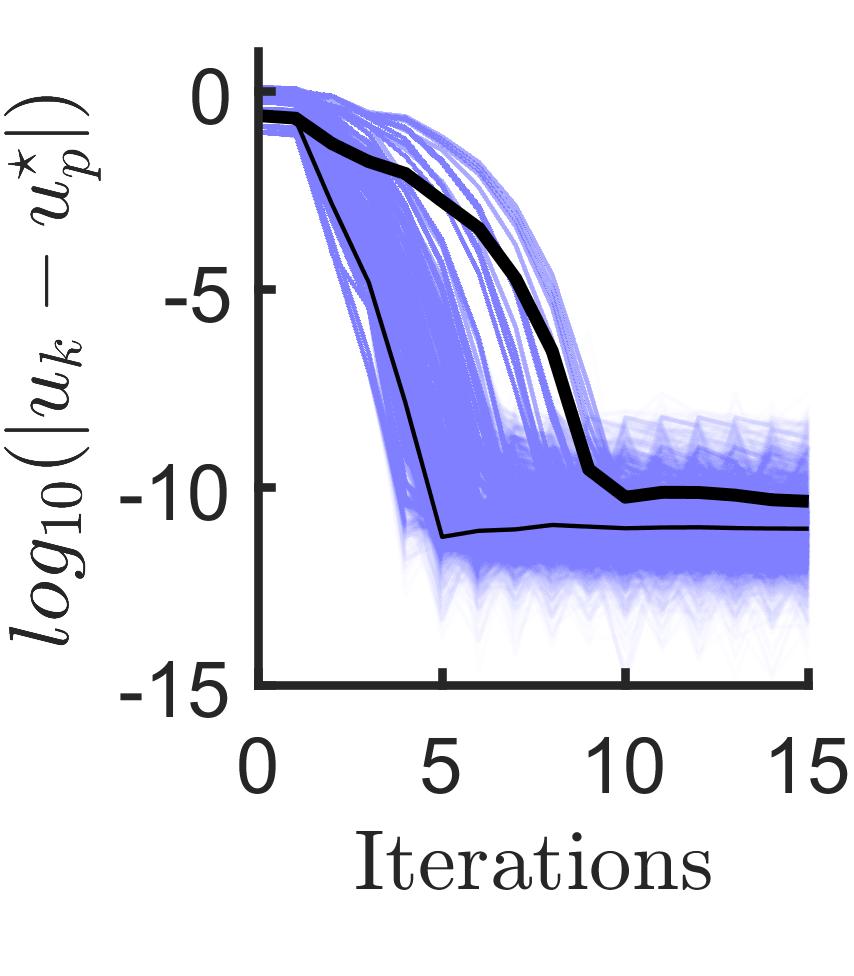}\hskip -0ex
		\includegraphics[width=3.45cm]{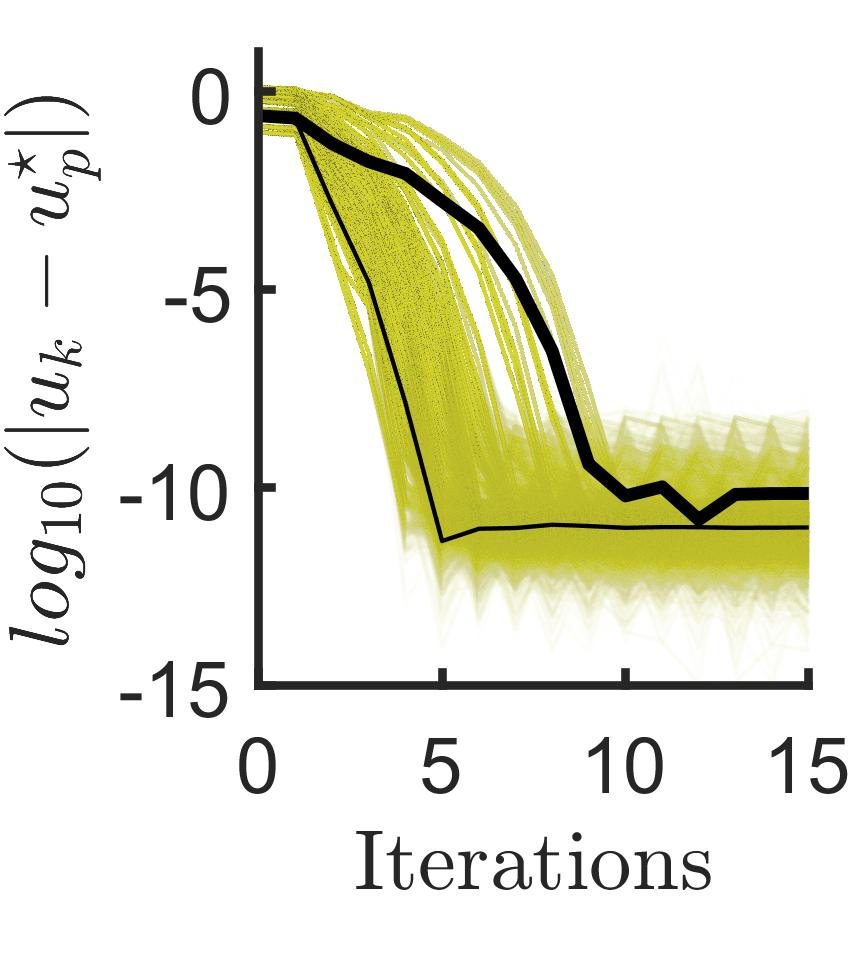}\hskip -0ex
		\includegraphics[width=3.45cm]{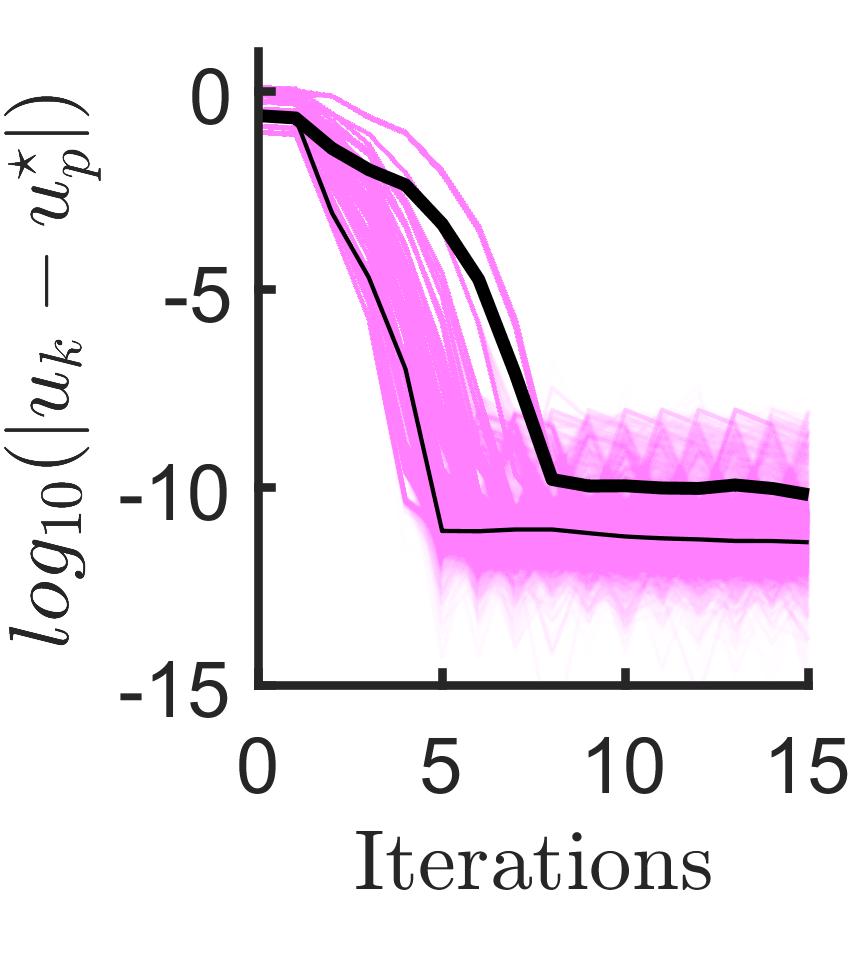}\\
	}
	\vspace{-4mm}
		b) \textbf{Sc.2} \textit{(Structures  D (\textcolor{gris_clair}{\raisebox{0.5mm}{\rule{0.3cm}{0.05cm}}}), I (\textcolor{blue}{\raisebox{0.5mm}{\rule{0.3cm}{0.05cm}}}), 
			A (\textcolor{gold}{\raisebox{0.5mm}{\rule{0.3cm}{0.05cm}}}),
			B (\textcolor{magenta}{\raisebox{0.5mm}{\rule{0.3cm}{0.05cm}}}), 
			mean  (\textcolor{black}{\raisebox{0.5mm}{\rule{0.3cm}{0.1cm}}}),
			median (\textcolor{black}{\raisebox{0.5mm}{\rule{0.3cm}{0.05cm}}}))}\\
	{\centering
		\includegraphics[width=3.45cm]{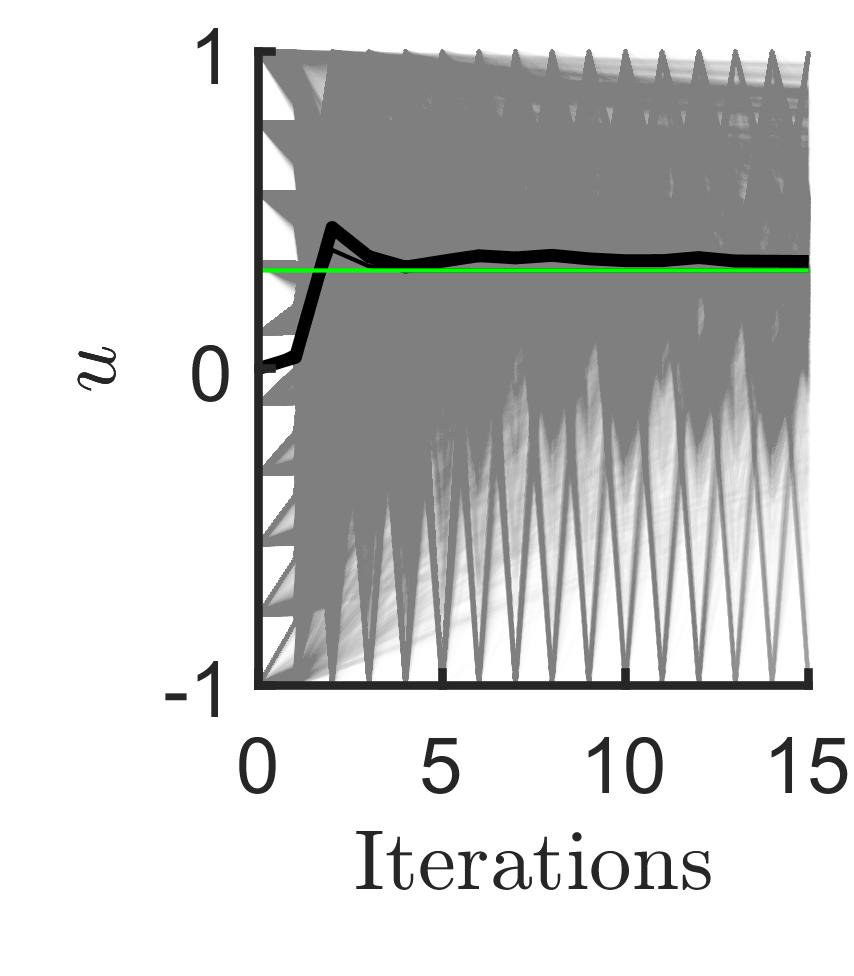}\hskip -0ex
		\includegraphics[width=3.45cm]{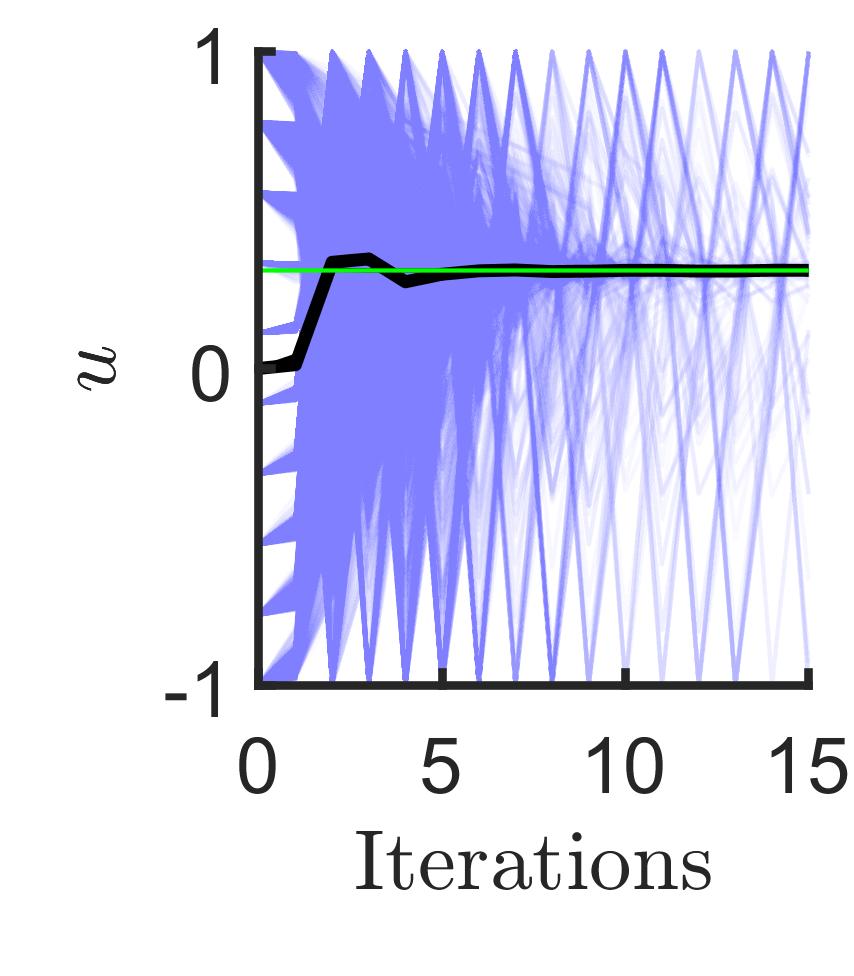}\hskip -0ex
		\includegraphics[width=3.45cm]{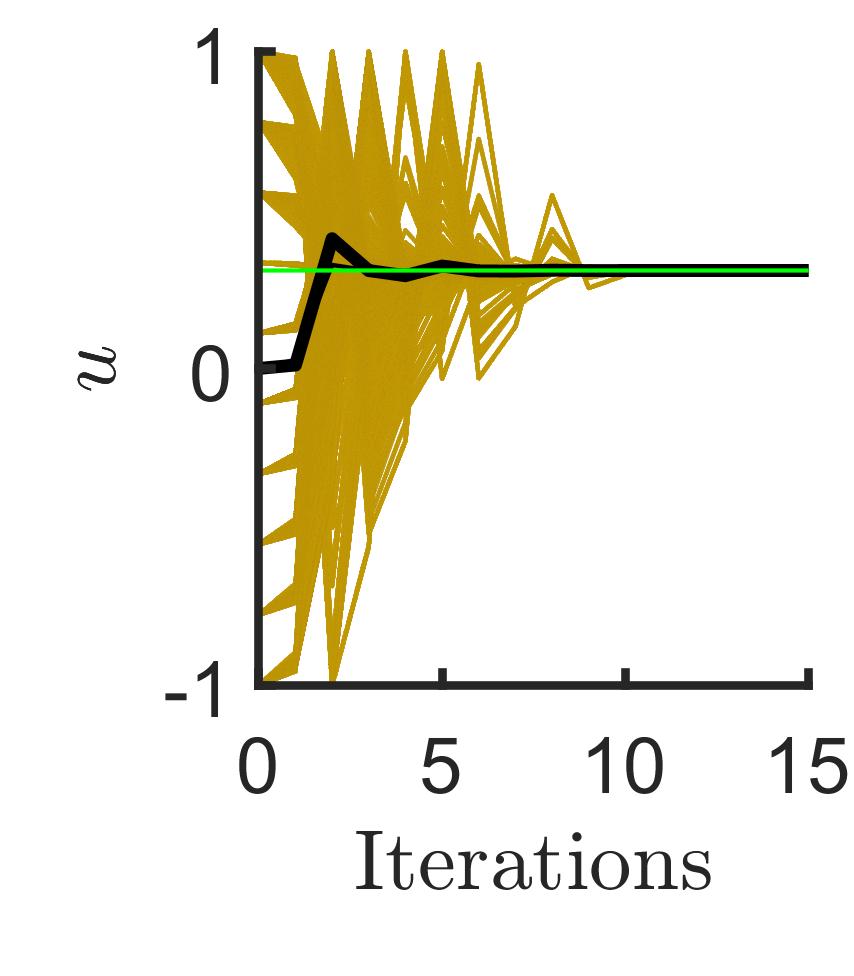}\hskip -0ex
		\includegraphics[width=3.45cm]{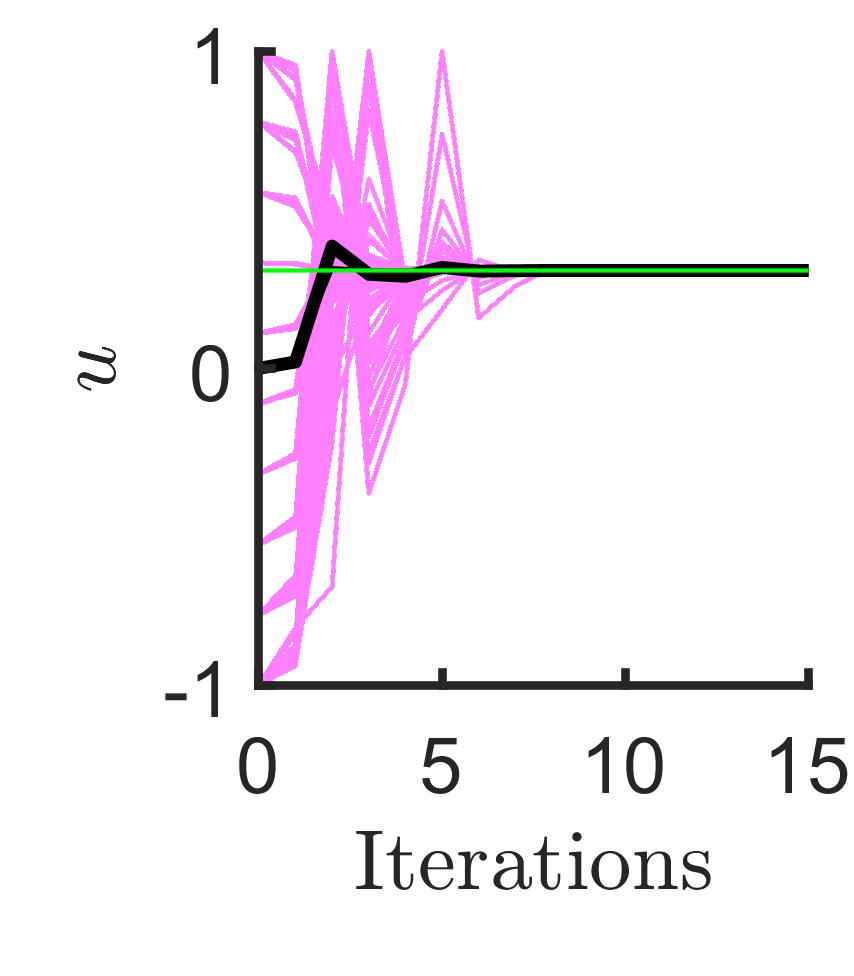} \\
		\vspace{-1.1cm}
		\includegraphics[width=3.45cm]{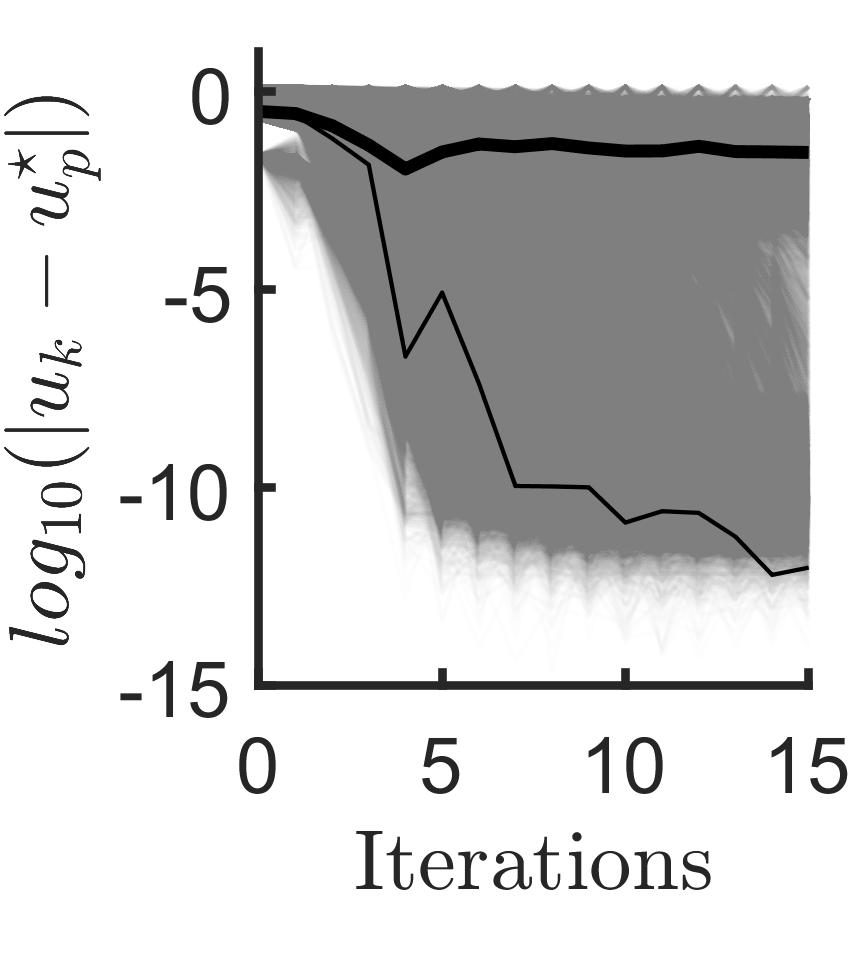}\hskip -0ex
		\includegraphics[width=3.45cm]{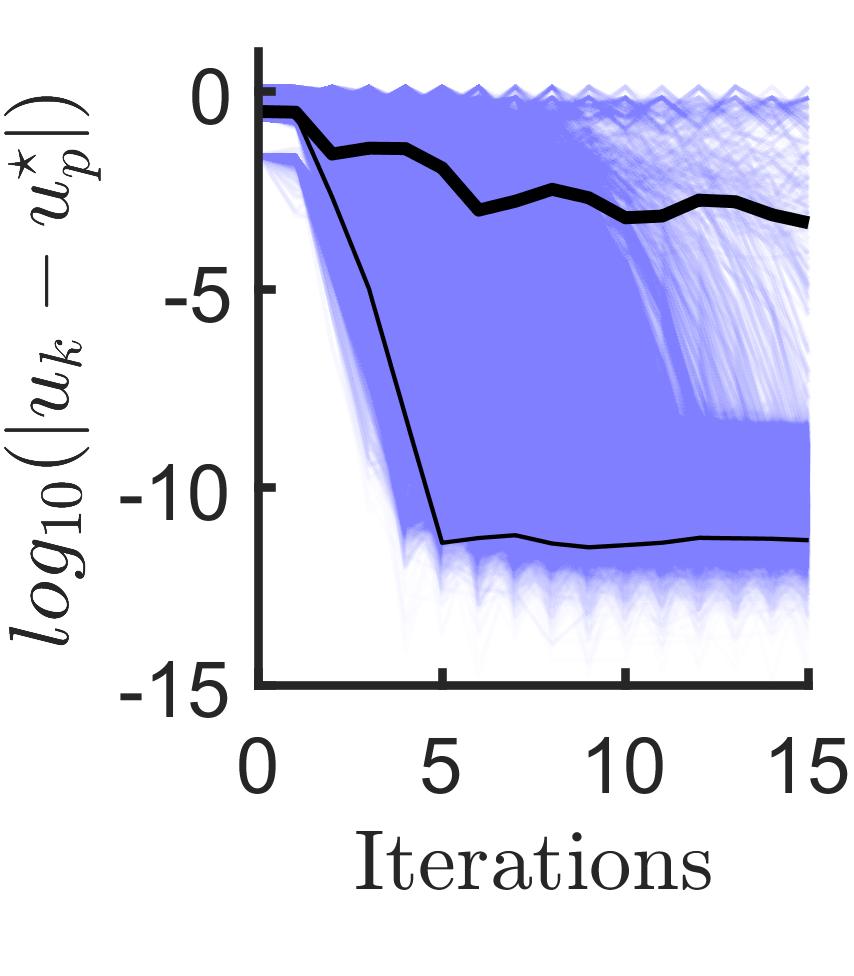}\hskip -0ex
		\includegraphics[width=3.45cm]{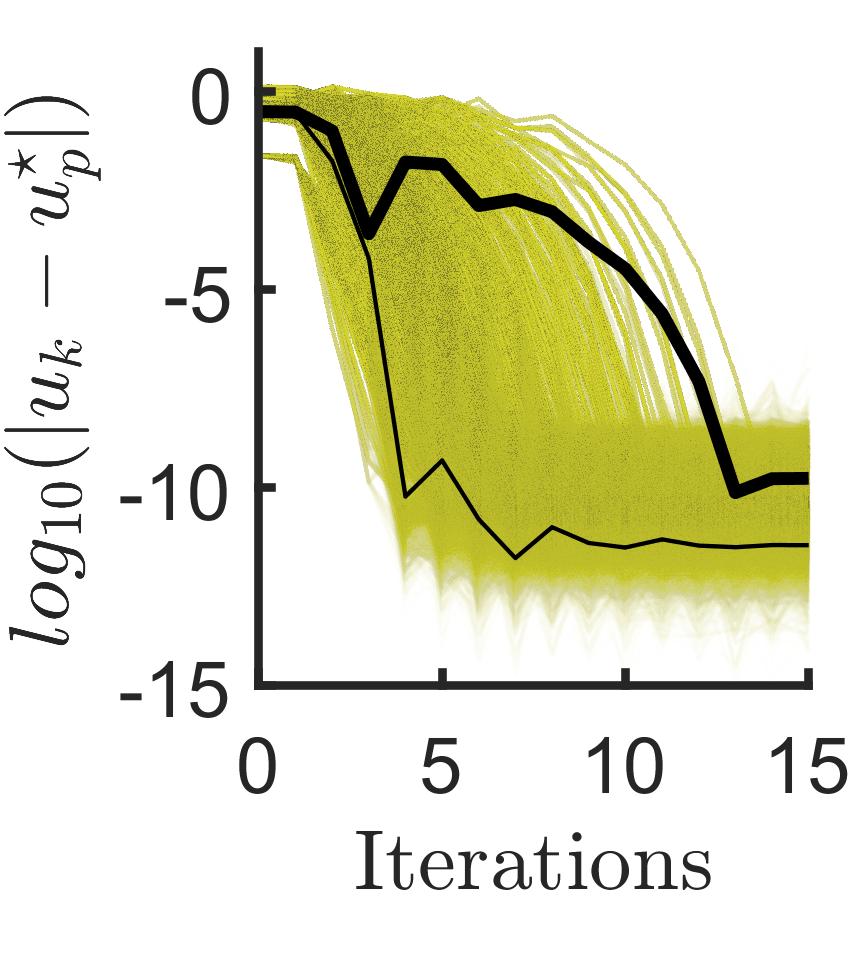}\hskip -0ex
		\includegraphics[width=3.45cm]{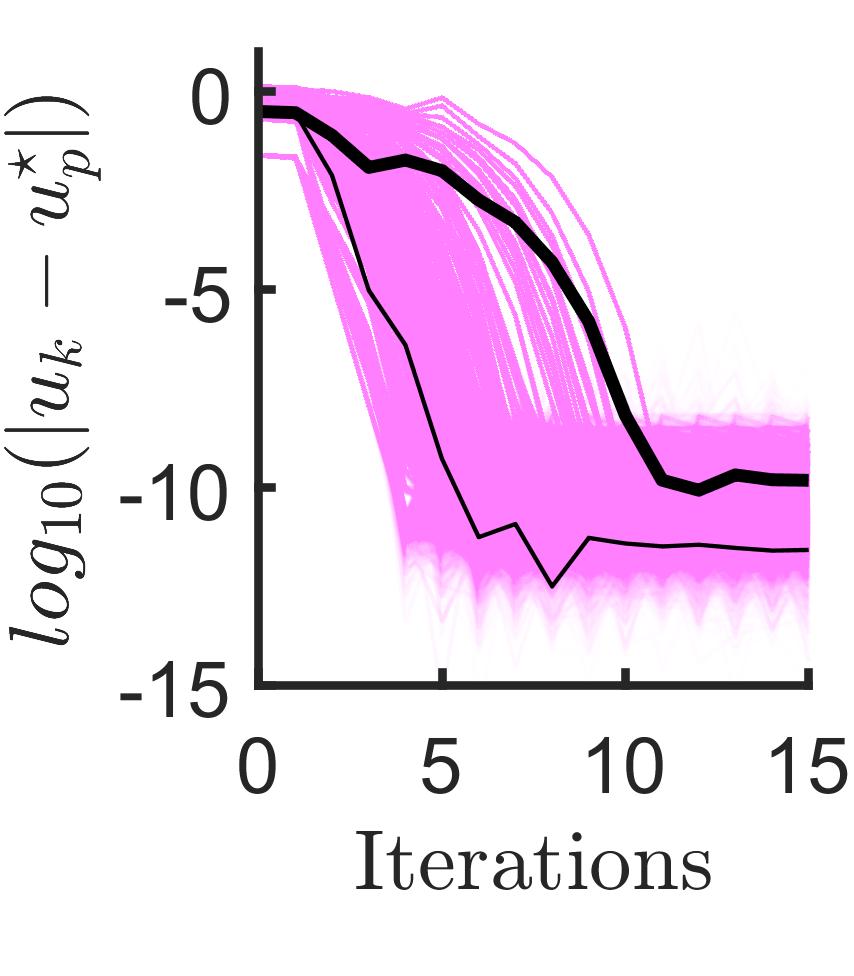}\\
	}
	\vspace{-4mm}	
	c) \textbf{Sc.3} \textit{(Structures  D (\textcolor{gris_clair}{\raisebox{0.5mm}{\rule{0.3cm}{0.05cm}}}), I (\textcolor{blue}{\raisebox{0.5mm}{\rule{0.3cm}{0.05cm}}}), 
		A (\textcolor{gold}{\raisebox{0.5mm}{\rule{0.3cm}{0.05cm}}}),
		B (\textcolor{magenta}{\raisebox{0.5mm}{\rule{0.3cm}{0.05cm}}}), 
		mean  (\textcolor{black}{\raisebox{0.5mm}{\rule{0.3cm}{0.1cm}}}),
		median (\textcolor{black}{\raisebox{0.5mm}{\rule{0.3cm}{0.05cm}}}))}
	\captionof{figure}{Simulation results for the structures D, I, A and B }
	\label{fig:5_17_Exemple_5_1_Results_3}
\end{minipage}
 
\ \\

\subsection{Study 2:  Systems in closed loop \& Consistent model}

\noindent
\begin{minipage}[h]{\linewidth}
	\includegraphics[width=14cm]{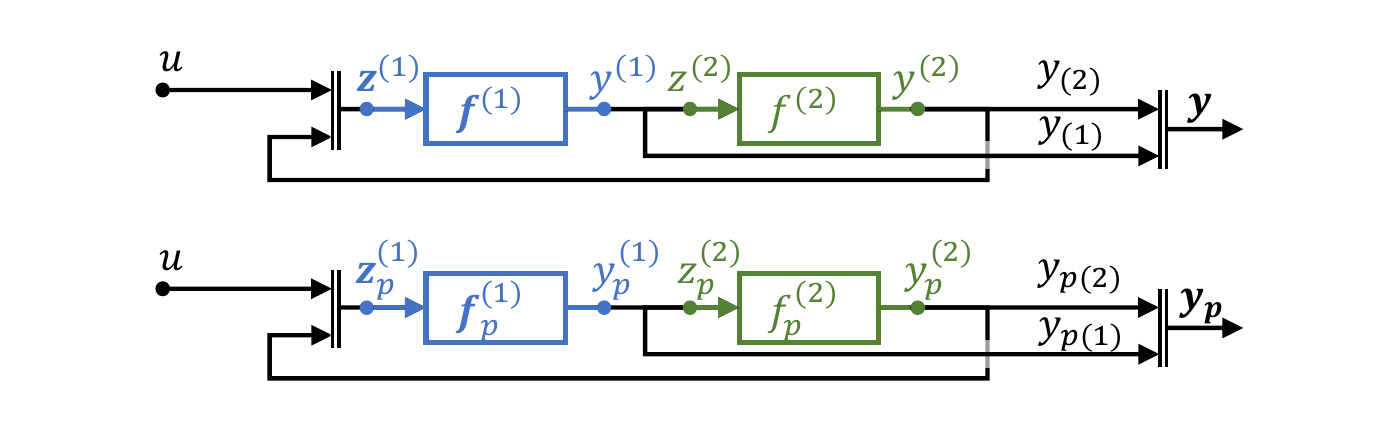}
	\captionof{figure}{Graphical description of the RTO problem of Study 2}
	\label{fig:5___19_Systemes_avec_retroaction}
\end{minipage}\\

One considers the theoretical RTO problem illustrated in Figure~\ref{fig:5___19_Systemes_avec_retroaction} and using the following functions: 
\begin{align*}
	f^{(1)} := \ &\theta_1 + 
	\left(\begin{array}{c}
		\theta_2 \\ 0.01
	\end{array}\right)^{\rm T}
	\bm{z}^{(1)} + \theta_3 
	\bm{z}^{(1)\rm T}
	\left(\begin{array}{cc}
		0.01 &  \theta_4 \\ \theta_4 & 0.002
	\end{array}\right)
	\bm{z}^{(1)}, \\
	f_p^{(1)} := \ & 1.5 + 
	\left(\begin{array}{c}
		-1 \\ 0.01
	\end{array}\right)^{\rm T}
	\bm{z}_p^{(1)} + 
	\theta_{p1} 
	\bm{z}_p^{(1)\rm T}
	\left(\begin{array}{cc}
		0.01 & 0.08\\ 0.08 & 0.002
	\end{array}\right)
	\bm{z}_p^{(1)},
\end{align*}
\vspace{-7mm}
\begin{align*}
	\bm{z}^{(1)} := \ & 
	\left(\begin{array}{c}
		u \\ y^{(2)}
	\end{array}\right),& 
	\bm{z}_p^{(1)} := \ & 
	\left(\begin{array}{c}
		u \\ y_p^{(2)}
	\end{array}\right), \\
	f^{(2)} := \ & \theta_5 + \theta_6z^{(2)} + \theta_7z^{(2)2}, & f_p^{(2)} := \ & 2 -2z_p^{(2)} + \theta_{p2}z_p^{(2)2}, \\
	z^{(2)} :=\ & y^{(1)}, & z_p^{(2)} := \ & y_p^{(1)}, \\
	\bm{y}  :=\ & \left[y^{(1)};y^{(2)}\right]^{\rm T}
\end{align*}
\begin{align*}
\phi_1(u,\bm{y}) := \ & y_{(2)},   & \phi_2(u,\bm{y}) := \ & uy_{(2)}, \\
\phi_3(u,\bm{y}) := \ & y_{(2)}^2, & \phi_4(u,\bm{y}) := \ & y_{(2)} + u^2 + 0.5uy_{(2)} + 0.05y_{(2)}^2.
\end{align*}
where $u\in[-1,1]$. 
The model is parameterized with the parameters $\{\theta_1,...,\theta_7\}$ and the plant with the parameters $\{\theta_{p1},\theta_{p2}\}$. In this example  3 scenarios are considered: 
\begin{itemize}
	\item \textit{Scenario 1 (A$\rightleftarrows$NL):}  The functions $\{\bm{f}^{(1)},\bm{f}_p^{(1)}\}$ are affine (A) and the functions $\{f^{(2)},f_p^{(2)}\}$ are non-linear (NL). 
	\item \textit{Scenario 2 (NL$\rightleftarrows$A):}  The functions $\{\bm{f}^{(1)},\bm{f}_p^{(1)}\}$ are NL and the functions $\{f^{(2)},f_p^{(2)}\}$ are A. 
	\item \textit{Scenario 3 (NL$\rightleftarrows$NL):}  The functions $\{\bm{f}^{(1)},\bm{f}_p^{(1)},f^{(2)},f_p^{(2)}\}$ are NL . 
\end{itemize}
The model and plant parameters considered for each of these scenarios are given in Table~\ref{tab:5_2_Exemple_5_2_Parametres_Modeles_Usine}. The Figures~\ref{fig:5_20_Exemple_5_2_sc1_Plant_and_Model}, \ref{fig:5_23_Exemple_5_2_sc2_Plant_and_Model} and \ref{fig:5_26_Exemple_5_2_sc3_Plant_and_Model} show the functions of the plant and of  the set of models that are  used. \\

\noindent
\begin{minipage}[h]{\linewidth}
	\centering
	\begin{tabular}{ccccccccccc}
		\toprule
		Scenario  & $\theta_1$          & $\theta_2$          & $\theta_3$   & $\theta_4$        
		          & $\theta_5$          & $\theta_6$          & $\theta_7$  
		          & $\theta_{p1}$ & $\theta_{p2}$ \\
				\midrule
				1 & $\mathbb{\Theta}$   & $-\mathbb{\Theta}$ & $0$ & $0$
				  & $1+\mathbb{\Theta}$ & $-1-\mathbb{\Theta}$  & $1+\mathbb{\Theta}$ 
				  & $0$ & $2$\\
				2 & $\mathbb{\Theta}$   & $-\mathbb{\Theta}$ & $1$ & $\mathbb{\Theta}^\prime$
			   	  & $1+\mathbb{\Theta}$ & $-1-\mathbb{\Theta}$  & $0$
				  & $1$  & $0$    \\
				3 & $\mathbb{\Theta}$   & $-\mathbb{\Theta}$ & $1$ & $\mathbb{\Theta}^\prime$
				  & $1+\mathbb{\Theta}$ & $-1-\mathbb{\Theta}$  & $1+\mathbb{\Theta}$
				  & $1$  & $2$    \\
		\bottomrule
	\end{tabular}%
	\captionof{table}{Definition of the models and of the plant used for each scenario. One defines: $\mathbb{\Theta}:=\{0.75,0.85, 0.95, 1.05, 1.15, 1.25\}$ and $\mathbb{\Theta}^{\prime}:=\{0.07, 0.074, 0.078, 0.082, 0.086, 0.09\}$.}
	\label{tab:5_2_Exemple_5_2_Parametres_Modeles_Usine}%
\end{minipage} \\

Like in the previous study, these three scenarios represent the majority of cases that can be encountered when two systems are in  closed loop. However, unlike the previous study, one can notice that the ``modeling error'' between the models and the plant are ``smaller''. The reason behind this choice is that when the parameters of a model are further separated from those of the plant the problem formed by the two equations $f^{(1)}$ and $f^{(2)}$ may no longer have a solution. However, in the framework of this study, one wants that for each correction structure, each point of the input space, and each choice of parameters, that this system has a solution.

Like in the previous study, for each scenario the ``local quality'' of the correction associated to each structure is quantified. To this end, one evaluates the distribution of the error on the prediction of the plant's Hessian for each model, correction structure, and cost function at 11 points uniformly distributed in the input space i.e. $u=\{-1,-0.8,...,0.8,1\}$. Then, by applying a statistical analysis to the obtained results, one gets the graphs of Figures~\ref{fig:5_21_Exemple_5_2_sc1_Results_1}, \ref{fig:5_24_Exemple_5_2_sc2_Results_1}, and \ref{fig:5_27_Exemple_5_2_sc3_Results_1}. In addition to these statistical analyses, on Figures~\ref{fig:5_22_Exemple_5_2_sc1_Results_2}, \ref{fig:5_25_Exemple_5_2_sc2_Results_2}, and  Figure~\ref{fig:5_28_Exemple_5_2_sc3_Results_2} one gives for each structure the percentage of cases for which no other correction structure brings better predictions of the Hessian of the plant.

What one can quite easily conclude from all these results is that if SM$i$ is in feedback with the rest of the SMs and the model is strictly consistent then the ranking of the correction structures from most to least efficient is: 
\begin{align}
	& B \geq A>I>D, & \text{if the systems are in feedback and $\forall i$ $\exists$ SP$i$ $\sim$ SM$i$.} 
\end{align}

As B almost systematically provides the best local corrections, let's check if this is again related to a better performance in terms of convergence speed on the optimum of the plant. One shows on Figure.~\ref{fig:5_29_Exemple_5_2_Results_3} the results of simulations initialized at 5 uniformly distributed point in the input space, i.e. $u_0 =\{-1,-0.5,0,0.5,1\}$, and using all the models of each scenario.  One can see that the methods that bring the best local corrections are those that converge the fastest on $\bm{u}_p^{\star}$. 

\noindent
\begin{minipage}[h]{\linewidth}
	\begin{center}
		\textbf{Scenario 1: A$\rightleftarrows$NL}
	\end{center}
	
	\vspace{-3mm}
	\begin{minipage}[h]{\linewidth}
		\vspace*{0pt}
		{\centering
			\begin{minipage}[t]{5cm}{\centering%
					\includegraphics[width=5cm]{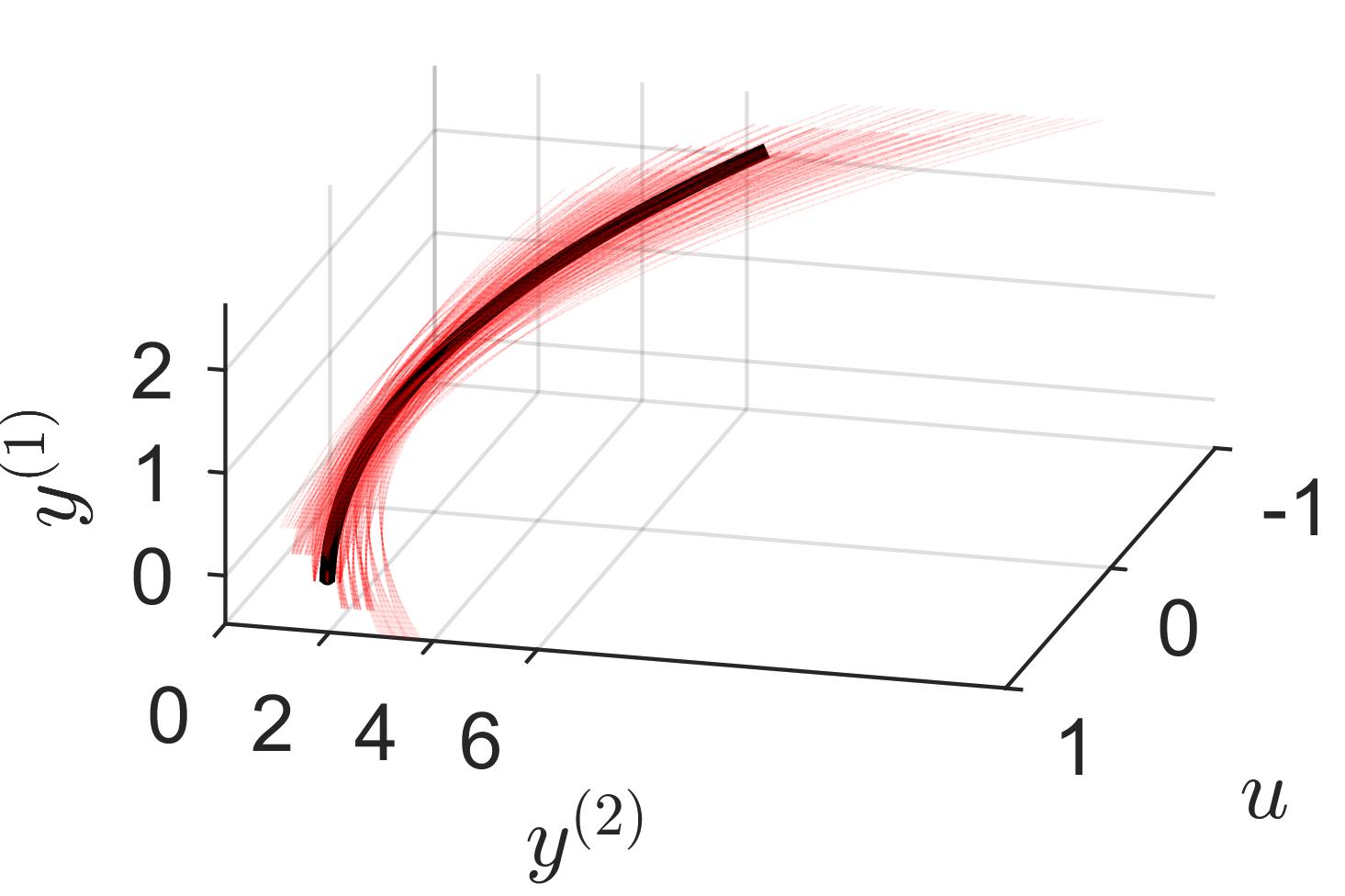}\\
					a) Functions $f^{(1)}$ and $f_p^{(1)}$}
			\end{minipage}\hskip -0ex
			\begin{minipage}[t]{5cm}\centering%
				\includegraphics[width=4.45cm]{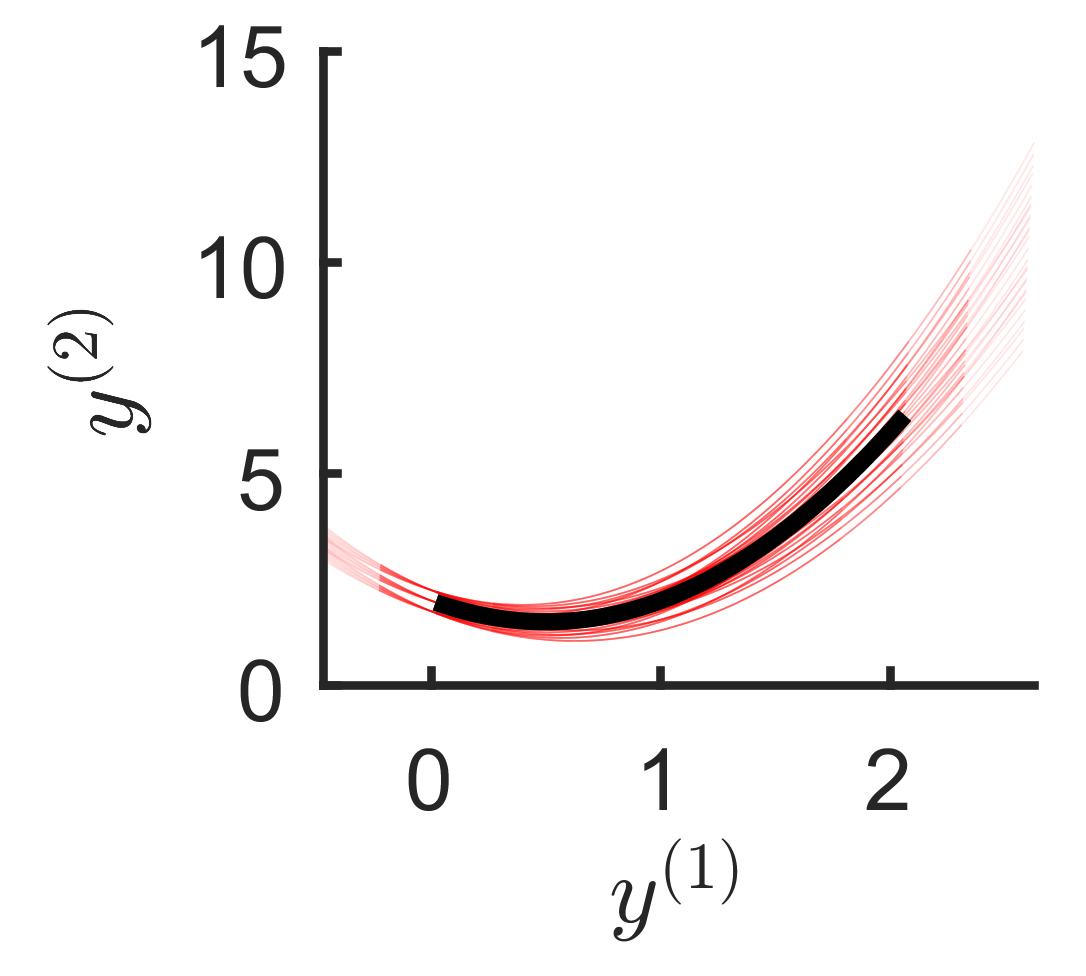}\\
				b) Functions $f^{(2)}$ and $f_p^{(2)}$
			\end{minipage}\hskip -0ex
			\begin{minipage}[t]{4.45cm}\centering%
				\includegraphics[width=4.45cm]{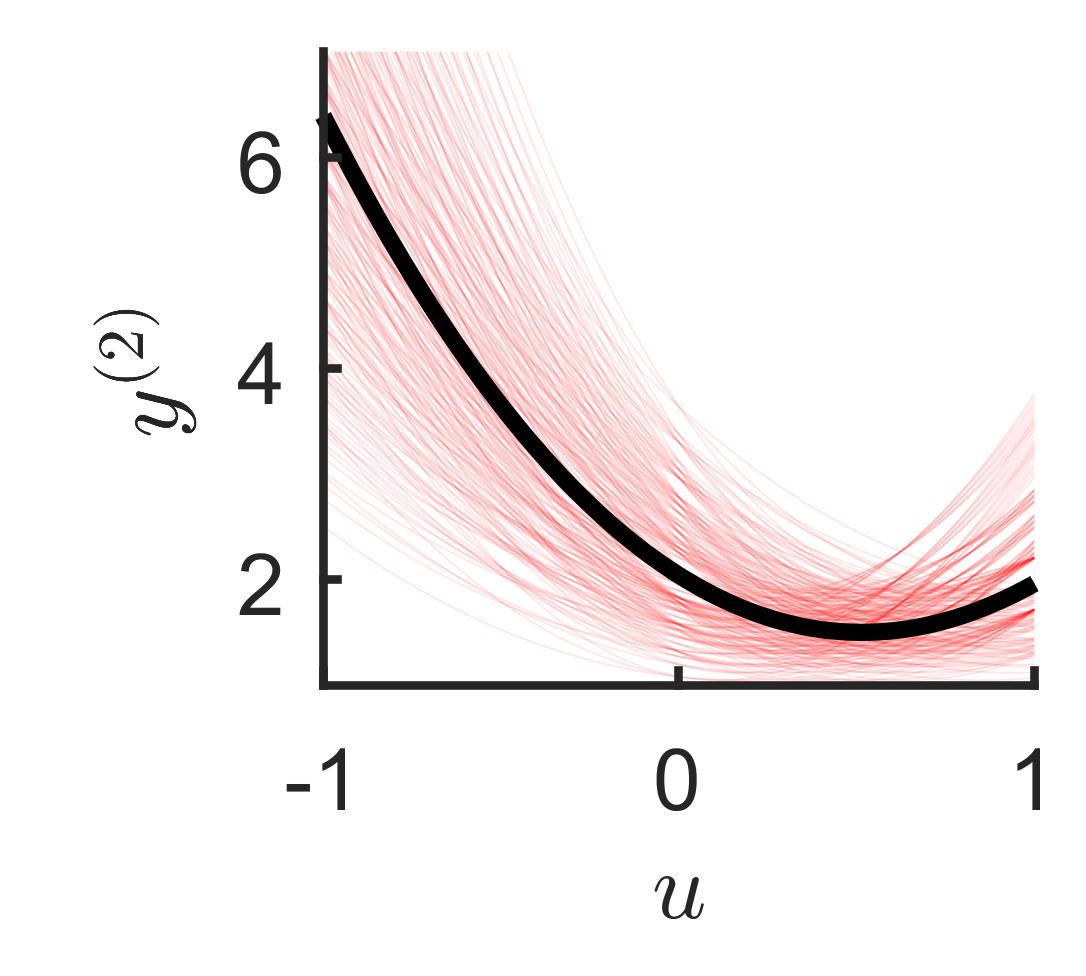}\\
				c) Functions $f$ and $f_p$
			\end{minipage} \\
			
			\medskip
			
			\begin{minipage}[h]{\linewidth} \centering
				\includegraphics[trim={0.1cm 0.2cm 0.2cm  0.1cm },clip,width=3.5cm]{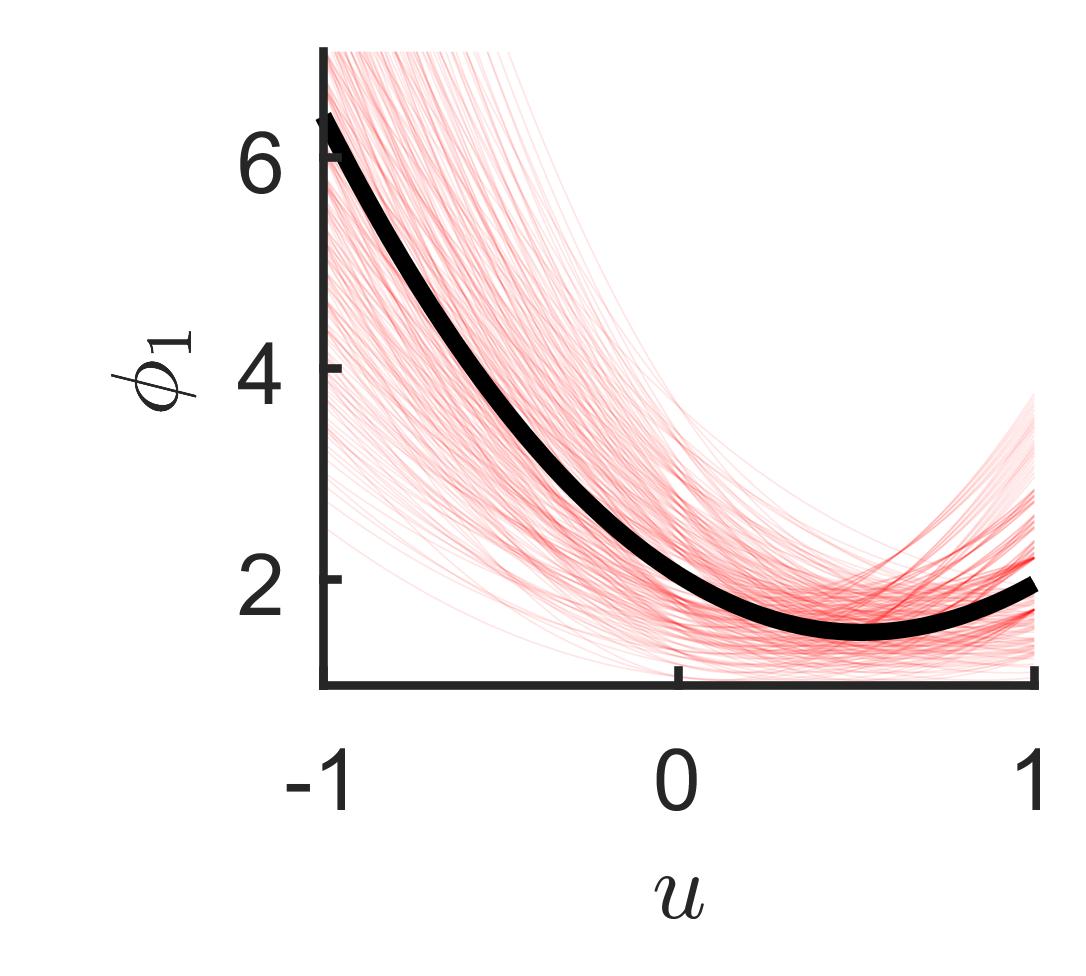}\hskip -0ex
				\includegraphics[trim={0.1cm 0.2cm 0.2cm  0.1cm },clip,width=3.5cm]{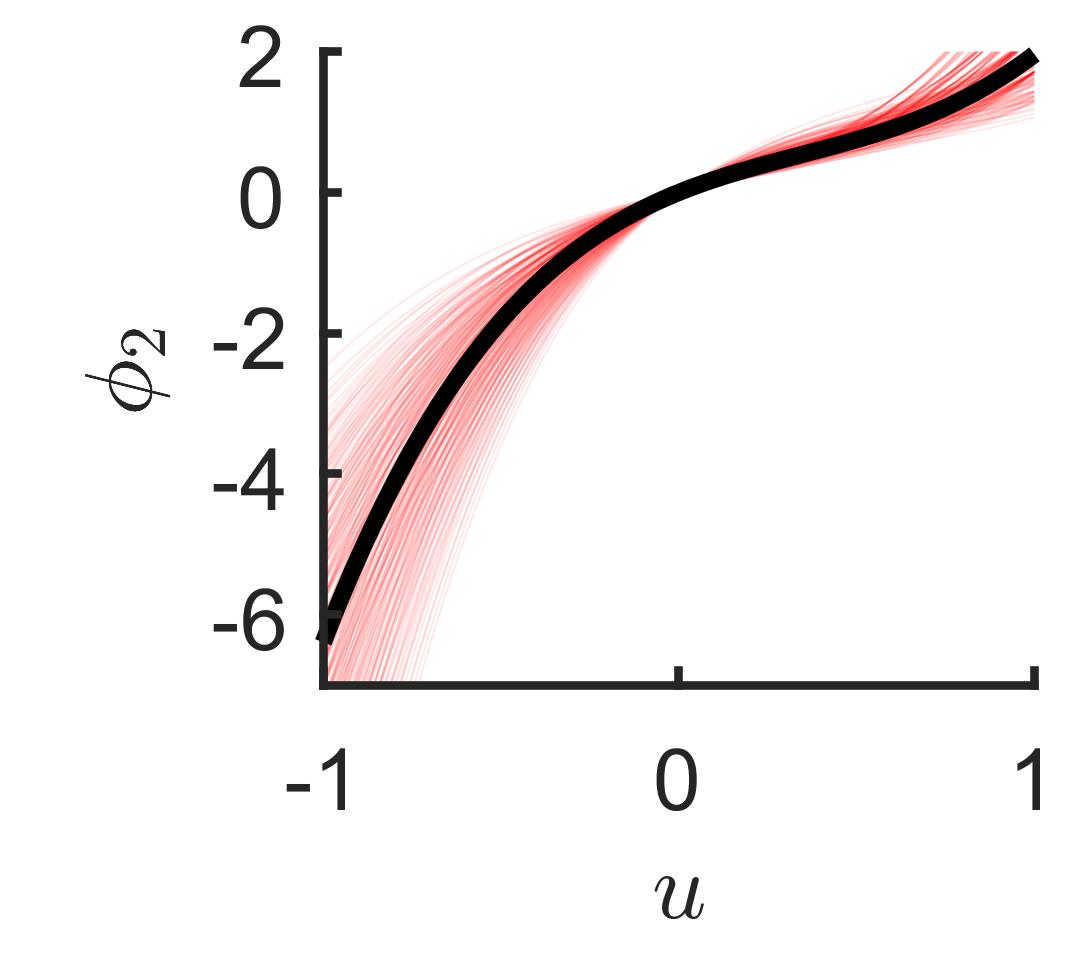}\hskip -0ex
				\includegraphics[trim={0.1cm 0.2cm 0.2cm  0.1cm },clip,width=3.5cm]{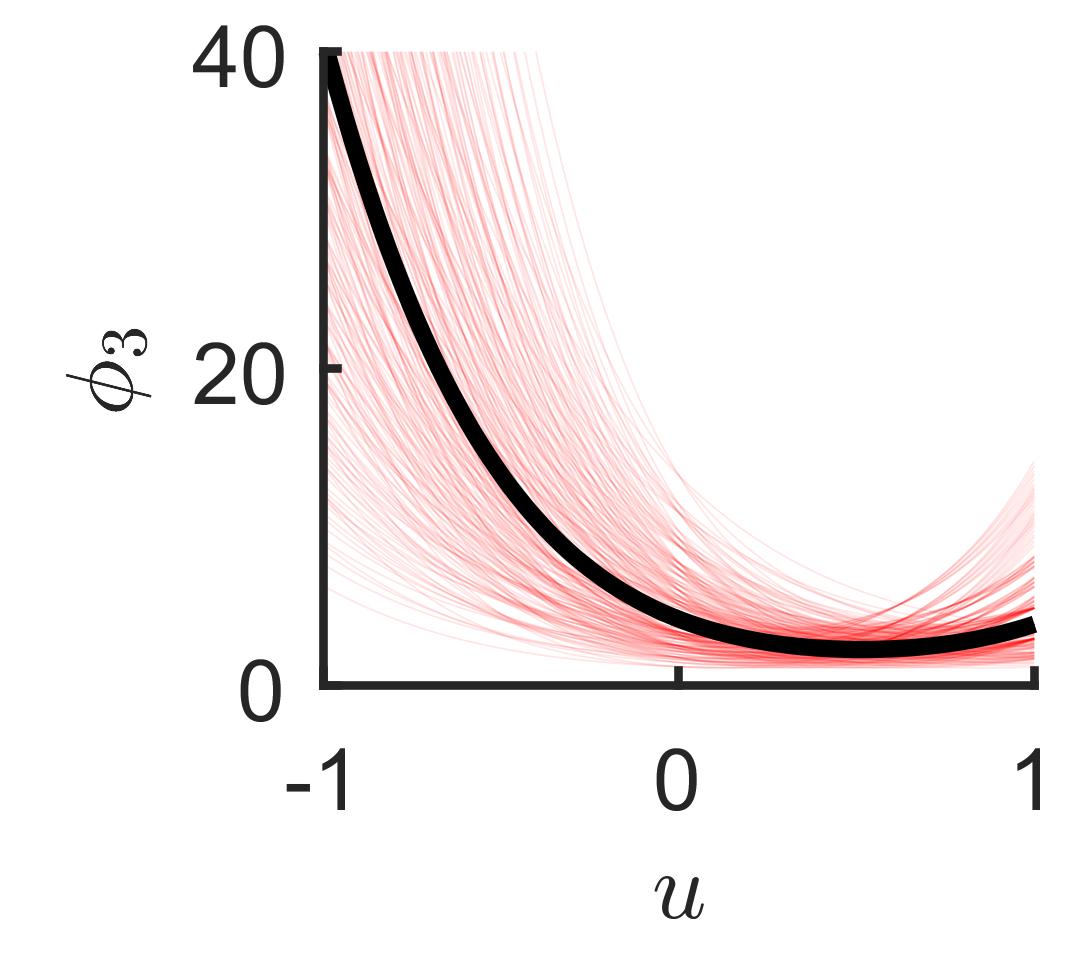}\hskip -0ex
				\includegraphics[trim={0.1cm 0.2cm 0.2cm  0.1cm },clip,width=3.5cm]{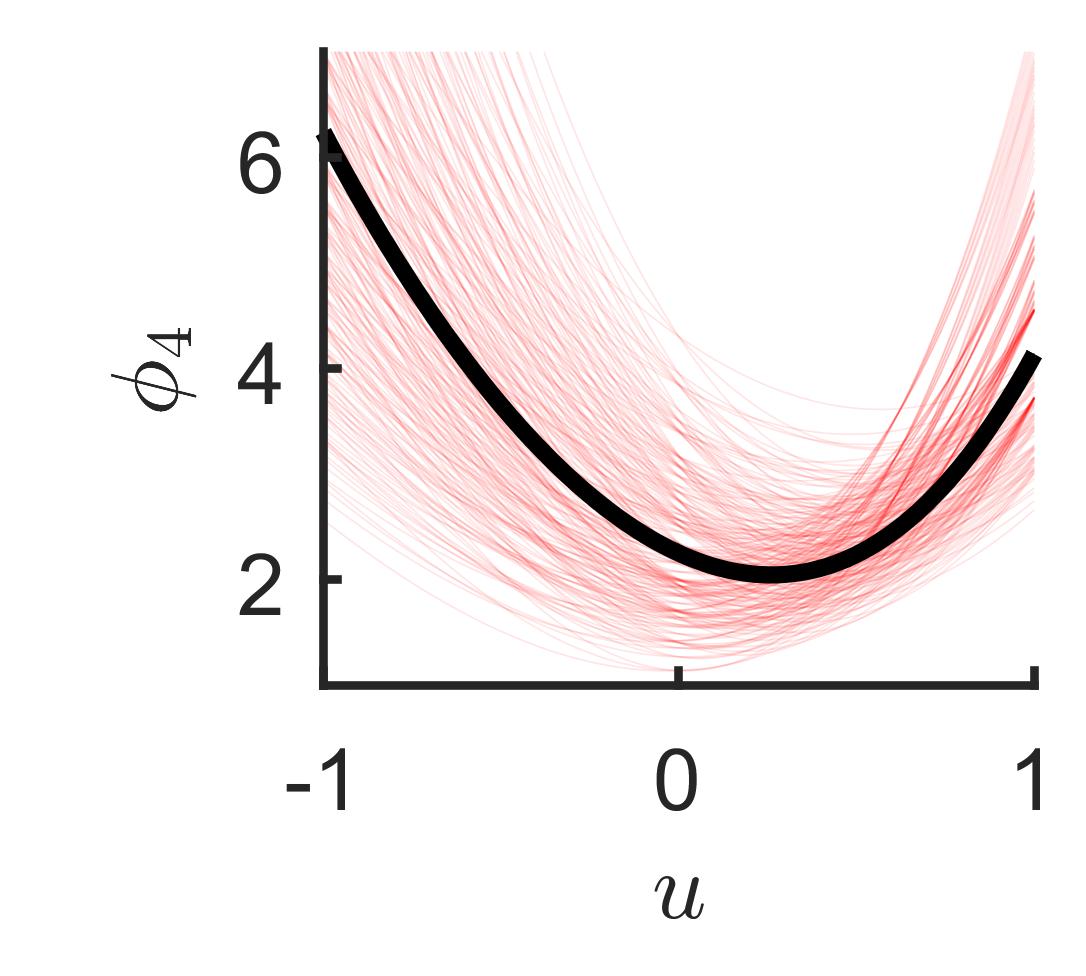}\\
				d) Functions $\phi_1$, $\phi_2$, $\phi_3$, and  $\phi_4$,
			\end{minipage} 
			\medskip
		}	
		\textcolor{red}{\raisebox{0.5mm}{\rule{0.5cm}{0.05cm}}}   : Model, 
		\textcolor{black}{\raisebox{0.5mm}{\rule{0.5cm}{0.1cm}}} : Plant.
		\vspace{-2mm}
		\captionof{figure}{Sc.1: Graphical description of the RTO problems}
		\label{fig:5_20_Exemple_5_2_sc1_Plant_and_Model}
	\end{minipage} \\
	
	\medskip
	
	\begin{minipage}[h]{\linewidth}
		\vspace*{0pt}
		{\centering
			\includegraphics[trim={2.75cm 0.2cm 0.2cm  0.2cm },clip,width=3.45cm]{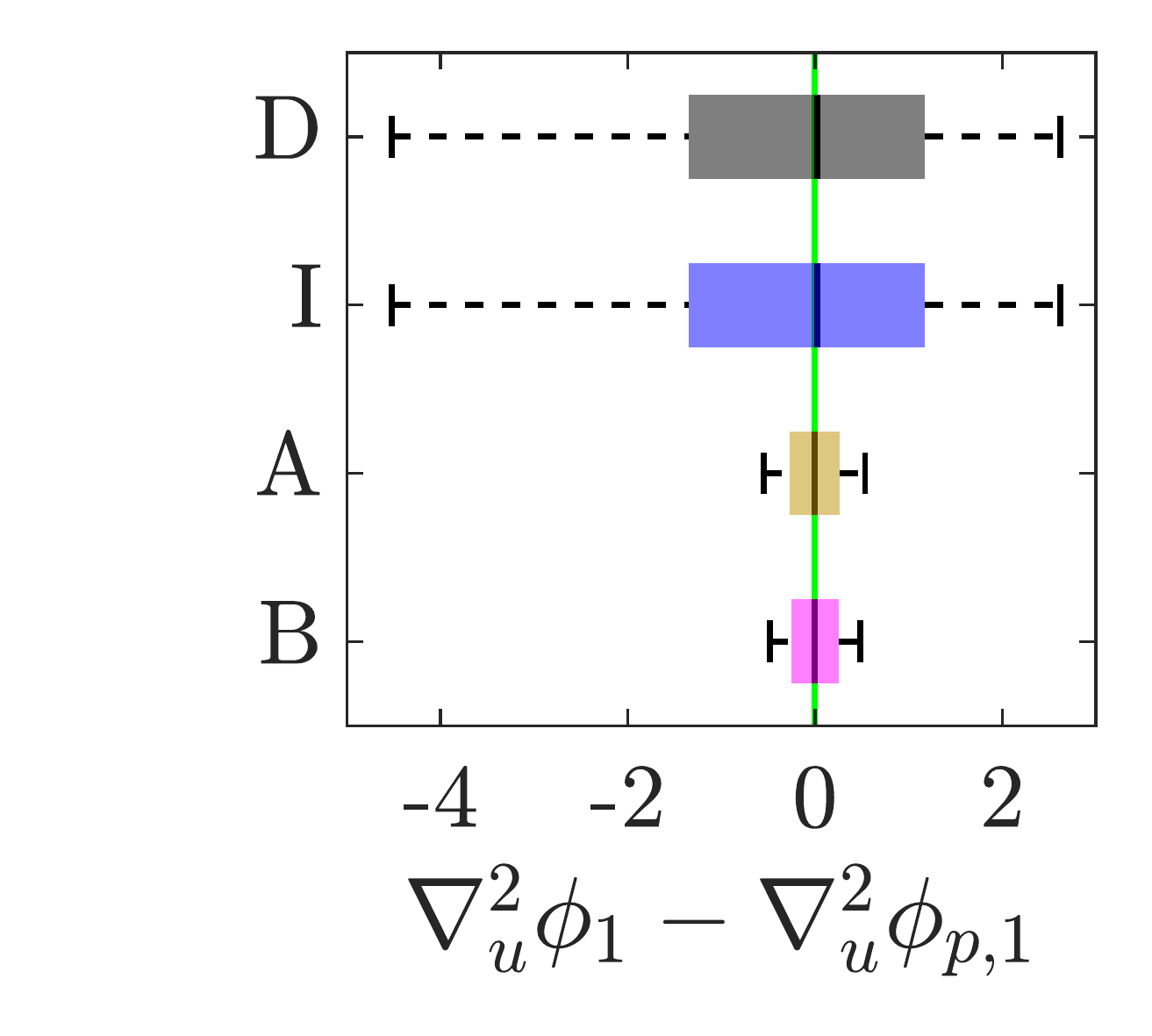}\hskip -0ex
			\includegraphics[trim={2.75cm 0.2cm 0.2cm  0.2cm },clip,width=3.45cm]{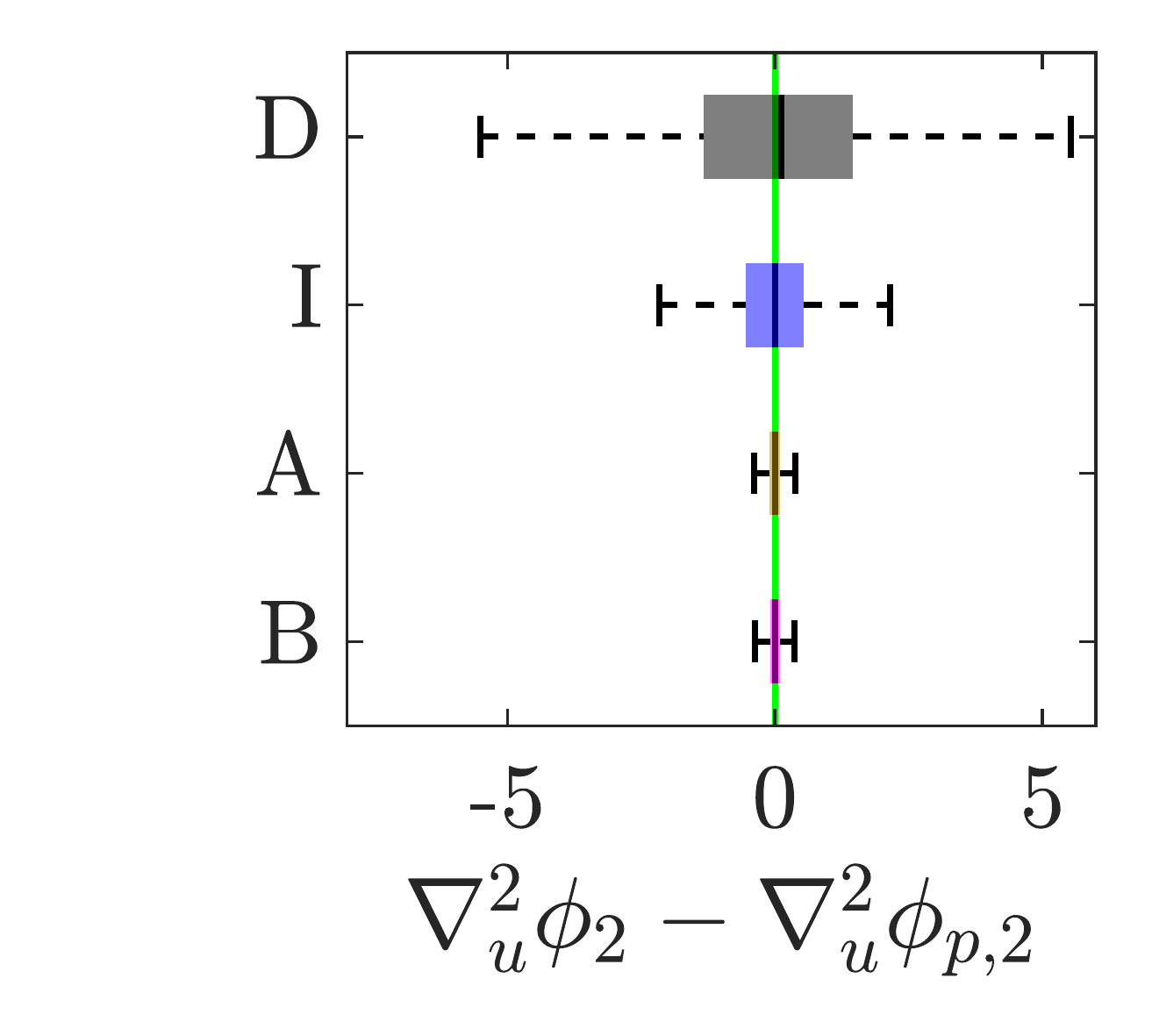}\hskip -0ex
			\includegraphics[trim={2.75cm 0.2cm 0.2cm  0.2cm },clip,width=3.45cm]{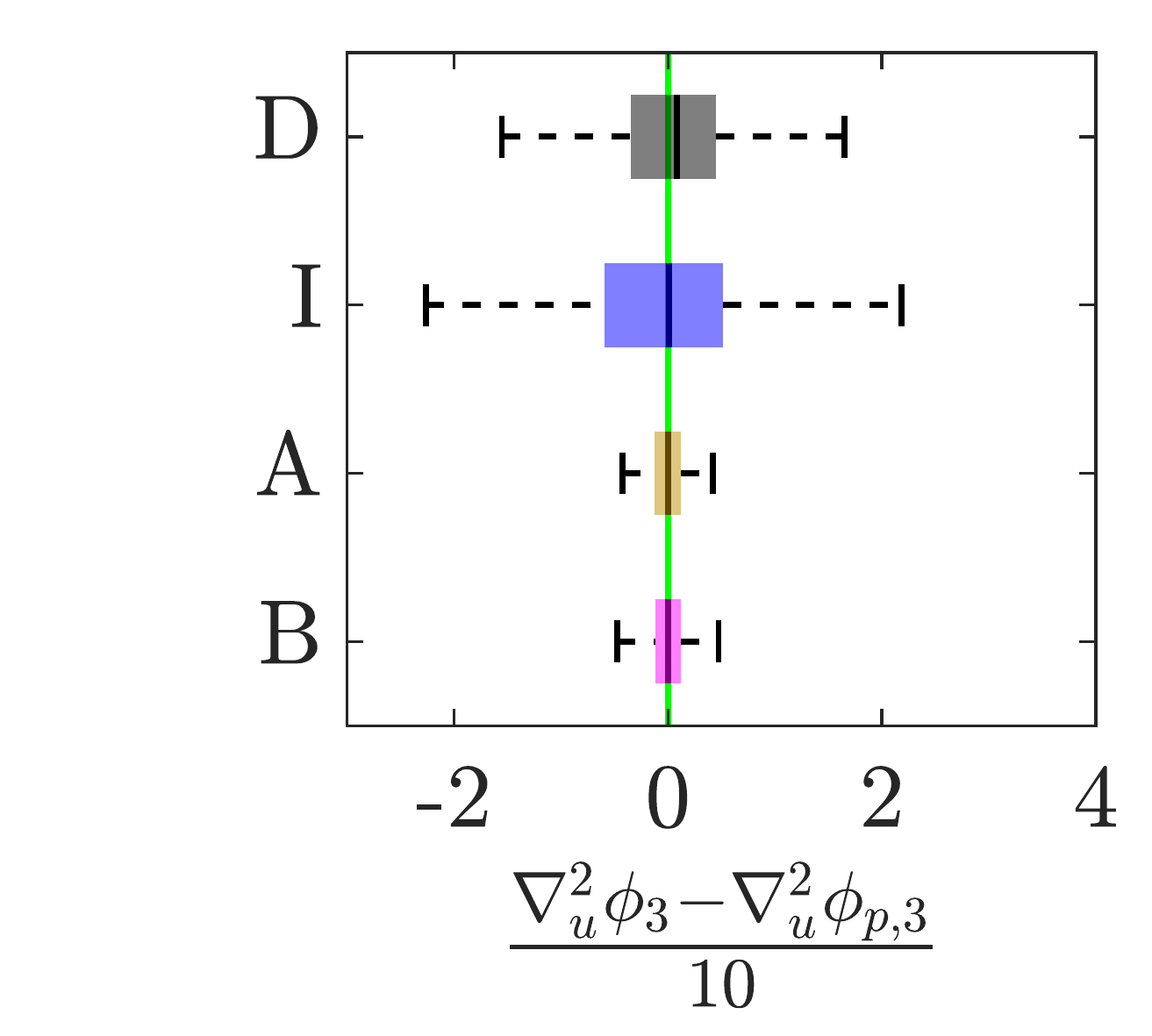}\hskip -0ex
			\includegraphics[trim={2.75cm 0.2cm 0.2cm  0.2cm },clip,width=3.45cm]{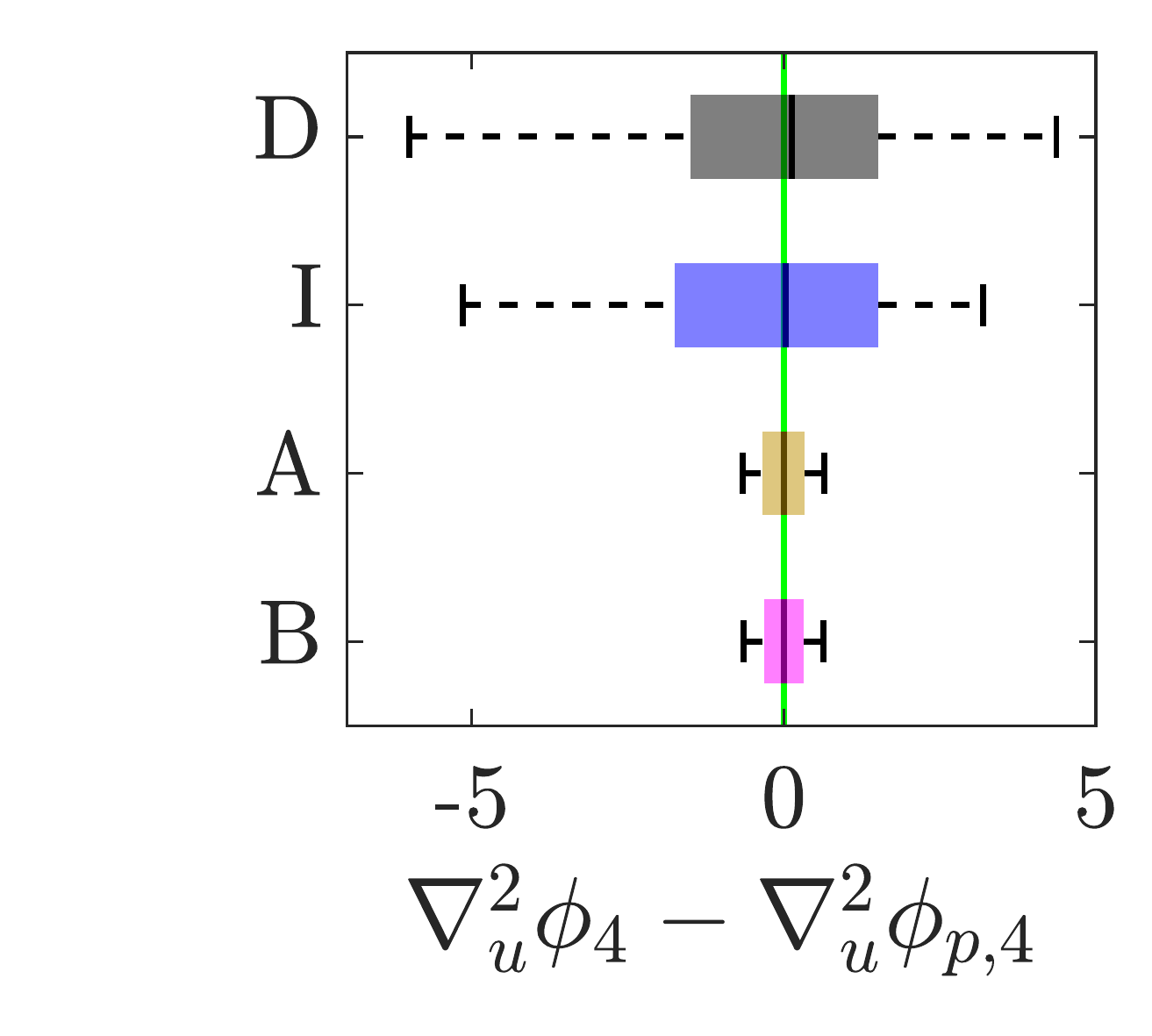}
		}
		\vspace{-2mm}
		\captionof{figure}{Sc.1: Statistical distributions of the prediction errors on the  Hessian of the plant's cost functions at the correction point for the structures D, I, A, and B.}
		\label{fig:5_21_Exemple_5_2_sc1_Results_1}
	\end{minipage} \\
	\begin{minipage}[h]{\linewidth}
		\vspace*{0pt}
		{\centering
			\includegraphics[trim={0.7cm 0.2cm 0.2cm  0.2cm },clip,width=3.45cm]{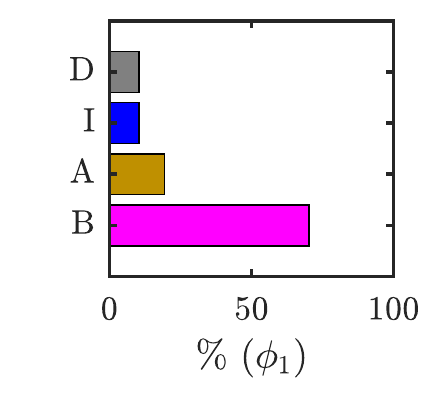}\hskip -0ex
			\includegraphics[trim={0.7cm 0.2cm 0.2cm  0.2cm },clip,width=3.45cm]{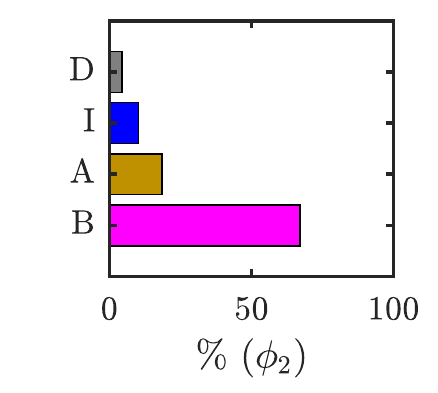}\hskip -0ex
			\includegraphics[trim={0.7cm 0.2cm 0.2cm  0.2cm },clip,width=3.45cm]{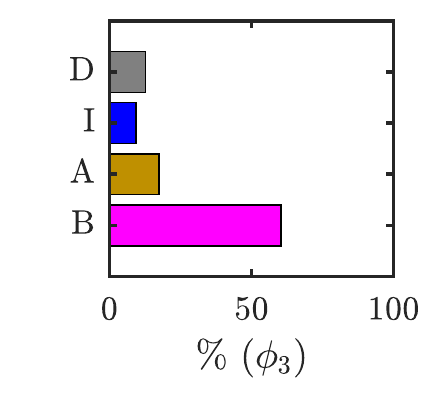}\hskip -0ex
			\includegraphics[trim={0.7cm 0.2cm 0.2cm  0.2cm },clip,width=3.45cm]{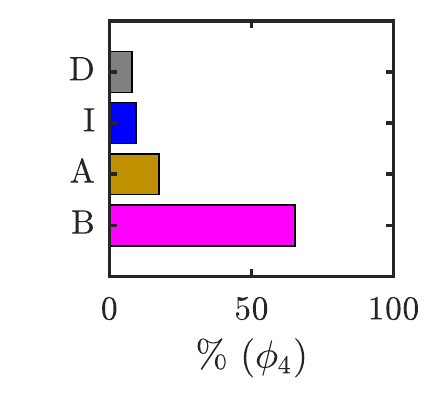} \\
		}
		\vspace{-2mm}
		\captionof{figure}{Sc.1: For each correction structure one gives here the percentage of cases for which no other structure provides better results  \textit{(if two structures provide the same best result then both take the point)}}
		\label{fig:5_22_Exemple_5_2_sc1_Results_2}
	\end{minipage} \\
\end{minipage}

\noindent
\begin{minipage}[h]{\linewidth}
	\begin{center}
		\textbf{Scenario 2: NL$\rightleftarrows$A}
	\end{center}
	
	\vspace{-3mm}
	\begin{minipage}[h]{\linewidth}
		\vspace*{0pt}
		{\centering
			\begin{minipage}[t]{5cm}{\centering%
					\includegraphics[width=5cm]{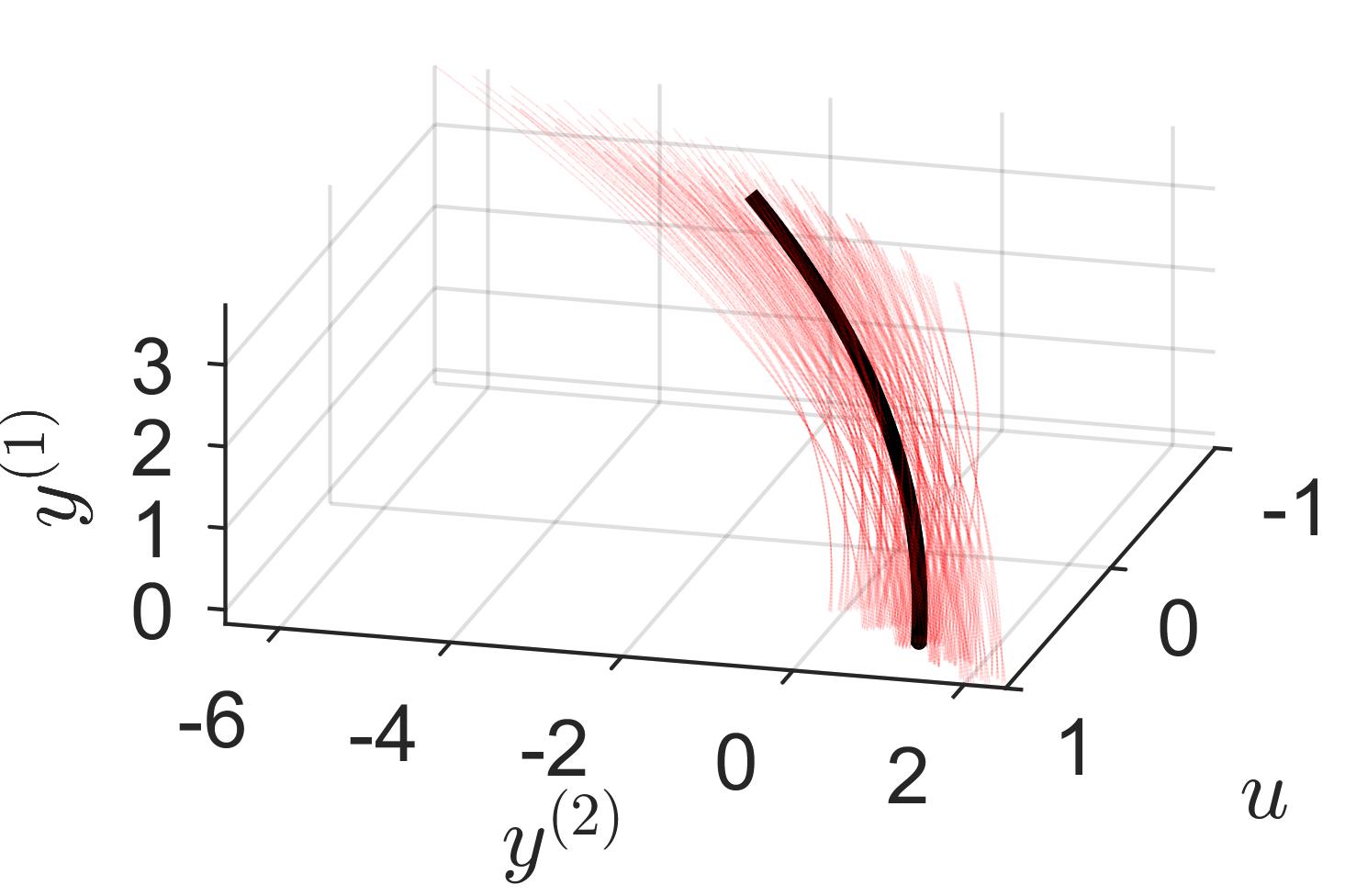}\\
					a) Fonction $f^{(1)}$ et $f_p^{(1)}$}
			\end{minipage}\hskip -0ex
			\begin{minipage}[t]{4.45cm}\centering%
				\includegraphics[width=4.45cm]{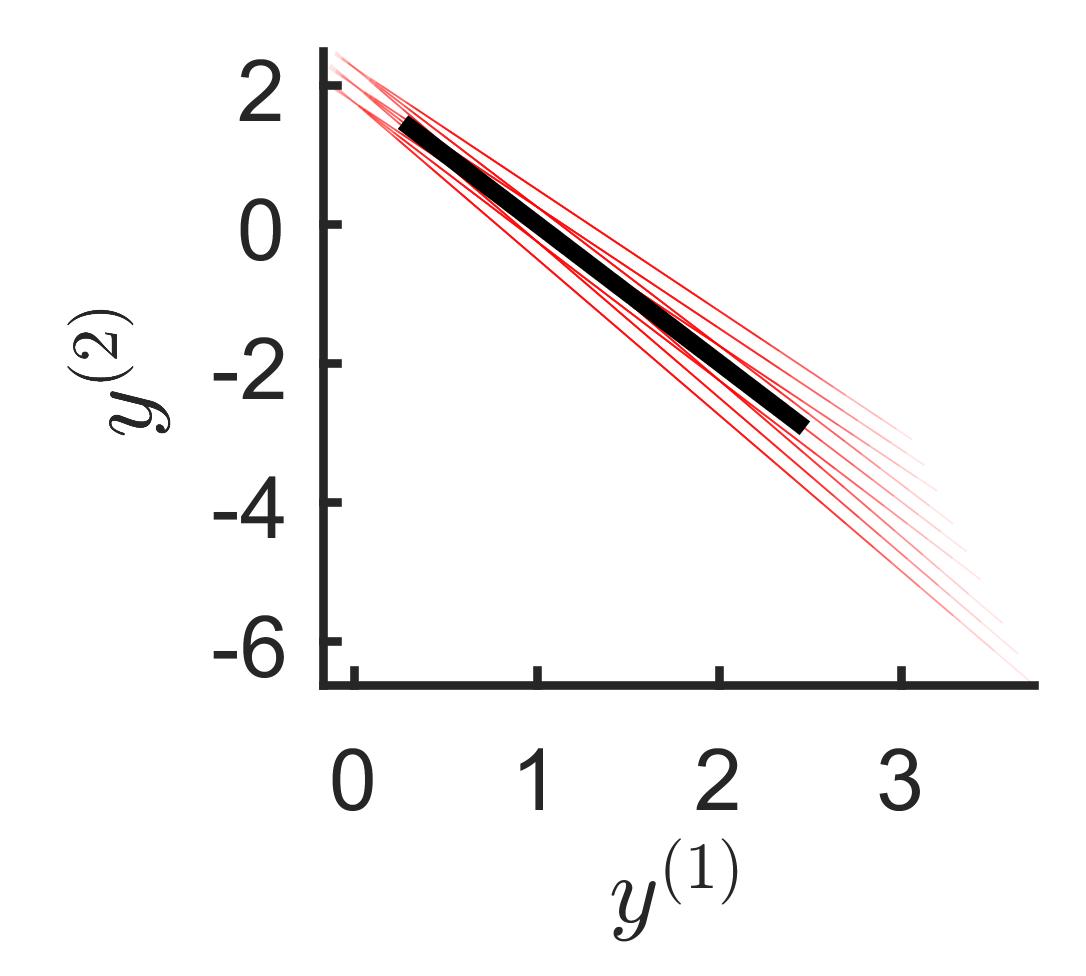}\\
				b) Fonction $f^{(2)}$ et $f_p^{(2)}$
			\end{minipage}\hskip -0ex
			\begin{minipage}[t]{4.45cm}\centering%
				\includegraphics[width=4.45cm]{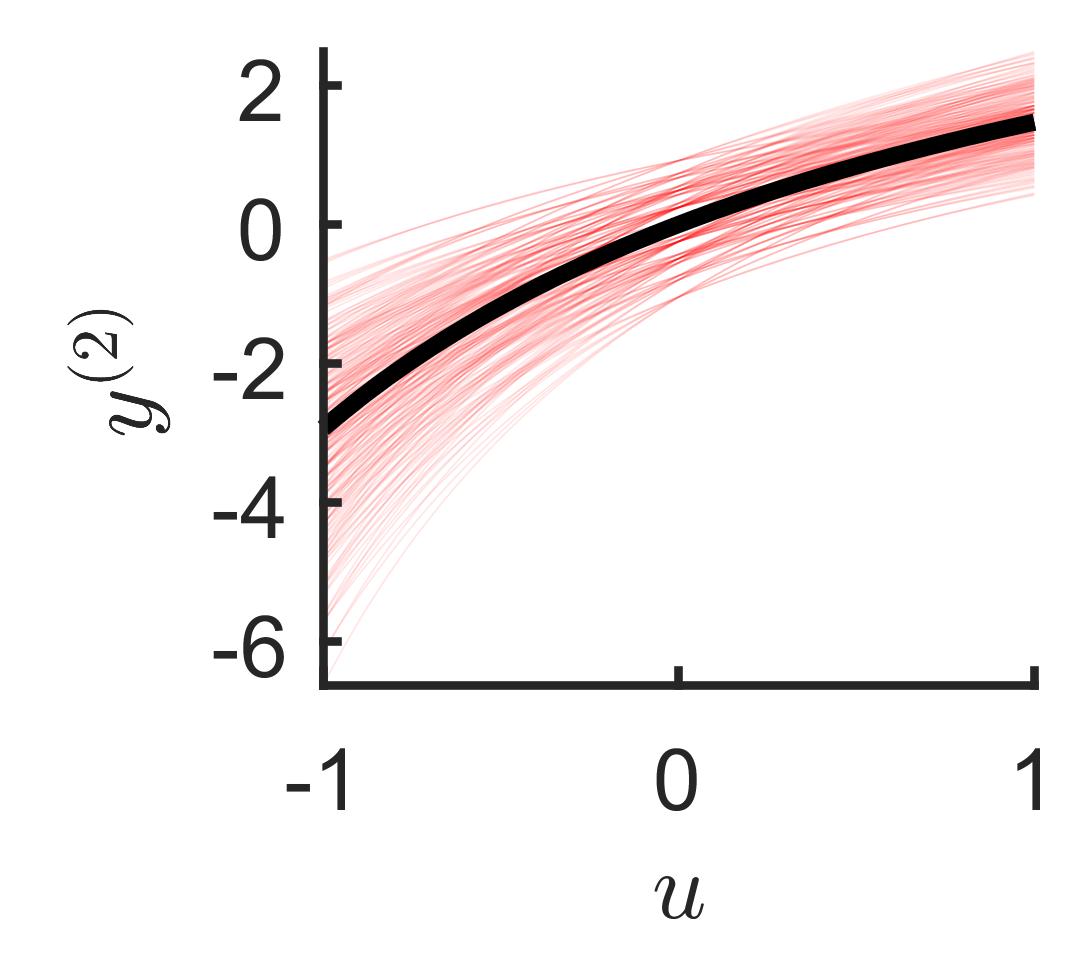}\\
				c) Fonction $f$ et $f_p$
			\end{minipage} \\
			
			\medskip
			
			\begin{minipage}[h]{\linewidth} \centering
				\includegraphics[trim={0.1cm 0.2cm 0.2cm  0.1cm },clip,width=3.5cm]{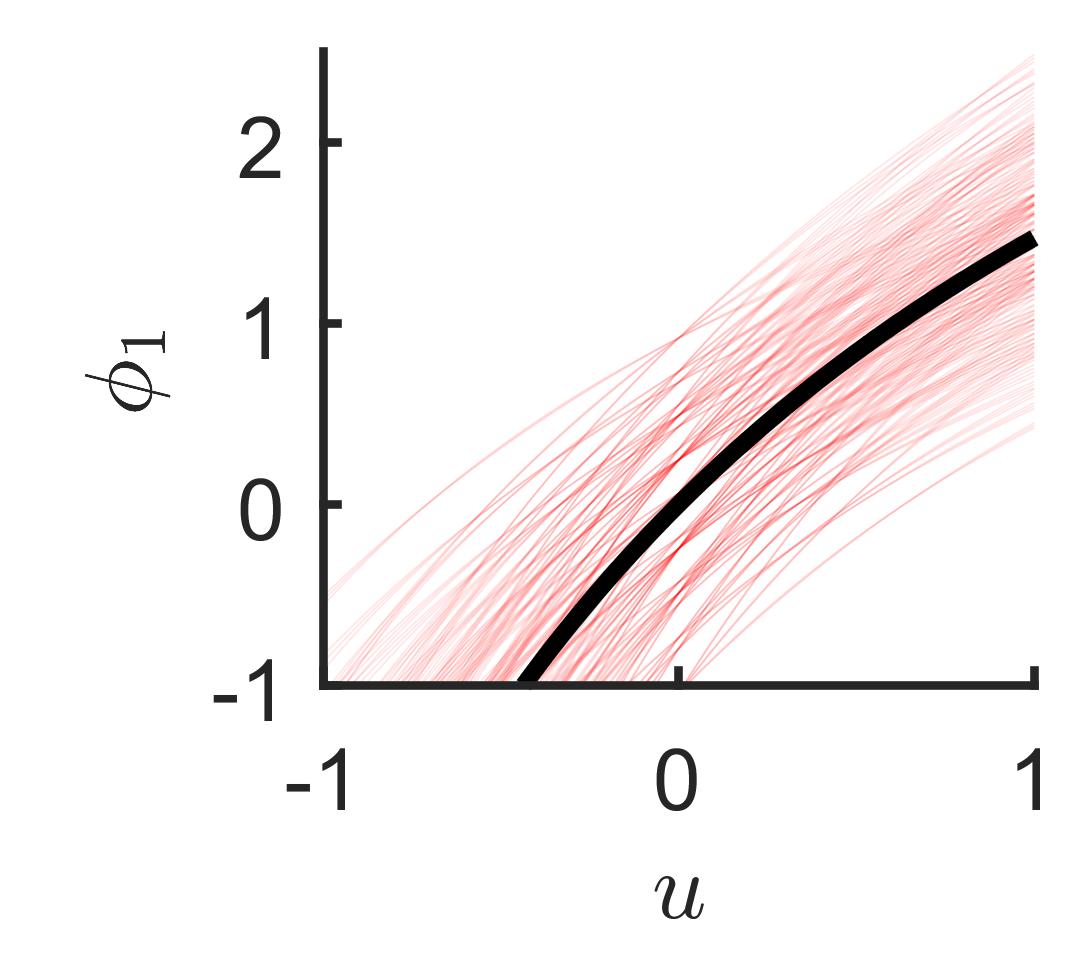}\hskip -0ex
				\includegraphics[trim={0.1cm 0.2cm 0.2cm  0.1cm },clip,width=3.5cm]{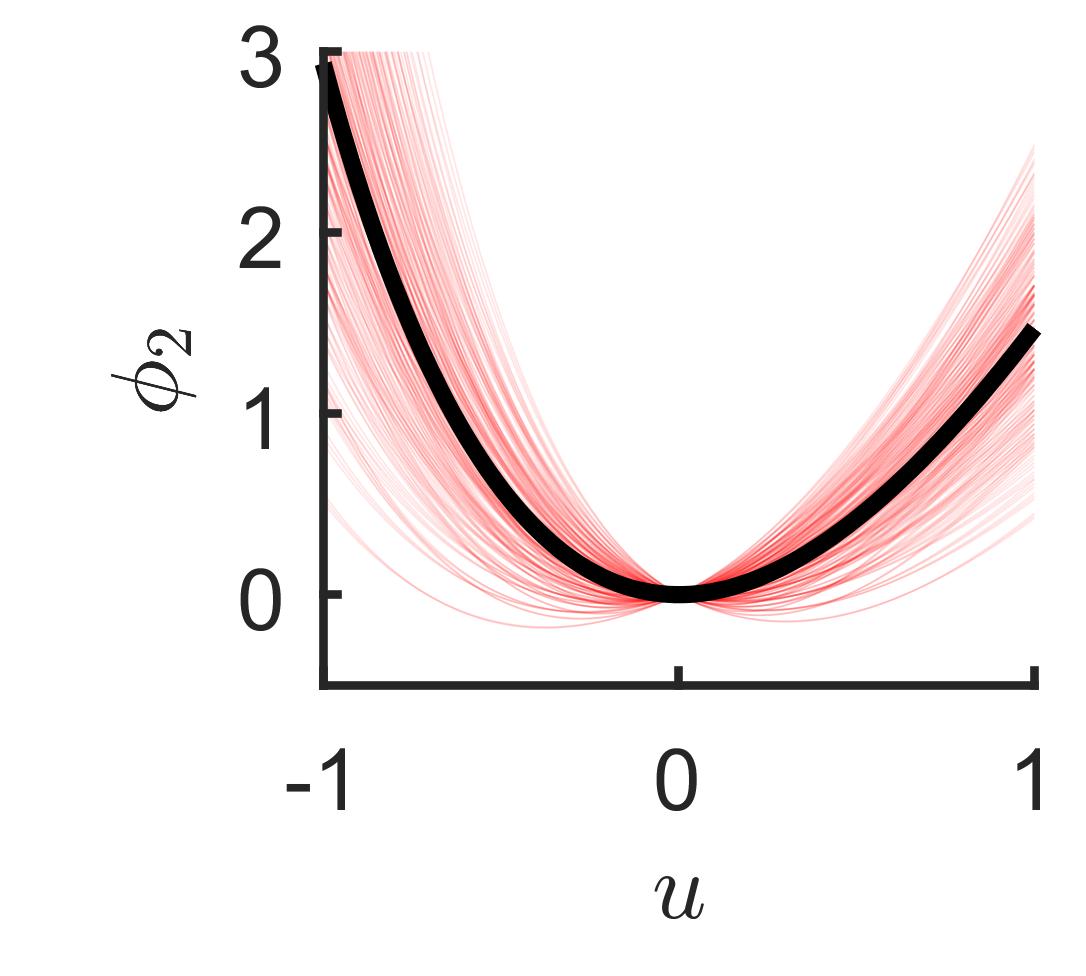}\hskip -0ex
				\includegraphics[trim={0.1cm 0.2cm 0.2cm  0.1cm },clip,width=3.5cm]{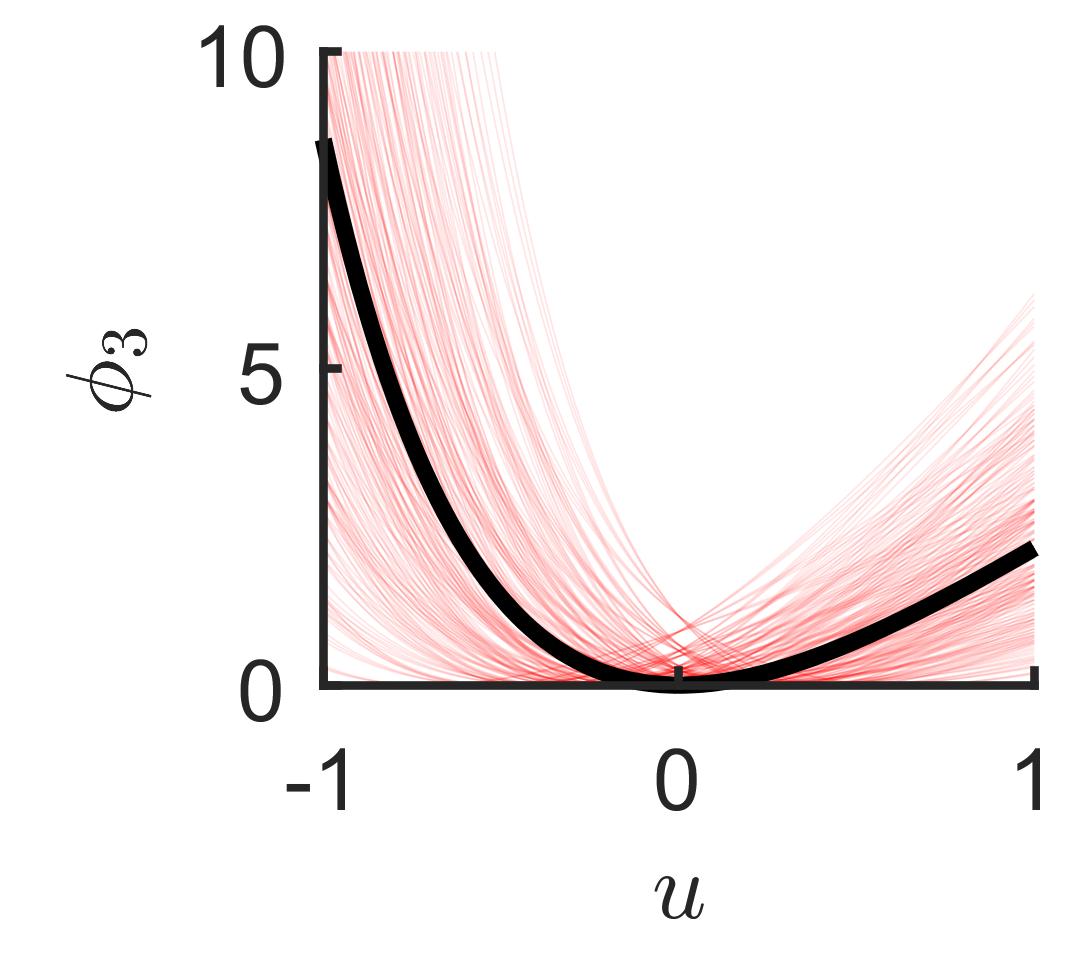}\hskip -0ex
				\includegraphics[trim={0.1cm 0.2cm 0.2cm  0.1cm },clip,width=3.5cm]{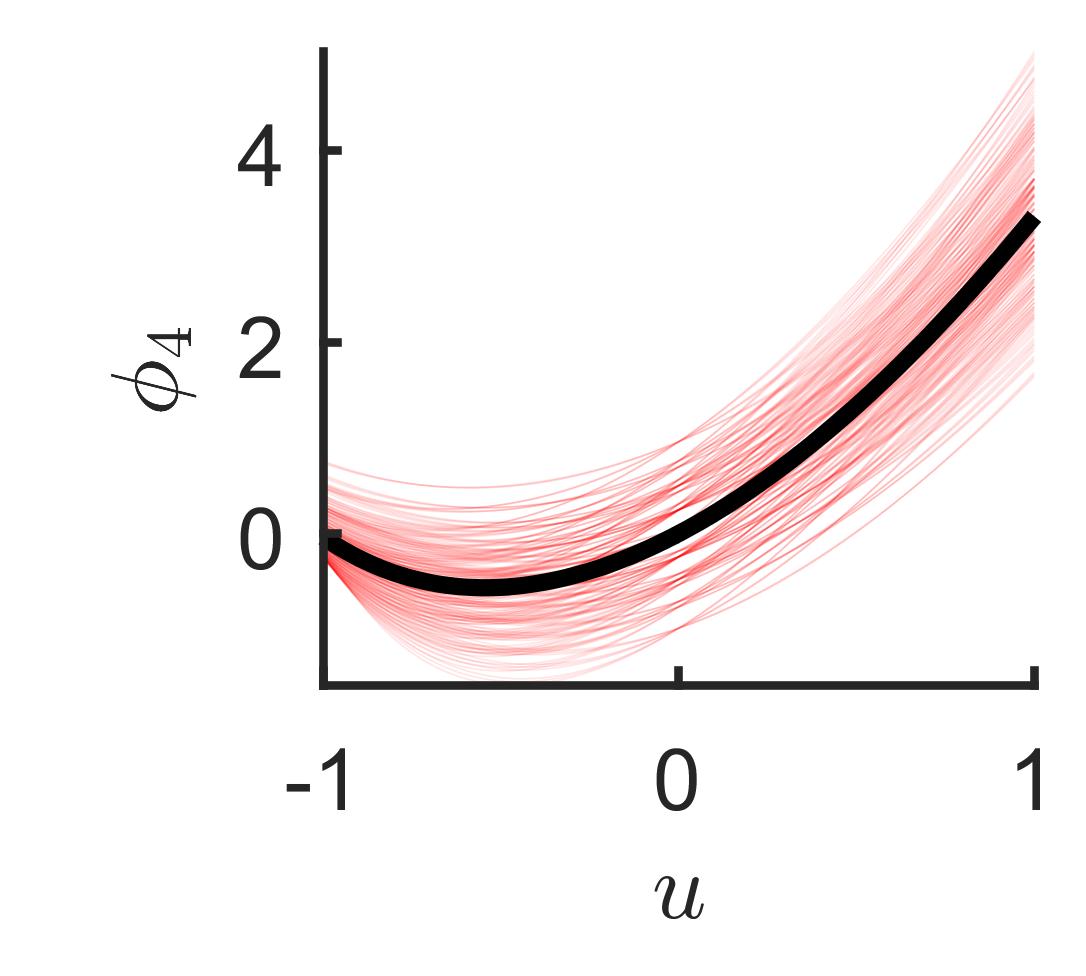}\\
				d) Fonctions $\phi_1$, $\phi_2$, $\phi_3$, et  $\phi_4$,
			\end{minipage} 
			\medskip
		}	
		\textcolor{red}{\raisebox{0.5mm}{\rule{0.5cm}{0.05cm}}}   : Modèles, 
		\textcolor{black}{\raisebox{0.5mm}{\rule{0.5cm}{0.1cm}}} : Usine.
		\vspace{-2mm}
		\captionof{figure}{Sc.2: Graphical description of the RTO problems}
		\label{fig:5_23_Exemple_5_2_sc2_Plant_and_Model}
	\end{minipage} \\
	
	\medskip
	
	\begin{minipage}[h]{\linewidth}
		\vspace*{0pt}
		{\centering
			\includegraphics[trim={2.75cm 0.2cm 0.2cm  0.2cm },clip,width=3.45cm]{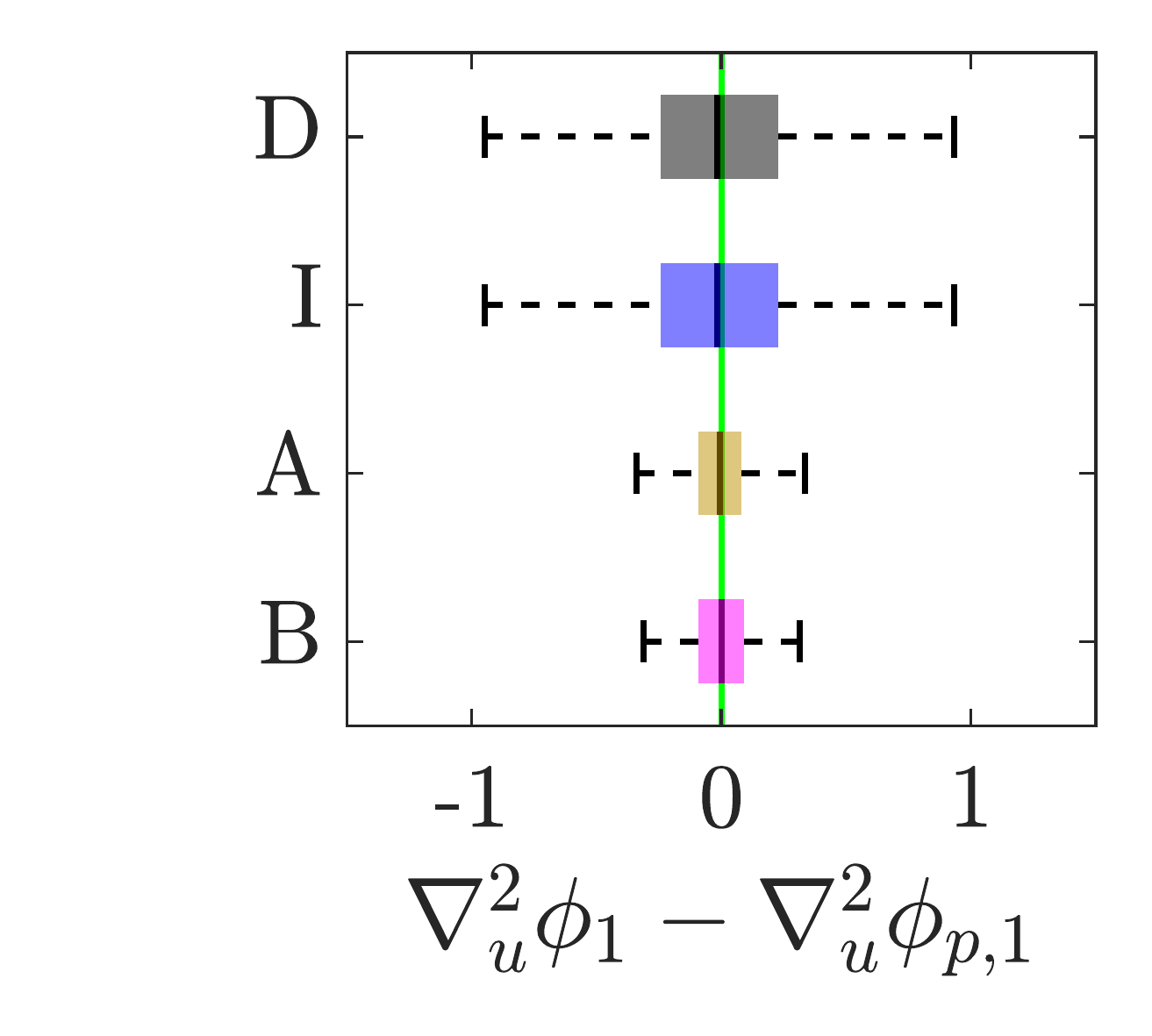}\hskip -0ex
			\includegraphics[trim={2.75cm 0.2cm 0.2cm  0.2cm },clip,width=3.45cm]{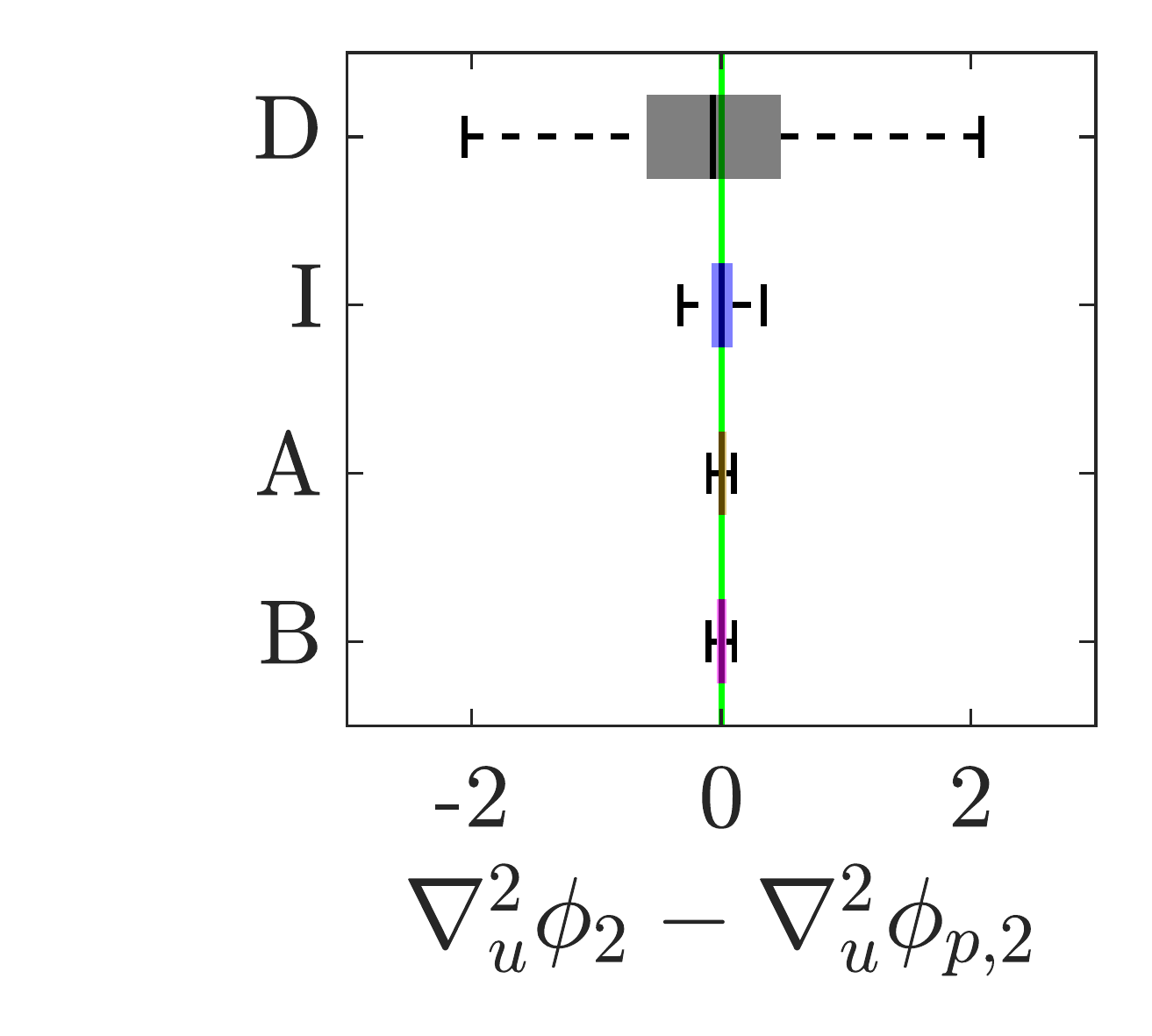}\hskip -0ex
			\includegraphics[trim={2.75cm 0.2cm 0.2cm  0.2cm },clip,width=3.45cm]{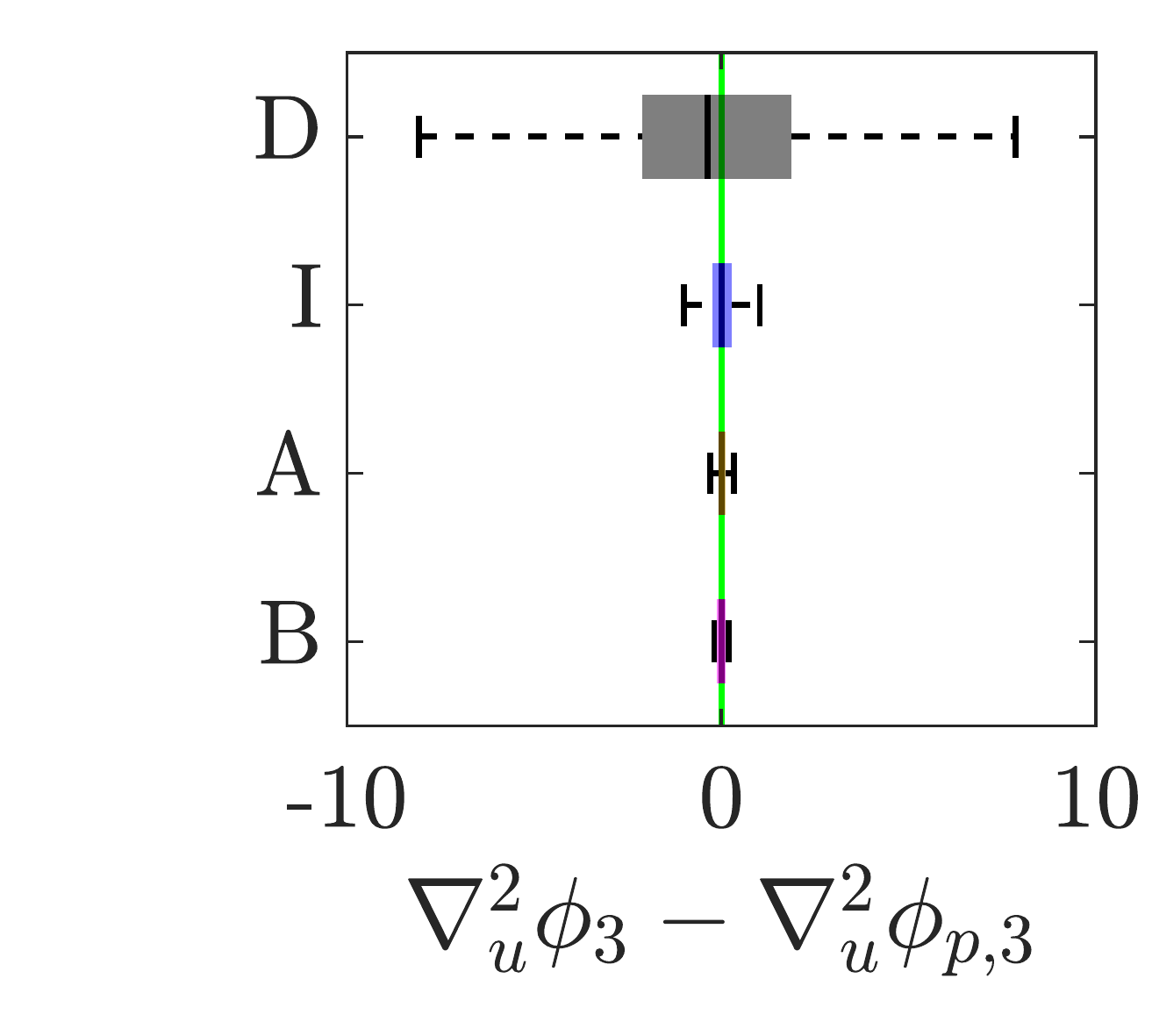}\hskip -0ex
			\includegraphics[trim={2.75cm 0.2cm 0.2cm  0.2cm },clip,width=3.45cm]{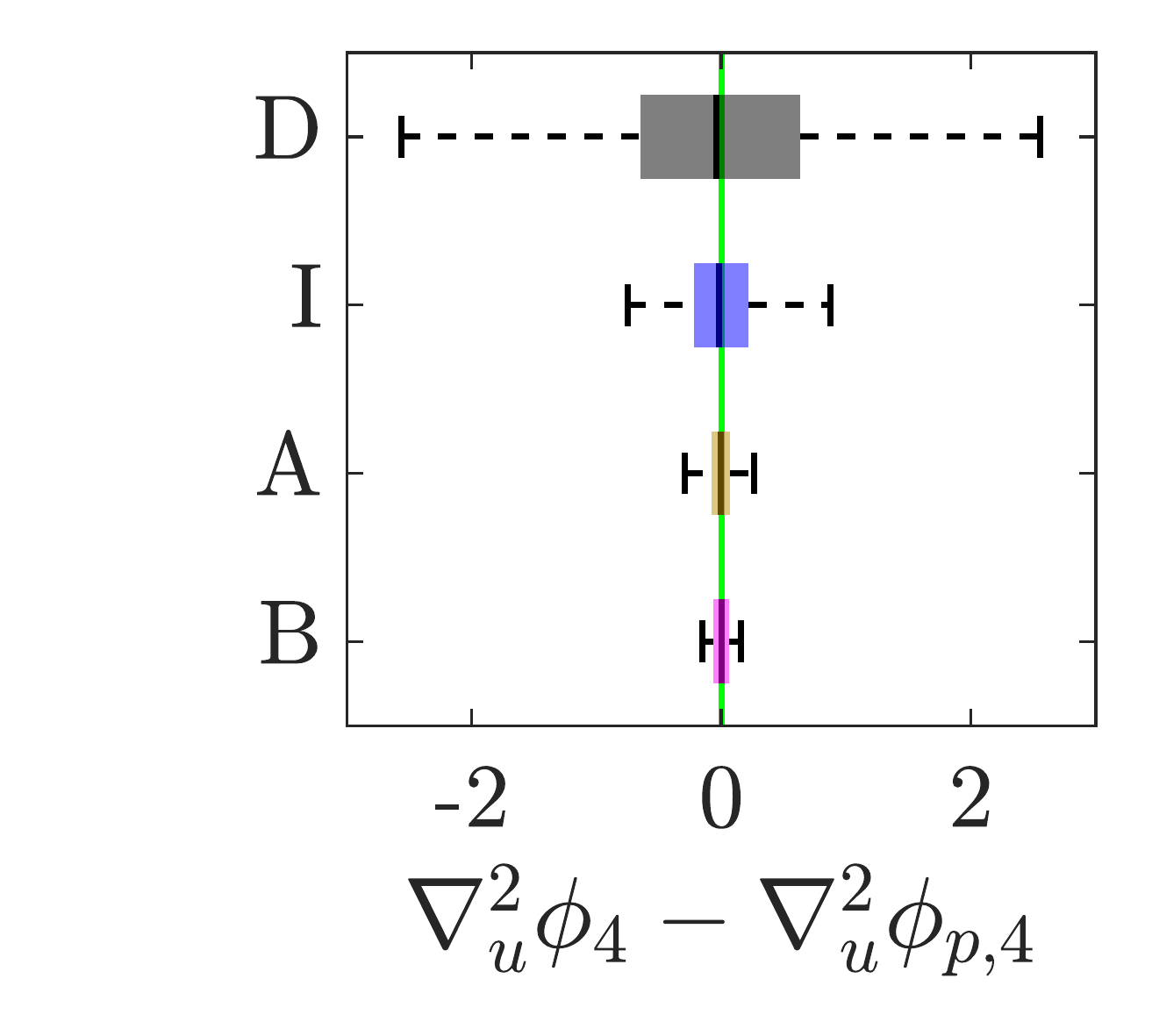}
		}
		\vspace{-2mm}
		\captionof{figure}{Sc.2: Statistical distributions of the prediction errors on the  Hessian of the plant's cost functions at the correction point for the structures D, I, A, and B.}
		\label{fig:5_24_Exemple_5_2_sc2_Results_1}
	\end{minipage} \\
	\begin{minipage}[h]{\linewidth}
		\vspace*{0pt}
		{\centering
			\includegraphics[trim={0.7cm 0.2cm 0.2cm  0.2cm },clip,width=3.45cm]{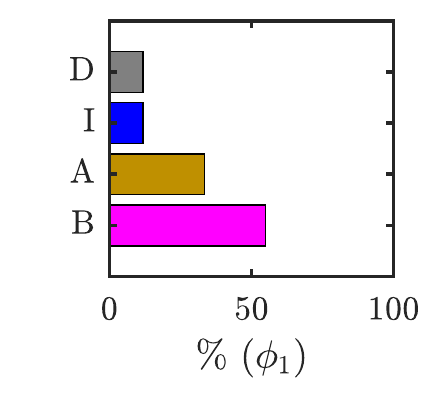}\hskip -0ex
			\includegraphics[trim={0.7cm 0.2cm 0.2cm  0.2cm },clip,width=3.45cm]{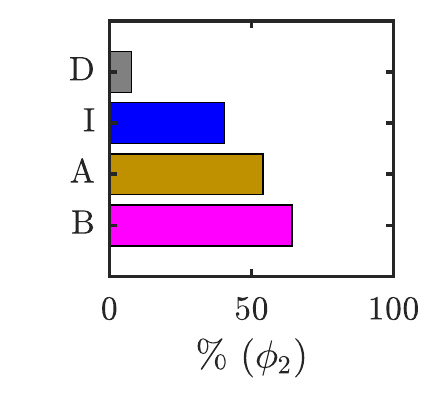}\hskip -0ex
			\includegraphics[trim={0.7cm 0.2cm 0.2cm  0.2cm },clip,width=3.45cm]{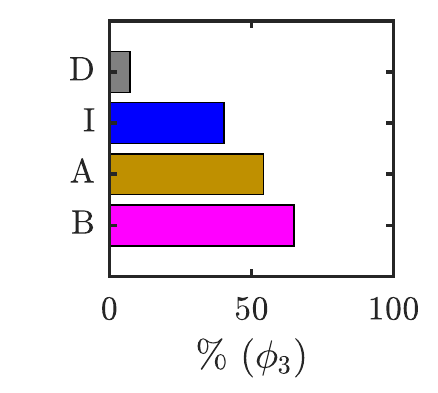}\hskip -0ex
			\includegraphics[trim={0.7cm 0.2cm 0.2cm  0.2cm },clip,width=3.45cm]{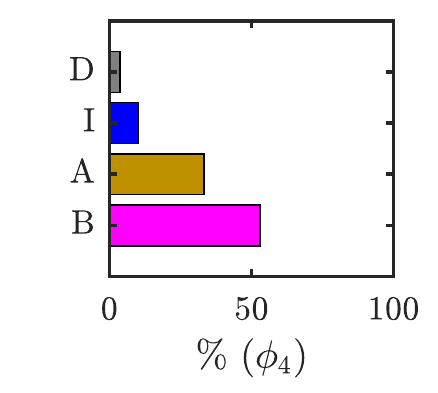} \\
		}
		\vspace{-2mm}
		\captionof{figure}{Sc.2: For each correction structure one gives here the percentage of cases for which no other structure provides better results  \textit{(if two structures provide the same best result then both take the point)}}
		\label{fig:5_25_Exemple_5_2_sc2_Results_2}
	\end{minipage} \\
\end{minipage}

\noindent
\begin{minipage}[h]{\linewidth}
	\begin{center}
		\textbf{Scenario 3: NL$\rightleftarrows$NL}
	\end{center}
	
	\vspace{-3mm}
	\begin{minipage}[h]{\linewidth}
		\vspace*{0pt}
		{\centering
			\begin{minipage}[t]{5cm}{\centering%
					\includegraphics[width=5cm]{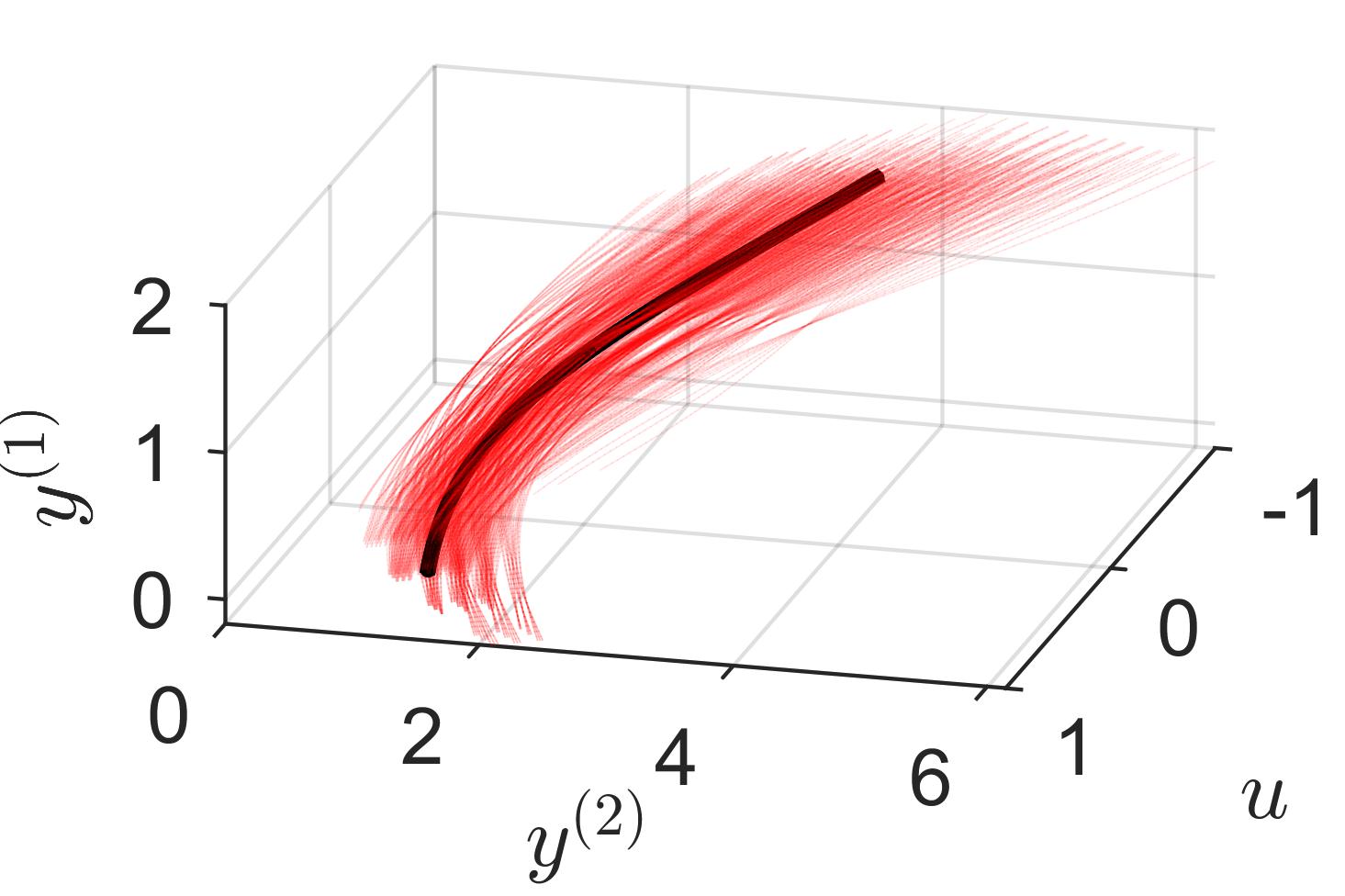}\\
					a) Function $f^{(1)}$ and $f_p^{(1)}$}
			\end{minipage}\hskip -0ex
			\begin{minipage}[t]{4.45cm}\centering%
				\includegraphics[width=4.45cm]{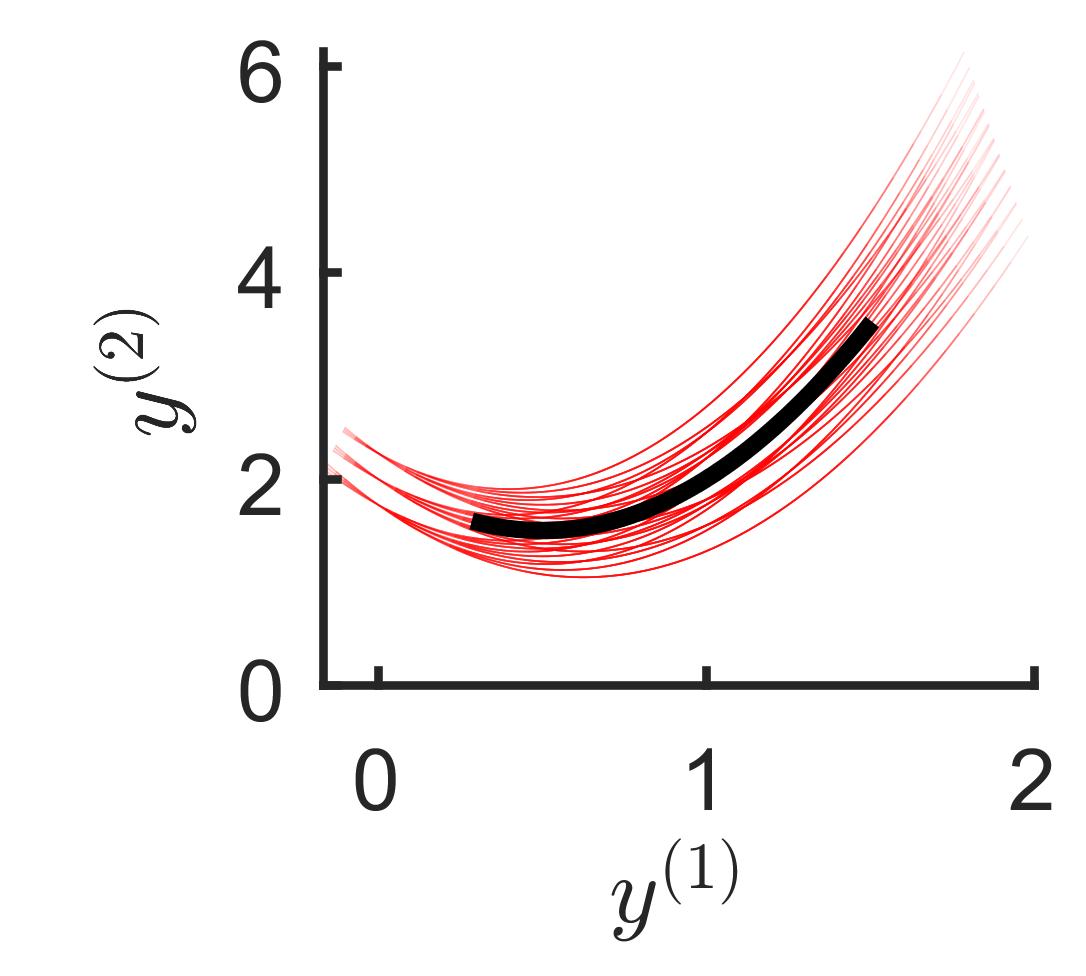}\\
				b) Function $f^{(2)}$ and $f_p^{(2)}$
			\end{minipage}\hskip -0ex
			\begin{minipage}[t]{4.45cm}\centering%
				\includegraphics[width=4.45cm]{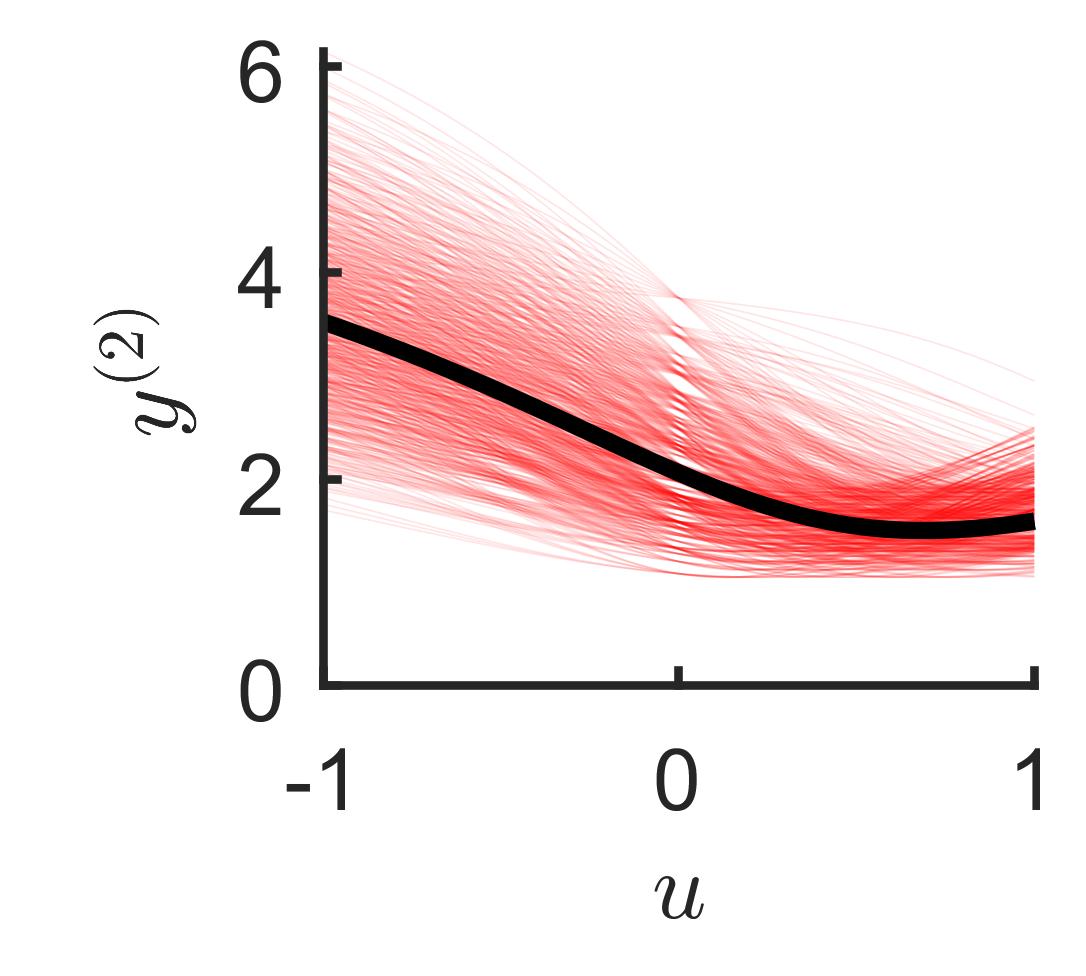}\\
				c) Function $f$ and $f_p$
			\end{minipage} \\
			
			\medskip
			
			\begin{minipage}[h]{\linewidth} \centering
				\includegraphics[trim={0.1cm 0.2cm 0.2cm  0.1cm },clip,width=3.5cm]{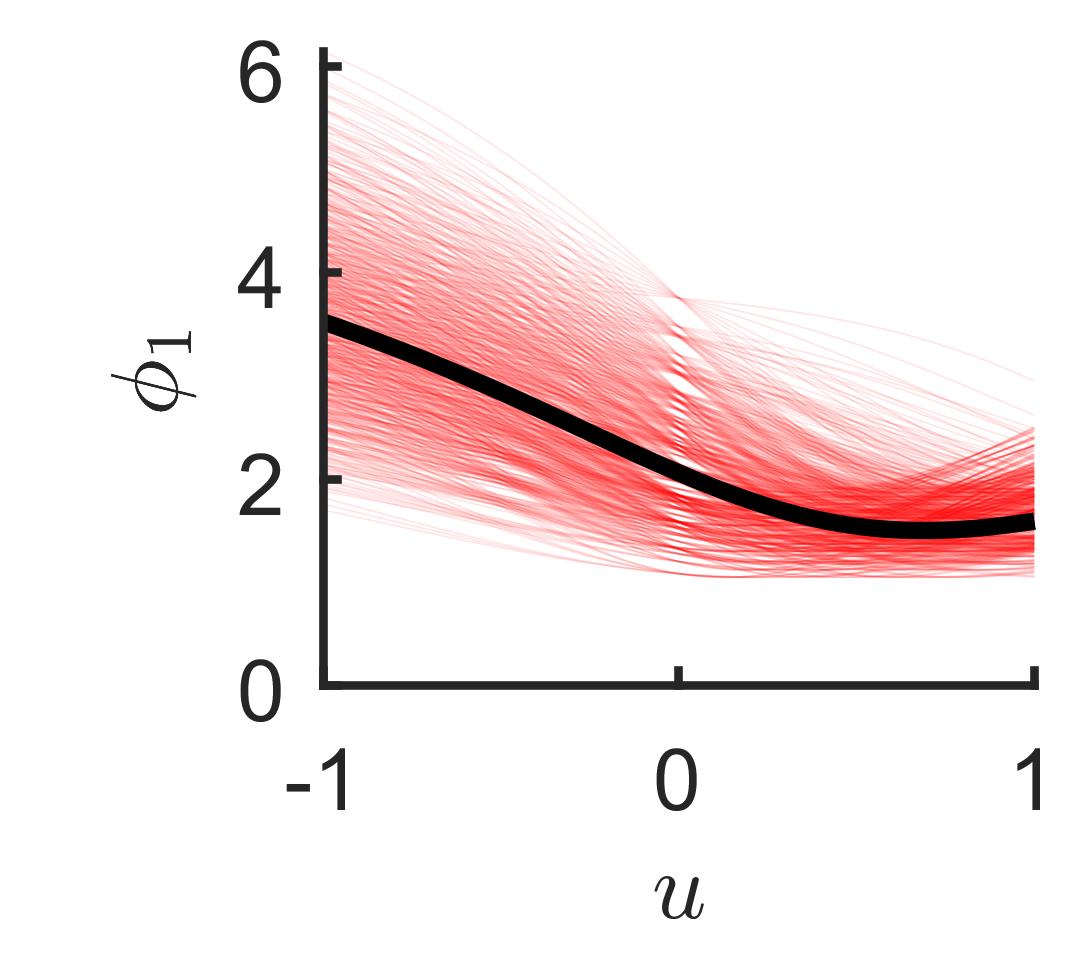}\hskip -0ex
				\includegraphics[trim={0.1cm 0.2cm 0.2cm  0.1cm },clip,width=3.5cm]{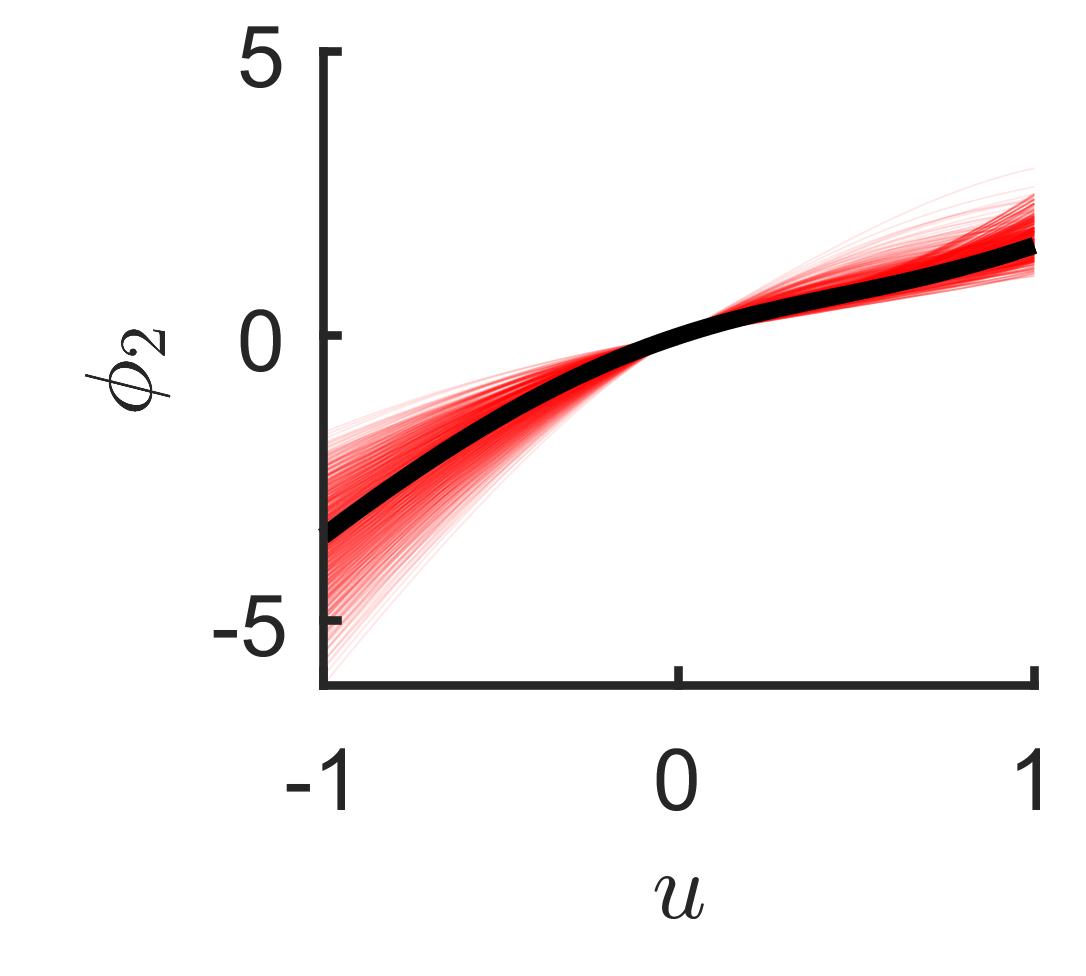}\hskip -0ex
				\includegraphics[trim={0.1cm 0.2cm 0.2cm  0.1cm },clip,width=3.5cm]{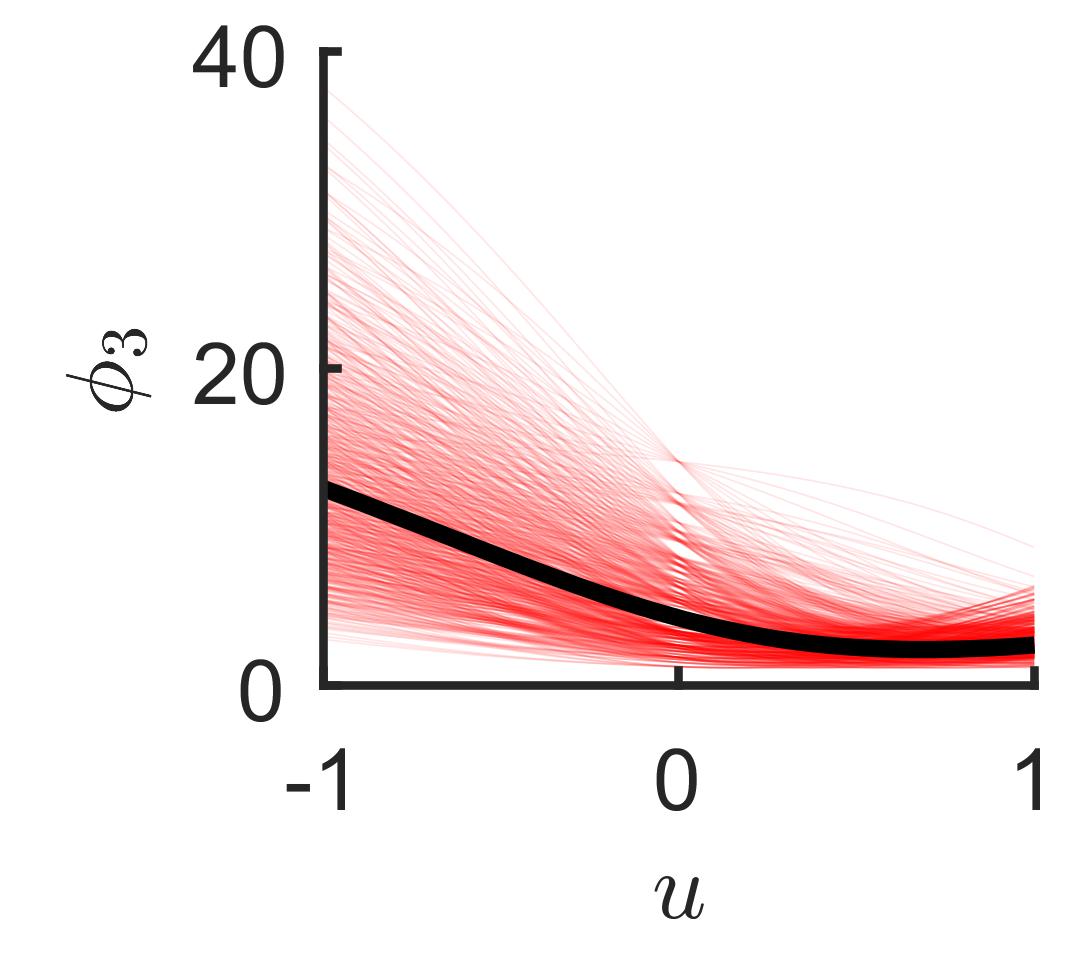}\hskip -0ex
				\includegraphics[trim={0.1cm 0.2cm 0.2cm  0.1cm },clip,width=3.5cm]{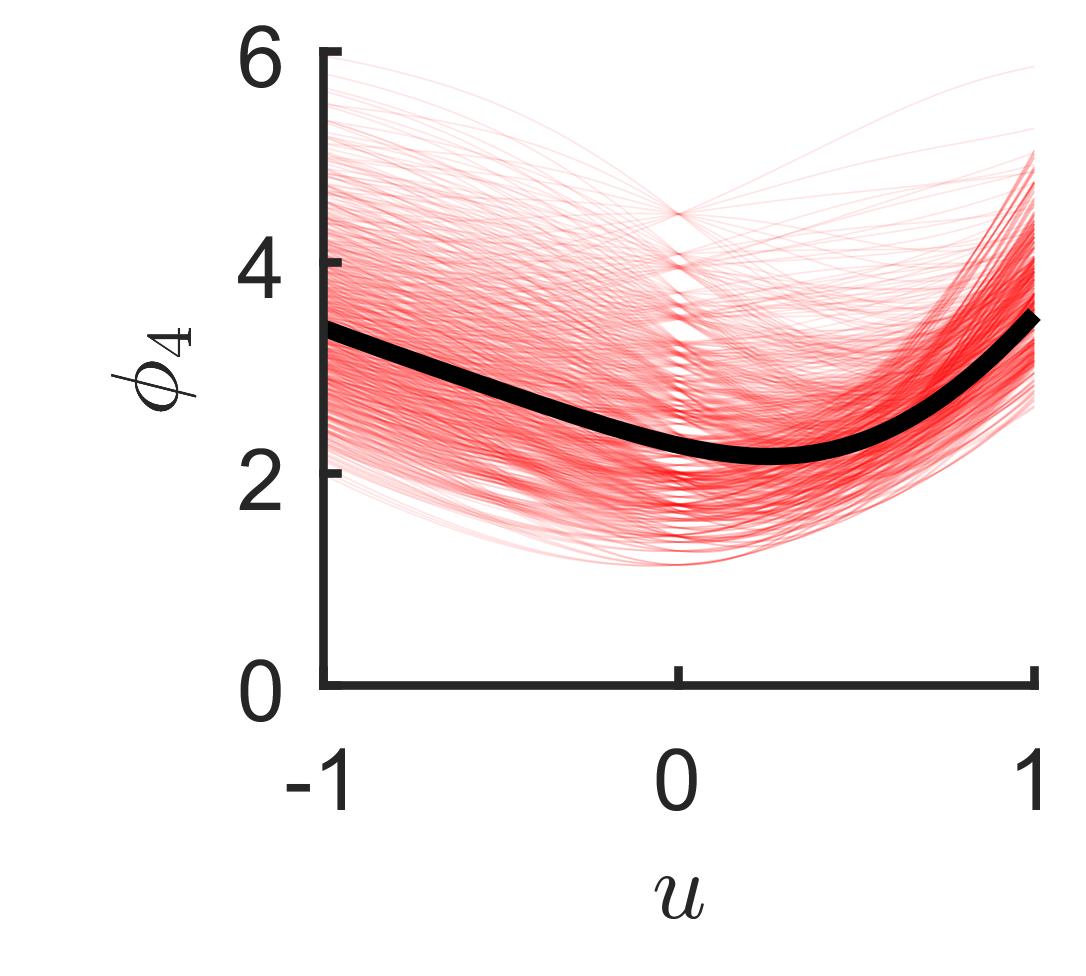}\\
				d) Function $\phi_1$, $\phi_2$, $\phi_3$, and  $\phi_4$,
			\end{minipage} 
			\medskip
		}	
		\textcolor{red}{\raisebox{0.5mm}{\rule{0.5cm}{0.05cm}}}   : Model, 
		\textcolor{black}{\raisebox{0.5mm}{\rule{0.5cm}{0.1cm}}} : Plant.
		\vspace{-2mm}
		\captionof{figure}{Sc.3: Graphical description of the RTO problems}
		\label{fig:5_26_Exemple_5_2_sc3_Plant_and_Model}
	\end{minipage} \\
	
	\medskip
	
	\begin{minipage}[h]{\linewidth}
		\vspace*{0pt}
		{\centering
			\includegraphics[trim={2.75cm 0.2cm 0.2cm  0.2cm },clip,width=3.45cm]{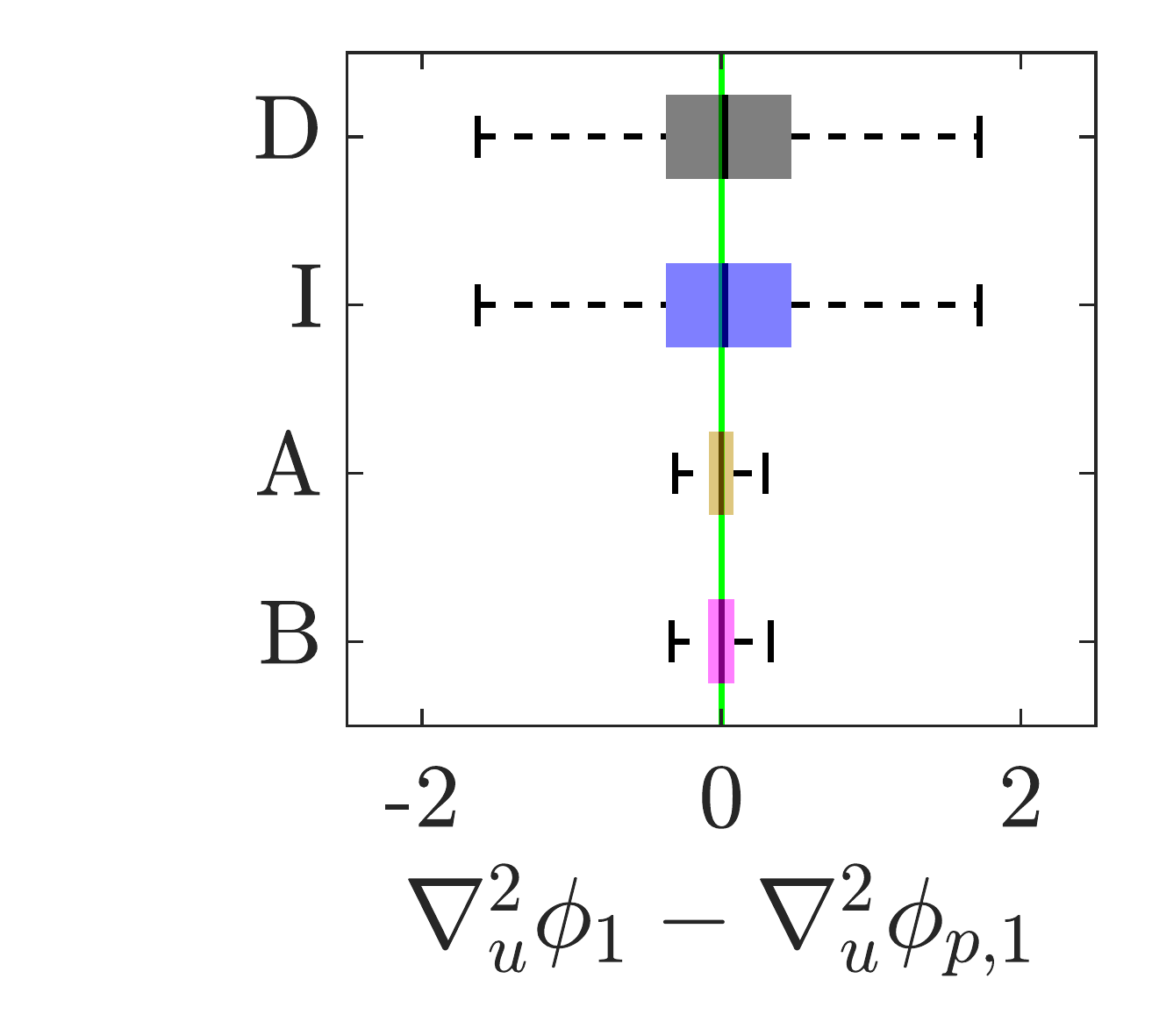}\hskip -0ex
			\includegraphics[trim={2.75cm 0.2cm 0.2cm  0.2cm },clip,width=3.45cm]{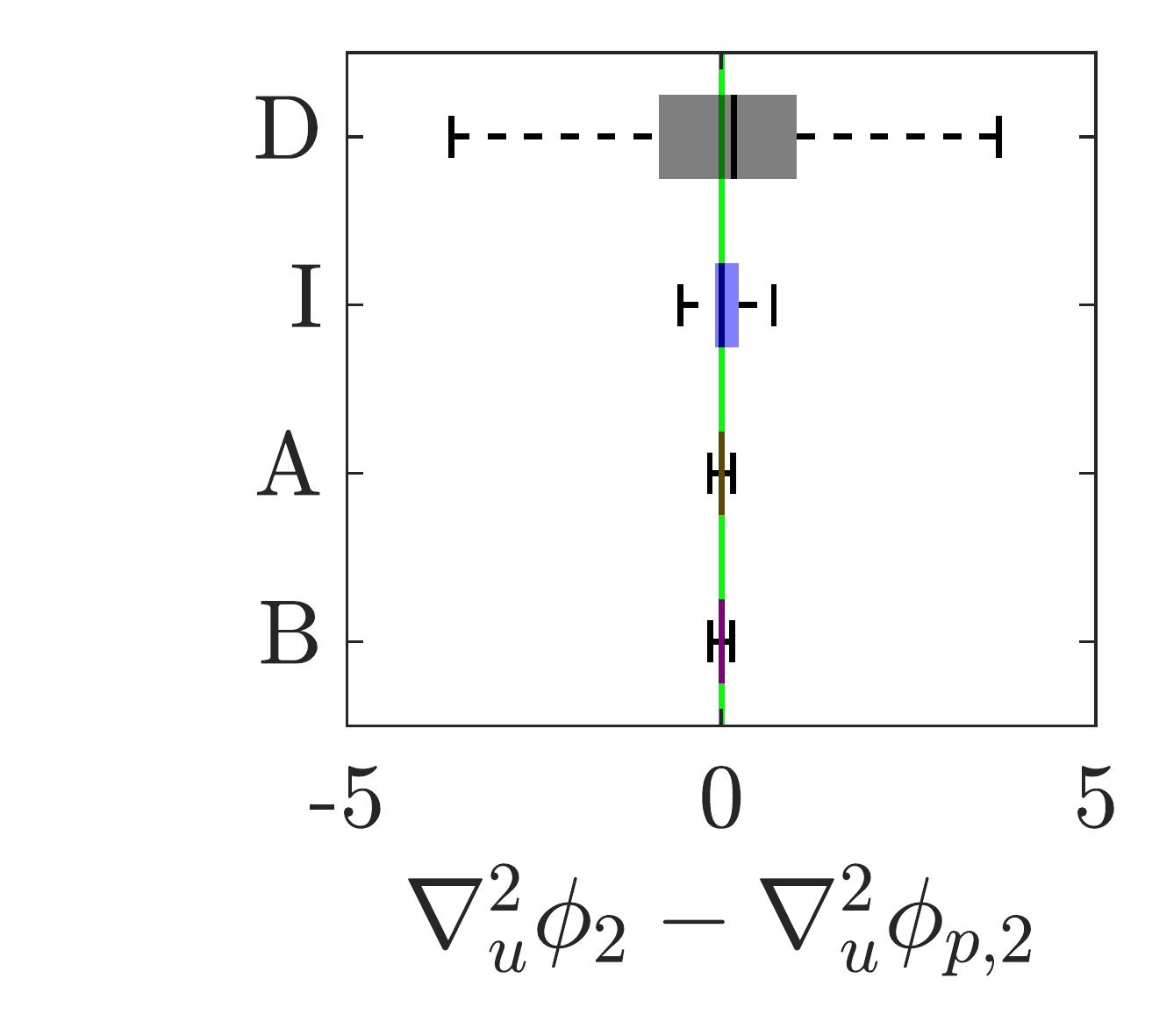}\hskip -0ex
			\includegraphics[trim={2.75cm 0.2cm 0.2cm  0.2cm },clip,width=3.45cm]{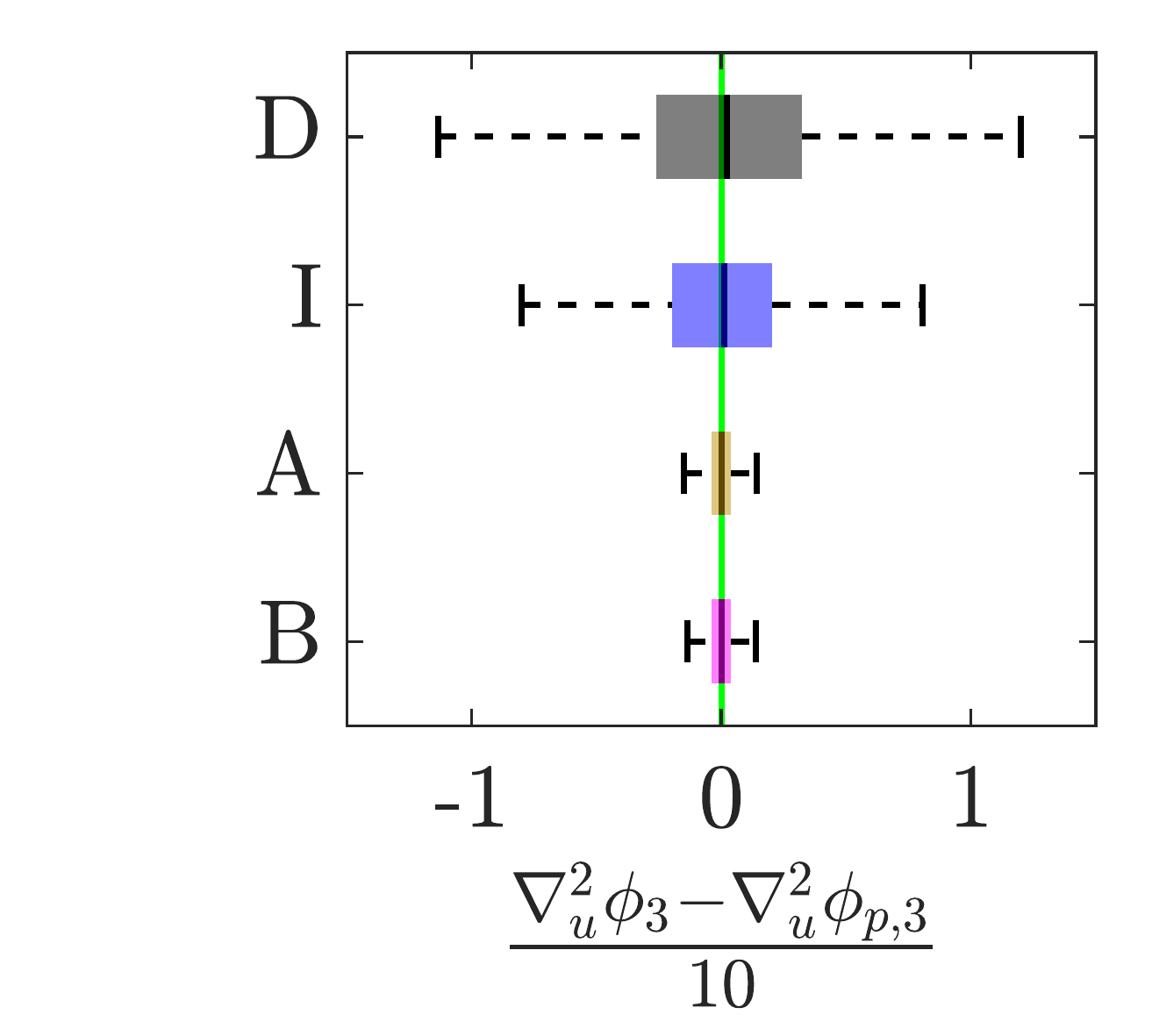}\hskip -0ex
			\includegraphics[trim={2.75cm 0.2cm 0.2cm  0.2cm },clip,width=3.45cm]{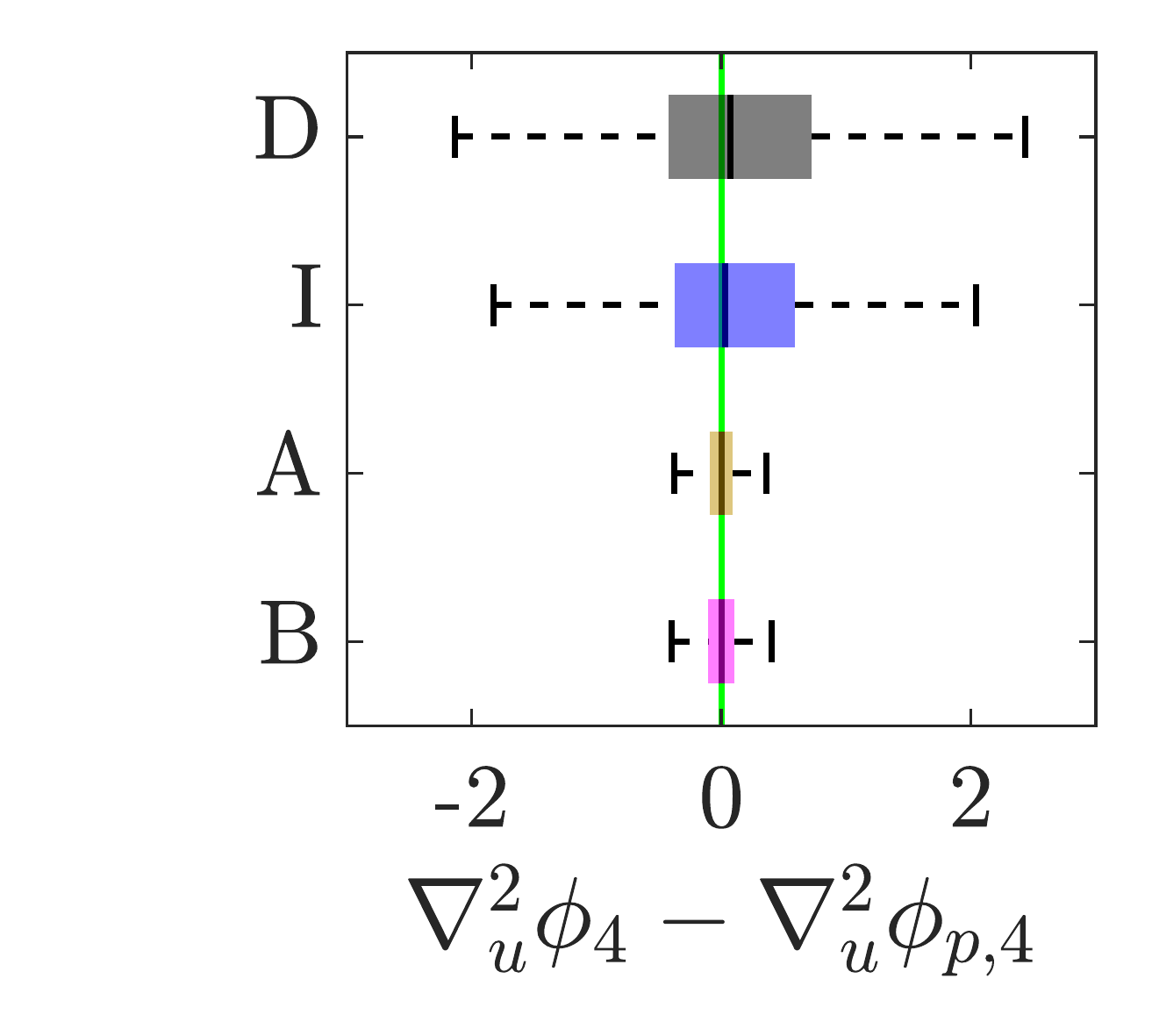}
		}
		\vspace{-2mm}
		\captionof{figure}{Sc.3:Statistical distributions of the prediction errors on the  Hessian of the plant's cost functions at the correction point for the structures D, I, A, and B.}
		\label{fig:5_27_Exemple_5_2_sc3_Results_1}
	\end{minipage} \\
	\begin{minipage}[h]{\linewidth}
		\vspace*{0pt}
		{\centering
			\includegraphics[trim={0.7cm 0.2cm 0.2cm  0.2cm },clip,width=3.45cm]{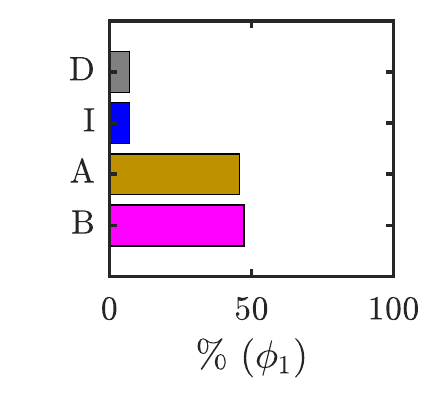}\hskip -0ex
			\includegraphics[trim={0.7cm 0.2cm 0.2cm  0.2cm },clip,width=3.45cm]{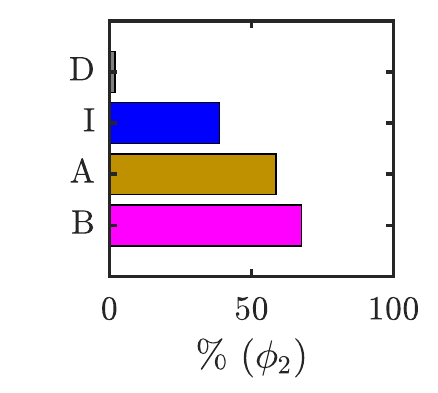}\hskip -0ex
			\includegraphics[trim={0.7cm 0.2cm 0.2cm  0.2cm },clip,width=3.45cm]{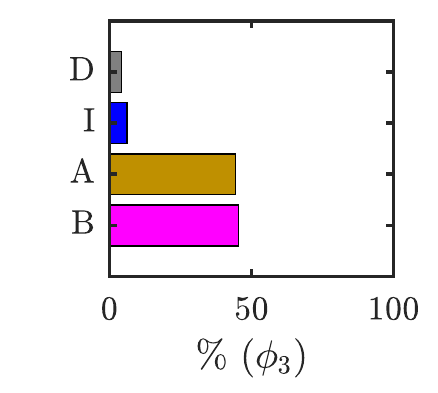}\hskip -0ex
			\includegraphics[trim={0.7cm 0.2cm 0.2cm  0.2cm },clip,width=3.45cm]{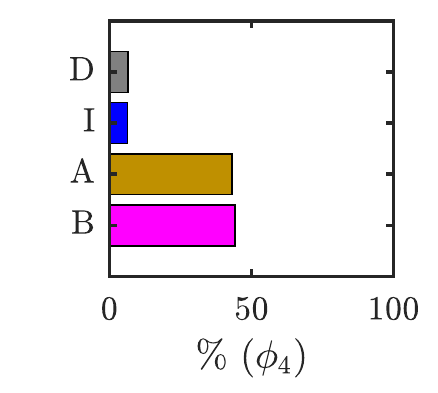} \\
		}
		\vspace{-2mm}
		\captionof{figure}{Sc.3: For each correction structure one gives here the percentage of cases for which no other structure provides better results  \textit{(if two structures provide the same best result then both take the point)}}
		\label{fig:5_28_Exemple_5_2_sc3_Results_2}
	\end{minipage} \\
\end{minipage} 
\noindent
\begin{minipage}[h]{\linewidth}
	\vspace*{0pt}
	{\centering
		\vspace{-5mm}
		\includegraphics[width=3.45cm]{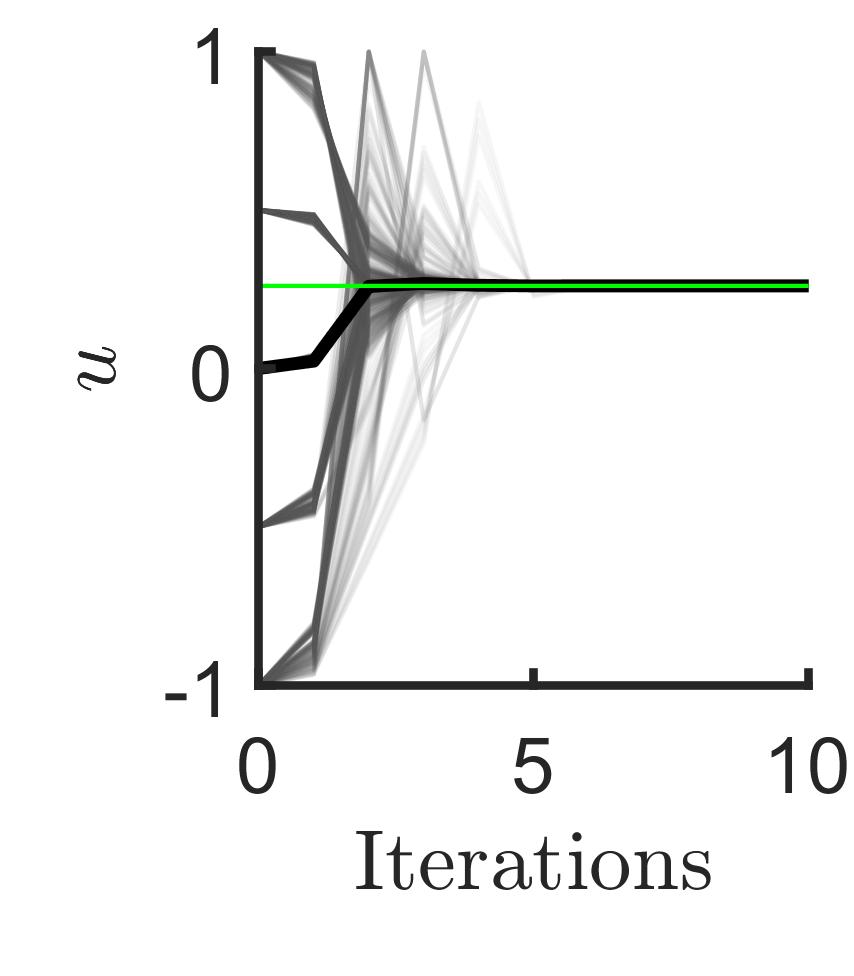}\hskip -0ex
		\includegraphics[width=3.45cm]{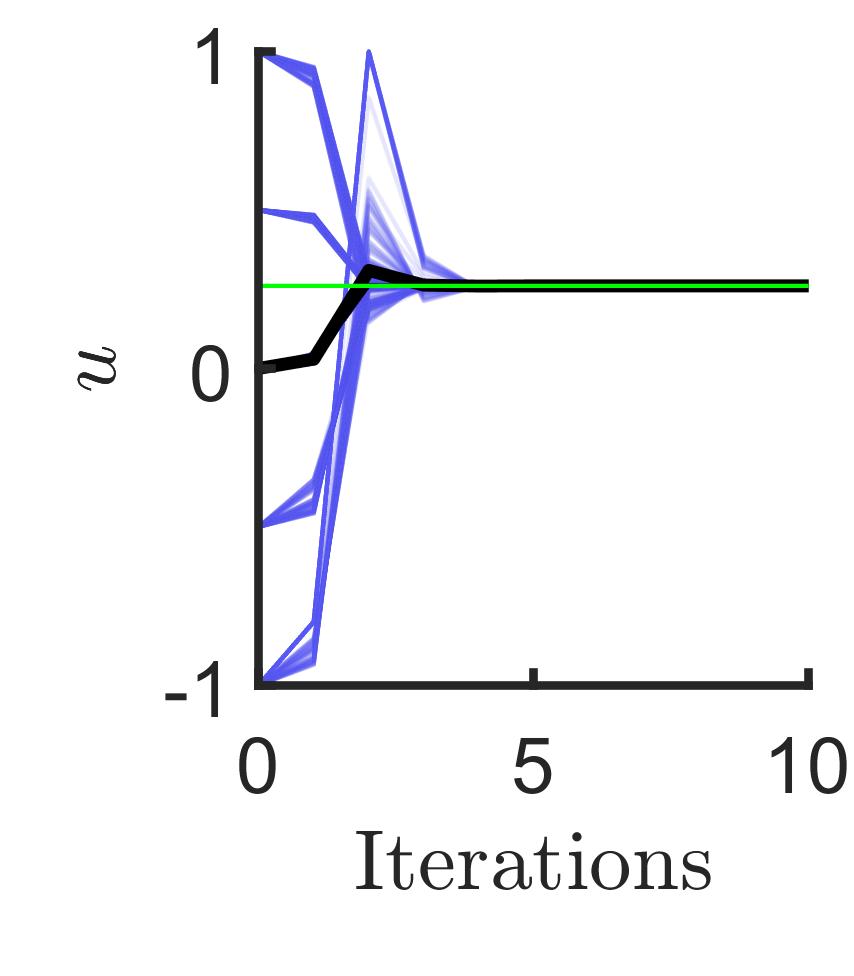}\hskip -0ex
		\includegraphics[width=3.45cm]{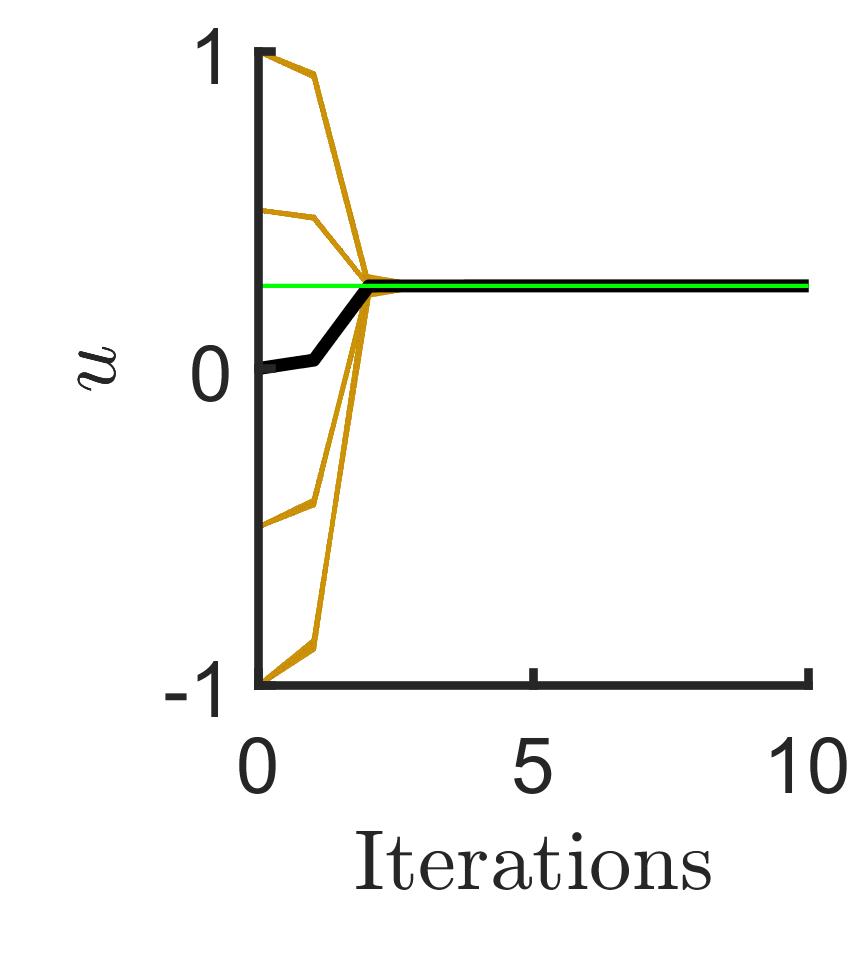}\hskip -0ex
		\includegraphics[width=3.45cm]{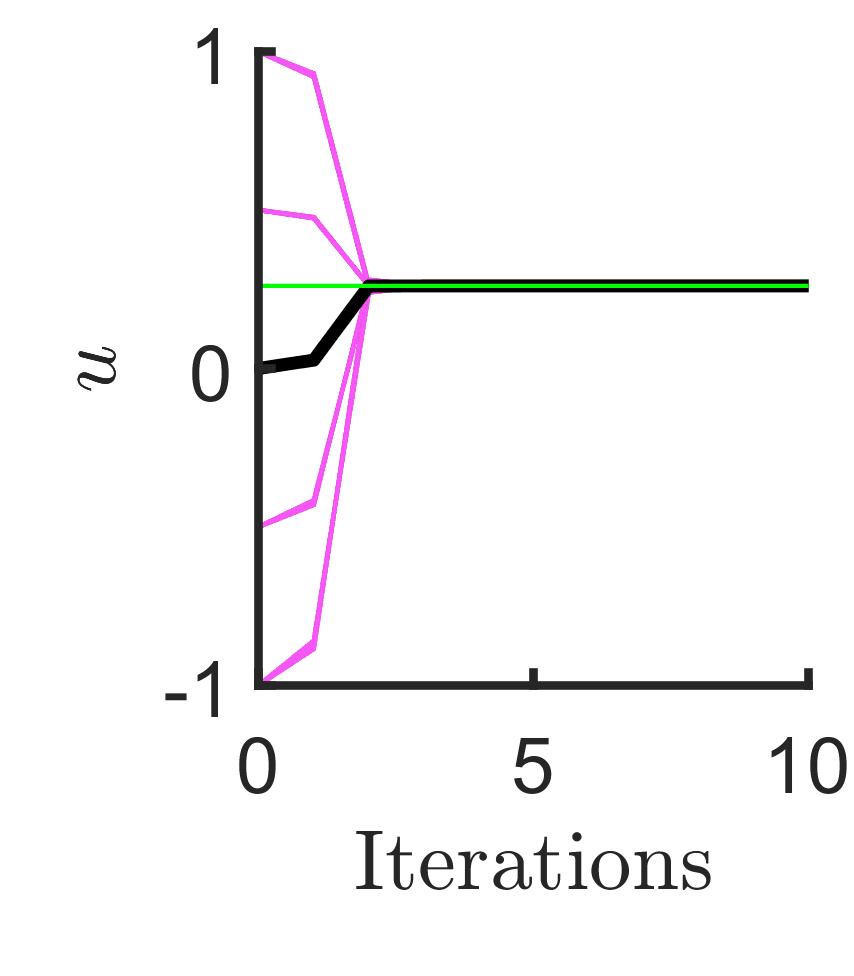} \\
		\vspace{-1.1cm}
		\includegraphics[width=3.45cm]{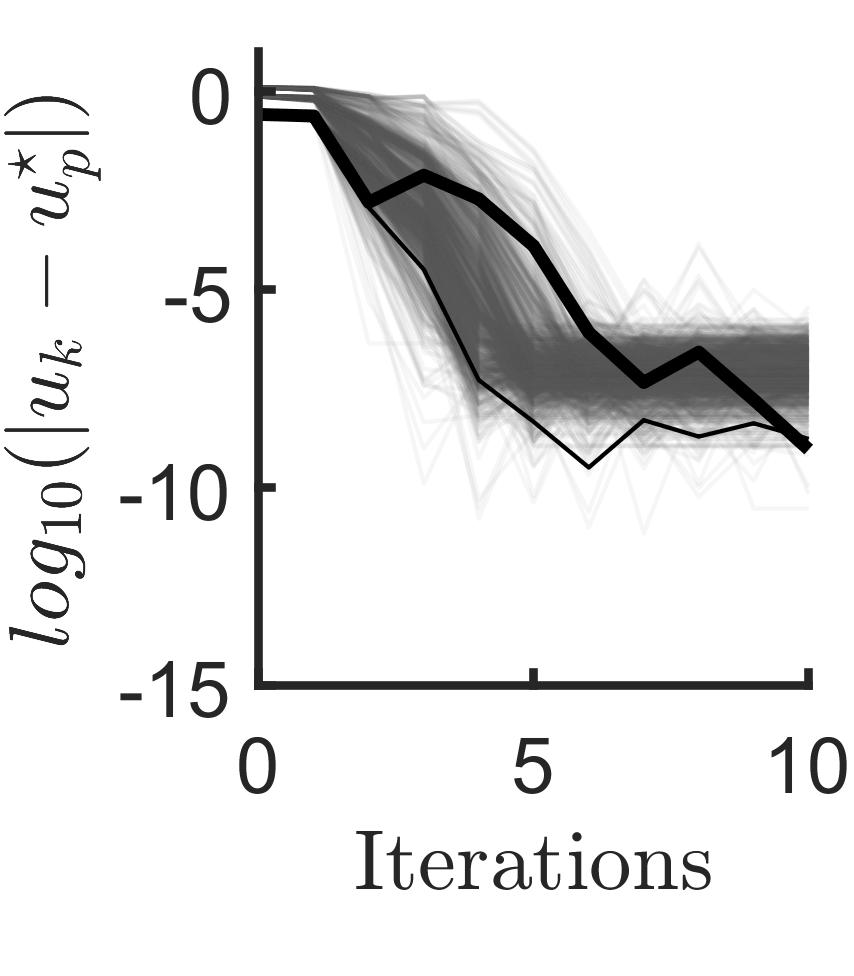}\hskip -0ex
		\includegraphics[width=3.45cm]{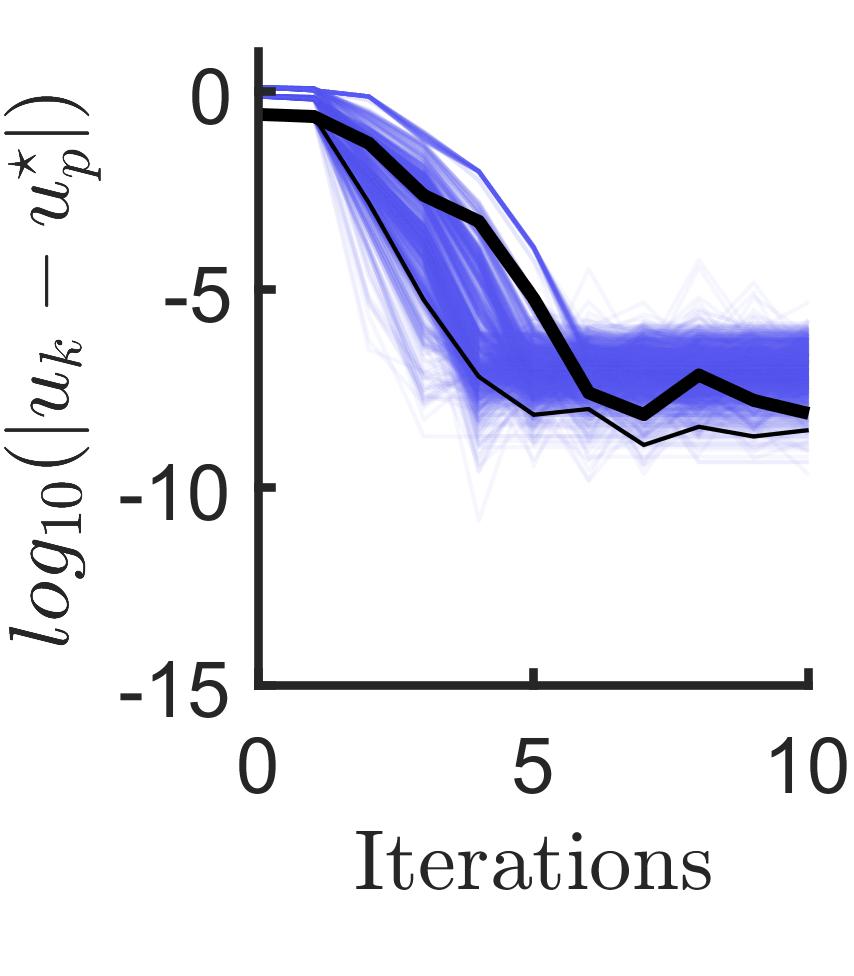}\hskip -0ex
		\includegraphics[width=3.45cm]{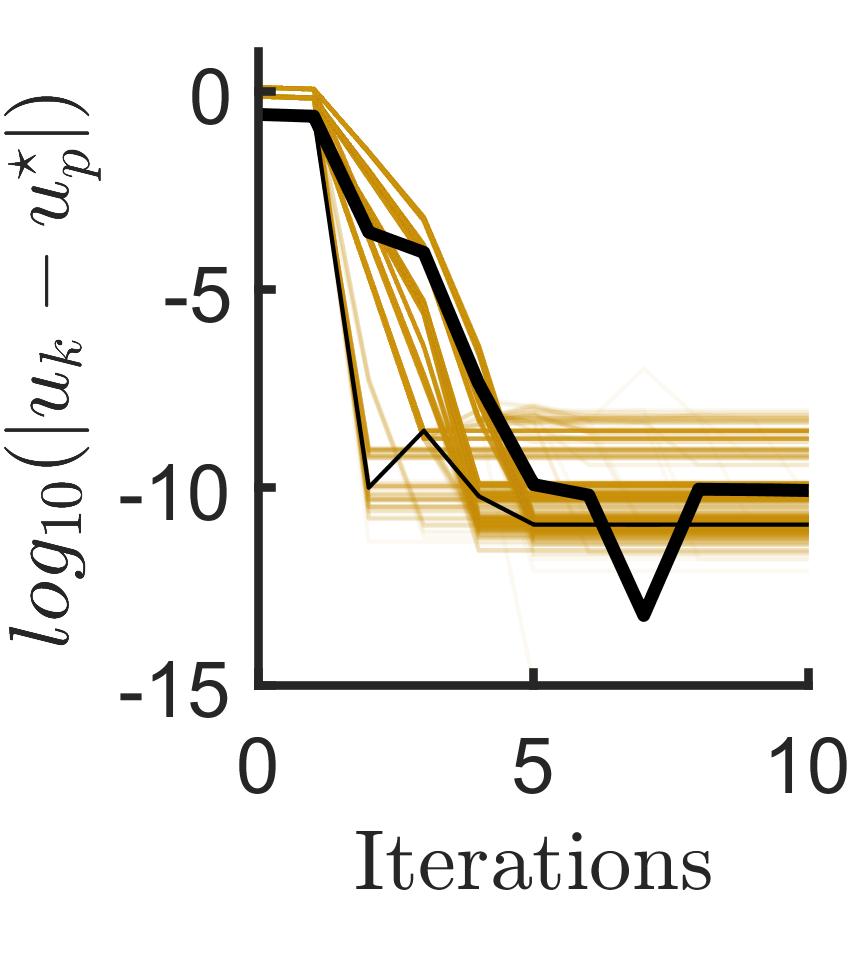}\hskip -0ex
		\includegraphics[width=3.45cm]{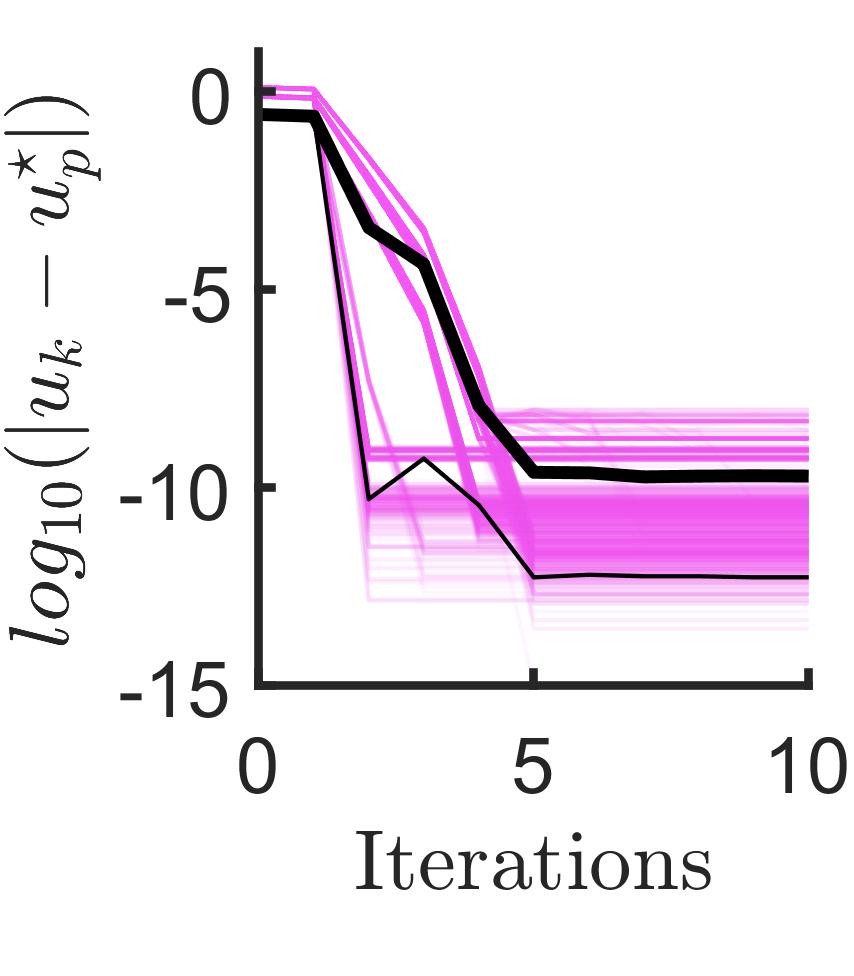}\\
	}
	\vspace{-4mm}
	a) \textbf{Sc.1} \textit{(Structure  D (\textcolor{gris_clair}{\raisebox{0.5mm}{\rule{0.3cm}{0.05cm}}}), I (\textcolor{blue}{\raisebox{0.5mm}{\rule{0.3cm}{0.05cm}}}), 
		A (\textcolor{gold}{\raisebox{0.5mm}{\rule{0.3cm}{0.05cm}}}),
		B (\textcolor{magenta}{\raisebox{0.5mm}{\rule{0.3cm}{0.05cm}}}), 
		mean  (\textcolor{black}{\raisebox{0.5mm}{\rule{0.3cm}{0.1cm}}}),
		median (\textcolor{black}{\raisebox{0.5mm}{\rule{0.3cm}{0.05cm}}}))}\\
	{\centering
		\includegraphics[width=3.45cm]{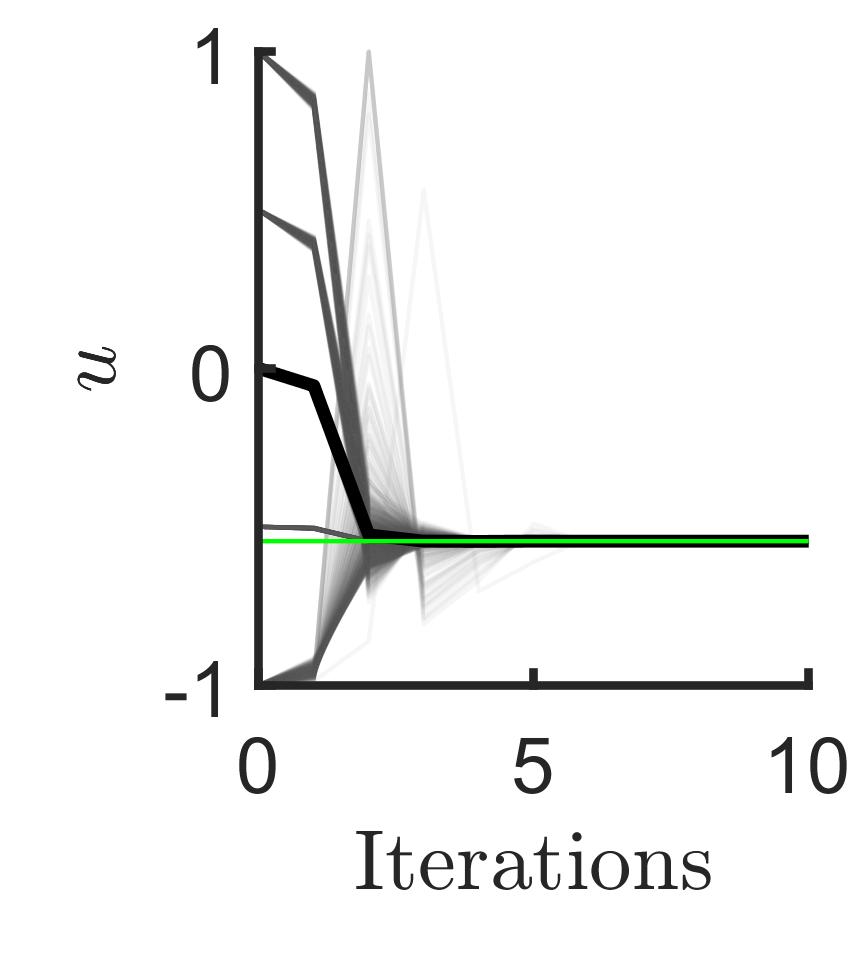}\hskip -0ex
		\includegraphics[width=3.45cm]{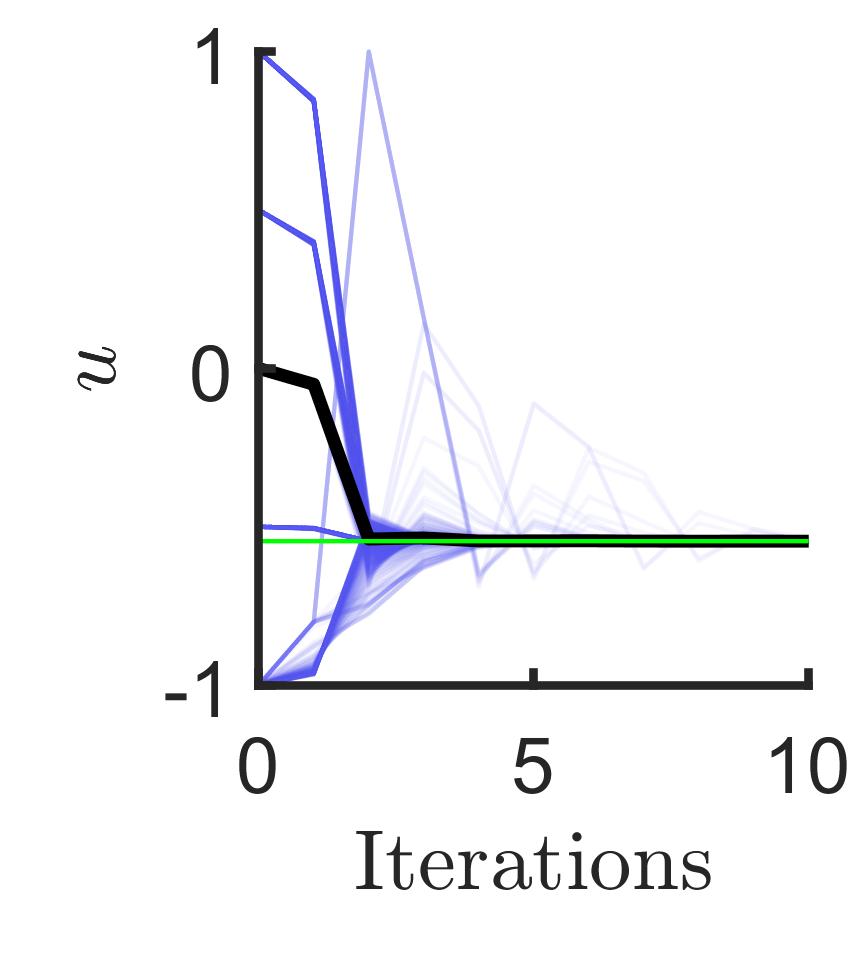}\hskip -0ex
		\includegraphics[width=3.45cm]{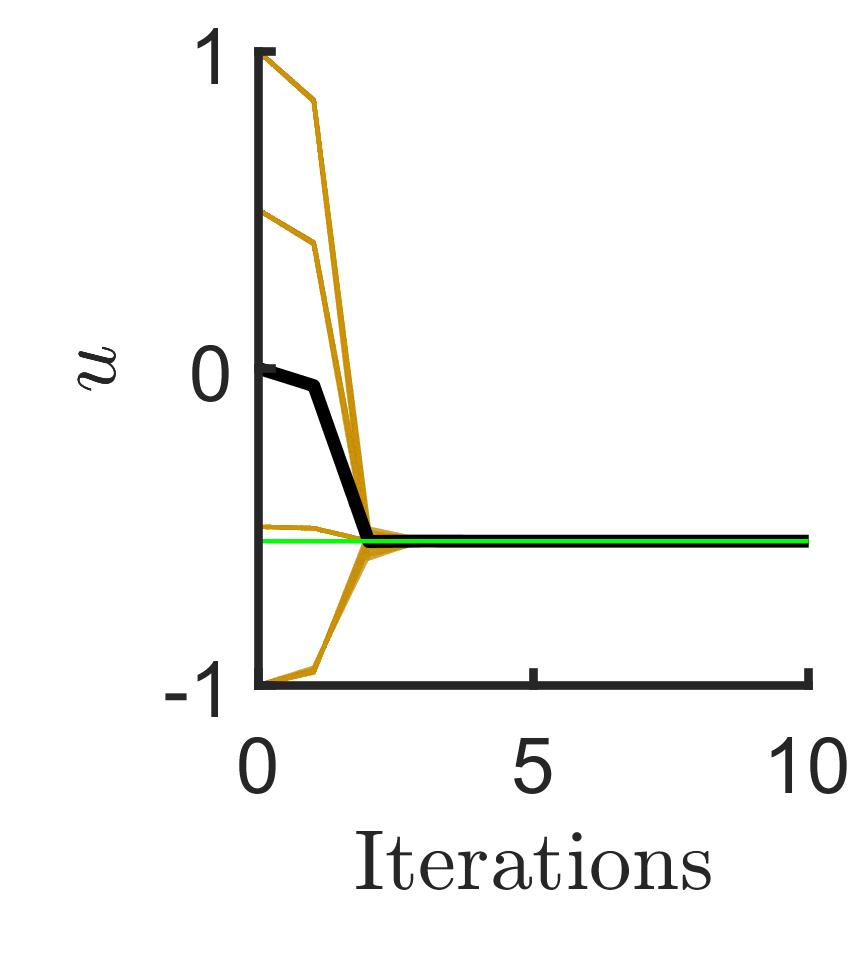}\hskip -0ex
		\includegraphics[width=3.45cm]{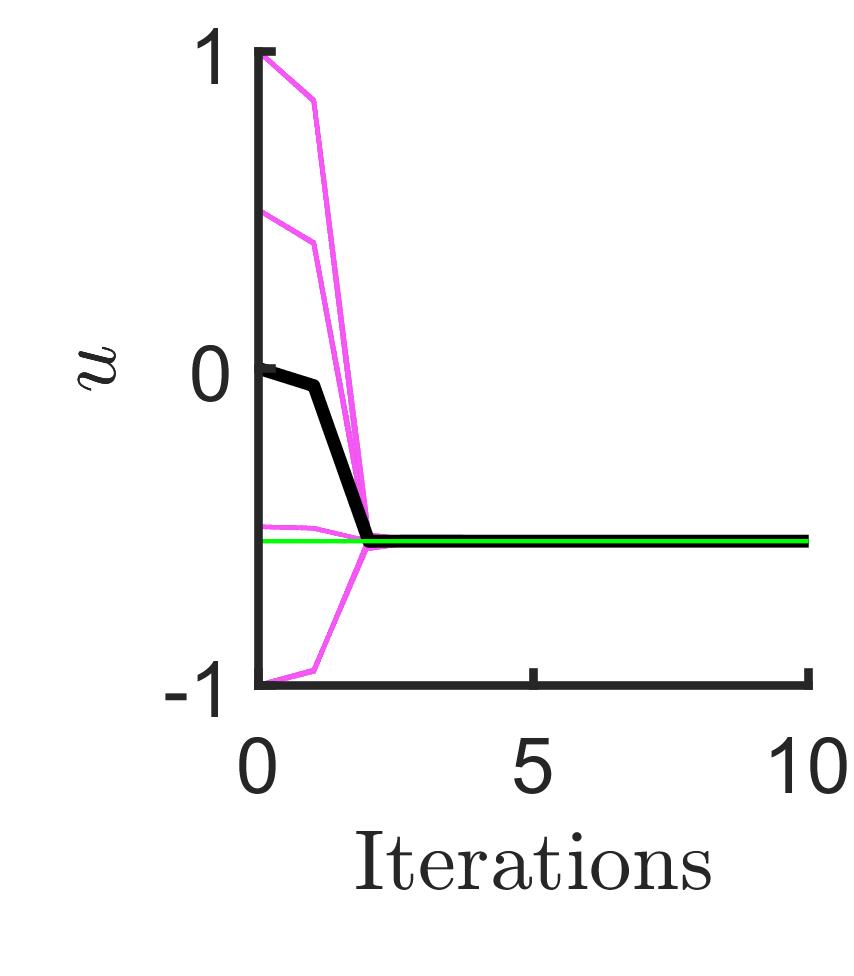} \\
		\vspace{-1.1cm}
		\includegraphics[width=3.45cm]{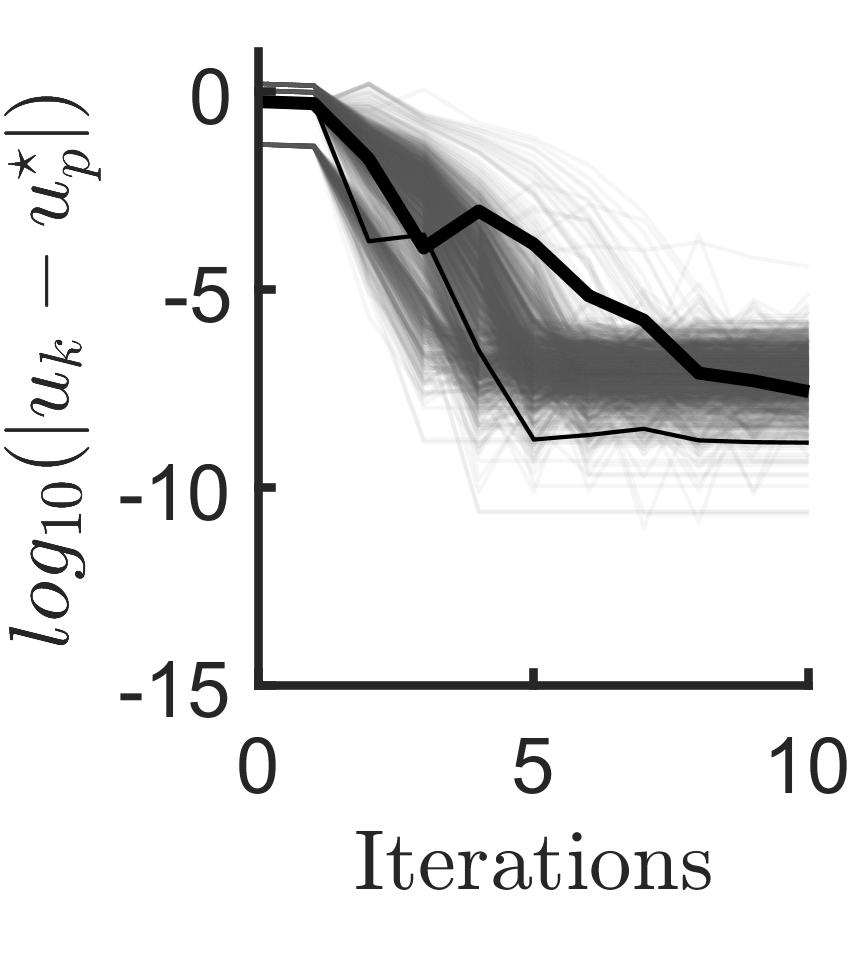}\hskip -0ex
		\includegraphics[width=3.45cm]{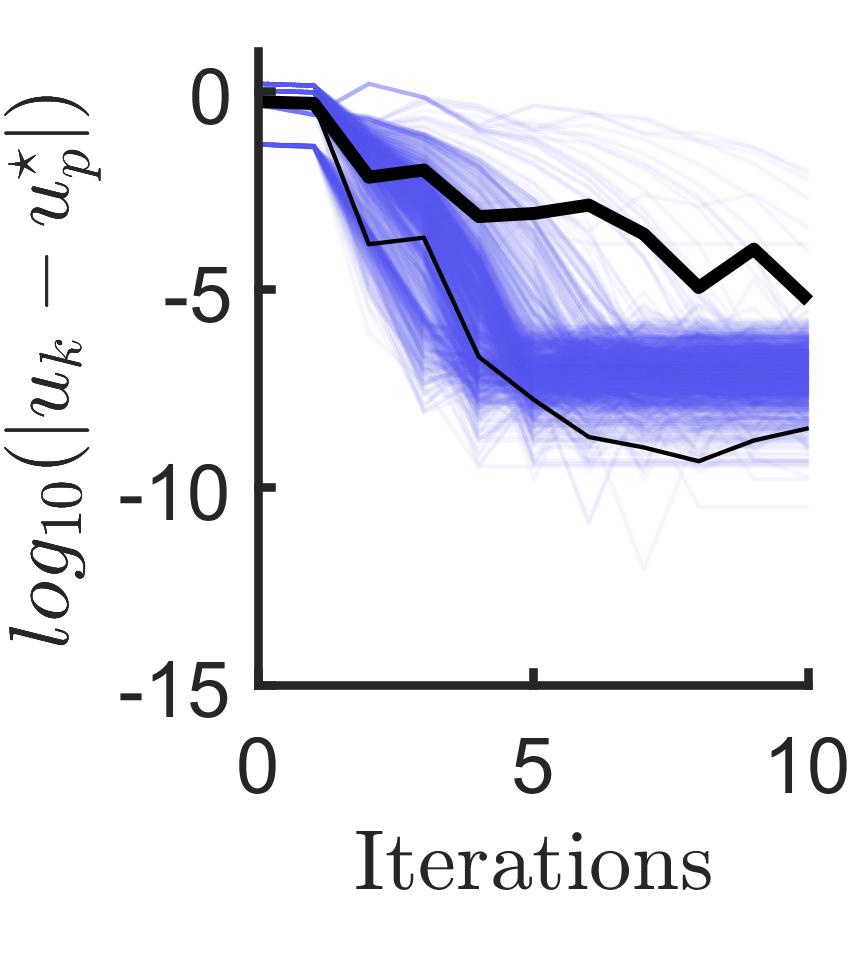}\hskip -0ex
		\includegraphics[width=3.45cm]{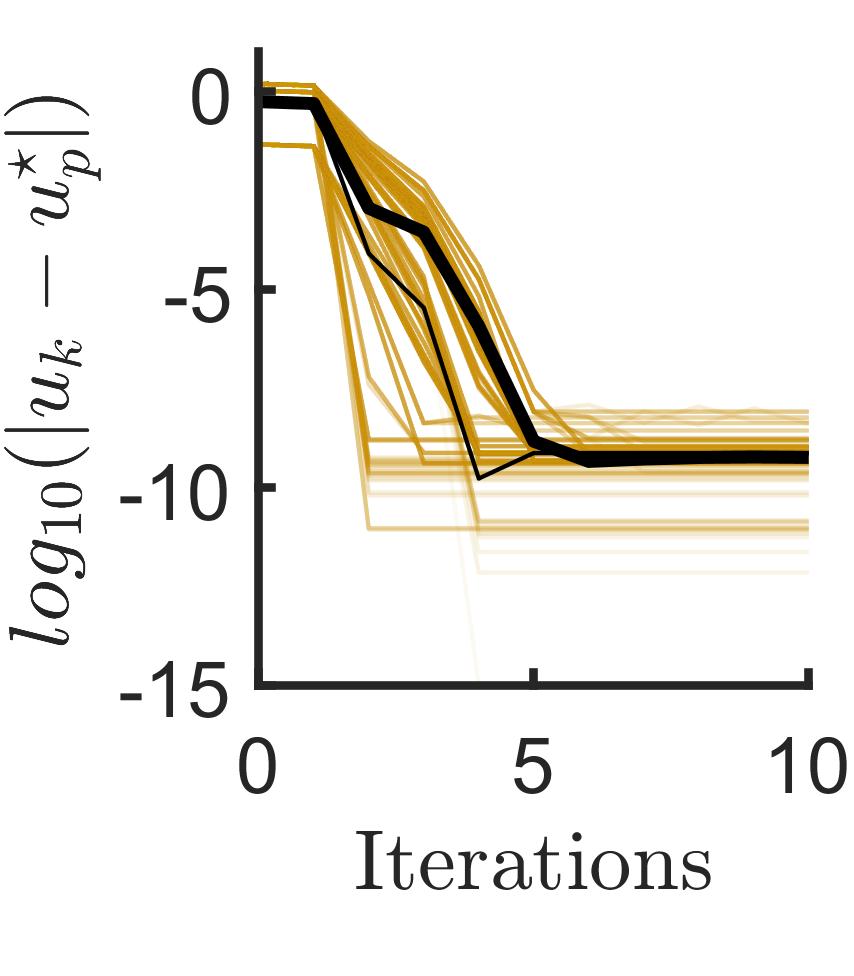}\hskip -0ex
		\includegraphics[width=3.45cm]{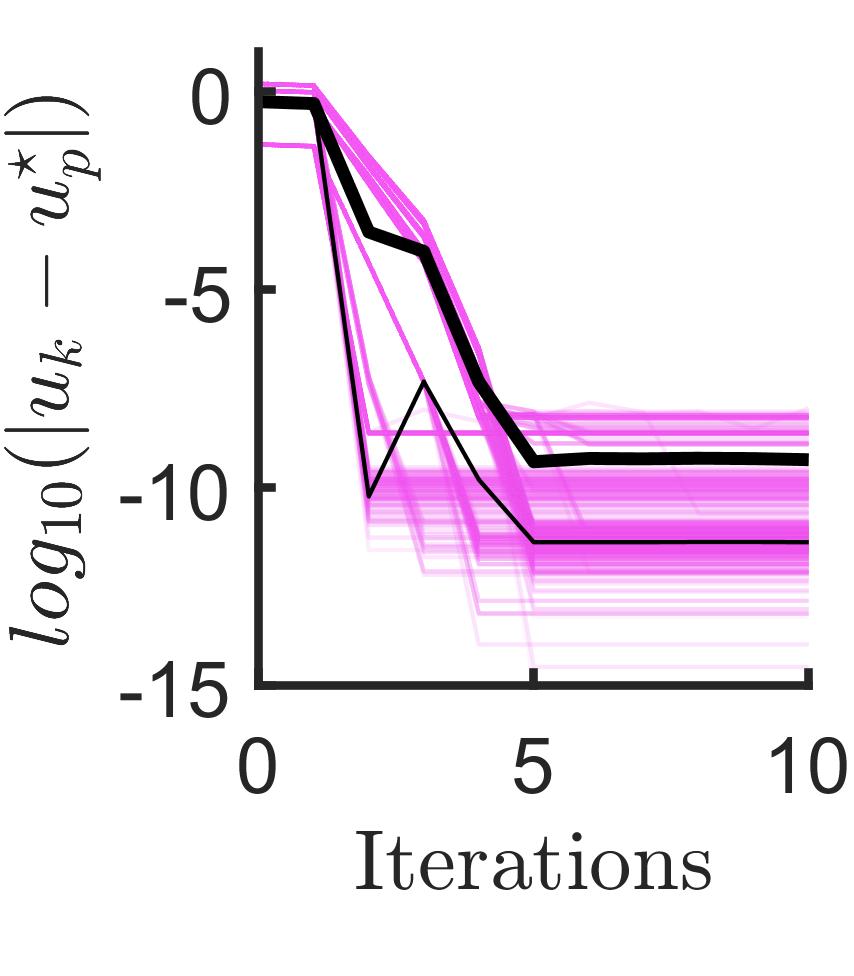}\\
	}
	\vspace{-4mm}
	b) \textbf{Sc.2} \textit{(Structure  D (\textcolor{gris_clair}{\raisebox{0.5mm}{\rule{0.3cm}{0.05cm}}}), I (\textcolor{blue}{\raisebox{0.5mm}{\rule{0.3cm}{0.05cm}}}), 
		A (\textcolor{gold}{\raisebox{0.5mm}{\rule{0.3cm}{0.05cm}}}),
		B (\textcolor{magenta}{\raisebox{0.5mm}{\rule{0.3cm}{0.05cm}}}), 
		mean  (\textcolor{black}{\raisebox{0.5mm}{\rule{0.3cm}{0.1cm}}}),
		median (\textcolor{black}{\raisebox{0.5mm}{\rule{0.3cm}{0.05cm}}}))}\\
	{\centering
		\includegraphics[width=3.45cm]{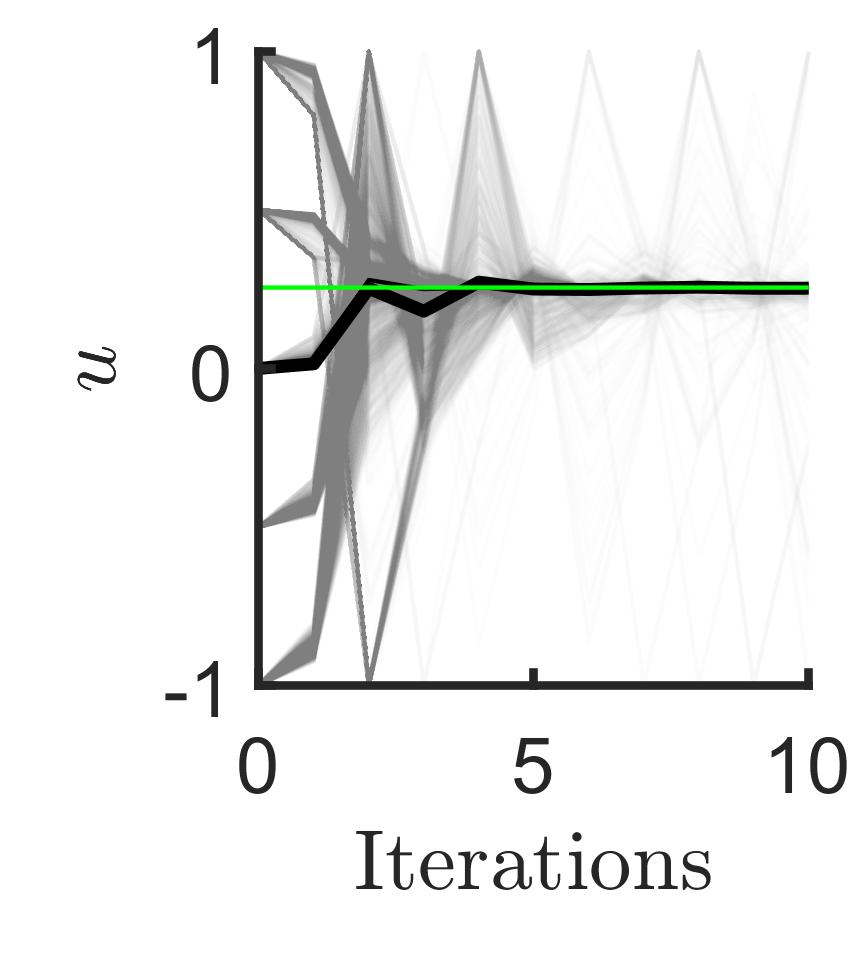}\hskip -0ex
		\includegraphics[width=3.45cm]{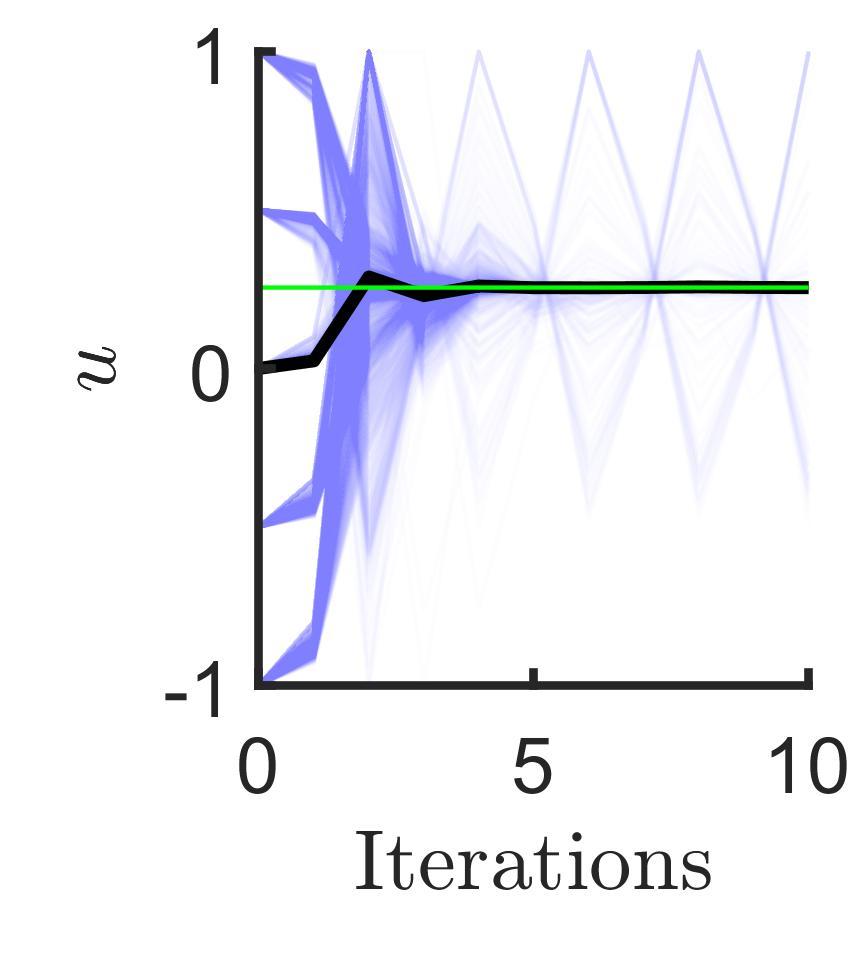}\hskip -0ex
		\includegraphics[width=3.45cm]{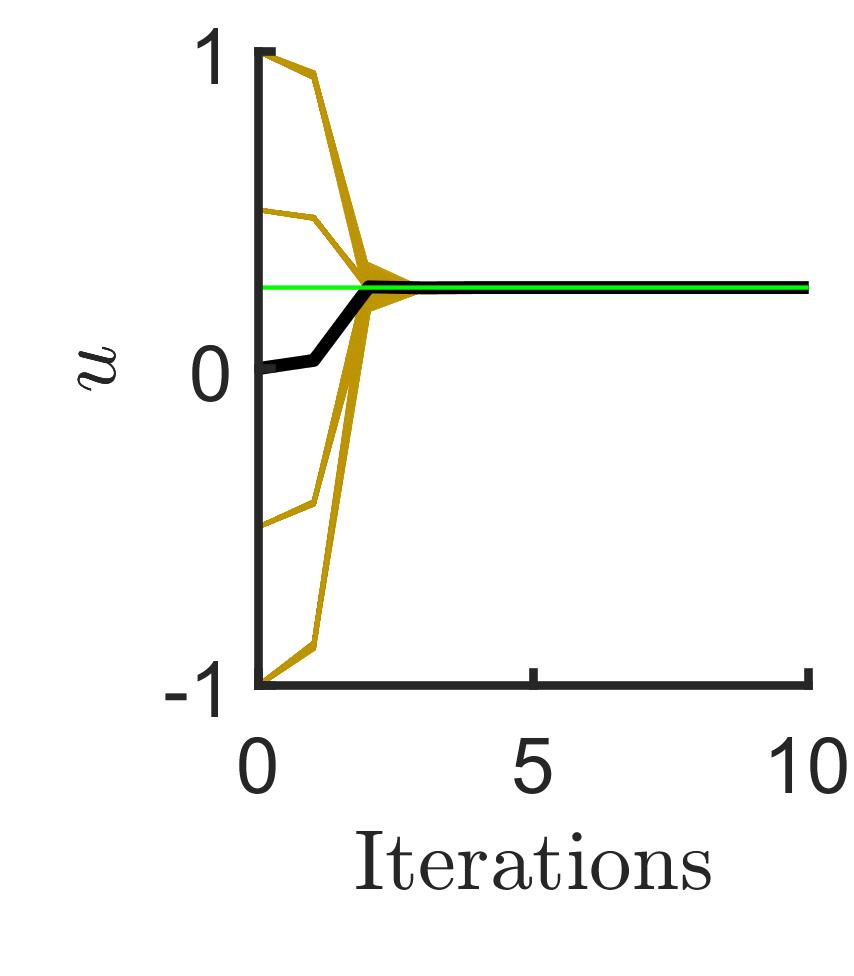}\hskip -0ex
		\includegraphics[width=3.45cm]{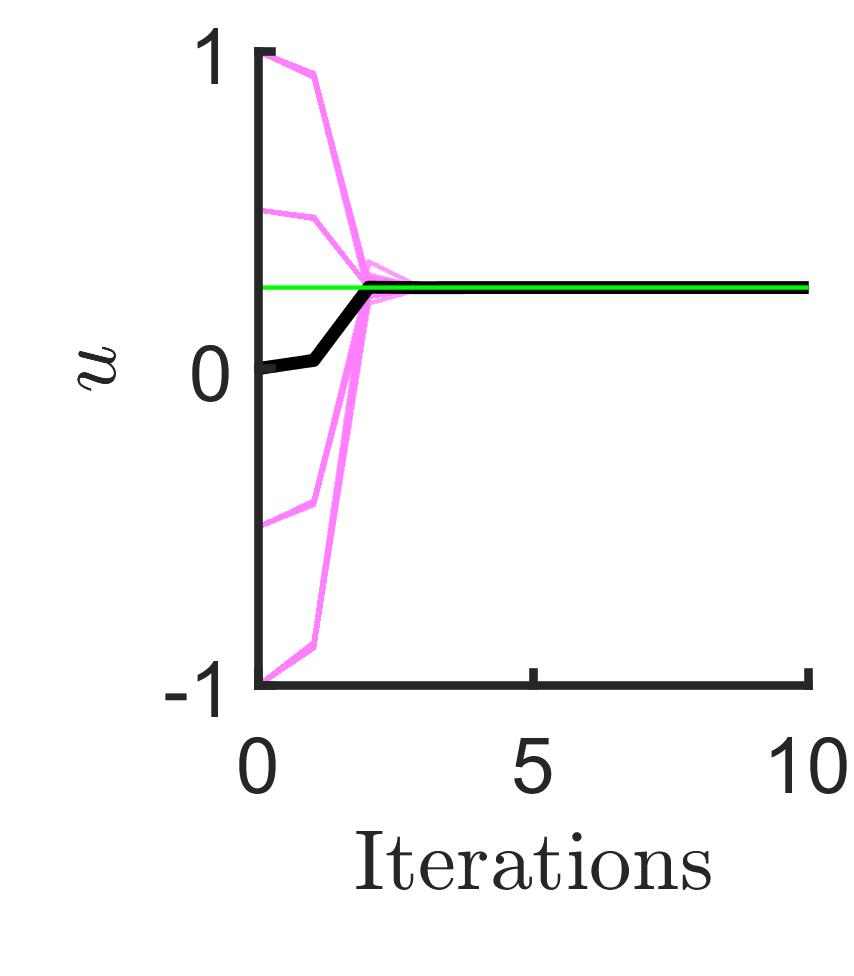} \\
		\vspace{-1.1cm}
		\includegraphics[width=3.45cm]{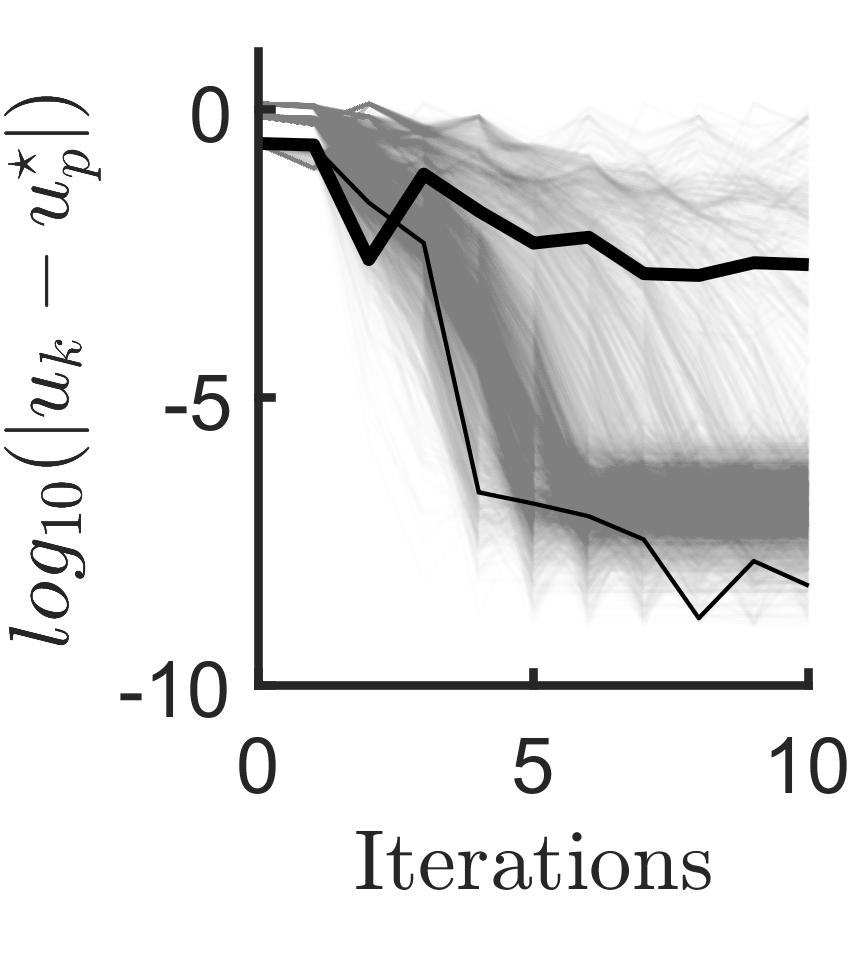}\hskip -0ex
		\includegraphics[width=3.45cm]{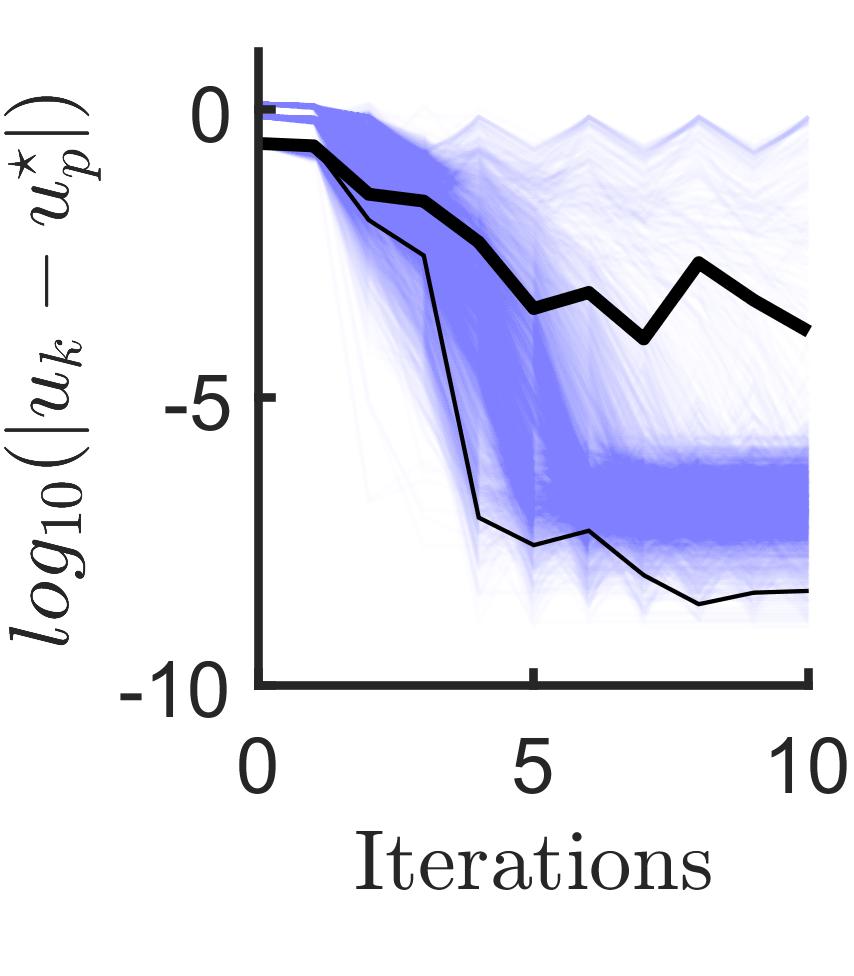}\hskip -0ex
		\includegraphics[width=3.45cm]{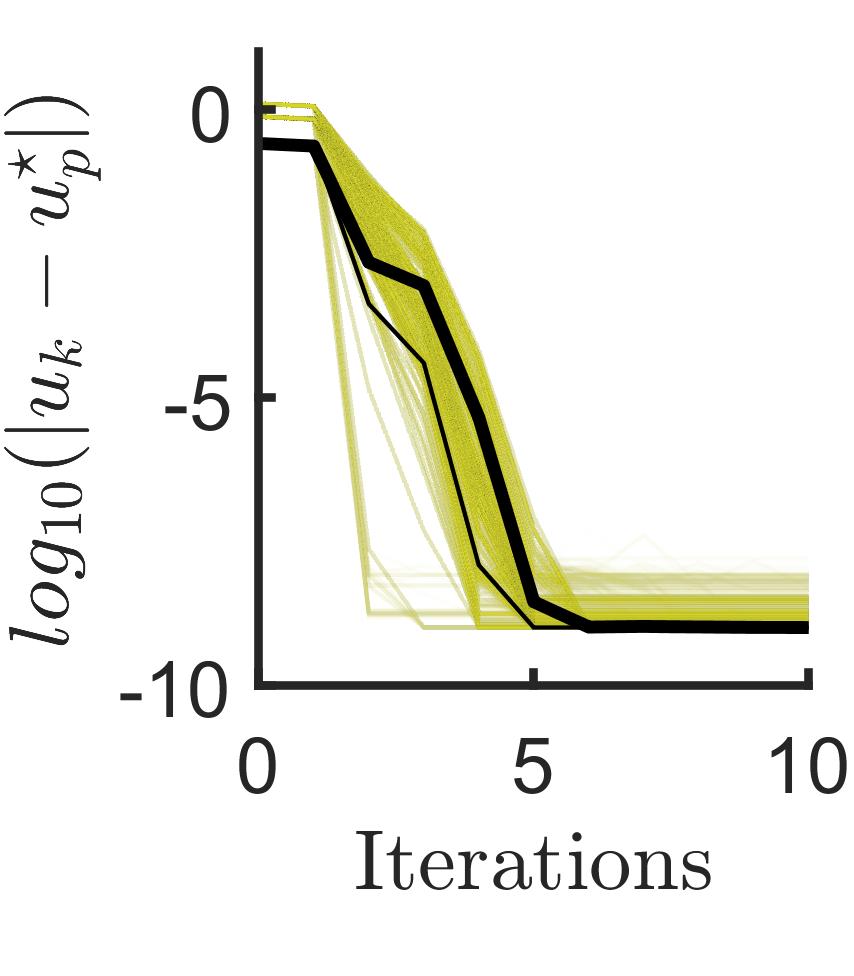}\hskip -0ex
		\includegraphics[width=3.45cm]{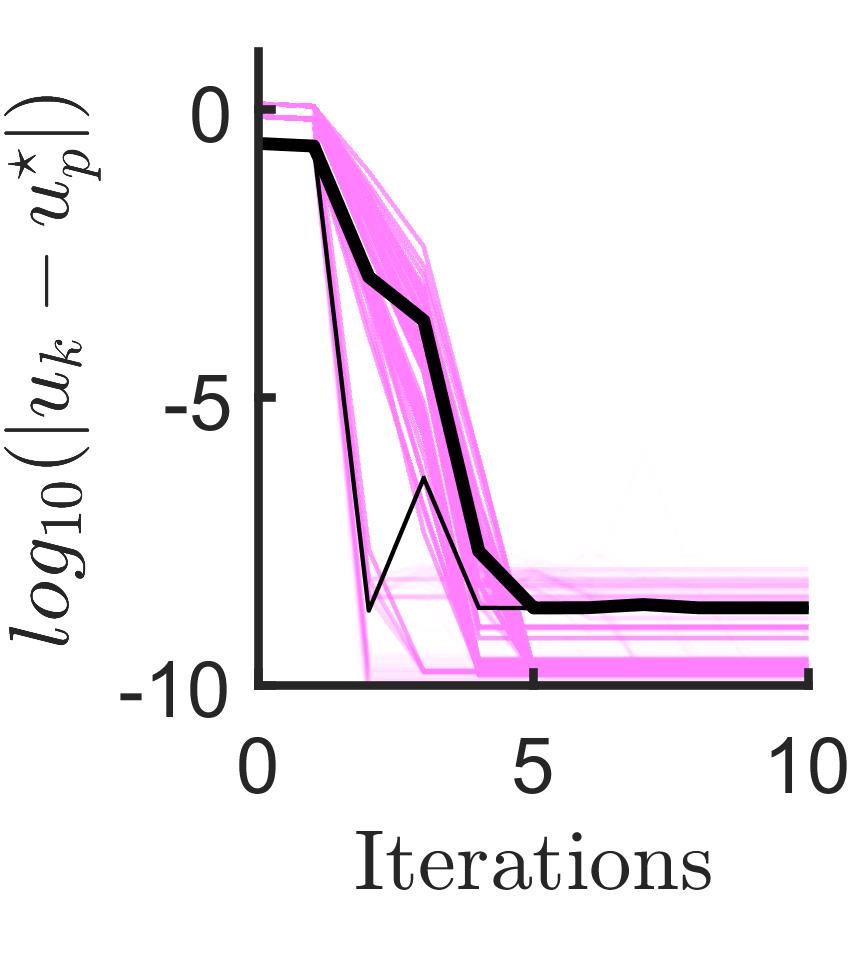}\\
	}
	\vspace{-4mm}	
	c) \textbf{Sc.3} \textit{(Structure  D (\textcolor{gris_clair}{\raisebox{0.5mm}{\rule{0.3cm}{0.05cm}}}), I (\textcolor{blue}{\raisebox{0.5mm}{\rule{0.3cm}{0.05cm}}}), 
		A (\textcolor{gold}{\raisebox{0.5mm}{\rule{0.3cm}{0.05cm}}}),
		B (\textcolor{magenta}{\raisebox{0.5mm}{\rule{0.3cm}{0.05cm}}}), 
		mean  (\textcolor{black}{\raisebox{0.5mm}{\rule{0.3cm}{0.1cm}}}),
		median (\textcolor{black}{\raisebox{0.5mm}{\rule{0.3cm}{0.05cm}}}))}
	\captionof{figure}{Simulation results for the structures D, I, A and B }
	\label{fig:5_29_Exemple_5_2_Results_3}
\end{minipage}
\clearpage

\subsection{Study 3: ``Wrong'' but consistent model}
\label{sec:5_5_3_Model_structurellement_Faux_mais_coherent}

In this section, one shows that even if the structure of a model is ``very wrong'' it can be consistent if the inputs of the SMs are \textit{sufficient} to predict the outputs of the SPs. In this case, the structure B can be used to correct the model.

One considers the theoretical RTO problem illustrated on the Figures~\ref{fig:5_30_ModelPlantCoherent} and \ref{fig:5_31_Exemple_5_3_Plant_and_Model} and using the following functions: 
\begin{align}
	f^{(1)}   := \ & \theta_1 - \theta_2z^{(1)} + \theta_3(z^{(1)})^2, &  z^{(1)} := \ & u, \nonumber  \\
	f_p^{(1)} := \ & 1 - \left(\begin{array}{c} 1 \\ 0.01 \end{array}\right)^{\rm T} \bm{z}_p^{(1)} + (\bm{z}_p^{(1)})^{\rm T}  \left(\begin{array}{cc} 0.01 & 0.08 \\ 0.08 & 0.002 \end{array}\right) \bm{z}_p^{(1)}, &  \bm{z}_p^{(1)} := \ & \left(\begin{array}{c} u \\ y^{(2)} \end{array}\right), \nonumber \\
	f^{(2)}   := \ & \theta_4 - \theta_5 z^{(2)}   + \theta_6 (z^{(2)})^2, &  z^{(2)}   := \ & y^{(1)},   \nonumber  \\
	f_p^{(2)} := \ & 2 - 2z_p^{(2)} + 2(z_p^{(2)})^2, &  z_p^{(2)} := \ & y_p^{(1)}, \nonumber 
\end{align}
\begin{align*}
	\phi_1(u,\bm{y}) := \ & y_{(2)},   & \phi_2(u,\bm{y}) := \ & uy_{(2)}, \\
	\phi_3(u,\bm{y}) := \ & y_{(2)}^2, & \phi_4(u,\bm{y}) := \ & y_{(2)} + u^2 + 0.5uy_{(2)} + 0.05y_{(2)}^2.
\end{align*}
where $u\in[-1,1]$, $\theta_1\in[0.75,1.25]$, $\theta_2\in[-1.25,-0.75]$, $\theta_3\in[0,0.02]$, $\theta_4\in[1.75,2.25]$, $\theta_5\in[-2.25,-1.75]$, and $\theta_6\in[1.8,2.2]$.

\noindent
\begin{minipage}[h]{\linewidth}
	\includegraphics[width=14cm]{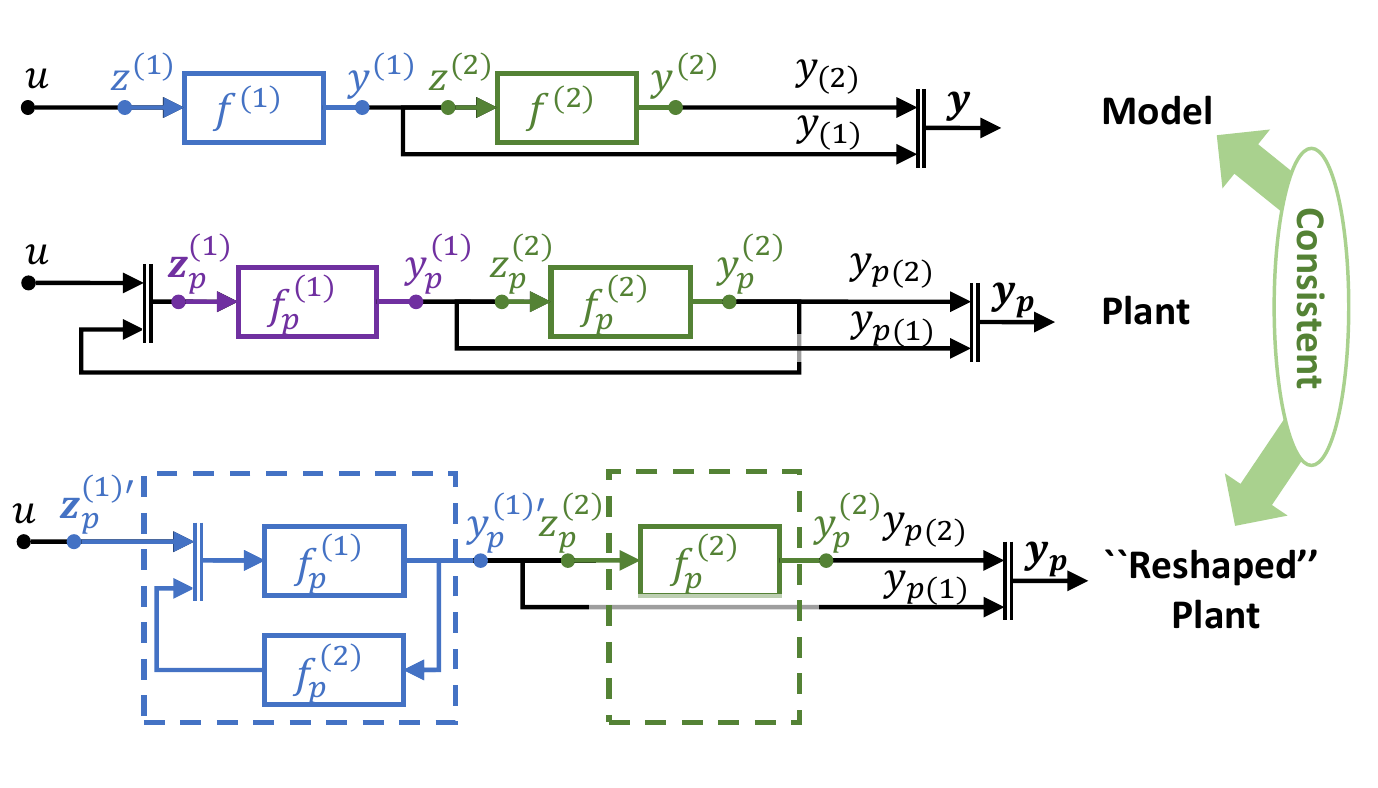}
	\captionof{figure}{Graphical description of the RTO problem of Study 3}
	\label{fig:5_30_ModelPlantCoherent}
\end{minipage}\\

If in the two previous studies the model and the plant had the same structures, here it is clearly not the case anymore. Indeed, SM1 and SP1 do not have the same inputs because the feedback of $y^{(2)}$ is not present. However, despite the presence of this structural error, it is possible to rearrange the SPs  to form groups of SPs whose interconnections are of the same nature and have the same structure as those of the model (see Figure~\ref{fig:5_30_ModelPlantCoherent}).

\noindent
	\begin{minipage}{\linewidth}
	\vspace*{0pt}
	{\centering
		\begin{minipage}[t]{6cm}\centering%
			\includegraphics[width=6cm]{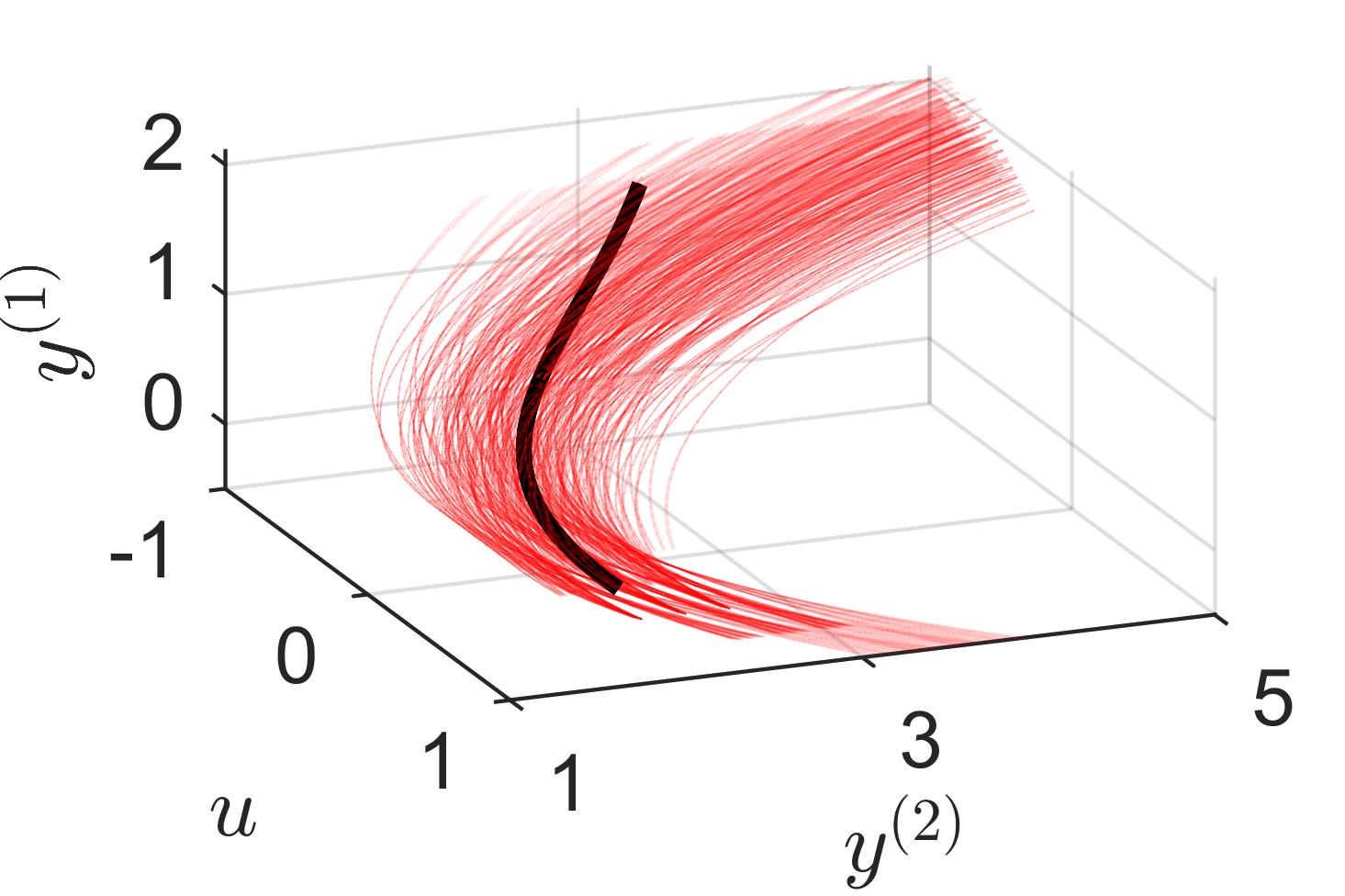}
			a) Functions $f^{(1)}$ and $f_p^{(1)}$
		\end{minipage}\hskip -0ex
		\begin{minipage}[t]{6cm}\centering%
			\includegraphics[trim={0.1cm 0.4cm 0.1cm  0.2cm },clip,width=6cm]{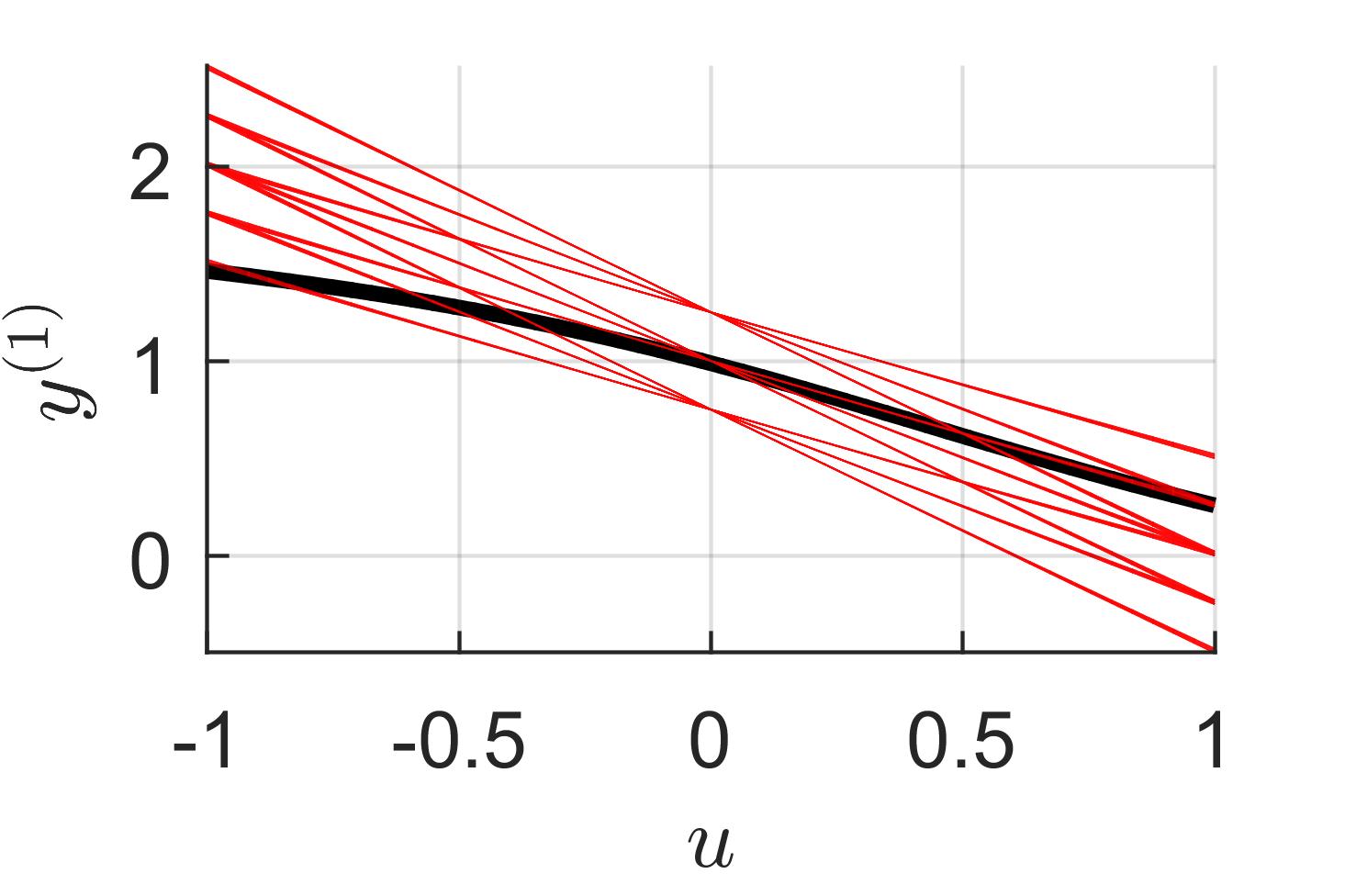}
			b) Functions $f^{(1)}$ and $f_p^{(1)}$
		\end{minipage} \\
		
		\medskip

		\begin{minipage}[t]{4.45cm}\centering%
			\includegraphics[trim={0.1cm 0.4cm 0.1cm  0.2cm },clip,width=4.45cm]{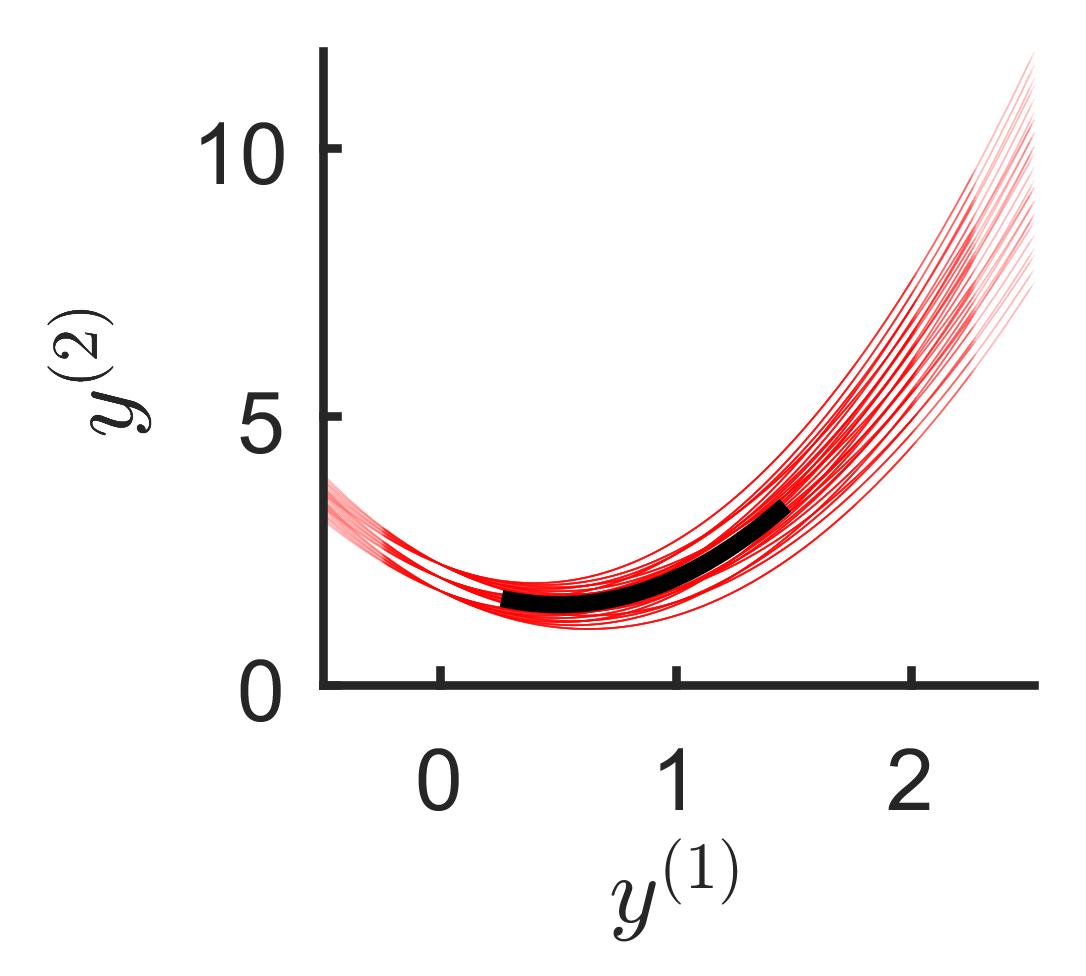}
			c) Functions $f^{(2)}$ and $f_p^{(2)}$
		\end{minipage} \hskip -0ex
		\begin{minipage}[t]{4.45cm}\centering%
			\includegraphics[trim={0.1cm 0.4cm 0.1cm  0.2cm },clip,width=4.45cm]{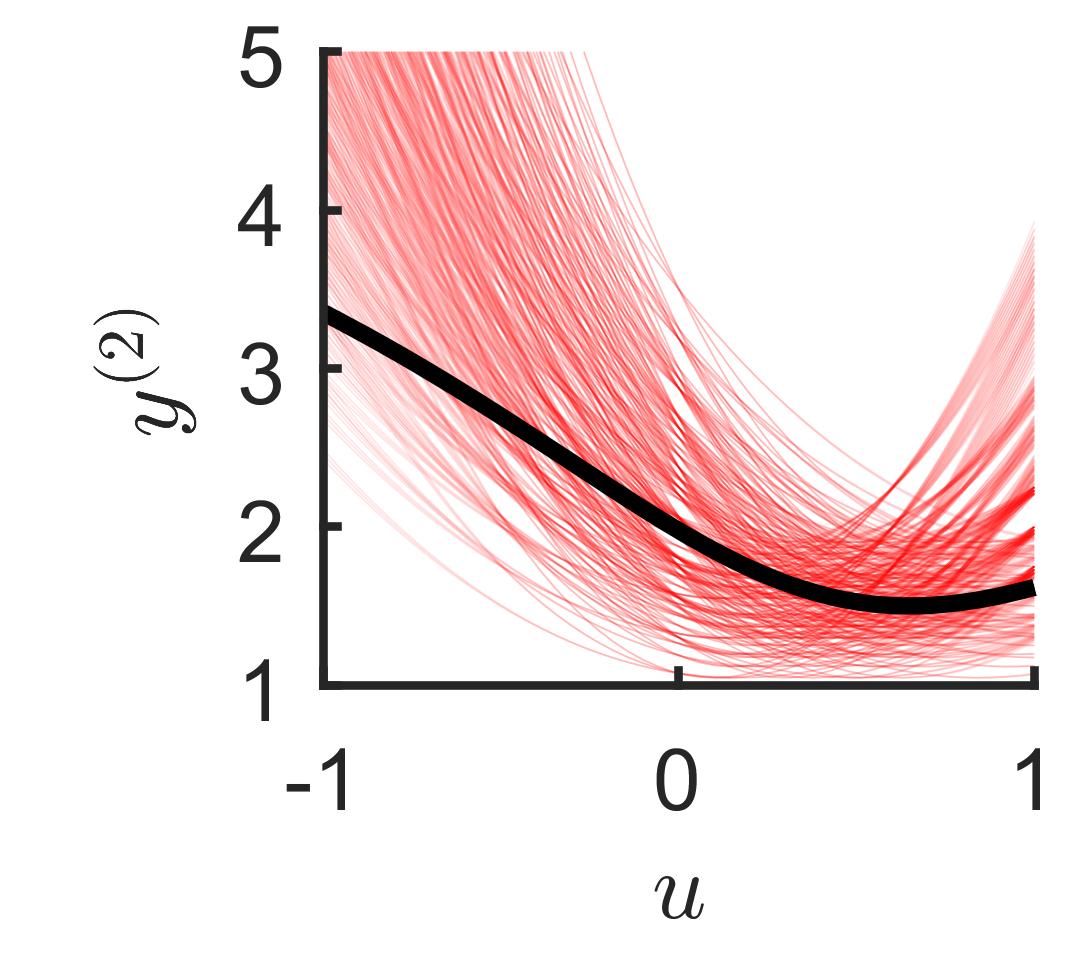}
			d) Functions $f$ and $f_p$
		\end{minipage} \\

		\medskip

		\begin{minipage}[h]{\linewidth}\centering%
			\includegraphics[trim={0.1cm 0.2cm 0.2cm  0.cm },clip,width=3.5cm]{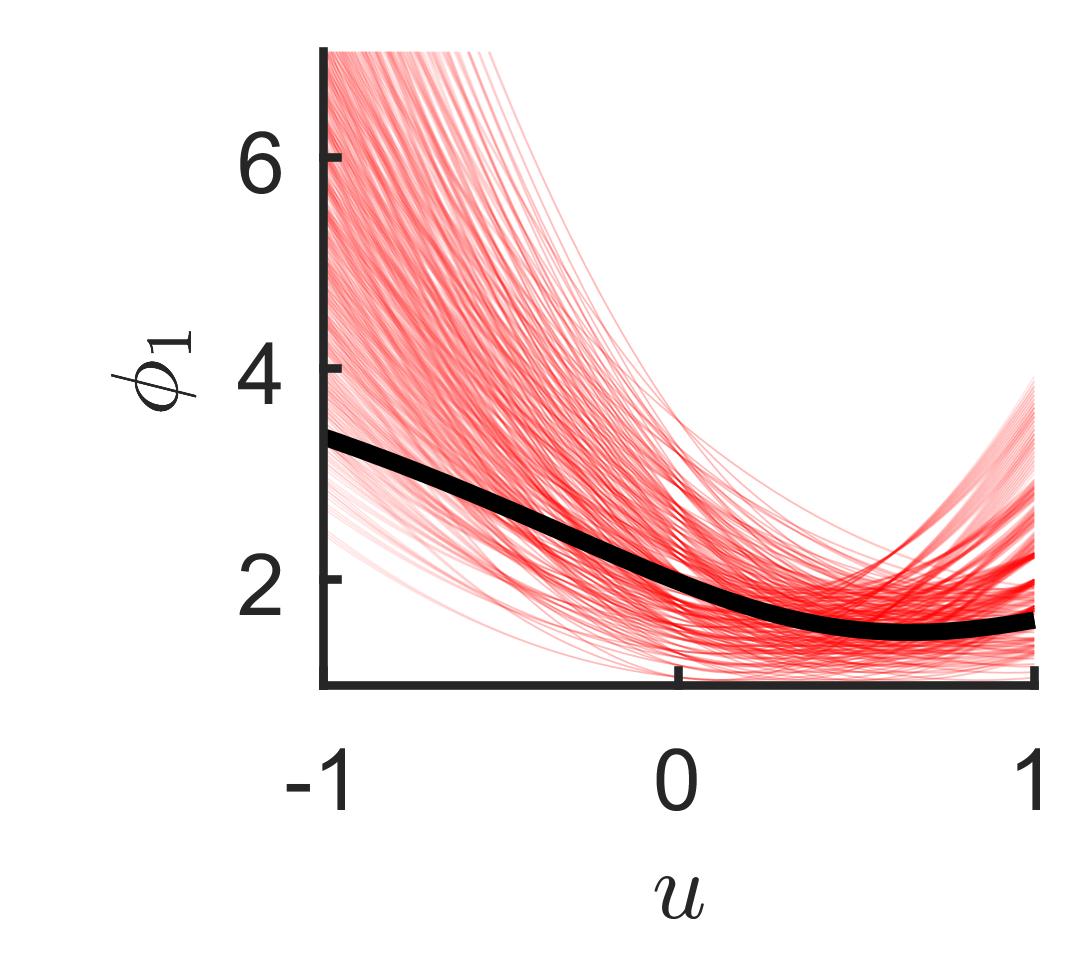}\hskip -0ex
			\includegraphics[trim={0.1cm 0.2cm 0.2cm  0.cm },clip,width=3.5cm]{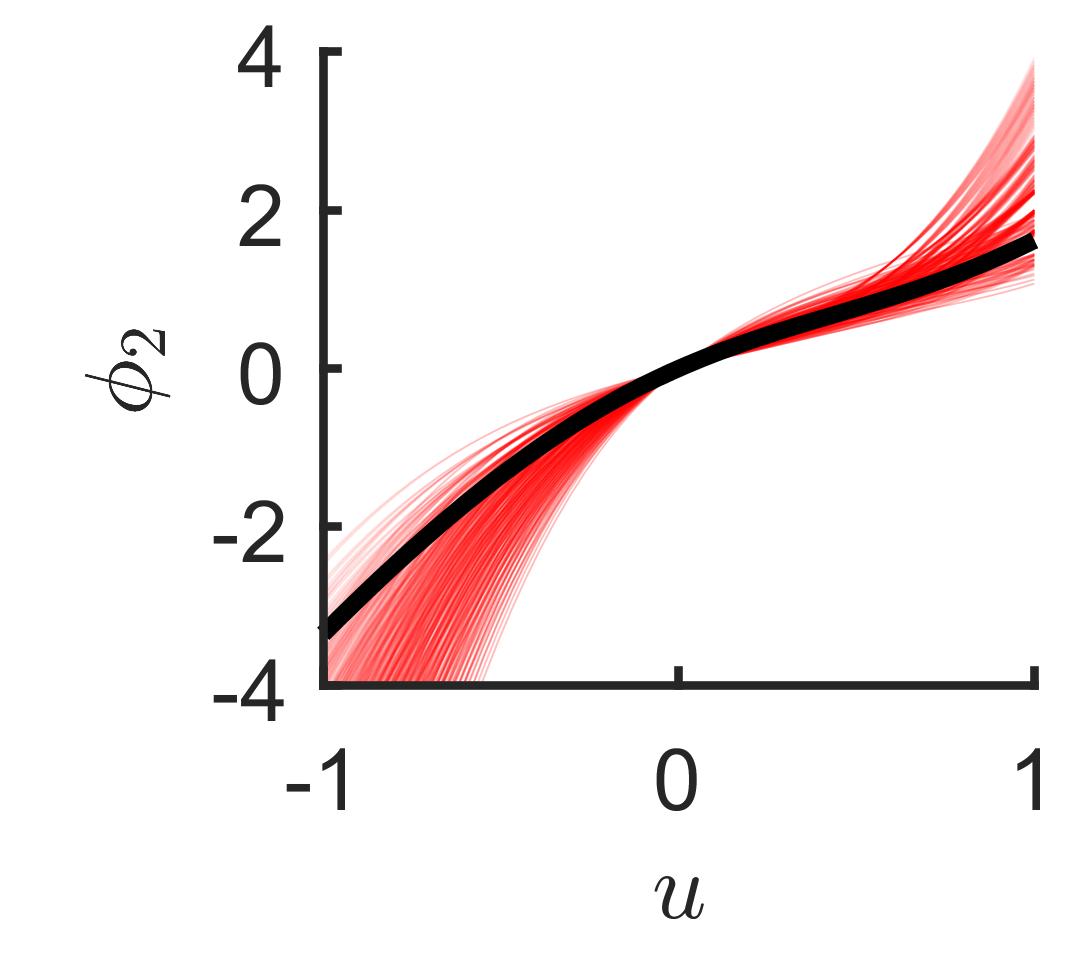}\hskip -0ex
			\includegraphics[trim={0.1cm 0.2cm 0.2cm  0.cm },clip,width=3.5cm]{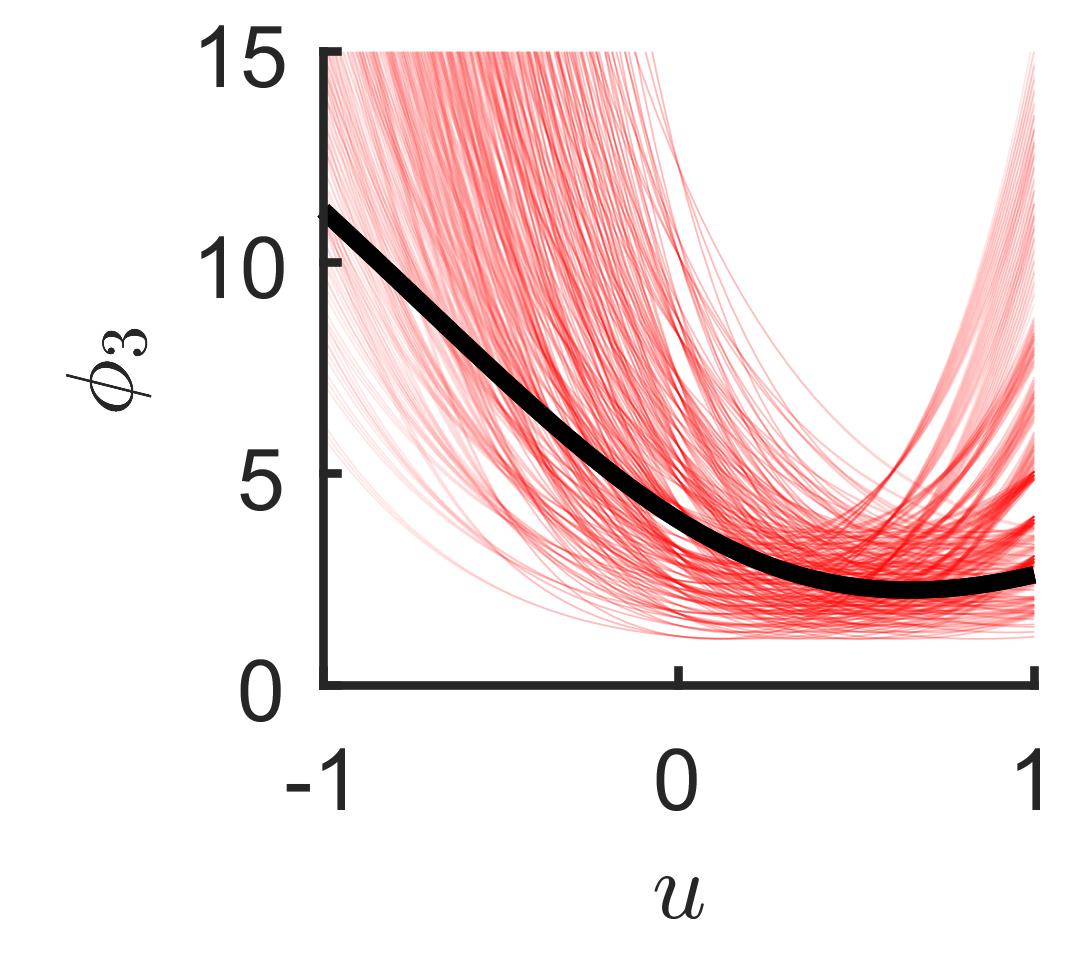}\hskip -0ex
			\includegraphics[trim={0.1cm 0.2cm 0.2cm  0.cm },clip,width=3.5cm]{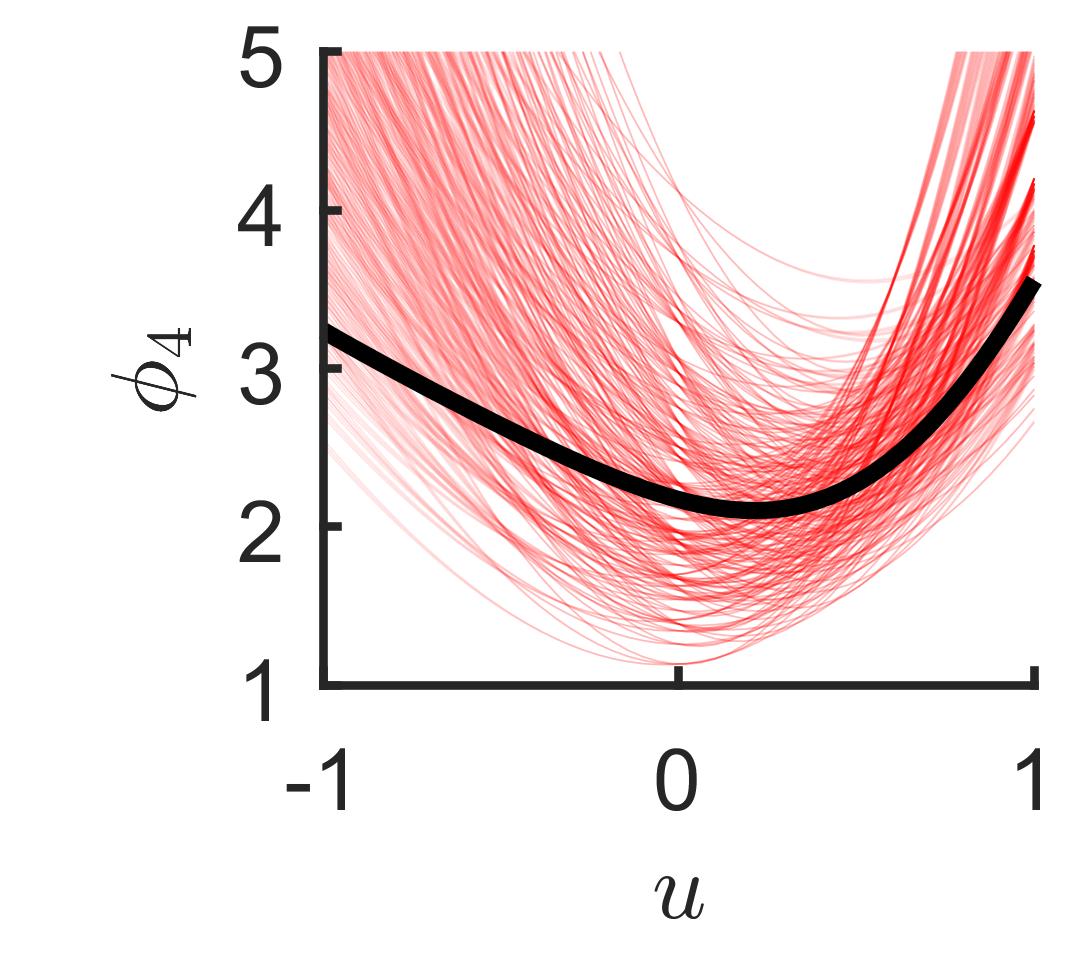}\\
			e) Fonctions $\phi_1$, $\phi_2$, $\phi_3$, et $\phi_4$, 
		\end{minipage}
		
	}	
	\textcolor{red}{\raisebox{0.5mm}{\rule{0.5cm}{0.05cm}}}   : Model, 
	\textcolor{black}{\raisebox{0.5mm}{\rule{0.5cm}{0.1cm}}} : Plant.\\
	\vspace{-3mm}
	\captionof{figure}{Graphical description of the RTO problems}
	\label{fig:5_31_Exemple_5_3_Plant_and_Model}
\end{minipage} \\

\noindent
\begin{minipage}{\linewidth}
	\vspace*{0pt}
	{\centering
		\includegraphics[trim={2.75cm 0.2cm 0.2cm  0.2cm },clip,width=3.45cm]{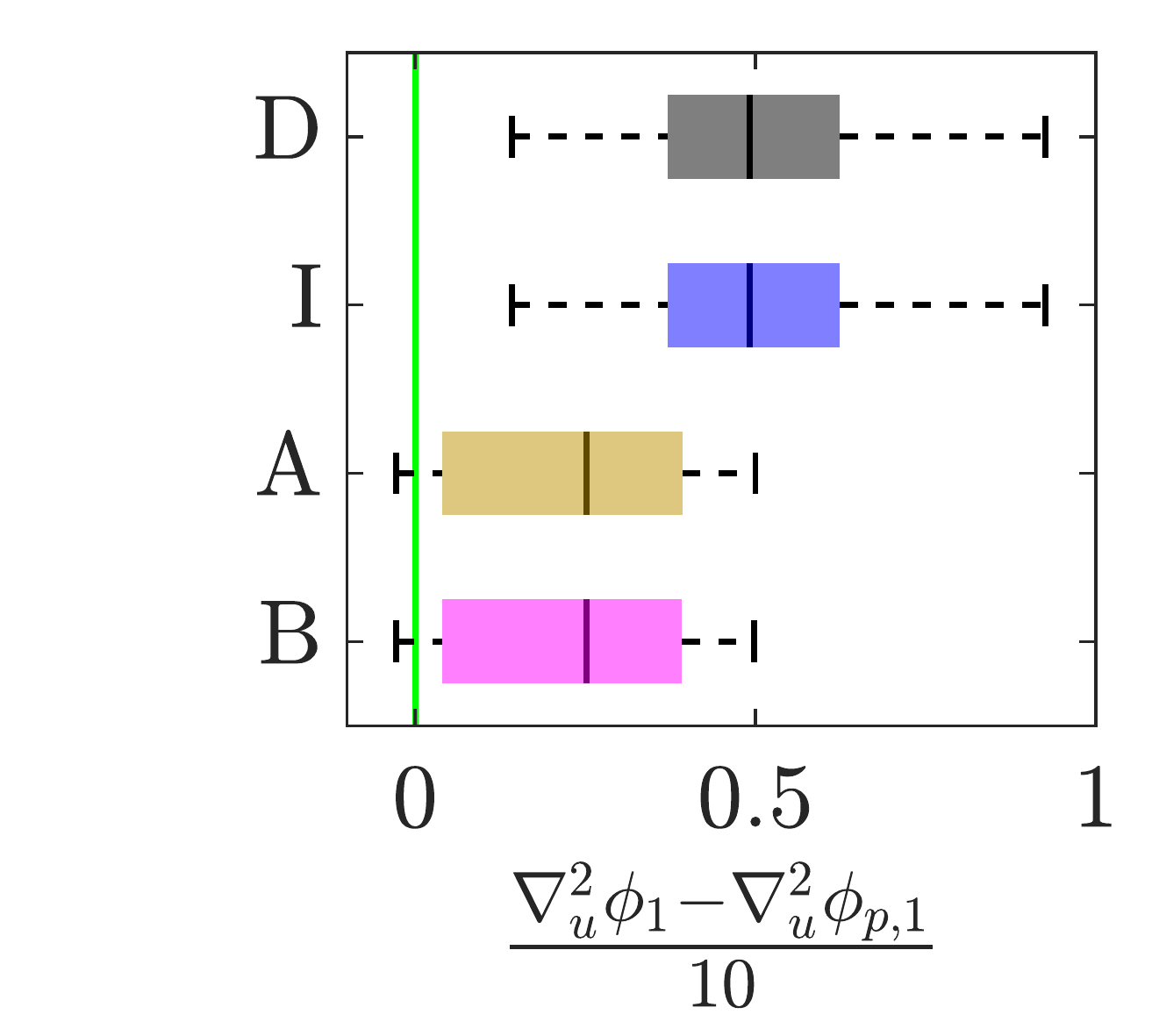}\hskip -0ex
		\includegraphics[trim={2.75cm 0.2cm 0.2cm  0.2cm },clip,width=3.45cm]{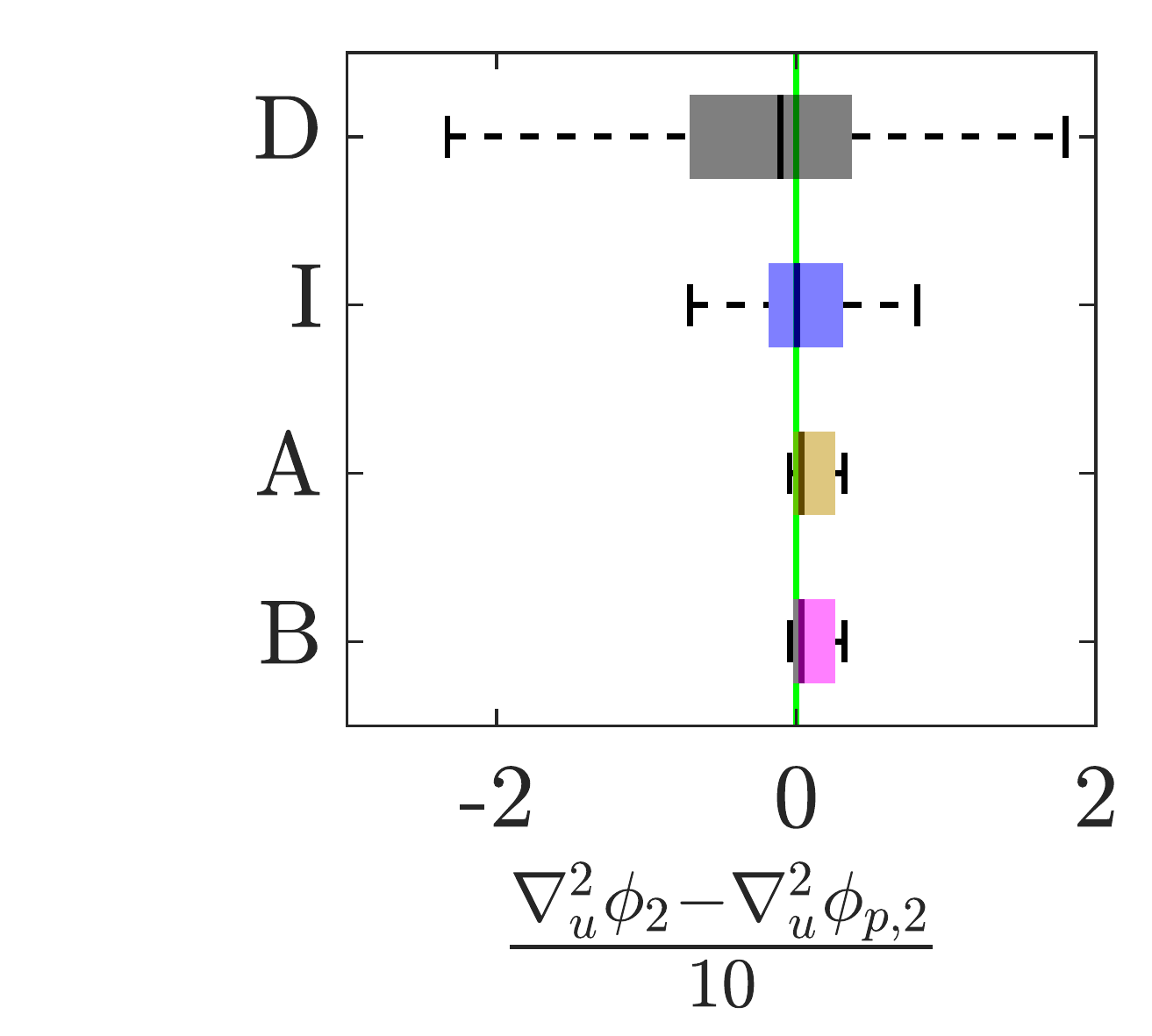}\hskip -0ex
		\includegraphics[trim={2.75cm 0.2cm 0.2cm  0.2cm },clip,width=3.45cm]{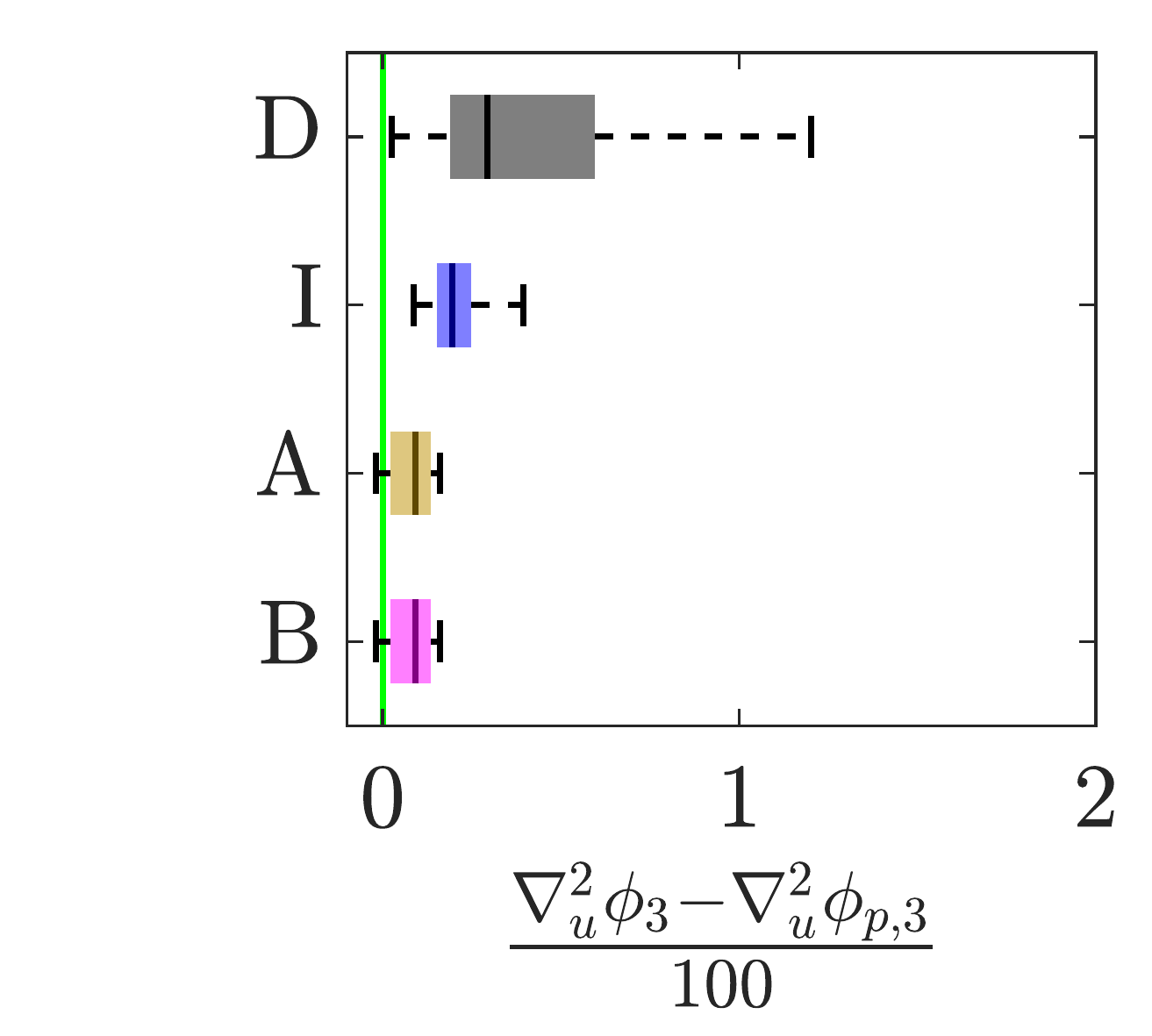}\hskip -0ex
		\includegraphics[trim={2.75cm 0.2cm 0.2cm  0.2cm },clip,width=3.45cm]{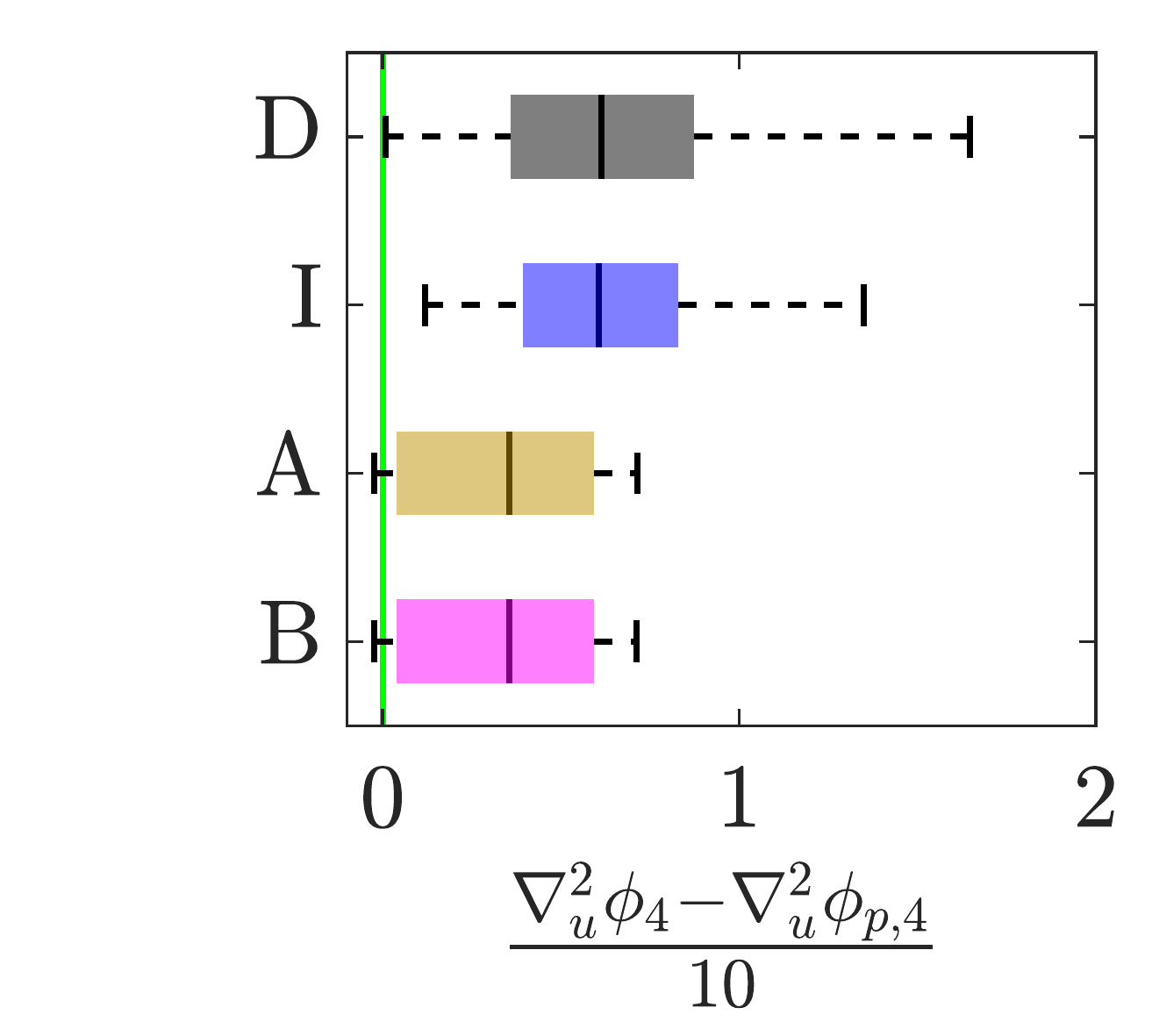}
	} 
	\captionof{figure}{Statistical distributions of the prediction errors on  Hessian of the plant's cost functions at the correction point for the structures D, I, A, and B.}
	\label{fig:5_32_Exemple_5_3_Results}
\end{minipage} 
\noindent
\begin{minipage}{\linewidth}
	\vspace*{0pt}
	{\centering
		\includegraphics[trim={0.7cm 0.2cm 0.2cm  0.2cm },clip,width=3.45cm]{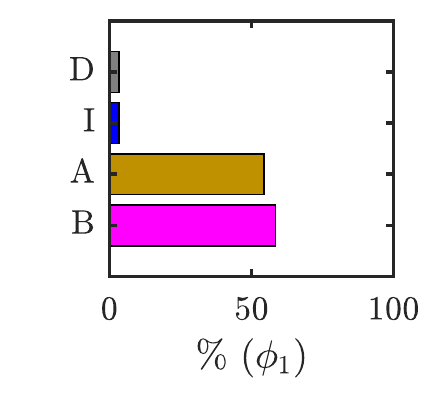}\hskip -0ex
		\includegraphics[trim={0.7cm 0.2cm 0.2cm  0.2cm },clip,width=3.45cm]{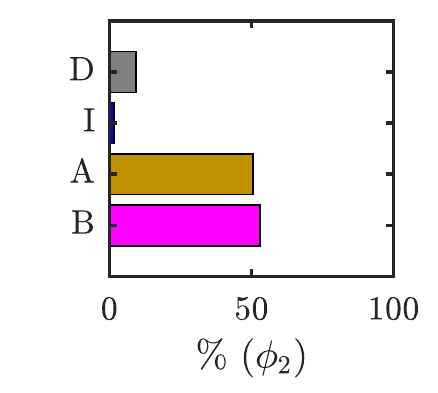}\hskip -0ex
		\includegraphics[trim={0.7cm 0.2cm 0.2cm  0.2cm },clip,width=3.45cm]{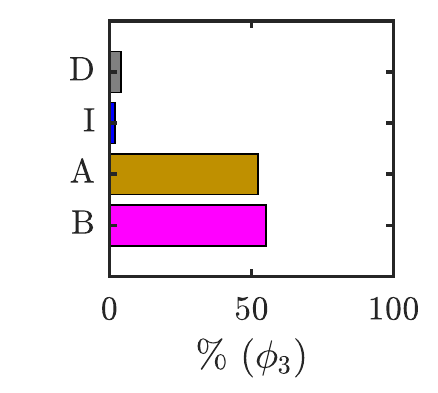}\hskip -0ex
		\includegraphics[trim={0.7cm 0.2cm 0.2cm  0.2cm },clip,width=3.45cm]{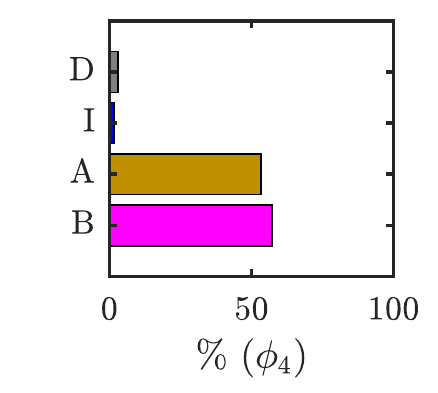} 
	}
	\vspace{-2mm}
	\captionof{figure}{Sc.3: For each correction structure one gives here the percentage of cases for which no other structure provides better results  \textit{(if two structures provide the same best result then both take the point)}}
	\label{fig:5_33_Exemple_5_2_sc3_Results_2}
\end{minipage} \\

As in the previous study, one quantifies for each scenario the ``local quality'' of the corrections with the prediction errors on the Hessian of the plant. It has been evaluated at  11 points uniformly distributed in the input space ($u=\{-1,-0.8,...,0.8,1\}$), for each model, and for each correction structure. Figure~\ref{fig:5_32_Exemple_5_3_Results} provides a statistical summary of those results, and Figure~\ref{fig:5_33_Exemple_5_2_sc3_Results_2} gives for each correction structure the percentage of cases for which no other correction structure brings better predictions of the Hessian.  One can observe that in the presence of significant structural error, structures A and B achieve a better correction of the model than structures D and I. 

(Contextualization) Although purely mathematical, the case that has been chosen for this study can correspond to (pieces of) a practical problems. For example, one can imagine that the first set of functions ($f^{(2)}$) is the set of equations describing the dynamics of a reaction at a given temperature and that the second set of equations ($f^{(1)}$) describes its thermodynamics. In this case, the variables that would interconnect these two systems would be the temperature and the concentrations of the mixture.  One could then imagine that a (simplified) model would have a fixed temperature $T_{set}$ breaking the feedback between these two systems, while the plant would ``contain'' it.  Therefore the RTO problem would have a structure similar to the mathematical problem discussed here. 

\subsection{Study 4: Inconsistent model}
\label{sec:5_5_4_Model_structurellement_Faux_et_incoherent}

In this section, an example is given of a case where a model is inconsistent and one shows that in this case the structure B is unable to provide affine corrections to the model's predictions. 

One considers the theoretical RTO problem illustrated on the Figures~\ref{fig:5_34_ModelPlantCoherent} and \ref{fig:5_35_Exemple_5_4_Plant_and_Model} and using the following functions:
\begingroup
 \allowdisplaybreaks
\begin{align}
	f^{(1)}   := \ & \theta_1 - \theta_2z^{(1)}  + \theta_3(z^{(1)})^2, &  z^{(1)} := \ & u_{(2)}, \nonumber  \\
	f_p^{(1)} := \ & 1 - z_p^{(1)} + 0.01 (z_p^{(1)})^2,                & z_p^{(1)} := \ & u_{(1)}, \nonumber \\
	f^{(2)}   := \ & \theta_4 - \theta_5 z_{(1)}^{(2)} + \theta_6\big( (z_{(1)}^{(2)})^2 + (z_{(2)}^{(2)})^2 \big), &  z^{(2)}   := \ & \left(\begin{array}{c} u_{(1)} \\ y^{(1)} \end{array}\right),   \nonumber  \\
	f_p^{(2)} := \ & 2 - 2z_{p(1)}^{(2)} + 2\big( (z_{p(1)}^{(2)})^2 + (z_{p(2)}^{(2)})^2 \big), &  \bm{z}_p^{(2)} := \ & \left(\begin{array}{c} u_{(2)} \\ y_p^{(1)} \end{array}\right),   \nonumber 
\end{align}
\begin{align*}
	\phi_1(u,\bm{y}) := \ & y_{(2)},   & \phi_2(u,\bm{y}) := \ & uy_{(2)}, \\
	\phi_3(u,\bm{y}) := \ & y_{(2)}^2, & \phi_4(u,\bm{y}) := \ & y_{(2)} + u^2 + 0.5uy_{(2)} + 0.05y_{(2)}^2.
\end{align*}
 \endgroup
where $u_{(1)}\in[-1,1]$, $u_{(2)}\in[-1,1]$, $\theta_1\in[0.75,1.25]$, $\theta_2\in[-1.25,-0.75]$, $\theta_3\in[0,0.02]$, $\theta_4\in[1.75,2.25]$, $\theta_5\in[-2.25,-1.75]$, and $\theta_6\in[1.8,2.2]$. 

The structural error of the model is clearly shown on Figure~\ref{fig:5_34_ModelPlantCoherent} where the inputs $u_{(1)}$ and $u_{(2)}$ are reversed. The effect of this inversion is that it is impossible to rearrange the SPs to form a structure similar to the one of the model. For instance, the model predicts a correlation between $u_{(2)}$ and $y_{(1)}$  and clearly such a correlation cannot exist in the plant. Hence the inconsistency of the model.   

To quantify the ``local quality'' of the correction associated to each structure, the distribution of the prediction errors of the \textit{gradient} and the Hessian of the plant for each model, correction structure, and cost function is evaluated at 6 points uniformly distributed in the input space ($u=\{-1,-0.6,-0.2,0.2,0.6,1\}$). Then, applying a statistical analysis to the obtained results, the graphs of Figures~\ref{fig:5_36_Exemple_5_4_Results} and \ref{fig:5_37_Exemple_5_4_Results} are obtained, 

Two observations can be made. First, it is clearly observable that the structure B is unable to correct the predictions of the plant's gradients, while all other structures do. Second, the Hessian predictions obtained with the structures A and B are generally worse  than those obtained with structure D which in turn are generally worse  than those obtained with structure I.

 This result is quite logical since it is obvious that to exploit the ``structural information'' of model that is inconsistent  is equivalent to introduce wrong information in the decision making process.  Therefore, if the structure A performs better than the structure B, it is clearly because it exploits ``less'' of this wrong information.

(Contextualization) Finally, although the inconsistency of a model is a critical weakness of the structures A and B,  one believes that in practice a model is unlikely to be that wrong.  

\clearpage

\begin{minipage}{\linewidth}
	\includegraphics[width=14cm]{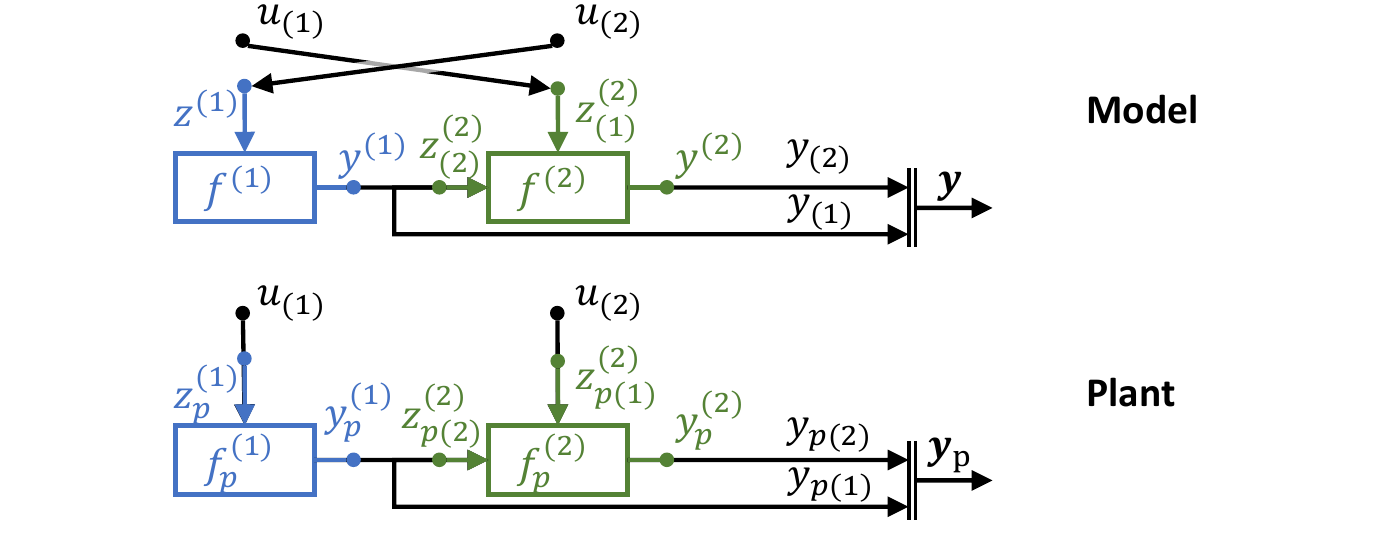}
	\captionof{figure}{Description of the RTO problem of Study 4}
	\label{fig:5_34_ModelPlantCoherent}
\end{minipage}\\

\begin{minipage}{\linewidth}
	\vspace*{0pt}
	{\centering
		\begin{minipage}[t]{6cm}\centering%
			\includegraphics[width=6cm]{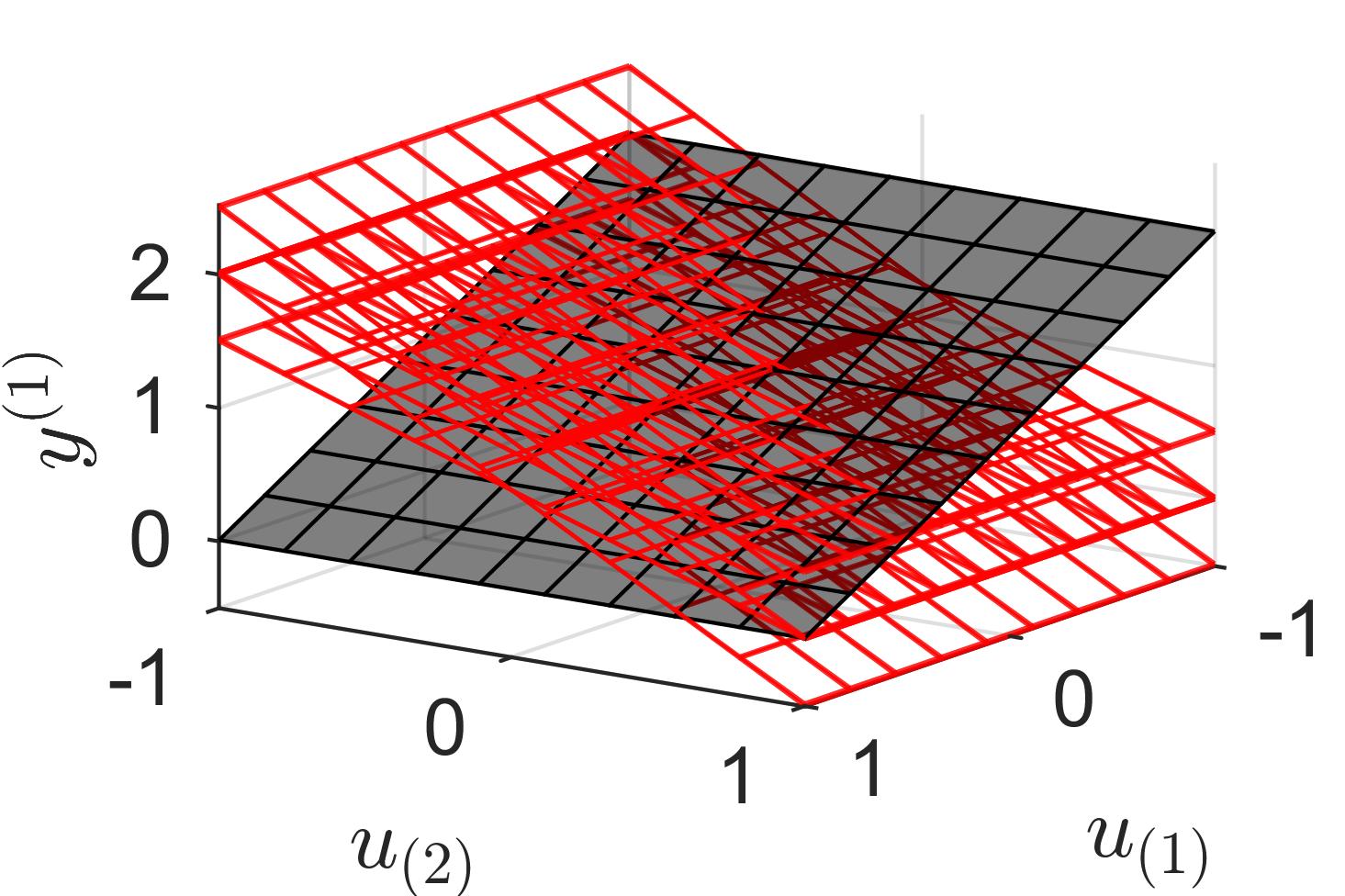}
			a) Functions $f^{(1)}$ and $f_p^{(1)}$
		\end{minipage}\hskip -0ex
		\begin{minipage}[t]{6cm}\centering%
			\includegraphics[trim={0.1cm 0.4cm 0.1cm  0.2cm },clip,width=6cm]{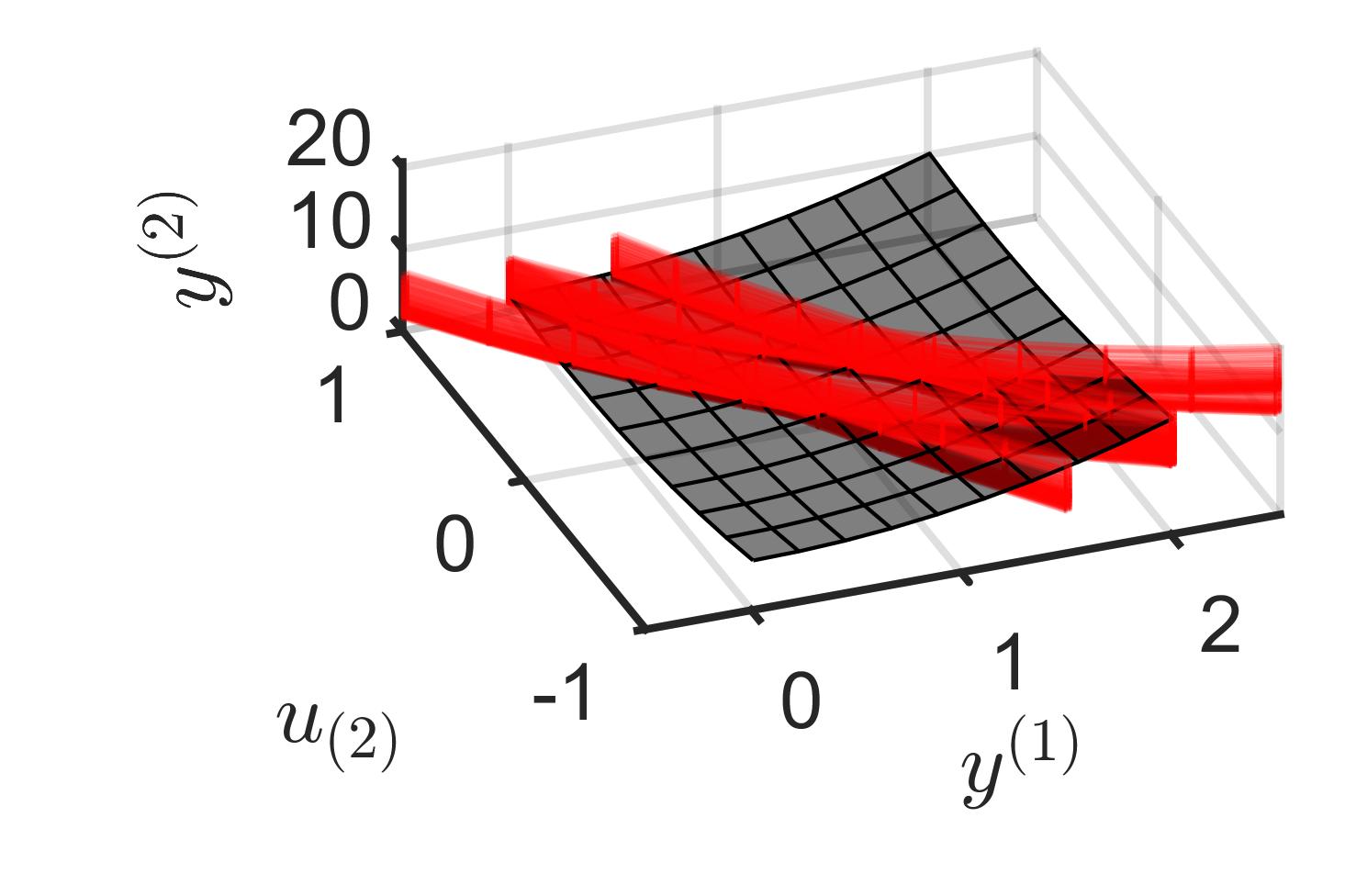}
			b) Functions $f^{(2)}$ and $f_p^{(2)}$
		\end{minipage} \\
		
		\medskip

		\begin{minipage}[t]{6cm}\centering%
			\includegraphics[trim={0.1cm 0.4cm 0.1cm  0.2cm },clip,width=6cm]{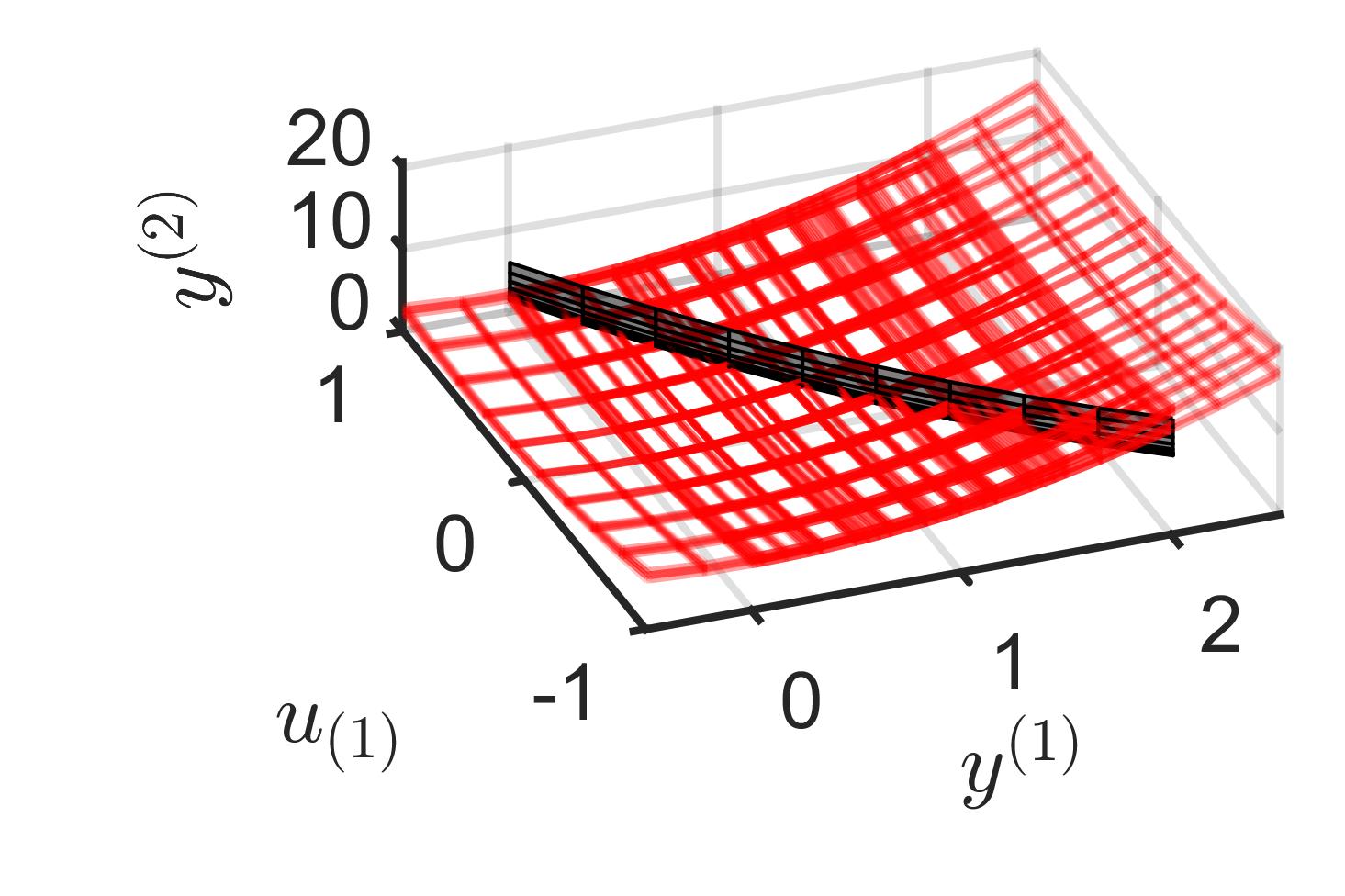}
			c) Functions $f^{(2)}$ and $f_p^{(2)}$
		\end{minipage} \hskip -0ex
		\begin{minipage}[t]{6cm}\centering%
			\includegraphics[trim={0.1cm 0.4cm 0.1cm  0.2cm },clip,width=6cm]{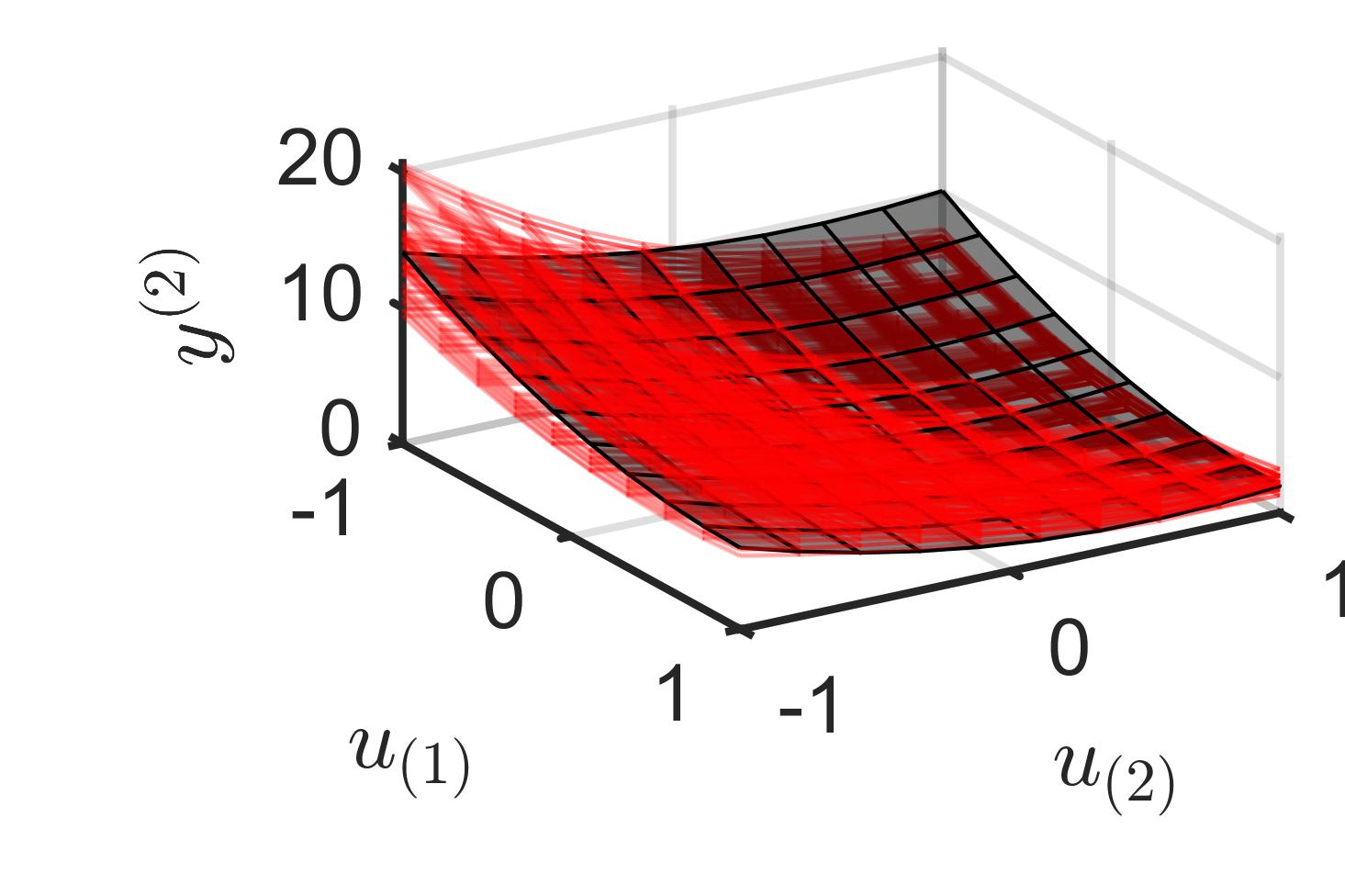}
			d) Functions $f$ and $f_p$
		\end{minipage} \\

		\medskip

		\begin{minipage}[h]{\linewidth}\centering%
			\includegraphics[trim={0.1cm 0.2cm 0.2cm  0.cm },clip,width=3.5cm]{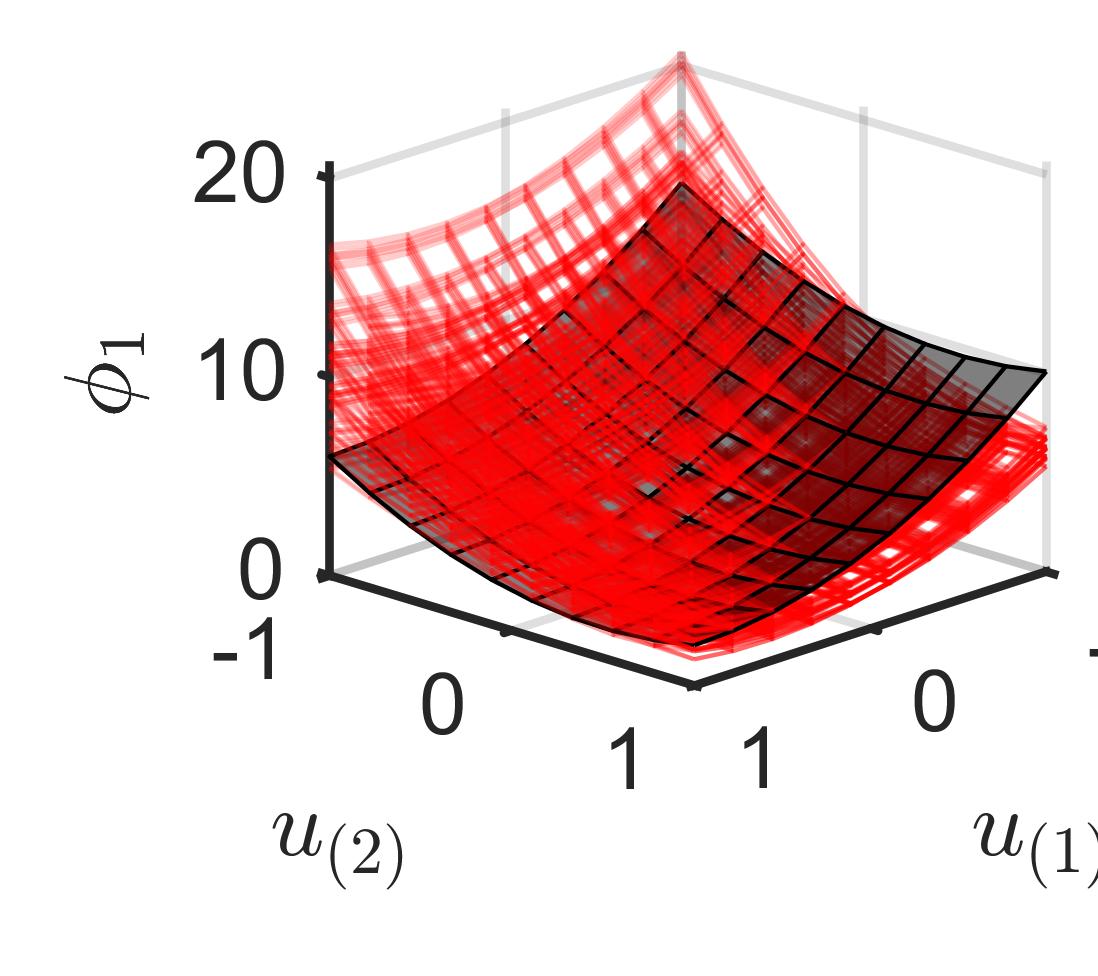}\hskip -0ex
			\includegraphics[trim={0.1cm 0.2cm 0.2cm  0.cm },clip,width=3.5cm]{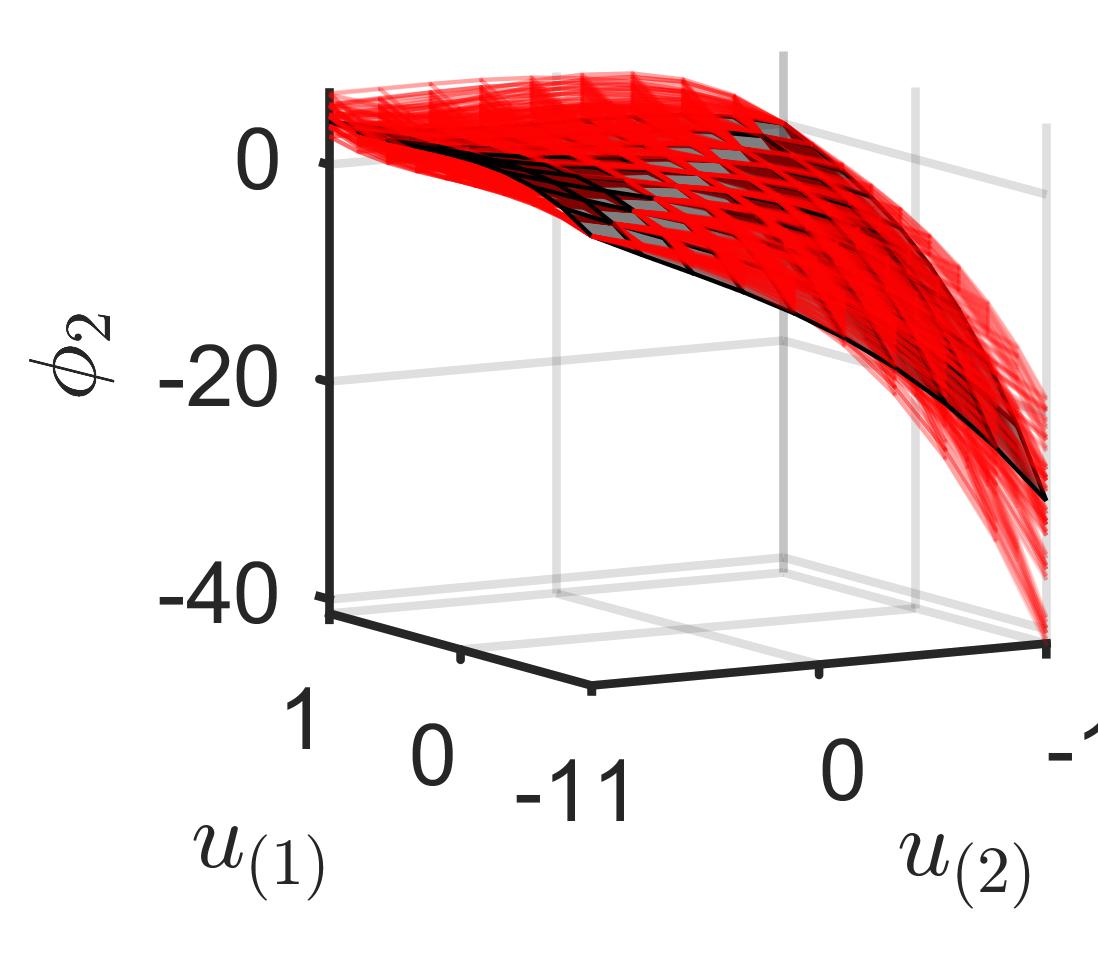}\hskip -0ex
			\includegraphics[trim={0.1cm 0.2cm 0.2cm  0.cm },clip,width=3.5cm]{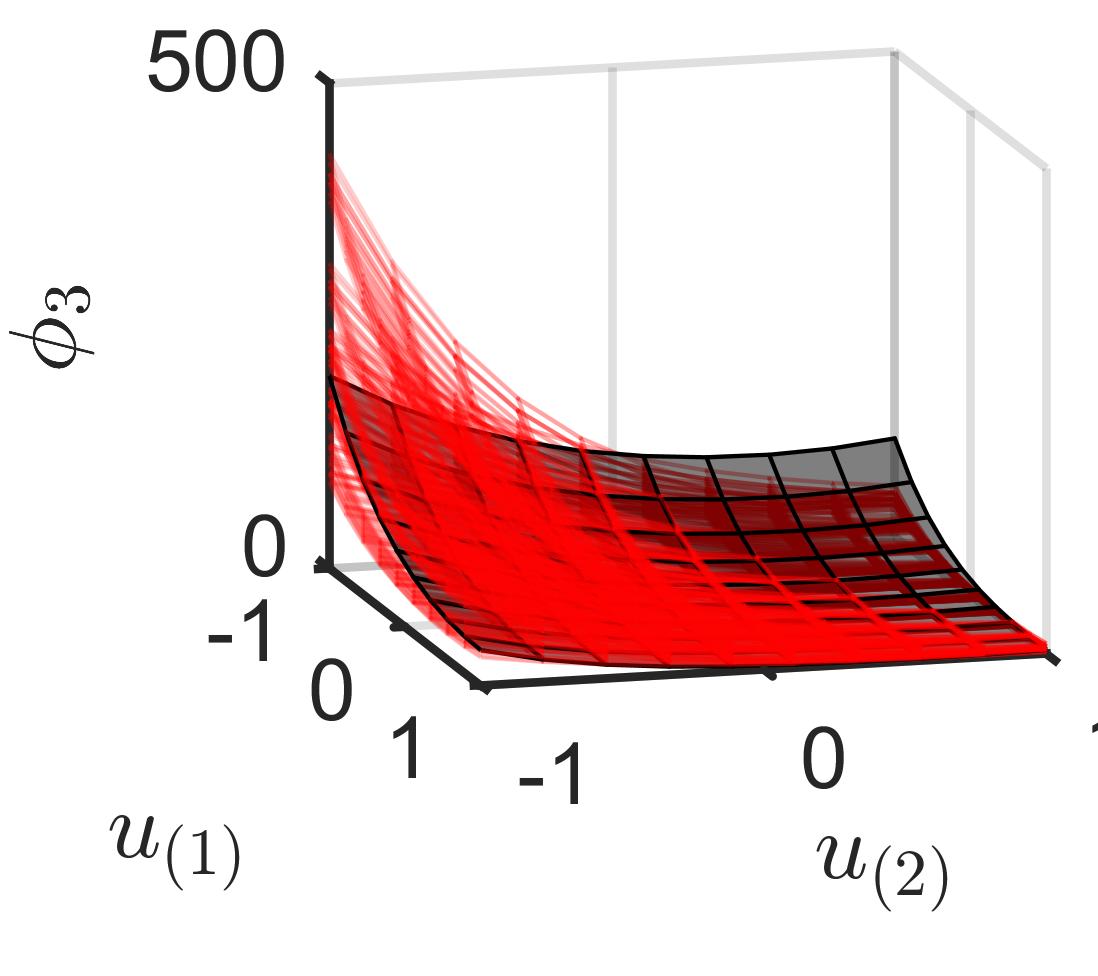}\hskip -0ex
			\includegraphics[trim={0.1cm 0.2cm 0.2cm  0.cm },clip,width=3.5cm]{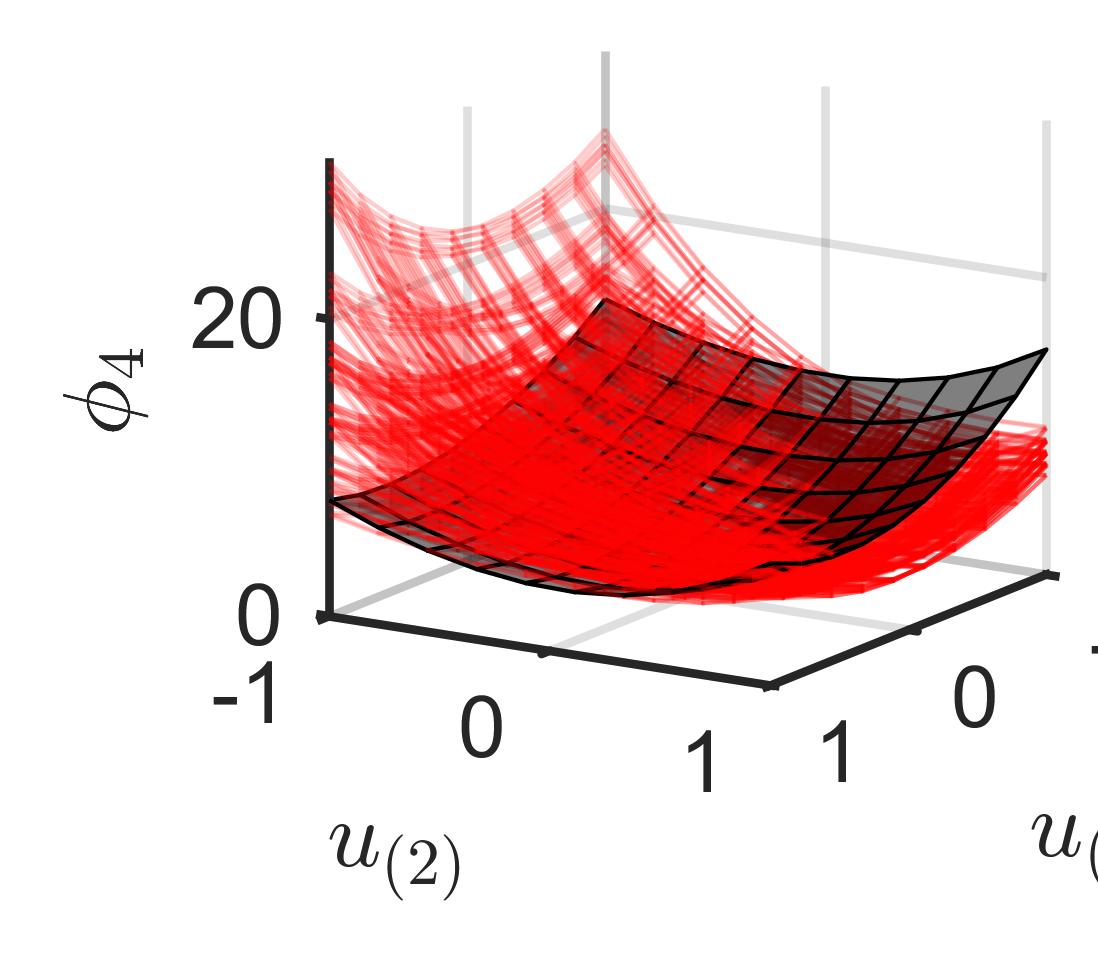}\\
			e) Functions $\phi_1$, $\phi_2$, $\phi_3$, and $\phi_4$, 
		\end{minipage}
		
	}	
	\textcolor{red}{\raisebox{0.5mm}{\rule{0.5cm}{0.05cm}}}   : Model, 
	\textcolor{black}{\raisebox{0.5mm}{\rule{0.5cm}{0.1cm}}} : Plant.\\
	\vspace{-3mm}
	\captionof{figure}{Graphical description of the RTO problems}
	\label{fig:5_35_Exemple_5_4_Plant_and_Model}
\end{minipage}

		\begin{minipage}[h]{\linewidth}
			\vspace*{0pt}
			{\centering
				\includegraphics[trim={2.75cm 0.2cm 0.2cm  0.2cm },clip,width=3.45cm]{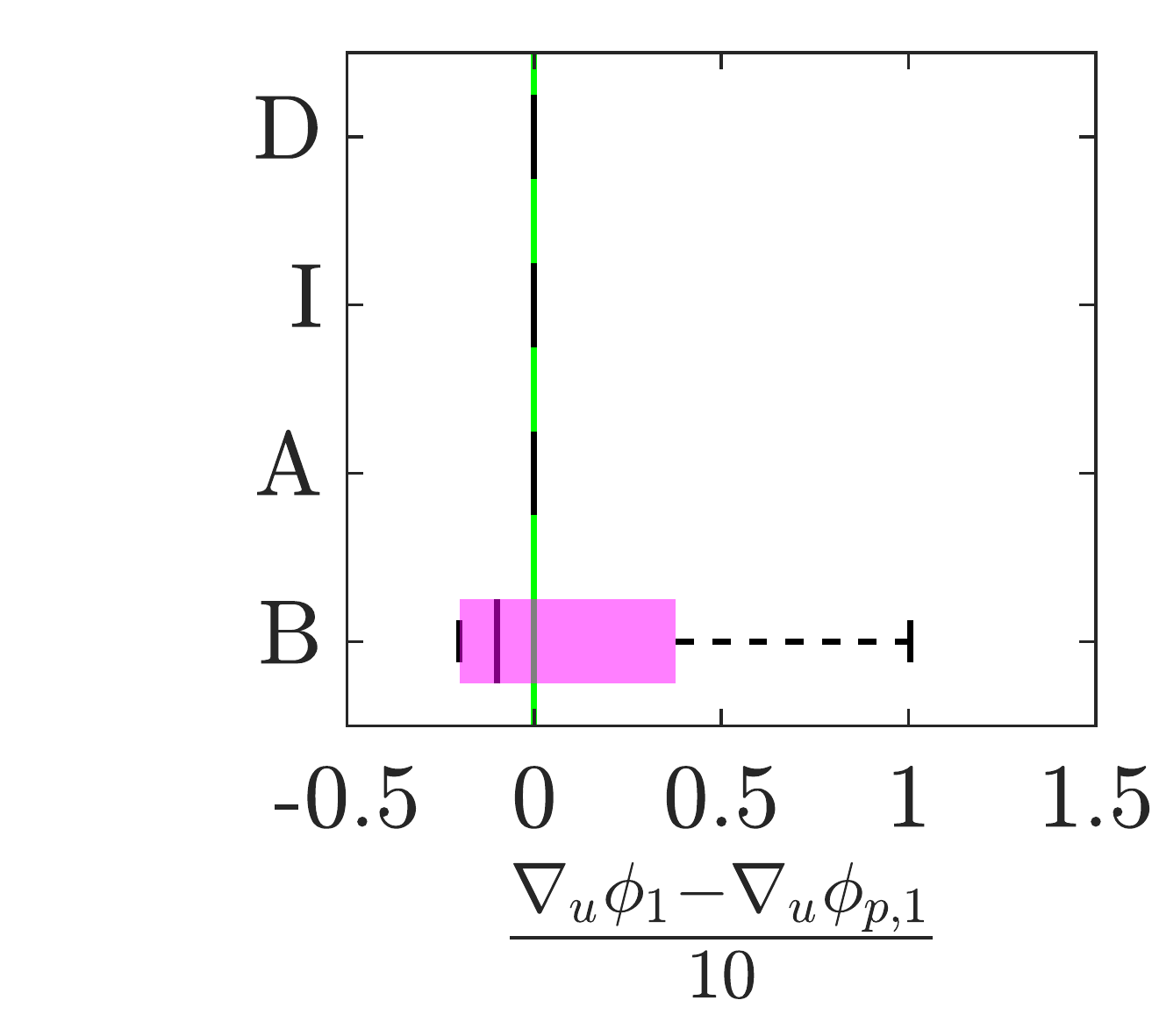}\hskip -0ex
				\includegraphics[trim={2.75cm 0.2cm 0.2cm  0.2cm },clip,width=3.45cm]{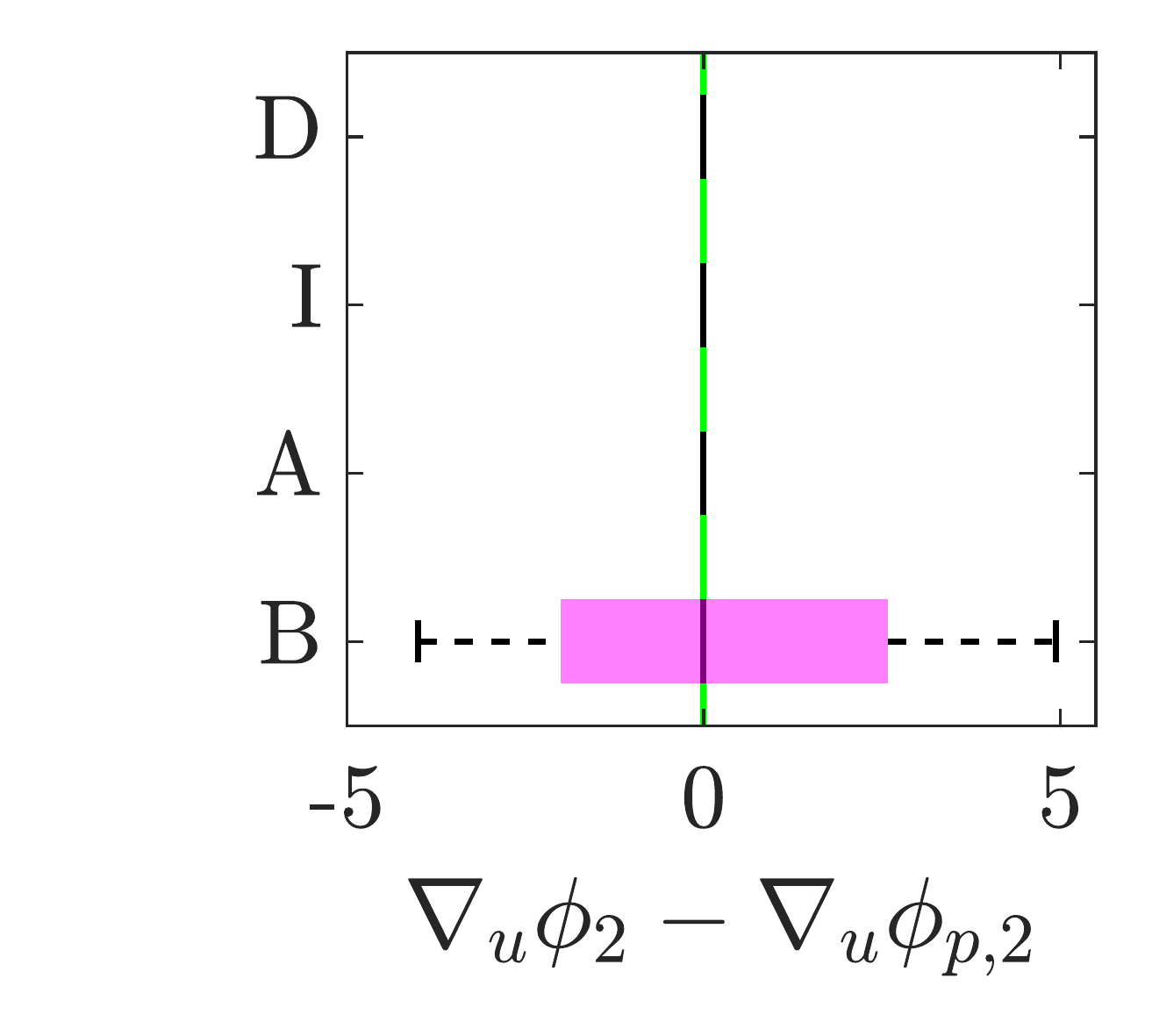}\hskip -0ex
				\includegraphics[trim={2.75cm 0.2cm 0.2cm  0.2cm },clip,width=3.45cm]{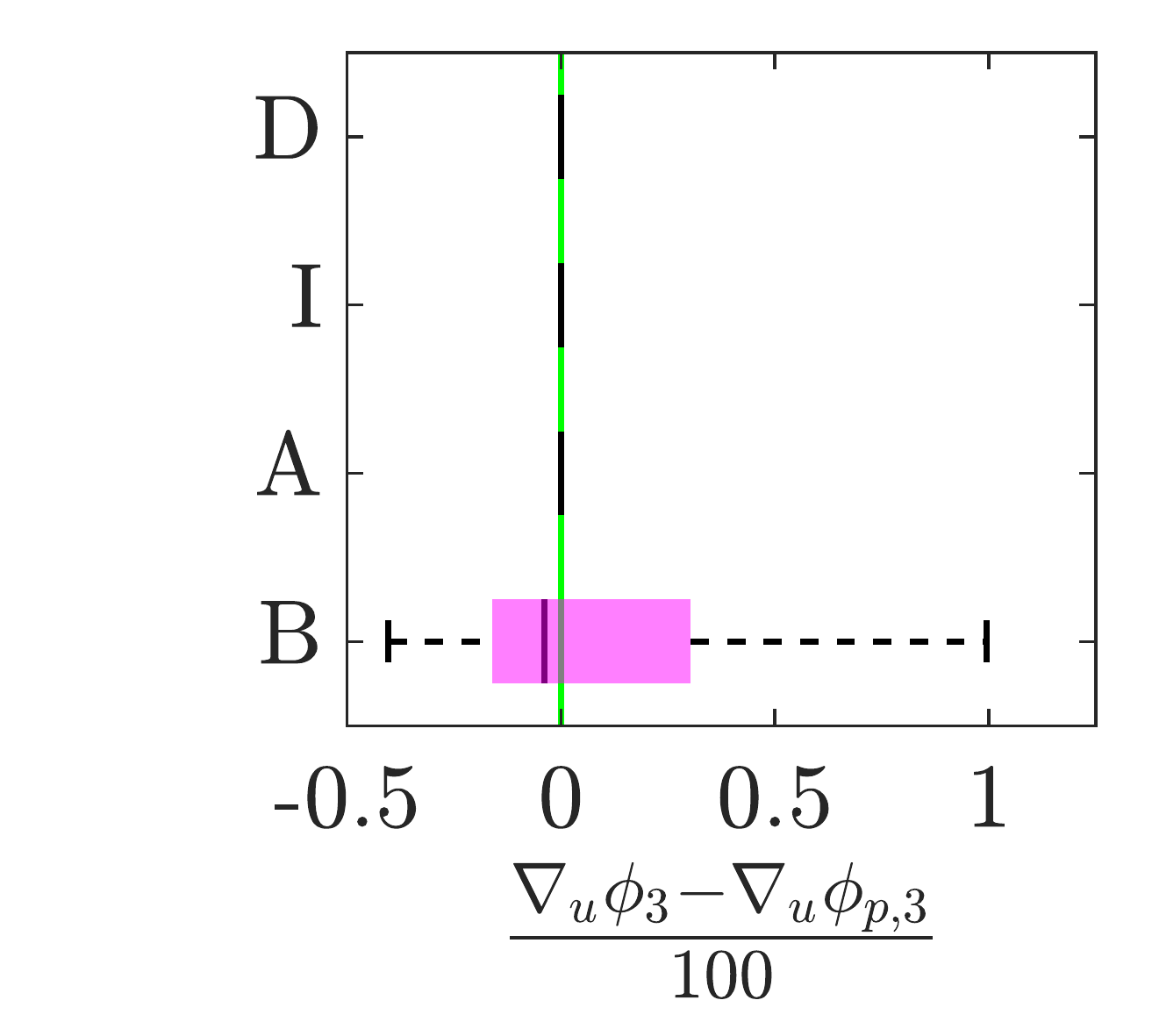}\hskip -0ex
				\includegraphics[trim={2.75cm 0.2cm 0.2cm  0.2cm },clip,width=3.45cm]{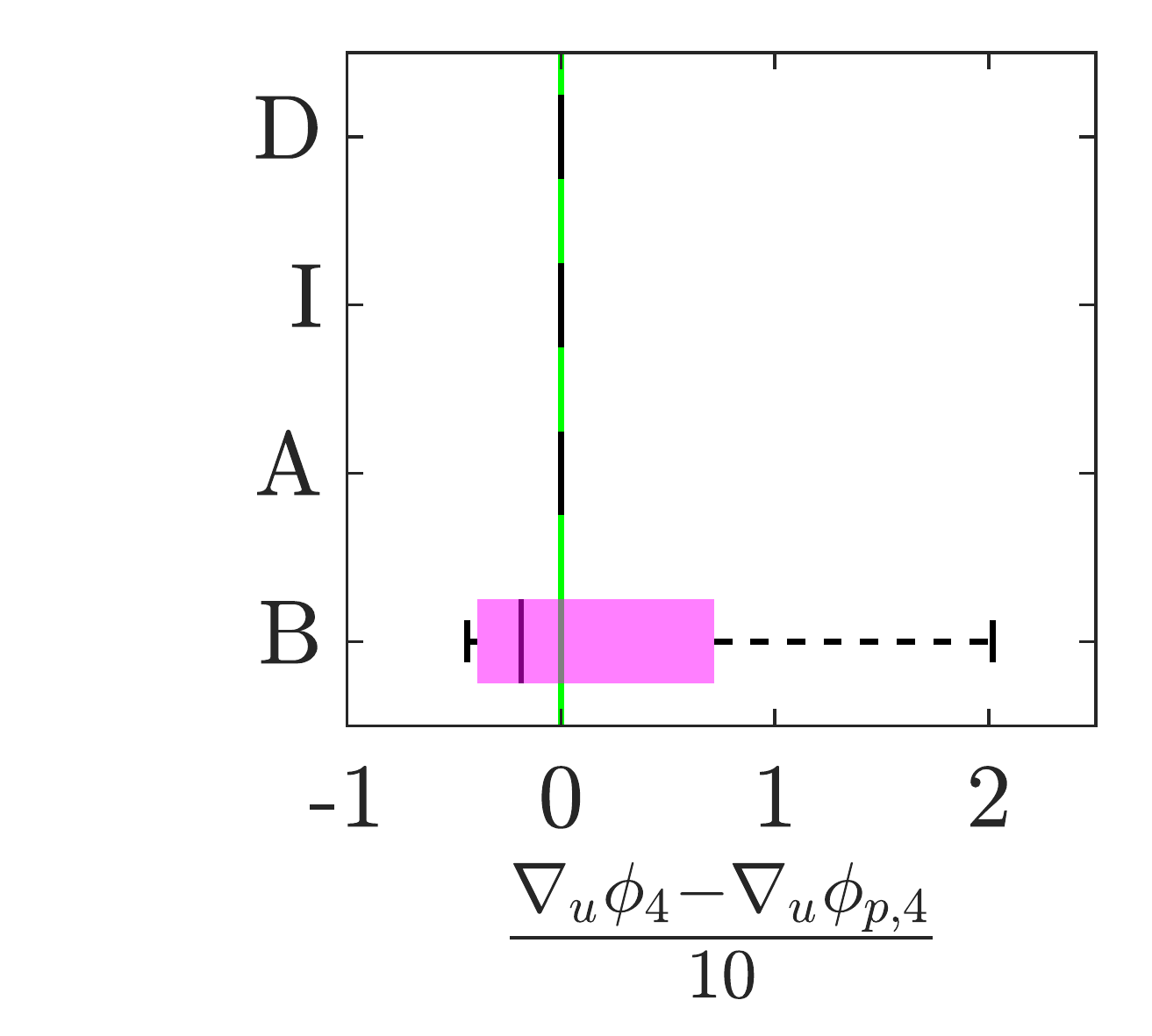}
			} 
			\captionof{figure}{
				Statistical distributions of the prediction errors on the gradient of the plant's cost at the correction point for the structures D,I,A, and B.}
			\label{fig:5_36_Exemple_5_4_Results}
		\end{minipage} 	
		
			\begin{minipage}[h]{\linewidth}
				\vspace*{0pt}
				{\centering	
				\includegraphics[trim={2.75cm 0.2cm 0.2cm  0.2cm },clip,width=3.45cm]{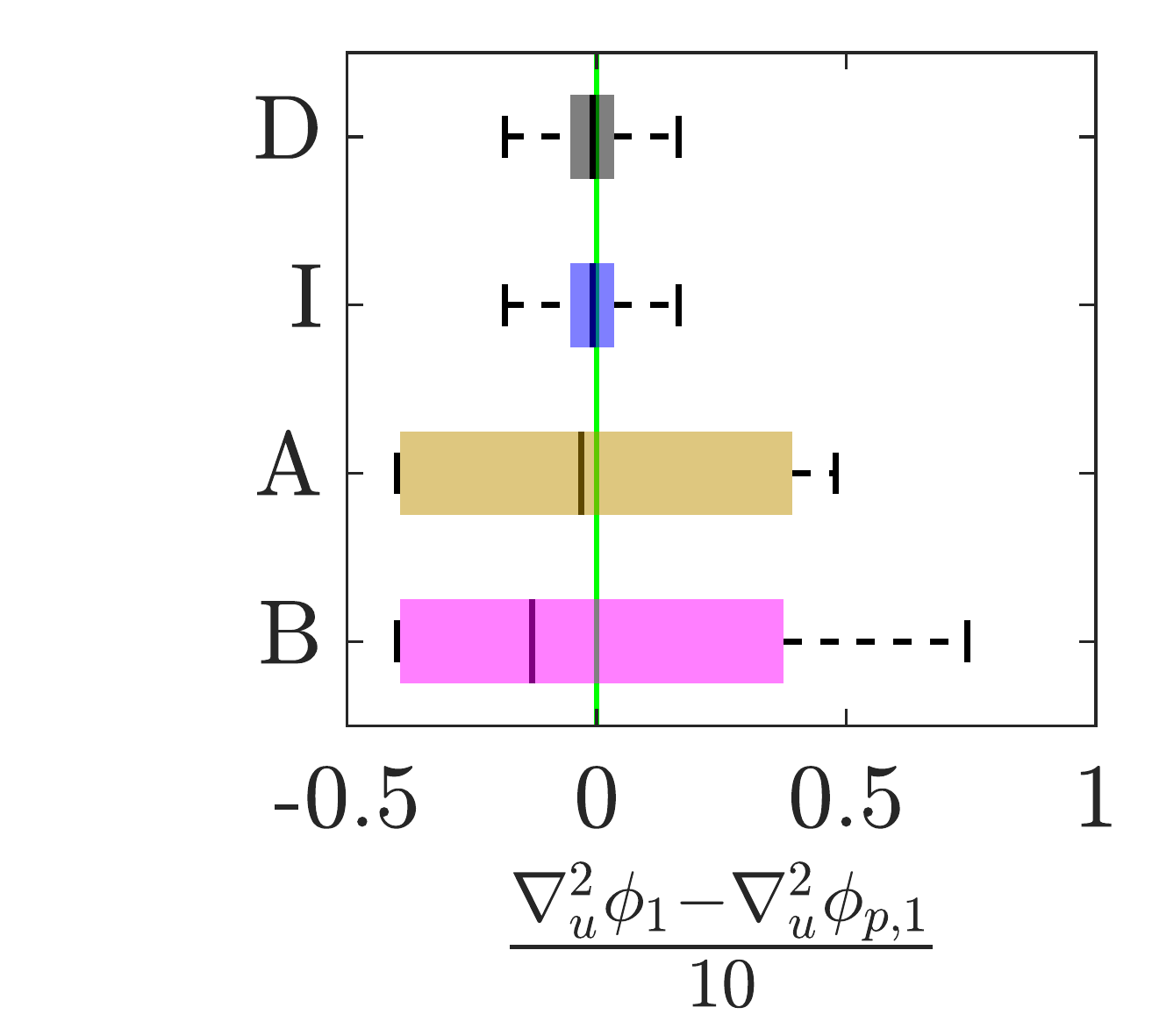}\hskip -0ex
				\includegraphics[trim={2.75cm 0.2cm 0.2cm  0.2cm },clip,width=3.45cm]{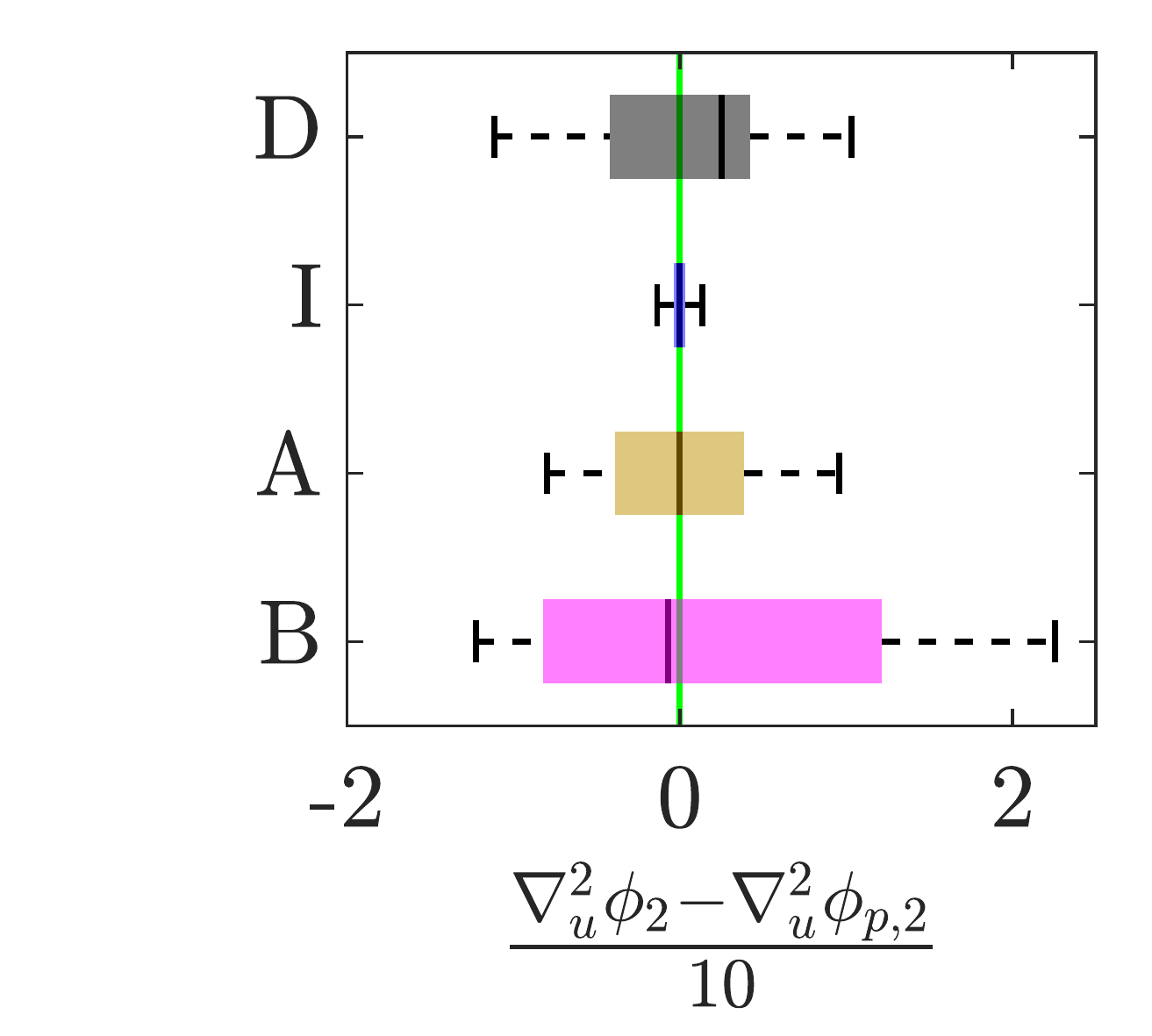}\hskip -0ex
				\includegraphics[trim={2.75cm 0.2cm 0.2cm  0.2cm },clip,width=3.45cm]{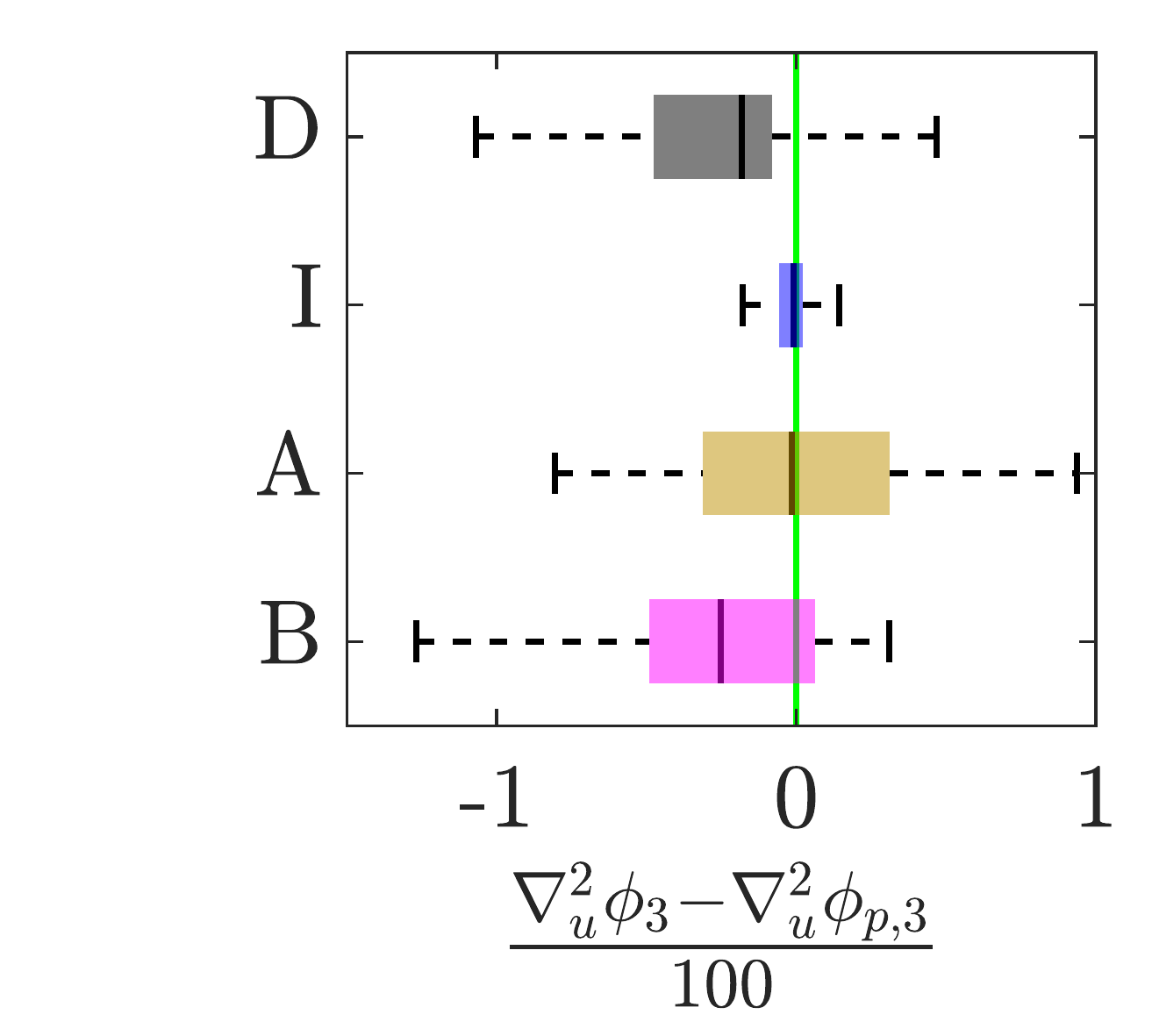}\hskip -0ex
				\includegraphics[trim={2.75cm 0.2cm 0.2cm  0.2cm },clip,width=3.45cm]{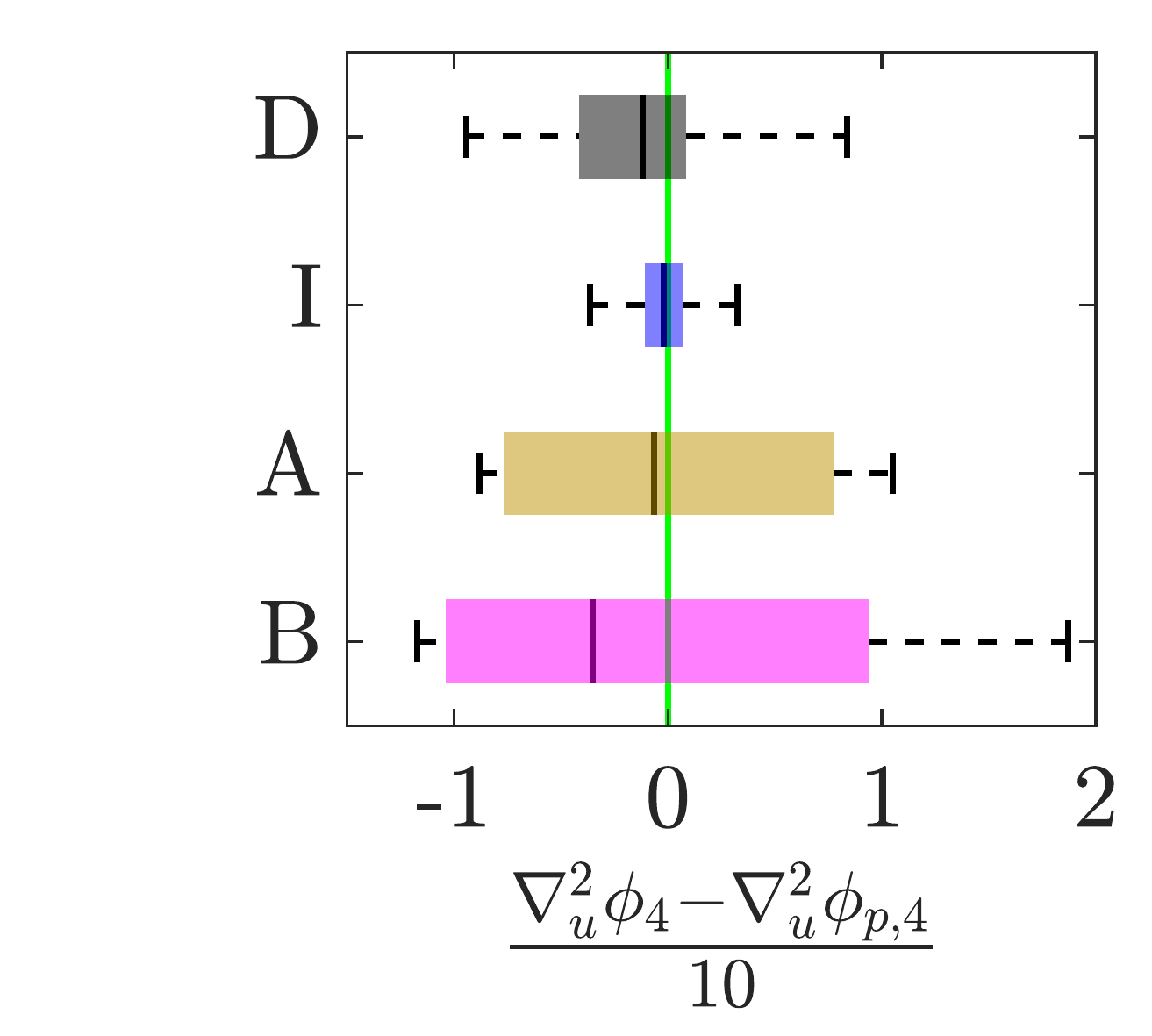} 
			}
			\captionof{figure}{Statistical distributions of the prediction errors on  Hessian of the plant's cost functions at the correction point for the structures D, I, A, and B.}
			\label{fig:5_37_Exemple_5_4_Results}
		\end{minipage}

\section{Case study: the Tennessee Eastman challenge problem}

\begin{minipage}[h]{\linewidth}
	\vspace*{0pt}
	{\centering	
		\includegraphics[width=14cm]{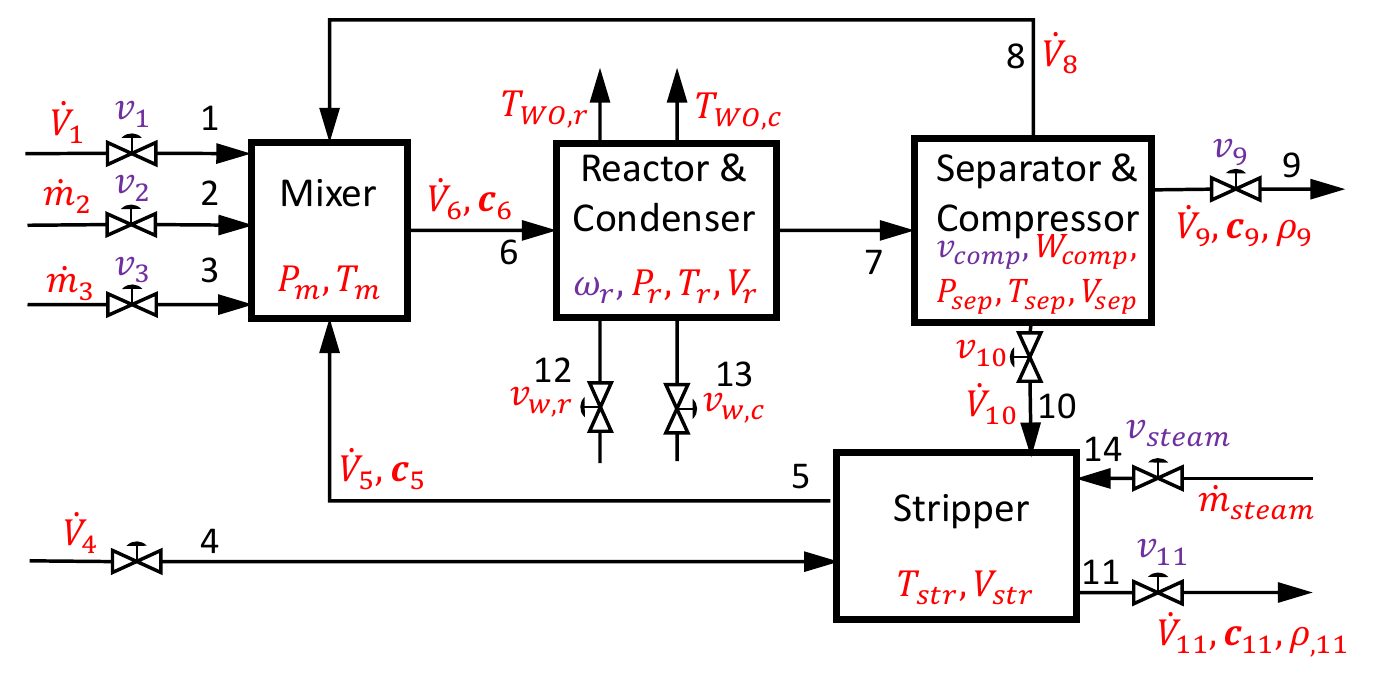}
	}
	\captionof{figure}{Unit based description of the Tennessee Eastman process. The variables manipulated by the controllers as well as the measured ones are in violet and red, respectively. More details about notations are available in the Appendix~\ref{app:B_TE}.}
	\label{fig:5_38_TE_challenge_process}
\end{minipage} 

\textit{(Disclaimer: In this study, only MA and two versions of \textit{internal modifier adaptation} (IMA) are implemented. IMA is basically KMFCaA without filter based constraints, filter adaptation, and curvature adaptation. Therefore, the only difference between MA and IMA is the correction structure. The two versions of IMA are named A and B to indicate which type of correction structure is used. Of course, the version A uses the structure A, and the structure B the B. Finally, this case study and its results were published in \cite{Papasavvas:2019c}).}\\

The Tennessee Eastman plant (TE) is simulated using a benchmark FORTRAN code (\cite{Bathelt:2015}). The plant is known to be open-loop unstable and is equipped with a decentralized control scheme (\cite{Ricker:1996}). 

The model at hand is a simplified version of a model available in the literature (\cite{Ricker:1995b}) sharing the same inputs as the plant of \cite{Bathelt:2015}. Nothing else than the plant information considered available in the Tennessee Eastman challenge problem statement\cite{Downs:1993}, is considered available in this work, e.g. for modeling purposes, that is:
\begin{itemize}
	\item Physical properties of the different components $A$, $B$, $C$, $D$, $E$, $F$, $G$, $H$, at 100 ${}^oC$. 
	\item An intentionally simplified set of reactions (as in \cite{Ricker:1995b}) to serve as the model and to implement structural plant-model mismatch (i.e.,  the simulated reality of the FORTRAN code considers more reactions): 
	\begin{align*}
	A_{(g)} + C_{(g)} + D_{(g)} \rightarrow \ & G_{(\ell)}, \\
	A_{(g)} + C_{(g)} + E_{(g)} \rightarrow \ & H_{(\ell)}, \\
	1/3A_{(g)} + D_{(g)} + 1/3E_{(g)} \rightarrow \ & F_{(\ell)}.
	\end{align*} 
	\item The compositions and temperatures of the inlet streams to the plant (i.e., streams 1, 2, 3, and 4 in Figure~\ref{fig:5_38_TE_challenge_process}).
	\item The sub-costs and -constraints of the plant,
	\item The measured variables. Densities $\rho_{9}$ and $\rho_{11}$ of the flows 9 and 11 are also known\footnote{These measurements are assumed to be available so that costs are functions of manipulated and measured variables only, for the reasons discussed in Section 3. Unlike in  \cite{Ricker:1995}, actual cost and constraints values are used. Also, the steady-state concentrations and volumetric flowrates of Stream 5, even though not directly measured, can be easily inferred from the other steady-state  measurements at the mixer.\label{footnote_exacte_fonctions}}. 
\end{itemize}
The only marginal differences with \cite{Ricker:1995b} are in order:
\begin{itemize}
	\item The original compositions of the streams 1, 2, 3, 4 of \cite{Downs:1993} are used,
	\item The temperature $T_6$ is not fixed. Instead, $F_6$ $F_9$, $\dot{m}_{steam}$, $P_m$, $T_{str}$, and $c_{E,6}$ are fixed\footnote{Details about these variables are available in Appendix~\ref{app:B_TE}.}, 
	\item The stripper is assumed to be ideal, i.e., it perfectly separates $G$ and $H$ from the other species, 
	\item The quantity of each component in each equipment is not modeled to infer the liquid volumes $V_r$, $V_{sep}$ and $V_{str}$. Instead, these volumes are directly manipulated through appropriate controller set-points.
\end{itemize}

The degrees of freedom for the controlled TE plant are indeed the set-points of the controllers, which are in order: $u_1 := \dot{m}_{11}^{sp}$ (mass flowrate of Stream 11 (kg/h)); $u_2 := c_{G,11}^{\%,sp}$ (molar concentration of $G$ in Stream 11 (\%)); $u_3 := c_{A}^{sp}$ (amount of $A$ relative the amount of $A+C$ in stream 6 (\%)); $u_4 := c_{AC}^{sp}$ (amount of $A+C$ in stream 6 (\%)); $u_5 := T_{r}^{sp}$ (reactor temperature (${}^{o}C$)) $u_6 := P_{r}^{sp}$ (reactor pressure (kPa)); $u_7 := V_{r}^{\%,sp}$ (reactor level (\%)); $u_8 := V_{sep}^{\%,sp}$ (separator level (\%)); $u_9 := V_{str}^{\%,sp}$ (stripper level (\%)); $u_{10} := v_{steam}$ (steam valve position (\%)); $u_{11} := \omega_{r}$ (reactor agitation speed (\%)).

According to \cite{Ricker:1995}, $u_6$ to $u_{11}$ can be fixed, and so are they in this work. More precisely, the optimal reactor liquid level and pressure are kept at their lower and upper bounds, respectively. To avoid constraint violations during the transients to steady state - something that has been observed during dynamical simulations -, the values of $P_r^{sp}$ and $V_{r}^{\%,sp}$ are set to $2800$ kPa and $65$\%, slightly backed-off from their ``real'' limit values of $2895$ kPa and $50$\%. Also, $V_{sep}^{\%,sp}$ and $V_{str}^{\%,sp}$ having negligible effects on the operating conditions of the plant, they are fixed to $50$\%, again far enough from their bounds. The energy consumed for stirring being not taken into account, the agitator speed is set to $100$\%. Finally, the steam valve position is set at $1$\%, since it has been proved in \cite{Ricker:1995} that steam has not effect on the steady state and costs money. 

The model pre-processing procedure has been applied and the resulting networks ($\mathcal{N}^{sm}$ and $\mathcal{N}^{sc}$), with detailed equations, are given in the Appendix~\ref{app:B_TE}. Notice that, contrary to the first case study, the structure of $\mathcal{N}^{sm}$ does not mimic the structure of the plant because \emph{not all} the process variables that connect the reactor, the condenser, the separator, the compressor, and the stripper are measured. Hence the difference between Figure~\ref{fig:5_38_TE_challenge_process} (which is a units-based description of the open-loop TE) and Figure~\ref{fig:B_TE_description_1} (in the Appendix~\ref{app:B_TE}), which depicts $\mathcal{N}^{sm}$.

Only the first of the six operating modes of \cite{Downs:1993} are considered. No disturbances are considered, since the purpose of this case study is to illustrate the implementation of IMA-A/B. Therefore, the model-based optimization problem considers four sub-costs and twelve sub-constraints, detailed in Appendix~\ref{app:B_TE}. Basically, the liquid levels in the reactor, the separator and in the stripper are bounded. The pressure and temperature of the reactor are upper bounded. The ratio between the flowrates of $G$ and $H$ in  stream~11 is 50/50 to satisfy quality requirements. There is an equality constraint on the production rate, i.e., the mass flowrate of $G$ and $H$ in stream 11 has to be equal to 14076 kg/h. The aggregated cost is composed of (i) the cost related to the power consumption of the compressor, (ii) the steam consumed by the stripper, (iii) the  flowrate in the purge, and (iv) the  flowrate in the production stream~11.

Figure~\ref{fig:5_39_SimResultsTE} depicts the simulation results for 30 RTO iterations with MA, IMA-A and -B starting from conservative initial inputs $\bm{u}_0=[10,30,60,50,120]^{\sf{T}}$. The filter gain used is $K=0.9$\footnote{With $K=1$, results are very similar, although more oscillations are observed before convergence for IMA-A and -B.}. One can see that MA does not converge to the plant optimum and, instead, oscillates around an infeasible point. On the other hand, IMA-A and -B lead to similar results and reach the plant optimum in about 5 iterations. Notice that the optimal steady-state reached here over performs slightly than the one from \cite{Ricker:1995}. This is due to the fact exact cost and constraint values are considered here (see footnote~\ref{footnote_exacte_fonctions}). Results are also different from \cite{Golshan:2000}, whereby the  two-step approach, i.e., online identification of the models parameters followed by re-optimization of the updated model, has been applied leading to convergence to a sub-optimal point, which happens to be quite luckily close from the true plant optimum. IMA converges to the true plant optimum, i.e. $m_{11}^{sp} = 22.797$  kg/h, $\% G^{sp} = 53.83$,  $c_A^{sp} = 63.22$, $c_{AC}^{sp} = 50.92$, and $T_r = 122.84 {}^oC$ for an operating cost of $113.53$\$/h, while satisfying all operating constraints.\\

\noindent
\begin{minipage}[h]{\linewidth}
	\vspace*{0pt}
	{\centering	
	\includegraphics[width=\linewidth]{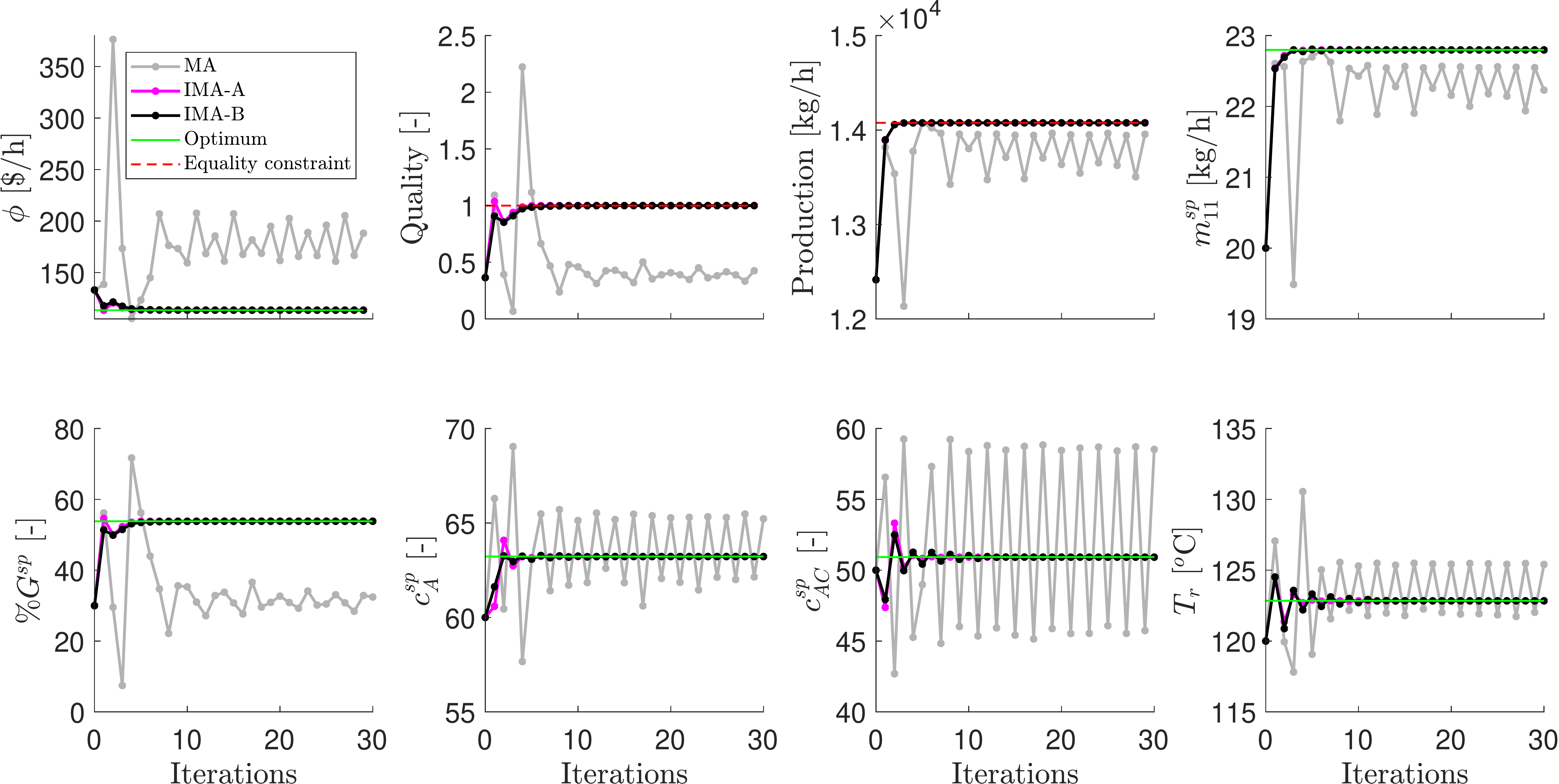}
	}
	\captionof{figure}{\textbf{TE:} Simulation results.}
	\label{fig:5_39_SimResultsTE}
\end{minipage}

\section{Conclusion}

In this chapter it has been shown that other correction structures than the classical direct and indirect ones exist. These new structures, called A and B, bring three improvements  in the following fields: 
\begin{itemize}
	\item Conceptual:  Because they enable the use of all the measurements of an experiment.
	\item Computational: Because they allow for parallel/distributed computing since each SM can be treated independently of he others.
	\item Efficiency: Because they generally produce better local corrections.
\end{itemize}
However, those pros come with the following cons:  
\begin{itemize}
	\item Complexity: Because their implementation comes with a great increase of complexity in terms of model pre-processing and implementation.  
\end{itemize}
In this chapter, it is seen how best to use the set of measurements available at the end of one experiment. In the next chapter, it will be discussed how the whole history of the experiments could be exploited.

	\chapter{An autopilot for steady-state processes}
	\label{Chap:6_ASP}

\section{Preliminaries}

In this chapter almost all the functions of S-ASP are reviewed and improved. The only ones that are not modified are both the signal purifier, for which no attempt is made to improve the SS detector; and the statistical converter, to which one believes not much improvements are possible. 
To simplify the developments in this chapter, it has been decided to work only with sub-functions $f_{(i)}$ and $f_{p(i)}$ of $\bm{f}$ and $\bm{f}_p$, and to simplify the notation, it has also been decided to remove the index $(\cdot)_{(i)}$: 
\begin{align}
	y_{p(i)} \Leftrightarrow \ & y_{p}, &
	f_{(i)} \Leftrightarrow  \ & f, &
	f_{p(i)} \Leftrightarrow \ & f_{p}.
\end{align}

Up to now the disturbances and their effects on the plant have been completely ignored; indeed, the focus has been mainly on the development of theoretical RTO. Taking into account the \textit{measured} disturbances in the decision making is, however, very simple since one just needs to consider them as virtual inputs $\bm{u}^\prime$ that are not manipulable and whose value is fixed by an equality constraint $\bm{u}^\prime=\bm{d}_p^{\prime}$. So any change on their part can be considered as an update of a constraint, i.e. should trigger the decision trigger as a change of constraint would. 

The real difficulties come from \textit{unmeasured} disturbances $\bm{d}_p^{\prime\prime}$. To manage them, the notion of \textit{consistent time period} (CTP) is introduced:
\begin{Definition} \label{def:6_PTC}
	\textbf{(consistent time period (CTP))}
	Over a CTP, the unmeasured disturbances are (almost) constant, i.e. $\bm{d}^{\prime\prime}\approx c^{te}$. As a result, the functions of the plant can be reduced to the simple relationship: 
	\begin{equation}
		y_p=f_p(\bm{x}),
	\end{equation}
	where $\bm{d}^{\prime\prime}$ would play the role of an unknown fixed parameter. Therefore, during a CTP the plant can be considered as a function whose inputs and outputs are the same as those of the model.
\end{Definition}

\section{Assumptions about the plant}

\subsection{Two assumptions}

\begin{figure}[H]
	\vspace*{0pt}
	{\centering
		\includegraphics[width=14cm]{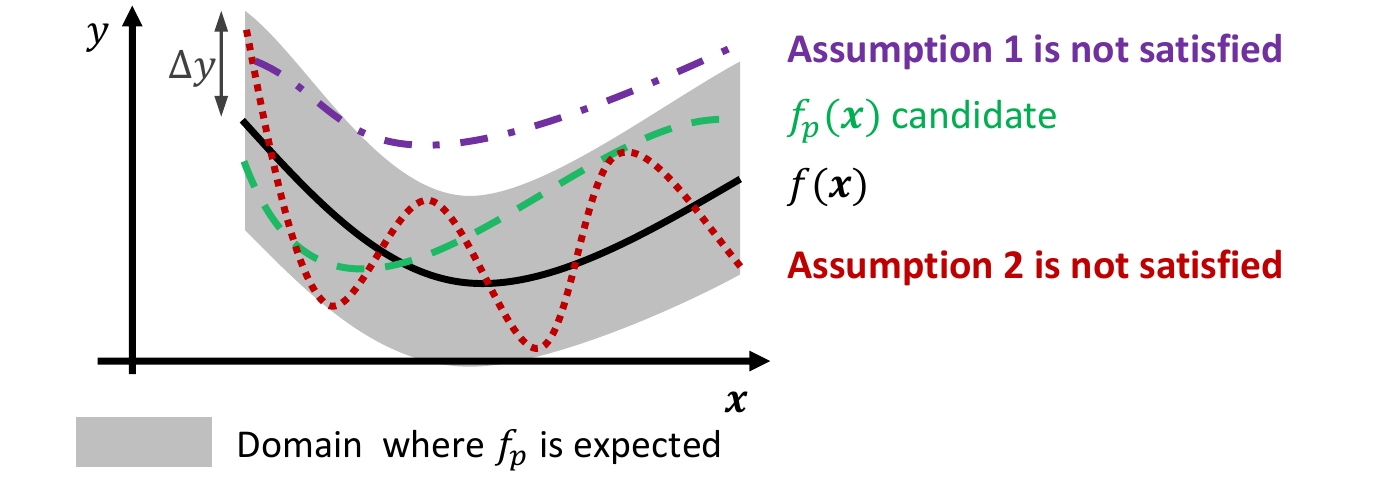} \\
	}
	\caption{Graphical illustration of the Assumptions~\ref{ass:6_1_Erreur_fp_f_bornee} and \ref{ass:6_2_Meme_Courbusre} on $f$ and $f_p$}
	\label{fig:AssumptionsFp}
\end{figure}

Even if $f_p$ is not known, the model at hand ($f$) is supposed to be a ``relevant'' approximation in the sense that it is supposed to represent the global trends of the plant's behavior at the level of detail that is appropriate for decision making.  
For example, if the model is supposed to be used to choose the temperature of a reactor with an accuracy of one degree, one does not need to know the details of its behavior variations to the thousandth of a degree. On the basis of this consideration, the following two hypotheses (illustrated on Figure~\ref{fig:AssumptionsFp}) are made.
\begin{Assumption} \label{ass:6_1_Erreur_fp_f_bornee}
	The error between the functions $f_p$ and $f$ is bounded by a known constant $\Delta y$: 
	\begin{align} \label{eq:6_3_Hypothese_Borne_fp_f}
		|f_p(\bm{x}) - f(\bm{x})| < \ & \Delta y,  & & \forall \bm{x}\in\amsmathbb{R}^{n_x}.
	\end{align}
\end{Assumption}

\begin{Assumption} \label{ass:6_2_Meme_Courbusre}
	The functions $f_p$ and $f$ have globally the same curvature\footnote{Simply put, the curvature at a given point is the quantification of the difference between a curve and a straight line, or between a surface and a plane, at that given point in space.}. More precisely, their minimum radii of curvature with respect to all directions in the input space are of the same order of magnitude. 
\end{Assumption}

\subsection{A Gaussian process}

These two assumptions imply that $f_p$ belongs to a set of functions whose smoothness and position w.r.t.  $f$ are a priori known.  This set of functions can be mathematically characterized by a Gaussian process (GP) with a squared-exponential with automatic relevance determination kernel (SE-ARD). 

Let us now introduce what a GP with a SE-ARD is. In the same way that a normal distribution characterizes the value of a random variable $\epsilon$ around a mean value $\epsilon_0$ with a variance  $\sigma_{\epsilon}^2$, i.e. $\epsilon\sim \mathcal{N}(\epsilon_0,\sigma_{\epsilon}^2)$,  a GP characterizes a distribution of random functions $f_p$ around a mean function $f$ with a kernel $\k_{\f\f}(\bm{x}_i,\bm{x}_j)$, i.e. $f_p\sim\mathcal{GP}(f,\k_{\f\f}(\bm{x}_i,\bm{x}_j))$. The kernel function characterizes the covariance of  $f_p$   with respect to itself at two locations  $\bm{x}_i$ and $\bm{x}_j$:
\begin{equation} \label{eq:6_4_Def_Generale_du_kernel}
	\k_{\f\f}(\bm{x}_i,\bm{x}_j) := \amsmathbb{E}[(f_p(\bm{x}_i)-f(\bm{x}_i)
	(f_p(\bm{x}_j)-f(\bm{x}_j))].
\end{equation}
By choosing the SE-ARD kernel, it is chosen that: 
\begin{equation}
	\k_{\f\f}(\bm{x}_i,\bm{x}_j):= \sigma_{f}^2
	\exp\big(-\frac{1}{2}(\bm{x}_i-\bm{x}_j)^{\rm T}\bm{\L}(\bm{x}_i-\bm{x}_j)\big),
\end{equation}
where $\sigma_f^2$ and $\bm{\L}$ are hyperparameters with the following meanings: 
\begin{itemize}
	\item $\sigma_f^2$ is the a priori variance of the values of $f_p(\bm{x}_i)$ around $f(\bm{x}_i)$. Indeed, it can easily be found that \eqref{eq:6_4_Def_Generale_du_kernel} evaluated at $\bm{x}_i=\bm{x}_j$ gives: 
	$$
	\k_{\f\f}(\bm{x}_{i},\bm{x}_{i}|\sigma_{f}^2,\bm{\L}) = 
	\amsmathbb{E}
	\left[
	\left(f_p(\bm{x}_{i}) - f(\bm{x}_{i})\right)^2
	\right]
	=
	\sigma_{f}^2. 
	$$
	So, a priori $f_p(\bm{x}_i)\sim \mathcal{N}(f(\bm{x}_i),\sigma_f^2)$ and 
	\begin{equation} \label{eq:6_6_confidence_domain}
		p(|f_p(\bm{x}_i)-f(\bm{x}_i)| \leq 2\sigma_f) = 95.45\%.
	\end{equation}
	\item $\bm{\L}$ is a diagonal matrix whose diagonal contains length scale parameters $\ell_1^{-2}$,..., $\ell_{n_x}^{-2}$ that characterize the curvature of $f_p$ w.r.t. its $n_x$ inputs. 
\end{itemize}
Therefore, any function $f_p$ whose location around a mean function $f$ is bounded, and whose curvatures w.r.t. all directions in its inputs space are known, can be defined with the following GP:
\begin{equation} \label{eq:6_7_GP_Fp}
	f_p \sim \mathcal{GP}(f,\k_{\f\f}(\bm{x}_i,\bm{x}_j|\sigma_{f}^2,\bm{\L})).
\end{equation}

The hyperparameters of this GP can be identified from the model, given Assumptions~\ref{ass:6_1_Erreur_fp_f_bornee} and \ref{ass:6_2_Meme_Courbusre}, as discussed next.

\subsection{Estimation of the plant's GP hyperparameters}

If the boundary on the error $|f_p(\bm{x})-f(\bm{x})|$ of \eqref{eq:6_3_Hypothese_Borne_fp_f} is equated to the confidence domain \eqref{eq:6_6_confidence_domain}, it can be determined that an appropriate choice of $\sigma_f^2$ would be: 
\begin{equation} \label{eq:6_8_choix_sigmaf}
	\sigma_{f}^2 = \left(\Delta y/2\right)^2.
\end{equation}

On the other hand, the choice of $\bm{\L}$ should be based on an analysis of the function  $f$, consisting of reconstructing the model with a GP using a SE-ARD kernel using a model selector  (i.e. an auto-tuning of the GP's hyperparameters, see Chapter~5 of \cite{Rasmussen:2006}). Once this approximation is constructed, the matrix $\bm{\L}$  identified to approximate the model should be,  according to the Hypothesis~\ref{ass:6_2_Meme_Courbusre}, a relevant choice of $\bm{\L}$ for \eqref{eq:6_7_GP_Fp}. 

\begin{figure}[h]
	\includegraphics[width=1\linewidth] {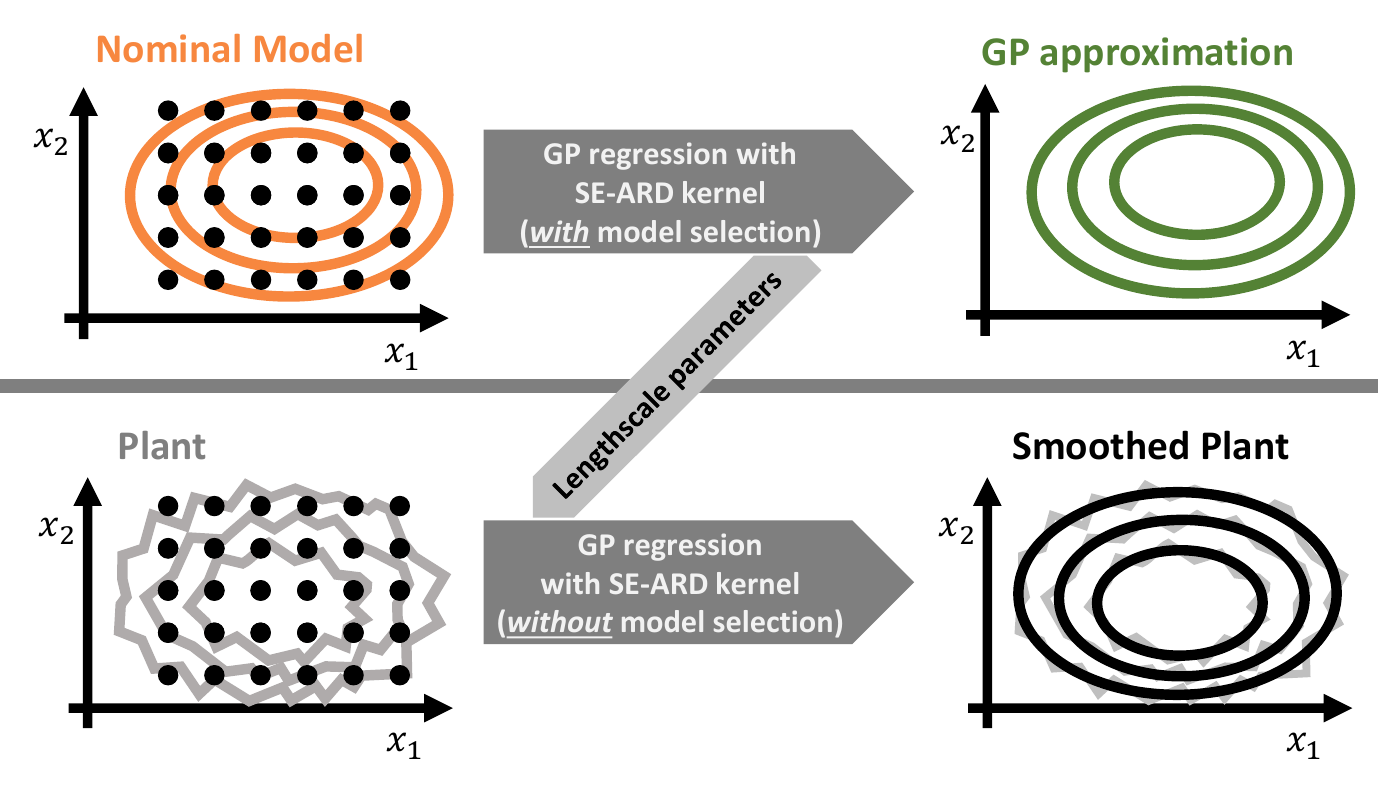}
	\caption{Smoothing effect introduced by the use of GP regression which would erase the least important local details of $f_p$ in order to keep only the trends that are useful for decision making.}
	\label{fig:6_2_SmoothPlant}
\end{figure}

\begin{Remark}
	\textbf{(Smoothing effect and deleting of ``irrelevant'' details)}
	It is clear that a function $f_p$ of a real plant usually has details that the model cannot contain. So, by approximating $f_p$ with a GP using the hyperparameters of the model, all these details would be erased. This can be interpreted as a loss of information since a simplified version of the outputs are being used, rather than the actual outputs of $f_p$. This can be interpreted as a simplification of $f_p$ which deletes all insignificant details and keeps only the main trends useful for decision making as shown on Figure~\ref{fig:6_2_SmoothPlant}. Of course, some of those details may be useful to find precisely where the actual  plant's optimum is. In such a case, the length scale hyperparameters of the model could be identified and assume the plant's ones are slightly smaller. In this later case, it would be up to the engineers to define what ``slightly smaller'' means. 
\end{Remark}

\section{The improved model-data combiner}
	In addition to knowing that $f_p$ is a GP of the form of \eqref{eq:6_7_GP_Fp}, one has access to the results of the  $\ell$$+$$1$ experiments already conducted on the plant: $\{(\widehat{\bm{x}}_{0}, \widehat{y}_{p,0}, \widehat{\sigma}_{y,0}^2
	, \widehat{\nu}_{y,0}^2,$ $\widehat{\bm{\Sigma}}_{x,0}
	, \widehat{\bm{N}}_{x,0}), ... ,$ $(\widehat{\bm{x}}_{\ell}, \widehat{y}_{p,\ell},\widehat{\sigma}_{y,\ell}^2,\widehat{\nu}_{y,\ell}^2,\widehat{\bm{\Sigma}}_{x,\ell}
	, \widehat{\bm{N}}_{x,\ell} )\}$, where  
	\begin{itemize}
		\item $\widehat{\bm{x}}_{i}$ is the estimated average of the plant inputs, 
		\item $\widehat{\bm{y}}_{p,i}$ is the estimated average of the plant outputs,
		\item $\widehat{\sigma}^2_{y,i}$ is the \textit{estimate} of the variance representing the uncertainty on the mean of the plant outputs,    
		\item $\widehat{\nu}^2_{y,i}$ is the \textit{estimate} of the variance representing the noise on the plant output,   
		\item $\widehat{\bm{\Sigma}}^2_{x,i}$ is the \textit{estimate} of the covariance matrix representing the uncertainty on the mean of the plant inputs, 
		\item $\widehat{\bm{N}}^2_{x,i}$ is the \textit{estimate} of the covariance matrix representing the noise on the plant inputs, 
	\end{itemize}
	of the $i$-th experiment. For the time being, the measurements of the noise variances ($\widehat{\nu}^2_{y,i}$, $\widehat{\bm{N}}^2_{x,i}$) is left to one side (in the conclusion one suggests an idea of what could be done with those measures), and the focus is  only on the other data. The links between these and the real values are the following: 
	\begin{align*}
		\widehat{\bm{x}}_{i} = \ & \bm{x}_{i} + \bm{\epsilon}_{x,i},&  \bm{\epsilon}_{x,i}\sim \ &   \mathcal{N}(\bm{0},\widehat{\bm{\Sigma}}_{x,i}), \\
		\widehat{y}_{p,i} = \ & y_{p,i} + \epsilon_{y,i},&  \epsilon_{y,i}\sim \ &   \mathcal{N}(0,\widehat{\sigma}_{y,i}^2).
	\end{align*}
	By definition, function  $f_p$ connects those measures in the following way:
	\begin{equation} \label{eq:6_9_yp_fp}
		\widehat{y}_{p,i} = f_p(\widehat{\bm{x}}_{i} - \bm{\epsilon}_{x,i}) + \epsilon_{y,i}.
	\end{equation}
	\eqref{eq:6_9_yp_fp} can be reformulated so as to separate the uncertain terms from the certain terms in order to be in the classical configuration of a regression problem and to be able to use standard regression methods to identify $f_p$.  To do this, one approximates \eqref{eq:6_9_yp_fp} with a Taylor series centered on  $\widehat{\bm{x}}_{i}$: 
	\begin{equation*}
		\widehat{y}_{p,i} =   f_p(\widehat{\bm{x}}_{i} ) + \epsilon_{y,i} - \bm{\epsilon}_{x,i}^{\rm T} \nabla_{\bm{x}} f_p(\widehat{\bm{x}}_{i}) + \mathcal{O}(\bm{\epsilon}_{x,i}^2),
	\end{equation*}
	and it is assumed that $\mathcal{O}(\bm{\epsilon}_{x,i}^2)$ is negligible:
	\begin{equation}  \label{eq:6_10_Link_Fp_Measurements}
		\widehat{y}_{p\{i\}} =  f_p(\widehat{\bm{x}}_{i}) + \epsilon_{xy,i}, \qquad \epsilon_{xy,i} \sim \mathcal{N}(0,\sigma_{xy,i}^{2}),
	\end{equation}
	where
	\begin{equation} \label{eq:6_11_MeasurementsUncertainty}
		\sigma_{xy\{i\}}^2 := \big| \nabla_{\bm{x}} f_p(\widehat{\bm{x}}_{i})\big|^{\rm T}\widehat{\bm{\Sigma}}_{x,i}\big|\nabla_{\bm{x}} f_p(\widehat{\bm{x}}_{i}) \big| + \widehat{\sigma}_{y,i}^2. 
	\end{equation}
	At this stage, $\widehat{y}_{p,i}$ has been separated into a sum of certain and uncertain terms.  Therefore, this is a standard form of regression problem where the goal is to identify the function  $f_p\sim \mathcal{GP}(f,\k_{\f\f}(\bm{x}_i,\bm{x}_j|\sigma_{f}^2,\bm{\L}))$ linking inputs $\widehat{\bm{x}}$ to outputs $\widehat{y}_p$, subject to measurement errors $\epsilon_{xy,i}\sim\mathcal{N}(0,\sigma_{xy,i}^2), \forall i$. The method for solving this type of regression problem when the database  contains  points and directional derivatives is given in Appendix~\ref{App:C_PartialDeriv}. 
	
	However, a technical difficulty arises. To determine the variance,	$\sigma_{xy,i}^2$, of $\epsilon_{xy,i}$,the gradient of $f_p$ at $\widehat{\bm{x}}_{i}$ need to be known.  However, to identify (regress) $f_p$ and thus estimate its gradient, $\sigma_{xy,i}^2$ is required. Therefore, the estimates of $f_p$ and $\sigma_{xy,i}^2$ are interdependent and should be calculated simultaneously, as such, the calculation is very costly. \cite{Mchutchon:2011} suggests to (i) assume that the gradient of $f_p$ is null to make the first estimates of $\sigma_{xy,i}^2$ and $f_p$ to then (ii) use this approximation of  $f_p$  to re-estimate its gradient and recompute  $\sigma_{xy,i}^2$ and $f_p$. If repeated until convergence, this procedure would converge on the true function  $f_p$ with the appropriate estimate of  $\sigma_{xy,i}^2$. Here, an even simpler (but less efficient) alternative is proposed consisting of replacing this unknown gradient of $f_p$ by an overestimation of its value, $\bm{b}(\widehat{\bm{x}}_{i})$. Of course, using an overestimation of this gradient implies that an overestimation of the measurement uncertainties $\widehat{\sigma}_{xy,i}$ is used. This has the effect that the measures are underused (because they are considered ``less reliable'' since they are subject to more uncertainties) in the regression and that the regression is therefore ``less efficient'' than it could be in an ideal case.  This overestimation can be derived from the hyperparameters and the model's gradients:
	 \begin{equation} \label{eq:6_12_UpperBound_NablaFp}
	 	\big| \nabla_{\bm{x}} f_p(\widehat{\bm{x}}_i)\big| \leq
	 	\bm{b}(\widehat{\bm{x}}_i) := \big| \nabla_{\bm{x}} f(\widehat{\bm{x}}_i)\big| + \frac{3}{\sqrt{2}}\sigma_f\bm{\ell}^{-1}
	 \end{equation}
	where $\bm{\ell}^{-1}:=[\ell_1^{-1},...,\ell_{n_x}^{-1}]^{\rm T}$. The mathematical foundations of \eqref{eq:6_12_UpperBound_NablaFp} are given in the Appendix~\ref{App:C_PartialDeriv} (equation \eqref{eq:C_13_wjhgvfrcje}) and are mainly based on the probability distribution of the derivative of a GP (SE-ARD) whose hyperparameters $\{\sigma_f^2,\bm{\L}\}$  are known.
	
	By combining the a priori distribution of $f_p$ \eqref{eq:6_7_GP_Fp} with measurements of values and gradients, e.g.  $\{\widehat{\bm{x}}_i, \widehat{y}_{p,i}, \sigma_{xy,i}^2\}$ and $\{\widehat{\bm{x}}_j, \widehat{\nabla_{\v_j} y}_{p,j}, \bm{\v}_{j},, \sigma_{\v,j}^2\}$, the conditional distribution of $f_p$ can be computed, given those data. Computing this a posteriori  distribution is equivalent to performing a GP regression. The theoretical basis for this type of regression is given in Chapter 2 of \cite{Rasmussen:2006} for the case where the data set is only made of points. For the case where the data set contains points \textit{and gradients}, \cite{Solak:2003} is recommended. The details of this regression procedure are given in Appendix~\ref{App:C_PartialDeriv}. 
	
	The result of this regression at a point $\bm{x}_*\in\amsmathbb{R}^{n_x}$ is  \eqref{eq:C_Fp_a_posteriori}:
	\begin{align} 
		& f_p(\bm{x}_*) \sim  \mathcal{N}(
		f(\bm{x}_*) + \bm{\cc}_{*}^{\rm T}\bm{\CC}^{-1}(\bm{y}_p-\bm{y})
		,
		\cc_{**} - \bm{\cc}_{*}^{\rm T}\bm{\CC}^{-1}\bm{\cc}_{*}
		)), \nonumber 
	\end{align}
	where all the terms constituting this equation are defined in Appendix~\ref{App:C_PartialDeriv}. This equation can be rewritten as a sum of  $f$ with a correction function\footnote{Where the index $(\cdot)_k$ indicates that it is a correction function established on the data of $k$  experiments conducted so far.} $\mu_k^f$, both evaluated at $\bm{x}_*$:
	\begin{align}  \label{eq:6_13_def_eta_kf}
		f_p(\bm{x}_*)    =  \ & \ f(\bm{x}_*) + \mu_k^f(\bm{x}_*), &
		\mu_k^f(\bm{x}_*)  \sim \ & 
		\mathcal{N}\big(\bm{\cc}_{*}^{\rm T}\bm{\CC}^{-1}(\bm{y}_p-\bm{y})
		,
		\cc_{**} - \bm{\cc}_{*}^{\rm T}\bm{\CC}^{-1}\bm{\cc}_{*}
		)\big).
	\end{align}
Clearly the correction function, $\mu_k^f$, is a correction function of an indirect correction structure (I). It should be noted that the latter has two particularities: (i) it is non-linear, and (ii) it is not only an estimate of the expectation of $f_p$, but also an estimate of the  uncertainty on this expectation.   From those two particularities, three important remarks can be made:
	
	\begin{Remark} \textbf{Concerning the particularity (i) -- Part 1:} The use of non-linear functions to correct the model was first proposed in \cite{Faulwasser:2014} where it is suggested to use quadratic correction functions (instead of affine). Of course, one criticism that could be made of this work is that if measuring the gradients of the plant is already a challenge in practice, then obtaining the Hessians is generally unthinkable. However, \textit{if the Hessians of the plant are not perfectly accessible, even their very approximate estimates are very interesting to exploit}. Indeed, it is their signs (defined positive, negative, or other) that play an important role in decision making because it is their signs that makes the difference between a minimum, a maximum and a saddle point. 
		The strategy that is used with KMFCaA consists of forcing the cost and active constraint functions to be positive definite at the correction point. This has the effect of enforcing the equilibrium condition at points that are not minimums of the plant as shown in Example~\ref{ex:2__4_PtsStable_Et_Instables}.
	If instead of blindly making these Hessians positive definite, only``their sign'' (and not necessarily their values) would be corrected, then maximums and saddle points would no longer be equilibrium points as shown in Example~\ref{ex:6_1_Correction_en_signe_du_Hessien}. 
	Therefore, one of the interests of using non-linear corrections could be to remove this small limitation of KMFCaA\footnote{Small because as it has been seen, even if these points are equilibrium points, they are unstable, so convergence on them is unlikely}. 
	\end{Remark}

	\begin{exbox} 
		\label{ex:6_1_Correction_en_signe_du_Hessien}
		\textbf{(Approximate Hessian corrections -- `` correcting the \textit{sign}, not necessarily the \textit{value})}	
		Let's take Example~\ref{ex:2__4_PtsStable_Et_Instables} and, in addition to the affine corrections of the model, add \textit{inaccurate} second order corrections. Then, the model updated at $u_k$ is:
		\begin{equation*}
			\phi_k(u) = \phi_p(u_k) + \nabla_u\phi_p|_{u_k}(u-u_k) + \frac{\widehat{\nabla_{uu}^2\phi_p}|_{u_k}}{2}(u-u_k)^2.
		\end{equation*} 
		where $\widehat{\nabla_{uu}^2\phi_p}|_{u_k}$ is the inaccurate estimate of the plant's Hessian $\nabla_{uu}^2\phi_p|_{u_k}$. It is defined as the following noisy finite difference approximation
		\begin{align*}
			\widehat{\nabla_{uu}\phi_p}|_{u_k}  = \ & \frac{\phi_p(u_k+h) + \varepsilon_1 - 2(\phi_p(u_k)+\varepsilon_2) + \phi_p(u_k-h) + \varepsilon_3 }{h^2}, \\
			= \ &\nabla_{uu}\phi_p|_{u_k} + \mathcal{O}(h^2) + \frac{\varepsilon_{123}}{h^2}, 
		\end{align*}  
		where $(\varepsilon_{1}, \varepsilon_{2}, \varepsilon_{3}) \sim  \mathcal{U}(-0.01, 0.01)$ are the measurement errors on the evaluations of $\phi_p(u_k+h)$, $\phi_p(u_k)$, and $\phi_p(u_k-h)$, respectively.  The combined effect of these errors is characterized by the variable  $\varepsilon_{123}\sim\mathcal{U}(-0.04, 0.04)$. Then, there are two possible choices of $h$:  
		
		If one chooses to use a \textit{large} $h$ then the estimate $\widehat{\nabla_{uu}^2\phi_p}|_{u_k}$ is: 
		\vspace{-\topsep}
		\begin{itemize}[noitemsep]
			\item only marginally affected by measurement uncertainties $\varepsilon_{123}$, since these are multiplied by $1/h^2$, 
			\item majorly affected by the term $\mathcal{O}(h^2)$. 
		\end{itemize}
	\vspace{-\topsep}
		If on this broad area the Hessian of the plant does not change, then it is clear that if:  
		\begin{align*}
			 \nabla_{uu}\phi_p|_{u_k} > \ & 0, & \Rightarrow \qquad \nabla_{uu}\phi_p|_{u_k} + \mathcal{O}(h^2)  > \ & 0, \\
			 \nabla_{uu}\phi_p|_{u_k} < \ & 0, & \Rightarrow \qquad \nabla_{uu}\phi_p|_{u_k} + \mathcal{O}(h^2)  < \ & 0.
		\end{align*}

		If one chooses to use a \textit{small} $h$, then the estimate $\widehat{\nabla_{uu}^2\phi_p}|_{u_k}$ is:
		\vspace{-\topsep}
		\begin{itemize}[noitemsep]
			\item majorly affected by measurement uncertainties $\varepsilon_{123}$ since these are multiplied by $1/h^2$, 
			\item only marginally affected by the term $\mathcal{O}(h^2)$. 
		\end{itemize}
	\vspace{-\topsep}
		And as the measurements are close to each other it should be clear that: 
		\begin{align*}
			\nabla_{uu}\phi_p|_{u_k} > \ & 0, & \not\Rightarrow \qquad \nabla_{uu}\phi_p|_{u_k} +  \frac{\varepsilon_{123}}{h^2} > \ & 0, \\
			\nabla_{uu}\phi_p|_{u_k} < \ & 0, & \not\Rightarrow \qquad \nabla_{uu}\phi_p|_{u_k} +  \frac{\varepsilon_{123}}{h^2}  < \ & 0.
		\end{align*}
		
			It is then clear that choosing a large $h$ is the option to choose if one wants a correction in sign and not necessarily in value of the Hessian of the model. 
			So, for the example one chooses to use $h=0.5$ which is a fairly large value  (according to the curves of Figure~\ref{fig:2___13_Exemple2_Fonctions}).
			Then, one evaluates $sol(u_k)$ for a large number of inputs by taking random numbers each time $(\varepsilon_{1}, \varepsilon_{2}, \varepsilon_{3}) \sim  \mathcal{U}(-0.01, 0.01)$. The blue and black curves of Figure~\ref{fig:2___16_Exemple2_LocalCone_perfectHessian} give the functions $sol(u_k)$  of the model with approximate and ideal corrections of its Hessian, respectively. \\
		\begin{minipage}[h]{\linewidth}
			\vspace*{0pt}
			\centering 
			\includegraphics[width=9cm]{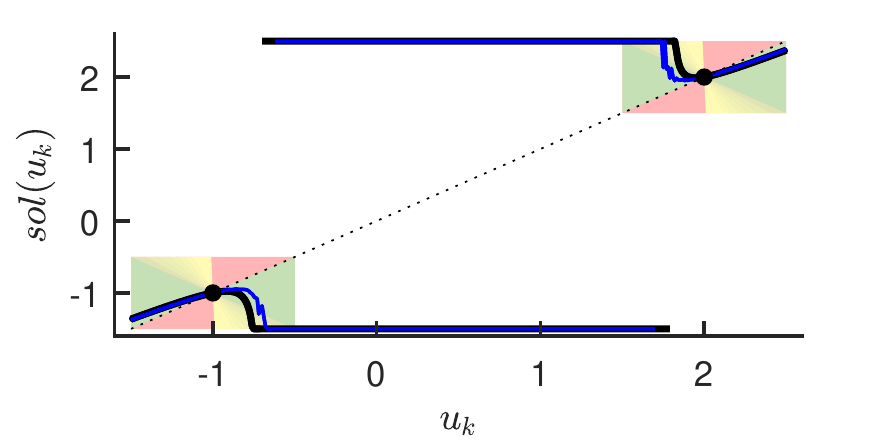}\hskip -0ex
			\includegraphics[width=4.35cm]{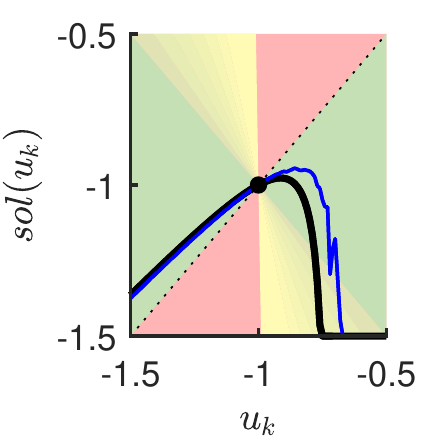}
			\captionof{figure}{Graphical analysis of convergence capabilities at different types of stationary points when using ideal and approximate Hessian sign corrections.}
			\label{fig:2___16_Exemple2_LocalCone_perfectHessian}
		\end{minipage}\\

		Several things can be observed: 
		\begin{itemize}
			\item When $u_k$ is between approximately $-0.8$ and $1.8$ the function $sol(u_k)$ has two solutions: $-1.5$ and $2.5$. This is due to the fact that the model updated at these points is concave. So, the solution on which the solver converges is induced by its initialization. 
			Here, the solutions obtained when the solver is systematically initialized with either $-1.5$, or $2.5$, are superimposed. Hence the presence of two solutions $\forall u_k\in [-0.8,1.8]$.
			One specifies that initializing the solver at the correction point (which is a common practice in RTO since $u_k$ is a priori a good candidate) is not a good idea since it can make a saddle-point stable despite second order corrections. 
			Indeed, around a saddle point the objective function of the optimization problem is such that all its derivatives are almost null (when evaluated by the solver) and the solver may consider this initialization as a minimum of the model and to stop immediately.  It is therefore recommended to systematically initialize the solver on a point which is not the correction point.  Any point relatively close to   $u_k$   should be suitable. 
			\item The blue and black curves are similar and  intersect the line  $sol(u_k)=u_k$ (dotted line) only at  the plant's minimums $u_k = \{-1, \ 2 \}$. Therefore, convergence is only possible on the plant's minimums, and the maximums and saddle points are no longer potential fixed points.    
			\item Without second order corrections it would be necessary to use a filter to make the plant's minimums stable, see Figure~\ref{fig:2___14_Exemple2_LocalCone} (the black curve is in the yellow area).  With an approximate second order correction, the use of filters is no longer necessary to ensure convergence, see Figure~\ref{fig:2___16_Exemple2_LocalCone_perfectHessian}. However, the filter remains useful to accelerate the convergence if it is chosen appropriately as explained in Chapter~\ref{Chap:2_Vers_Une_meilleure_Convergence}.  
		\end{itemize}
	\end{exbox}

\begin{figure}[h]
		\includegraphics[width=1\linewidth] {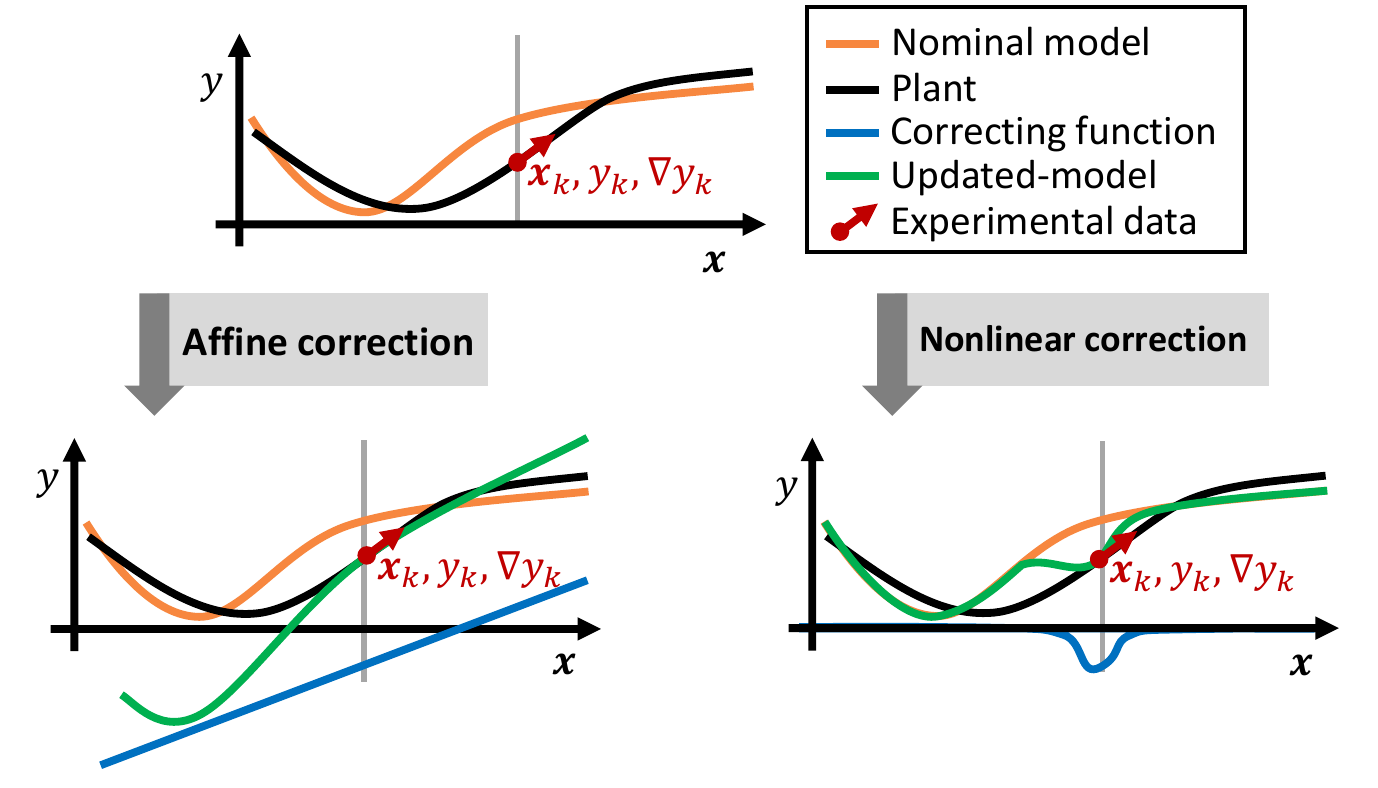}
		\caption{Conceptual illustration of the idea that the model and experimental data are information that is global and local respectively, and should be used as such.}
		\label{fig:6_3_LocalVsGlobalCorrection}
	\end{figure}

	\begin{Remark} \label{rem:6_4_Particularite_i_part_2}
		\textbf{Concerning the particularity (i) -- Part 2:}
		A criticism of simple RTO methods (such as KMFCaA) might be that they tend to use experimental data collected at a specific point in the input space to make global corrections via affine correction functions. When one makes the assumption that the model is not a significantly incorrect representation of the plant (Assumptions~\ref{ass:6_1_Erreur_fp_f_bornee} and \ref{ass:6_2_Meme_Courbusre}), these global effects of local data can globally reduce the predictive quality of the model as illustrated in Figure~\ref{fig:6_3_LocalVsGlobalCorrection}. The use of non-linear correction functions can avoid this issue. 
	\end{Remark}

	\begin{Remark} \label{rem:6_5_Particularie_ii}
		\textbf{(Concerning the particularity (ii))} 
		The ability to simultaneously provide an estimate of the expectation and uncertainty of the position of $f_p$ for a given $\bm{x}$ enables two things that were not previously considered:
		\begin{itemize} 
			\item \textit{To quantify the risks of a decision.} Indeed, knowing the uncertainty on $f_p$ for a given  $\bm{x}$ can be used to evaluate the uncertainty on the violation of the plant constraints $\bm{g}(\bm{x})$. Then, real benefit-risk analyses for each decisions can be made. This new possibility allows us to improve the experiment selector (see section~\ref{sec:6_7_Selecteur_d_exprience}).
			\item  \textit{To quantify the surprise associated with the observation of an experimental result.} 
			Once the experiment has been performed, one can  evaluate the performance of the updated model by comparing the confidence domains of the \textit{expected} (by the model) and \textit{actual} (taken on the plant) observations.  For example, if the experimental result is in a domain that the model considered very unlikely, then one can legitimately wonder if the model is actually up to date. Therefore, the observation of unlikely results can be used to detect CTP variations, and thus improve the
			consistency monitor (see section~\ref{sec:6_6_Contrôoleur_de_coherence}). 
		\end{itemize} 
	\end{Remark}
	
	Before moving on to improving the experiment selector and the consistency monitor, note that to compute $\mu_k^f(\bm{x})$ with \eqref{eq:6_13_def_eta_kf}  one needs to invert the matrix $\bm{\CC}$ which is potentially very large since it represents the whole history of the results of experiments. It is known that inverting a large matrix can be difficult or even impossible from a computational point of view, and it is therefore necessary to guarantee that this matrix is not ``too big'' (compared to the available computational capacities). One way to reduce its size, \textit{i.e. to compress the information it contains}, would be to replace the results of experiments ``close to each other in the input space'' by local averages and trends.  In the following section an explanation is given of how such a compression procedure can be carried out, and an improved version of the compressor is proposed.
		
\section{The improved compressor}
\label{sec:6_5_Compresseur}
	
	\begin{figure}[h]
		\centering
		\includegraphics[width=11cm] {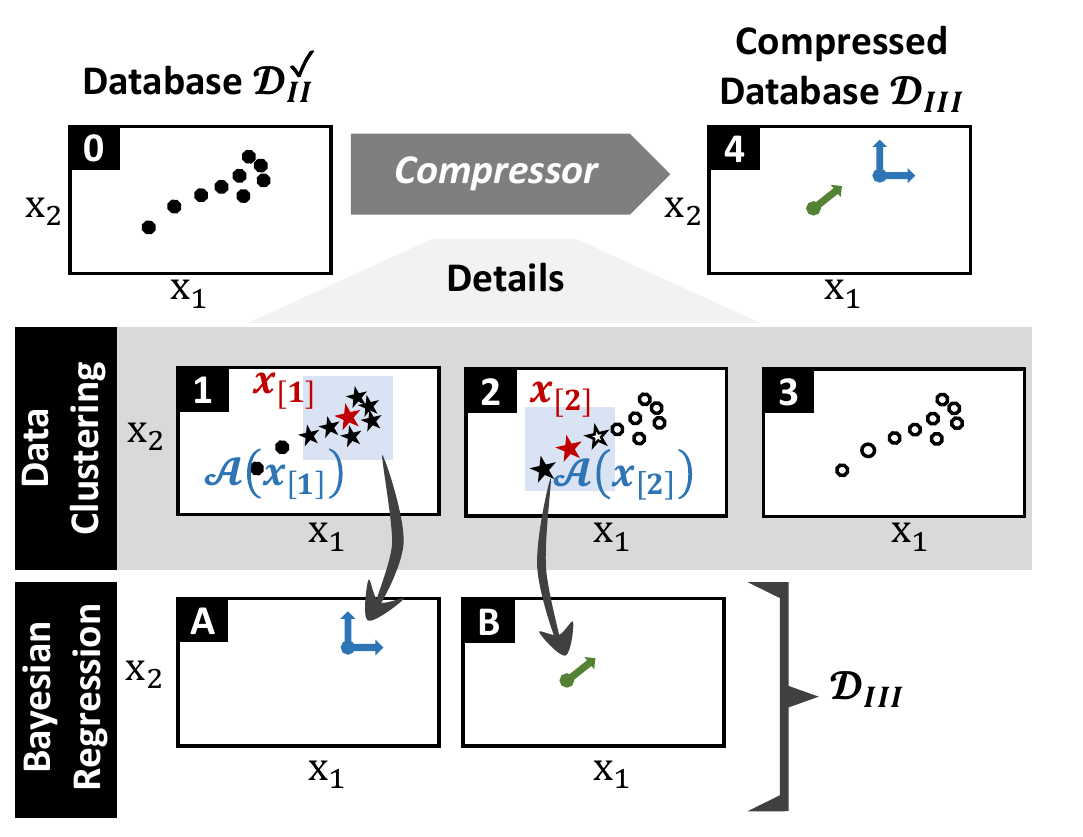}
		\caption{Graphical explanation of the functioning of the improved compressor.}
		\label{fig:6_5_ClusteringRawDataset}
	\end{figure} 
	The main ideas at the heart of the functioning of the new improved compressor are illustrated in Figure~\ref{fig:6_5_ClusteringRawDataset}. It can be seen that the compression mechanism is divided into two actions:  
	\vspace{-\topsep}
	\begin{itemize}[noitemsep]
		\item a data clustering action, 
		\item an action of transforming a group of data into a mean and \textit{some} trends, i.e. some directional derivatives. 
	\end{itemize}
	Let's start by explaining the data clustering strategy:
	
	\subsection{Data clustering -- a simple strategy} 
	
	 The idea is to extract from $\mathcal{D}_{II,k}^{\checkmark}$ the input sequence:
	\begin{equation*}
		X_{p,k} = \{\widehat{\bm{x}}_{p,0}, \hdots, \widehat{\bm{x}}_{p,k}\},
	\end{equation*}
	of which a copy 	$X^{\prime} = X_{p,k}$ is made which is called \textit{set of candidates}. Then, one randomly selects a point in this set of candidates that is called $\bm{x}_{[1]}$, and one calculates which points of $X_{p,k}$ are in the rectangular domain  $\mathcal{A}(\bm{x}_{[1]})$ around $\bm{x}_{[1]}$ in which it is assumed that $f_p$ has a quasi-affine behavior.  For the sake of simplicity, it is chosen to always define $\mathcal{A}(\bm{x}_{[1]})$  as a rectangle centered on $\bm{x}_{[1]}$, whose ``main directions'' are those of the axes of the input space, and which is defined as:
	\begin{align} \label{eq:6_14_definition_A_x}
		\mathcal{A}(\bm{x}_{[1]}) := \ & \{ \bm{x}\in\mathbb{R}^{n_x} | -\bm{a}/2 \leq \bm{x}-\bm{x}_{[1]} \leq \bm{a}/2 \}, &
		\bm{a} :=  \ & (a_1, ..., a_{n_x})^{\rm T},
	\end{align} 
	where all the $a_i$ are predefined (i.e. independent of the location of the $\bm{x}_{[1]}$) and defines the size of the rectangle in the direction of the $i$-th axis  of the input space. The appropriate choice of $a_i$, $\forall i\in\{1,..,n_x\}$, is explained in Remark~\ref{rem:6_3_Choix_ai}. All the points of  $X_{p,k}$ that belong to $\mathcal{A}(\bm{x}_{[1]})$ constitute the cluster $X_{[1]}:=\{\bm{x}\in X_{p,k} | \bm{x}\in \mathcal{A}(\bm{x}_{[1]})\}$. In Figure~\ref{fig:6_5_ClusteringRawDataset} - \textcolor{white}{\colorbox{black}{\textbf{\textbf{1}}}} one can find 
	\vspace{-\topsep}
	\begin{itemize}[noitemsep]
		\item $\bm{x}_{[1]}$ represented with a red star;
		\item $\mathcal{A}(\bm{x}_{[1]})$ represented with the blue rectangle;
		\item $X_{[1]}$ represented with black or red stars;
		\item $X^{\prime}$ represented with filled stars or filled dots (all dots are filled since it is the first step $X^{\prime}:= X_{p,k}$);
	\end{itemize}
	Then, the set of point in  $X^\prime$ that are in  $X_{[1]}$ are removed from $X^\prime$. The idea here is to maintain in $X^\prime$ only points that do not belong to any cluster. Then, this process is repeated by selecting   $\bm{x}_{[2]}$ randomly in $X^{\prime}$ to then build a new cluster $X_{[2]}$, etc. In  Figure~\ref{fig:6_5_ClusteringRawDataset} - \textcolor{white}{\colorbox{black}{\textbf{\textbf{2}}}} one can find: 
	\vspace{-\topsep}
	\begin{itemize}[noitemsep]
		\item $\bm{x}_{[2]}$ represented with a red star;
		\item $\mathcal{A}(\bm{x}_{[2]})$ represented with the blue rectangle;
		\item $X_{[2]}$ represented with black or red stars (filled or unfilled);
		\item $X^{\prime}$ represented with filled stars or filled dots ;
	\end{itemize}
	By repeating this procedure until the set of candidates $X^{\prime}$ is completely emptied, one can place all the data in groups from which one can extract averages and trends.

	\begin{Remark} \label{rem:6_3_Choix_ai}
		\textbf{(How to choose $\bm{a_i}$)} By definition the domain $\mathcal{A}(\bm{x})$ is a domain containing $\bm{x}$ and within which $f_p$ is quasi-affine. So, if one chooses to fix the shape of $\mathcal{A}(\bm{x})$ to that of the rectangle \eqref{eq:6_14_definition_A_x} an easy choice of $a_i$ would be to take it as small as possible. Indeed, by doing this, it is certain that the assumption that  $f_p$ is quasi-affine in $\mathcal{A}(\bm{x})$ irrespective of $\bm{x}\in\amsmathbb{R}^{n_x}$ is true. However, the negative effect of such a choice is to increase the number of clusters and thus to indirectly reduce the compression effect. There is therefore a trade-off between (i) ensuring the validity of the assumption that $f_p$  is quasi-affine in $\mathcal{A}(\bm{x})$, and (ii) achieving an effective compression. To find the right balance between these two points one proposes to start  from Assumptions~\ref{ass:6_1_Erreur_fp_f_bornee} and \ref{ass:6_2_Meme_Courbusre} which can be summarized with equation~\eqref{eq:6_7_GP_Fp}.  Given these assumptions, it is possible to estimate the relative approximation error (given by $error$ in [\%]) associated with the affine approximation of $f_p$ over a domain w.r.t. its size. This study is done in Appendix~\ref{App:AffineLikeDomain} where one finds that the relative error associated to an affine approximation of $f_p$ in a direction $\bm{\v}_i$  and at a distance $a_i$ is most likely bounded by:
		\begin{align} \label{eq:6_15_QuantifiedLoss}
			& error \leq  100 \frac{3\sqrt{3\pi}}{4}  \left( \frac{a_1}{\ell_{\bm{\v}_i}} \right)^2,
		\end{align} 
		\begin{align} 
			\text{where} \qquad &  \ell_{\bm{\v}_i}^{-2} :=  \bm{\v}_i^{\rm T} \bm{\L} \bm{\v}_i,
			\qquad  
			\bm{\v}_i := 
			\left(\begin{array}{c}
				\delta_{i1} \\
				\vdots \\
				\delta_{in_x}
			\end{array}\right),
			\quad \text{and} \quad 
			\delta_{ij} := \left\{
			\begin{array}{ll}
				1, & \text{if } i=j, \\
				0, & \text{otherwise}.
			\end{array}
			\right.
			. \nonumber
		\end{align} 
	From this inequality, one can compute the upper bound on the relative error of approximation of $f_p$ when linearized over a domain of size $a_i$ in the $\bm{\v}_i$ direction,  as a function of the values of  $a_i$, $\forall i=1,..,n_x$. Table~\ref{tab:6_1_Loss_wrt_ai} gives this relative error for 4 choices of $a_i$.  \textit{This error should be interpreted as the likely maximal percentage of information lost during the compression of the database $\mathcal{D}_{II}^{\checkmark}$, i.e. not as wrong information added to the database.} According to the case studies (see  ~\ref{sec:6_11_Cas_Etude_WOplant}),  a choice that seems appropriate would be  $a_i = \ell_{\bm{\v}_i}/6$. 
		
		\begin{table}[H]
			\centering
			\begin{tabular}{lccccc}
				\toprule
				$a_i$ & = & $\ell_{\bm{\v}_i}/2$   & $\ell_{\bm{\v}_i}/4$ & $\ell_{\bm{\v}_i}/6$ & $\ell_{\bm{\v}_i}/8$  \\
				\midrule
				$error$  &$\leq$     & $57.6\%$    & $14.4\%$    & $7.4\%$ & $3.6\%$ \\
				\bottomrule
			\end{tabular}%
			\captionsetup{justification=centering,width=.75\linewidth}
			\caption{Relative approximation error w.r.t. $a_i$.}
			\label{tab:6_1_Loss_wrt_ai}%
		\end{table}%
	\end{Remark}
	
\begin{figure}[h]
		\centering
		\includegraphics[width=10cm] {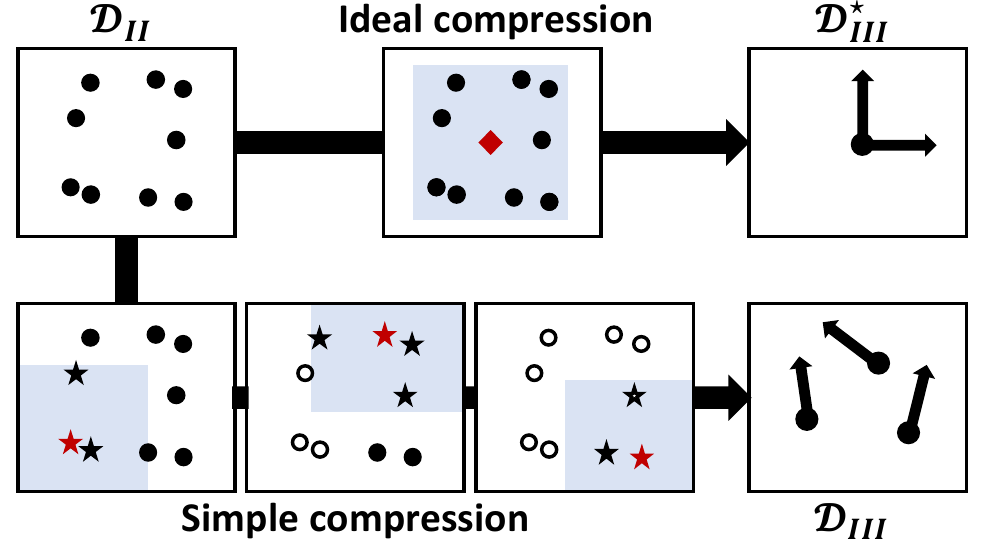}
		\caption{\textbf{Remark~\ref{rem:6_4_clustering performant}:}
			Schematic comparison of a simple and an ``ideal'' clustering method.}
		\label{fig:6_6_Ideal_Clustering}
	\end{figure}

	\begin{Remark} \label{rem:6_4_clustering performant}
		\textbf{(A more efficient clustering strategy)}
		If the idea of randomly selecting points $\bm{x}_{[i]}$ in the set of candidates $X^{\prime}$ is to guarantee that all the data of the $\mathcal{D}_{II}^{\checkmark}$ are allocated to clusters in a finite number of iterations of the clustering algorithm. It does not guarantee an optimal grouping of the data, i.e. the final number of groups will not necessarily be the smallest possible. A potential improvement of what is proposed would be to use an algorithm inspired by K-coverage\footnote{ These algorithms are for the moment mainly used in research on the placement and optimization of communication networks, see for example \cite{chakrabarty:2002,wang:2003}.} algorithms consisting of optimizing the placement of the domains $\mathcal{A}$  in the input space in order to cover all data as efficiently as possible.
		Figure~\ref{fig:6_6_Ideal_Clustering} illustrates this point for a relatively simple case. One can see that with a smart placing of the center of  $\mathcal{A}$ (the red diamond), all the data can be gathered within a unique cluster. On the other hand, the simple approach leads to three clusters, and therefore to a less efficient compression. 
	\end{Remark}
	
	\subsection{Bayesian regression -- Transform a group of points into an average and trends} 
	
	 For the sake of simplicity, the way to implement a Bayesian regression is explained only for the cluster $X_{[1]}$ -- the other clusters are treated in exactly the same way. The idea is to transform the $m$ points of $X_{[1]}$: $(\widehat{\bm{x}}_{[1,1]},\widehat{y}_{p[1,1]},\sigma_{xy[1,1]}^2), ... , (\widehat{\bm{x}}_{[1,m]},\widehat{y}_{p[1,m]},$ $\sigma_{xy[1,m]}^2)$, into a mean, trends, and the uncertainties on them. Since Assumptions~\ref{ass:6_1_Erreur_fp_f_bornee} and \ref{ass:6_2_Meme_Courbusre} are assumed to hold, and since they can be summarized with \eqref{eq:6_7_GP_Fp}, these mean, trends, and uncertainties can be found by simply evaluating the value and gradient of a GP regression based on those points and evaluated at $\bm{x}_{[1]}$, as shown in Appendix~\ref{App:Predict_Val_and_Grad}. Basically, if one considers equations  \eqref{eq:C_App_CovMat_F_DF_1} and \eqref{eq:C_App_CovMat_F_DF} with the matrices $\bm{\CC}_{**}$, $\bm{\CC}_{*}$, and $\bm{\CC}$ appropriately defined w.r.t. the cloud of points of $X_{[1]}$ and definitions 	\eqref{eq:C_8_verjrhcjdsx}, \eqref{eq:A_iywbvwe}, then one gets: 
	 \begin{align}
		 \left(\begin{array}{c}
		 	\alpha \\  \bm{\beta}
		 \end{array}
		 \right) :=
	 	\amsmathbb{E}
	 	\left[\left.\left(\begin{array}{c}
	 		f_{p} \\  \nabla_{\bm{x}} f_{p}
	 	\end{array}
	 	\right)\right| X_{[1]} \right]= \ &
	 	\text{\eqref{eq:C_App_CovMat_F_DF_1}}|X_{[1]} \\
	 	\left(\begin{array}{cc}
	 		\sigma_{\alpha}^2 & \cdot \\ \cdot & \bm{\Sigma}_{\bm{\beta}}
	 	\end{array}
	 	\right) := \amsmathbb{V}
	 	\left[\left.\left(\begin{array}{c}
	 		f_{p} \\  \nabla_{\bm{x}} f_{p}
	 	\end{array}
	 	\right)\right| X_{[1]} \right] = \ & \text{\eqref{eq:C_App_CovMat_F_DF}}|X_{[1]},
	 \end{align}
	 
	 Now that one has converted a cloud of data into a mean and trends, let's extract the subpart of these new data that actually contain useful information, i.e. \textit{the data of interest}.  
	 
	\subsection{Selection of the data of interest}
	
	One considers the average point of cluster,
		\begin{align}
			& 	(\overline{\bm{x}}_{[1]}, \overline{y}_{[1]},\overline{\sigma}_{[1]}^2), & 
			\text{where} \qquad  & \overline{y}_{[1]} = \overline{\alpha}, \quad  \text{and}  \quad  \overline{\sigma}_{[1]}^2 := \sigma_{\alpha}^2
		\end{align}
	is in any case, an interesting data to record in $\mathcal{D}_{III}$. In fact the only case where this data could be of no interest is when the measurement uncertainties $\sigma_{xy}^2$ are so large that the uncertainty of the location $f_p$ w.r.t.  $f$  (which has been called $\sigma_f^2$, see Assumption~\ref{ass:6_1_Erreur_fp_f_bornee})  is significantly smaller than $\sigma_{\alpha}^2$. In this case, as the union of measurements would be more uncertain than the model itself ($\sigma_{\alpha}^2>\sigma_f^2$), it is clear that it is of no interest for decision making.\footnote{If this occurs , then the  idea of using measurements to guide the decision is inappropriate and RTO is probably not a solution to your problem. You should rather consider robust optimization methods that base decision making mainly on the model, see \cite{Srinivasan:03b,Monnigmann:2003}.}
	
	Concerning the gradient $\bm{\beta}$,  one decides to preserve only the components (the directional derivatives) which present ``little uncertainty''. To decide if  $\bm{\beta}$ presents an interest w.r.t. a direction $\bm{\v}$, one suggest to compare its \textit{a priori} and \textit{a posteriori} distributions, i.e. if 
  	\begin{align} \label{eq:6_20_DirectionOfInterst}
  		\gamma \bm{\v}^T \bm{\Sigma}_{\beta} \bm{\v} \leq  \ &\sigma_f^2 \bm{\v}^T \bm{\L} \bm{\v}, & 
  		\|\bm{\v}\|_2 = \ & 1,
  	\end{align} 
 		then the directional derivative $\bm{\v}^{\rm T}\bm{\beta}$ is relevant to decision making. The variable $\gamma \geq 1$ is a scalar to be chosen that quantifies the ratio $\left(\sigma_f^2 \bm{\v}^T \bm{\L} \bm{\v}\right)/\left(
  		\bm{\v}^T \bm{\Sigma}_{\beta} \bm{\v} \right)$from which the \textit{a priori} variance is considered as significantly greater than the \textit{a posteriori} one. Thus, the larger $\gamma$ is, the more precise a directional derivative must be to be considered interesting. A simple manipulation of  equation~\eqref{eq:6_20_DirectionOfInterst}allows to rewrite it as: 
  		\begin{equation}
  			\bm{\v}^T \left(
  			\sigma_f^2 \bm{\L}- \gamma\bm{\Sigma}_{\beta}
  			\right)\bm{\v} \geq 0.
  		\end{equation}
  		One can then evaluate \eqref{eq:6_20_DirectionOfInterst} in all directions $\bm{\v}$ simultaneously. To do this, it is sufficient to perform a spectral analysis of $\sigma_f^2 \bm{\L}- \gamma\bm{\Sigma}_{\beta}$. All eigenvectors associated with \textit{positive} eigenvalues are the directions to be kept for the decision making.  Let's assume that there is $n_{\v}\leq n_x$ eigenvectors associated with positive eigenvalues and let's name them  $\{\bm{v}_{[1,1]},...,\bm{v}_{[1,n_{\v}]}\}$. In this case the directional derivatives to be saved in $\mathcal{D}_{III}$ are:
  		\begin{align}
  			& (\overline{\bm{x}}_{[1]}, \nabla_{\bm{\v}} y_{[1,1]},\bm{\v}_{[1,1]}, \sigma^2_{\v[1,1]}), \ ..., \ (\overline{\bm{x}}_{[1]}, \nabla_{\bm{\v}} y_{[1,n_{\v}]},\bm{\v}_{[1,n_{\v}]}, \sigma^2_{\v[1,n_{\v}]}), \nonumber \\ 
  			\text{where} \quad  & \forall i=1,...,n_{\v}: 
  			\qquad 	
  			\nabla_{\bm{\v}} y_{[1,i]} :=   \bm{\v}_{[1,n_{\v}]}^{\rm T}\overline{\bm{\beta}}, 
  			\qquad 
  			\sigma^2_{\v[1,i]} :=  \bm{\v}_{[1,i]}^{\rm T} \bm{\Sigma}_{\beta} \bm{\v}_{[1,i]}. 
		\end{align}
		
		\subsection{The compressed database}
  		
  		After having done this compression work on each cluster, one ends up with the data of interest of all the clusters that must be saved in $\mathcal{D}_{III}$. This saving can be done by simply stacking the data:  
  		\begin{align}
  			\bm{X}  := \ & (\overline{\bm{x}}_{[1]}, \overline{\bm{x}}_{[1]}, ...,\overline{\bm{x}}_{[1]}, \quad \overline{\bm{x}}_{[2]}, ...)^{\rm T}, &
  			\bm{y}_p := \ & (\overline{y}_{[1]}, \nabla_{\bm{\v}} y_{[1,1]}, ..., \nabla_{\bm{\v}} y_{[1,n_v]}, \quad\overline{y}_{[2]},...
  			)^{\rm T}, \nonumber \\
  			\bm{\V} := \ & (\emptyset, \bm{\v}_{[1,1]}, ..., \bm{\v}_{[1,n_v]},\quad \emptyset, ...)^{\rm T}, & 
  			\bm{\sigma} :=  \ & (\overline{\sigma}_{[1]}^2, \sigma^2_{\v[1,1]}, ..., \sigma^2_{\v[1,n_v]}, \quad \overline{\sigma}_{[2]}^2)^{\rm T}.
  		\end{align}
  		It is then quite simple to make the link between $\mathcal{D}_{III} := \{\bm{X},\bm{y}_p,\bm{\V},\bm{\sigma}\}$ and the equivalent data matrices (using lighter notations) used in Appendix~\ref{App:C_PartialDeriv} \eqref{eq:C_6_Matrices_de_donnees} that are: 
	  	\begin{align}
	  		\bm{X}  := \ & (\bm{x}_1, \bm{x}_2, ..., \bm{x}_n)^{\rm T}, &
  			\bm{y}_p := \ & (\widehat{y}_{p,1}, \widehat{\nabla_{\v_j} y}_{p,j}, ..., \widehat{y}_{p,n})^{\rm T}, \nonumber \\
  			\bm{\V} := \ & (\emptyset, \bm{\v}_2, ..., \emptyset)^{\rm T}, & 
  			\bm{\sigma} :=  \ & (\sigma_{xy,1}^2, \sigma_{\v,2}^2, ..., \sigma_{xy,n}^2)^{\rm T}.
  		\end{align}
  		
	\subsection{Implementation advice}
	\label{sec:6_5_5_Conseil_implementation}
  		
  		As the database $\mathcal{D}_{II}$ gets bigger at the end of each experiment, it is necessary to repeat the compression regularly. However, it is not necessary to repeat the entire procedure.   Instead, when a point $(\widehat{\bm{x}}_{{k+1}},\widehat{y}_{p,k+1},\widehat{\sigma}^2_{xy,k+1})$ is added to $\mathcal{D}_{II}$, one can first check if it belongs to an existing domain $\mathcal{A}(\bm{x}_{[i]})$, $\forall i$, resulting from a previous compression.
  		\begin{itemize} 
  			\item If it belongs to such a set, then this point could be integrated to the associated cluster. Then, only the bayesian regression of this cluster needs to be repeated, and two options are possible: (\textit{A -- if there are few points)} The full regression procedure is repeated. \textit{(B -- if there are many points)} 
  			One replaces the cloud of points of the cluster (minus the last point) by the value and the directional derivatives $\alpha$ and $\bm{\beta}$, and performs the regression with the new point and those aggregated values.  
  			\item If it does not belong to any cluster $\mathcal{A}(\bm{x}_{[i]})$, $\forall i$, then it is sufficient to continue the compression algorithm as if $\widehat{\bm{x}}_{{k+1}}$ was the last item in the set of candidates $X^{\prime}$.
		\end{itemize}
		With this strategy, it is believed that regardless of the size of the size of $\mathcal{D}_{II}^{\checkmark}$, the progressive addition of new data allows the compression calculations to be carried out in an equally progressive way.  This avoids a computational blocking associated with a regression or a clustering on too much data.
  		 
\section{The improved  consistency monitor}
\label{sec:6_6_Contrôoleur_de_coherence}

	\subsection{Data expiration factors}

	The role of the consistency monitor is to ensure that only data from $\mathcal{D}_{II}$ that are up to date are sent to the compressor. The subpart of $\mathcal{D}_{II}$ that passes through the consistency monitor  is named the validated database ($\mathcal{D}_{II}^{\checkmark}$). 
	
	It is considered that the only two factors that can lead to data expiration are:
	\begin{itemize}
		\item \textbf{Expiration factor 1:} Unmeasured and progressive disturbances, such as plant aging. 
		\item \textbf{Expiration factor 2:} Unmeasured and instantaneous disturbances, such as a step in the composition of one of the processed raw products or an instantaneous degradation of a property of a plant component.  
	\end{itemize}

	\subsection{Managing the expiration factor 1}
	
		To manage the effects of unmeasured progressive disturbances, it is suggested to consider the time $t$ (i.e. the age of the plant) as a measured disturbance whose effect on the plant outputs is unknown. Since this disturbance is not a priori part of the model, one cannot use the model to identify its hyperparameter $\ell_t$  for each output of the plant. Therefore, a different attention needs to be given to this new input. One possibility would be to let the engineers set $\ell_t$. In fact, they should be able to assess the order of magnitude of the time it takes for the progressive unmeasured disturbances to have significant effects on the plant output. This order of magnitude is the hyperparameter $\ell_t$ to be used.    Another possibility would be to select the $\ell_t$ that maximizes the likelihood of the measurements already made via a model selection strategy (i.e. auto-tuning of the  hyperparameters, see Chapter 5 of \cite{Rasmussen:2006}) specifically focused on the choice of this hyperparameter with constraints on its order of magnitude so that it does not capture unmeasured high frequency perturbations (not progressive). 
	
\begin{figure}[h]
	 	\centering
	 	\includegraphics[width=10cm]{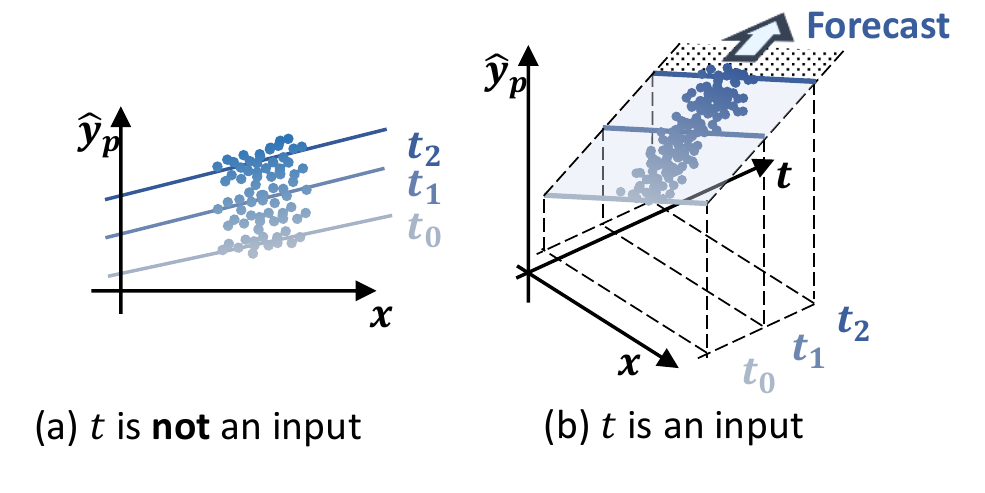}
	 	\caption{To manage the expiration factor 1 , one can consider the time $t$ as a measured disturbance}
	 	\label{fig:6_7_TimeIsKey}
	 \end{figure}

	As illustrated in Figure~\ref{fig:6_7_TimeIsKey}, the idea of considering $t$ as a measured disturbance could enable the identification of trends in temporal variations of  $f_p$. Therefore, data too old to represent the current behavior of the plant could still be used to improve current (and future\footnote{Which could be the basis for a feedback iteration between the RTO and the planning/scheduling layers of the management of a plant.}) predictions through this time trend of $f_p$.  
		 
	\subsection{Managing the expiration factor 2}
	
\begin{figure}[h]
		\centering
		\includegraphics[width=14cm]{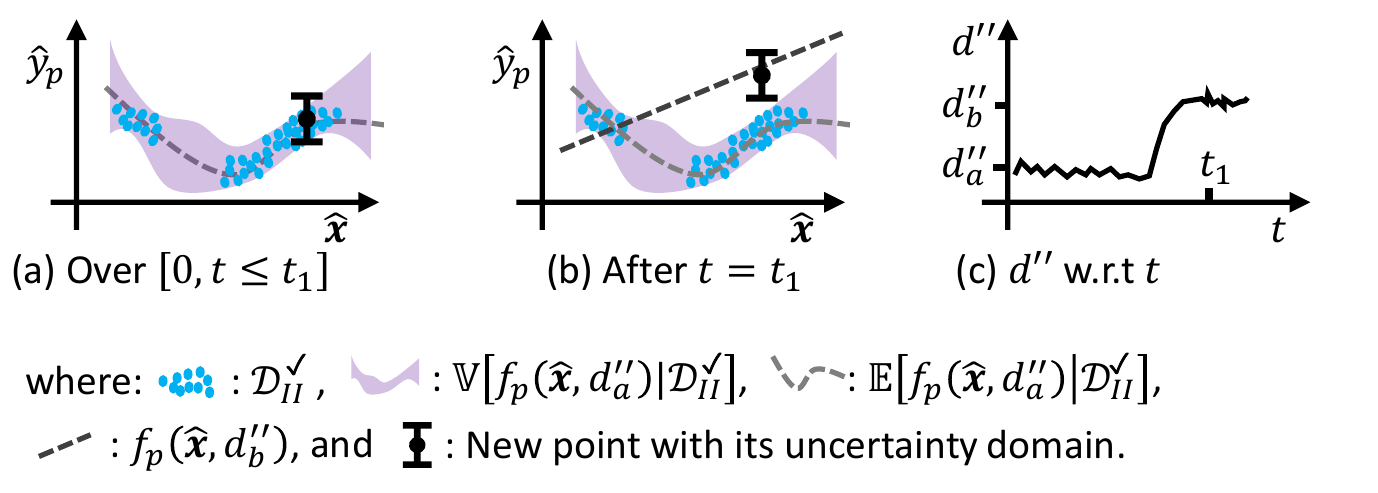}
		\caption{Graphic illustration of the proposed CTP transition detector}
		\label{fig:6_8_ManageUnmeasuredDisturbance}
	\end{figure}

	An instantaneous unmeasured disturbance has the effect of instantaneously changing the function $f_p$.  In other words, such a disturbance has the effect of ending a CTP (see Definition~\ref{def:6_PTC}). To manage such a situation one proposes to empty the databases $\mathcal{D}_I$ and $\mathcal{D}_{II}$ when a CTP transition is detected. Two cases are possible: 
	\begin{itemize}
		\item \textbf{Case 1 -- the simple case:} When the plant is at steady state and its inputs  $\bm{u}$  and $\bm{d}^{\prime}$ are not perturbed.  If in such a case there is a sudden change of an unmeasured disturbance $\bm{d}^{\prime\prime}$ which has significant effects on the output $\widehat{y}_{p}$. Then it can be detected with a classical fault detection method (FDM) as proposed in \cite{Ye:2018,Mukkula:2019}. 
		\item \textbf{Case 2 -- the complex case:} When the plant inputs are being manipulated, e.g. the optimum plant is being sought.  In this case, it is less easy to dissociate the effects of input manipulations from those of variations of unmeasured disturbances. To dissociate these two effects, one proposes to quantify the surprise associated to the observation of the results of a new experiment. Indeed, if this observation is \textit{too far} from the others, then one could suspect that this difference is due to a change of $f_p$. To determine whether an observation $(\widehat{\bm{x}}_{k+1}, \widehat{y}_{p,k+1}, \sigma_{xy,k+1}^2)$ is \textit{too far}  from the previous observations,  one proposes to compare it to the domain of expectation of this observation based on the previous observations:
		\begin{align*}
			\widetilde{y}_{p,k+1} := \ & \amsmathbb{E}[f_p(\widehat{\bm{x}}_{k+1})|\widehat{\bm{x}}_{k+1},\mathcal{D}_{III}], &
			\widetilde{\sigma}_{k+1}^2 := \ & \amsmathbb{V}[f_p(\widehat{\bm{x}}_{k+1})|\widehat{\bm{x}}_{k+1},\mathcal{D}_{III}].
		\end{align*}
		If the intersection of the two segments is such that 
		\begin{equation} 
			[\widehat{y}_{p,k+1}-2\sigma_{k+1}, \  \widehat{y}_{p,k+1}+2\sigma_{k+1}] \cap  [\widetilde{y}_{p,k+1}-2\widetilde{\sigma}_{k+1}, \ 
			\widetilde{y}_{p,k+1}+2\widetilde{\sigma}_{k+1}]=\emptyset, \label{eq:Condition_end_CPT}
		\end{equation} 
		then one considers that the new observation is surprising. If this is the case, then one supposes that a CTP transition has occurred and one reacts accordingly, i.e. one empties the databases $\mathcal{D}_{I}$ and $\mathcal{D}_{II}$. \textit{Note that in this case the results of the last experiment can be kept since they belong to the new CTP.}
		
Finally, Figure~\ref{fig:6_8_ManageUnmeasuredDisturbance} gives a graphical interpretation of this strategy. One can see, among other things, the condition of intersection of the two segments which is represented by the overlap of the purple domain with the black segment, and the fact that the purple domain is influenced by the results of previous experiments (the blue points).
	\end{itemize}

At this stage it has been defined (i) how to combine the data and the nominal model, (ii) how to compress the database to make it usable, and (iii) how to keep the database up to date. So the whole data processing part is over. Now let's see how the decisions should be made, starting with the experiment selector.
	
\section{The improved experiment selector}
\label{sec:6_7_Selecteur_d_exprience}

\subsection{Definition of the point to go to}

The role of the experiment selector is to define the operating point around which the next experiments should be performed. As noted in the Remark~\ref{rem:6_4_Particularite_i_part_2}, the corrected model provides estimates of $\amsmathbb{E}[f_p(\bm{x})]$ and of $\amsmathbb{V}[f_p(\bm{x})]$ that can be used to quantify the risks associated with a decision.   Clearly, what one is looking for is $\bm{u}_p^{\star}$ such that: 
\begin{equation}
	\bm{u}_p^{\star} :=  \operatorname{arg}
	\underset{\bm{u}}{\operatorname{min}} \quad  \phi(\bm{u},\bm{f}_p(\bm{x}_p)), \quad 
	\text{s.t.} \quad  
	\bm{g}(\bm{u},\bm{f}_p(\bm{x}_p)) \leq \bm{0},
	\label{eq:6_28_PB_Opt_usine} 
\end{equation}
and one wants to reach it with a series of $\bm{u}_k$ that does not violate the plant constraints. Therefore, the following approach is proposed:
\begin{equation}
	\bm{u}_{k+1}:=  \operatorname{arg}
	\underset{\bm{u}}{\operatorname{min}} \quad \varphi(\bm{u},\bm{f}_k(\widehat{\bm{x}})),\quad 
	\text{s.t.} \quad 
	\bm{h}(\bm{u},\bm{f}_k(\widehat{\bm{x}})) \leq \bm{0}, \quad \widehat{\bm{x}} := [\bm{u}^{\rm T}, \bm{d}^{\prime \rm T}]^{\rm T},
	\label{eq:6_29_Model_based_PB} 
\end{equation}
where $\varphi$ is the expectancy of the plant's cost function, and $h_{(i)}$ is the estimated upper bound of the plant's constraint function $g_{(i)}$, $\forall $$i$$=$$ 1,$$...,$$n_g$:
\begin{align}
	\varphi(\bm{u},\bm{f}_k(\widehat{\bm{x}})) := \ & \amsmathbb{E}[\phi(\bm{u},\bm{f}_k(\widehat{\bm{x}}))], \label{eq:6_29_varphi} \\
	h_{(i)}(\bm{u},\bm{f}_k(\widehat{\bm{x}})) := \ & \amsmathbb{E}[g_{(i)}(\bm{u},\bm{f}_k(\widehat{\bm{x}}))] + 2 \sqrt{\amsmathbb{V}[g_{(i)}(\bm{u},\bm{f}_k(\widehat{\bm{x}}))]}. \label{eq:6_30_h}
\end{align}

Using the expected upper bound of the constraint instead of its expectation manages the risk associated to an experiment. \textit{(In fact, using this expected upper bound is only useful on constraints that must be satisfied, often called hard constraints. For the constraints that can be temporarily violated without much consequence (usually named  soft constraints),  one does not recommend the use of such an upper bound; rather the use of the expectancy. ) }

Definitions \eqref{eq:6_29_varphi} and \eqref{eq:6_30_h} make it possible to ensure that the solution of	\eqref{eq:6_29_Model_based_PB} is a relevant point to test, but evaluating them can be complicated. Therefore, in the  next subsection a  method (which works well when the uncertainty on $\bm{f}_k$ is small) is proposed to evaluate them.

\subsection{Evaluation method of functions \eqref{eq:6_29_varphi} and \eqref{eq:6_30_h}}
\label{sec:6_6_2_Uncertaintyspreads}

The objective of this section is to simplify expressions  \eqref{eq:6_29_varphi} and \eqref{eq:6_30_h} to make them easily solvable. To do this, let's start writing the Taylor series of $\phi$ and $g_{(i)}$, $\forall i=1,..,n_g$, centered on the evaluation point $(\bm{u}_*,\bm{y}_*=\amsmathbb{E}[\bm{f}_k(\bm{x}_*)])$: 
\begin{align}
	\phi(\bm{u},\bm{y}) = \ &  \phi|_{\bm{u}_*,\bm{y}_*} + 
	\sum_{j=1}^{n_u} (u_{(j)} - u_{*(j)}) \partial_{u_{(j)}} \phi|_{\bm{u}_*,\bm{y}_*} + 
	\sum_{j=1}^{n_y} (y_{(j)} - y_{*(j)}) \partial_{y_{(j)}}\phi|_{\bm{u}_*,\bm{y}_*} + ...  \nonumber \\
	& \mathcal{O}(\|\bm{y} - \bm{y}_{*}\|^2) + \mathcal{O}(\|\bm{u} - \bm{u}_{*}\|^2),  \nonumber \\
	g_{(i)}(\bm{u},\bm{y}) = \ &  g_{(i)}|_{\bm{u}_*,\bm{y}_*} + 
	\sum_{j=1}^{n_u} (u_{(j)} - u_{*(j)}) \partial_{u_{(j)}} g_{(i)}|_{\bm{u}_*,\bm{y}_*} + ...  \nonumber \\
	& \sum_{j=1}^{n_y} (y_{(j)} - y_{*(j)}) \partial_{y_{(j)}}g_{(i)}|_{\bm{u}_*,\bm{y}_*} + \mathcal{O}(\|\bm{y} - \bm{y}_{*}\|^2) + \mathcal{O}(\|\bm{u} - \bm{u}_{*}\|^2). \nonumber 
\end{align}
The effect of uncertainties on the outputs $\bm{y}_*$ can be captured by evaluating these functions at the point $(\bm{u}_*,\bm{y}_*+\bm{\epsilon}_y)$ with  $\epsilon_{y(j)} \sim\mathcal{N}(0,\amsmathbb{V}[f_{k(j)}(\bm{x}_*)])$, $\forall j=1,...,n_y$:
\begin{align}
	\phi(\bm{u}_*,\bm{y}_*+\bm{\epsilon}_y) = \ &  \phi|_{\bm{u}_*,\bm{y}_*} + 
	\sum_{j=1}^{n_y} \epsilon_{y(j)} \partial_{y_{(j)}}\phi|_{\bm{u}_*,\bm{y}_*} +  \mathcal{O}(\|\bm{\epsilon}_y\|^2),  \label{eq:6_31_yfuuycewrds}\\
	g_{(i)}(\bm{u}_*,\bm{y}_*+\bm{\epsilon}_y) = \ &  g_{(i)}|_{\bm{u}_*,\bm{y}_*} + 
	\sum_{j=1}^{n_y} \epsilon_{y(j)} \partial_{y_{(j)}}g_{(i)}|_{\bm{u}_*,\bm{y}_*} +  \mathcal{O}(\|\bm{\epsilon}_y\|^2), \label{eq:6_32_yfuvcjwds}
\end{align}
If one supposes that:
\begin{Assumption} \label{ass:6_3_Petites_variances}
	$\forall j = 1,...,n_y$ and $\forall \bm{x}_*\in\amsmathbb{R}^{n_x}$ the variances $\amsmathbb{V}[f_{k(j)}](\bm{x}_*)$ are small enough so that the terms  $ \mathcal{O}(\|\bm{\epsilon}_y\|^2)$ of \eqref{eq:6_31_yfuuycewrds} and \eqref{eq:6_32_yfuvcjwds} are negligible. 
\end{Assumption}
\noindent
Then \eqref{eq:6_31_yfuuycewrds} and \eqref{eq:6_32_yfuvcjwds} reduce to: 
\begin{align*}
	\phi(\bm{u}_*,\bm{y}_*+\bm{\epsilon}_y) = \ &  \phi|_{\bm{u}_*,\bm{y}_*} + 
	\sum_{j=1}^{n_y} \epsilon_{y(j)} \partial_{y_{(j)}}\phi|_{\bm{u}_*,\bm{y}_*},  \\
	g_{(i)}(\bm{u}_*,\bm{y}_*+\bm{\epsilon}_y) = \ &  g_{(i)}|_{\bm{u}_*,\bm{y}_*} + 
	\sum_{j=1}^{n_y} \epsilon_{y(j)} \partial_{y_{(j)}}g_{(i)}|_{\bm{u}_*,\bm{y}_*}.
\end{align*}
So, if one makes Assumption~\ref{ass:6_3_Petites_variances} the following equalities are obvious: 
\begin{align}
	\amsmathbb{E}\left[\phi(\bm{u},\bm{f}_k(\widehat{\bm{x}}))\right] = \ & \amsmathbb{E}\left[ \phi|_{\bm{u}_*,\bm{y}_*} + 
	\sum_{j=1}^{n_y} \epsilon_{y(j)} \nabla_{y_{(j)}}\phi|_{\bm{u}_*,\bm{y}_*} \right] &
	= \ & 	\phi|_{\bm{u}_*,\bm{y}_*},
	\label{eq:6_33_E_phi} \\
	\amsmathbb{E}\left[g_{(i)}(\bm{u},\bm{f}_k(\widehat{\bm{x}}))\right] = \ & 
	 \amsmathbb{E}\left[ g_{(i)}|_{\bm{u}_*,\bm{y}_*} + 
	\sum_{j=1}^{n_y} \epsilon_{y(j)} \partial_{y_{(j)}}g_{(i)}|_{\bm{u}_*,\bm{y}_*} \right] &
	= \ &  g_{(i)}|_{\bm{u}_*,\bm{y}_*}, 
	\label{eq:6_34_E_g} \\ 
	\amsmathbb{V}\left[ g_{(i)}(\bm{u},\bm{f}_k(\widehat{\bm{x}}))\right] = \ & 
	\amsmathbb{V}\left[
	g_{(i)}|_{\bm{u}_*,\bm{y}_*} + 
	\sum_{j=1}^{n_y} \epsilon_{y(j)} \partial_{y_{(j)}}g_{(i)}|_{\bm{u}_*,\bm{y}_*}
	\right], \nonumber  \\
	= \ &  \sum_{j=1}^{n_y} \big| \partial_{y_{(j)}}g_{(i)}|_{\bm{u}_*,\bm{y}_*} \big|^2  \amsmathbb{V}[ f_{k(j)}(\widehat{\bm{x}}_{*}) ]
	\label{eq:6_35_V_g}. 
\end{align}
One specifies that to obtain equation~\eqref{eq:6_35_V_g} it is assumed that the uncertainties on the outputs, $\bm{y}_p$, are independent since each subfunction, $f_{k(i)}$, is corrected individually. Hence the absence of covariance terms, $\amsmathbb{C}\big(f_{k(j)}(\widehat{\bm{x}}_{*}), f_{k(j^\prime)}(\widehat{\bm{x}}_{*})\big)$ with $j\neq j^\prime$. One deduces from \eqref{eq:6_33_E_phi}, \eqref{eq:6_34_E_g}, and \eqref{eq:6_35_V_g} the following estimates of \eqref{eq:6_29_varphi} and \eqref{eq:6_30_h}:
\begin{align}
	\varphi(\bm{u},\bm{f}_k(\widehat{\bm{x}})) \approx \ &
	\phi(\bm{u},\amsmathbb{E}[\bm{f}_{k}(\widehat{\bm{x}})]),
	 \\
	h_{(i)}(\bm{u},\bm{f}_k(\widehat{\bm{x}})) \approx \ &
	g_{(i)}(\bm{u},\amsmathbb{E}[\bm{f}_{k}(\widehat{\bm{x}})])
	 + 2 \sqrt{
	 \sum_{j=1}^{n_y} \big| \nabla_{y_{(j)}}g_{(i)}|_{\bm{u},\amsmathbb{E}(\bm{f}_{k}(\widehat{\bm{x}}))}\big|^2  \amsmathbb{V}[f_{k(j)}(\widehat{\bm{x}})]
 	}. 
\end{align}

\begin{figure}[h]
	\centering
	\includegraphics[width=14cm]{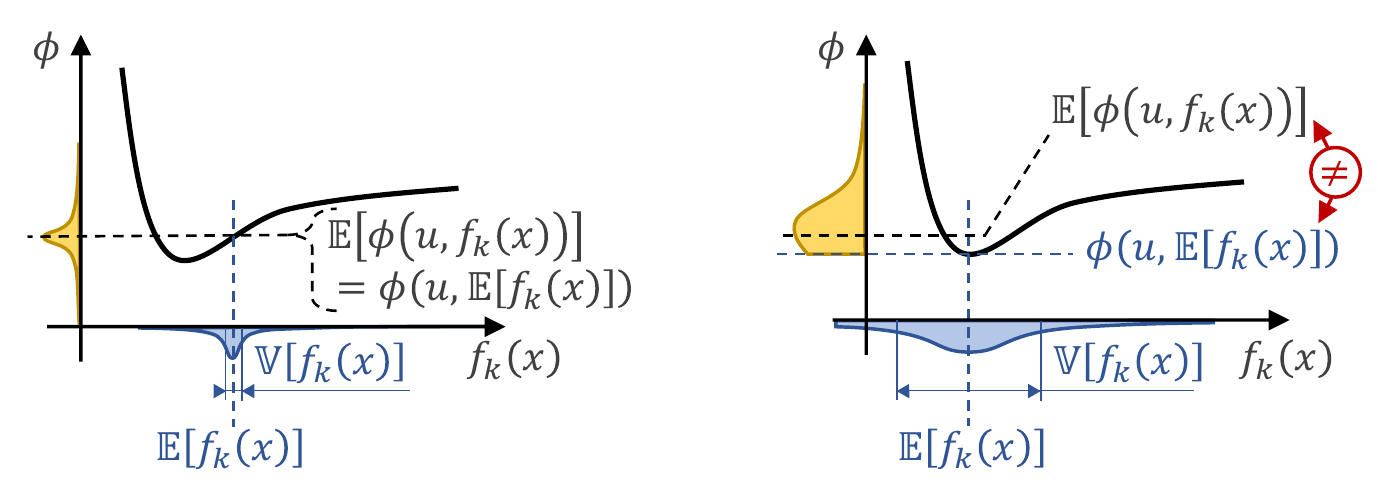}
	\begin{minipage}{7cm}
		\centering
		a) Small uncertainty on $f_k(\bm{x})$
	\end{minipage}
	\begin{minipage}{7cm}
		\centering
		b) Large uncertainty on $f_k(\bm{x})$
	\end{minipage}
	\caption{Illustration of the limit of Assumption~\ref{ass:6_3_Petites_variances}.}
	\label{fig:6_9_Limites_Hypothese_petite_variance_fk}
\end{figure}

\begin{Remark}
	\textbf{(Limit of Assumption~\ref{ass:6_3_Petites_variances})}  Assumption~\ref{ass:6_3_Petites_variances} is valid when the confidence domain of $f_{k}$ is small enough such that functions  $\phi$ and $g_{(i)}$ are  quasi-affine, see
	Figures~\ref{fig:6_9_Limites_Hypothese_petite_variance_fk} a) and b). In case  a), the confidence domain of $f_k$ is small and within it \eqref{eq:6_29_varphi} and \eqref{eq:6_30_h} can be approximated by \eqref{eq:6_31_yfuuycewrds} and \eqref{eq:6_32_yfuvcjwds}. In case b), this is clearly no longer true. To handle such cases, one should either increase the order of the Taylor series approximation, or use a Monte Carlo method, i.e. evaluate $\phi$ and $g_{(i)}$ $\forall i$ for a ``large number'' of inputs following the distribution: 
	\begin{equation*}
		\left(\begin{array}{c}
			\bm{u} \\ \bm{y}
		\end{array}\right) \sim 
	\mathcal{N}\left(
	\left( \begin{array}{c}
		\bm{u} \\ \amsmathbb{E}[f_{k(1)}(\widehat{\bm{x}})] \\ \vdots \\
		\amsmathbb{E}[f_{k(n_y)}(\widehat{\bm{x}})]
	\end{array} \right), \ 
	\left( \begin{array}{cccc}
		\bm{0} &                                           &        &  \\
		       & \amsmathbb{V}[f_{k(1)}(\widehat{\bm{x}})] &        &  \\
		       &                                           & \ddots &  \\
		       &                                           &        & \amsmathbb{V}[f_{k(n_y)}(\widehat{\bm{x}})]
	\end{array} \right) 
	\right)
	\end{equation*}
	which can be computationally expensive to perform at each iteration of the solver (loop 1 of Figure~\ref{fig:2___1_Solver_Model_Plant}), i.e. not each RTO iteration (loop 2 of Figure~\ref{fig:2___1_Solver_Model_Plant}).
\end{Remark}

Finally, thanks to the experiment designer, Assumption~\ref{ass:6_3_Petites_variances} can be forced to be true at the convergence point as explained in the following section. 

\section{The improved experiment designer}
\label{sec:6_7_The_improved_experiment_designer}

The role of the  experiment designer is to convert a target  $\bm{u}_{k+1}$ into a series of experiments to be conducted on the plant. For the S-ASP  it consisted in converting $\bm{u}_{k+1}$  into a series of at least $n_u+1$  experiments structured in a predefined way so that estimates of plant's values and gradients can easily be extracted from them.   

In the context of the ASP, there is no need for the experiments to follow a predefined structure since the compressor is advanced enough to accept ``unstructured'' data, see section~\ref{sec:6_5_Compresseur}.  

 If the collection of data does not need to be structured anymore, it is still interesting to force the  appearance of some particular structures when the autopilot converges. For example, to try to force the conditions of application of Theorem~\ref{thm:2___1_AffineCorrection_Impl_KKTmathcing} to always be true. In other words, with ASP the experiment designer's purpose is to \textit{validate} the optimality of a point on which ASP is converging.

One considers that the autopilot has converged if the minimum  $\widehat{\bm{x}}_{k+1}$ of the model corrected at $\widehat{\bm{x}}_{k}$  is such that it lies in the rectangular domain centered on  $\widehat{\bm{x}}_{k}$ and on which one considers $f_p$ to be almost affine. In more mathematical terms, one would say that the autopilot converges if: 
\begin{equation} \label{eq:6_38_Critere_Convergence}
	\widehat{\bm{x}}_{k+1}
	\in
	\mathcal{A}(\widehat{\bm{x}}_{k}),	
\end{equation}
where $\mathcal{A}(\widehat{\bm{x}}_{k}):=$ \eqref{eq:6_14_definition_A_x}. 

If at an iteration $k$ the criterion \eqref{eq:6_38_Critere_Convergence} is satisfied, then one considers that the convergence is almost reached and one must make sure that the next experiments bring information on the first order optimality condition of the plant, \textit{if the latter is too uncertain}. To do so, the three following indicators are used together  with an acceptable loss ($loss_{acc}$) and a scaling matrix ($\bm{\Delta}u_{scal}$) that must be defined by the engineers:

\begin{Definition} (The acceptable loss: $loss_{acc}$; and the scaling matrix: $\bm{\Delta u_{scal}}$) In theory, one considers that one is at the optimum of a plant if the projection of derivative of the cost function  on the null space of the active constraints is null. In practice, this  never happens since no perfectly accurate measurements are available. So, one needs to replace it with something else. One suggests to use a scaling matrix to quantify what a ``small'' step in the inputs is:
	\begin{equation}
		\bm{\Delta u_{scal}} = \left(
		\begin{array}{ccc}
			\Delta u_{scal,1} &        & \\
			& \ddots & \\
			&        & \Delta u_{scal,n_u}
		\end{array}
		\right),
	\end{equation}
where $\bm{\Delta u_{scal,i}}$ is what the engineers managing the plant consider as being a ``small'' variation of the $i$-th input. These small variations being defined, one considers that the plant is at an acceptable (i.e. quasi-optimal) operating condition if a small feasible step $\bm{\Delta u_{scal}}\bm{v}$, where $\bm{v}\in\amsmathbb{R}^{n_u}$ and $\|\bm{v}\|=1$ cannot imply a ``great'' improvement of the cost. It is also the engineers managing the plant that have to define what such a  ``great'' improvement is with a parameter $loss_{acc}$. Basically, all cost improvement greater than $loss_{acc}$ are ``great'', and the others are acceptable. 
\end{Definition}

 
\begin{itemize}
	\item \textbf{Indicator A:} This indicator is focused on the uncertainty on the \textit{gradients} of the cost function at $\widehat{\bm{x}}_{k+1}$ \textit{in the unconstrained directions}.
	Its objective is to check whether the operating conditions of the plant can be improved by moving it in unconstrained directions. 
	To build this indicator one needs $\bm{g}^a_k$: the vector of the active constraints of the updated model at $\widehat{\bm{x}}_{k+1}$. Then, to evaluate if improvements are possible \textit{in the unconstrained directions}, one just need to check whether $\exists \bm{v} \in \amsmathbb{R}^{n_u}$ such that:
	\begin{align}
		2\sqrt{	\bm{v}_{scal,proj}^{\rm T} \amsmathbb{V}[\nabla_{\bm{u}}\phi_k(\bm{u}_{k+1},\bm{f}_k(\widehat{\bm{x}}_{k+1}))] 	\bm{v}_{scal,proj}} \geq loss_{acc}.
	\end{align}
	where 
	\begin{align}
		\bm{v}_{scal,proj} :=  \bm{\Delta u_{scal}}\bm{v} - \sum_{i=1}^{n_a}  \frac{\amsmathbb{E}[\nabla_{\bm{u}}g^{a}_{k(i)}]^{\rm T} \bm{\Delta u_{scal}} \bm{v}
		}{\|\amsmathbb{E}[\nabla_{\bm{u}}g^{a}_{k(i)}]\|^2} \amsmathbb{E}[\nabla_{\bm{u}}g^{a}_{k(i)}]. \label{eq:6:39_v_scal_proj}
	\end{align}
	If such a $\bm{v}$ exists, then the actual optimality of what the ASP considers as being the plant optimum is questionable and one should increase the knowledge one has on the plant's gradients along this direction $\bm{v}$. 
	\item \textbf{Indicator B:} This indicator is focused on the uncertainty on the \textit{values} of the active constraint at $\widehat{\bm{x}}_{k+1}$. 
	Its objective is to check whether the operating condition of the plant can be improved thanks to better knowledge on what is believed to be the active constraints of the plant. 
	To build this indicator, one uses one of the common interpretations of Lagrange multipliers, that they measure the marginal benefit that can be obtained from the relaxation of the constraint with which it is associated \cite{Bodenhorn:1974}. 
	In the context of RTO, these multipliers have an analogous meaning because they allow us to associate a cost to the uncertainty on an active constraints.  Indeed, if one separates the uncertainties from the expectations on the constraints and consider the following two optimization problems: 
	\begin{align}
		\phi_A := \ & \underset{\bm{u}}{\operatorname{min}} \  \varphi (\bm{u},\bm{f}_k(\widehat{\bm{x}})) \  s.t. \   \amsmathbb{E}[g^a_{k(i)}(\bm{u},\bm{f}_k(\widehat{\bm{x}}))] = 0, \ \forall i \ \  active,  \nonumber  \\
		\phi_B := \ & \underset{\bm{u}}{\operatorname{min}} \  \varphi (\bm{u},\bm{f}_k(\widehat{\bm{x}})) \  s.t. \  \amsmathbb{E}[g^a_{k(i)}(\bm{u},\bm{f}_k(\widehat{\bm{x}}))] +\epsilon_{g_{(i)}} = 0, \ \forall i \ \  active, \nonumber \\
		\epsilon_{g_{(i)}} \sim \ & \mathcal{N}(0,\amsmathbb{V}[g^a_{k(i)}(\bm{u},\bm{f}_k(\widehat{\bm{x}}))]), \forall i, \nonumber
	\end{align}
	Then the difference $\phi_B-\phi_A$ can be bounded (with a high probability) as follows:
	\begin{align}
		|\phi_B  - \phi_A |=  \sum_i \lambda_{(i)} \epsilon_{g_{(i)}} \leq   2 \sum_i \lambda_{(i)} \sqrt{\amsmathbb{V}[g^a_{k(i)}(\bm{u},\bm{f}_k(\widehat{\bm{x}}))] },
	\end{align}
	(were one assumes that $\epsilon_{g_{(i)}}$ are independent $\forall i$, hence the absence of covariance matrices). 
	If one associates this loss with what has been defined as the matrix of tolerable loss,  $\bm{\Delta} \phi_{min}$, the following criterion can be used to assess whether additional experiments are required to validate that $\widehat{\bm{x}}_{k+1}$ is indeed an acceptable operating point to stay at:
	\begin{equation} 
		2 \sum_{i} \big( \lambda_{(i)} \sqrt{ \amsmathbb{V}[g^a_{k(i)} (\bm{u}_{k+1},\bm{f}_k(\widehat{\bm{x}}_{k+1}))]} \big) \leq loss_{acc}
	\end{equation}
	If this criterion is not met, then more experiments should be conducted in the vicinity of $\widehat{\bm{x}}_{k+1}$ to reduce the values of $\amsmathbb{V}[g^a_{k(i)} (\bm{u}_{k+1},\bm{f}_k(\widehat{\bm{x}}_{k+1}))]$ $\forall i$. 
	\item \textbf{Indicator C:} This indicator is focused on the  \textit{sign} of the gradients of the cost function $\widehat{\bm{x}}_{k+1}$ \textit{in the constrained directions}.  
	Its objective is to verify if what the model predicts as active constraints can be, with additional experiments, deactivated. 
	To carry out this verification it is enough to check whether the \textit{sign} of the projection of the model cost gradient  at $\widehat{\bm{x}}_{k+1}$  on the constrained directions can change. In other words, if $\exists \bm{v} \in \amsmathbb{R}^{n_u}$  such that 
	\begin{align}
		& 	\big|
		\bm{v}^{\prime}_{scal,proj}{}^{\rm T}
		\amsmathbb{E}[\nabla_{\bm{u}}\phi_k(\bm{u}_{k+1},\bm{f}_k(\widehat{\bm{x}}_{k+1}))]\big| < ... \\
		&\qquad \qquad 2 \sqrt{
			\bm{v}^{\prime}_{scal,proj}{}^{\rm T}
			\amsmathbb{V}[\nabla_{\bm{u}}\phi_k(\bm{u}_{k+1},\bm{f}_k(\widehat{\bm{x}}_{k+1}))]
		\bm{v}^{\prime}_{scal,proj}
		} . \label{eq:6_42_IsActive}
	\end{align}
	where 
	\begin{align}
		\bm{v}^{\prime}_{scal,proj} :=  \bm{\Delta u_{scal}}\bm{v} -\bm{v}_{scal,proj},
	\end{align}
	then there is an non-negligible probability that additional experiments may lead to changes of the sign of the updated model gradients in the constrained directions, i.e. to the deactivation of active constraints. 
	Indeed, the term on the right gives the maximum probable (at 96\%) variation of $ \nabla_{\bm{u}}\phi_k(\bm{u}_{k+1},\bm{f}_k(\widehat{\bm{x}}_{k+1}))$ around  $\amsmathbb{E}[\nabla_{\bm{u}}\phi_k(\bm{u}_{k+1},\bm{f}_k(\widehat{\bm{x}}_{k+1}))]$. So, if this value is greater than the absolute value of $\amsmathbb{E}[\nabla_{\bm{u}}\phi_k(\bm{u}_{k+1},\bm{f}_k(\widehat{\bm{x}}_{k+1}))]$, then it is clear that the sign of $ \nabla_{\bm{u}}\phi_k(\bm{u}_{k+1},\bm{f}_k(\widehat{\bm{x}}_{k+1}))$ could change if more data were collected. Also, if this sign changes, then what was an active constraint of the updated model could be deactivated.  In such a case, it would be possible for indicators A and B to be misused because of poor evaluation of active constraints due to too much uncertainty about the gradients of the plant.
\end{itemize}
A graphical interpretation of these three indicators is proposed below:
\begin{Side}
	\begin{GraphicalInterpretation}   \label{gi:Indicators}
	Let's consider a problem where only one variable can be manipulated. In this case, the optimum predicted by the model can correspond to one of the three cases of Figure~\ref{fig:6_10_Graph_Interpretation}, that are:
	\begin{itemize}
		\item The optimum does not activate any constraint. 
		\item The optimum activates a  constraint subject to uncertainties.
		\item The optimum activates a  constraint subject to no uncertainties.
	\end{itemize}
\begin{minipage}[h]{\linewidth}
\smallskip
\centering
\includegraphics[width=14cm]{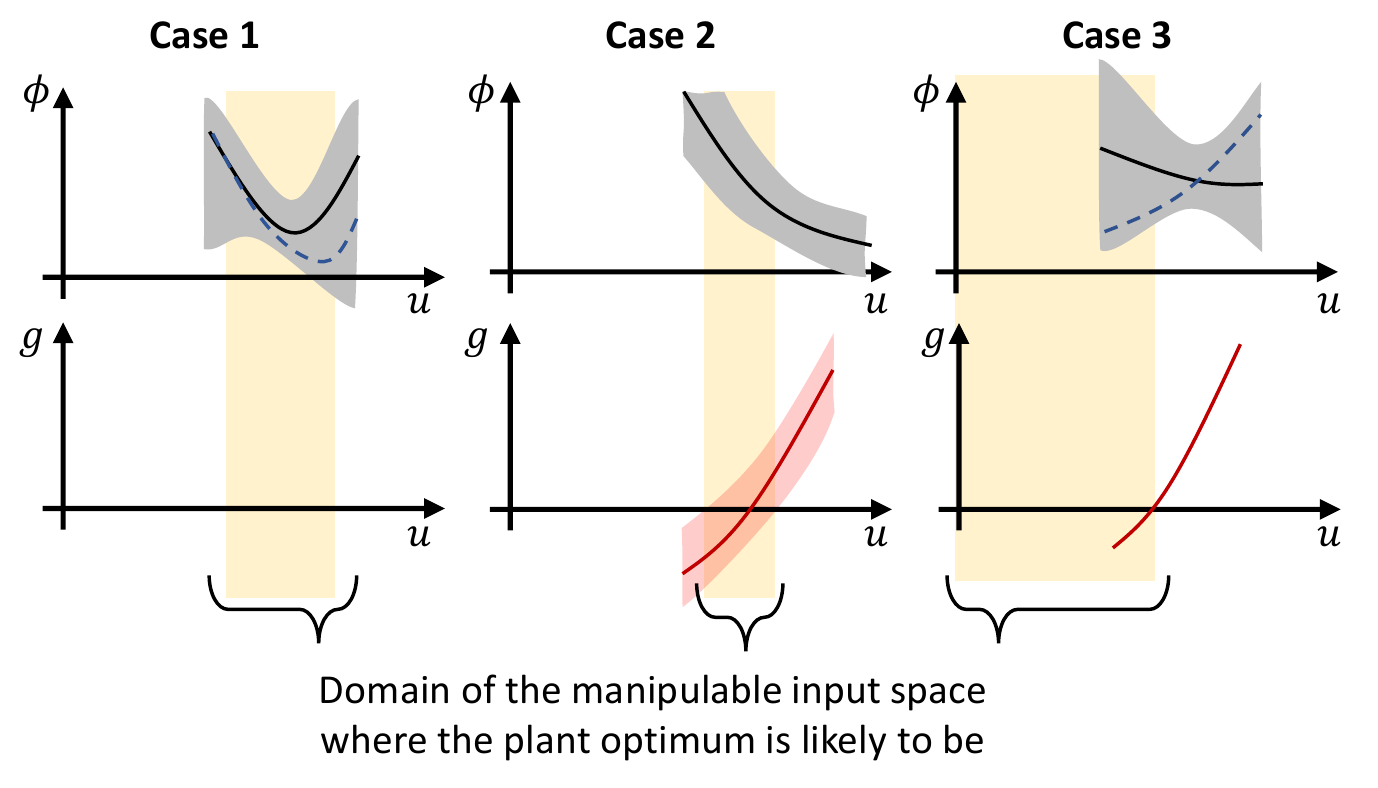}
\captionof{figure}{Three cases}
\label{fig:6_10_Graph_Interpretation} 	
\smallskip
\end{minipage} 
	As it will be seen now, the three indicators manage these three situations correctly. Let's discuss them one by one:
	\begin{itemize}
		\item In case 1, since no constraint is active, indicators B and C are zero while indicator A gives an estimate of what the probable loss  due to uncertainty on the cost is. For example, if the cost expectation is the black curve, the blue curve presents an alternative whose  likelihood  is quantifiable and whose minimum may be different. Indicator A quantifies the probable losses associated with these alternatives and thus motivates additional experiments around the minimum of the model. 
		
		\item In case 2, since a constraint is active and since there is only one input, the indicator A is necessarily zero. Moreover, as there is little uncertainty about the cost function, it is sure that the active constraint of the model is an active constraint of the plant. This should be confirmed by indicator C. Therefore, only the uncertainty on the constraint function affects the decisions. To know if one needs to reduce this uncertainty, one can link this uncertainty to the cost function to evaluate the probable losses associated with it. This is what indicator B does.   
		
		\item In case 3, since a constraint is active and since there is only one input, the indicator A is necessarily null. Moreover, since there is no uncertainty about the constraint, indicator B is also null. However, as it can be observed, there are probable alternatives (blue dotted curve) to the expected cost function (black curve), which if they corresponded to the cost function of the plant would cause the deactivation of this constraint. There is therefore a doubt about the activation of this constraint. This doubt is quantified by the indicator C. 
	\end{itemize} 
	These three cases form almost all the possibilities that can be encountered when $n_u=1$. The only case they do not contain is the case where the activation of an uncertain constraint is not sure (a combination of cases 2 and 3).  

	\end{GraphicalInterpretation}
\end{Side}
These indicators should be used as following:

First, the indicator C should be checked. If it identifies that some active constraints could in fact be inactive, i.e. if \eqref{eq:6_42_IsActive} holds, then the next experiment should be such that the uncertainty  $\amsmathbb{V}[\nabla_{\bm{u}}\phi_k(\bm{u}_{k+1},\bm{f}_k(\widehat{\bm{x}}_{k+1}))]$ decreases the most. For instance one could take the greatest eigenvalue of $\amsmathbb{V}[\nabla_{\bm{u}}\phi_k(\bm{u}_{k+1},\bm{f}_k(\widehat{\bm{x}}_{k+1}))]$ and make the next step in the direction of its associated eigenvector (in the feasible side of the input space). 

If the indicator C confirms that the active constraints seems to be right, then the indicators A and B should be checked. Since both of them quantify the potential loss associated to the updated-model's uncertainties, they can be combined to assess total loss associated to them w.r.t. a direction $\bm{v}\in\amsmathbb{R}^{n_u}$ ($loss_{tot}(\bm{v})$): 
\begin{align}
	loss_{tot}(\bm{v}) := \ & 2\sqrt{	\bm{v}_{scal,proj}^{\rm T} \amsmathbb{V}[\nabla_{\bm{u}}\phi_k(\bm{u}_{k+1},\bm{f}_k(\widehat{\bm{x}}_{k+1}))] 	\bm{v}_{scal,proj}} + ... \nonumber \\
		& \quad 2 \sum_{i} \big( \lambda_{(i)} \sqrt{ \amsmathbb{V}[g^a_{k(i)} (\bm{u}_{k+1},\bm{f}_k(\widehat{\bm{x}}_{k+1}))]} \big), \nonumber \\
		& s.t. \ \ \bm{v}_{scal,proj} =  \eqref{eq:6:39_v_scal_proj}.
\end{align}
From this definition one can define the maximal potential loss due to uncertainties ($loss_{max} $) and the search direction ($\bm{v}_{search}$) associated to this maximal uncertainty: 
\begin{align}
	loss_{max}      := \ & \underset{\bm{v}\in\amsmathbb{R}^{n_u}, \|\bm{v}\|=1}{\operatorname{min}} loss_{tot}(\bm{v}), \\
	\bm{v}_{search} := \ & \operatorname{arg} \ \underset{\bm{v}\in\amsmathbb{R}^{n_u}, \|\bm{v}\|=1}{\operatorname{min}} loss_{tot}(\bm{v}),
\end{align}
If $loss_{max} > loss_{acc}$, then additional experiments are required to validate (or invalidate) $\widehat{\bm{x}}_{k+1}$ as an acceptable operating condition to stay on. Ideally, the next experiment should be performed in the feasible neighborhood of  $\widehat{\bm{x}}_{k+1}$ and in the direction of $\bm{v}_{search}$, which is the direction that is the most affected by the model uncertainties.

\begin{Remark} \textbf{(How to compute $\amsmathbb{E}[\nabla_{\bm{u}}\phi_k(\bm{u}_{k+1},\bm{f}_k(\widehat{\bm{x}}_{k+1}))]$ and  $\amsmathbb{V}[\nabla_{\bm{u}}\phi_k(\bm{u}_{k+1},$ \\ $\bm{f}_k(\widehat{\bm{x}}_{k+1}))]$)}
	
	On sets $(\epsilon_j,\bm{\tau}_j)$ $\forall j=1,...,n_y$, the uncertainties on the values and the gradients of $f_{p(j)}$, i.e. $\forall j = 1,...,n_y$:
	\begin{align}
		f_{k(j)}                 = \ & f_{p(j)} + \epsilon_j, & 
		\nabla_{\bm{u}} f_{k(j)} = \ & \nabla_{\bm{u}} f_{p(j)} + \bm{\tau}_j, &
		\left(\begin{array}{c}
			\epsilon_j \\ 
			\bm{\tau}_j
		\end{array}\right) \sim \ & \mathcal{N}\left( \bm{0},
		\left(\begin{array}{cc}
			\sigma_j^2        & \bm{\varsigma}_j^{\rm T} \\
			\bm{\varsigma}_j  & \bm{\Sigma}_j
		\end{array}\right) \right),
	\end{align}
	where $\sigma_{j}^2\in\amsmathbb{R}$ is the variance of  $\epsilon_j$, 
	$\bm{\Sigma}_{j}\in\amsmathbb{R}^{n_u\times n_u}$  is the variance of    $\bm{\tau}_j$, and $\bm{\varsigma}_{j} \in\amsmathbb{R}^{n_u}$ is the covariance of  $(\epsilon_j,\bm{\tau}_j)$. These variances and covariances can be computed on the basis of the model and the database as explained in Appendix~\ref{App:Predict_Val_and_Grad} with equation~\eqref{eq:C_17_Variances_et_Covariances}. Let's write the Taylor series of function $\phi$ around $(\bm{u}_k,\bm{f}_{k}(\widehat{\bm{x}}_k))$:
	\begingroup
	\allowdisplaybreaks
	\begin{align*}
		\phi(\bm{u},\bm{f}_p(\bm{x})) = \ & 
		\phi\left(
		\bm{u}_k, \left[
		\begin{array}{c}
			f_{k(1)}(\widehat{\bm{x}}_k)+\epsilon_1+ \bm{\tau}_1^{\rm T}\delta\bm{u} \\
			\vdots \\
			f_{k(n_y)}(\widehat{\bm{x}}_k) +\epsilon_{n_y}+\bm{\tau}_{n_y}^{\rm T}\delta\bm{u}
		\end{array}
		\right]
		\right), \\
		= \ & 
		\phi(\bm{u}_k,\bm{f}_k(\widehat{\bm{x}}_k))
		+ 
		\sum_{j=1}^{n_y} \Big[
		\partial_{y_{(j)}}\phi(\bm{u}_k,\bm{f}_k(\widehat{\bm{x}}_k)) (\epsilon_j +  \bm{\tau}_j^{\rm T}\delta\bm{u})
		\Big]. 
	\end{align*}
	\endgroup
	By taking the first derivative:
	\begin{align}
		& \nabla_{\bm{u}}\phi(\bm{u},\bm{f}_p(\bm{x})) =   \nabla_{\bm{u}} \phi(\bm{u}_k,\bm{f}_k(\widehat{\bm{x}}_k))  + 
		\sum_{j=1}^{n_y} \Big[ 
		\big[ 
		\partial_{\bm{u}}\partial_{y_{(j)}}\phi(\bm{u}_k,\bm{f}_k(\widehat{\bm{x}}_k)) 
		+ 
		\nabla_{\bm{u}} \bm{f}_k(\widehat{\bm{x}}_k)^{\rm T} ...
		\nonumber  \\ 
		&  \qquad 
		\partial_{\bm{y}}\partial_{y_{(j)}}\phi(\bm{u}_k,\bm{f}_k(\widehat{\bm{x}}_k))
		\big] \big(\epsilon_j + \sum_{j=1}^{n_y}\big[ \bm{\tau}_j^{\rm T}\big](\bm{u}-\bm{u}_k)\big) +
		\bm{\tau}_j 
		\partial_{y_{(j)}}\phi(\bm{u}_k,\bm{f}_k(\widehat{\bm{x}}_k))
		\Big],\label{eq:TaylorApprox_gi}
	\end{align}
	one deduces that:
	\begingroup
	\allowdisplaybreaks
	\begin{align}
		\amsmathbb{E}[\nabla_{\bm{u}}\phi(\bm{u},\bm{f}_p(\bm{x}))] = \ & \nabla_{\bm{u}} \phi(\bm{u}_k,\bm{f}_k(\widehat{\bm{x}}_k))  \\
		\amsmathbb{V}[\nabla_{\bm{u}}\phi(\bm{u},\bm{f}_p(\bm{x}))] = \ & 
		 \sum_{j=1}^{n_y} 
		\amsmathbb{V}[
		\epsilon_j \bm{l}_{j} + 
		\bm{\tau}_j l_{j}^{\prime}], \nonumber \\
		=  \ & 
		 \sum_{j=1}^{n_y}  
		\sigma_{j}^2 \bm{l}_{j}\bm{l}_{j}^{\rm T} + 
		l_{j}^{\prime 2} \bm{\Sigma}_{j} + \bm{l}_{j} \bm{\varsigma}_{j}^{\rm T} l_{j}^{\prime} +    l_{j}^{\prime} \bm{\varsigma}_{j}\bm{l}_{j}^{\rm T}. \label{eq:ModelOpt_var_PlantOpt}
	\end{align}
	\endgroup
	where:
	\begingroup
	\allowdisplaybreaks
	\begin{align}
		\bm{l}_{j} := \ & \big( 
		\partial_{\bm{u}} \partial_{y_{(j)}} \phi (\bm{u}_k^{\star},\bm{f}_k(\widehat{\bm{x}}_k^{\star})) + 
		\nabla_{\bm{u}} \bm{f}_k(\widehat{\bm{x}}_k^{\star})^{\rm T} \partial_{\bm{y}} \partial_{y_{(j)}} \phi (\bm{u}_k^{\star},\bm{f}_k(\widehat{\bm{x}}_k^{\star}))
		\big), \\
		l_{j}^{\prime} := \ & \partial_{y_{(j)}} \phi(\bm{u}_k^{\star},\bm{f}_k(\widehat{\bm{x}}_k^{\star})).
	\end{align}
	\endgroup
	Which concludes this remark.
\end{Remark}

\section{The improved decision trigger}
\label{sec:improved_decision_trigger}

The ASP's decision trigger reacts to the six events listed in section~\ref{sec:4_4_1_Declencheur_de_decision} in the following way: 
\begin{itemize}
	\item (E1) The change of a cost or a constraint is easy to detect since it is a change that must be done on the autopilot's software and it is obviously simple to detect such a change.  When this happens, the autopilot must stop what it is doing and use the updated model to revise the last decision.
	\item (E2) The availability of a new model is easy to detect. Obviously, when this event is detected, a decision must be made. 
	\item (E3) One decides not to consider this event in the ASP. Indeed, if the engineers are not satisfied with an experiment in progress because, for example, it violates constraints, they will have to force the autopilot to change its decision. What is recommended is to let the experiment finish in order to measure the output of the plant (even if it is not satisfactory). Once these measurements are obtained, the Autopilot will realize that this area of the input space is dangerous or uninteresting based on real data. If the engineers really don't want to finalize an experiment, they have to restart a decision while modifying the database the autopilot uses so that it doesn't make the same decision as before. There are then three options available: 
	\begin{itemize}
		\item Estimate what would have been the experimental results of the interrupted experiment at the risk of creating inconsistencies in the database. This could lead to an erroneous detection of the CTP transition and an unnecessary deletion of the database.
		\item Modify or add constraints, at the risk of losing the optimality upon  convergence guarantees ASP provides. 
		\item Increase the a priori uncertainty on the location of $f_p$ by increasing $\sigma_f^2$, see equation~\eqref{eq:6_7_GP_Fp}. This should have the effect of reducing the size of the iterations to the unvisited areas of the input space resulting in a more secure but slower progression. 
	\end{itemize}
	\item (E4) Discussed later.
	\item (E5) An anomaly in the plant's behavior can be detected with the consistency monitor to which one can associate a fault detector. 
	\item (E6) If the engineers want to relaunch the decision making, then this can be very easily integrated into the autopilot software. However, if they don't want the autopilot to restart the decision they interrupted, they will either have to modify the databases, modify the nominal model, or modify (or add) the cost and constraint functions.   
\end{itemize}

The event (E4) is not discussed in this list because it will be given special attention. Indeed, the event (E4) is associated with variations of the measured disturbances. $\widehat{\bm{d}}^{\prime}$. Ideally, a decision should be made at each variation of $\widehat{\bm{d}}^{\prime}$. However, it is clear that the measures of $\widehat{\bm{d}}^{\prime}$ vary all the time, because of measurement noise, among other things.  It is therefore not practical to trigger a decision at each variation of $\widehat{\bm{d}}^{\prime}$. However, it is possible to approach this ideal behavior with an appropriate revision of the \textit{RTO-plant interface} managing the small variations of $\widehat{\bm{d}}^{\prime}$. Then, the \textit{decision trigger} only has to react to ``significant'' changes in  $\widehat{\bm{d}}^{\prime}$, see Remark~\ref{rem:6_10_SignificantDisturbanceVar}. Now let's see how the \textit{RTO-plant interface} is revised. 

\section{The improved RTO-plant interface}
\label{sec:VLC_VLC_VLC}
To maintain the plant at its optimal operating conditions despite the presence of high-frequency disturbances $\widehat{\bm{d}}^{\prime}$, it would be necessary to be able to solve the optimization problem \eqref{eq:6_29_Model_based_PB}  at the same high-frequency. Doing so would turn the autopilot into a kind of controller (supposed to make high frequency decisions). However, the time required to compute \eqref{eq:6_29_Model_based_PB} can make the reaction to disturbances ineffective or even harmful. To avoid these problems, one proposes to turn the RTO-plant interface into a \textit{virtual control layer} (VCL) that would be placed upstream of the plant control system.  This VCL would contain the explicit solution of a multi-parametric quadratic program (MP-QP) approximating  \eqref{eq:6_29_Model_based_PB} around $\widehat{\bm{x}}_{k+1}$ and whose parameters are the measured disturbances. In short, one proposes to use the following quadratic and affine  approximations $\{\widetilde{\varphi},\widetilde{\bm{h}}^a\}$ of $\{\varphi,\bm{h}^a\}$:
\begin{align}
	\widetilde{\varphi}(\bm{u},\bm{f}_k(\widehat{\bm{x}})) := \ & 
	\varphi(\bm{u}_{k+1},\bm{f}_k(\widehat{\bm{x}}_{k+1}))  + 
	\nabla_{\bm{u}} \varphi|_{\bm{u}_{k+1},\bm{f}_k(\widehat{\bm{x}}_{k+1})}
	\delta\bm{u} 
	+  ... \nonumber \\
	& 
	\nabla_{\bm{d}^{\prime}} \varphi|_{\bm{u}_{k+1},\bm{f}_k(\widehat{\bm{x}}_{k+1})} 
	\delta\bm{d}^{\prime}    
	+ 
	\frac{1}{2}
	\delta\bm{u}^{\rm T} 
		\nabla_{\bm{u}\bm{u}} \varphi|_{\bm{u}_{k+1},\bm{f}_k(\widehat{\bm{x}}_{k+1})}
	\delta\bm{u} 
	+ ... \nonumber \\
	& \
	\frac{1}{2}
	\delta\bm{d}^{\prime \rm T}
		\nabla_{\bm{d}^{\prime}\bm{d}^{\prime}} \varphi|_{\bm{u}_{k+1},\bm{f}_k(\widehat{\bm{x}}_{k+1})}
	\delta\bm{d}^{\prime}
	+ 
	\delta\bm{d}^{\prime \rm T}
		\nabla_{\bm{d}^{\prime}\bm{u}} \varphi|_{\bm{u}_{k+1},\bm{f}_k(\widehat{\bm{x}}_{k+1})}
	\delta\bm{u}. \\ 
	\widetilde{\bm{h}}^a(\bm{u},\bm{f}_k(\widehat{\bm{x}})) := \ &  
	\bm{h}^a(\bm{u}_{k+1},\bm{f}_k(\widehat{\bm{x}}_{k+1})) + 
	\nabla_{\bm{u}} \bm{h}^a|_{\bm{u}_{k+1},\bm{f}_k(\widehat{\bm{x}}_{k+1})}
	\delta\bm{u} 
	+  
	... \nonumber \\
	& 
	\nabla_{\bm{d}^{\prime}} \bm{h}^a|_{\bm{u}_{k+1},\bm{f}_k(\widehat{\bm{x}}_{k+1})} 
	\delta\bm{d}^{\prime}, 
\end{align}
where $\bm{h}^a$ are the  active constraints at $\widehat{\bm{x}}_{k+1}$, and 
\begin{align*}
	\delta\bm{u}                     := \ & \bm{u}-\bm{u}_{k+1}, &
	\delta\widehat{\bm{d}}^{\prime}  := \ & \widehat{\bm{d}}^{\prime}-\widehat{\bm{d}}_{p,k}^{\prime}.
\end{align*}
To simplify the following reasoning, let's rewrite these expressions as:
\begin{align}
	\widetilde{\varphi}(\bm{u},\bm{f}_k(\widehat{\bm{x}})) := \ &
	\varphi_{\star}(\delta\bm{d}^{\prime}) + \bm{q}(\delta\bm{d}^{\prime})^{\rm T}\delta\bm{u} + \frac{1}{2} \delta\bm{u}^{\rm T} \bm{Q} \delta\bm{u}, \\
	\widetilde{\bm{h}}^a(\bm{u},\bm{f}_k(\widehat{\bm{x}})) := \ &  \bm{h}^a_{\star}(\delta\bm{d}^{\prime}) + \bm{H}\delta\bm{u},
\end{align}
where: 
\begingroup
\allowdisplaybreaks
\begin{align*}
	\varphi_{\star}(\delta\bm{d}^{\prime}) := \ &\varphi(\bm{u}_{k+1},\bm{f}_k(\widehat{\bm{x}}_{k+1})) +  \nabla_{\bm{d}^\prime} \varphi|_{\bm{u}_{k+1},\bm{f}_k(\widehat{\bm{x}}_{k+1})}^{\rm T} \delta\widehat{\bm{d}}^{\prime} + ... \qquad \qquad \qquad \qquad \qquad \qquad  \\ 
	&  \frac{1}{2} (\delta\widehat{\bm{d}}^{\prime} )^{\rm T}\nabla_{\bm{d}^\prime\bm{d}^\prime}\varphi|_{\bm{u}_{k+1},\bm{f}_k(\widehat{\bm{x}}_{k+1})} \delta\widehat{\bm{d}}^{\prime} , \\
	\bm{q}(\delta\bm{d}^{\prime})^{\rm T} := \ 
	& \nabla_{\bm{u}} \varphi|_{\bm{u}_{k+1},\bm{f}_k(\widehat{\bm{x}}_{k+1})} +  
	\nabla_{\bm{ud}^\prime} \varphi|_{\bm{u}_{k+1},\bm{f}_k(\widehat{\bm{x}}_{k+1})}	\delta\widehat{\bm{d}}^{\prime} , \\
	\bm{Q} := \ 
	& \nabla_{\bm{uu}} \varphi|_{\bm{u}_{k+1},\bm{f}_k(\widehat{\bm{x}}_{k+1})}, \\
	\bm{h}_{\star}(\delta\bm{d}^{\prime}) := \ 
	&  \bm{h}^a(\bm{u}_{k+1},\bm{f}_k(\widehat{\bm{x}}_{k+1})) + \nabla_{\bm{d}^\prime}\bm{h}^a|_{\bm{u}_{k+1},\bm{f}_k(\widehat{\bm{x}}_{k+1})}\delta\widehat{\bm{d}}^{\prime} , \\
	\bm{H} := \ 
	& \nabla_{\bm{u}}\bm{h}^a|_{\bm{u}_{k+1},\bm{f}_k(\widehat{\bm{x}}_{k+1})}.
\end{align*}
\endgroup
If one replaces $\{\varphi,\bm{h}\}$ with $\{\widetilde{\varphi},\widetilde{\bm{h}}^a\}$ in 
\eqref{eq:6_29_Model_based_PB}, the the following   optimization problem is obtained:
\begin{align}
	\delta\bm{u}^{\star}:=  \operatorname{arg}
	\underset{\delta\bm{u}}{\operatorname{min}} \  &  \bm{q}(\delta\bm{d}^{\prime})^{\rm T}\delta\bm{u} + \frac{1}{2} \delta\bm{u}^{\rm T} \bm{Q} \delta\bm{u} \nonumber \\
	\text{s.t.} \quad & 
	\bm{h}_{\star}(\delta\bm{d}^{\prime})+ \bm{H}\delta\bm{u} = \bm{0},
	\label{eq:6_58_MPQP_VCL} 
\end{align}
whose solution w.r.t.  $\delta\bm{d}^\prime$ can be explicitly calculated using the general solution formula of a QP (see Appendix~\ref{sec:A_2_1_Resoudre_QP_Contraintes_Egalite}, equation~\eqref{eq:A_3_sol_Analytique}):  
\begin{equation}
	\left(\begin{array}{c}
		\delta \bm{u}^{\star}\\
		\bm{\lambda}
	\end{array}
	\right)
	= 
	-
	\left(\begin{array}{cc}
		\bm{Q} & \bm{H}^{\rm T} \\
		\bm{H} & \bm{0}
	\end{array}
	\right)^{-1}
	\left(\begin{array}{c}
		\bm{q}(\delta\bm{d}^{\prime}) \\
		\bm{h}_{\star}(\delta\bm{d}^{\prime})
	\end{array}
	\right),
\end{equation}
where $\bm{\lambda}$ is the vector of Lagrange multipliers associated with the active constraints.  If one defines:
\begin{equation*}
	\bm{\Psi} := \left(\begin{array}{cc}
		\bm{\Psi}_{11} & \bm{\Psi}_{12} \\
		\bm{\Psi}_{21} & \bm{\Psi}_{22}
	\end{array}
	\right)
	:=
	\left(\begin{array}{cc}
		\bm{Q} & \bm{H}^{a\rm T} \\
		\bm{H}^a & \bm{0}
	\end{array}
	\right)^{-1}
\end{equation*}
where $\bm{\Psi}_{11}\in\amsmathbb{R}^{n_u\times n_u}$ and $\bm{\Psi}_{12}\in\amsmathbb{R}^{n_u\times a}$ are the two matrices constituting the first rows of  $\bm{\Psi}$. Then the optimal variations $\delta\bm{u}^{\star}$ to be added to  $\bm{u}_{k+1}$ when $\widehat{\bm{d}}^{\prime}_k$ is subject to variations $\delta\bm{d}^{\prime}$ can be computed:
\begin{equation}
	\delta \bm{u}^{\star} = 
	\bm{o}^{\star}_k  +\bm{G}^{\star}_{k} \delta\widehat{\bm{d}}^{\prime},
\end{equation}
where
\begin{align*}
	\bm{o}^{\star}_{k} := \ & -\bm{\Psi}_{11} \nabla_{\bm{u}} \varphi|_{\bm{u}_{k+1},\bm{f}_k(\widehat{\bm{x}}_{k+1})} -   \bm{\Psi}_{12} \bm{h}^a|_{\bm{u}_{k+1},\bm{f}_k(\widehat{\bm{x}}_{k+1})} , \nonumber \\
	\bm{G}^{\star}_{k} := \ & -\bm{\Psi}_{11} \nabla_{\bm{ud}} \varphi|_{\bm{u}_{k+1},\bm{f}_k(\widehat{\bm{x}}_{k+1})} - 
	\bm{\Psi}_{12} \nabla_{\bm{d}}\bm{h}^a|_{\bm{u}_{k+1},\bm{f}_k(\widehat{\bm{x}}_{k+1})}.
\end{align*}
Normally, $\bm{o}^{\star}_{k}=\bm{0}$ since by definition $\bm{u}_{k+1}$ is the minimum of \eqref{eq:6_29_Model_based_PB} when $\bm{d}^{\prime} = \widehat{\bm{d}}^{\prime}_k$. Therefore, one can conclude that the optimal variations of the inputs w.r.t. $\delta\bm{d}^{\prime}$  are: 
\begin{equation} \label{eq:6_61VCL_Gains}
	\delta \bm{u}^{\star}= 
	\bm{G}^{\star}_{k}\delta\widehat{\bm{d}}^{\prime}.
\end{equation}

Finally, based on the updated model, the best RTO-plant interface would be the following VCL: 
\begin{equation}	\label{eq:6_62_MBOpt_plus_VCL}
	\bm{u}_{k+1} = \bm{u}_{k+1}^\star + \bm{G}_{k}^{\star}\delta\widehat{\bm{d}}^{\prime}.
\end{equation}

\begin{Remark} \label{rem:6_9_ImprovedVCL}
	The explicit solution of a MP-QP is a piecewise affine function of its parameters (\cite{Tondel:2003}). In this section, the focus is on the piece of this function containing the point $(\bm{u}_{k+1}^{\star},\widehat{\bm{d}}_{k}^{\prime})$ since it is indirectly assumed that the variations  $\delta\bm{d}^{\prime}$ of $\widehat{\bm{d}}_{k}^{\prime}$ are small. Therefore, the optimal variation of $\bm{u}$ w.r.t.  $\delta\bm{d}^{\prime}$ that is proposed with \eqref{eq:6_62_MBOpt_plus_VCL} is only valid when $\widehat{\bm{d}}_{k}^{\prime}$ is small enough not to cause any change in the set of active constraints of the model. However, if this is the case, then one should replace \eqref{eq:6_62_MBOpt_plus_VCL} by a piecewise affine function following a method like the one proposed in \cite{Tondel:2003}. 
\end{Remark}

\begin{Remark} \label{rem:6_10_SignificantDisturbanceVar}
	\textbf{(Significant disturbance variation)}
	It has been seen that a part of the measured disturbances is managed by the RTO-plant interface so that the decision trigger is not permanently triggered. 
	In this section it is explained how the RTO-plant interface can handle the small $\delta\bm{d}^\prime$. It only remains to define what a ``large''  $\delta\bm{d}^\prime$  is to have the triggering condition of (E5). Two options are proposed. (Option 1) One could say that a  $\delta\bm{d}^{\prime}$ is large if \eqref{eq:6_58_MPQP_VCL}  is a bad estimate of \eqref{eq:6_29_Model_based_PB}, i.e. if:
	\begin{equation} \label{eq:6_63_Condition_Deltad_is_too_large}
		\left|
		\begin{array}{ccc}
			\left(
			\begin{array}{l}
				\operatorname{arg}
				\underset{\bm{u}}{\operatorname{min}} \ \varphi(\bm{u},\bm{f}_k(\widehat{\bm{x}}))  \\
				\text{s.t.} \quad  
				\bm{h}(\bm{u},\bm{f}_k(\widehat{\bm{x}})) \leq \bm{0}
			\end{array}
			\right)
			& 
			-
			&
			\left(\bm{u}^{\star}_{k+1} + \bm{G}_{k}^{\star}\delta\widehat{\bm{d}}^{\prime}\right)
		\end{array}
		\right| 
		\leq \bm{\zeta},
	\end{equation}
	where $\bm{\zeta}\in\amsmathbb{R}^{n_u}$ is the vector defining the acceptable estimation error on $\bm{u}_p^{\star}$  that engineers must define. This criterion can be tested for a large number of  $\delta\bm{d}^{\prime}$ at the time of VCL design to determine a domain that, if left, would trigger a VCL reconstruction. (Option 2) The other option would be to let the engineers choose what this domain is. This second option is probably less optimal but is much easier to implement. 
\end{Remark}

\section{Discussion}

\subsection{ASP and self-optimizing control}

If the plant is perpetually disturbed in such a way that it never really reaches a steady state, then no SS data from the plant could be collected and the RTO-plant interface would keep the plant around the minimum of the nominal (un-updated) model $\bm{u}_0^{\star}$  with adjustments for the effects of the measured disturbances $\bm{G}_0^{\star}\delta\widehat{\bm{d}}^{\prime}$  also based on the nominal model.

One can then notice that the ASP would have a behavior very similar to that of a self-optimizing control system (SOC -- \cite{Skogestad:00}) which consists of (i) finding the optimal operating conditions of the nominal model, (ii) extracting active constraints from these conditions and install controllers to monitor them, and finally (iii) finding the linear combination of control variables that when held constant at their optimal values (based on the nominal model) minimize the effects of perturbations on the optimality of the model. There are two notable differences between SOC and ASP: 
\begin{itemize}
	\item ASP does not require any particular control structure, whereas  SOC requires that the constraints ``to be activated'' are controlled.
	\item ASP, unlike SOC, learns from the results of the $k$ experiments that have been conducted and can improve both its estimate of the plant's optimum position  ($\bm{u}_0^\star \rightarrow \bm{u}_k^{\star}$), as well as  the small disturbance compensation strategy  ($\bm{G}_0^{\star} \rightarrow \bm{G}_k^{\star}$).
\end{itemize}

\subsection{ASP and the practical challenges of RTO}

The following is a list of practical RTO challenges associated with ASP's response to them. This list can be seen as a kind of summary of the contributions of this chapter. 
\begin{itemize}
	\item \textbf{Challenge 1 (Optimality guarantee  in case of convergence):} ASP guarantees the optimality of the plant upon convergence  thanks to the validation procedures encapsulated in the experiment designer (see section~\ref{sec:6_7_The_improved_experiment_designer}). 
	
	\item \textbf{Challenge 2 (Convergence guarantee):}  ASP guarantees that a convergence will be observed after a finite number of experiments. This rather trivial result is the subject of Theorem~\ref{thm:6_1_Garantie_de_convergence} given after this enumeration. 
	
	\item \textbf{Challenge 3 (Ensure compliance with the constraints): }If Assumptions~\ref{ass:6_1_Erreur_fp_f_bornee}, \ref{ass:6_2_Meme_Courbusre} and \ref{ass:6_3_Petites_variances} are true, then the ASP guarantees with a probability of (approximately) 95.4\%  that the plant's constraints will not be violated (see  \eqref{eq:6_30_h}). 
	
	\item \textbf{Challenge 4 (Management of measurement uncertainty):} Since ASP accounts for measurement uncertainty in its model correction strategy, and since measurements of plant values and gradients are not assumed to be ideal, these uncertainties are not a particular problem.
	
	\item \textbf{Challenge 5 (Speed):} ASP is ``fast'' since all superficial decisions such as measuring the gradients of the plant at each iteration have been removed. In a way, one can consider that ASP combines the speed of  constraints adaptation (CA - \cite{Chachuat:08}), the parsimonious exploration of  directional modifier adaptation \cite{Costello:2016}   (with the advantage that the modeling error of the model does not need to be parametric and that the parametric error does not need to be quantified with a set of parameters containing ``those of the plant''), with the optimality upon convergence guarantees of ISO.

	\item \textbf{Challenge 6 (Disturbances management):} ASP is able to handle high frequency and low amplitude measured disturbances thanks to the VCL integrated in the RTO-plant interface (see section~\ref{sec:VLC_VLC_VLC}). Unmeasured, slow and progressive disturbances are taken into account through a time variable (age of the plant) which is considered as a measured disturbance (as illustrated on Figure~\ref{fig:6_7_TimeIsKey}). Unmeasured and instantaneous disturbances are indirectly managed through the maintenance of the database in the consistency monitor (as illustrated on Figure~\ref{fig:6_8_ManageUnmeasuredDisturbance}). 
	
	\item \textbf{Challenge 7 (Automatic reaction to major events):}  Thanks to the the improved decision trigger (section~\ref{sec:improved_decision_trigger}) ASP reacts automatically to most (all) major events that one (the author) could imagine.
	
	\item \textbf{Challenge 8 (Maintenance of the database):} Thanks to the improved  consistency monitor (see section\ref{sec:6_6_Contrôoleur_de_coherence}) ASP is able to identify significant variations in plant behavior and is therefore able to keep the experimental database up to date, i.e. it automatically deletes any data that is no longer useful for modeling the current behavior of the plant.  It should be noted that the method used by ASP to maintain the database  makes it possible to preserve a \textit{maximum} of \textit{up-to-date}  data whereas more simple  methods, e.g. where 
	\vspace{-\topsep}
	\begin{itemize}[noitemsep]
		\item All data are kept,
		\item Only the minimum amount of the most recent data is kept, 
	\end{itemize}
	\vspace{-\topsep}
	do not guarantee that only up-to-date data is used.
\end{itemize}

\begin{thmbox} \label{thm:6_1_Garantie_de_convergence}
	If Assumptions~\ref{ass:6_1_Erreur_fp_f_bornee} and \ref{ass:6_2_Meme_Courbusre} are true and the plant is not subject to CTP transitions, then it is certain that ASP will converge after a finite number of experiments. 
\end{thmbox}
\begin{proofbox}
	Any RTO method can  exhibit only three types of behavior: It can  (i) diverge, (ii) enter a limit cycle, or (iii) converge. Let's analyze whether behaviors (i) and (ii) can appear when using ASP. 
	
	Concerning (i): Since the ASP uses all of the data history and since there is no CTP transition, it is certainly not going to revisit points that it has identified as suboptimal. Also, each domain $\mathcal{A}$ centered on these sub-optimal points are areas that ASP avoids. So, if ASP is continuously exploring new operating points and if the input space of the plant is finite (which is always the case with real processes), then it will eventually have explored the whole space. In this case, as the whole space would be explored, the updated model could only be almost equal to the plant. Therefore, the minimum of the updated model would be almost that of the plant and the ASP would converge. It has been proved by reduction that ASP cannot diverge (indefinitely). 
	
	Concerning (ii): In light of the discussion associated with the point (i), it should be clear that ASP can in no way enter a limit cycle since the domain   $\mathcal{A}$ around a point that has been identified as suboptimal cannot be revisited. So a cycle, i.e. returning to a point already visited, is structurally not possible. 
	
	Finally, the only possible behavior is (iii), which concludes this proof.
\end{proofbox}

\clearpage

\section{Case study -- the Williams-Otto plant}
\label{sec:6_11_Cas_Etude_WOplant}

\subsection{Problem description}

\begin{figure}[b]
	\centering
	\includegraphics[width=14cm]{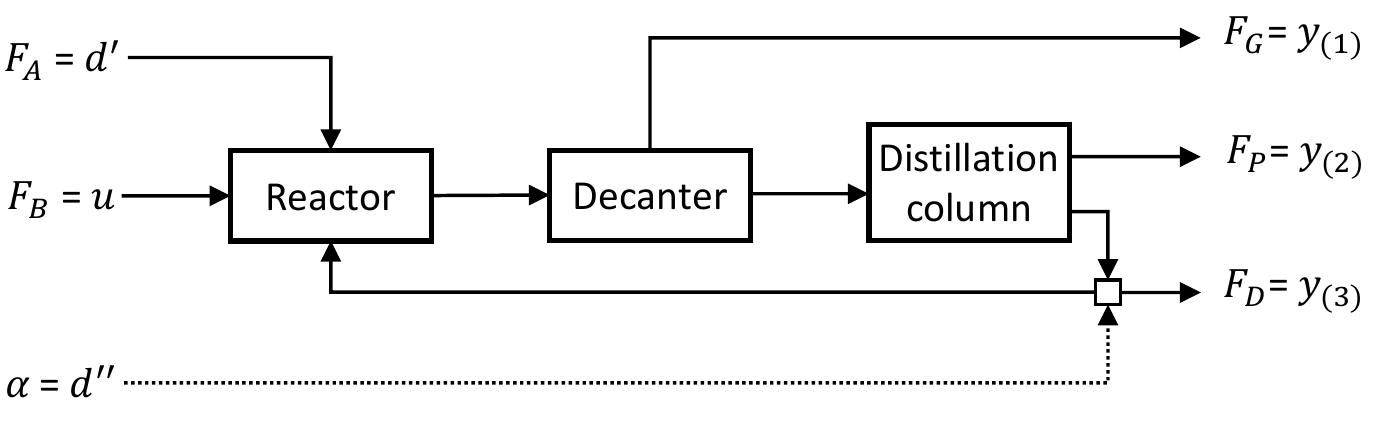}
	\caption{The Williams-Otto plant}
	\label{fig:6_10_WO_plant_description}
\end{figure}

ASP is illustrated on a simplified version of the Williams-Otto plant presented in Figure~\ref{fig:6_10_WO_plant_description} \cite{Williams:60}. The simplification consist in replacing (a) the heat exchanger of the reactor, and (b) the heat exchanger between the reactor and the decanter by thermal control systems. This means that, at steady state, the temperatures in the reactor and of the flow entering the decanter correspond to the fixed temperature setpoints (a) 355.72~K and (b) 310.93~K. The model of the simplified Williams-Otto plant is obtained in the following way: 
\vspace{-\topsep}
\begin{itemize}[noitemsep]
	\item The reactor is modeled with the equations (9)-(15) and (18)-(20) of \cite{Williams:60}, with the reactor temperature fixed at 355.72~K.
	\item The decanter is modeled with the equations (30)-(36), with the decanter temperature is fixed at 310.93~K.
	\item The distillation column is modeled with the equations (38)-(48).
	\item The recycle is modeled with the equations (49)-(55).
\end{itemize}
The reactor is a continuous stirred-tank reactor (CSTR), where the following three reactions take place:
\begin{align}
	\text{A} + \text{B} & \overset{k_1}{\longrightarrow}   \text{C},           &  k_1 = A_{1}e^{-B_1/355.72} \label{eq:reaction1} \\
	\text{C} + \text{B} & \overset{k_2}{\longrightarrow}   \text{P} + \text{E},&  k_2 = A_{2}e^{-B_2/355.72} \\
	\text{P} + \text{C} & \overset{k_3}{\longrightarrow}   \text{G},           &  k_3 = A_{3}e^{-B_3/355.72} 
\end{align}
The plant inputs are the inflows $\{F_A,F_B\}$ of materials \{A,B\} and the recycle ratio $\alpha$. Among those inputs: 
\vspace{-\topsep}
\begin{itemize}[noitemsep]
	\item $F_A$ is a mean flowrate $\overline{F_A}=6.577$ t/h affected by a measured disturbance ($d^{\prime}$) which depends on systems upstream of this process and whose value belongs to the domain $ [-1000,1000]$ kg/h and is subject to step type changes about every five days. 
	\item $F_B$ is the variable that can be manipulated ($u$) whose value must be chosen by the ASP. 
	\item $\alpha$ is an unmeasured disturbance ($d^{\prime\prime}$) which can be handled by the engineers operating the plant, e.g., during maintenance phases. When the plant is operating at its nominal state $\alpha = 0.5499$. However, during maintenance phases engineers may manipulate the recycling  valve and this value may change for short or long periods (depending on the duration of the maintenance). 
\end{itemize}
The outputs of the plant $\bm{y}^{\rm T} := [F_G, \ F_P, \ F_D, \ F_A]$ correspond to the mass flows of product, waste, by-product (fuel), and material A, respectively. The functions linking the inputs to the outputs are named as follows: 
\begin{align}
	F_P := \ & f_{p(1)}(u,\bm{d}), & F_G := \ & f_{p(2)}(u,\bm{d}), &	F_D := \ & f_{p(3)}(u,\bm{d}),
\end{align}
where there is no function for $F_A$ since it is a measured disturbances. 
Material G is insoluble in the reactant stream, and it is assumed that the waste G is perfectly extracted through the mass Each plant output measurement is subject to measurement uncertainty corresponding to an unbiased normal noise whose variance is 5 kg/h (which means that the error on \textit{one measurement} of a flowrate is with an expectancy of more than 95\% within the domain [-10,10] kg/h). This measurement uncertainty is passively managed by the collection of several measures whenever a steady-state (SS) is reached. More precisely, one always waits 30 min when a SS is detected in order to collect 60 measures of the same flowrate to reduce the measurement uncertainty (a measurement of a flowrate is available each 30s).

The model corresponds to the plant where these two following structural plant-model mismatches are added:
\vspace{-\topsep}
\begin{itemize}[noitemsep]
	\item The reactor is modeled with an approximation of the three-reaction scheme of the plant by a two-reaction scheme \cite{Marchetti:09}. The two reactions of the model reads:
	\begin{align}
		\text{A} + 2\text{B}           & \overset{k_1^{*}}{\longrightarrow}   \text{P} + \text{E},  &  k_1^{*} = A_{1}^{*}e^{-B_1^{*}/355.72} \\
		\text{A} + \text{B} + \text{P} & \overset{k_2^{*}}{\longrightarrow}   \text{P} + \text{E},  &  k_2^{*} = A_{2}^{*}e^{-B_2^{*}/355.72}
	\end{align}
	\item The distillation column has a different efficiency. In fact, the efficiency $\eta_c$ of the distillation column affects the production flowrate of material P according to the following equation:
	\begin{equation} \label{eq:ColumnEff}
		F_P = (X_P^M - \eta_c X_E^M)F_M
	\end{equation}
	where $F_M$ is the total mass inflow in the distillation column and $X_P^M$ and $X_E^M$ are the weight ratios of species P and E of the flow $F_M$,  respectively.
\end{itemize}

The operating cost in \$/h is given by the following function (where the flowrates should be in ton (metric) per hour (t/h)):
\begin{align}
	\phi(u,\bm{y}) := \ &  1000(0.02F_A + 0.03F_B + 0.01 F_G- 0.3F_P - 0.0068F_D)/0.45,\\
	               := \ &  1000(0.02y_{(4)} + 0.03u+ 0.01 y_{(1)}- 0.3y_{(2)}- 0.0068y_{(3)}/0.45.
\end{align}

The operating constraints are: 
\vspace{-\topsep}
\begin{itemize}[noitemsep]
	\item Boundaries on the manipulated variables: $u\in[10,20]$ t/h,
	\item An upper bound on the production rate characterizing the maximum production that the market can absorb: $F_P \leq 2.37$ t/h. It should be noted that this is not a physical constraint that should never be exceeded. The objective of a good management of the plant is to maintain, \textit{in the long term}, the plant below this production limit. Violating this constraint would imply a storage of the production which is not necessarily problematic in the short term. 
\end{itemize}

Of course, since one considers that no measurement is perfect, the goal of the ASP is not to find the operating condition that both satisfies the constraints and minimizes the operating cost. The objective is rather to find a satisfactory operating condition which is defined by the following two parameters: 
\begin{align}
	\Delta u_{scal}    := \ & 0.01 \ \text{t/h},  & 
	loss_{acc} := \ & 0.1  \ \text{\$/h}. \label{eq:6_CaseStudy_satisfaction_parameters}
\end{align}
Basically, one considers that an operating point is satisfactory if a (feasible) modification of the manipulable inputs  $\pm \Delta u_{scal}$ can imply an improvement of the operating conditions of less than  $loss_{acc}$. Notice that  $loss_{acc}$ is rather small compared to the value of the operating cost magnitude ($\sim$400 \$/h, i.e. an improvement of 0.1 corresponds to an improvement of 0.025\%). In actual implementation, a larger $loss_{acc}$ would be sufficient, but since here one wants ``nice'' and ``clean'' results, one voluntarily chooses a very small value.

The Figures~\ref{fig:6_10_WO_Maps} a) and b) show what are the operating cost (the contour lines), the feasibility range (the red domain violates the production's constraints), and the position of the  optimal $u$ w.r.t. measured disturbance $d^{\prime}$ (the green line) in nominal operation ($d^{\prime\prime}=\alpha=0.5499$) for the plant and the model. 
It can be observed that for the high values of $d^{\prime}$ the optimum of the plant activates the constraint on the production and thus be defined by it, whereas for low values of $d^{\prime}$ the optimum of the plant does not activate any constraint. 
Figures~\ref{fig:6_10_WO_Maps} c)--k) show what functions  $f_{p(1)}$,.., $f_{p(3)}$, $f_{(1)}$, .., $f_{(3)}$ are as well as what the plant-model mismatches are (which are clearly not linear). Also, on the Figure~\ref{fig:6_10_WO_Maps}  e)--k) it is clearly visible that the plant model mismatch is non-linear. 

Finally, Table~\ref{tab:FixedParam} summarizes all the plant and model's parameters and gives the hyperparameters that are used in all simulation. 

\clearpage

\noindent
\begin{minipage}[h]{\linewidth}
	\vspace*{0pt}
	{\centering	
		\noindent
		\begin{minipage}[h]{8cm}
			\arraycolsep=1pt\def\arraystretch{1}
				\centering
				\captionof{table}{Plant and model parameters.}
				\begin{tabular}{lll}
					\toprule
					Parameter & Value & Unit \\
					\midrule 
					{$A_1$}      & {$5.975 \times 10^9$}     & {$h^{-1}$} \\
					{$A_2$}      & {$2.596 \times 10^{12}$}  & {$h^{-1}$} \\
					$A_3$      & $9.628 \times 10^{15}$  & $h^{-1}$ \\
					$B_1$      & $6.666 \times 10^{3}$   & $K$ \\
					$B_2$      & $8.333 \times 10^{3}$   & $K$ \\
					$B_3$      & $1.111 \times 10^{4}$   & $K$ \\
					$\eta_c$   & $0.1$                    &  - \\
					$A_1^*$    & $7.8812 \times 10^{11}$  & $h^{-1}$ \\
					$A_2^*$    & $1.5515 \times 10^{17}$  & $h^{-1}$ \\
					$B_1^*$    & $8.0776 \times 10^{3}$   & $K$ \\
					$B_2^*$    & $12.438 \times 10^{3}$   & $K$ \\
					$\eta_c^*$ & $0.105$                   & - \\
					\bottomrule
				\end{tabular}
				\label{tab:FixedParam}
		\end{minipage} \hskip -0ex
		\begin{minipage}[h]{7cm}
			\centering	
				\centering
				\captionof{table}{Hyperparameters}
				\begin{tabular}{lll}
					\toprule
					Parameter & Value & Unit \\
					\midrule 
					\multirow{2}[0]{*}{$f_{(1)}$}       & {$\ell_1$} & {$4$ t/h} \\
					& {$\ell_2$} & {$0.5$ t/h} \\
					\midrule 
					\multirow{2}[0]{*}{$f_{(2)}$}       & {$\ell_1$} & {$4$ t/h} \\
					& {$\ell_2$} & {$0.5$ t/h}\\
					\midrule 
					\multirow{2}[0]{*}{$f_{(3)}$}       & {$\ell_1$} & {$4$ t/h} \\
					& {$\ell_2$} & {$0.5$ t/h}\\
					\bottomrule
				\end{tabular}
				\label{tab:HyperparametersParam}
		\end{minipage} 
	}
\end{minipage}

\clearpage

\noindent
\begin{minipage}[h]{\linewidth}
	\vspace*{0pt}
	{\centering	
		\begin{minipage}[h]{6.675cm}
			\centering	
			\includegraphics[width=6.675cm]{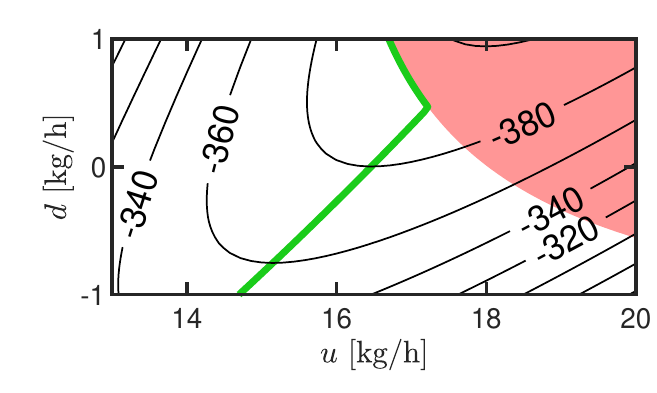} \\
			\ \ \ \ a) Plant's cost and constraint
		\end{minipage} \hskip -0ex
		\begin{minipage}[h]{6.675cm}
			\centering	
			\includegraphics[width=6.675cm]{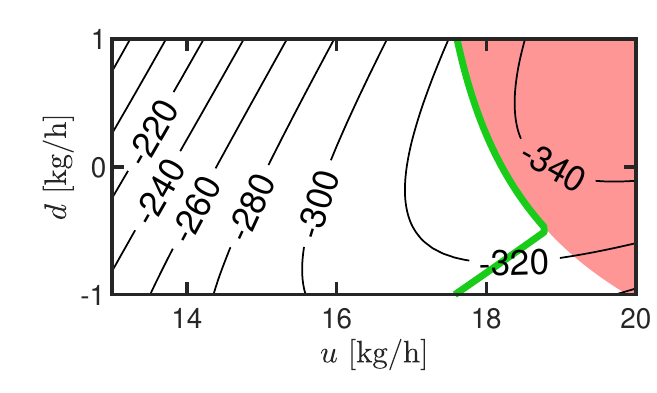}\\
			\ \ \ \ b) Model's cost and constraint
		\end{minipage} 
	
		
		\begin{minipage}[h]{4.45cm}
			\centering	
			\includegraphics[width=4.45cm]{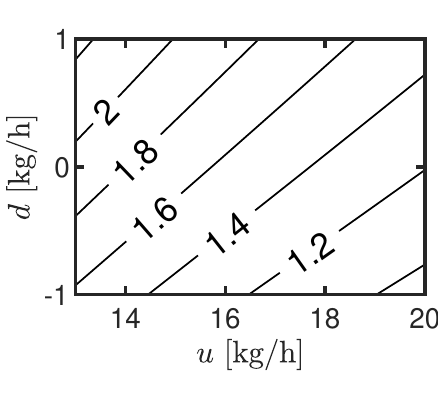} \\
			\ \ \ \ c) $f_{p(1)}(u,d)$ 
		\end{minipage} \hskip -0ex
		\begin{minipage}[h]{4.45cm}
			\centering	
			\includegraphics[width=4.45cm]{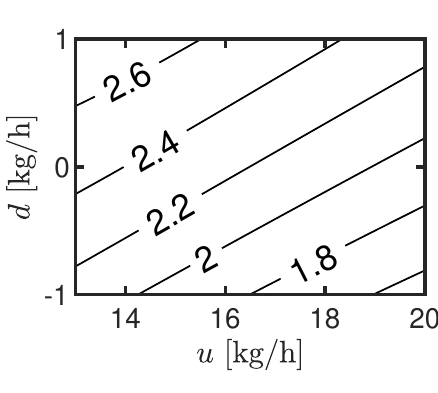} \\
			\ \ \ \ d) $f_{(1)}(u,d)$ 
		\end{minipage} \hskip -0ex
		\begin{minipage}[h]{4.45cm}
			\centering	
			\includegraphics[width=4.45cm]{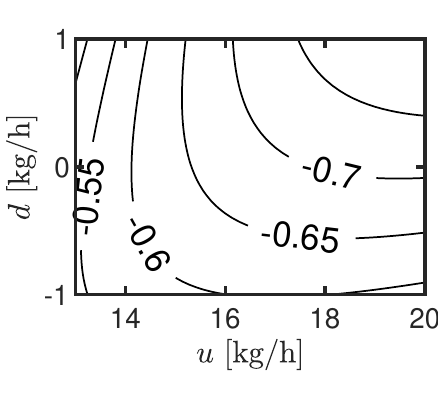} \\
			\ \  e)  $f_{p(1)}(u,d) - f_{(1)}(u,d)$ 
		\end{minipage} 
	
		
		\begin{minipage}[h]{4.45cm}
			\centering	
			\includegraphics[width=4.45cm]{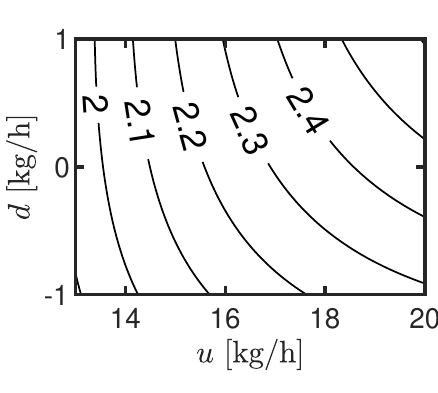} \\
			\ \ \ \ f) $f_{p(2)}(u,d)$ 
		\end{minipage} \hskip -0ex
		\begin{minipage}[h]{4.45cm}
			\centering	
			\includegraphics[width=4.45cm]{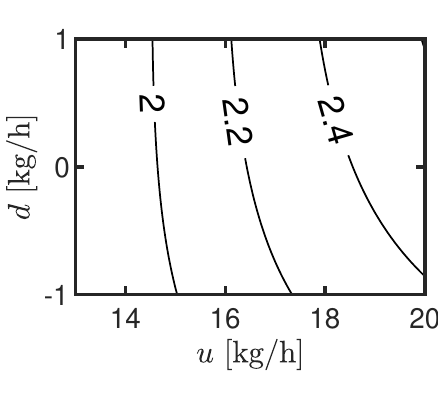} \\
			\ \ \ \ g) $f_{(2)}(u,d)$ 
		\end{minipage} \hskip -0ex
		\begin{minipage}[h]{4.45cm}
			\centering	
			\includegraphics[width=4.45cm]{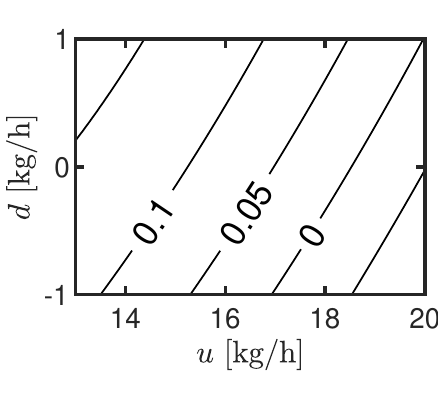} \\
			  h)  $f_{p(2)}(u,d) - f_{(2)}(u,d)$ 
		\end{minipage} 
	
		
		\begin{minipage}[h]{4.45cm}
			\centering	
			\includegraphics[width=4.45cm]{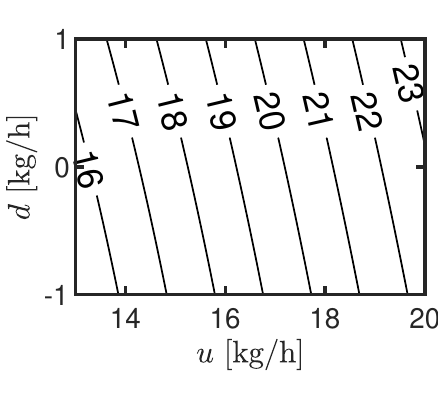} \\
			\ \ \ \ i) $f_{p(3)}(u,d)$ 
		\end{minipage} \hskip -0ex
		\begin{minipage}[h]{4.45cm}
			\centering	
			\includegraphics[width=4.45cm]{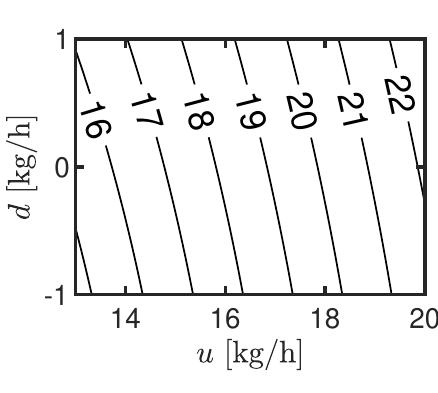} \\
			\ \ \ \ j) $f_{(3)}(u,d)$ 
		\end{minipage} \hskip -0ex
		\begin{minipage}[h]{4.45cm}
			\centering	
			\includegraphics[width=4.45cm]{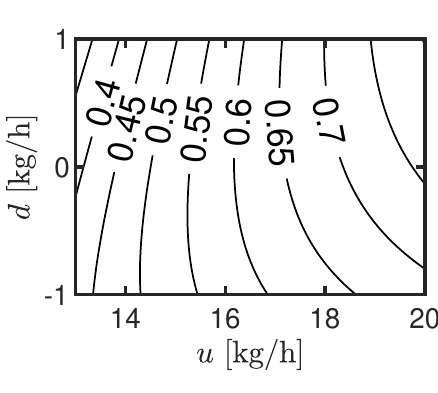} \\
			  k)  $f_{p(3)}(u,d) - f_{(3)}(u,d)$ 
		\end{minipage} 
		
	}
	\captionof{figure}{Williams-Otto process -- Problem description}
	\label{fig:6_10_WO_Maps}
\end{minipage}

\clearpage

\subsection{Simulations}

In order to illustrate all the contributions of this chapter, three studies are proposed: 
\begin{itemize}
	\item \textbf{Study 1 (A simple optimization):} This study consists in illustrating the effects of the validation system proposed in section~\ref{sec:6_7_The_improved_experiment_designer}. To do so, two scenarios are considered: 
	\begin{itemize}[noitemsep]
		\item \textbf{Scenario 1:} Standard ASP (the one proposed in this chapter) is used, i.e. the validation method, the stand-by system, etc, everything is implemented. 
		\item \textbf{Scenario 2:} Standard ASP without the validation process of the experience designer is used. Notice that without this validation the ASP cannot enter stand-by. 	
	\end{itemize}
	The simulation results obtained for those two scenarios are given in Figures~\ref{fig:6_11_ShortSim} to \ref{fig:6_14_ShortSim}. Clearly, the main observation is that without the validation procedure the ASP is unable to converge to the plant optimum. It is only capable of converging to a feasible operating point.   In addition to that, it can also be seen that ASP needs only 4 iterations to identify an operating point satisfying the optimality conditions \eqref{eq:6_CaseStudy_satisfaction_parameters}, and that after these 4 iterations the ASP is in stand-by. 
		\item \textbf{Study 2 (Reaction to unmeasured disturbances)} This study consists in illustrating the functioning of the improved couple \{decision trigger, consistency monitor\}.  To do so, three scenarios of the variation of the unmeasured disturbance are considered:
		\begin{itemize}[noitemsep]
			\item \textbf{Scenario 1:} A large step on the unmeasured constraint during a stand-by phase. The simulation results are given in Figure~\ref{fig:6_study_2_sc1__1_inputs} to \ref{fig:6_study_2_sc1__4_database}.
			\item \textbf{Scenario 2:} A very short step on the unmeasured constraint during a stand-by phase. The simulation results are given oin Figure~\ref{fig:6_study_2_sc2__1_inputs} to \ref{fig:6_study_2_sc2__3_database}.
			\item \textbf{Scenario 3:} A long step on the unmeasured constraint during the research phase (before ASP converges to an optimum). The simulation results are given in Figure~\ref{fig:6_study_2_sc3__1_inputs} to \ref{fig:6_study_2_sc3__4_database}.
		\end{itemize}
		In scenario 1, the ASP is in stand-by mode when the unmeasured disturbance is subject to a step variation. About eight hours after this step, the decision trigger observes that the plant is (a) in stand-by mode, (b) at SS, and (c) that this SS is far from the confidence domain based on the previous measurements. The union of those three observations \{(a),(b),(c)\} is sufficient to trigger the reset of the database (one keeps only the last measurement) and the interruption of stand-by mode. As one can see in Figures~\ref{fig:6_study_2_sc1__3_details} and \ref{fig:6_study_2_sc1__4_database}, the database is indeed reset. Figure~\ref{fig:6_study_2_sc1__2_decision} shows that the exploration of the plant restarts to discover the new optimal operating point.
		
		In the scenario 2, the ASP is also in stand-by mode when the unmeasured disturbance is subject to a step variation. However, less than eight hours after this event, the unmeasured disturbance is again subject to an unmeasured disturbance that brings it back to its initial value. This time, when the plan finally reaches SS, it matches with the expectation one may have based on  the previous observations. As a result, everything works as if nothing had happened, i.e. the database is not reset, the inputs are changed, and the stand-by mode is not left. 
		
		In the scenario 3, the ASP is performing experiments to identify the optimal operating condition of the plant when the unmeasured disturbance is  subject to a step variation. When the experiment during which this step occurred is finished, the consistency monitor sees that the SS of this experiment does not match with what one could expect according to the previous experimental results (see Figure~\ref{fig:6_study_2_sc3__3_details}, at iteration 4 the new observation is far from the confidence domain of the iteration 3). Therefore, it is automatically decided to reset the database (see Figure~\ref{fig:6_study_2_sc3__4_database}) and to make new decisions to restart the identification of the plant's optimum (see Figure~\ref{fig:6_study_2_sc3__2_decision}). One can clearly see in Figure~\ref{fig:6_study_2_sc3__1_inputs} that this event is perfectly managed and that ASP is able to efficiently reach the plant optimum in four iterations after this event.

		
		\item \textbf{Study 3 (A long simulation):} This study consists in illustrating the ability of the data management features to sustain long term implementations. To do so, the WO process is simulated over a year (8760h) without unmeasured disturbances variations (so that the database is not reset to maximize its size). The simulation results are shown on Figures~\ref{fig:6_study_3__1_inputs} to \ref{fig:6_study_3__6_computationnal}. 
		
		A first observation could be that the ASP is able to maintain the WO process at its optimal operating conditions almost all year, see Figure~\ref{fig:6_study_3__1_inputs}.  
		
		A second observation could be that the ASP seems to present two phases: a training phase and a trained phase. The training phase would correspond to approximately the first 1000h during which  the ASP needs several iterations to identify a plant optimum. This is clearly visible on Figure~\ref{fig:6_study_3__3_training} where several oscillations around the optimal cost are observed $\forall t < 1000$ and are not observed anymore after. It can also be seen on Figure~\ref{fig:6_study_3__4_decisions} where the decisions and type w.r.t the time are shown for two time periods. One corresponding to the first hours of the plant under the ASP's control, and the other one   corresponding to what happens after 4000h. One can see that initially ASP needs several decisions, i.e. it needs to explore to identify the plant optimum, whereas after many hours it just go to what it believes is the plant optimum (first decision) and decides to stay there without requiring any validation (second decision). Hence the couples of red points on the lower plot of Figure~\ref{fig:6_study_3__4_decisions}. 
		
		A third observation could be that the database size stays relatively small. There are two reasons for that. First, the stand-by mode enables to limit the quantity of data to be used. For instance, if iterative method such as KMFCaA where used, then one would never stop collecting data and the  x-axis of Figure~\ref{fig:6_study_3__5_training}   would be much longer, and therefore the black line would get much higher.  Second, the compression method is efficient enough to combine similar experiments and extract from them only the ``useful'' information.  
		
		Finally, Figure~\ref{fig:6_study_3__6_computationnal} shows that the computational time of ASP does not increase significantly over time. This computational time gathers everything the ASP does at each iteration.	(\textit{It also take into account some ``data management'' related to the construction of the plots presented in this study. Also, since the code is written in MATLAB (which is clearly not the most efficient language if one seeks speed) and the code is absolutely not optimized for speed. So, here what one should observe is rather the trend than the value, and the trend is that the computational time does not increase significantly over time.}) 
\end{itemize}

\clearpage 

\noindent
\begin{minipage}[h]{\linewidth}
\begin{minipage}[h]{\linewidth}
	\vspace*{0pt}
	{\centering	
		
		\includegraphics[width=6.675cm]{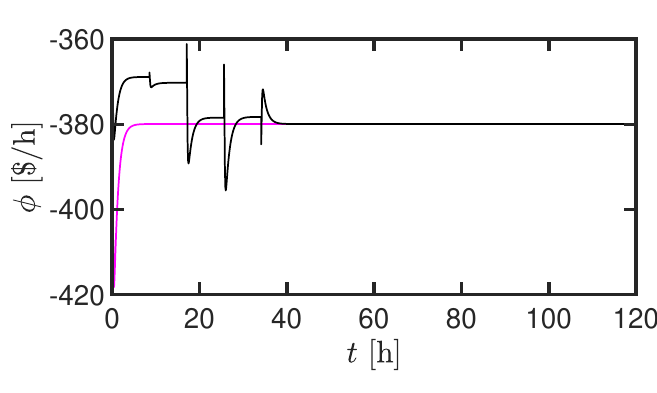}    \hskip -0ex
		\includegraphics[width=6.675cm]{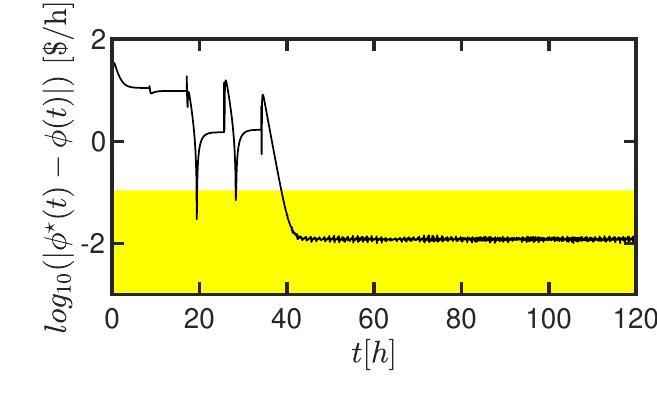} \hskip -0ex
		\includegraphics[width=6.675cm]{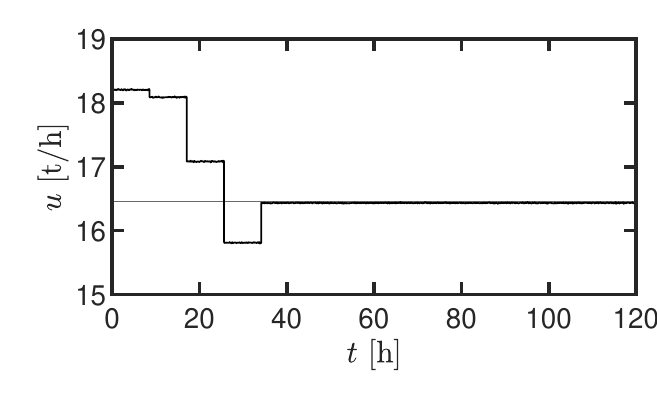}   \hskip -0ex
		\includegraphics[width=6.675cm]{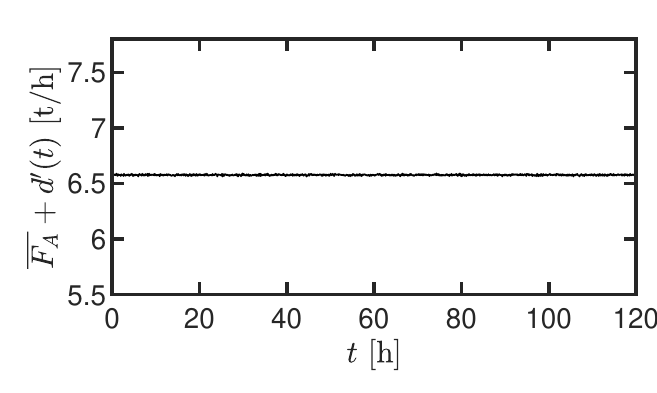} \hskip -0ex
		\includegraphics[width=6.675cm]{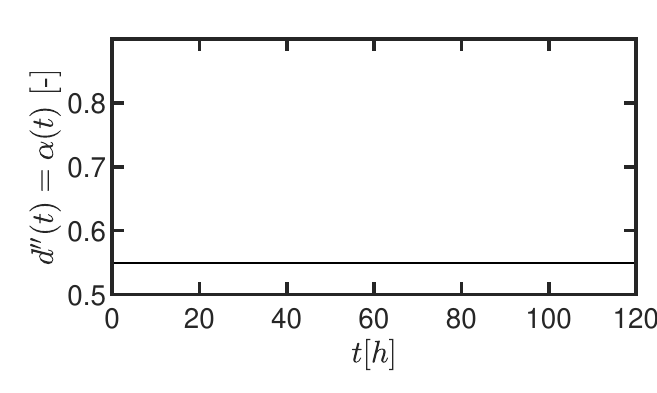}   \hskip -0ex
		\includegraphics[width=6.675cm]{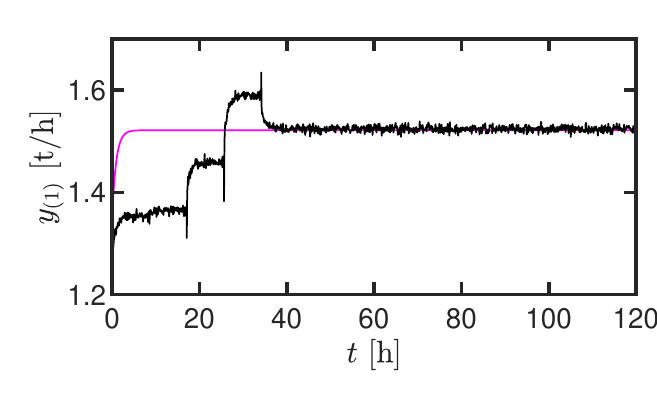}   \hskip -0ex
		\includegraphics[width=6.675cm]{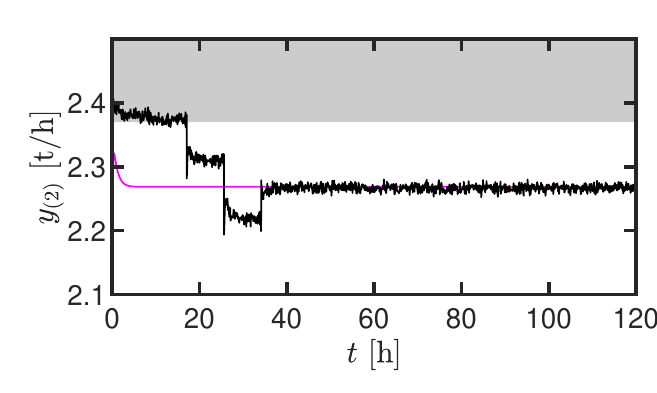} \hskip -0ex
		\includegraphics[width=6.675cm]{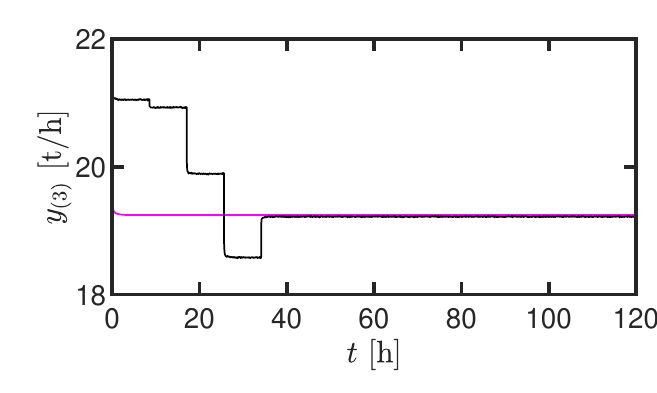}   \hskip -0ex
		
	}
	
	\textcolor{black}{\raisebox{1mm}{\rule{0.5cm}{0.05cm}}}: results obtained with the model and inaccurate measurements, \\
	\textcolor{magenta}{\raisebox{1mm}{\rule{0.5cm}{0.05cm}}}: results obtained with perfect model and measurements,\\
	\textcolor{gray}{\raisebox{-0.5mm}{\rule{0.5cm}{0.3cm}}}: constrained area, 
	\textcolor{yellow}{\raisebox{-0.5mm}{\rule{0.5cm}{0.3cm}}}: acceptable loss.
	\captionof{figure}{\textbf{Study 1 -- Scenario 1:} ASP is able to converge to an optimum not defined by constraints}
	\label{fig:6_11_ShortSim}
\end{minipage} 

\noindent
\begin{minipage}[h]{\linewidth}
	\vspace*{0pt}
	{\centering	
		
		\includegraphics[trim={1.5cm 0cm 1cm  0cm},clip,width=13.35cm]{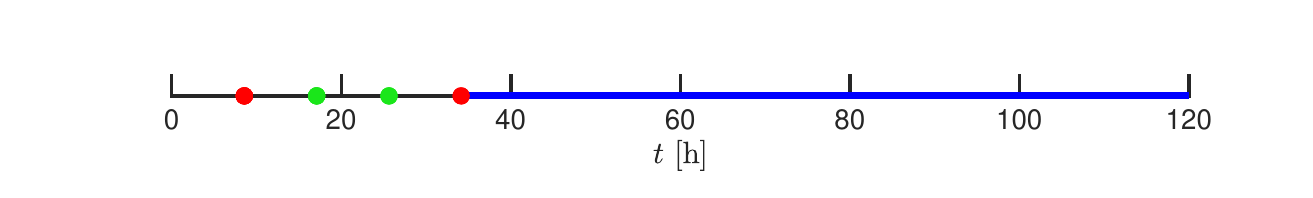}   
		
	}
	\textcolor{red}{$\bullet$}: ``normal'' decision, 
	\textcolor{green}{$\bullet$}: decision related to the validation, 
	\textcolor{blue}{\raisebox{1mm}{\rule{0.5cm}{0.05cm}}}: stand-by mode.  
	\captionof{figure}{\textbf{Study 1 -- Scenario 1:} ASP Decision dates \& stand-by mode}
	\label{fig:6_12_ShortSim}
\end{minipage} 
\end{minipage}

\noindent
\begin{minipage}[h]{\linewidth}
\begin{minipage}[h]{\linewidth}
	\vspace*{0pt}
	{\centering	
		
		\includegraphics[width=6.675cm]{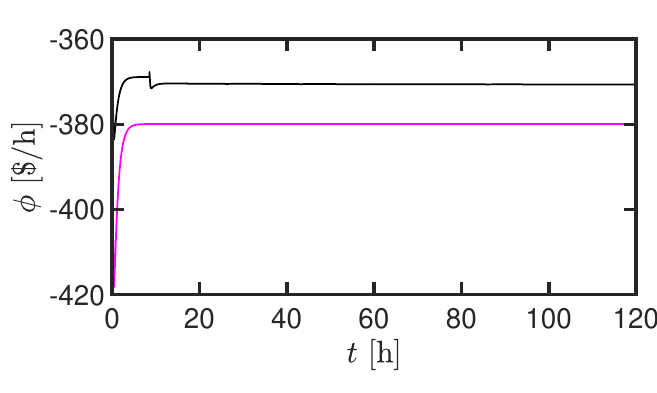}    \hskip -0ex
		\includegraphics[width=6.675cm]{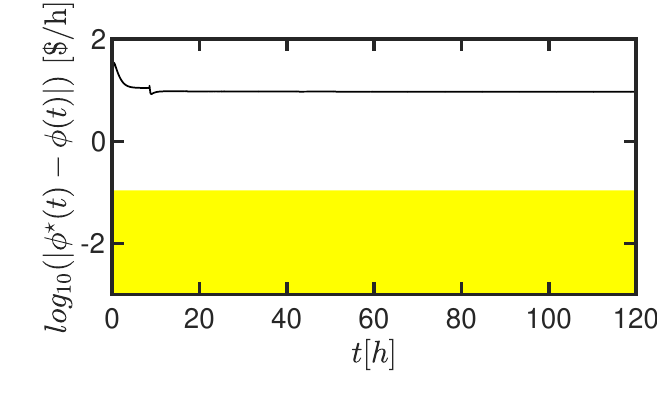} \hskip -0ex
		\includegraphics[width=6.675cm]{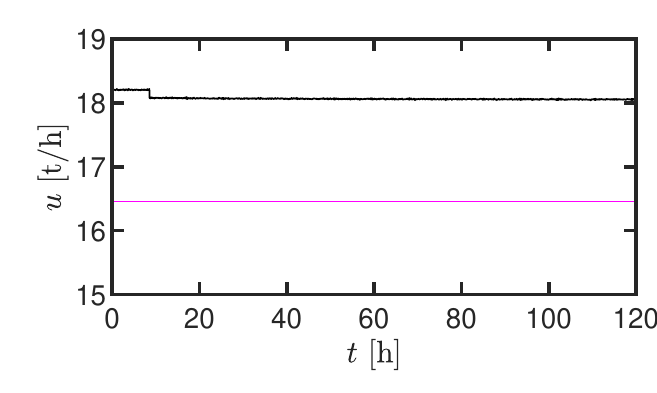}   \hskip -0ex
		\includegraphics[width=6.675cm]{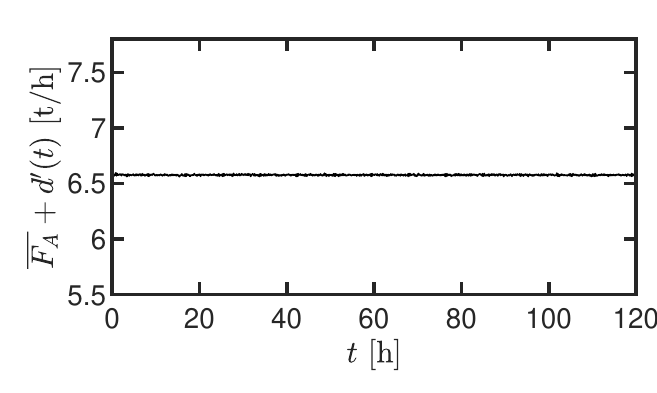} \hskip -0ex
		\includegraphics[width=6.675cm]{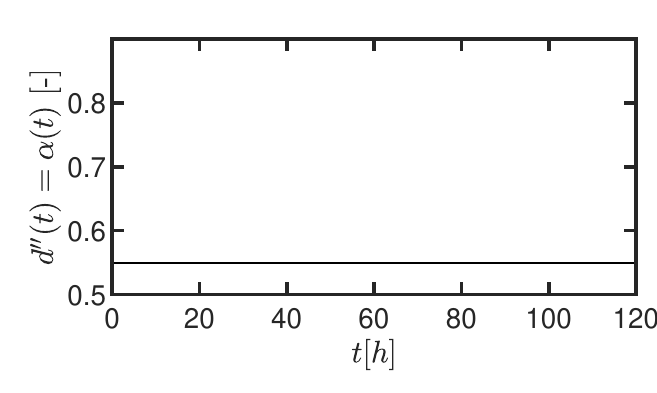}   \hskip -0ex
		\includegraphics[width=6.675cm]{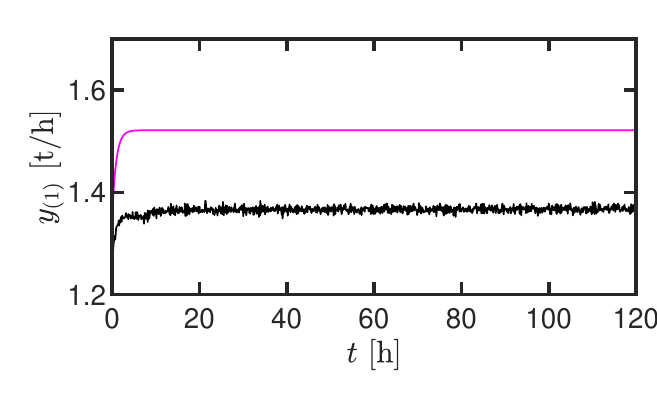}   \hskip -0ex
		\includegraphics[width=6.675cm]{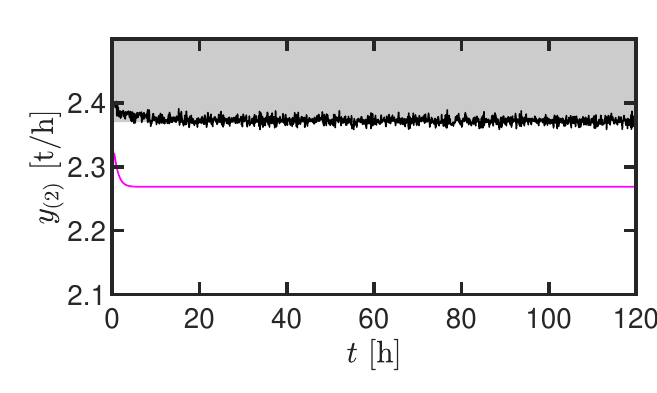} \hskip -0ex
		\includegraphics[width=6.675cm]{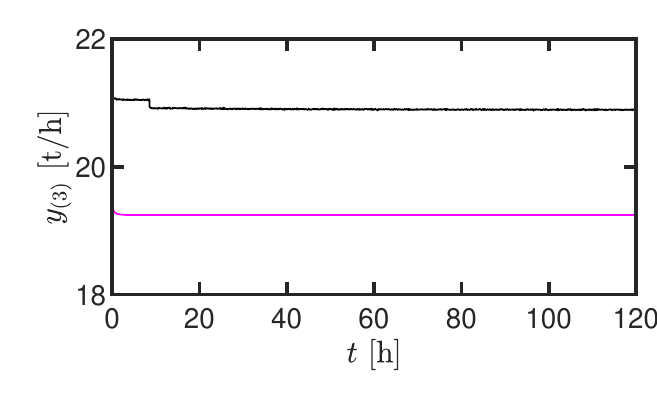}   \hskip -0ex
		
	}
	\textcolor{black}{\raisebox{1mm}{\rule{0.5cm}{0.05cm}}}: results obtained with the model and inaccurate measurements, \\
	\textcolor{magenta}{\raisebox{1mm}{\rule{0.5cm}{0.05cm}}}: results obtained with perfect model and measurements,\\
	\textcolor{gray}{\raisebox{-0.5mm}{\rule{0.5cm}{0.3cm}}}: constrained area, 
	\textcolor{yellow}{\raisebox{-0.5mm}{\rule{0.5cm}{0.3cm}}}: acceptable loss.
	\captionof{figure}{\textbf{Study 1 -- Scenario 2:} Without proper validation the ASP is \textbf{not} able to converge to an optimum not defined by constraints}
	\label{fig:6_13_ShortSimes}
\end{minipage}

\noindent
\begin{minipage}[h]{\linewidth}
	\vspace*{0pt}
	{\centering	
		
		\includegraphics[trim={1.5cm 0cm 1cm  0cm},clip,width=13.35cm]{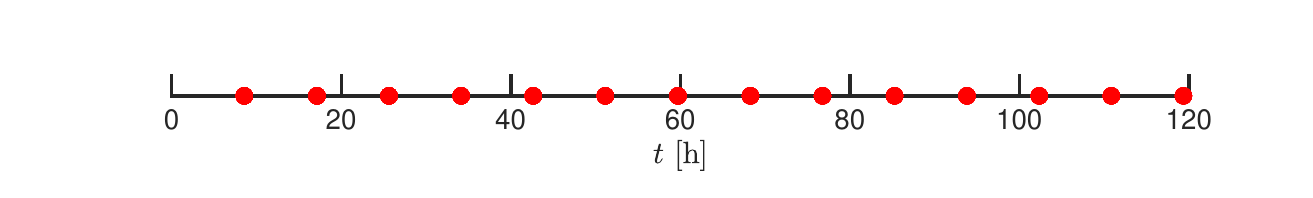}   
		
	}
	\textcolor{red}{$\bullet$}: ``normal'' decision, 
	\textcolor{green}{$\bullet$}: decision related to the validation, 
	\textcolor{blue}{\raisebox{1mm}{\rule{0.5cm}{0.05cm}}}: stand-by mode.  
	\captionof{figure}{\textbf{Study 1 -- Scenario 2:} ASP-without-validation decision dates}
	\label{fig:6_14_ShortSim}
\end{minipage} 
\end{minipage}

\noindent
\begin{minipage}[h]{\linewidth}
	\begingroup
	\fontsize{10pt}{12pt}\selectfont
	\vspace*{0pt}
	{\centering	
		\begin{minipage}[t]{4.3cm}
			\vspace{0pt}
			\includegraphics[trim={13.5cm 4.5cm 5.8cm  0.4cm},clip,height=4.5cm]{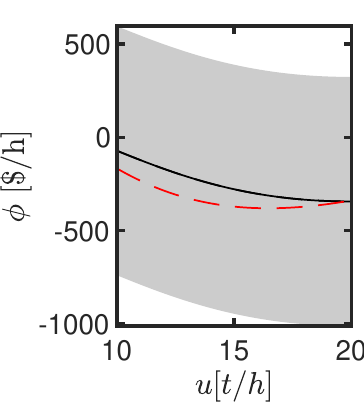}
		\end{minipage} \hskip -0ex  
		\begin{minipage}[t]{3cm}
			\vspace{0pt}
			\begin{overpic}[trim={14.65cm 4.5cm 5.8cm  0.4cm},clip,height=4.5cm]{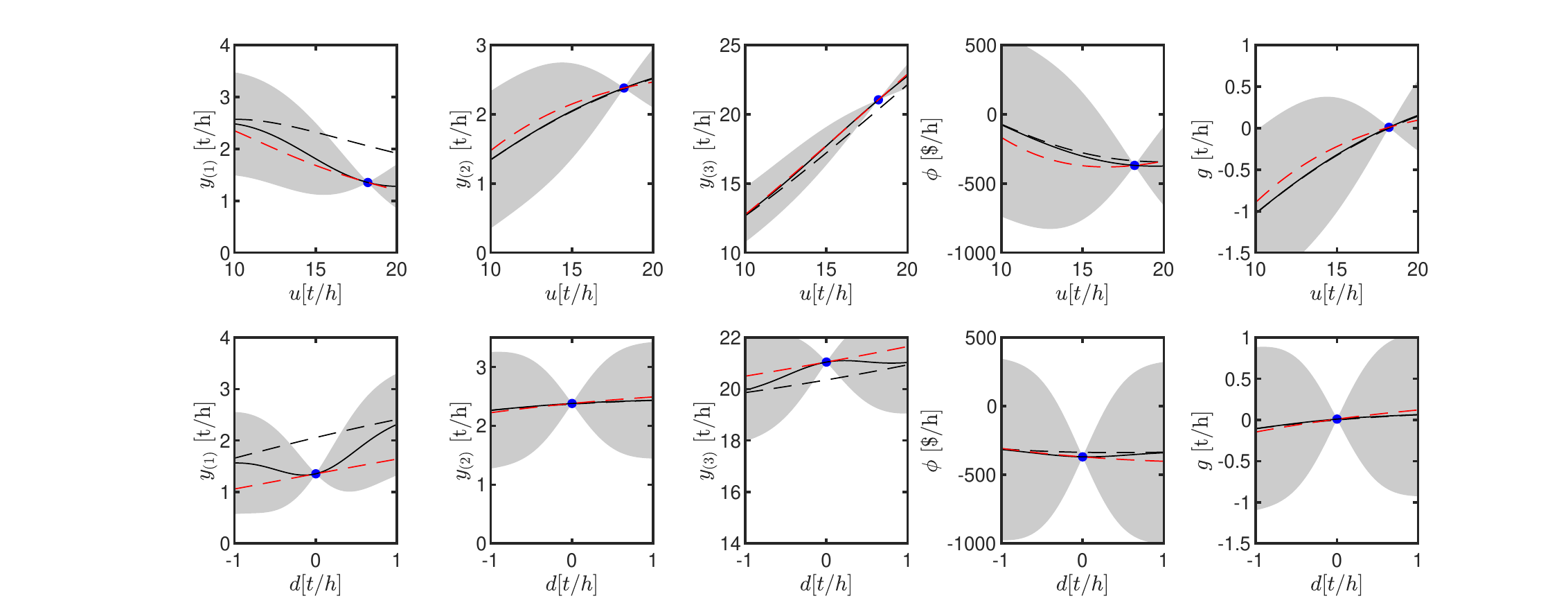}
				\put (-2,10.6) {\textcolor{black}{\colorbox{white}{10}}}
				\put (24,10.6) {\textcolor{black}{\colorbox{white}{15}}}
				\put (50,10.6) {\textcolor{black}{\colorbox{white}{20}}}
			\end{overpic}
		\end{minipage} \hskip -0ex   
		\begin{minipage}[t]{3cm}
			\vspace{0pt}
			\begin{overpic}[trim={14.65cm 4.5cm 5.8cm  0.4cm},clip,height=4.5cm]{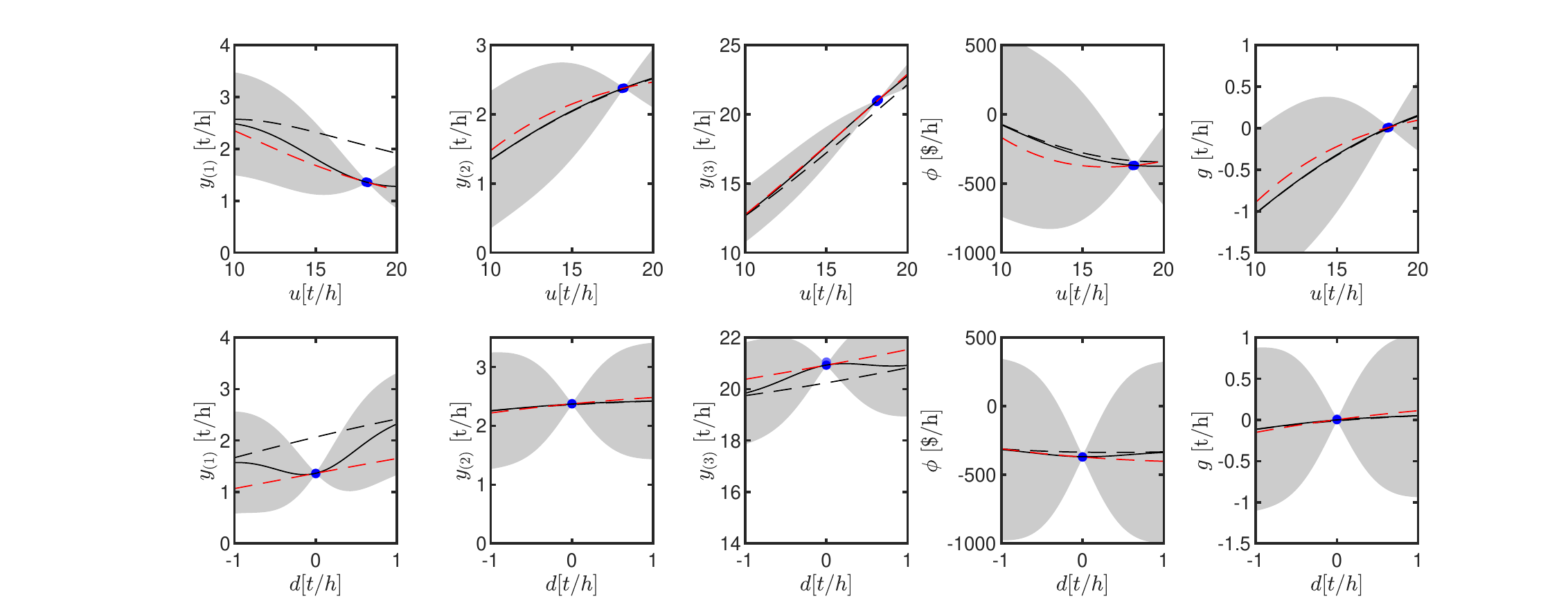}
				\put (-2,10.6) {\textcolor{black}{\colorbox{white}{10}}}
				\put (24,10.6) {\textcolor{black}{\colorbox{white}{15}}}
				\put (50,10.6) {\textcolor{black}{\colorbox{white}{20}}}
			\end{overpic}
		\end{minipage} \hskip -0ex   
		\begin{minipage}[t]{3cm}
			\vspace{0pt}
			\begin{overpic}[trim={14.65cm 4.5cm 5.8cm  0.4cm},clip,height=4.5cm]{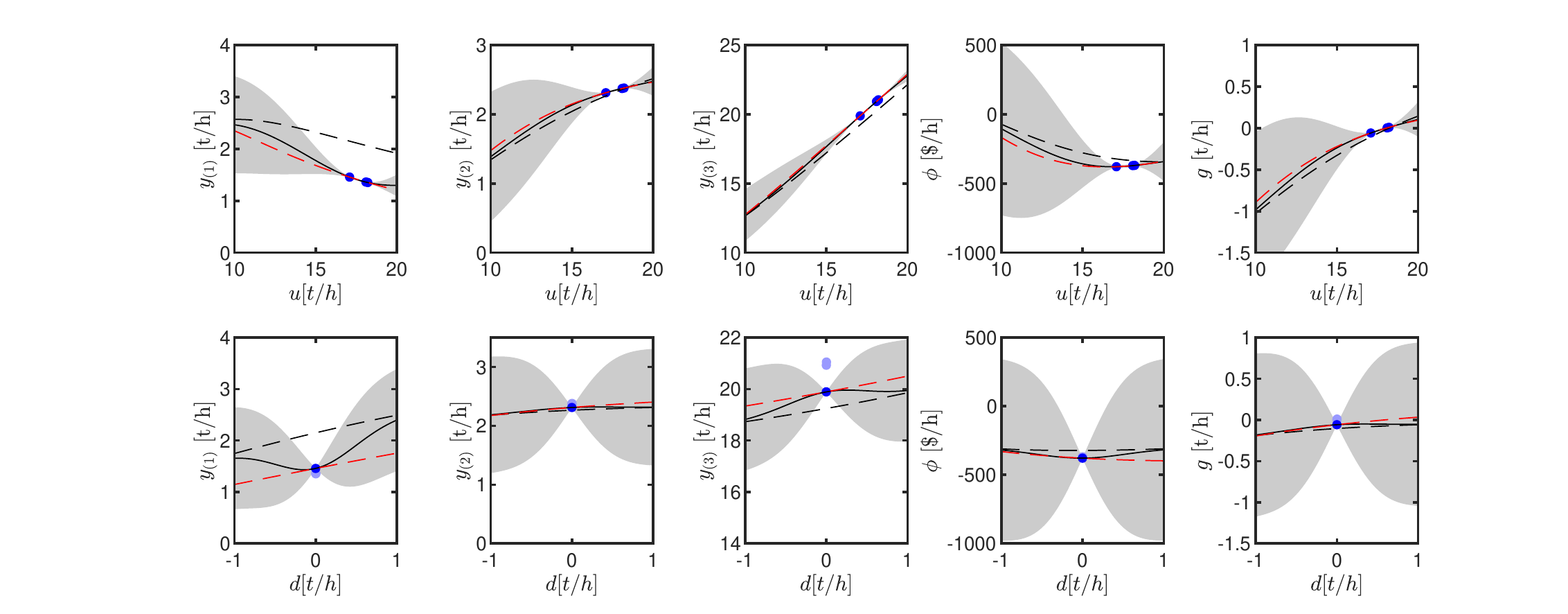}
				\put (-2,10.6) {\textcolor{black}{\colorbox{white}{10}}}
				\put (24,10.6) {\textcolor{black}{\colorbox{white}{15}}}
				\put (50,10.6) {\textcolor{black}{\colorbox{white}{20}}}
			\end{overpic}
		\end{minipage} \hskip -0ex  
		
		\begin{minipage}[t]{4.3cm}
			\vspace{0pt}
			\includegraphics[trim={13.5cm 4.5cm 5.8cm  0.4cm},clip,height=4.5cm]{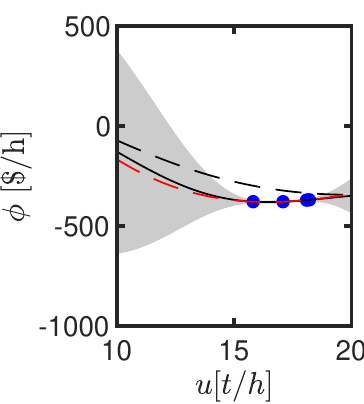}
		\end{minipage} \hskip -0ex  
		\begin{minipage}[t]{3cm}
			\vspace{0pt}
			\begin{overpic}[trim={14.65cm 4.5cm 5.8cm  0.4cm},clip,height=4.5cm]{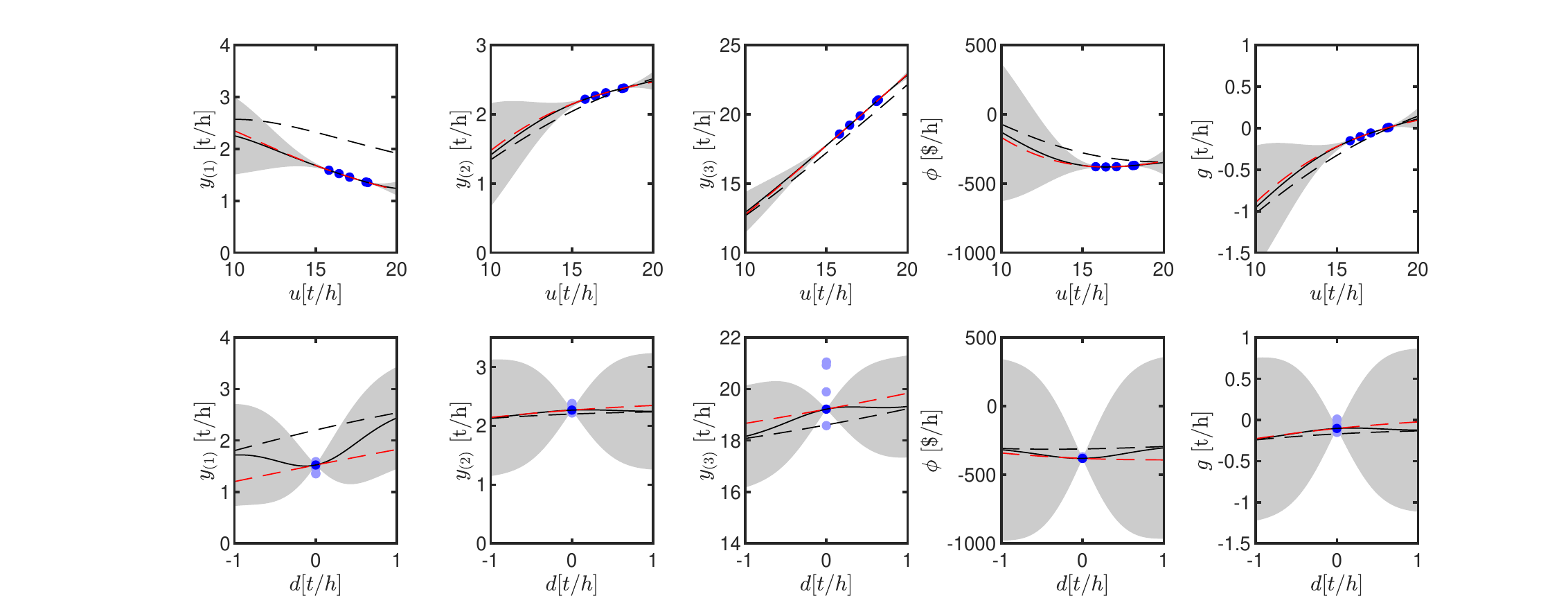}
				\put (-2,10.6) {\textcolor{black}{\colorbox{white}{10}}}
				\put (24,10.6) {\textcolor{black}{\colorbox{white}{15}}}
				\put (50,10.6) {\textcolor{black}{\colorbox{white}{20}}}
			\end{overpic}
		\end{minipage} \hskip -0ex   
		\begin{minipage}[t]{3cm}
			\vspace{0pt} 
			\begin{overpic}[trim={14.65cm 4.5cm 5.8cm  0.4cm},clip,height=4.5cm]{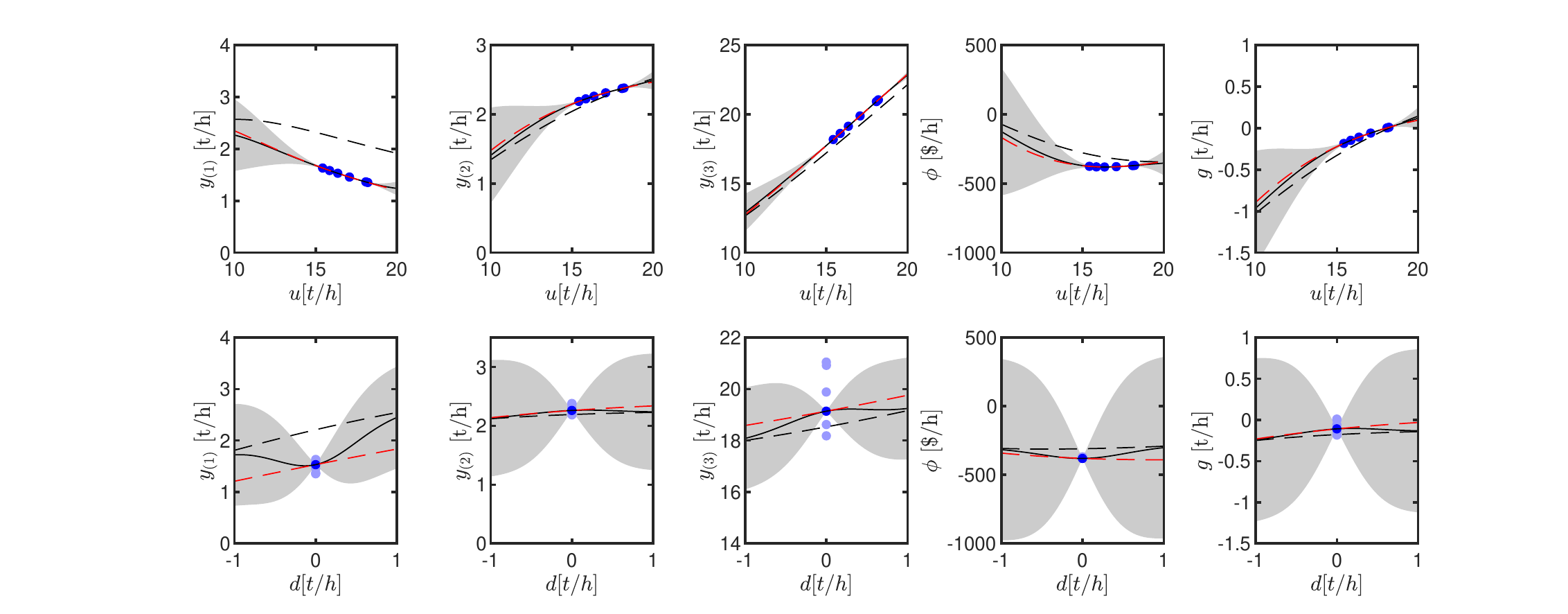}
				\put (-2,10.6) {\textcolor{black}{\colorbox{white}{10}}}
				\put (24,10.6) {\textcolor{black}{\colorbox{white}{15}}}
				\put (50,10.6) {\textcolor{black}{\colorbox{white}{20}}}
			\end{overpic}
		\end{minipage} \hskip -0ex   
		
	}
	\textcolor{blue}{$\bullet$}: data, 
	\textcolor{gray}{\raisebox{-0.5mm}{\rule{0.5cm}{0.3cm}}}: uncertainty on $\phi$, 
	\textcolor{black}{\raisebox{1mm}{\rule{0.2cm}{0.05cm}\hspace{0.1cm}\rule{0.2cm}{0.05cm}}}: model,
	\textcolor{red}{\raisebox{1mm}{\rule{0.2cm}{0.05cm}\hspace{0.1cm}\rule{0.2cm}{0.05cm}}}: plant,
	\textcolor{black}{\raisebox{1mm}{\rule{0.5cm}{0.05cm}}}: updated model.
	\endgroup
	\captionof{figure}{\textbf{Study 1 -- Scenario 1:} Some details about the iterations}
	\label{fig:6_13_ShortSim}
%
%
%
%
	\begingroup
	\fontsize{10pt}{12pt}\selectfont
	\vspace*{0pt}
	{\centering	
		\begin{minipage}[t]{4.3cm}
			\vspace{0pt}
			\includegraphics[trim={13.5cm 4.5cm 5.8cm  0.4cm},clip,height=4.5cm]{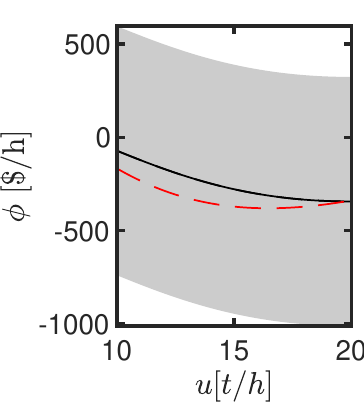}
		\end{minipage} \hskip -0ex  
		\begin{minipage}[t]{3cm}
			\vspace{0pt}
			\begin{overpic}[trim={14.65cm 4.5cm 5.8cm  0.4cm},clip,height=4.5cm]{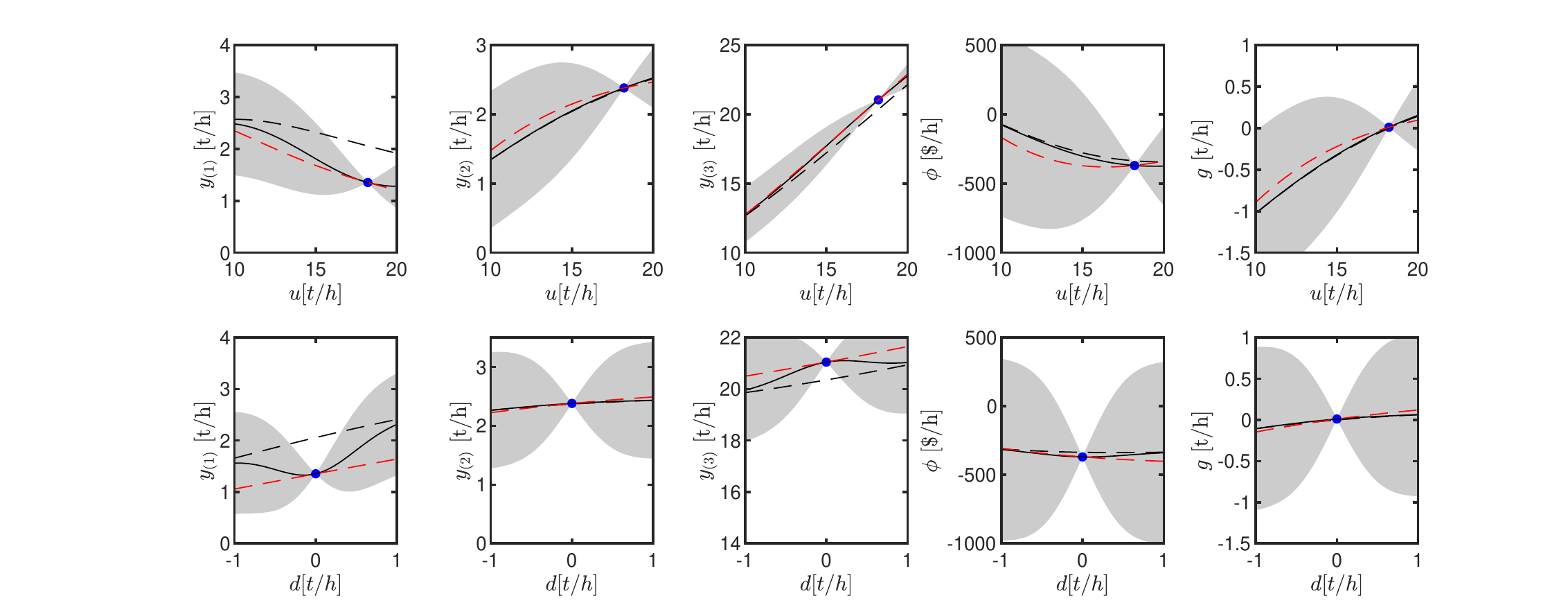}
				\put (-2,10.6) {\textcolor{black}{\colorbox{white}{10}}}
				\put (24,10.6) {\textcolor{black}{\colorbox{white}{15}}}
				\put (50,10.6) {\textcolor{black}{\colorbox{white}{20}}}
			\end{overpic}
		\end{minipage} \hskip -0ex   
		\begin{minipage}[t]{3cm}
			\vspace{0pt}
			\begin{overpic}[trim={14.65cm 4.5cm 5.8cm  0.4cm},clip,height=4.5cm]{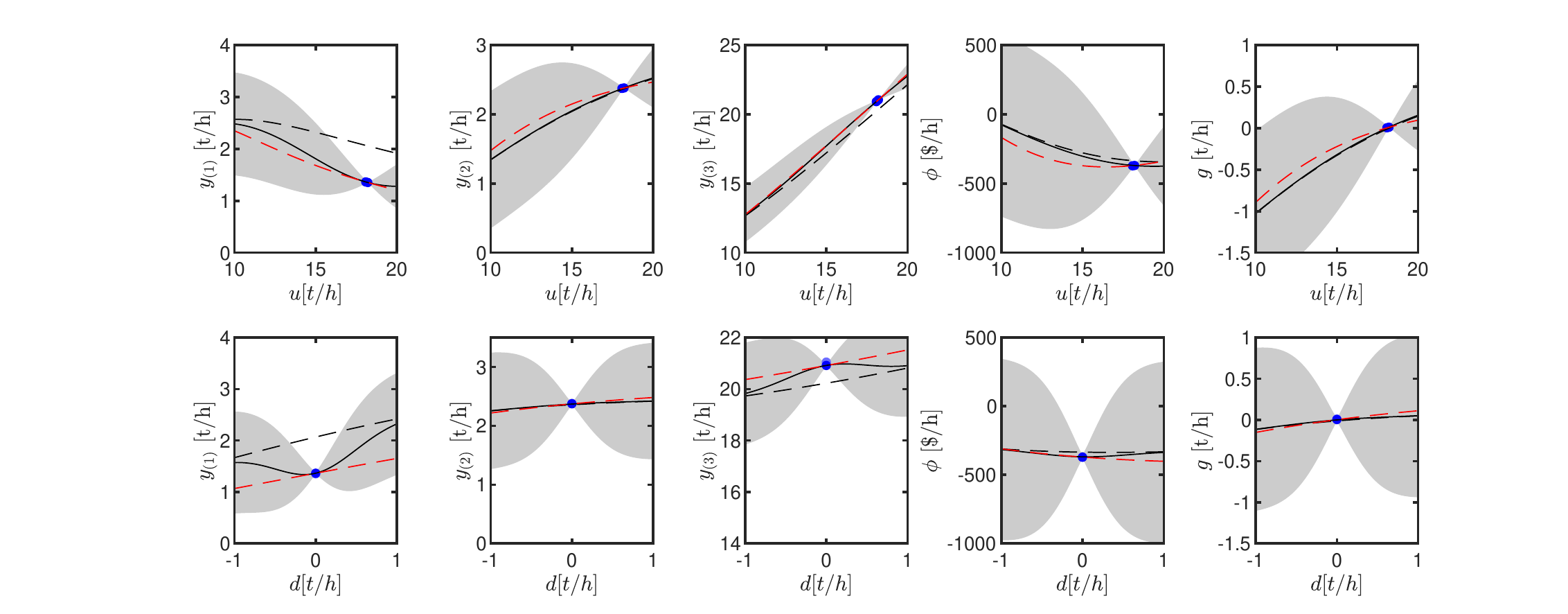}
				\put (-2,10.6) {\textcolor{black}{\colorbox{white}{10}}}
				\put (24,10.6) {\textcolor{black}{\colorbox{white}{15}}}
				\put (50,10.6) {\textcolor{black}{\colorbox{white}{20}}}
			\end{overpic}
		\end{minipage} \hskip -0ex   
		\begin{minipage}[t]{3cm}
			\vspace{0pt}
			\begin{overpic}[trim={14.65cm 4.5cm 5.8cm  0.4cm},clip,height=4.5cm]{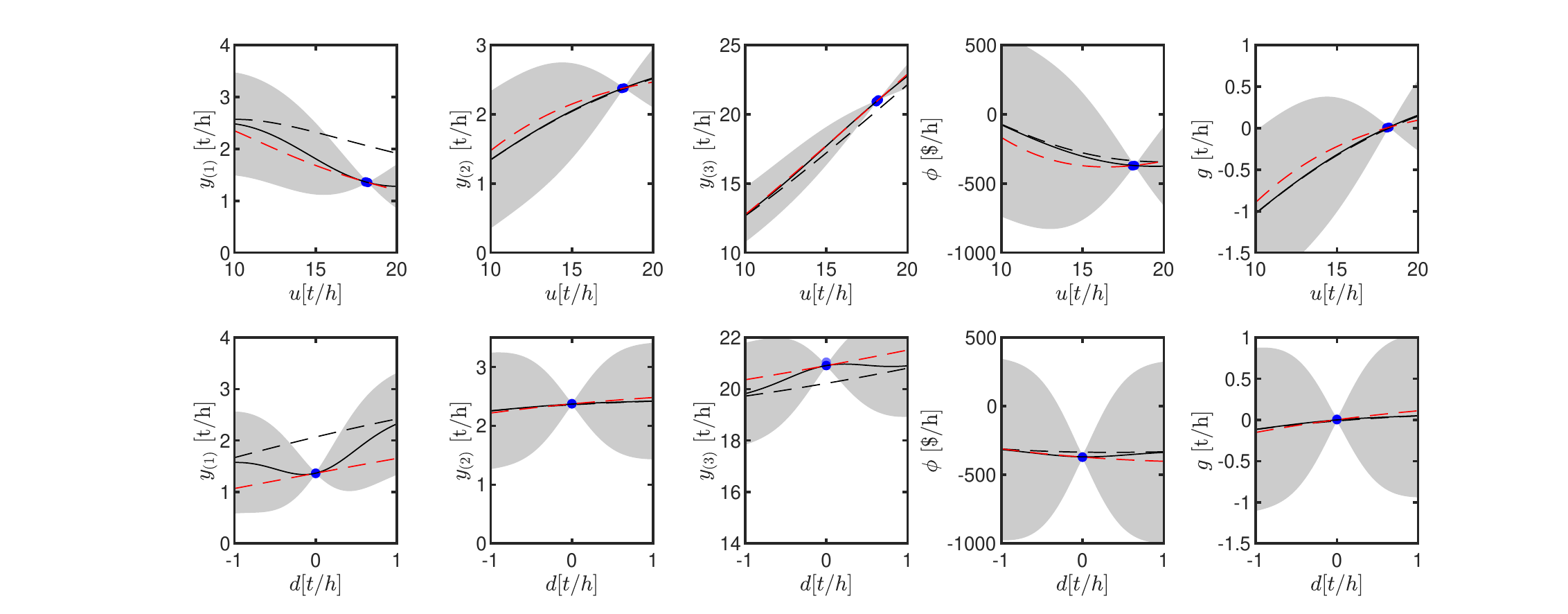}
				\put (-2,10.6) {\textcolor{black}{\colorbox{white}{10}}}
				\put (24,10.6) {\textcolor{black}{\colorbox{white}{15}}}
				\put (50,10.6) {\textcolor{black}{\colorbox{white}{20}}}
			\end{overpic}
		\end{minipage} \hskip -0ex  
		
	}
	
	\textcolor{blue}{$\bullet$}: data, 
	\textcolor{gray}{\raisebox{-0.5mm}{\rule{0.5cm}{0.3cm}}}: uncertainty on $\phi$, 
	\textcolor{black}{\raisebox{1mm}{\rule{0.2cm}{0.05cm}\hspace{0.1cm}\rule{0.2cm}{0.05cm}}}: model,
	\textcolor{red}{\raisebox{1mm}{\rule{0.2cm}{0.05cm}\hspace{0.1cm}\rule{0.2cm}{0.05cm}}}: plant,
	\textcolor{black}{\raisebox{1mm}{\rule{0.5cm}{0.05cm}}}: updated model.
	\endgroup
	\captionof{figure}{\textbf{Study 1 -- Scenario 2:} Some details about the iterations}
	\label{fig:6_13_ShortSim3456y}
\end{minipage}

\clearpage

\checkoddpage
\ifoddpage (Page deliberately left blank) \clearpage \else   
\fi


\noindent
\begin{minipage}[h]{\linewidth}
	
	\begin{minipage}[h]{\linewidth}
		\vspace*{0pt}
		{\centering	
			
			\includegraphics[width=6.675cm]{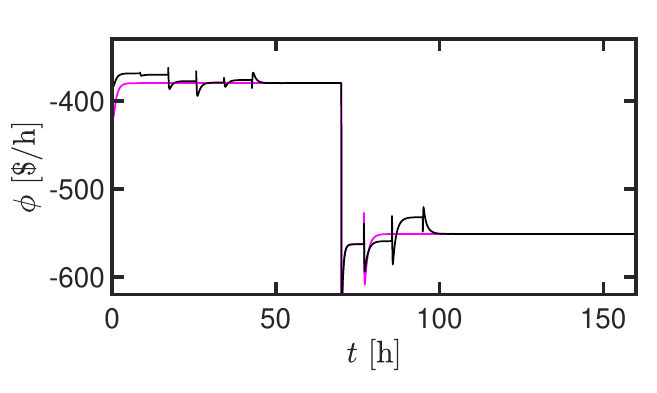}    \hskip -0ex
			\includegraphics[width=6.675cm]{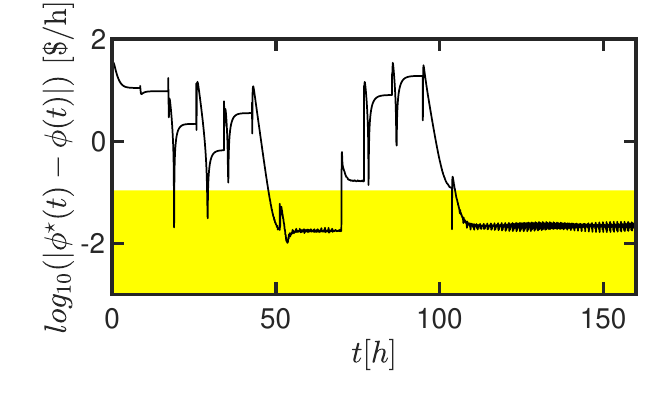} \\
			\vspace{-3mm}
			\includegraphics[width=6.675cm]{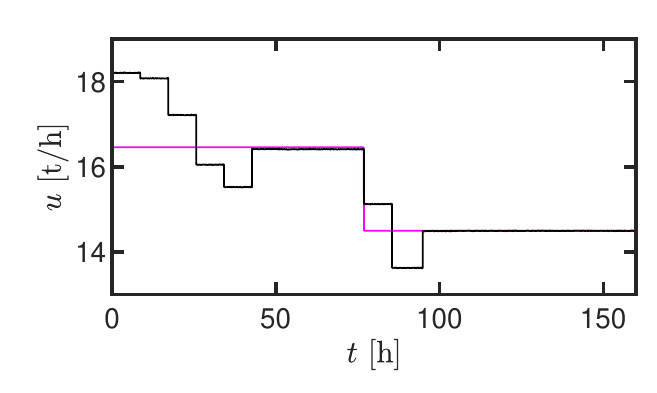}   \hskip -0ex
			\includegraphics[width=6.675cm]{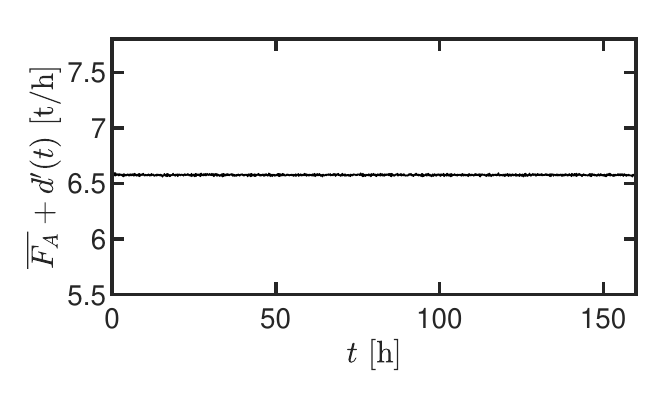} \\
			\vspace{-3mm}
			\includegraphics[width=6.675cm]{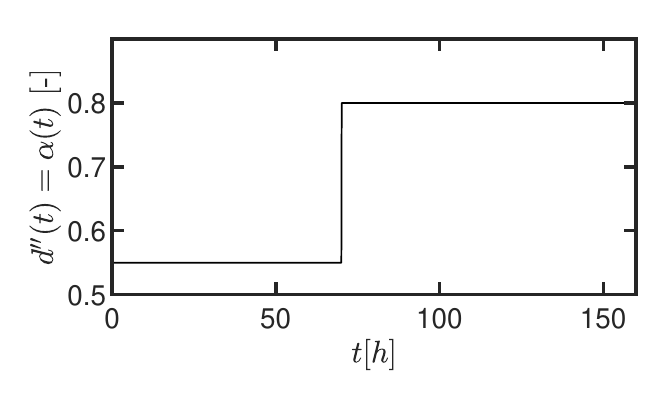}   \hskip -0ex
			\includegraphics[width=6.675cm]{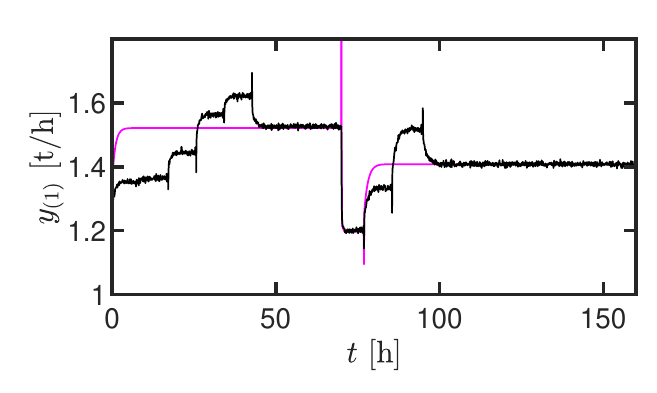}   \\
			\vspace{-3mm}
			\includegraphics[width=6.675cm]{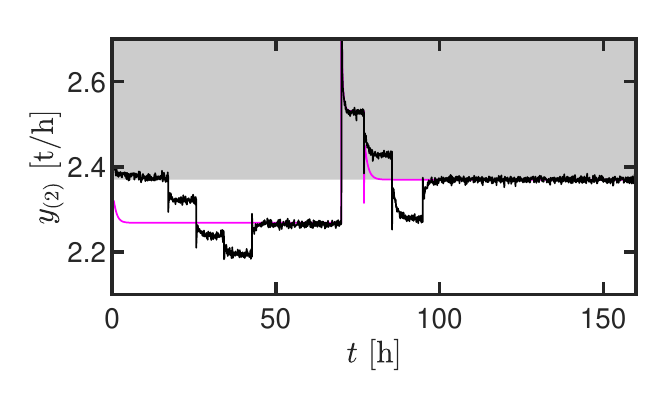} \hskip -0ex
			\includegraphics[width=6.675cm]{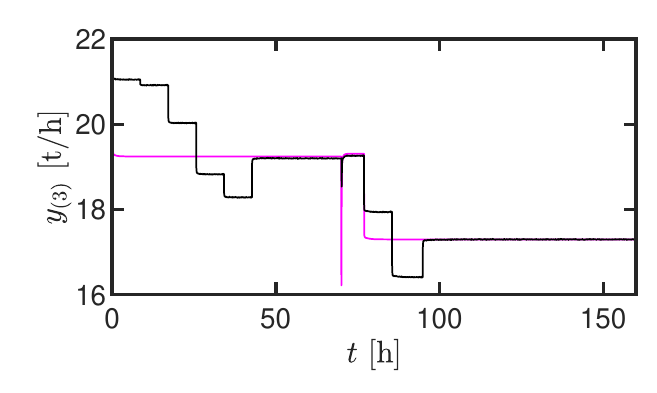} 
			
		}
	\textcolor{black}{\raisebox{1mm}{\rule{0.5cm}{0.05cm}}}: results obtained with the model and inaccurate measurements, \\
	\textcolor{magenta}{\raisebox{1mm}{\rule{0.5cm}{0.05cm}}}: results obtained with perfect model and measurements,\\
	\textcolor{gray}{\raisebox{-0.5mm}{\rule{0.5cm}{0.3cm}}}: constrained area, 
	\textcolor{yellow}{\raisebox{-0.5mm}{\rule{0.5cm}{0.3cm}}}: acceptable loss.
		\captionof{figure}{\textbf{Study 2 -- Scenario 1:} Inputs, outputs, convergence criterion, and cost w.r.t. time}
		\label{fig:6_study_2_sc1__1_inputs}
	\end{minipage} 
	\begin{minipage}[h]{\linewidth}
		\vspace*{0pt}
		{\centering	
			
			\includegraphics[trim={1.5cm 0cm 1cm  0cm},clip,width=13.35cm]{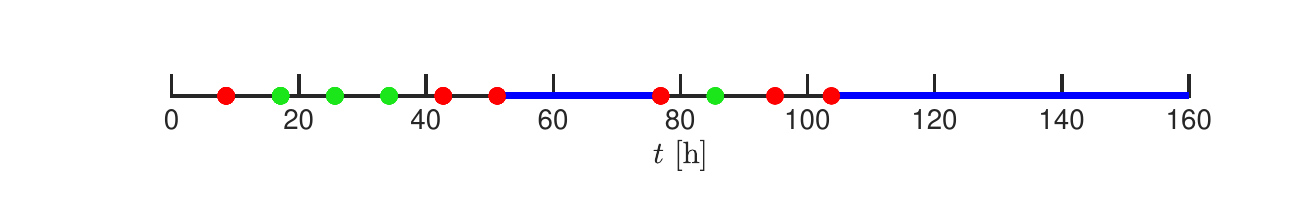}   
			
		}
	\textcolor{red}{$\bullet$}: ``normal'' decision, 
	\textcolor{green}{$\bullet$}: decision related to the validation, 
	\textcolor{blue}{\raisebox{1mm}{\rule{0.5cm}{0.05cm}}}: stand-by mode.  
		\captionof{figure}{\textbf{Study 2 -- Scenario 1:} Decision dates \& stand-by mode}
		\label{fig:6_study_2_sc1__2_decision}
	\end{minipage}
\end{minipage}

\noindent
\begin{minipage}[h]{\linewidth}
	\begingroup
	\fontsize{10pt}{12pt}\selectfont
	\vspace*{0pt}
	{\centering	
		\begin{minipage}[t]{4.3cm}
			\vspace{0pt}
			\begin{overpic}[trim={2.6cm 4.5cm 17cm  0.4cm},clip,height=4cm]{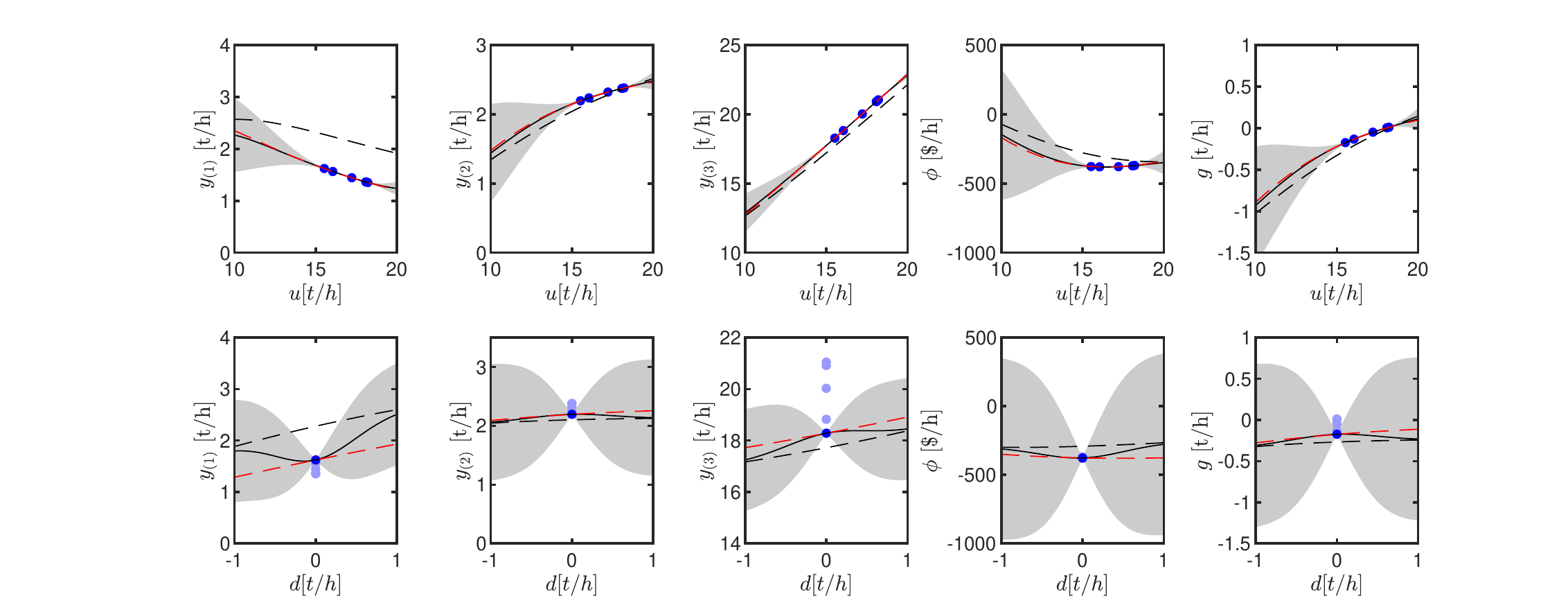}
				\put (20,22) {\textcolor{white}{\colorbox{black}{\textbf{\textbf{5}}}} }
			\end{overpic}
		\end{minipage} \hskip -3ex  
		\begin{minipage}[t]{3cm}
			\vspace{0pt}
			\begin{overpic}[trim={3.4cm 4.5cm 17cm 0.4cm},clip,height=4cm]{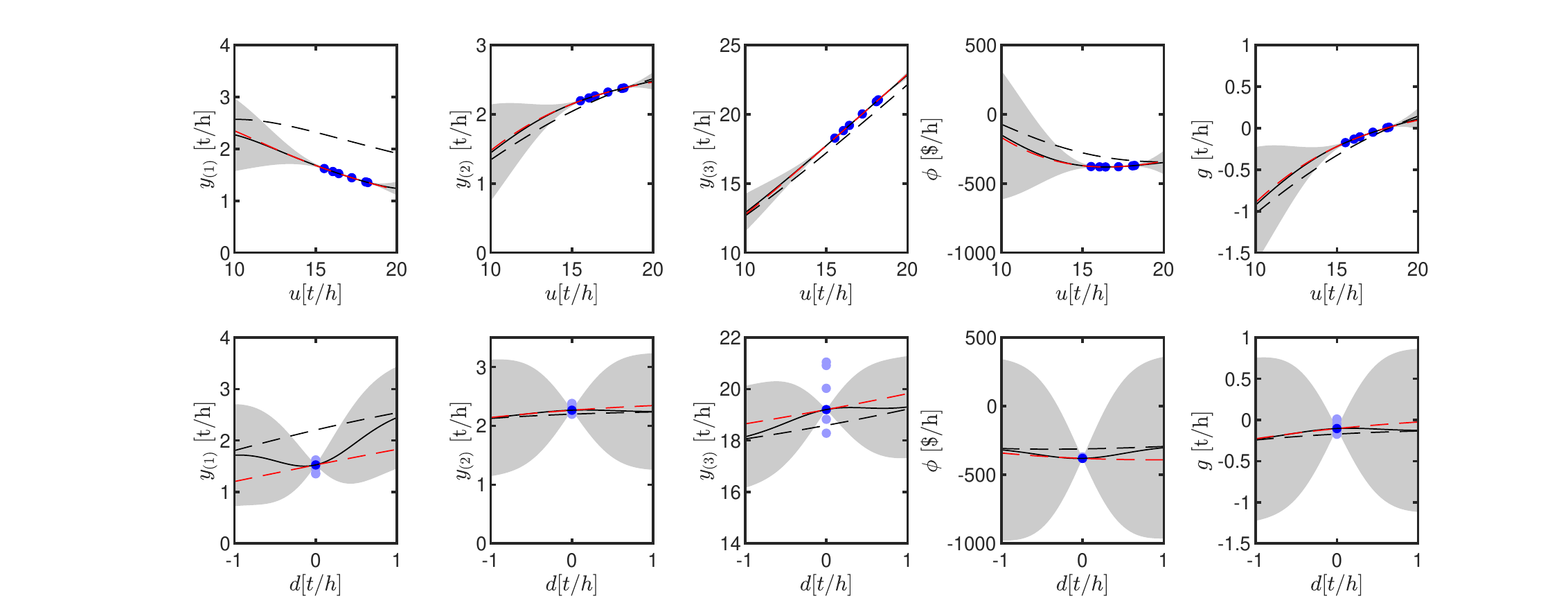}
				\put (-2,10.6) {\textcolor{black}{\colorbox{white}{10}}}
				\put (24,10.6) {\textcolor{black}{\colorbox{white}{15}}}
				\put (50,10.6) {\textcolor{black}{\colorbox{white}{20}}}
				\put (1,22) {\textcolor{white}{\colorbox{black}{\textbf{\textbf{6}}}} }
			\end{overpic}
		\end{minipage} \hskip -0ex   
		\begin{minipage}[t]{3cm}
			\vspace{0pt}
			\begin{overpic}[trim={3.4cm 4.5cm 17cm 0.4cm},clip,height=4cm]{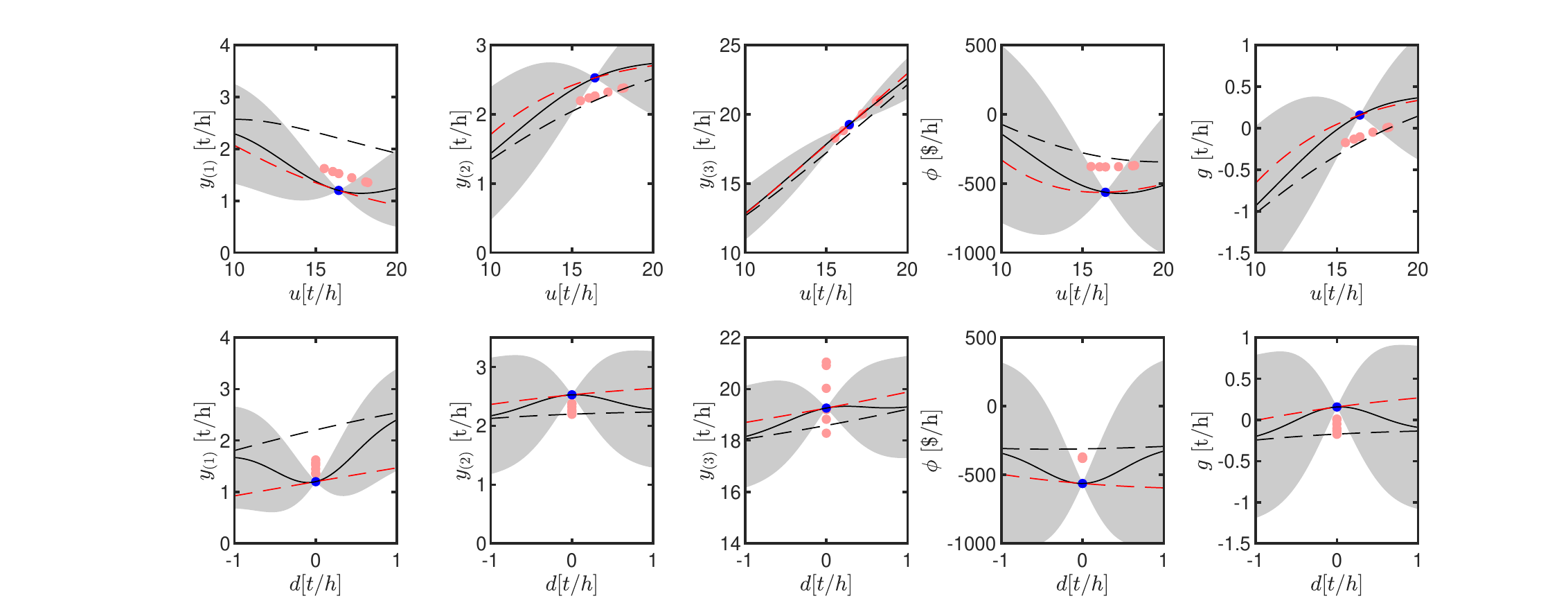}
				\put (-2,10.6) {\textcolor{black}{\colorbox{white}{10}}}
				\put (24,10.6) {\textcolor{black}{\colorbox{white}{15}}}
				\put (50,10.6) {\textcolor{black}{\colorbox{white}{20}}}
				\put (1,22) {\textcolor{white}{\colorbox{black}{\textbf{\textbf{7}}}} }
			\end{overpic}
		\end{minipage} \hskip -0ex   
		\begin{minipage}[t]{3cm}
			\vspace{0pt}
			\begin{overpic}[trim={3.4cm 4.5cm 17cm 0.4cm},clip,height=4cm]{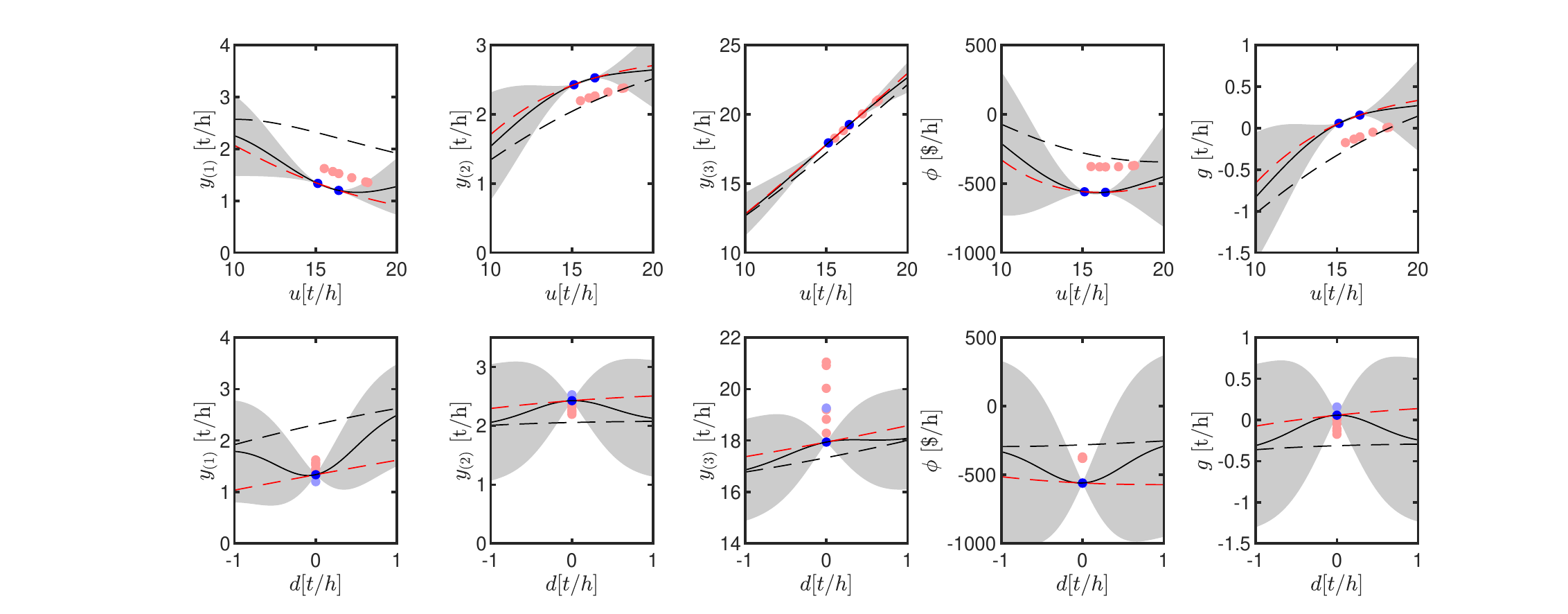}
				\put (-2,10.6) {\textcolor{black}{\colorbox{white}{10}}}
				\put (24,10.6) {\textcolor{black}{\colorbox{white}{15}}}
				\put (50,10.6) {\textcolor{black}{\colorbox{white}{20}}}
				\put (1,22) {\textcolor{white}{\colorbox{black}{\textbf{\textbf{8}}}} }
			\end{overpic}
		\end{minipage} \\ 
	\vspace{3mm}
	a) Updated versions of $f_{(1)}$\\	
	}
	\endgroup
	\begingroup
	\fontsize{10pt}{12pt}\selectfont
	\vspace*{-3pt}
	{\centering	
		\begin{minipage}[t]{4.3cm}
			\vspace{0pt}
			\begin{overpic}[trim={6.37cm 4.5cm 13.1cm  0.4cm},clip,height=4cm]{images/ch6/WO/Study_2/WO_Opt2_Sc1_GP_view_k_5.pdf}
				\put (20,22) {\textcolor{white}{\colorbox{black}{\textbf{\textbf{5}}}} }
			\end{overpic}
		\end{minipage} \hskip -3ex  
		\begin{minipage}[t]{3cm}
			\vspace{0pt}
			\begin{overpic}[trim={7.15cm 4.5cm 13.3cm  0.4cm},clip,height=4cm]{images/ch6/WO/Study_2/WO_Opt2_Sc1_GP_view_k_6.pdf}
				\put (-2,10.6) {\textcolor{black}{\colorbox{white}{10}}}
				\put (24,10.6) {\textcolor{black}{\colorbox{white}{15}}}
				\put (50,10.6) {\textcolor{black}{\colorbox{white}{20}}}
				\put (1,22) {\textcolor{white}{\colorbox{black}{\textbf{\textbf{6}}}} }
			\end{overpic}
		\end{minipage} \hskip -0ex   
		\begin{minipage}[t]{3cm}
			\vspace{0pt}
			\begin{overpic}[trim={7.15cm 4.5cm 13.3cm  0.4cm},clip,height=4cm]{images/ch6/WO/Study_2/WO_Opt2_Sc1_GP_view_k_7.pdf}
				\put (-2,10.6) {\textcolor{black}{\colorbox{white}{10}}}
				\put (24,10.6) {\textcolor{black}{\colorbox{white}{15}}}
				\put (50,10.6) {\textcolor{black}{\colorbox{white}{20}}}
				\put (1,22) {\textcolor{white}{\colorbox{black}{\textbf{\textbf{7}}}} }
			\end{overpic}
		\end{minipage} \hskip -0ex   
		\begin{minipage}[t]{3cm}
			\vspace{0pt}
			\begin{overpic}[trim={7.15cm 4.5cm 13.3cm  0.4cm},clip,height=4cm]{images/ch6/WO/Study_2/WO_Opt2_Sc1_GP_view_k_8.pdf}
				\put (-2,10.6) {\textcolor{black}{\colorbox{white}{10}}}
				\put (24,10.6) {\textcolor{black}{\colorbox{white}{15}}}
				\put (50,10.6) {\textcolor{black}{\colorbox{white}{20}}}
				\put (1,22) {\textcolor{white}{\colorbox{black}{\textbf{\textbf{8}}}} }
			\end{overpic}
		\end{minipage} \\ 
	\vspace{3mm}
	b) Updated versions of $f_{(2)}$	\\
	}
	\endgroup
	\begingroup
	\fontsize{10pt}{12pt}\selectfont
	\vspace*{-3pt}
	{\centering	
		\begin{minipage}[t]{4.3cm}
			\vspace{0pt}
			\begin{overpic}[trim={10.1cm 4.5cm 9.5cm  0.4cm},clip,height=4cm]{images/ch6/WO/Study_2/WO_Opt2_Sc1_GP_view_k_5.pdf}
				\put (20,22) {\textcolor{white}{\colorbox{black}{\textbf{\textbf{5}}}} }
			\end{overpic}
		\end{minipage} \hskip -3ex  
		\begin{minipage}[t]{3cm}
			\vspace{0pt}
			\begin{overpic}[trim={10.9cm 4.5cm 9.5cm  0.4cm},clip,height=4cm]{images/ch6/WO/Study_2/WO_Opt2_Sc1_GP_view_k_6.pdf}
				\put (-2,10.6) {\textcolor{black}{\colorbox{white}{10}}}
				\put (24,10.6) {\textcolor{black}{\colorbox{white}{15}}}
				\put (50,10.6) {\textcolor{black}{\colorbox{white}{20}}}
				\put (1,22) {\textcolor{white}{\colorbox{black}{\textbf{\textbf{6}}}} }
			\end{overpic}
		\end{minipage} \hskip -0ex   
		\begin{minipage}[t]{3cm}
			\vspace{0pt}
			\begin{overpic}[trim={10.9cm 4.5cm 9.5cm  0.4cm},clip,height=4cm]{images/ch6/WO/Study_2/WO_Opt2_Sc1_GP_view_k_7.pdf}
				\put (-2,10.6) {\textcolor{black}{\colorbox{white}{10}}}
				\put (24,10.6) {\textcolor{black}{\colorbox{white}{15}}}
				\put (50,10.6) {\textcolor{black}{\colorbox{white}{20}}}
				\put (1,22) {\textcolor{white}{\colorbox{black}{\textbf{\textbf{7}}}} }
			\end{overpic}
		\end{minipage} \hskip -0ex   
		\begin{minipage}[t]{3cm}
			\vspace{0pt}
			\begin{overpic}[trim={10.9cm 4.5cm 9.5cm  0.4cm},clip,height=4cm]{images/ch6/WO/Study_2/WO_Opt2_Sc1_GP_view_k_8.pdf}
				\put (-2,10.6) {\textcolor{black}{\colorbox{white}{10}}}
				\put (24,10.6) {\textcolor{black}{\colorbox{white}{15}}}
				\put (50,10.6) {\textcolor{black}{\colorbox{white}{20}}}
				\put (1,22) {\textcolor{white}{\colorbox{black}{\textbf{\textbf{8}}}} }
			\end{overpic}
		\end{minipage}  \\ 
	\vspace{3mm}
	c) Updated versions of $f_{(3)}$\\
	}
	\endgroup
	\textcolor{blue}{$\bullet$}: \textit{used} data, 
	\textcolor{red!50!white}{$\bullet$}: \textit{rejected} data, 
	\textcolor{black}{\raisebox{1mm}{\rule{0.2cm}{0.05cm}\hspace{0.1cm}\rule{0.2cm}{0.05cm}}}: model,
	\textcolor{red}{\raisebox{1mm}{\rule{0.2cm}{0.05cm}\hspace{0.1cm}\rule{0.2cm}{0.05cm}}}: plant,
	\textcolor{black}{\raisebox{1mm}{\rule{0.5cm}{0.05cm}}}: updated model,\\
	\textcolor{gray}{\raisebox{-0.5mm}{\rule{0.5cm}{0.3cm}}}: confidence domain. 
	
	\captionof{figure}{\textbf{Study 2 -- Scenario 1:} The updated model of $f_{p(1)}$, $f_{p(2)}$, and $f_{p(3)}$,  at the \textcolor{white}{\colorbox{black}{\textbf{\textbf{k-th}}}} iteration of ASP. The  datum obtained at the 7th iteration is far from the confidence area (the gray domain). 
	}  
	\label{fig:6_study_2_sc1__3_details}
{\centering	
	\includegraphics[width=6.675cm]{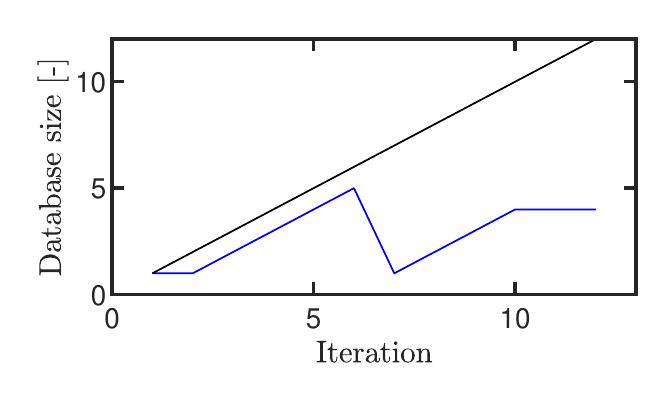} \\ 
}

\textcolor{black}{\raisebox{1mm}{\rule{0.5cm}{0.05cm}}}: amount of collected data, 
\textcolor{blue}{\raisebox{1mm}{\rule{0.5cm}{0.05cm}}}: amount of data in $\mathcal{D}_{III}$. 
\captionof{figure}{\textbf{Study 2 -- Scenario 1:} Database size w.r.t. time}
\label{fig:6_study_2_sc1__4_database}
\end{minipage}


\noindent
\begin{minipage}[h]{\linewidth}
	
	\begin{minipage}[h]{\linewidth}
		\vspace*{0pt}
		{\centering	
			
			\includegraphics[width=6.675cm]{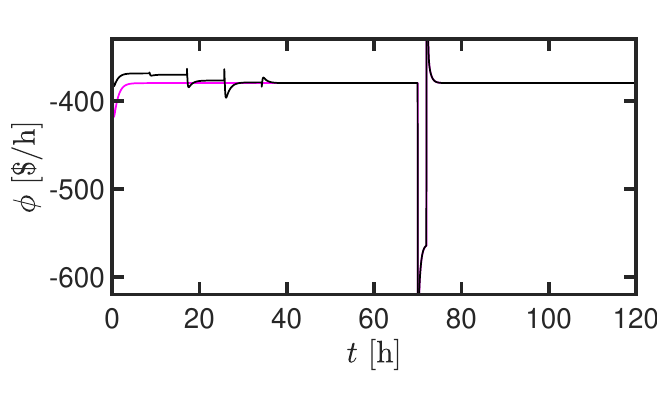}    \hskip -0ex
			\includegraphics[width=6.675cm]{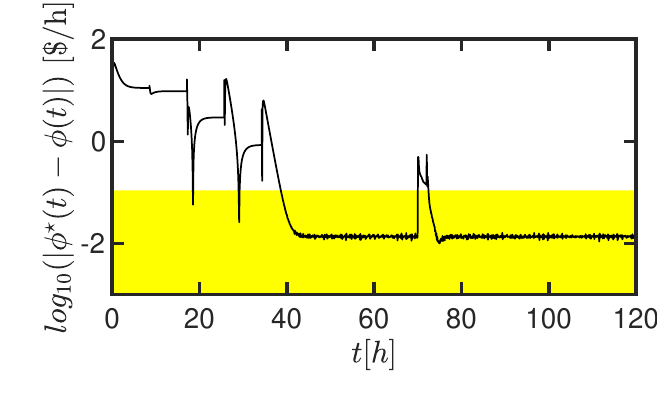} \hskip -0ex
			\includegraphics[width=6.675cm]{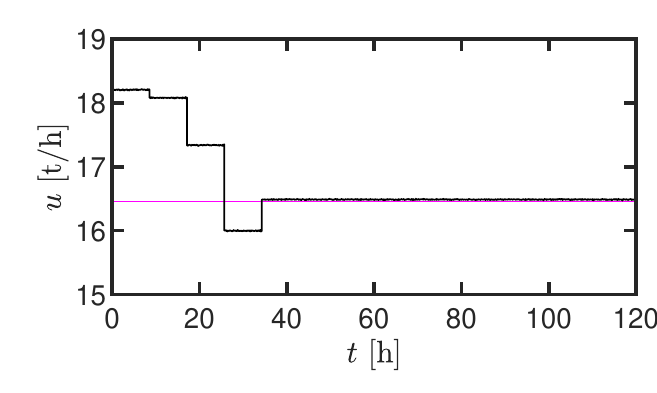}   \hskip -0ex
			\includegraphics[width=6.675cm]{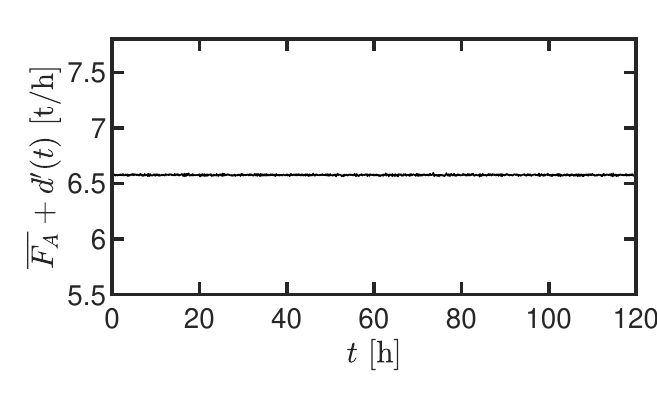} \hskip -0ex
			\includegraphics[width=6.675cm]{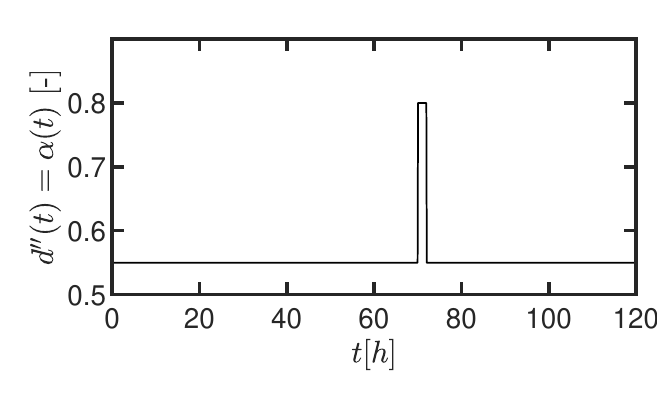}   \hskip -0ex
			\includegraphics[width=6.675cm]{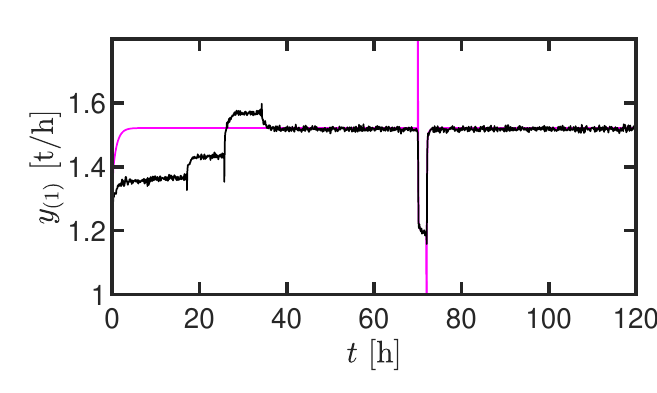}   \hskip -0ex
			\includegraphics[width=6.675cm]{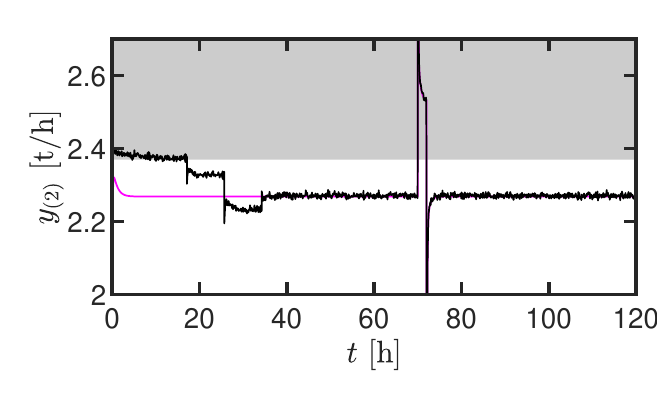} \hskip -0ex
			\includegraphics[width=6.675cm]{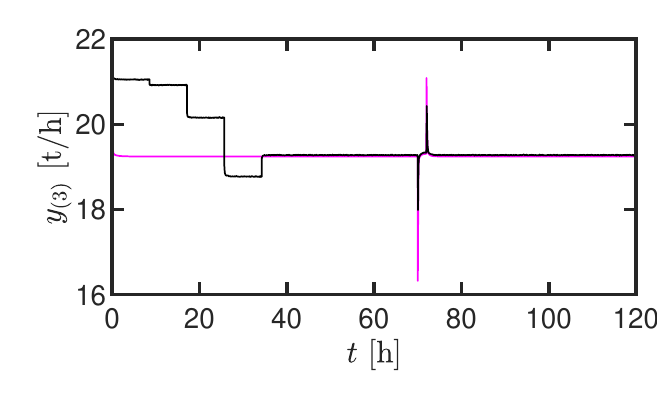}   \hskip -0ex
			
		}
		\textcolor{black}{\raisebox{1mm}{\rule{0.5cm}{0.05cm}}}: results obtained with the model and inaccurate measurements, \\
		\textcolor{magenta}{\raisebox{1mm}{\rule{0.5cm}{0.05cm}}}: results obtained with perfect model and measurements,\\
		\textcolor{gray}{\raisebox{-0.5mm}{\rule{0.5cm}{0.3cm}}}: constrained area, 
		\textcolor{yellow}{\raisebox{-0.5mm}{\rule{0.5cm}{0.3cm}}}: acceptable loss.
		\captionof{figure}{\textbf{Study 2 -- Scenario 2:} Inputs, outputs, convergence criterion, and cost w.r.t. time}
		\label{fig:6_study_2_sc2__1_inputs}
	\end{minipage} 
	\begin{minipage}[h]{\linewidth}
		\vspace*{0pt}
		{\centering	
			
			\includegraphics[trim={1.5cm 0cm 1cm  0cm},clip,width=13.35cm]{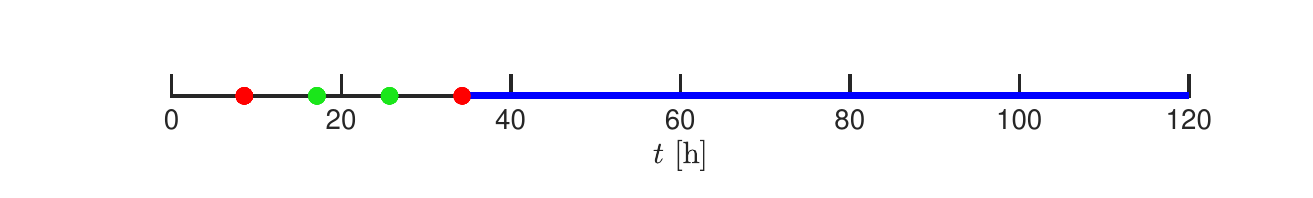}   
			
		}
		\textcolor{red}{$\bullet$}: ``normal'' decision, 
		\textcolor{green}{$\bullet$}: decision related to the validation, 
		\textcolor{blue}{\raisebox{1mm}{\rule{0.5cm}{0.05cm}}}: stand-by mode.  
		\captionof{figure}{\textbf{Study 2 -- Scenario 2:} Decision dates \& stand-by mode}
		\label{fig:6_study_2_sc2__2_decision}
	\end{minipage}
\end{minipage}

\noindent
\begin{minipage}[h]{\linewidth}
	{\centering	
		\includegraphics[width=6.675cm]{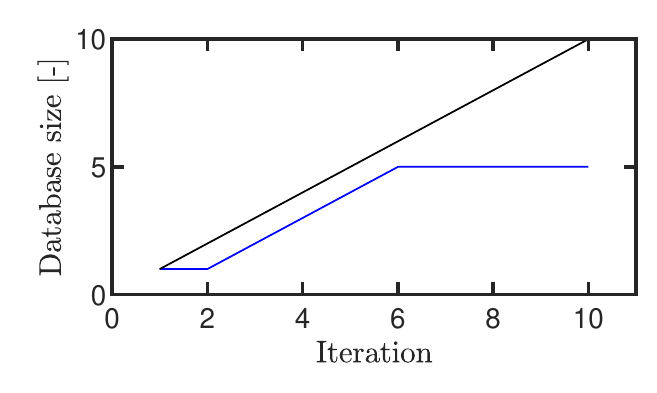} \\ 
	}
	\textcolor{black}{\raisebox{1mm}{\rule{0.5cm}{0.05cm}}}: amount of collected data, 
	\textcolor{blue}{\raisebox{1mm}{\rule{0.5cm}{0.05cm}}}: amount of data in $\mathcal{D}_{III}$. 
	\captionof{figure}{\textbf{Study 2 -- Scenario 2:} Database size w.r.t. time}
	\label{fig:6_study_2_sc2__3_database}
\end{minipage}


\noindent
\begin{minipage}[h]{\linewidth}
	
	\begin{minipage}[h]{\linewidth}
		\vspace*{0pt}
		{\centering	
			
			\includegraphics[width=6.675cm]{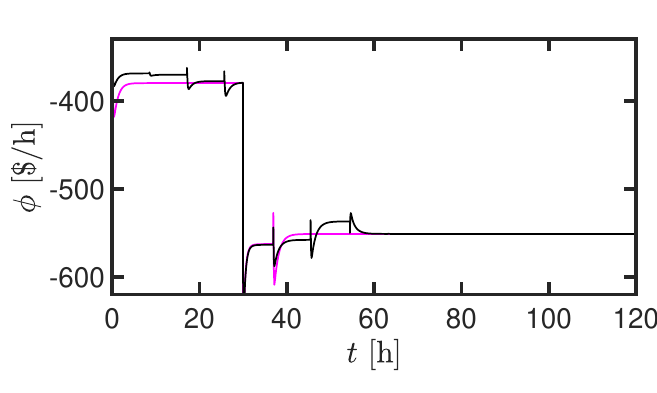}    \hskip -0ex
			\includegraphics[width=6.675cm]{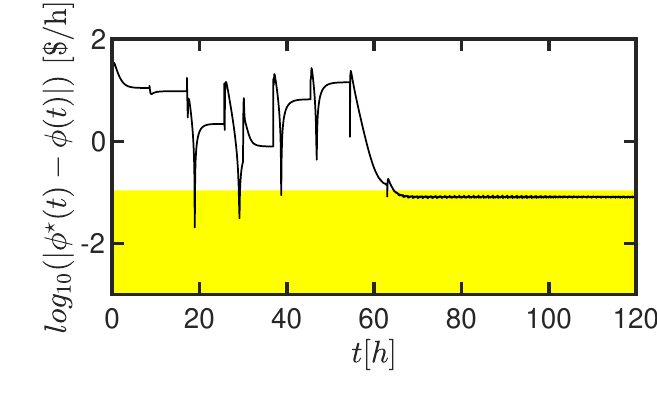} \hskip -0ex
			\includegraphics[width=6.675cm]{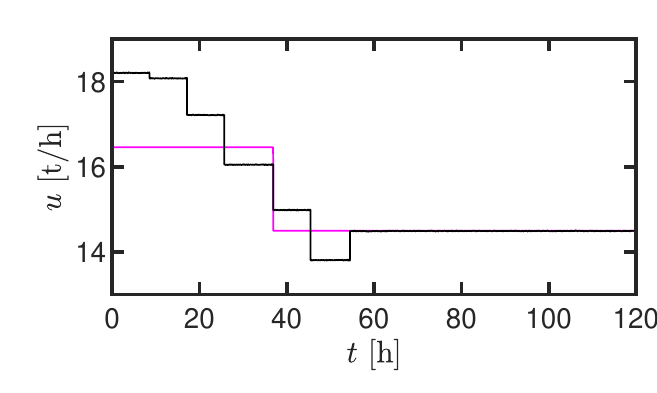}   \hskip -0ex
			\includegraphics[width=6.675cm]{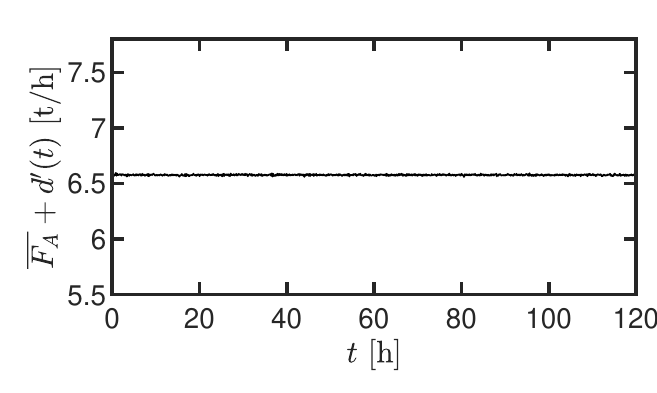} \hskip -0ex
			\includegraphics[width=6.675cm]{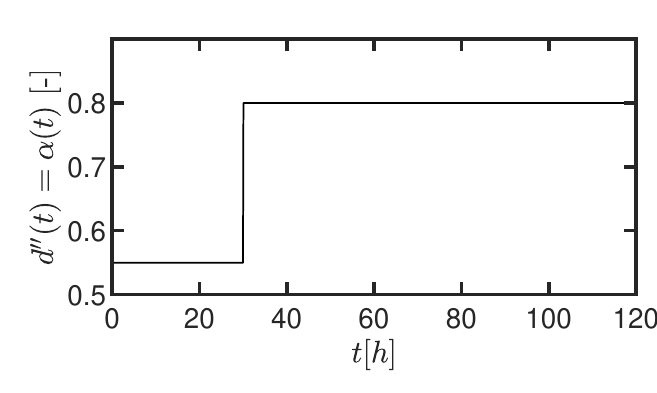}   \hskip -0ex
			\includegraphics[width=6.675cm]{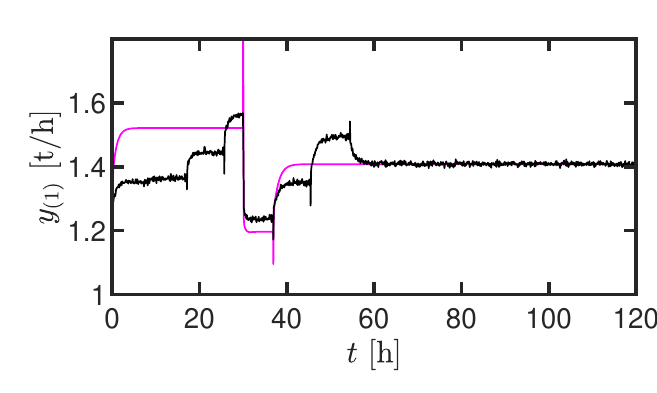}   \hskip -0ex
			\includegraphics[width=6.675cm]{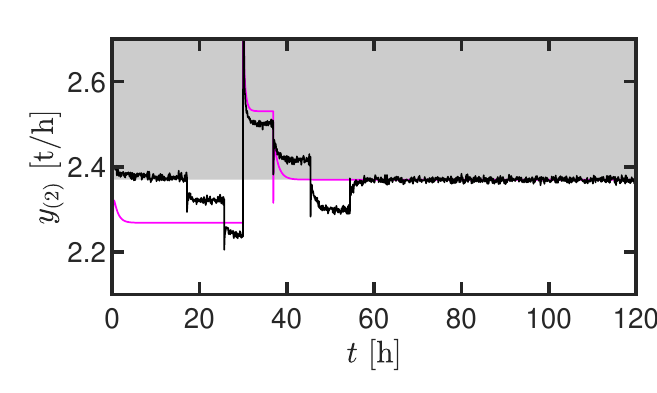} \hskip -0ex
			\includegraphics[width=6.675cm]{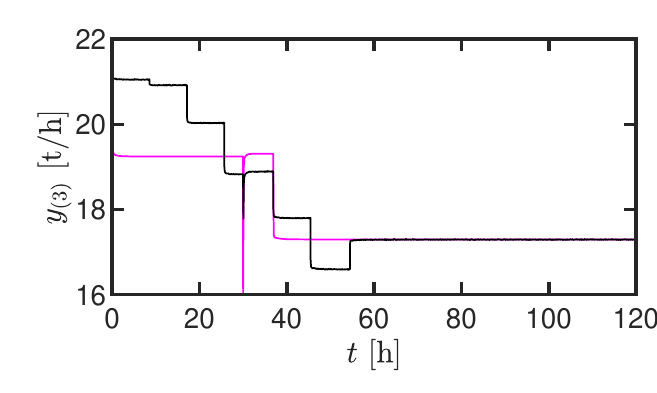}   \hskip -0ex
			
		}
		\textcolor{black}{\raisebox{1mm}{\rule{0.5cm}{0.05cm}}}: results obtained with the model and inaccurate measurements, \\
		\textcolor{magenta}{\raisebox{1mm}{\rule{0.5cm}{0.05cm}}}: results obtained with perfect model and measurements,\\
		\textcolor{gray}{\raisebox{-0.5mm}{\rule{0.5cm}{0.3cm}}}: constrained area, 
		\textcolor{yellow}{\raisebox{-0.5mm}{\rule{0.5cm}{0.3cm}}}: acceptable loss.
		\captionof{figure}{\textbf{Study 2 -- Scenario 3:} Inputs, outputs, convergence criterion, and cost w.r.t. time}
		\label{fig:6_study_2_sc3__1_inputs}
	\end{minipage} 
	\begin{minipage}[h]{\linewidth}
		\vspace*{0pt}
		{\centering	
			
			\includegraphics[trim={1.5cm 0cm 1cm  0cm},clip,width=13.35cm]{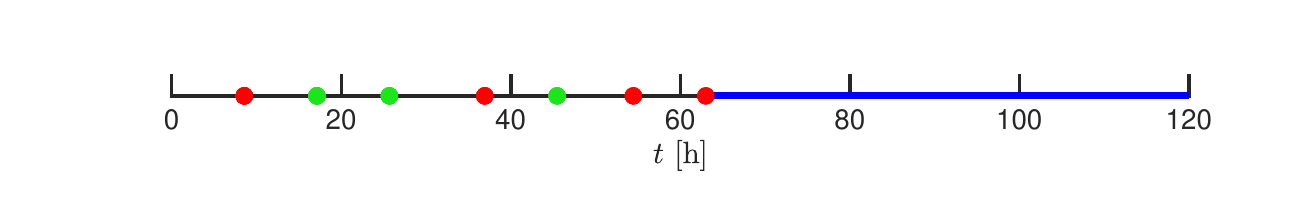}   
			
		}
		\textcolor{red}{$\bullet$}: ``normal'' decision, 
		\textcolor{green}{$\bullet$}: decision related to the validation, 
		\textcolor{blue}{\raisebox{1mm}{\rule{0.5cm}{0.05cm}}}: stand-by mode.  
		\captionof{figure}{\textbf{Study 2 -- Scenario 3:} Decision dates \& stand-by mode}
		\label{fig:6_study_2_sc3__2_decision}
	\end{minipage}
\end{minipage}

\noindent
\begin{minipage}[h]{\linewidth}
	\begingroup
	\fontsize{10pt}{12pt}\selectfont
	\vspace*{0pt}
	{\centering	
		\begin{minipage}[t]{4.3cm}
			\vspace{0pt}
			\begin{overpic}[trim={2.6cm 4.5cm 17cm  0.4cm},clip,height=4cm]{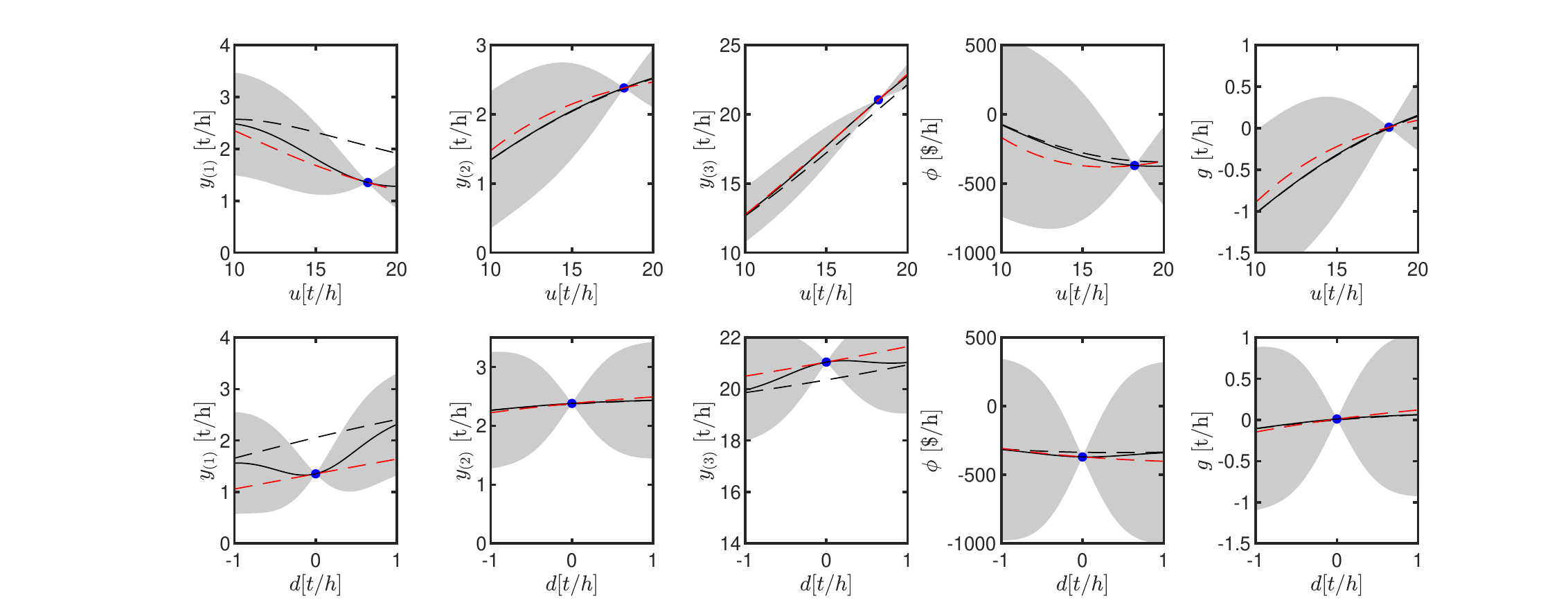}
				\put (20,22) {\textcolor{white}{\colorbox{black}{\textbf{\textbf{1}}}} }
			\end{overpic}
		\end{minipage} \hskip -3ex  
		\begin{minipage}[t]{3cm}
			\vspace{0pt}
			\begin{overpic}[trim={3.4cm 4.5cm 17cm 0.4cm},clip,height=4cm]{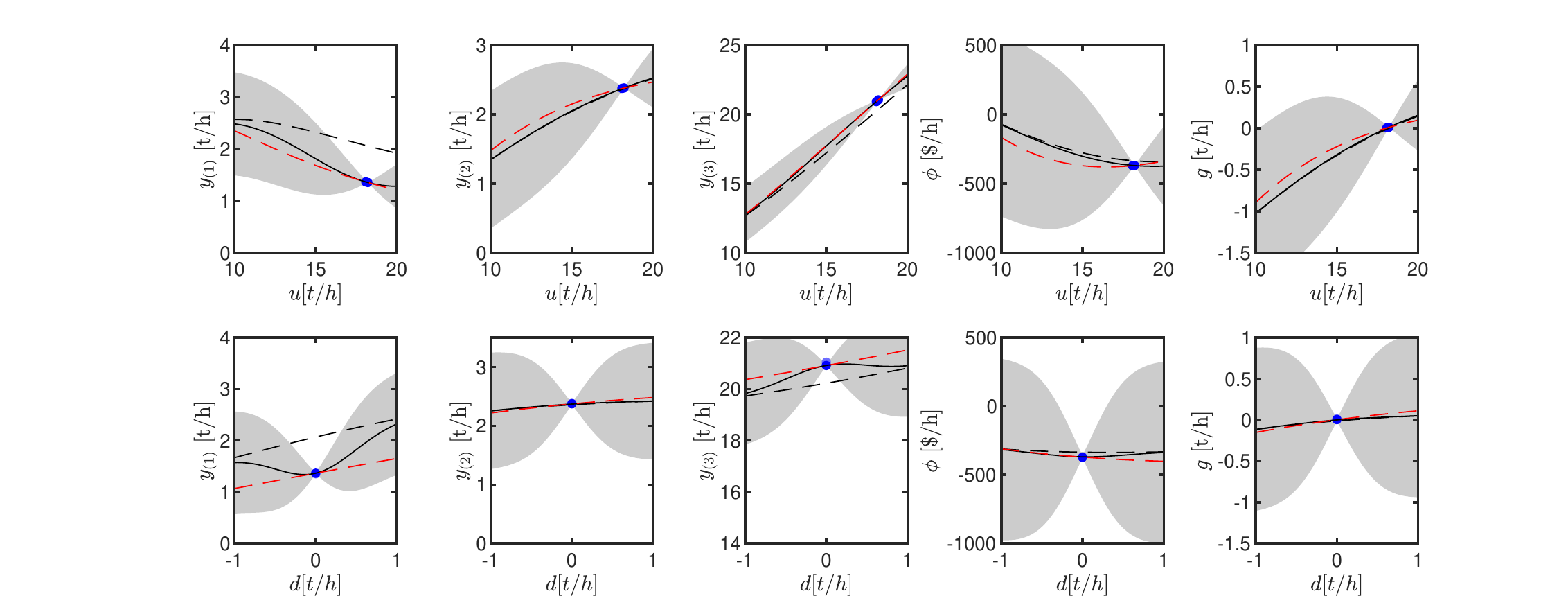}
				\put (-2,10.6) {\textcolor{black}{\colorbox{white}{10}}}
				\put (24,10.6) {\textcolor{black}{\colorbox{white}{15}}}
				\put (50,10.6) {\textcolor{black}{\colorbox{white}{20}}}
				\put (1,22) {\textcolor{white}{\colorbox{black}{\textbf{\textbf{2}}}} }
			\end{overpic}
		\end{minipage} \hskip -0ex   
		\begin{minipage}[t]{3cm}
			\vspace{0pt}
			\begin{overpic}[trim={3.4cm 4.5cm 17cm 0.4cm},clip,height=4cm]{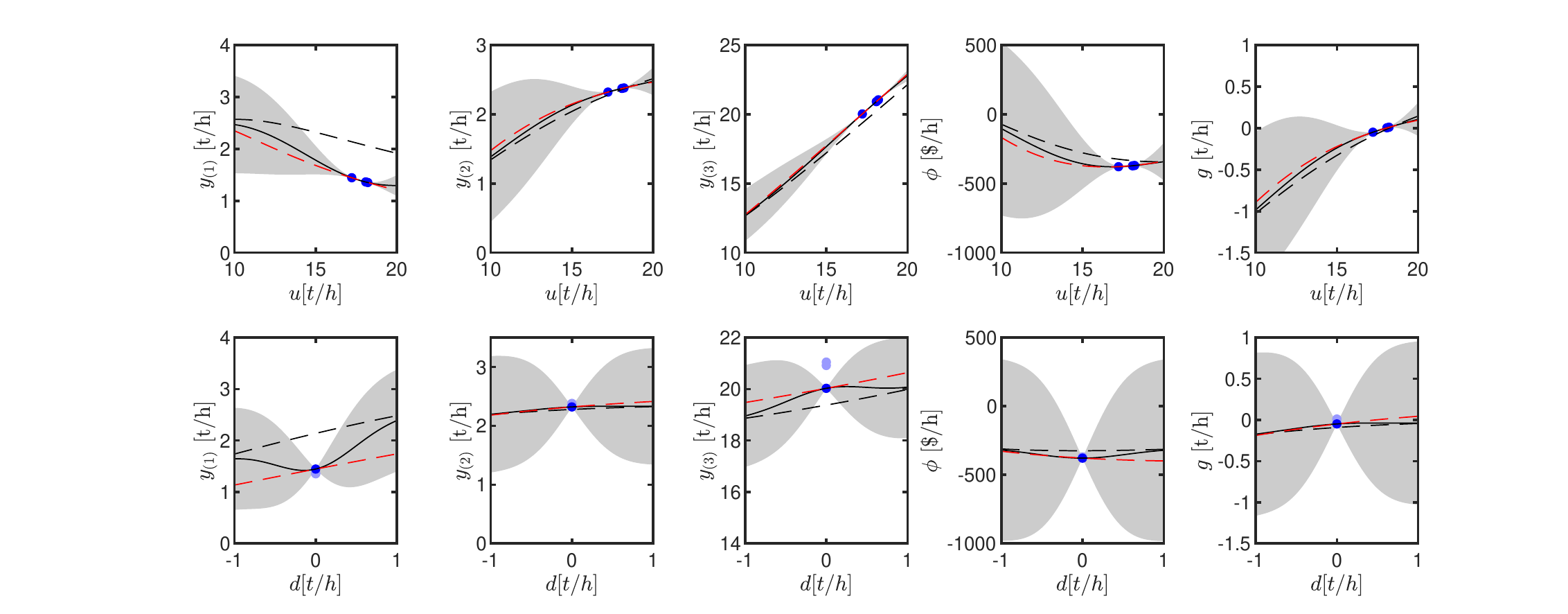}
				\put (-2,10.6) {\textcolor{black}{\colorbox{white}{10}}}
				\put (24,10.6) {\textcolor{black}{\colorbox{white}{15}}}
				\put (50,10.6) {\textcolor{black}{\colorbox{white}{20}}}
				\put (1,22) {\textcolor{white}{\colorbox{black}{\textbf{\textbf{3}}}} }
			\end{overpic}
		\end{minipage} \hskip -0ex   
		\begin{minipage}[t]{3cm}
			\vspace{0pt}
			\begin{overpic}[trim={3.4cm 4.5cm 17cm 0.4cm},clip,height=4cm]{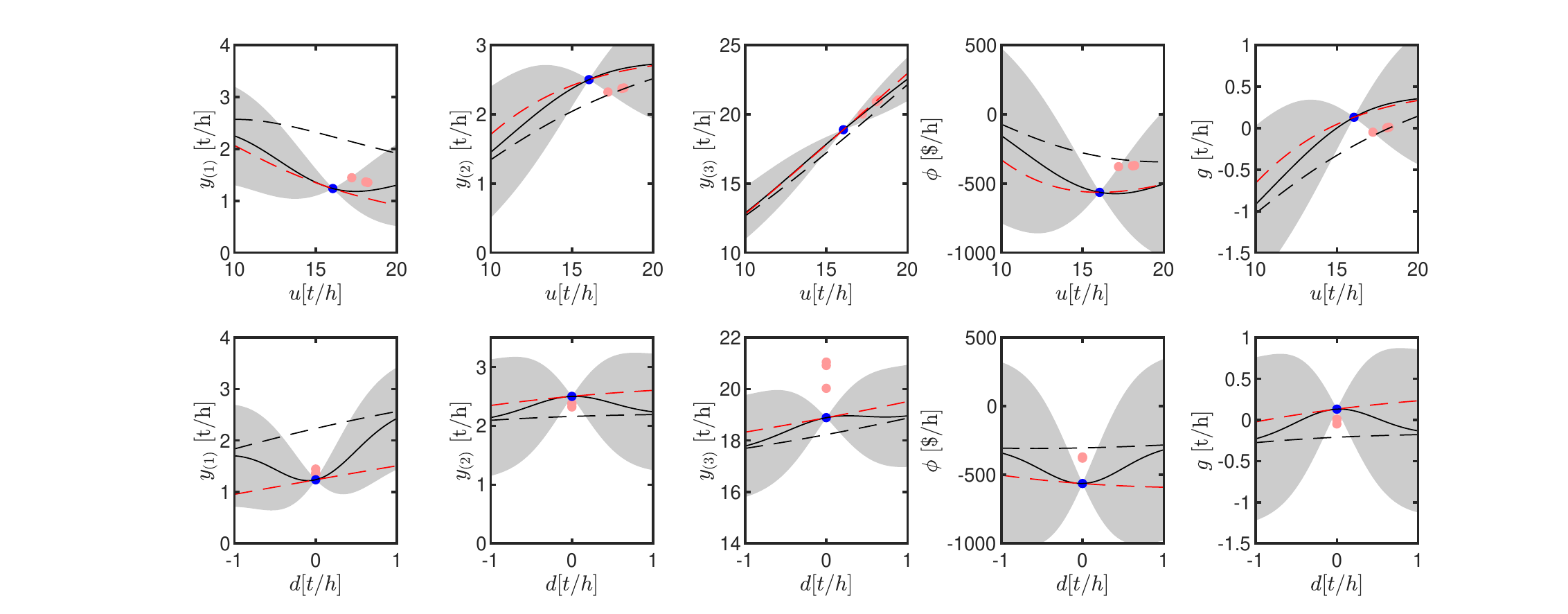}
				\put (-2,10.6) {\textcolor{black}{\colorbox{white}{10}}}
				\put (24,10.6) {\textcolor{black}{\colorbox{white}{15}}}
				\put (50,10.6) {\textcolor{black}{\colorbox{white}{20}}}
				\put (1,22) {\textcolor{white}{\colorbox{black}{\textbf{\textbf{4}}}} }
			\end{overpic}
		\end{minipage} \\ 
		\vspace{3mm}
		a) Updated versions of $f_{(1)}$\\	
	}
	\endgroup
	\begingroup
	\fontsize{10pt}{12pt}\selectfont
	\vspace*{-3pt}
	{\centering	
		\begin{minipage}[t]{4.3cm}
			\vspace{0pt}
			\begin{overpic}[trim={6.37cm 4.5cm 13.1cm  0.4cm},clip,height=4cm]{images/ch6/WO/Study_2/WO_Opt2_Sc3_GP_view_k_1.pdf}
				\put (20,22) {\textcolor{white}{\colorbox{black}{\textbf{\textbf{1}}}} }
			\end{overpic}
		\end{minipage} \hskip -3ex  
		\begin{minipage}[t]{3cm}
			\vspace{0pt}
			\begin{overpic}[trim={7.15cm 4.5cm 13.3cm  0.4cm},clip,height=4cm]{images/ch6/WO/Study_2/WO_Opt2_Sc3_GP_view_k_2.pdf}
				\put (-2,10.6) {\textcolor{black}{\colorbox{white}{10}}}
				\put (24,10.6) {\textcolor{black}{\colorbox{white}{15}}}
				\put (50,10.6) {\textcolor{black}{\colorbox{white}{20}}}
				\put (1,22) {\textcolor{white}{\colorbox{black}{\textbf{\textbf{2}}}} }
			\end{overpic}
		\end{minipage} \hskip -0ex   
		\begin{minipage}[t]{3cm}
			\vspace{0pt}
			\begin{overpic}[trim={7.15cm 4.5cm 13.3cm  0.4cm},clip,height=4cm]{images/ch6/WO/Study_2/WO_Opt2_Sc3_GP_view_k_3.pdf}
				\put (-2,10.6) {\textcolor{black}{\colorbox{white}{10}}}
				\put (24,10.6) {\textcolor{black}{\colorbox{white}{15}}}
				\put (50,10.6) {\textcolor{black}{\colorbox{white}{20}}}
				\put (1,22) {\textcolor{white}{\colorbox{black}{\textbf{\textbf{3}}}} }
			\end{overpic}
		\end{minipage} \hskip -0ex   
		\begin{minipage}[t]{3cm}
			\vspace{0pt}
			\begin{overpic}[trim={7.15cm 4.5cm 13.3cm  0.4cm},clip,height=4cm]{images/ch6/WO/Study_2/WO_Opt2_Sc3_GP_view_k_4.pdf}
				\put (-2,10.6) {\textcolor{black}{\colorbox{white}{10}}}
				\put (24,10.6) {\textcolor{black}{\colorbox{white}{15}}}
				\put (50,10.6) {\textcolor{black}{\colorbox{white}{20}}}
				\put (1,22) {\textcolor{white}{\colorbox{black}{\textbf{\textbf{4}}}} }
			\end{overpic}
		\end{minipage} \\ 
		\vspace{3mm}
		b) Updated versions of $f_{(2)}$	\\
	}
	\endgroup
	\begingroup
	\fontsize{10pt}{12pt}\selectfont
	\vspace*{-3pt}
	{\centering	
		\begin{minipage}[t]{4.3cm}
			\vspace{0pt}
			\begin{overpic}[trim={10.1cm 4.5cm 9.5cm  0.4cm},clip,height=4cm]{images/ch6/WO/Study_2/WO_Opt2_Sc3_GP_view_k_1.pdf}
				\put (20,22) {\textcolor{white}{\colorbox{black}{\textbf{\textbf{1}}}} }
			\end{overpic}
		\end{minipage} \hskip -3ex  
		\begin{minipage}[t]{3cm}
			\vspace{0pt}
			\begin{overpic}[trim={10.9cm 4.5cm 9.5cm  0.4cm},clip,height=4cm]{images/ch6/WO/Study_2/WO_Opt2_Sc3_GP_view_k_2.pdf}
				\put (-2,10.6) {\textcolor{black}{\colorbox{white}{10}}}
				\put (24,10.6) {\textcolor{black}{\colorbox{white}{15}}}
				\put (50,10.6) {\textcolor{black}{\colorbox{white}{20}}}
				\put (1,22) {\textcolor{white}{\colorbox{black}{\textbf{\textbf{2}}}} }
			\end{overpic}
		\end{minipage} \hskip -0ex   
		\begin{minipage}[t]{3cm}
			\vspace{0pt}
			\begin{overpic}[trim={10.9cm 4.5cm 9.5cm  0.4cm},clip,height=4cm]{images/ch6/WO/Study_2/WO_Opt2_Sc3_GP_view_k_3.pdf}
				\put (-2,10.6) {\textcolor{black}{\colorbox{white}{10}}}
				\put (24,10.6) {\textcolor{black}{\colorbox{white}{15}}}
				\put (50,10.6) {\textcolor{black}{\colorbox{white}{20}}}
				\put (1,22) {\textcolor{white}{\colorbox{black}{\textbf{\textbf{3}}}} }
			\end{overpic}
		\end{minipage} \hskip -0ex   
		\begin{minipage}[t]{3cm}
			\vspace{0pt}
			\begin{overpic}[trim={10.9cm 4.5cm 9.5cm  0.4cm},clip,height=4cm]{images/ch6/WO/Study_2/WO_Opt2_Sc3_GP_view_k_4.pdf}
				\put (-2,10.6) {\textcolor{black}{\colorbox{white}{10}}}
				\put (24,10.6) {\textcolor{black}{\colorbox{white}{15}}}
				\put (50,10.6) {\textcolor{black}{\colorbox{white}{20}}}
				\put (1,22) {\textcolor{white}{\colorbox{black}{\textbf{\textbf{4}}}} }
			\end{overpic}
		\end{minipage}  \\ 
		\vspace{3mm}
		c) Updated versions of $f_{(3)}$\\
	}
	\endgroup
	\textcolor{blue}{$\bullet$}: \textit{used} data, 
	\textcolor{red!50!white}{$\bullet$}: \textit{rejected} data, 
	\textcolor{black}{\raisebox{1mm}{\rule{0.2cm}{0.05cm}\hspace{0.1cm}\rule{0.2cm}{0.05cm}}}: model,
	\textcolor{red}{\raisebox{1mm}{\rule{0.2cm}{0.05cm}\hspace{0.1cm}\rule{0.2cm}{0.05cm}}}: plant,
	\textcolor{black}{\raisebox{1mm}{\rule{0.5cm}{0.05cm}}}: updated model,
	\textcolor{gray}{\raisebox{-0.5mm}{\rule{0.5cm}{0.3cm}}}: confidence domain. 
	
	\captionof{figure}{\textbf{Study 2 -- Scenario 3:} The updated model of $f_{p(1)}$, $f_{p(2)}$, and $f_{p(3)}$,  at the \textcolor{white}{\colorbox{black}{\textbf{\textbf{k-th}}}} iteration of ASP. The  datum obtained at the 4th iteration is far from the confidence area (the gray domain). 
	}  
	\label{fig:6_study_2_sc3__3_details}
	{\centering	
		\includegraphics[width=6.675cm]{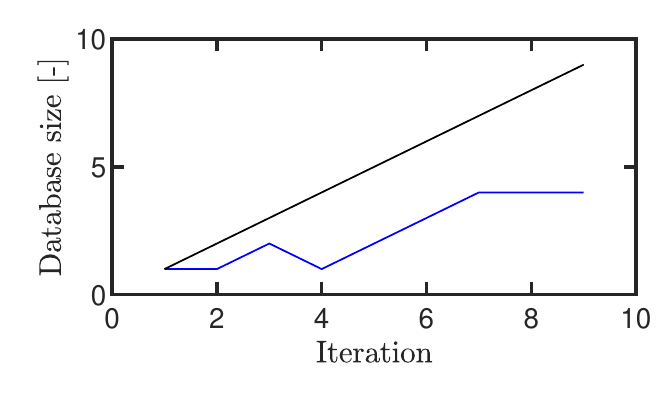} \\ 
	}
	
	\textcolor{black}{\raisebox{1mm}{\rule{0.5cm}{0.05cm}}}: amount of collected data, 
	\textcolor{blue}{\raisebox{1mm}{\rule{0.5cm}{0.05cm}}}: amount of data in $\mathcal{D}_{III}$. 
	\captionof{figure}{\textbf{Study 2 -- Scenario 2:} Database size w.r.t. time}
	\label{fig:6_study_2_sc3__4_database}
\end{minipage}

\noindent
\begin{minipage}[h]{\linewidth}
	
	\begin{minipage}[h]{\linewidth}
		\vspace*{0pt}
		{\centering	
			
			\includegraphics[width=6.675cm]{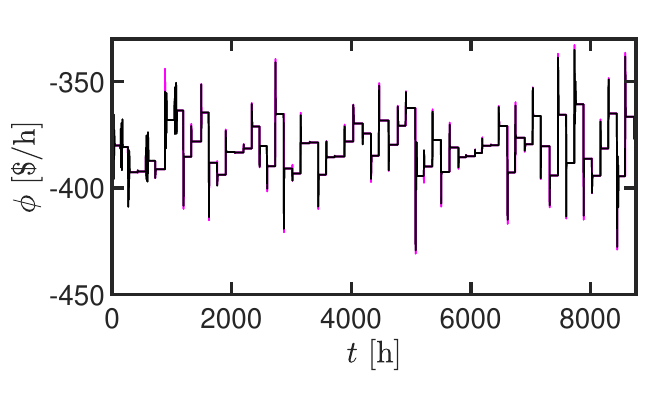}    \hskip -0ex
			\includegraphics[width=6.675cm]{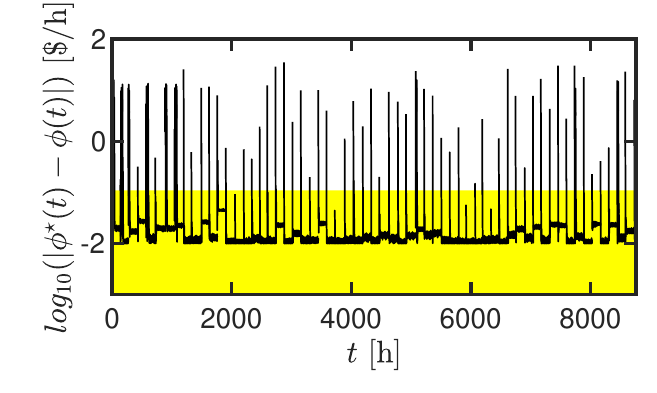} \hskip -0ex
			\includegraphics[width=6.675cm]{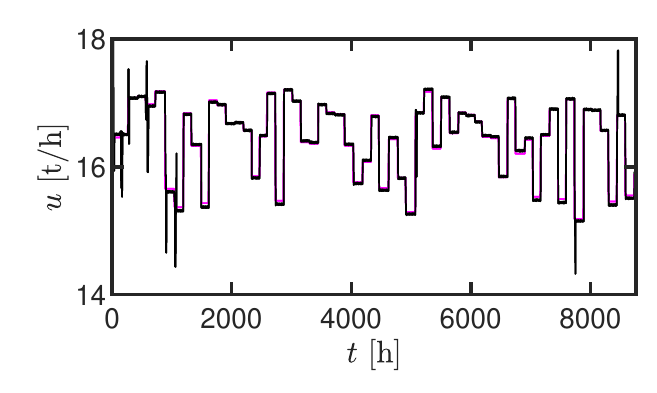}   \hskip -0ex
			\includegraphics[width=6.675cm]{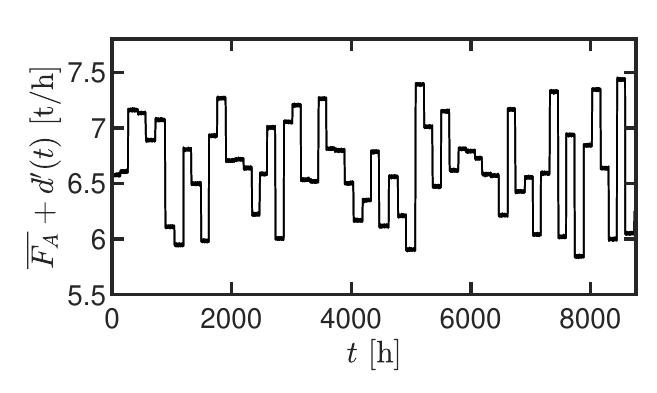} \hskip -0ex
			\includegraphics[width=6.675cm]{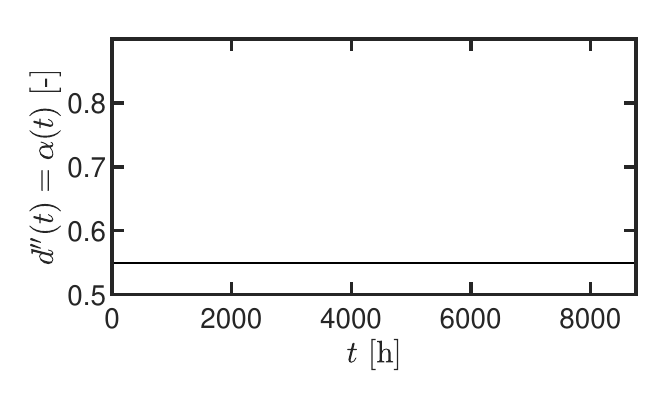}   \hskip -0ex
			\includegraphics[width=6.675cm]{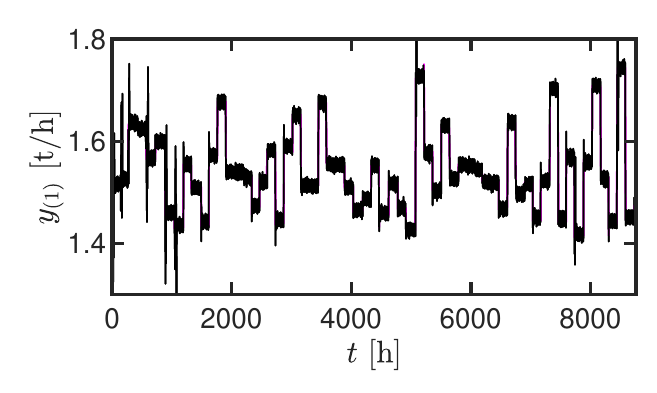}   \hskip -0ex
			\includegraphics[width=6.675cm]{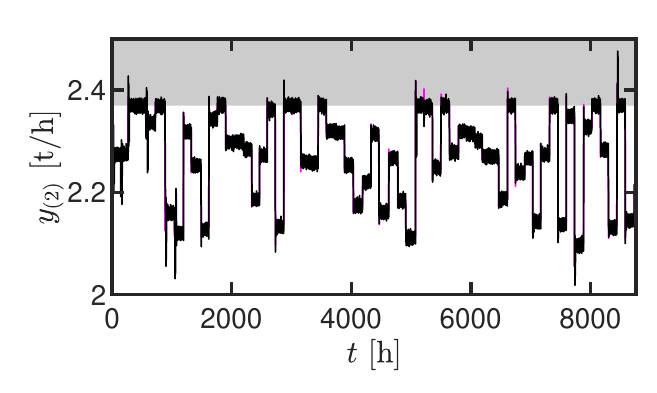} \hskip -0ex
			\includegraphics[width=6.675cm]{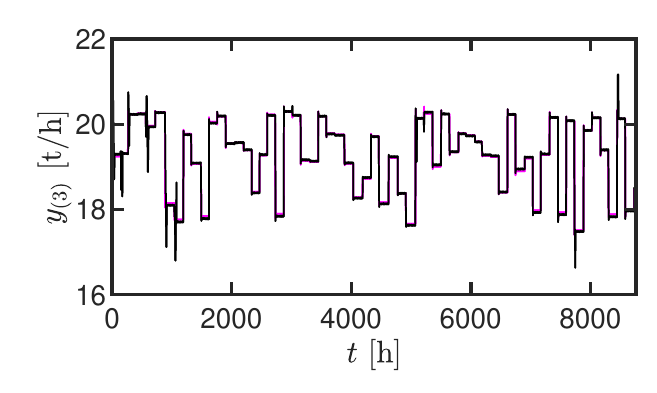}   \hskip -0ex
			
		}
		\textcolor{black}{\raisebox{1mm}{\rule{0.5cm}{0.05cm}}}: results obtained with the model and inaccurate measurements, \\
		\textcolor{magenta}{\raisebox{1mm}{\rule{0.5cm}{0.05cm}}}: results obtained with perfect model and measurements,\\
		\textcolor{gray}{\raisebox{-0.5mm}{\rule{0.5cm}{0.3cm}}}: constrained area, 
		\textcolor{yellow}{\raisebox{-0.5mm}{\rule{0.5cm}{0.3cm}}}: acceptable loss.
		\captionof{figure}{\textbf{Study 3:} Inputs, outputs, convergence criterion, and cost w.r.t. time}
		\label{fig:6_study_3__1_inputs}
	\end{minipage} 
	\begin{minipage}[h]{\linewidth}
		\vspace*{0pt}
		{\centering	
			
			\includegraphics[trim={1.5cm 0cm 1cm  0cm},clip,width=13.35cm]{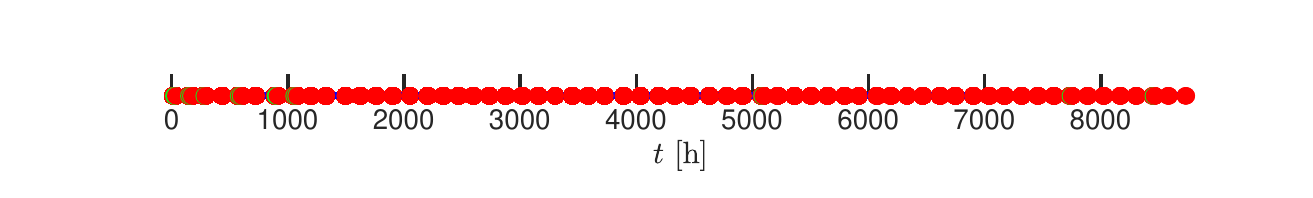}   
			
		}
		\textcolor{red}{$\bullet$}: ``normal'' decision, 
		\textcolor{green}{$\bullet$}: decision related to the validation, 
		\textcolor{blue}{\raisebox{1mm}{\rule{0.5cm}{0.05cm}}}: stand-by mode.  
		\captionof{figure}{\textbf{Study 3:} Decision dates \& stand-by mode}
		\label{fig:6_study_3__2_decision}
	\end{minipage}
\end{minipage} 

\noindent
\begin{minipage}[h]{\linewidth}
	\vspace*{0pt}
	{\centering	
		\includegraphics[width=6.675cm]{images/ch6/WO/Study_3/WO_LongSim_phi.pdf} 
		\includegraphics[width=6.675cm]{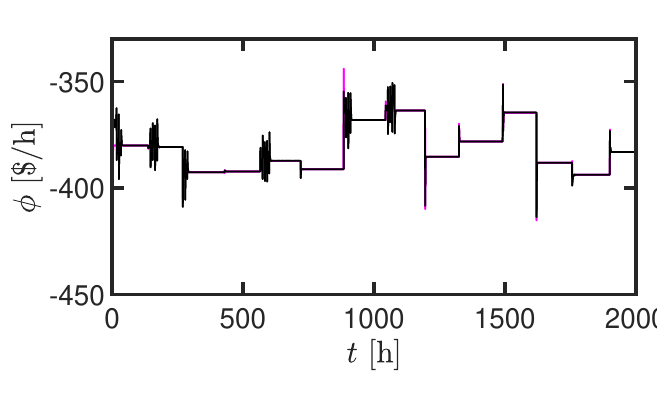} \\
	} 
	\textcolor{black}{\raisebox{1mm}{\rule{0.5cm}{0.05cm}}}: results obtained with the model and inaccurate measurements, \\
	\textcolor{magenta}{\raisebox{1mm}{\rule{0.5cm}{0.05cm}}}: results obtained with perfect model and measurements.
	\captionof{figure}{\textbf{Study 3:} Training \& trained phases}
	\label{fig:6_study_3__3_training}
	
	\medskip 
	
{\centering	
	\includegraphics[trim={1.5cm 0cm 1cm  0.8cm},clip,width=13.35cm]{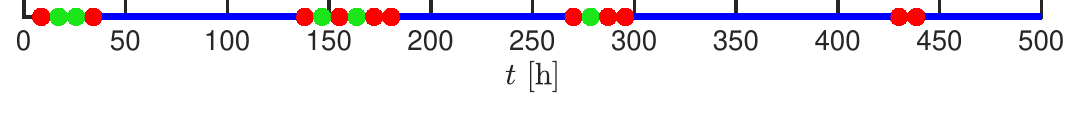}  \\
	\includegraphics[trim={1.5cm 0cm 1cm  0.8cm},clip,width=13.35cm]{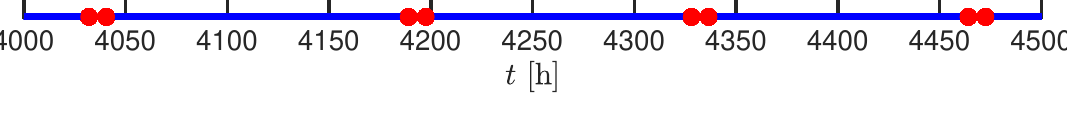}  \\ 
} 
\textcolor{red}{$\bullet$}: ``normal'' decision, 
\textcolor{green}{$\bullet$}: decision related to the validation, 
\textcolor{blue}{\raisebox{1mm}{\rule{0.5cm}{0.05cm}}}: stand-by mode.  
\captionof{figure}{\textbf{Study 3:} Decisions making  over the training \& trained phases}
\label{fig:6_study_3__4_decisions}
%
	\vspace*{0pt}
	{\centering	
		\includegraphics[width=6.675cm]{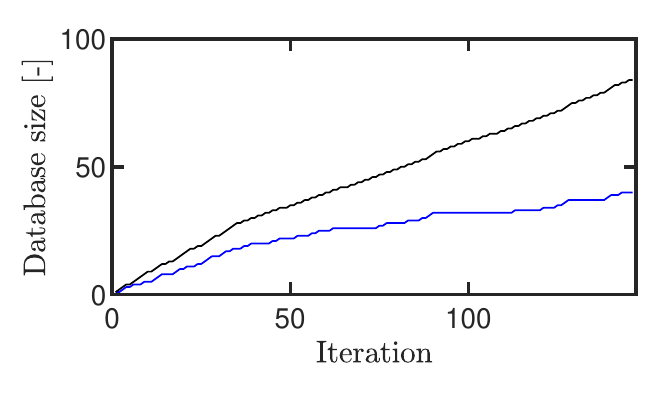}  \\ 
	} 

	\textcolor{black}{\raisebox{1mm}{\rule{0.5cm}{0.05cm}}}: amount of collected data, 
	\textcolor{blue}{\raisebox{1mm}{\rule{0.5cm}{0.05cm}}}: amount of data in $\mathcal{D}_{III}$.
	\captionof{figure}{\textbf{Study 3:} Database size w.r.t time \& trained phases}
	\label{fig:6_study_3__5_training}
%
	\vspace*{0pt}
	{\centering	
		\includegraphics[width=6.675cm]{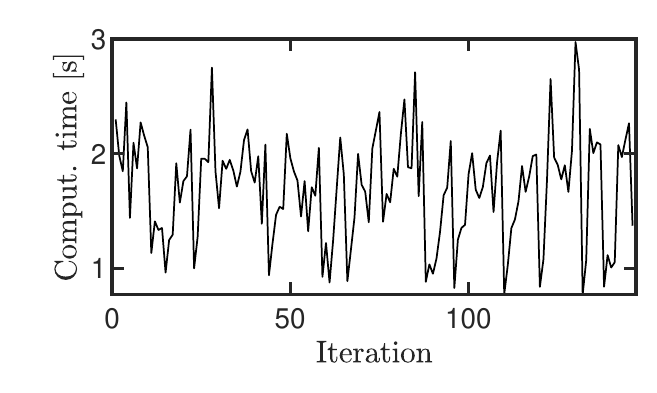}  \\   
	} 
	
	\captionof{figure}{\textbf{Study 3:} Computationnal time. It gathers everything ASP does at the end of an experiment from the statistical analysis up to the validation procedure passing by the database compression, the coherency check and the model-based optimisation solving.  
	}
	\label{fig:6_study_3__6_computationnal}
\end{minipage}

\section{Conclusion}

In this chapter improvements are proposed for each of the autopilot functions introduced in chapter~\ref{Chap:4_S_ASP} except for the signal purifier and the statistical analysis. Most of these improvements are illustrated on a relatively simple example, illustrating in a comprehensible way the functionalities of the new algorithm that has been built. There is indeed only one manipulated variable, one measured perturbation, and one un-measured perturbation. The non-illustrated features are the RTO-plant interface and the use of time as a measured perturbation. However, one believed that these two features are relatively well explained, easy to understand, and what they offer is not particularly subject to debate. Indeed, the RTO-plant interface corresponds to a pre-computation of what the RTO loop would do if it could be computed at a very high speed. So it only brings speed.  A preliminary  study on a fuel cell has shown that this interface brings significant improvements in the presence of continuous low amplitude disturbances. However, implementing the whole ASP on this case study was not done due to lack of time (this work will certainly appear in a future article together with everything that has been proposed in this chapter).

	\chapter{Conclusion}

\section{Summary}

In this work the challenges related to the monitoring of plants designed to operate at steady-state have been presented and it has been explained how a large number of tasks associated with this monitoring can be automated. One believes to have significantly improved a large part of the state of the art. 
In particular, chapters 2, 3 and 5 give clear answers to problems that until then were either unresolved or unidentified, be it the choice of the filtering gain, the model-based feasibility of the decisions that are made, or the correction structure that is used. 
In addition, chapters 4 and 6 give a list of functions that need to be fulfilled in order to run a plant efficiently and suggests ways to build these functions and interconnect them.
Many of these functions are completely new such as the ASP's RTO-plant interface, the experiment designer, the consistency monitor, the compressor, etc. introduced in chapter 6.

However, while this work provides many answers, it also raises many questions and opens research directions in RTO and in other fields related to optimization.

\section{Potentially interesting research directions}

\subsection{In RTO}

Here are seven potentially interesting research directions: 
\begin{enumerate}
	\item ASP with a correction structures of type A or B. One has started to study this research direction but it turns out that it is a very big piece of work that would probably require a complete PhD to be properly treated. Nevertheless, based on the preliminary study one has done, one can already indicate what the advantages and disadvantages would be of such an autopilot as well as the potential difficulties associated to its design. 
	
	Probable advantages: 
	\begin{itemize}
		\item One can expect a better database management. Indeed, the CTP would no longer be defined at the level of the whole plant, but rather at the level of the SMs. Therefore, when the plant's behavior changes, only the databases of the  SMs associated to the SPs whose behaviors have actually changed would be reset, whereas the other databases could be kept.
		\item One can expect better corrections. As it has been shown in chapter 5, the correction structures A and B provide better corrections than the D or I ones.   
		\item It opens the possibility to adapt the correction strategy to the plant-model mismatch of each SMs. Since the model and the plant are split into pieces that are SMs and SPs, one could consider   developing methods to identify the local nature  of the plant-model mismatch, i.e. at the level of the SMs and SPs. Then, for instance, one could implement parametric corrections where the SMs are only subject to parametric plant-model mismatch. 
	\end{itemize}

	Likely disadvantages/difficulties:  
	\begin{itemize}
		\item Implementation is clearly one of the main issue of distributed methods (using correction structures of type A of B). Indeed, on the one hand, the structures D and I only require one set of correction functions that are always ``plugged'' to the model the same way (their inputs are the plant's inputs, and their outputs are connected to either the plant's cost and constraints or the plant's outputs). While on the other hand, the structures A and B are  case dependent and may be complex networks of sets of functions that are not always easy to manage when it comes to implement them, i.e. writing the code.   
		\item The initialization may be less intuitive since Assumptions~\ref{ass:6_1_Erreur_fp_f_bornee} and \ref{ass:6_2_Meme_Courbusre} would have be made at the level of each SMs. 
		\item Computational difficulties may rise from the fact that if one evaluates the uncertainty on the outputs of each SMs, then to evaluate how those uncertainties diffuse and interact through the model may require a lot of computations. For instance, in the section~\ref{sec:6_6_2_Uncertaintyspreads} one shows how the uncertainty on the outputs of the model are affected when they pass through the cost function. Implementing the structures A or B means doing such computations for each SM simultaneously, which is a priori not and easy task. 
	\end{itemize}

	
	\item One could use  the measurement of the noise variance $\bm{N}_y$ to adjust the plant controller.  In this work one does not use $\bm{N}_y$, see Figure~\ref{fig:2_Convertisseur_Statistique}. However, when the plant is at steady-state, the quality of two different sets of settings of the controllers could be assessed by their efficiency at rejecting disturbances and measurement noise. Therefore, the better configuration would maintain the outputs of the plant closer to their set-point, i.e. the covariance matrix of the measurement $\bm{N}_y$ would have smaller eigenvalues. So, $\bm{N}_y$  could be used as an indicator of the quality of the control structure, and it could be used to maintain or even update it. 
	
	
	\item The problem of ``implicit'' constraints.  A steady state  model of a plant is generally capable of providing a prediction of the steady state of the plant for any input. However, one of the things that has been observed with the TE case study of chapter~\ref{Chap:5_IMA} is that for some tested inputs, the plant does not have a stable state. Either it oscillates or it diverges slowly until it activates an alarm causing the plant to stop. At the moment, no RTO algorithm using SS models (except perhaps fast RTO) can satisfactorily handle such cases.  One thinks that it would be important to find a simple and efficient protocol to manage this type of problem via, for example, the addition of a constraint associated to the existence of SS. 
	
	
	\item Towards a network of ASPs. The sharing of databases of several similar systems could be useful to  improve the ASP of each system, identify misunderstood behaviors on a large number of systems (e.g. degradation processes related to the use of equipment over the long term), etc. 
	
	
	 \item The analysis of chapter~\ref{Chap:5_IMA} showed that the correction structures A and B can, thanks to a measurement-based splitting of a model, enable better/deeper corrections. However, this result could also be used to motivate the placement of new sensors on the plant in order to improve its monitoring by enabling better splitting of its subparts.  
	 
	 
	 \item Sometimes the plant is affected by known but unmeasured disturbances. It could be very interesting to observe them to improve the management of databases (i.e. reset them less often), and to be able to manage their effects in real time with an improved VCL. However, given that the model changes at each iteration and that unknown unmeasured disturbances may be present, the approach to observe them does not seem so obvious. The design of such an observer is an open problem. 
	 	 
	 
	 \item The waiting time that is used at the end of each experiment to passively reduce the effects of measurement uncertainty could be chosen together with  $\bm{u}_{k+1}$, or even revised in real time. In fact, it is proposed to set this time to a sufficiently large value so that each statistical data on the outputs of the plant is accurate enough. However, nothing says that all experiments must provide accurate results. It is possible that having a rough idea of the output of the plant is enough to motivate a new decision. For instance, if an experiment violates a constraint, then knowing exactly to what extent it is violated may not be of the greatest interest, whereas immediate reaction to bring the plant back to a feasible point might be better. The addition of such a function to  ASP (i.e. which would be a revision of the experiment designer) would be interesting since it could accelerate its convergence on  $\bm{u}_p^{\star}$.
\end{enumerate}

\subsection{In other research areas}

\subsubsection{Theoretical RTO or Sequential Optimization}

One can notice that the mechanism of SQP is very similar to that of the theoretical RTO algorithms that have been built throughout this thesis. We claim that the filter selection mechanism introduced in Chapter 2 can be used to link SQP iterations together via parameter $a^*$ of  \eqref{eq:A_6_SQP_parametre_a_asteri} which seems to play a role similar to the one of KMFCaA's filter. So replacing the merit functions that are usually used to define $a^*$  by the ``simple'' filter update strategy of KMFCaA could be, in some particular cases to be defined, a pertinent idea.

\subsubsection{Theoretical RTO or ROM-Optimization}
\label{sec:7_x_OTR_theotique_et_ROM_Opt}

The problems that theoretical RTO  (see section~\ref{sec:2_1_Definition_PB_ORT_Theorique}) and  reduced-order-model based optimization (ROM-optimization)  try to solve are almost identical.   The only notable difference between these two research areas is of a semantic nature. Indeed, if their objectives are written, their similarities are obvious: 
\begin{itemize}
	\item The objective of the theoretical RTO is to identify the minimum of a plant using a minimal number of experiments, because each experiment has a cost in terms of time and money that one wants to limit as much as possible. To achieve this goal, an easy-to-evaluate model that can be called a low-fidelity model can be used to guide and improve the exploration of the hard-to-evaluate system that is the plant. 
	\item ROM-Optimization aims at identifying the minimum of a high-fidelity model by evaluating it in a minimal number of points because each simulation has a cost in terms of time that one wants to limit as much as possible.  To achieve this goal, an easy-to-evaluate model that can be called a low-fidelity model can be used to guide and improve the exploration of the hard-to-evaluate system that is the high-fidelity model.
\end{itemize}
Also, an additional objective of theoretical RTO is to maintain the plant within feasible operation state whereas ROM-Optimization seems to be less focus on this problematic (since its working with simulations and not real processes that can break).  

Due to the similarity between these two research areas it is not surprising that the methods developed to solve these problems are similar. It is therefore not surprising that an equivalent of ISO  \cite{Gao:05} developed by the RTO community in 2005 had already been proposed by the ROM-Optimization community 7 years earlier in \cite{Alexandrov:98}. This also means that any theoretical RTO method can be applied to ROM-Optimization.  It is therefore possible to use KMFCaA to identify the minimum of a high fidelity model.  

\begin{figure}[H]
	\centering
	\includegraphics[width=15cm]{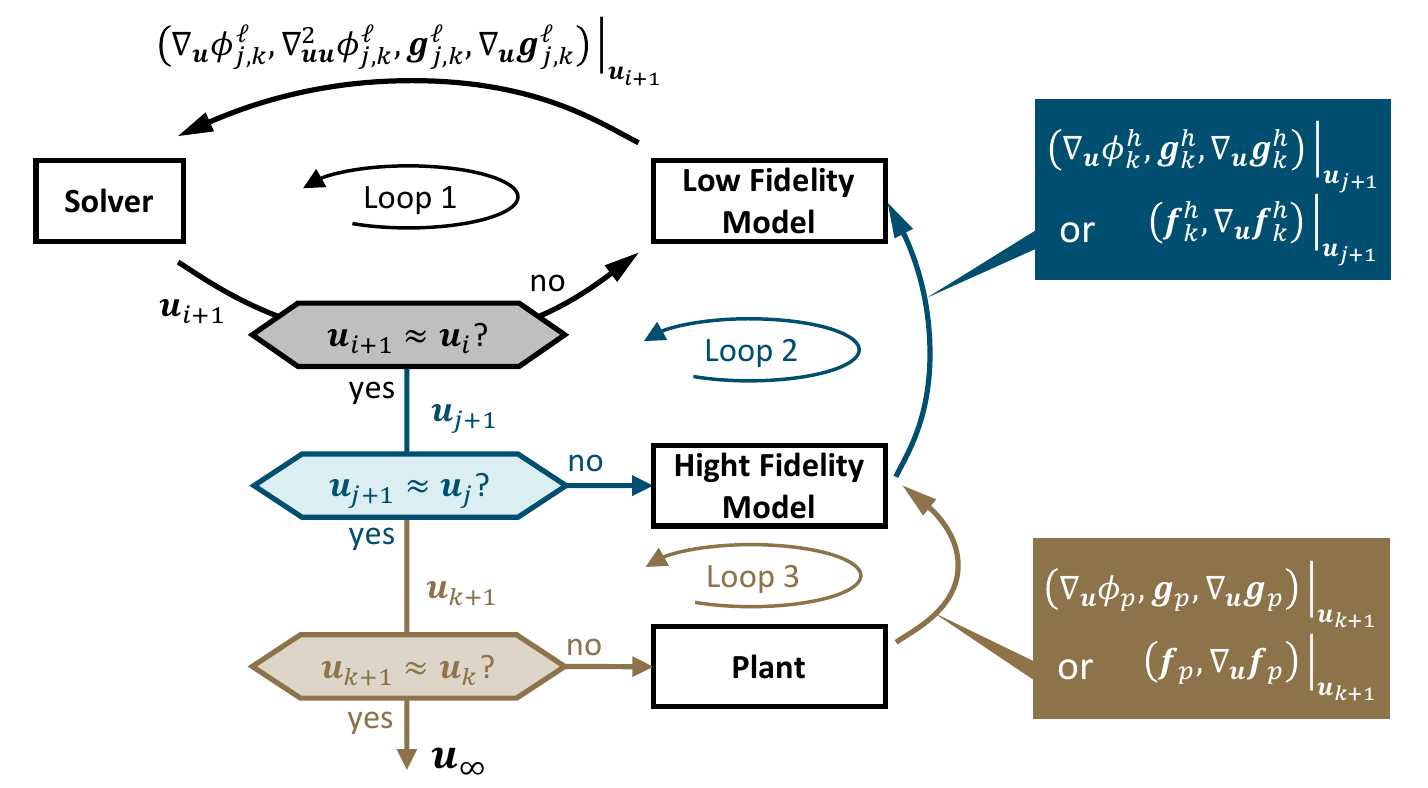}
	\caption{RTO with high- and low-fidelity models. The functions $(\phi^{h}_{k},\bm{g}^{h}_{k},\bm{f}^{h}_{k})$ are the functions of the updated high fidelity model $\bm{u}_k$ on the basis of experiments carried out on the plant.  The functions $(\phi^{\ell}_{j,k},\bm{g}^{\ell}_{j,k},\bm{f}^{\ell}_{j,k})$ are the functions of the low-fidelity model updated at $\bm{u}_j$ on the basis of simulations on the high-fidelity model updated at $\bm{u}_k$, i.e. $(\phi^{h}_{k},\bm{g}^{h}_{k},\bm{f}^{h}_{k})$.}
	\label{fig:rvfd}
\end{figure}

The concept can even be taken a step further and one can consider optimizing a plant using a high fidelity model. 
One could consider the implementation of a double KMFCaA loop to first optimize the high fidelity model by using a low fidelity model easy to evaluate. To then optimize the plant with a KMFA loop, where KMFA corresponds to a KMFCaA method from which one would remove the convexification step. Indeed, to perform this convexification step, one would have to evaluate the Hessian of the high fidelity model. For the same reasons that motivate us not to do it for the plant, one does not wish to do it for the high-fidelity model (requires too many evaluations of functions expensive to evaluate). Figure~\ref{fig:rvfd} shows what such a structure looks like. 

One might ask what is the point of such a structure when using a simple low fidelity model is (in theory) sufficient to find the optimum of the plant.  The interest could be that interposing a high-fidelity model between the low-fidelity model and the plant could reduce the number of experiments to be performed by making them more relevant. Indeed, the evaluations of the high-fidelity model and the plant are time-consuming, but only the evaluations of the plant can entail costs and real risks.

	\bibliographystyle{acm}
	\bibliography{Thesis} 
	\chapter{Nomenclature}

\noindent

\begin{tcolorbox}[
	blanker,
	breakable]
	\textbf{Acronyms (all):}
	\vspace{-\topsep}
	\begin{longtable}{p{2cm}p{13cm}}
		A      & A type of correction structure (see Figure~\ref{fig:5___3_Observation_Un_SM}),\\
		ASP    & Autopilot for steady-state process (designed to operate at SS), \\
		B      & A type of correction structure (see Figure~\ref{fig:5___3_Observation_Un_SM}),\\
		CSTR   & Continuous stirred-tank reactor, \\
		CTP    & Consistent time period, \\
		D      & A type of correction structure (see Figure~\ref{fig:5___3_Observation_Un_SM}),\\
		E$i$   & Event $i$, see section \ref{sec:4_4_1_Declencheur_de_decision}, \\
		FDM    & Fault-detection method, \\
		GP     & Gaussian process, \\
		h.o.t. & High order terms, \\
		KMFCA  & MFCA with filter-based constraints,\\
		KMFCaA & Modifier-filter-(active)curvature adaptation with filter-based constraints,\\
		I      & A type of correction structure (see Figure~\ref{fig:5___3_Observation_Un_SM}),\\
		ISO    & Iterative setpoint optimization, \\
		ISOPE  & Integrated system optimization and parameter estimation, \\
		KKT    & Karush–Kuhn–Tucker,\\
		MA     & Modifier adaptation, \\
		MAy    & Output modifier adaptation, \\
		MFCA   & Modifier-Filter-Curvature Adaptation, \\
		MFCQ   & Mangasarian-Fromovitz constraint qualification, \\
		MP-QP  & Multi-parametric quadratic program, \\
		NLP    & Nonlinear programming, \\
		LICQ   & Linear independence constraint qualification, \\
		QP     & Quadratic program,\\
		ROM    & Reduced-order-model (based optimization),\\
		RTO    & Real time optimization,\\
		SE-ARD & Squared-exponential with automatic relevance determination (kernel),\\
		SM     & Submodel, \\
		SOC    & Self optimizing control, \\
		SP     & Subplant, \\
		SQP    & Sequential quadratic programming, \\
		SS     & Steady-state, \\
		S-ASP  & Simple autopilot for steady-state process (designed to operate at SS), \\
		TS     & Two-step approach,\\ 
		VLC    & Virtual control layer,\\
		WO     & Williams-Otto,\\
	\end{longtable}
\end{tcolorbox}

\bigskip

\begin{tcolorbox}[
	blanker,
	breakable]
	\textbf{Acronyms (new algorithms only):}
	\vspace{-\topsep}
	\begin{longtable}{p{2cm}p{13cm}}
		MFCA      & Theoretical RTO method -- Modifier-Filter-Curvature Adaptation, (see chapter~\ref{Chap:2_Vers_Une_meilleure_Convergence}) --  MFCA guarantees the stability conditions hold (see Figure~\ref{fig:2___2_Domaines}),\\
		KMFCA     & Theoretical RTO method -- MFCA with filter-based constraints (see chapter~\ref{Chap:3_KMA}) --  KMFCA provides what MFCA provides plus better guarantees of constraints satisfaction,\\
		KMFCaA    & Theoretical RTO method -- Modifier-filter-(active)curvature adaptation with filter-based constraints (see chapter~\ref{Chap:3_KMA}) --  KMFCaA provides what KMFCA provides with a speed boost,\\
		S-ASP     & Practical RTO method -- Simple autopilot for steady-state processes (see chapter~\ref{Chap:4_S_ASP}) --  S-ASP is based on KMFCaA and has several additional features to handle some practical problems,\\
		ASP       & Practical RTO method -- Autopilot for steady-state processes (see chapter~\ref{Chap:6_ASP}) --   ASP is, to the author's knowledge, the first attempt to fully automate the supervision of a plant.
	\end{longtable}
\vspace{-\topsep}
\textit{Most of the time these algorithms acronyms are followed by a letter indicating which correction structure is employed. The four structures that are used in this thesis are referred to the letters D, I, A and B. Those structures are illustrated and compared on Figure~\ref{fig:5___3_Observation_Un_SM}.}
\end{tcolorbox}
\bigskip 

\noindent
\begin{tcolorbox}[
	blanker,
	breakable]
	\textbf{Typography (bold, capitals, etc):}
	\vspace{-\topsep}
	\begin{longtable}{p{1.5cm}p{13cm}}
		$x$         & $x$       is a scalar, \\
		$\bm{x}$    & $\bm{x}$  is a vector, \\
		$\bm{X}$    & $\bm{X}$  is a matrix, \\
		$X$         & $X$       is a sequence of \{scalars, gradients, or matrices\}, or a cluster,\\
	\end{longtable}
\end{tcolorbox}

\bigskip 

\noindent
\begin{tcolorbox}[
	blanker,
	breakable]
	\textbf{Superscripts \& Subscripts \& Accents:}
	\vspace{-\topsep}
	\setlength{\tabcolsep}{0.7em} 
	{\renewcommand{\arraystretch}{1.5}
		\begin{longtable}{p{1.5cm}p{13cm}}
		$(\cdot)_k$           & The variable $(\cdot)$ at iteration $k$, \\
		$(\cdot)_p$           & The variable $(\cdot)$ of the plant, \\
		$(\cdot)_{p,k}$       & The variable $(\cdot)$ of the plant at iteration $k$, \\
		$\widehat{(\cdot)}_p$ & The measurement of the plant's variable $(\cdot)$, \\
		$(\cdot)_{(i)}$       & $i$-th element of a the vector $(\cdot)$, \\
		$(\cdot)_{\{\ell\}}$  & The sequence of measurements of $(\cdot)$ associated to the $\ell$th experiment, \\
		$(\cdot)_{\{\ell|i\}}$  & The $i$th element of the sequence of measurements of $(\cdot)$ associated to the $\ell$th experiment, \\
		$(\cdot)_{[i]}$       & Refers to the cluster $i$, \\
		$\overline{(\cdot)}$  & Upper bound of $(\cdot)$, \\
		$\underline{(\cdot)}$ & Lower bound of $(\cdot)$, \\
		$\utilde{\text{\hskip 0.1ex $(\cdot)$}}$ & Estimated lower bound of $(\cdot)$,\\
		$\widetilde{g}$       & Generalized version of function $g$, see \eqref{eq:2_Generalized_version_function},\\
		$(\cdot)^{c}$         & Convex approximation of $(\cdot)$, \\
		$\bm{g}^{a}$          & The vector of active constraints, \\
		$(\cdot)^{h}$          & Function $(\cdot)$ associated to a high-fidelity model, \\
		$(\cdot)^{\ell}$          & Function $(\cdot)$ associated to a low-fidelity model, \\
		$u^{\star}$           & A local minimum, \\
		$u^{\bullet}$         & A KKT point, \\
	\end{longtable}
	}
\end{tcolorbox}

\bigskip 

\noindent
\begin{tcolorbox}[
	blanker,
	breakable]
	\textbf{Roman Symbols:}
	\vspace{-\topsep}
	\begin{longtable}{p{1.5cm}p{13cm}}
		$a$                      & Parameters characterizing a domain $\mathcal{A}$,\\
		$\bm{a}$                 & Convergence criterion \eqref{eq:4___27_Critere_arret_S_ASP},  \\
		$\bm{d}$                 & Disturbances vector, \\
		$\bm{d}^{\prime}$        & Measured disturbances vector, \\
		$\bm{d}^{\prime\prime}$  & Unmeasured disturbances vector, \\
		$\bm{f}$                 & Inputs-outputs mappings,\\ 
		$filt_k$                 & Function defined there: \eqref{eq:3___46_fonction_filt},\\
		$\bm{g}$                 & Constraint functions, \\
		$\bm{G}_k^{\star}$       & VCL optimal gain w.r.t. the updated model at iteration $k$, \\
		$h$                      & Constraint function of ASP, see \eqref{eq:6_30_h}, \\
		$k$                      & Usually refers to the $k$-th RTO iteration, \\ 
		$K$                      & Filter gain, \\
		$\k_{\f\f}$              & Kernel function characterizing the covariance of a function $\f$ with itself at two locations $\bm{x}_i$ and $\bm{x}_j$, \\
		$\k_{\f\f^1}$            & Kernel function characterizing the covariance of a function $\f$ with its derivative at two locations $\bm{x}_i$ and $\bm{x}_j$, \\
		$\k_{\f^1\f^1}$           & Kernel function characterizing the covariance of the derivative of a function $\f$ with itself at two locations $\bm{x}_i$ and $\bm{x}_j$, \\
		$\ell$                   & Can be the index of an experiment, or a length scale parameter (an hyperparameter of a GP), \\
		$\bm{\L}$                & Matrix of length scale parameters, \\
		$n_{(\cdot)}$            & Size of vector $(\cdot)$, \\
		$nab_k$                  & Function defined there: \eqref{eq:3___48_fonction_nab},\\
		$loss_{acc}$             & Acceptable loss\\,
		$\bm{N}_{y\{\ell\}}$     & Covariance matrix of the noise on the $\ell$th sequence of measurements of  $y$, \\
		$p(a)$                   & Probability of $a$, \\
		$p(a|b)$                 & Probability of $a$ given  $b$, \\
		$r$                      & Radius, \\
		$\bm{S}_{y\{\ell\}}$     & Covariance matrix  representing the uncertainty on the mean value of $\ell$th sequence of measurements of  $y$, \\
		$\bm{sol}$               & Solution of an optimization problem, \\
		$t$                      & The time (the age of the plant), \\
		$\bm{u}$                 & Manipulated variables (or decision variables),\\
		$\bm{v}$                 & Copy of the manipulated variables, or directional vector,\\
		$\bm{x}$                 & Inputs of the plant, \\
		$\bm{y}$                 & Outputs vector (or measured variables vector),\\
	\end{longtable}
\end{tcolorbox}

\bigskip 

\noindent
\begin{tcolorbox}[
	blanker,
	breakable]
	\textbf{Greek Symbols:}
	\vspace{-\topsep}
	\begin{longtable}{p{1.5cm}p{13cm}}
		$\bm{\Delta u_{scal}}$   & Scaling matrix, \\
		$\Delta y$               & Boundary on the modeling error, \\
		$\bm{\theta}$            & Parameters, \\
		$\bm{\lambda}$           & Vector of Lagrange multipliers, \\
		$\mu_k^{\bm{f}}$         & Correction function updated at the $k$th RTO iteration  and applied to the model's function $\bm{f}$,\\ 
		$\widehat{\nu}^2_{y,i}$  & The \textit{estimate} of the variance representing the noise on the plant output,  ,\\
		$\sigma_{xy,i}^2$           & The combined uncertainty on the plant outputs, see \eqref{eq:6_11_MeasurementsUncertainty},\\
		$\widehat{\sigma}_{y,i}^2$  & The \textit{estimate} of the variance representing the uncertainty on the mean of the plant output,\\           
		$\sigma_f^2$               & A priori variance of $f_p(\bm{x}_i$ around $f(\bm{x}_i)$ (hyperparameter of a GP), \\
		$\widehat{\bm{\Sigma}}^2_{x,i}$  & The \textit{estimate} of the covariance matrix representing the uncertainty on the mean of the plant inputs,\\ 
 		$\phi$                   & Cost, or objective,  function, \\
 		$\varphi$                & Objective  function of ASP, see \eqref{eq:6_29_varphi}, \\
	\end{longtable}
\end{tcolorbox}

\bigskip 

\noindent
\begin{tcolorbox}[
	blanker,
	breakable]
	\textbf{Calligraphic symbols:}
	\vspace{-\topsep}
	\begin{longtable}{p{1.5cm}p{13cm}}
		$\mathcal{A}$              & Domain over which one assumes a function is quasi-affine, \\
		$\mathcal{C}^0$            & Set of continuous functions, \\
		$\mathcal{C}^1$            & Set of continuous functions whose derivatives exist and are continuous, \\
		$\mathcal{C}^2$            &  Set of continuous functions whose derivatives and second derivatives exist and are continuous, \\
		$\mathcal{D}$              & A database, \\
		$\mathcal{D}_I$            & The \textit{raw} database containing all the SS measurements,\\
		$\mathcal{D}_{II}$         & The \textit{usable} database containing the results of statistical analyses of the measurement sequences associated to each experimented SS,\\
		$\mathcal{D}_{II}^{\bm{\checkmark}}$         
		                           & The \textit{validated} database which is the subset of $\mathcal{D}_{II}$ that is supposed to contain only up to date data,\\
		$\mathcal{D}_{III}$        & The \textit{compressed} database which is supposed to ``summarize'' $\mathcal{D}_{II}^{\bm{\checkmark}}$,\\
		$\mathcal{GP}$             & Gaussian process,\\
		$\mathcal{L}$              & A Lagrangian, \\
		$\mathcal{N}(a,b)$         & A normal distribution of mean $a$ and variance $b$, \\
		$\mathcal{O}$              & High order terms function, \\
	\end{longtable}
\end{tcolorbox}

 \bigskip 
 
\begin{table}[htbp]
	\textbf{Custom symbols}:\\
	\vspace{-\topsep}
	\begin{tabular}{p{0.05\textwidth}p{0.95\textwidth}}
		$\nabla^S$    & Directional derivative of a function $\bm{sol}$, \\
		$\Delta^S$   & Estimate of a directional derivative of a function $\bm{sol}$, \\
	\end{tabular}
\end{table}

\bigskip 

\noindent
\begin{table}[htbp]
	\textbf{Sets:}\\
	\vspace{-\topsep}
	\begin{tabular}{p{0.05\textwidth}p{0.95\textwidth}}
		$\amsmathbb{R}$        & Set of real numbers, \\
		$\mathbb{\Theta}$   & Set of parameters, \\
	\end{tabular}
\end{table}

\bigskip 

\noindent
\begin{table}[htbp]
	\textbf{Operators:}\\
	\vspace{-\topsep}
	\begin{tabular}{p{0.1\textwidth}p{0.9\textwidth}}
		$\nabla_x f$             & Gradient of $f$ w.r.t. $x$, \\
		$\nabla^2_{xx} f$        & Hessian of $f$ w.r.t. $x$, \\
		$\nabla_x f|_{x=x_i}$    & Gradient of $f$ w.r.t. $u$ evaluated at $x=x_i$, \\
		$\partial_x f$           & Partial derivative of $f$ w.r.t. $x$, \\
		$\partial^2_{xx} f$      & Partial second derivative of $f$ w.r.t. $x$, \\
		$\partial_x f|_{x=x_i}$  & Partial derivative of $f$ w.r.t. $u$ evaluated at $x=x_i$, \\
		$\amsmathbb{E}[(\cdot)]$ & Expectancy of $(\cdot)$, \\
		$\amsmathbb{V}[(\cdot)]$ & Variance of $(\cdot)$, \\
		$\amsmathbb{V}[(\cdot),(\star)]$ & Covariance of $(\cdot)$ and $(\star)$, \\
	\end{tabular}
\end{table}

	\appendix
	\renewcommand\chaptername{Appendix}

\chapter{Technical complements to Chapter 2}

\section{Correction of a composed  function}
\label{App:A___CorrectionFonctionComposite}

Let there be three $\mathcal{C}^2$ functions:
\begin{align*}
	\phi:\amsmathbb{R}^{n_u},\amsmathbb{R}^{n_y}\rightarrow \ &  \amsmathbb{R}, & 
	\bm{f}_p:\amsmathbb{R}^{n_u}\rightarrow \ &  \amsmathbb{R}^{n_y}, & 
	\bm{f}_{k+1}:\amsmathbb{R}^{n_u}\rightarrow \ &  \amsmathbb{R}^{n_y},
\end{align*}
and let's state that at the point $\bm{u}_{k}\in\amsmathbb{R}^{n_u}$ the function $\bm{f}_{k+1}$ is:
\begin{align} \label{eq:A___1_AffineCorrection}
	\mathcolor{green!50!black}{
		\bm{y}_k := \bm{f}_{k+1}|_{\bm{u}_{k}} :=
	} \ & 
	\mathcolor{green!50!black}{
		\bm{f}_p|_{\bm{u}_{k}} := \bm{y}_{p,k},
	} &
	\mathcolor{blue!50!black}{
		\nabla_{\bm{u}} \bm{f}_{k+1}|_{\bm{u}_{k}} := \ 
	}
	& 
	\mathcolor{blue!50!black}{
		\nabla_{\bm{u}} \bm{f}_p|_{\bm{u}_{k}}},
\end{align}
where the variables  $\bm{y}_k $ and $\bm{y}_{p,k}$ are introduced to make the following developments more readable.  Then, one can observe that the composed functions $\phi \circ \bm{f}_p$ and $\phi \circ \bm{f}_{k+1}$ have similar Taylor series expansions around $\bm{u}_{k}$: 
\begin{equation*}
	\begin{array}{r@{\ }l@{}c@{}c@{}l}
		\phi(\bm{u},\bm{f}_{k+1}(\bm{u})) = & \
		\phi|_{\bm{u}_{k}, \mathcolor{green!50!black}{ \bm{y}_{k} }} & \
		+ & \
		\big(
		\partial_{\bm{u}} \phi|_{\bm{u}_{k},\mathcolor{green!50!black}{
				\bm{y}_{k}} } + \partial_{\bm{y}}\phi|_{\bm{u}_{k},
			\mathcolor{green!50!black}{\bm{y}_{k}} }
		\mathcolor{blue!50!black}{
			\nabla_{\bm{u}} \bm{f}_{k+1}|_{\bm{u}_{k}}
		}
		\big)
		\left(
		\bm{u}-\bm{u}_{k}
		\right) & \ + 
		\mathcolor{red!50!black}{h.o.t.}, \\
		\phi(\bm{u},\bm{f}_{p}(\bm{u})) = & \
		\phi|_{\bm{u}_{k}, \mathcolor{green!50!black}{ \bm{y}_{p,k} }} & \
		+ & \
		\big(
		\partial_{\bm{u}} \phi|_{\bm{u}_{k},\mathcolor{green!50!black}{
				\bm{y}_{p,k}} } + \partial_{\bm{y}}\phi|_{\bm{u}_{k},
			\mathcolor{green!50!black}{\bm{y}_{p,k}} }
		\mathcolor{blue!50!black}{
			\nabla_{\bm{u}} \bm{f}_{p}|_{\bm{u}_{k}}
		}
		\big)
		\left(
		\bm{u}-\bm{u}_{k}
		\right) & \ +
		\mathcolor{red!50!black}{ h.o.t.}, \\
	\end{array}
\end{equation*}
where the effects of the  definitions \eqref{eq:A___1_AffineCorrection} are illustrated through the correspondence of the color code used (green and blue). If differences exist between $\phi \circ \bm{f}_p$ and $\phi \circ \bm{f}_{k+1}$,  then they reside in the higher order terms of this development (``$\mathcolor{red!50!black}{ h.o.t.}$'' -- where the color red is used to mark a source of the differences).  It goes without saying that if the functions $\bm{f}_{k+1}$ and $\bm{f}_p$  are affine in $\bm{u}$ and $\bm{y}$, then $\bm{f}_{k+1}=\bm{f}_p$ and $(\phi \circ \bm{f}_p) = (\phi \circ \bm{f}_{k+1}) $. 

\begin{Remark}
	\textbf{(The deep effects of the indirect approach)}
	One can write the Hessian of $(\phi \circ \bm{f}_{k+1})$ at $\bm{u}_{k}$ in the following way: 
	\begin{align*}
		& \nabla_{\bm{uu}}(\phi\circ \bm{f}_{k+1})|_{\bm{u_k},\mathcolor{green!50!black}{\bm{y}_k} } = ...\\
		&  \qquad \Big[
		\partial_{\bm{u}}\phi(\partial_{\bm{u}}\phi)^{\rm T} + 2\partial_{\bm{uy}}\phi \mathcolor{blue!50!black}{\nabla_{\bm{u}} \bm{f}_{k+1}} + 
		(
		\mathcolor{blue!50!black}{\nabla_{\bm{u}} \bm{f}_{k+1}}
		)^{\rm T}
		\partial_{\bm{yy}}\phi
		\mathcolor{blue!50!black}{\nabla_{\bm{u}} \bm{f}_{k+1}} + ... \\
		& \qquad  \sum_{i=1}^{n_y} \left[
		\partial_{y_{(i)}}\phi \mathcolor{red!50!black}{\nabla_{\bm{uu}}\bm{f}_{k+1}}
		\right] 
		\Big]\Big|_{\bm{u}_k, \mathcolor{green!50!black}{\bm{y}_k}}.
	\end{align*}
	In the case of direct approaches, such that the equalities \eqref{eq:A___1_AffineCorrection} are not guaranteed, all green and blue terms in this equation are wrong. Specifically, several things could be sources of error in the estimation of $(\nabla_{\bm{uu}}(\phi\circ \bm{f}_{p})|_{\bm{u_k},\bm{y}_k})$ with $(\nabla_{\bm{uu}}(\phi\circ \bm{f}_{k+1})|_{\bm{u_k},\bm{y}_k})$:
	\begin{itemize}
		\item The functions $\nabla_{\bm{u}}\bm{f}_{k+1}$ would be structurally wrong, therefore \textit{a priori} $$\mathcolor{red!50!black}{\nabla_{\bm{u}}\bm{f}_{k+1} \neq \nabla_{\bm{u}}\bm{f}_{p}} ,\qquad  \forall \bm{u}\in\amsmathbb{R}^{n_u}$$ with a few potential exceptions. 
		\item The derivatives and Hessians of the cost function, although structurally correct, would not be evaluated at the right point if $\mathcolor{red!50!black}{\bm{y}_k\neq\bm{y}_{p,k}}$:
		\begin{align*}
			\partial_{\bm{uu}}\phi|_{\bm{u}_k,\mathcolor{red!50!black}{\bm{y}_k}} \neq \ & \partial_{\bm{uu}}\phi|_{\bm{u}_k,\mathcolor{red!50!black}{\bm{y}_{p,k}}}, & 
			\partial_{\bm{uy}}\phi|_{\bm{u}_k,\mathcolor{red!50!black}{\bm{y}_k}} \neq \ & \partial_{\bm{uy}}\phi|_{\bm{u}_k,\mathcolor{red!50!black}{\bm{y}_{p,k}}}, \\ 
			\partial_{\bm{yy}}\phi|_{\bm{u}_k,\mathcolor{red!50!black}{\bm{y}_k}} \neq \ & \partial_{\bm{yy}}\phi|_{\bm{u}_k,\mathcolor{red!50!black}{\bm{y}_{p,k}}}, &
			\partial_{y_{(i)}}\phi |_{\bm{u}_k,\mathcolor{red!50!black}{\bm{y}_k}} \neq \ & \partial_{y_{(i)}}\phi |_{\bm{u}_k,\mathcolor{red!50!black}{\bm{y}_{p,k}}}, \ \forall i =1,...,n_y.
		\end{align*}
	\end{itemize}
	It is therefore clear that choosing an indirect approach that brings the equalities \eqref{eq:A___1_AffineCorrection} allows to partially correct the high order terms (``$\mathcolor{red!50!black}{ h.o.t.}$''). However, it is less clear that these partial corrections ultimately provide a real benefit. A realistic case study (the WO reactor) and a statistical mathematical analysis show that a slight improvement can be expected \cite{Papasavvas:2019a}.
\end{Remark}

\section{Solution of the quadratic optimization problem}

\subsection{Solving a quadratic problem with \textit{equality} constraints}
\label{sec:A_2_1_Resoudre_QP_Contraintes_Egalite}

Consider a quadratic optimization problem with equality constraints: 
\begin{align} 
	\bm{u}_{k+1} = \operatorname{arg}
	\underset{\bm{u}}{\operatorname{min}} \quad &  q_0 + \bm{q}^{\rm T} \bm{u} + \frac{1}{2} \bm{u}^{\rm T} \bm{Q} \bm{u}, \nonumber \\
	\text{s.t.} \quad & \bm{h}^a + \bm{H}^a \bm{u} = \bm{0}.  \label{eq:A_2_QP_egalite}
\end{align}
This problem has an analytical solution which is (see \cite{Nocedal:2006} page 451): 
\begin{equation} \label{eq:A_3_sol_Analytique}
	\left( 	\begin{array}{c}
		\bm{u}_{k+1} \\ 
		\bm{\lambda}
	\end{array} \right) 
	= 
	-
	\left(\begin{array}{cc}
		\bm{Q}      & \bm{H}^{a\rm T} \\ 
		\bm{H}^{a}  & \bm{0}
	\end{array}\right)^{-1}
	\left(\begin{array}{c}
		\bm{q} \\ \bm{h}^a
	\end{array}\right),
\end{equation}
where $\bm{\lambda}$ is the vector of Lagrange multipliers associated with each of the equality constraints of \eqref{eq:A_2_QP_egalite}. Starting from this result, let's see how one can solve the cases where the constraints are inequalities instead of equalities.

\subsection{Solving a quadratic problem with \textit{inequality} constraints}
\label{sec:A_2_2_Resoudre_QP_Contraintes_Inegalite}

Consider a quadratic optimization problem with inequality constraints:
\begin{align} 
	\bm{u}_{k+1} = \operatorname{arg}
	\underset{\bm{u}}{\operatorname{min}} \quad &  q_0 + \bm{q}^{\rm T} \bm{u} + \frac{1}{2} \bm{u}^{\rm T} \bm{Q} \bm{u}, \nonumber \\
	\text{s.t.} \quad & \bm{g}(\bm{u})  = \bm{h} + \bm{H} \bm{u} \leq \bm{0}.  \label{eq:A_2_QP_inegalite}
\end{align}
One idea to solve this problem is to identify the constraints that are active at  $\bm{u}_{k+1}$ to replace the inequality constraints of \eqref{eq:A_2_QP_inegalite} by equality constraints on, and only on, these active constraints. 

This work of identifying active constraints can be done in the following way:

\textbf{(Step 0 -- Initialization}) First, a point $\bm{u}_{0}^{test}$ must be chosen. Then, the active constraints at $\bm{u}_0^{test}$ are identified and used to define a set of active constraints: $\mathcal{A}^{test}$.  Then one enters in a cycle of update of this set  $\mathcal{A}^{test}$. One initializes $j=0$ and goes to Step 1: 

\textbf{(Step 1 -- Computation of $\bm{u}^{test}$)} Consider a quadratic optimization problem with inequality constraints:
\begin{align} 
	\bm{u}_{j+1}^{test} = \operatorname{arg}
	\underset{\bm{u}}{\operatorname{min}} \quad &  q_0 + \bm{q}^{\rm T} \bm{u} + \frac{1}{2} \bm{u}^{\rm T} \bm{Q} \bm{u}, \nonumber \\
	\text{s.t.} \quad & \bm{g}^{test}(\bm{u})  = \bm{h}^{test} + \bm{H}^{test} \bm{u} = \bm{0},  \label{eq:A_2_QP_inegalite_test}
\end{align}
where  $\bm{g}^{test} = [g_{(i)}]_{i\in\mathcal{A}^{test}}$ a vector concatenating the constraints that one believes are active at $\bm{u}_{k+1}$. The solution of this problem is "easy" to calculate with the analytical solution \eqref{eq:A_3_sol_Analytique}. 

\textbf{(Step 2 -- Test 1)} Check if at the point $\bm{u}_{j+1}^{test}$ the set of constraints $\bm{g}\leq 0$, i.e. not only $\bm{g}^{test}$. If this test is passed, then the constraints that are active at $\bm{u}_{k+1}$ are part of $\mathcal{A}^{test}$ and one can move on to Step 3. If this is not the case, then at the point  $\bm{u}^{test}$ one violates a constraint that was supposedly inactive according to $\bm{u}_j^{test}$. One will therefore explore the segment linking $\bm{u}_j^{test}$ to $\bm{u}_{j+1}^{test}$ to find the closest point to $\bm{u}_{j+1}^{test}$ that does not violate the constraints. Then, one takes the active constraints at this point to redefine  $\mathcal{A}^{test}$, increments the counter $j$, and returns to Step 1.

\textbf{(Step 3 -- Test 2)} If one reaches this step it is that a priori $\mathcal{A}^{test}$  contains \textit{enough} constraints so that $\bm{u}_{j+1}^{test}$ is feasible. Now let's check if $\mathcal{A}^{test}$  contains \textit{too many} constraints. To carry out this verification, it is sufficient to check whether the vector of Lagrange multipliers calculated at the same time as $\bm{u}_{j+1}^{test}$ with \eqref{eq:A_3_sol_Analytique} are all $\geq 0$. If this is not the case, then one must remove from  $\mathcal{A}^{test}$ all constraints associated with Lagrange multipliers that are $<0$, increment the counter $j$, and return to Step 1. Otherwise, one moves to step 4. 

\textbf{(Step 4 -- Finalization)} If both tests are passed, then $\bm{u}^{test}$ is the solution we are looking for, i.e. $\bm{u}_{j+1}^{test}=\bm{u}_{k+1}$.

Finally, the real challenge associated with solving a quadratic problem with inequality constraints is the identification of the active constraints. 

\begin{Remark} \label{rem:A_Risque_viole_contrainte_1}
	(Risk of constraint violation part 1)
	It is clear that in order to identify these active constraints one has to implement a strategy that may involve the repetitive and/or frequent violation of constraints, especially during Step 2 with the exploration of the segment $[\bm{u}_j^{test},\bm{u}_{j+1}^{test}]$.
\end{Remark}

On the basis of this strategy, let us introduce the operating principles of the SQP (sequential quadratic programming) nonlinear optimization problem solver.

\subsection{Solving a non-linear optimization problems (SQP method)}

Consider a nonlinear optimization problem with inequality constraints: 
\begin{align} 
	\bm{u}_{k+1} = \operatorname{arg}
	\underset{\bm{u}}{\operatorname{min}} \quad  \phi(\bm{u}), \quad 
	\text{s.t.} \quad  \bm{g}(\bm{u})  \leq \bm{0}.  \label{eq:A_4_NL_inegalite}
\end{align}
One idea to solve this problem would be to choose a candidate $\bm{u}_0^{*}$ that one hopes to be close to $\bm{u}_{k+1}$, and to repeat   $\forall \ell\rightarrow \infty$ the following steps:
\begin{itemize}
	\item \textbf{(Step A)} Replace the functions $\phi$ and  $\bm{g}$  by quadratic and linear approximations of themselves at $\bm{u}_\ell^*$,
	\item \textbf{(Step B)} Solve the quadratic problem with inequality constraints based on these approximations as seen in the previous section to find $\bm{u}_{\ell+1}$
	\item \textbf{(Step C)} Select a scalar $a^*$ so that the functions $\phi$ and $\bm{g}$ are significantly ``better'' at
	\begin{equation} \label{eq:A_6_SQP_parametre_a_asteri}
		\bm{u}_{\ell+1}^{*}=\bm{u}_\ell^{*} + a^*(\bm{u}_{\ell+1}-\bm{u}_\ell^{*})
	\end{equation}
	than at $\bm{u}_{\ell}^{*}$. If $\bm{u}_{\ell+1}^{*}= \bm{u}_\ell^{*}$ stop the procedure and return $\bm{u}_{k+1} = \bm{u}_\ell^{*}$. Otherwise, increment  $\ell$ , and return to \textbf{Step A}.
\end{itemize}

The challenge of SQP is to ensure (i) progress towards the desired solution $\bm{u}_{k+1}$ at each iteration, (ii) that the convergence on   $\bm{u}_{k+1}$ will be observed in a limited time (or limited number of iterations), as well as (iii) the set of challenges associated with solving a quadratic problem with inequality constraints.

\begin{Remark} \label{rem:A_Risque_viole_contrainte_2}
	(Risk of constraint violation part 2)
	As at each step SQP uses an approximation of the NLP to progress. The intermediate points ($\bm{u}_{\ell}^*$) only satisfies the \textit{approximated} constraints of this NLP, and this satisfaction obviously does not imply the satisfaction of the \textit{actual} constraints.  Constraints violations can then be observed during these intermediate steps and/or during the choice of $a^*$.
\end{Remark}

The Remarks \ref{rem:A_Risque_viole_contrainte_1} and \ref{rem:A_Risque_viole_contrainte_2} show that a solver like SQP presents risks of constraint violation at different levels of its exploration strategy towards the optimum of an optimization problem. Moreover, at each Step A of SQP, one has to obtain the Hessian of the cost function of the problem one wants to solve. This is in addition to the potentially large number of evaluations required for the calculation of $a^*$  and for the convergence on $\bm{u}_{k+1}$, which ultimately results in SQP being not necessarily implementable (as such) on a real plant.

\section{Proof of Theorem \ref{thm:2___3_Proprietes_Sol}} 
\label{App:A_3_Preuve_THM2_3}

\subsection{Statement A: Global shape of \textit{sol}}
\label{App:A___Fonction_globalement_C0_et_Cn_par_morceaux}

Let's rewrite the optimization problem \eqref{eq:2___27a_ModelBasedPB} in the following way: 
\begin{align} \label{eq:A___2_OptPBExtended}
	\bm{u}^{*}_{k+1}:= \bm{sol}(\bm{v}_k) = \operatorname{arg}
	\underset{\bm{u}}{\operatorname{min}} \quad &  \phi_k(\bm{u}) := \widetilde{g}_{(0)}(\bm{u},\bm{v}_k), \\ 
	\text{s.t.} \quad & \bm{g}_k(\bm{u}) := \widetilde{\bm{g}}(\bm{u},\bm{v}_k ) \leq \bm{0}.   \nonumber   
\end{align}	
where the variable $\bm{v}\in\amsmathbb{R}^{n_u}$ is introduced as a copy of $\bm{u}$ indicating the point at which the model is updated. The functions of this problem can be approximated around  $(\bm{u}^*_{k+1},\bm{v}_k)$ via the following Taylor series $\forall i = 0, ..., n_g$:
\begin{align*}
	\widetilde{g}_{(i)}(\bm{u},\bm{v})  = \ &  \widetilde{g}_{(i)}|_{*} + 
	\partial_{\bm{v}}\widetilde{g}_{(i)}|_{*}^{\rm T}(\bm{v}-\bm{v}_k)  + 
	\partial_{\bm{u}}\widetilde{g}_{(i)}|_{*}^{\rm T}(\bm{u}-\bm{u}_{k+1}^{*}) +  ... \nonumber \\ 
	& \frac{1}{2} (\bm{v}-\bm{v}_k)^{\rm T}
	\partial_{\bm{vv}}\widetilde{g}_{(i)}|_{*}
	(\bm{v}-\bm{v}_k) + 
	\frac{1}{2} (\bm{u}-\bm{u}_{k+1}^{*})^{\rm T}
	\partial_{\bm{uu}}\widetilde{g}_{(i)}|_{*}
	(\bm{u}-\bm{u}_{k+1}^{*}) + ... \nonumber \\
	&
	(\bm{u}-\bm{u}_{k+1}^{*})^{\rm T}
	\partial_{\bm{uv}}\widetilde{g}_{(i)}|_{*}
	(\bm{v}-\bm{v}_k) + h.o.t..
\end{align*}
Which can be rewritten as: 
\begin{align} \label{eq:A___3_FonctionsParametriseesEnV}
	\widetilde{g}_{(i)}(\bm{u},\bm{v}) = \ &  
	\big[ a_i + \bm{b}_i^{\rm T} \bm{v} + \frac{1}{2}\bm{v}^{\rm T} \bm{C}_i \bm{v} \big] +  \big[\bm{d}_i^{\rm T} + \bm{v}^{\rm T}\bm{E}_i\big] \bm{u} + \frac{1}{2}\bm{u}^{\rm T} \bm{F}_i \bm{u}.
\end{align}
where the details are
\begin{align*}
	a_i := \ & 
	\widetilde{g}_{(i)}|_{*} -
	\partial_{\bm{v}}\widetilde{g}_{(i)}|_{*}^{\rm T}\bm{v}_{k}  -
	\partial_{\bm{u}}\widetilde{g}_{(i)}|_{*}^{\rm T}\bm{u}_{k+1}^{*} +
	\frac{1}{2} \bm{v}_k^{\rm T}\partial_{\bm{vv}}\widetilde{g}_{(i)}|_{*}\bm{v}_k + ... \\
	&
	\frac{1}{2} (\bm{u}_{k+1}^{*})^{\rm T}\partial_{\bm{uu}}\widetilde{g}_{(i)}|_{*}\bm{u}_{k+1}^{*} + 
	(\bm{u}_{k+1}^{*})^{\rm T}\partial_{\bm{uv}}\widetilde{g}_{(i)}|_{*}\bm{v}_{k}, \\
	\bm{b}_i := \ & 
	\partial_{\bm{v}}\widetilde{g}_{(i)}|_{*}^{\rm T} - 
	\partial_{\bm{vv}}\widetilde{g}_{(i)}|_{*}\bm{v}_k - \partial_{\bm{uv}}\widetilde{g}_{(i)}|_{*}\bm{u}_{k+1}^*, \\
	\bm{C}_i := \ & \partial_{\bm{vv}}\widetilde{g}_{(i)}|_{*}, \\
	\bm{d}_i := \ & \partial_{\bm{u}}\widetilde{g}_{(i)}|_{*} -  \partial_{\bm{uu}}\widetilde{g}_{(i)}|_{*} \bm{u}_{k+1}^* - \partial_{\bm{uv}}\widetilde{g}_{(i)}|_{*}\bm{v}_{k}, \\
	\bm{E}_i  := \ & \partial_{\bm{vu}}\widetilde{g}_{(i)}|_{*}, \\
	\bm{F}_i  := \ & \partial_{\bm{uu}}\widetilde{g}_{(i)}|_{*}.
\end{align*}
If one replaces the functions of \eqref{eq:A___2_OptPBExtended} by their convex approximations \eqref{eq:A___3_FonctionsParametriseesEnV}, then one obtains the following multi-parametric quadratic program (mp-QP):
\begin{align} \label{eq:A___4_mpQP}
	\bm{u}^{*}= \bm{sol}(\bm{v}) = \operatorname{arg} \ &
	\underset{\bm{u}}{\operatorname{min}} \quad  \big[ \bm{d}_0^{\rm T} + \bm{v}^{\rm T}\bm{E}_0 \big]  \bm{u} + \frac{1}{2}\bm{u}^{\rm T} \bm{F}_0 \bm{u}  \\ 
	\text{s.t.} \quad & \big[ a_i + \bm{b}_i^{\rm T} \bm{v} + \frac{1}{2}\bm{v}^{\rm T} \bm{C}_i \bm{v} \big] +  \big[\bm{d}_i^{\rm T} + \bm{v}^{\rm T}\bm{E}_i\big] \bm{u}  = \bm{0}, \quad  \forall i \in \mathcal{A}(\bm{v}), \nonumber    
\end{align}	
where $\mathcal{A}(\bm{v})$ is the subset of the constraints of \eqref{eq:A___4_mpQP} that are active at $\bm{u}^*$ when the model is updated at $\bm{v}$. By definition,  \eqref{eq:A___4_mpQP} correctly approximates the solution of \eqref{eq:A___2_OptPBExtended} when $\bm{v}$ is close to $\bm{v}_k$. The explicit solution of this problem is (according to \eqref{eq:A_3_sol_Analytique}): 
\begin{equation} \label{eq:A___5_SolExplicite}
	-\left(\begin{array}{cc}
		\bm{F}_0 & \bm{DE}_{\mathcal{A}(\bm{v})}^{\rm T} \\
		\bm{DE}_{\mathcal{A}(\bm{v})} & \bm{0}
	\end{array}\right) 
	\left(\begin{array}{c}
		\bm{sol}(\bm{v}) \\
		\bm{\lambda}(\bm{v})
	\end{array}\right) =
	\left(\begin{array}{c}
		(\bm{d}_0^{\rm T} + \bm{v}^{\rm T}\bm{E}_0) \\
		 \bm{A}_{\mathcal{A}(\bm{v})}
	\end{array}\right)
\end{equation}
where: 
\begin{itemize}
	\item $\bm{\lambda}(\bm{v})$ is a function that returns the Lagrange multiplier vector of the active constraints at the minimum of the model updated at $\bm{v}$,
	\item $\bm{DE}_{\mathcal{A}(\bm{v})} :=  [\bm{d}_i^{\rm T} + \bm{v}^{\rm T}\bm{E}_i]_{\forall i \in \mathcal{A}(\bm{v})},$ and 
	$\bm{A}_{\mathcal{A}(\bm{v})} :=  \big[ a_i + \bm{b}_i^{\rm T} \bm{v} + \frac{1}{2}\bm{v}^{\rm T} \bm{C}_i \bm{v}
	\big]_{\forall i \in \mathcal{A}(\bm{v})}$, where $[(\cdot)_i]_{\forall i \in \mathcal{A}(\bm{v})}$ is the ordered concatenation of $(\cdot)_i$ whose $i\in \mathcal{A}(\bm{v})$. For example, if only the constraints $i=1,2,4$ are active at $\bm{u}^*$ then:
	\begin{equation*}
		\bm{DE}_{\mathcal{A}(\bm{v})} =  [\bm{d}_i^{\rm T} + \bm{v}^{\rm T}\bm{E}_i]_{\forall i \in \mathcal{A}(\bm{v})} = \left( \begin{array}{c}
			\bm{d}_1^{\rm T} + \bm{v}^{\rm T}\bm{E}_1 \\
			\bm{d}_2^{\rm T} + \bm{v}^{\rm T}\bm{E}_2\\
			\bm{d}_4^{\rm T} + \bm{v}^{\rm T}\bm{E}_4
		\end{array}\right).
	\end{equation*}
\end{itemize}

What one can learn from \eqref{eq:A___5_SolExplicite} is that $\bm{sol}$ and $\bm{\lambda}$ are functions of $\bm{v}$ which are globally continuous ($\mathcal{C}^0$) and at least piecewise-$\mathcal{C}^1$, where each piece is defined by a different set of active constraints $\mathcal{A}(\bm{v})$. 
This result is illustrated in the following example:

\begin{exbox}
		One considers the following RTO problem: 
		\begin{align*}   
			\phi(u,y) := \ & y, \ & 
			y_p = f_p(u) := \ & \frac{(u-1)^2}{3} + \frac{(u-1)^4}{60} ,  \    & 
			y = f(u)  :=  \ & \frac{u^2}{4}.
		\end{align*}
		Subject to the following two cases of constraints:
		\begin{align*}  
			& \text{\textbf{Scenario 1}: $-5 \leq u \leq 5$, i.e.:} &  & g_{(1)}  = -5-u,&  & g_{(2)}  = u-5, \\
			& \text{\textbf{Scenario 2}: $\phantom{-}1 \leq u \leq 5$, i.e.:} &  & g_{(1)}  = \phantom{-}1-u,&  & g_{(2)}  = u-5.
		\end{align*}
		If one implements ISO-D/I, then the modified optimization problem at a point $v$ is according to the notations of \eqref{eq:A___2_OptPBExtended}: 
		\begin{align*}
			sol(v) = \operatorname{arg}
			\underset{u}{\operatorname{min}} \quad &  \widetilde{g}_{(0)}(u,v), \\ 
			\text{s.t.} \quad &  \widetilde{\bm{g}}(u,v) \leq \bm{0}.   \nonumber    
		\end{align*}
		where: 
		\begin{align*}
			\widetilde{g}_{(0)}(u,v) := \ & \frac{(v-1)^2}{3} + \frac{(v-1)^4}{60} + \left( \frac{2(v-1)}{3} + \frac{(v-1)^3}{15}
			\right) u +  \frac{u^2}{2}, \\
			\widetilde{g}_{(1)}(u,v) := \ &  g_{(1)}(u), \\
			\widetilde{g}_{(2)}(u,v)  :=  \ &  g_{(2)}(u). 
		\end{align*}
		One will not try to explicitly compute $sol(v)$ or the values of the Lagrange multipliers $\lambda_{(1)}(v)$ and  $\lambda_{(2)}(v)$ associated to the constraints $\widetilde{g}_{(1)}$ and $\widetilde{g}_{(2)}$. Instead one uses a solver to evaluate their values at a large number of points. The results obtained for each scenario are shown in Figure~\ref{fig:A___1_SolU_LambdaU_CActive}. These results confirm and illustrate that the functions  $(sol,\lambda_{(i)} \ \forall i = 1, ..., n_g)$  are globally $\mathcal{C}^0$ and at least piecewise-$\mathcal{C}^1$, where each piece is associated with a different set of active constraints.  Scenario 1 shows a case where these functions are at least $\mathcal{C}^1$ at  $v=u_p^{\star} =1$ and the Scenario 2 shows a case where these functions are only $\mathcal{C}^0$ at $v=u_p^{\star} =1$. 
		
		\begin{minipage}[h]{\linewidth}
			
			\vspace*{0pt}
			\includegraphics[width=4.45cm]{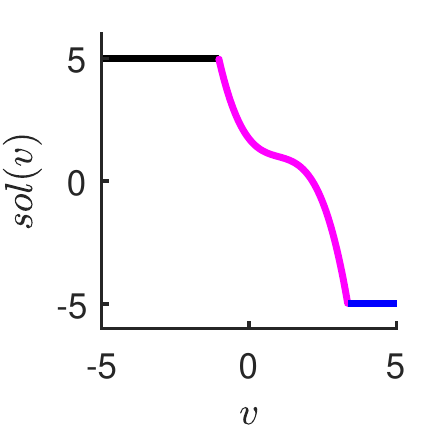}\hskip -0ex
			\includegraphics[width=4.45cm]{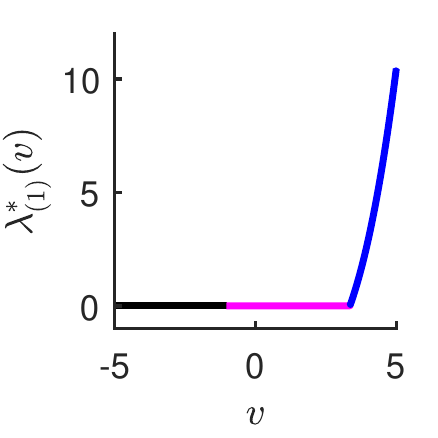}\hskip -0ex
			\includegraphics[width=4.45cm]{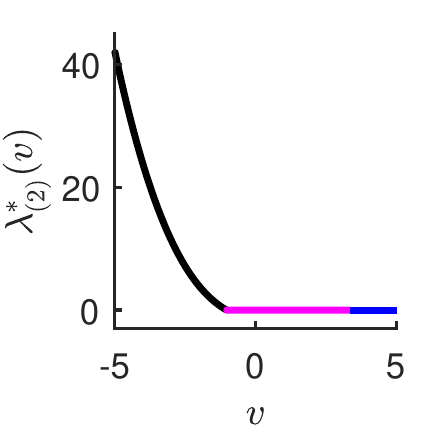}
		\begin{center}
			a) Scenario 1
		\end{center}

			\vspace*{0pt}
			\includegraphics[width=4.45cm]{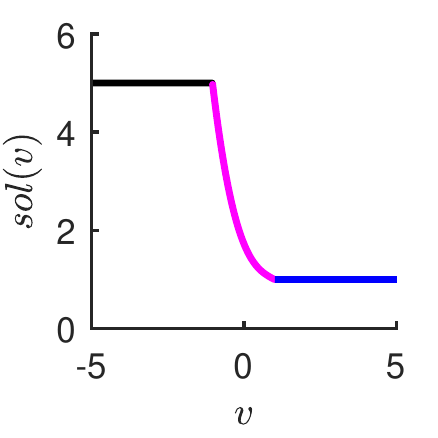}\hskip -0ex
			\includegraphics[width=4.45cm]{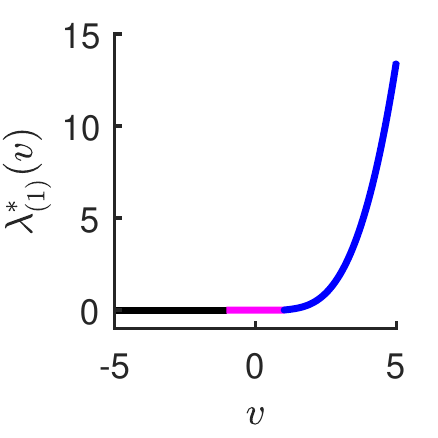}\hskip -0ex
			\includegraphics[width=4.45cm]{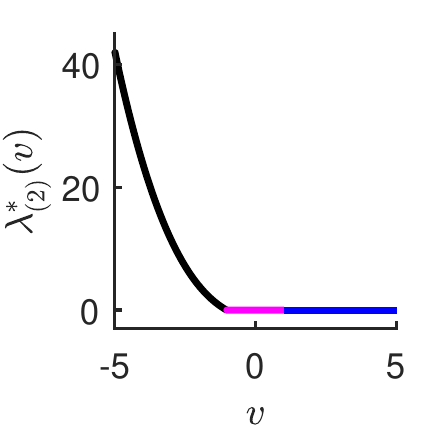}
			\begin{center}
				b) Scenario 2
			\end{center}

			\captionof{figure}{The functions $sol$, $\lambda_{(1)}$ and $\lambda_{(2)}$ are globally $\mathcal{C}^0$ and at least piecewise-$\mathcal{C}^1$. Each ``piece'' is distinguished by a a color: respectively  \textcolor{magenta}{$\bm{\mathcal{A}(v)=\emptyset}$} ,  \textcolor{blue}{$\bm{\mathcal{A}(v)=1}$}, and $\bm{\mathcal{A}(v)=2}$.}
			\label{fig:A___1_SolU_LambdaU_CActive}
		\end{minipage}\\
\end{exbox}

\subsection{Statements B and C: \textit{sol} in the vicinity of the plant's stationary points}
\label{an:A___2_3_Differents_types_de_pts_stationnaires}

If a point $(\bm{u}^{\bullet},\bm{\lambda}_p)$ is a KKT point of the plant, then the derivative of the Lagrangian of the plant $\mathcal{L}_{p}(\bm{u}) = \sum_{i=0}^{n_g} \lambda_{p(i)} g_{p(i)}|_{\bm{u}^{\bullet}},$ satisfies: 
\begin{equation} \label{eq:A_____2132_ercbgfjgfhj}
	\nabla_{\bm{u}} \mathcal{L}_{p}(\bm{u}) = \sum_{i=0}^{n_g} \lambda_{p(i)} \nabla_{\bm{u}} g_{p(i)}|_{\bm{u}^{\bullet}} = 0. 
\end{equation}
In addition, if $(\bm{u}^{\bullet},\bm{\lambda}_p)$ is a KKT point of the plant then it is also a KKT point of the model modified with ISO-D or -I, or any other method satisfying the applicability condition of Theorem~\ref{thm:2___2_ConditionNiveau1}. So: 
\begin{align}
	\bm{sol}(\bm{v}^{\bullet}) = \ &\bm{v}^{\bullet}, & 
	\bm{\lambda}(\bm{v}^{\bullet}) = \bm{\lambda}_p.
\end{align}

Now, let's analyze how  $\bm{sol}$ behaves in the neighborhood of $\bm{v}^{\bullet}$. To do this, one approximates the Lagrangian of the model around this stationary point in the following way:
\begin{equation*}
	\begin{array}{l@{\ }l}
		\mathcal{L}(\bm{u},\bm{v}) & :=   \sum_{i=0}^{n_g} \lambda_{(i)}|_{\bm{v}} \widetilde{g}_{(i)}|_{\bm{u},\bm{v}},\\
		\widetilde{g}_{(i)}(\bm{u},\bm{v}) & := g_{p(i)}|_{\bm{v}} + \nabla_{\bm{v}}g_{p(i)}|_{\bm{v}}^{\rm T} (\bm{u}-\bm{v}) + \frac{1}{2} (\bm{u}-\bm{v})^{\rm T}  \nabla_{\bm{uu}} \widetilde{g}_{(i)}|_{\bm{v},\bm{v}}(\bm{u}-\bm{v}), \\
		\lambda_{(0)}(\bm{v})  & := 1. 
	\end{array}
\end{equation*}
Indeed, since the correction point and the minimum of the updated model are the same, the 0th and 1st order terms of the functions $\widetilde{g}_{(i)}$ are those of the plant evaluated at the correction point while the 2nd order term is the Hessian of the updated model evaluated at $\bm{v}^{\bullet}$. Of course, this approximation is only valid for $\forall \bm{v},\bm{u}\in\mathcal{B}(\bm{u}^{\bullet},r\rightarrow 0)$.

The functions  $\bm{sol}(\bm{v})$ and $\bm{\lambda}(\bm{v})$ can be linked within the same expression in the following way: 
\begin{equation} \label{eq:A___6_Lagrangien}
	\partial_{\bm{u}} \mathcal{L}(\bm{sol}(\bm{v}),\bm{v}) =
	\sum_{i=0}^{n_g} \lambda_{(i)}|_{\bm{v}}\left(
	\nabla_{\bm{v}}g_{p(i)}|_{\bm{v}} + \nabla_{\bm{uu}} \widetilde{g}_{(i)}|_{\bm{v},\bm{v}}(\bm{sol}(\bm{v})-\bm{v})
	\right) = \bm{0}
\end{equation}
By insulating  $\bm{sol}(\bm{v})$ on the left side of the equation one finds:
\begin{align}
	\bm{sol}(\bm{v}) = \ &  \bm{v} - 
	\left( \sum_{i=0}^{n_g}  \lambda_{(i)}|_{\bm{v}}   \nabla_{\bm{uu}}\widetilde{g}_{(i)}|_{\bm{v},\bm{v}}
	\right)^{-1}
	\left(
	\sum_{i=0}^{n_g}  \lambda_{(i)}|_{\bm{v}}   \nabla_{\bm{v}}g_{p(i)}|_{\bm{v}}
	\right). 
\end{align}
\textit{(One notices that one can recognize that $\bm{sol}(\bm{v}^{\bullet}) =  \bm{v}^{\bullet}$ because for $\bm{v}=\bm{v}^{\bullet}$ there is a $\bm{\lambda}_p\in\amsmathbb{R}^{n_g +}$ such that \eqref{eq:A_____2132_ercbgfjgfhj} holds.)} From this expression one finds that: 
\begingroup
 \allowdisplaybreaks
\begin{align}
	 \nabla_{\bm{v}} \bm{sol}(\bm{v})|_{\bm{v}^{\bullet}} = \ & 
	 \bm{I} - \Big(
	 \sum_{i=0}^{n_g}
	 \lambda_{(i)}|_{\bm{v}^{\bullet}}
	 \nabla_{\bm{uu}} 	\widetilde{g}_{(i)}|_{\bm{v}^{\bullet},\bm{v}^{\bullet}}\Big)
	 ^{-1} ... \nonumber \\
	 & 
	  \Big(
	 	\sum_{i=0}^{n_g} 
	 \lambda_{(i)}|_{\bm{v}^{\bullet}} 
	 \nabla_{\bm{vv}}g_{p(i)}|_{\bm{v}^{\bullet}} 
	 +
	 \nabla_{\bm{v}}g_{p(i)}|_{\bm{v}^{\bullet}}
	 \nabla_{\bm{v}} \lambda_{(i)}|_{\bm{v}^{\bullet}}^{\rm T}
	 \Big), \nonumber \\
	 = \ & 
	 \bm{I} -
	 \left(
	 \nabla_{\bm{uu}}\mathcal{L}|_{\bm{v}^{\bullet},\bm{v}^{\bullet}}
	 \right)^{-1}
	 \nabla_{\bm{vv}}\mathcal{L}_p|_{\bm{v}^{\bullet}} 
	 -
	 \left(
	 \nabla_{\bm{uu}}\mathcal{L}|_{\bm{v}^{\bullet},\bm{v}^{\bullet}}
	 \right)^{-1}
	 \sum_{i=0}^{n_g}
	 \nabla_{\bm{v}}g_{p(i)}|_{\bm{v}^{\bullet}}
	 \nabla_{\bm{v}} \lambda_{(i)}|_{\bm{v}^{\bullet}}^{\rm T} \label{eq:A_merde}
\end{align}
\endgroup
Let's focus on the term:
\begin{equation}\label{eq:A___8_rbhuceeeeegc}
	\left(
	\nabla_{\bm{uu}}\mathcal{L}|_{\bm{v}^{\bullet},\bm{v}^{\bullet}}
	\right)^{-1} \sum_{i=0}^{n_g}  \nabla_{\bm{v}}g_{p(i)}|_{\bm{v}^{\bullet}}
	\nabla_{\bm{v}} \lambda_{(i)}|_{\bm{v}^{\bullet}}^{\rm T}.
\end{equation}
which can be simplified in the following three ways:  
\begin{itemize}
	\item 
	\textbf{Concerning the cost function ($i=0$):} Starting from the definition of $\lambda_{(0)} $:
	\begin{align}
		& \forall \bm{v}\in\amsmathbb{R}^{n_u}, \ \lambda_{(0)}|_{\bm{v}}:= \  1,  &
		\Rightarrow \ \ & 
		\nabla_{\bm{v}} \lambda_{(0)}|_{\bm{v}^{\bullet}}
		= \bm{0}, \nonumber \\
		&&
		\Rightarrow \ \ &
		\nabla_{\bm{v}}g_{p(0)}|_{\bm{v}^{\bullet}}
		\nabla_{\bm{v}} \lambda_{(0)}|_{\bm{v}^{\bullet}}^{\rm T}
		=   \bm{0}.  \label{eq:A___9_kebcgwjnhv}
	\end{align}
	\item 
	\textbf{Concerning the \textit{inactive} constraints  ($i\not\in\mathcal{A}(\bm{v}^{\bullet})$):} 
	For all constraints that are not active at $\bm{v}^{\bullet}$, i.e.  $\forall i\not\in\mathcal{A}(\bm{v}^{\bullet})$, it is clear that:
	\begin{align}
		g_{p(i)}|_{\bm{v}^{\bullet}} < 0, \  
		\& \  
		g_{p(i)} \in \mathcal{C}^0  \quad  \Rightarrow \quad & 
		\forall \bm{v} \in \mathcal{B}(\bm{v}^{\bullet}, r\rightarrow 0): \  g_{p(i)}|_{\bm{v}} < 0, \nonumber \\
		\Rightarrow \quad & 
		\forall \bm{v} \in \mathcal{B}(\bm{v}^{\bullet}, r\rightarrow 0): \  \lambda_{(i)}|_{\bm{v}} = 0, \nonumber \\
		\Rightarrow \quad &  
		\forall i\not\in\mathcal{A}(\bm{v}^{\bullet}): \ 
		\nabla_{\bm{v}} \lambda_{(i)}|_{\bm{v}^{\bullet}} = \bm{0}, \nonumber \\
		\Rightarrow  \quad &  
		\sum_{i\not\in\mathcal{A}(\bm{v}^{\bullet})}
		\nabla_{\bm{v}}g_{p(i)}|_{\bm{v}^{\bullet}}
		\nabla_{\bm{v}} \lambda_{(i)}|_{\bm{v}^{\bullet}}^{\rm T}
		= \bm{0}.  \label{eq:A___10_sdafukygj}
	\end{align}
	\item 
	\textbf{Concerning the \textit{weakly active} constraints according to the direction $\delta\bm{v}$ ($i\not\in \mathcal{A}(\bm{v}^{\bullet},\delta\bm{v})$):} One says that a constraint $i$ is weakly active in a direction  $\delta\bm{v}$ if: 
	\begin{align*}
		g_{p(i)}|_{\bm{v}^{\bullet}} = \ & 0, & 
		\underset{\substack{h\rightarrow 0 \\ h > 0}}{\operatorname{lim}} \  g_{p(i)}|_{\bm{v}^{\bullet} + h\delta\bm{v}} < \ & 0.
	\end{align*}
	If this is the case, as $\forall h$ close to $0$ and greater than  $0$ one has $g_{p(i)}|_{\bm{v}^{\bullet} + h\delta\bm{v}} <  0$, it is then clear that $\lambda_{(i)}|_{\bm{v}^{\bullet} + h\delta\bm{v}} = 0$. And since $\lambda_{(i)}$ is a  $\mathcal{C}^0$ function of $\bm{v}$: 
	\begin{align}
		\lambda_{(i)}|_{\bm{v}^{\bullet}}= \ & 0, &
		\nabla_{\bm{v}}\lambda_{(i)}|_{\bm{v}^{\bullet}} \delta\bm{v} =\  & 0, & \text{if $i \not\in \mathcal{A}(\bm{v}^{\bullet},\delta\bm{v})$}.
	\end{align}
\end{itemize}

From these three observations one can reduce \eqref{eq:A___8_rbhuceeeeegc} to: 
\begin{equation} \label{eq:A___14_cjherghj}
	\text{\eqref{eq:A___8_rbhuceeeeegc}}
	= 
	\left(
	\nabla_{\bm{uu}}\mathcal{L}|_{\bm{v}^{\bullet},\bm{v}^{\bullet}}
	\right)^{-1}
	\sum_{i\in\mathcal{A}(\bm{v}^{\bullet},\delta\bm{v})} \nabla_{\bm{v}}g_{p(i)}|_{\bm{v}^{\bullet}}
	\nabla_{\bm{v}} \lambda_{(i)}|_{\bm{v}^{\bullet}}^{\rm T}.
\end{equation}
So only the active constraints at $\bm{v}^{\bullet}$ remain active when one moves in the direction of $\delta\bm{v}$ (which is linked to $\nabla_{\bm{v}}\bm{sol}$). By definition the constraints are such that:   
\begin{align} \label{eq:A_kucherg}
	\bm{g}^{a}(\bm{sol}(\bm{v})) = \ & 0, &
	\nabla_{\bm{v}}\bm{g}^{a}(\bm{sol}(\bm{v})) = \ & 0,
\end{align}
where
\begin{align*}
	\bm{g}^{a}_{p} := \ & [g_{p(i)}]_{i\in\mathcal{A}(\bm{v}^{\bullet},\delta\bm{v})}.
\end{align*}
One approximates them around $\bm{v}^{\bullet}$ with a Taylor series:
\begin{align}
	\bm{g}^{a}(\bm{sol}(\bm{v})) = \ & 
	\bm{g}^a_{p}|_{\bm{v}^{\bullet}} + 
	\nabla_{\bm{v}}\bm{g}^a_{p}|_{\bm{v}^{\bullet}} 
	(\bm{sol}(\bm{v})-\bm{v})|_{\bm{v}^{\bullet}}, \nonumber\\
	\Rightarrow \ \	\nabla_{\bm{v}}\bm{g}^{a}(\bm{sol}(\bm{v})) = \ &
	\nabla_{\bm{v}} \bm{g}^a_{p}|_{\bm{v}^{\bullet}} +
	\nabla_{\bm{vv}}\bm{g}^a_{p}|_{\bm{v}^{\bullet}} 
	(\bm{sol}(\bm{v})-\bm{v})
	+ 
	\nabla_{\bm{v}}\bm{g}^a_{p}|_{\bm{v}^{\bullet}} 
	(\nabla_{\bm{v}}\bm{sol}(\bm{v})-\bm{I})^{\rm T},  \nonumber \\
	= \ &
	\nabla_{\bm{v}} \bm{g}^a_{p}|_{\bm{v}^{\bullet}} + 	\nabla_{\bm{v}}\bm{g}^a_{p}|_{\bm{v}^{\bullet}} 
	(\nabla_{\bm{v}}\bm{sol}(\bm{v})-\bm{I})^{\rm T}, \label{eq:A_kgfvejrc}
\end{align}
One combines \eqref{eq:A_kucherg} and \eqref{eq:A_kgfvejrc}:
\begin{align*}
	\nabla_{\bm{v}} \bm{g}^a_{p}|_{\bm{v}^{\bullet}}	(\bm{I}-\nabla_{\bm{v}}\bm{sol}(\bm{v}))^{\rm T} =  \ &
	\nabla_{\bm{v}}\bm{g}^a_{p}|_{\bm{v}^{\bullet}}, 
\end{align*}
One multiplies on each side by $\nabla_{\bm{v}} \bm{g}^a_{p}|_{\bm{v}^{\bullet}}^{\rm T}
\left(
\nabla_{\bm{v}} \bm{g}^a_{p}|_{\bm{v}^{\bullet}}
\nabla_{\bm{v}} \bm{g}^a_{p}|_{\bm{v}^{\bullet}}^{\rm T}
\right)^{-1}$:
\begin{align*}
	\bm{M}
	(\bm{I}-\nabla_{\bm{v}}\bm{sol}(\bm{v}))^{\rm T}
	=  \ & \bm{M}, 
\end{align*}
where 
\begin{align*}
	\bm{M} := \ & 
	\nabla_{\bm{v}} \bm{g}^a_{p}|_{\bm{v}^{\bullet}}^{\rm T}
	\left(
	\nabla_{\bm{v}} \bm{g}^a_{p}|_{\bm{v}^{\bullet}}
	\nabla_{\bm{v}} \bm{g}^a_{p}|_{\bm{v}^{\bullet}}^{\rm T}
	\right)^{-1}
	\nabla_{\bm{v}} \bm{g}^a_{p}|_{\bm{v}^{\bullet}}.
\end{align*}
Then, one carries out a series of "simple" manipulations:	
\begin{align}
	\bm{M}\nabla_{\bm{v}}\bm{sol}(\bm{v})^{\rm T}
	=  \ & \bm{0}, \nonumber  \\
	\nabla_{\bm{v}}\bm{sol}(\bm{v})\bm{M}
	=  \ & \bm{0}, \nonumber \\
	\nabla_{\bm{v}}\bm{sol}(\bm{v})\bm{M} \delta\bm{v}
	=  \ & \bm{0}, \nonumber \\
	\nabla_{\bm{v}}\bm{sol}(\bm{v})(\bm{M} \delta\bm{v} - \delta\bm{v}  + \bm{I}\delta\bm{v})
	=  \ & \bm{0}, \nonumber \\
	\nabla_{\bm{v}}\bm{sol}(\bm{v})\delta\bm{v} 
	=  \ & \nabla_{\bm{v}}\bm{sol}(\bm{v})(\bm{I}-\bm{M}) \delta\bm{v}. \label{eq:A___flberbwhfer}
\end{align}
Finally, one combines \eqref{eq:A_merde} and \eqref{eq:A___flberbwhfer}:
\begin{align*}
	 & \delta\bm{v}^{\rm T} \nabla_{\bm{v}} \bm{sol}(\bm{v})|_{\bm{v}^{\bullet}}  \delta\bm{v} =
	\delta\bm{v}^{\rm T}  \left(\bm{I} -
	\left(
	\nabla_{\bm{uu}}\mathcal{L}|_{\bm{v}^{\bullet},\bm{v}^{\bullet}}
	\right)^{-1}
	\nabla_{\bm{vv}}\mathcal{L}_p|_{\bm{v}^{\bullet}} \right)
	\left(
	\bm{I} - \bm{M}
	\right) \delta\bm{v}
	- ... 	\\ 
	& \qquad \qquad 
	\delta\bm{v}^{\rm T} 
	\big(
	\left(
	\nabla_{\bm{uu}}\mathcal{L}|_{\bm{v}^{\bullet},\bm{v}^{\bullet}}
	\right)^{-1}
	\sum_{i\in\mathcal{A}(\bm{v}^{\bullet},\delta\bm{v})}
	\nabla_{\bm{v}} \lambda_{(i)}|_{\bm{v}^{\bullet}}^{\rm T}
	\underbrace{
	\nabla_{\bm{v}}g_{p(i)}|_{\bm{v}^{\bullet}}
	\big)
	\big(
	\bm{I} - \bm{M}
	\big) \delta\bm{v}
	}_{=0},
\end{align*}
where indeed: 
\begin{align*}
	\nabla_{\bm{v}} \bm{g}^a_{p}|_{\bm{v}^{\bullet}}(\bm{I}-\bm{M}) = \ &  \nabla_{\bm{v}} \bm{g}^a_{p}|_{\bm{v}^{\bullet}} - \nabla_{\bm{v}} \bm{g}^a_{p}|_{\bm{v}^{\bullet}}\nabla_{\bm{v}} \bm{g}^a_{p}|_{\bm{v}^{\bullet}}^{\rm T}
	\left(
	\nabla_{\bm{v}} \bm{g}^a_{p}|_{\bm{v}^{\bullet}}
	\nabla_{\bm{v}} \bm{g}^a_{p}|_{\bm{v}^{\bullet}}^{\rm T}
	\right)^{-1}
	\nabla_{\bm{v}} \bm{g}^a_{p}|_{\bm{v}^{\bullet}}, \\
	= \ & \nabla_{\bm{v}} \bm{g}^a_{p}|_{\bm{v}^{\bullet}} - \nabla_{\bm{v}} \bm{g}^a_{p}|_{\bm{v}^{\bullet}}, \\
	= \ & \bm{0}.
\end{align*}
This only leaves:
\begin{align*}
	\delta\bm{v}^{\rm T} 
	\nabla_{\bm{v}} \bm{sol}(\bm{v})|_{\bm{v}^{\bullet}}  \delta\bm{v} = \ & 
	\delta\bm{v}^{\rm T} 
	\left(\bm{I} -
	\left(
	\nabla_{\bm{uu}}\mathcal{L}|_{\bm{v}^{\bullet},\bm{v}^{\bullet}}
	\right)^{-1}
	\nabla_{\bm{vv}}\mathcal{L}_p|_{\bm{v}^{\bullet}} \right)
	\left(
	\bm{I} - \bm{M}
	\right) \delta\bm{v}, \\
	= \ & \delta\bm{v}^{\rm T}\bm{A} \delta\bm{v} - \delta\bm{v}^{\rm T}\bm{B} \delta\bm{v}.
\end{align*}
where
\begin{align}
	\bm{A} := \ & \bm{I}-\bm{M}, &
	\bm{B} := \ & 
	\left(\left(
	\nabla_{\bm{uu}}\mathcal{L}|_{\bm{v}^{\bullet},\bm{v}^{\bullet}}
	\right)^{-1}
	\nabla_{\bm{vv}}\mathcal{L}_p|_{\bm{v}^{\bullet}} \right)
	\left(
	\bm{I} - \bm{M}
	\right). 
\end{align}
Let us now analyze the terms $\delta\bm{v}^{\rm T}\bm{A}\delta\bm{v}$ and  $\delta\bm{v}^{\rm T}\bm{B}\delta\bm{v}$. Let's start with the first one: 
\begin{equation}
	\frac{\delta\bm{v}^{\rm T} \bm{A}\delta\bm{v}}{\|\delta\bm{v}\|^2} 
	=  
	1 - \frac{\delta\bm{v}^{\rm T}\bm{M}\delta\bm{v}  }{\|\delta\bm{v}\|^2},  \qquad \qquad 
	\frac{\delta\bm{v}^{\rm T} \bm{A}\delta\bm{v}}{\|\delta\bm{v}\|^2}  
	\left\{
	\begin{array}{ll}
		= 1, & \text{ if $\bm{M}\delta\bm{v} = 0$}, \\
		< 1, & \text{ if $\bm{M}\delta\bm{v} \neq 0$}. \\
	\end{array}
	\right. \label{eq:A___17_dAd}
\end{equation}	
As the model is systematically convex at the correction point, see \eqref{eq:2_29_uywgkfyurwkc}:
\begin{equation}
	\nabla_{\bm{uu}} \mathcal{L}|_{\bm{v}^{\bullet},\bm{v}^{\bullet}} > 0,
\end{equation}
\begin{align}
		& \bullet\text{ If $\bm{v}^{\bullet}$ is a minimum of the plant: } 
		&\Rightarrow &  
		\nabla_{\bm{vv}} \mathcal{L}_p|_{\bm{v}^{\bullet}} > 0 , \nonumber
		\\
		&& \Rightarrow &  
		(\nabla_{\bm{uu}} \mathcal{L}|_{\bm{v}^{\bullet},\bm{v}^{\bullet}})^{-1}
		\nabla_{\bm{vv}} \mathcal{L}_p|_{\bm{v}^{\bullet}}
		>0, \qquad \  \label{eq:A___19_ertnjh}\\
		&& \Rightarrow&   
		(\nabla_{\bm{uu}} \mathcal{L}|_{\bm{v}^{\bullet},\bm{v}^{\bullet}})^{-1}
		\nabla_{\bm{vv}} \mathcal{L}_p|_{\bm{v}^{\bullet}}(\bm{I}-\bm{M})
		\geq 0, \nonumber\\
		&& \Rightarrow&   
		\delta\bm{v}^{\rm T} \bm{B} \delta\bm{v}
		\left\{
		\begin{array}{ll}
			= 0, & \text{if $(\bm{I}-\bm{M})\delta\bm{v} = 0$}, \\
			> 0, & \text{if $(\bm{I}-\bm{M})\delta\bm{v}\neq 0$}.
		\end{array}
		\right. \nonumber \\
		& \bullet\text{ If $\bm{v}^{\bullet}$ is a maximum of the plant: } 
		& \Rightarrow& 
		\delta\bm{v}^{\rm T} \bm{B} \delta\bm{v}
		\left\{
		\begin{array}{ll}
			= 0, & \text{if $(\bm{I}-\bm{M})\delta\bm{v} = 0$}, \\
			< 0, & \text{if $(\bm{I}-\bm{M})\delta\bm{v}\neq 0$}.
		\end{array}
		\right. \nonumber \\
		&\bullet\text{ If $\bm{v}^{\bullet}$ is a saddle-point of the plant: } 
		&\Rightarrow & 
		\delta\bm{v}^{\rm T} \bm{B} \delta\bm{v}
		\left\{
		\begin{array}{ll}
			= 0, & \text{if $(\bm{I}-\bm{M})\delta\bm{v} = 0$}, \\
			\gtrless 0, & \text{if $(\bm{I}-\bm{M})\delta\bm{v}\neq 0$}.
		\end{array}
		\right. \nonumber
\end{align}
And by combining these results with   \eqref{eq:A___17_dAd} one finds that $\forall \delta\bm{v}\in\amsmathbb{R}^{n_u}\backslash \bm{0}$:
\begin{align}
	& \bullet\text{ If $\bm{v}^{\bullet}$ is a minimum of the plant: } 
	& & \Rightarrow 
	\frac{\delta\bm{v}^{\rm T}\bm{A} \delta\bm{v} - \delta\bm{v}^{\rm T}\bm{B} \delta\bm{v}}{\|\delta\bm{v}\|^2} < 1 , \nonumber
	\\
	& & & \Rightarrow 
	\frac{ \delta\bm{v}^{\rm T} 
	\nabla_{\bm{v}} \bm{sol}(\bm{v})|_{\bm{v}^{\bullet}}  \delta\bm{v}}{\|\delta\bm{v}\|^2} < 1 . \nonumber \\
	& & & \Rightarrow 
	\nabla^S(\bm{v}^{\bullet},\delta\bm{v}) < 1 . \phantom{\underbrace{a}_{a}} \nonumber \\
	& \bullet\text{ If $\bm{v}^{\bullet}$ is a maximum of the plant:} 
	& & \Rightarrow
	\nabla^S(\bm{v}^{\bullet},\delta\bm{v}) > 1, \nonumber \\
	&\bullet\text{ If  $\bm{v}^{\bullet}$ is a saddle-point of the plant: } 
	& &\Rightarrow
	\nabla^S(\bm{v}^{\bullet},\delta\bm{v})  \gtrless 1. \nonumber
\end{align}
This concludes the proof of statement B.

Finally, one can reuse the results of the previous analysis to prove statement C: 
\begin{align}
 & & \nabla^S(\bm{v},\delta\bm{v}) < \ & -1,  \label{eq:A_erhbherbhgebrhv} \\
 & & \frac{\delta\bm{v}^{\rm T} \bm{A} \delta\bm{v} - \delta\bm{v}^{\rm T}  \bm{B} \delta\bm{v} }{\|\delta\bm{v}\|^2} < \ & -1, \nonumber \\
 & & 1 - \frac{\delta\bm{v}^{\rm T}  \bm{B} \delta\bm{v}}{\|\delta\bm{v}\|^2} < \ & -1, \nonumber \\
 & & 1 - \frac{\delta\bm{v}^{\rm T}  
 	\left((\nabla_{\bm{uu}}\mathcal{L}|_{\bm{v}^{\bullet},\bm{v}^{\bullet}})^{-1}
 	\nabla_{\bm{uu}}\mathcal{L}_p|_{\bm{v}^{\bullet}}\right)
 	\left(
 		\bm{I}-\bm{M}
 	\right) \delta\bm{v}}{\|\delta\bm{v}\|^2} < \ & -1, \nonumber \\ 
& & 1 - \frac{\delta\bm{v}^{\rm T}  
	\left((\nabla_{\bm{uu}}\mathcal{L}|_{\bm{v}^{\bullet},\bm{v}^{\bullet}})^{-1}
	\nabla_{\bm{uu}}\mathcal{L}_p|_{\bm{v}^{\bullet}}\right) \delta\bm{v}}{\|\delta\bm{v}\|^2} < \ & -1, \nonumber \\ 
& & \frac{\delta\bm{v}^{\rm T}  
\left((\nabla_{\bm{uu}}\mathcal{L}|_{\bm{v}^{\bullet},\bm{v}^{\bullet}})^{-1}
\nabla_{\bm{uu}}\mathcal{L}_p|_{\bm{v}^{\bullet}}\right) \delta\bm{v}}{\|\delta\bm{v}\|^2} > \ & 2. \label{eq:A_erhbherbhgebrhvb}
\end{align}
Where one mainly uses the fact that as $\delta\bm{v}$ is a vector orthogonal to the set of active constraints of the plant at $\bm{v}^{\bullet}$:
\begin{equation*} 
	\bm{M}\delta\bm{v} = \nabla_{\bm{v}} \bm{g}^a_{p}|_{\bm{v}^{\bullet}}^{\rm T}
	\left(
	\nabla_{\bm{v}} \bm{g}^a_{p}|_{\bm{v}^{\bullet}}
	\nabla_{\bm{v}} \bm{g}^a_{p}|_{\bm{v}^{\bullet}}^{\rm T}
	\right)^{-1}
	\underbrace{
		\nabla_{\bm{v}} \bm{g}^a_{p}|_{\bm{v}^{\bullet}}\delta\bm{v}
	}_{=\bm{0}} = \bm{0}.
\end{equation*} 
The equivalence \eqref{eq:A_erhbherbhgebrhv} $\Leftrightarrow$  \eqref{eq:A_erhbherbhgebrhvb} demonstrates   \eqref{eq:2___32_Condition_NablaS_pluspetitque_m1}. Which concludes the proof of the Theorem~\ref{thm:2___3_Proprietes_Sol}.


\chapter{The Tennessee Eastman challenge process }
\label{app:B_TE}

\begin{figure}
	\centering	
	\includegraphics[width=1\linewidth,trim={0mm  430mm  10mm 25mm},clip]{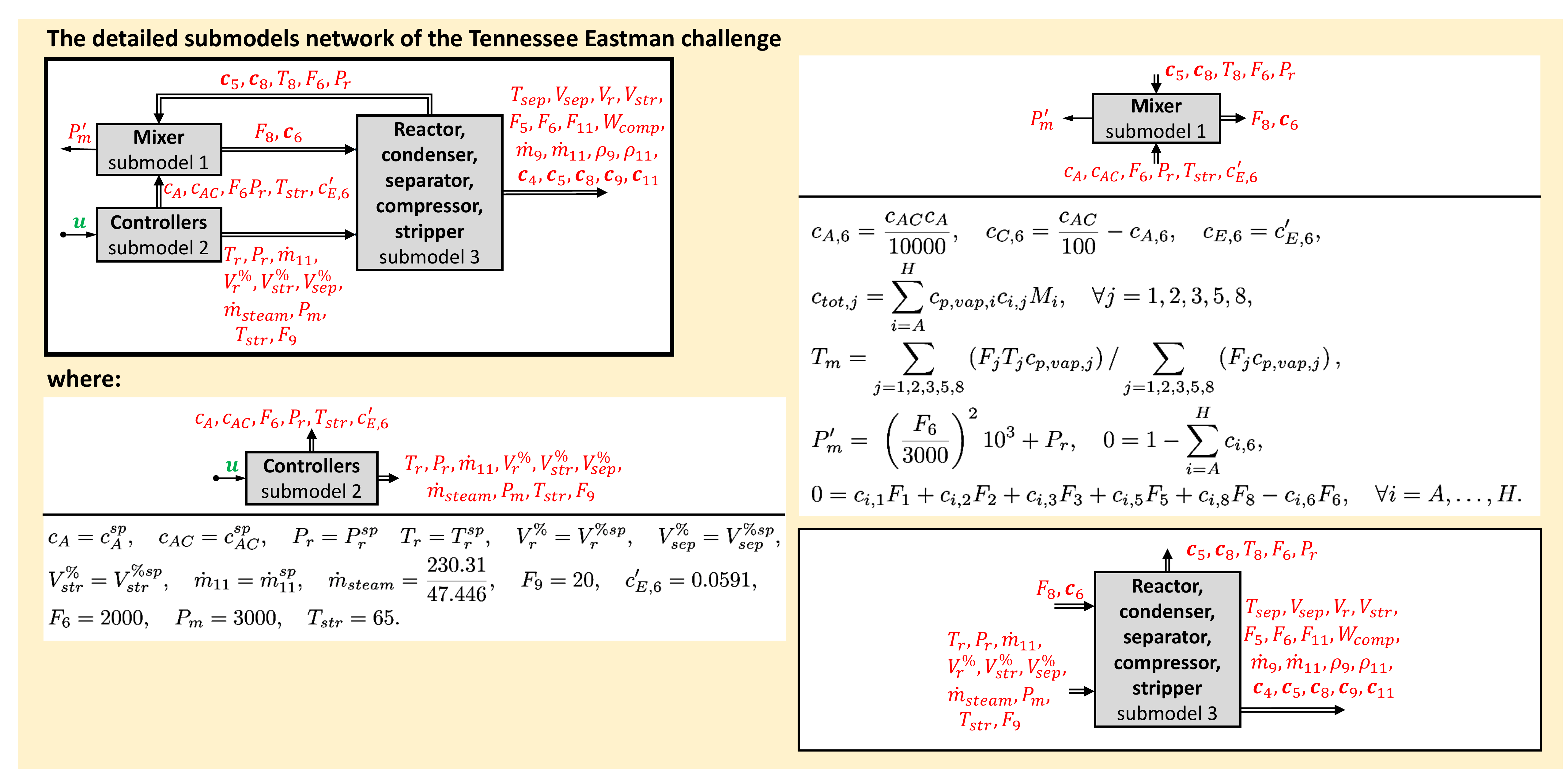}
	\caption{Description of the submodels network obtained for the TE case study with the constitutive equations of sub-models 1 and 2. } 
	\label{fig:B_TE_description_1} 	
\end{figure}

\begin{figure}
	\centering	
	\includegraphics[width=1\linewidth,trim={0mm  45mm  10mm 290mm},clip]{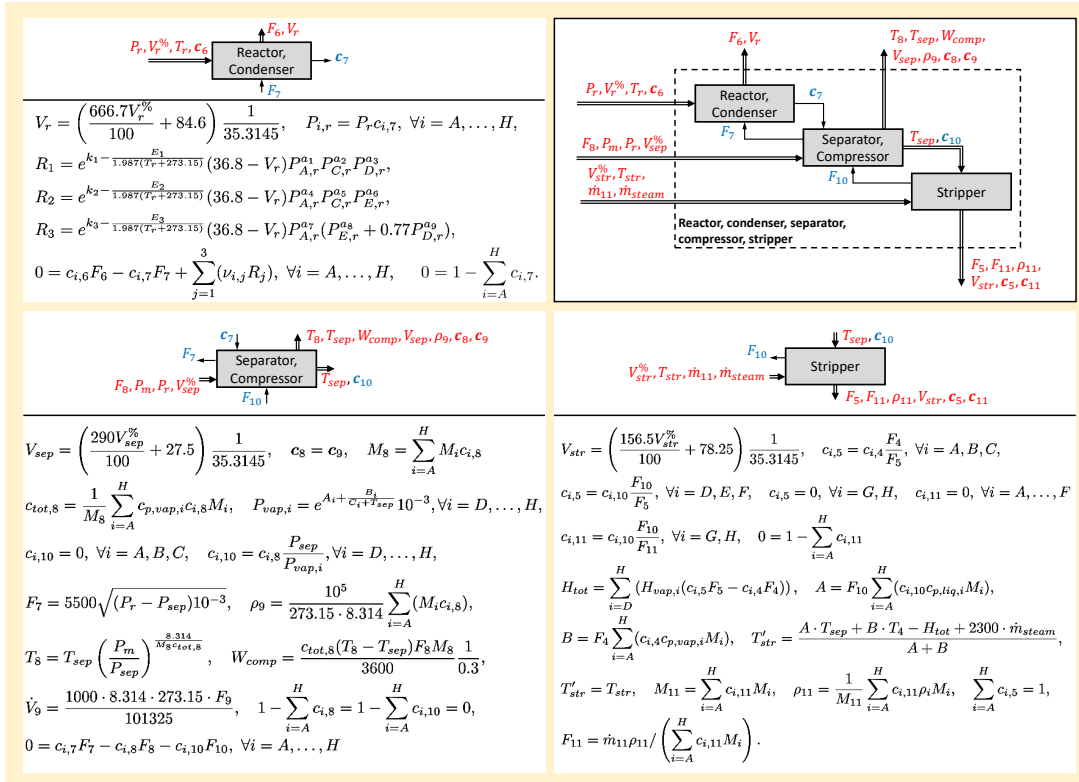}
	\caption{Equations of submodel 3 for the TE case study.} 
	\label{fig:B_TE_description_2} 	
\end{figure}

The TE model is described in Figures~\ref{fig:B_TE_description_1} and \ref{fig:B_TE_description_2} and the notation is explained in Table~\ref{table:B_1_TE_SimplifiedModelNomenclature}. The four sub-costs are: 
\begin{align}
	& \phi_{(1)}=   0.0536 \ W_{comp}, & & \phi_{(2)}=   0.0318\  \dot{m}_{steam}, \nonumber \\
	& \phi_{(3)} =  \frac{\dot{V}_9\rho_9}{M_9} \sum_{\substack{i=A \\ i\neq B}}^H C_ic_{i,9}, & &
	\phi_{(4)} =  \frac{\dot{V}_{11}\rho_{11}}{M_{11}} \sum_{i=D}^F C_ic_{i,11}, \nonumber
\end{align}
where $C_i$ are the costs of the components $i=A,\hdots,H$ listed in Table~\ref{tab:B2_CostsComponents}.  The inequality constraints are: 
\begin{align}
	\bm{g} = \ \Big[ & g_{(1)} := T_r - 150, 
	\quad g_{(2)} := P_r - 2895,
	\quad g_{(3)} := V_r - 21.3, 
	\quad g_{(4)} := 11.8 - V_r, \nonumber\\
	&  g_{(5)} := V_{sep} - 9,
	\quad g_{(6)} := 3.3 - V_{sep}, 
	\quad g_{(7)} := V_{str} - 6.6, 
	\quad g_{(8)} := 3.5 - V_{str},\nonumber \\
	&  g_{(9)} := \frac{c_{G,11}M_G}{c_{H,11}M_H}-1, 
	\quad g_{(10)} := 1 - \frac{c_{G,11}M_G}{c_{H,11}M_H}, \nonumber \\
	& g_{(11)} := \frac{\dot{V}_{11}\rho_{11} }{M_{11}}(M_Gc_{G,11} + M_Hc_{H,11}) - 14076,\nonumber \\
	&  g_{(12)} := 14076 - \frac{\dot{V}_{11}\rho_{11} }{M_{11}}(M_Gc_{G,11} + M_Hc_{H,11}) \Big]^{\rm T} \leq \bm{0},  \nonumber
\end{align}
where $\{g_{(9)}, g_{(10)}\}$ and $\{g_{(11)}, g_{(12)}\}$ are quality and production constraints, respectively.

\begin{Remark} (An alternative to affine corrections) \label{remark:NL_modifiers}
	For the implementation of IMA-A and -B, the way some of the measured variables have been corrected has been slightly changed. Indeed the (expected) affine-in-input correction of $\bm{c}_5$, $\bm{c}_6$, $\bm{c}_8$, $\bm{c_9}$, and $\bm{c}_{11}$ for IMA-A was implemented instead as follows, $ \forall i =A,\hdots,H$, and $j= 5, \, 6,\,8,\,9,\,11$:
	\begin{equation}
	c_{i,j,m}(\bm{u}) = \frac{c_{i,j}(\bm{u}) + \bm{\varepsilon}_k^{c_{i,j}} + \bm{\lambda}_k^{c_{i,j}}(\bm{u}-\bm{u}_k)}{\sum_{\ell=A}^H c_{\ell,j}(\bm{u}) + \bm{\varepsilon}_k^{c_{\ell,j}} + \bm{\lambda}_k^{c_{\ell,j}}(\bm{u}-\bm{u}_k)}.
	\end{equation}
	Doing so allows maintaining $c_{i,j,m}\in[0,1]$ and avoids numerical problems due to unrealistic linear corrections of the molar fractions (i.e. outside $[0,\ 1]$), which led  simulations to fail. Still, it is easy to verify that doing so preserved the affine correction of $c_{i,j}$ at $\bm{u}_k$:
	\begin{align*}
	c_{i,j,m}\big|_{\bm{u}_k} = \ & \frac{c_{i,j,p}\big|_{\bm{u}_k}}{\sum_{\ell=A}^H c_{\ell,j,p}\big|_{\bm{u}_k}}, \\
	\nabla_u c_{i,j,m}\big|_{\bm{u}_k} = \ & \Big(
	\nabla_u c_{i,j,p}\big|_{\bm{u}_k} \sum_{\ell=A}^H c_{\ell,j,p}\big|_{\bm{u}_k} +  
	c_{i,j,p}\big|_{\bm{u}_k} \sum_{\ell=A}^H \nabla_u c_{\ell,j,p}\big|_{\bm{u}_k}
	\Big) \Big(\sum_{\ell=A}^H c_{\ell,j,p}\big|_{\bm{u}_k}\Big)^{ -2} .
	\end{align*}
	Since $\sum_{\ell=A}^H c_{\ell,j,p}\big|_{\bm{u}_k} = 1$ and  $\sum_{\ell=A}^H \nabla_u c_{\ell,j,p}\big|_{\bm{u}_k}=0$, these equations can be rewritten:
	\begin{align*}
	c_{i,j,m}\big|_{\bm{u}_k} = \ & \frac{c_{i,j,p}\big|_{\bm{u}_k}}{1} = c_{i,j,p}, \\
	\nabla_u c_{i,j,m}\big|_{\bm{u}_k} = \ & \frac{
		\nabla_u c_{i,j,p}\big|_{\bm{u}_k} + 0}{1^{2}} =  \nabla_u c_{i,j,p}\big|_{\bm{u}_k},
	\end{align*}
	which proves that the aim of the input-affine correction is preserved. To the best of the authors knowledge, this is the first time that a nonlinear correction is performed (for outputs and/or cost and constraints) in the context of MA or MAy, which still allows the reconciliation of zeroth and first order behavior of the corrected and plant quantity.
\end{Remark}


\begin{figure}[b]
	\centering
	\begin{tabular}{p{1.1cm} p{3cm} p{9cm}}
		\toprule[1.5pt]
		\textbf{Symbols} & \textbf{\ \ \ Unit} & \textbf{Description} \\
		\midrule
		$c^{sp}_A$        & [\% . kmol. A / kmol. (A+C)] & ratio between the molar concentration of A and A+C in the mixer\\
		$c^{sp}_{AC}$     & [\% . kmol .(A+C) / kmol]    &  molar fraction of A and C in the mixer  \\
		$c_{i,j}$         & [kmol . \text{$i$} / kmol]   & molar fraction of  $i$ in  $j$ \\
		$\bm{c}_{j}$      & [kmol . \text{$i$} / kmol]   & vector of the molar concentrations $[c_{A,j},\hdots,c_{H,j}]^{\sf{T}}$ \\
		$c_{p,liq,i}$     & [kJ/(kg.${}^oC$)]            & specific heat capacity of $i$ in liquid phase \\
		$c_{p,vap,i}$     & [kJ/(kg.${}^oC$)]            & specific heat capacity of $i$ in vapor phase \\
		$F_j$             & [kmol / h]                   & molar flow rate of $j$ \\
		$H_{vap,i}$       & [kJ / kg]                    & enthalpy of vaporisation of $i$ \\
		$\dot{m}_j$       & [kg/ h]                      & mass flow rate of stream $j$ \\
		$M_i$             & [kg/kmol . \text{$i$}]       & molecular weight of component or stream $i$ \\
		$P_{\ell}$        & [MPa]                        & total pressure in  $\ell$ \\
		$P_{i,\ell}$      & [MPa]                        & partial pressure of  $i$ in $\ell$ \\
		$P_{vap,i}$       & [MPa]                        & vapor pressure of $i$ in the separator \\
		$R_i$             & [kmol . \text{$i$} / h]      & reaction conversion rate for reaction $i$ \\
		$T_{\ell}$        & [${}^oC$]                    & temperature in $\ell$  \\
		$T_{j}$           & [${}^oC$]                    & temperature in $j$  \\
		$T_{WO,\ell}$     & [${}^oC$]                    & temperature of the water living $\ell$  \\
		$v_j$             & [\%]                         & valve position in $j$ \\
		$V_{\ell}$        & [$m^3$]                      & liquid volume in $\ell$ \\
		$\dot{V}_j$       & [$m^3/h$]                    & volume flow rate in $j$ \\
		$V_{\ell}^{\%}$   & [\%]                         & liquid level in $\ell$ in percent of the maximal level \\
		$W_{comp}$        & [kW]                         & power consumption of the compressor \\
		$\rho_{i}$        & [kg/$m^3$]                   & liquid density of component or stream $i$ at 100${}^oC$ \\
		$\omega_{r}$      & [\%]                         & agitator speed in the reactor \\
		\bottomrule 
	\end{tabular}
	\captionof{table}{Notations used for the Tennessee Eastman challenge process. The values of $c_{p,liq,i},$ $c_{p,vap,i},$ $H_{vap,i},$ $M_i,$ and $\rho_{i}$ can be found in \cite{Downs:1993}. The subscript $i$ is used to denote components (A,...,H), while $j$ refers to streams ($j = 1,\dots,11$) and $\ell$ to units, e.g., reactor (r), separator (sep) or stripper (str).}
	\label{table:B_1_TE_SimplifiedModelNomenclature}
\end{figure}

\begin{figure}[b]
	\centering
	\begin{tabular}{ccccccc}
		\toprule[1.5pt]
		$C_A$    & $C_C$    & $C_D$ & $C_E$    & $C_F$    & $C_G$    & $C_H$ \\
		\midrule
		$2.206$ & $6.177$ & $22.06$ & $14.56$  & $17.89$ & $30.44$ & $22.94$\\
		\bottomrule[1.5pt]
	\end{tabular}
	\captionof{table}{Cost associates to each components in \$/kmol.}
	\label{tab:B2_CostsComponents}
\end{figure}

\chapter{Technical complements to Chapter 6}

\section{Regression of a Gaussian process (SE-ARD-WN) based on points and directional derivatives} 
\label{App:C_PartialDeriv}

Consider the regression problem consisting of identifying the function  $f_p\sim \mathcal{GP}\big(f,$ $ \k_{\f\f}(\bm{x}_i,\bm{x}_j|\sigma_f^2,\bm{\L})\big)$ linking inputs $\widehat{\bm{x}}$ to measured outputs  $\widehat{y}_{p}$ subject to measurement uncertainties $\epsilon_{xy} \sim \mathcal{N}(0,\sigma_{xy}^{2})$. Thus, each experiment result  $\{\widehat{\bm{x}}_{\{i\}}, \widehat{y}_{p\{i\}}, \epsilon_{xy\{i\}}\}$ are linked by the following relationship:
\begin{equation}  \label{eq:C_1_PB_regression}
	\widehat{y}_{p\{i\}} =  f_p(\widehat{\bm{x}}_{\{i\}}) + \epsilon_{xy\{i\}}, \qquad \epsilon_{xy\{i\}} \sim \mathcal{N}(0,\sigma_{xy\{i\}}^{2}).
\end{equation}
In this section one:
\vspace{-\topsep}
\begin{itemize}[noitemsep] 
	\item[] \textbf{(\S~\ref{App:DirectionalDerivatives})} expresses the functions $\nabla_{\bm{x}_{\{i\}}}^\n
	\nabla_{\bm{x}_{\{j\}}}^\m
	\k_{\f\f}(\bm{x}_{\{i\}},\bm{x}_{\{j\}})$, 
	\item[] \textbf{(\S~\ref{App:C2_GP_regression_Par_deriv})} show how these functions can be used in a GP regression context with points \textit{and directional derivatives},
	\item[] \textbf{(\S~\ref{App:C3})} gives the distribution of the derivatives of a GP (SE-ARD) w.r.t. its hyperparameters, 
	\item[] \textbf{(\S~\ref{App:C4})} give the upper bounds on the derivatives of a GP (SE-ARD) w.r.t. its hyperparameters, 
	\item[] \textbf{(\S~\ref{App:Predict_Val_and_Grad})}  shows how the gradient of a regressed function can be predicted with a GP (SE-ARD). 
\end{itemize}

\subsection{Derivatives of the kernel SE-ARD}
\label{App:DirectionalDerivatives}

The directional derivative of $f_p$ w.r.t. the direction $\bm{\v}\in\amsmathbb{R}^{n_x}$ is defined as: 
\begin{equation}
	\nabla_{\bm{\v}}f_p(\bm{x}) := \lim_{\delta\rightarrow 0} \frac{f_p(\bm{x}+\delta\bm{\v}) - f_p(\bm{x})}{\delta}. 
\end{equation} 
The covariance functions $\k_{\f\f^1}(\bm{x}_i,\bm{x}_j,\bm{\v}_j)$,  $\k_{\f^1\f}(\bm{x}_i,\bm{x}_j,\bm{\v}_i)$, and  $\k_{\f^1\f^1}(\bm{x}_i,$ $\bm{x}_j,\bm{\v}_i,\bm{\v}_j)$ that characterize the covariance of the pairs $(f_p(\bm{x}_i)$ $,$ $\nabla_{\bm{\v}_j}f_p(\bm{x}_j))$,  $\big(\nabla_{\bm{\v}_i}f_p(\bm{x}_i))$ $,$ $f_p(\bm{x}_j)\big)$, and 
$\big(\nabla_{\bm{\v}_i}f_p(\bm{x}_i))$ $,$ $\nabla_{\bm{\v}_j}f_p(\bm{x}_j))$ respectively are given in  \cite{Solak:2003}\footnote{In fact \cite{Solak:2003} gives the covariance functions between points and ``full'' derivatives, i.e. not directional ones. It is however quite simple to deduce the covariance functions between points and ``directional'' derivatives, which is done here.} and are the following: 
\begingroup
\allowdisplaybreaks
\begin{align}
	\k_{\f\f^{1}}(\bm{x}_i,\bm{x}_j,\bm{\v}_j) :=  \ & \nabla_{\bm{\v}_j}\k_{\f\f}(\bm{x}_i,\bm{x}_j), \nonumber \\
	=  \ &   \bm{\v}_j^{\rm T} \bm{\L}(\bm{x}_i-\bm{x}_j)\k_{\f\f}(\bm{x}_i,\bm{x}_j), \label{eq:kff1} \\
	\k_{\f^{1}\f}(\bm{x}_i,\bm{x}_j,\bm{\v}_i) :=  \ &
	\nabla_{\bm{\v}_i}\k_{\f\f}(\bm{x}_i,\bm{x}_j) ,
	\nonumber \\
	=  \ & 
	-\bm{\v}_i^{\rm T} \bm{\L}(\bm{x}_i-\bm{x}_j)
	\k_{\f\f}(\bm{x}_i,\bm{x}_j). \label{eq:kf1f} \\
	\k_{\f^{1}\f^{1}}(\bm{x}_i,\bm{x}_j,\bm{\v}_i, \bm{\v}_j) := \	& \nabla_{\bm{\v}_j}\nabla_{\bm{\v}_i}\widehat{\k}_{\f\f}(\bm{x}_i,\bm{x}_j), \nonumber \\
	= \ &\bm{\v}_j^{\rm T}
	\big(
	\bm{\L}(\bm{I}_{n_x} - (\bm{x}_i-\bm{x}_j)(\bm{x}_i-\bm{x}_j)^{\rm T}
	\bm{\L})
	\big)\bm{\v}_i
	\k_{\f\f}(\bm{x}_i,\bm{x}_j). \label{eq:App_Kf1f1}
\end{align}
\endgroup
And in the presence of measurement noise:
\begin{align*}
	\widehat{\k}_{\f\f}(\bm{x}_i,\bm{x}_j) =\ & \k_{\f\f}(\bm{x}_i,\bm{x}_j) + \delta_{ij} \sigma_{xy,i}^2, \\ 
	\widehat{\k}_{\f\f^{1}}(\bm{x}_i,\bm{x}_j,\bm{\v}_j) =\ & \k_{\f\f^{1}}(\bm{x}_i,\bm{x}_j,\bm{\v}_j), \\ 
	\widehat{\k}_{\f^{1}\f}(\bm{x}_i,\bm{x}_j,\bm{\v}_j) =\ & \k_{\f^{1}\f}(\bm{x}_i,\bm{x}_j,\bm{\v}_j), \\ 
	\widehat{\k}_{\f^{1}\f^{1}}(\bm{x}_i,\bm{x}_j,\bm{\v}_i, \bm{\v}_j) = \ &
	\k_{\f^{1}\f^{1}}(\bm{x}_i,\bm{x}_j,\bm{\v}_i, \bm{\v}_j) + \delta_{ij}\sigma_{\v,i}^2,
\end{align*}
where $\sigma_{\v,i}^2$ is the variance of the uncertainty on the measurements of $\nabla_{\bm{\v}_i}f_p(\bm{x}_i)$.

\subsection{Regression of a GP with points and directional derivatives}
\label{App:C2_GP_regression_Par_deriv}
Let's assemble the available data, i.e. points $\{\widehat{\bm{x}}_i, \widehat{y}_{p,i}, \sigma_{xy,i}^2\}$ and directional derivatives $\{\widehat{\bm{x}}_j, \widehat{\nabla_{\v_j} y}_{p,j},$ $ \bm{\v}_{j}, \sigma_{\v,j}^2\}$, in the following way: 
\begin{align}
	\bm{X}  := \ & (\bm{x}_1, \bm{x}_2, ..., \bm{x}_n)^{\rm T}, &
	\bm{y}_p := \ & (\widehat{y}_{p,1}, \widehat{\nabla_{\v_2} y}_{p,2}, ..., \widehat{y}_{p,n})^{\rm T}, \nonumber \\
	\bm{\V} := \ & (\emptyset, \bm{\v}_2, ..., \emptyset)^{\rm T}, & 
	\bm{\sigma} :=  \ & (\sigma_{xy,1}^2, \sigma_{\v,2}^2, ..., \sigma_{xy,n}^2)^{\rm T}. \label{eq:C_6_Matrices_de_donnees}
\end{align}
Then the location of the curve $f_p$ regressed at a point $\bm{x}_*\in\amsmathbb{R}^{n_x}$ is: 
\begin{align} 
	\amsmathbb{E}[f_p(\bm{x}_*)|\bm{x}_*,\bm{X},\bm{y}_p,\bm{\V},\bm{\sigma}] = \ & f(\bm{x}_*) + \bm{\cc}_{*}^{\rm T}\bm{\CC}^{-1}(\bm{y}_p-\bm{y}), \nonumber \\
	\amsmathbb{V}[f_p(\bm{x}_*)|\bm{x}_*,\bm{X},\bm{y}_p,\bm{\V},\bm{\sigma}] = \ &
	\cc_{**} - \bm{\cc}_{*}^{\rm T}\bm{\CC}^{-1}\bm{\cc}_{*}. \label{eq:ConditionnalDistributionGP_Gradient}
\end{align}
where:
\begin{align}
	\bm{\CC} := \ &
	\left(
	\begin{array}{ccc}
		\cc(\bm{x}_1,\bm{x}_1|\bm{\V},\bm{\sigma}) & \hdots &  \cc(\bm{x}_1,\bm{x}_n|\bm{\V},\bm{\sigma}) \\
		\vdots & \ddots & \vdots \\
		\cc(\bm{x}_n,\bm{x}_1|\bm{\V},\bm{\sigma})                  &           \hdots         & \cc(\bm{x}_n,\bm{x}_n|\bm{\V},\bm{\sigma})
	\end{array}
	\right), & 
	\bm{\cc}_{*} := \ & \left(
	\begin{array}{c}
		\cc(\bm{x}_*,\bm{x}_1|\bm{\V},\bm{\sigma}) \\
		\vdots \\
		\cc(\bm{x}_*,\bm{x}_n|\bm{\V},\bm{\sigma}) 
	\end{array}
	\right), \nonumber  \\
	\bm{y} := \ & (f(\bm{x}_{1}), \nabla_{\bm{\v}_2} f(\bm{x}_{2}), ..., f(\bm{x}_{n}))^{\rm T}, &
	\cc_{**} := \ & \cc(\bm{x}_*,\bm{x}_*), 
	\label{eq:C_8_verjrhcjdsx}
\end{align}
with: 
\begin{equation}
	\cc (\bm{x}_i,\bm{x}_j|\bm{\D},\bm{\sigma}) := 
	\left\{
	\begin{array}{l l}
		\widehat{\k}_{\f\f}(\bm{x}_i,\bm{x}_j), & \text{if case 1,} \\
		\widehat{\k}_{\f\f^{1}}(\bm{x}_i,\bm{x}_j,\bm{\v}_j), & \text{if case 2,} \\
		\widehat{\k}_{\f^{1}\f}(\bm{x}_i,\bm{x}_j,\bm{\v}_i), & \text{if case 3,} \\
		\widehat{\k}_{\f^{1}\f^{1}}(\bm{x}_i,\bm{x}_j,\bm{\v}_i,\bm{\v}_j), & \text{if case 4.} 
	\end{array}
	\right.
\end{equation}
where each case corresponds to:
\vspace{-\topsep}
\begin{itemize}[noitemsep]
	\item[] \textbf{Case 1:} The $i$-th and $j$-th data are points;
	\item[] \textbf{Case 2:} The $i$-th datum is a point and the $j$-th datum is a directional derivative; 
	\item[] \textbf{Case 2:} The $i$-th datum is a directional derivative and the $j$-th datum is a point;
	\item[] \textbf{Case 4: } The $i$-th and $j$-th data are directional derivatives
\end{itemize}
In other words: 
\begin{equation} \label{eq:C_Fp_a_posteriori}
	f_P(\bm{x}_*) \sim \mathcal{N}(
	f(\bm{x}_* + \bm{\cc}_{*}^{\rm T}\bm{\CC}^{-1}(\bm{y}_p-\bm{y})
	,
	\cc_{**} - \bm{\cc}_{*}^{\rm T}\bm{\CC}^{-1}\bm{\cc}_{*}
	)).
\end{equation}

\subsection{A priori expectation of the value, gradients, and Hessians of the plant}
\label{App:C3}

According to  equation~\eqref{eq:App_Kf1f1}, in the absence of data the variance of the directional derivative of $f_p$ in the direction of $\bm{\v}_i$ is:
\begin{equation} \label{eq:C_10_var_F1}
	\k_{\f^{1}\f^{1}}(\bm{x}_i,\bm{x}_j=\bm{x}_i,\bm{\v}_i, \bm{\v}_j=\bm{\v}_i) = \sigma_f^2/\ell_{\bm{\v}_i}^2,
\end{equation}
where
\begin{equation*}
	\ell_{\bm{\v}_i}^{-2} := \bm{\v}_i^{\rm T} \bm{\L} \bm{\v}_i.
\end{equation*}
By similar reasoning, the covariance of the couple  $\big(\nabla^2_{\bm{\v}_i}f_p(\bm{x}_i))$ $,$ $\nabla^2_{\bm{\v}_j}f_p(\bm{x}_j)\big)$ is:  
\begin{align*}
	\k_{\f^{2}\f^{2}}(\bm{x}_i,\bm{x}_j,\bm{\v}_i, \bm{\v}_j) :=  
	\nabla^2_{\bm{\v}_j}\nabla^2_{\bm{\v}_i}\k_{\f\f}(\bm{x}_i,\bm{x}_j). 
\end{align*}
Therefore, the a priori variance of  $\nabla^2_{\bm{\v}_i}f_p(\bm{x}_i)$ is:
\begin{align} \label{eq:C_11_var_F2}
	\k_{\f^{2}\f^{2}}(\bm{x}_i,\bm{x}_j=\bm{x}_i,\bm{\v}_i, \bm{\v}_j=\bm{\v}_i) =  3 \sigma_f^2/\ell_{\bm{\v}_i}^4
\end{align}
From \eqref{eq:C_10_var_F1} and \eqref{eq:C_11_var_F2}, one can derive the a priori distributions of the gradients and Hessians of the plant   $\forall \bm{x}_i\in\amsmathbb{R}^{n_x}$. 
\begin{equation}
	\begin{array}{l@{}l l@{}l}
	\amsmathbb{E}[\nabla_{\bm{\v}_i}f_p(\bm{x}_i)] = \ & \nabla^1_{\bm{\v}_i}f(\bm{x}_i), & \qquad \amsmathbb{V}[\nabla_{\bm{\v}_i}f_p(\bm{x}_i)] = \ &  \sigma_f^2/\ell_{\bm{\v}_i}^2, \\
	\amsmathbb{E}[\nabla^2_{\bm{\v}_i}f_p(\bm{x}_i)] = \ & \nabla^2_{\bm{\v}_i}f(\bm{x}_i), 
	& 
	\qquad \amsmathbb{V}[\nabla^2_{\bm{\v}_i}f_p(\bm{x}_i)] = \ & 3 \sigma_f^2/\ell_{\bm{\v}_i}^4.
	\end{array} \label{App:C_ExpectedValuesGrad}
\end{equation}

\subsection{A priori upper bounds of the plant gradients}
\label{App:C4}
Let's rewrite \eqref{App:C_ExpectedValuesGrad} as:
\begin{equation*}
	\nabla_{\bm{\v}_i}f_p(\bm{x}) \sim \mathcal{N}(\nabla_{\bm{\v}_i}f(\bm{x}), 
	\sigma_f^2/\ell_{\bm{\v}_i}^{2}
	).
\end{equation*}
So the absolute value of $\nabla_{\bm{\v}_i}f_p(\bm{x})$ satisfies the following (triangular) inequality:
\begin{align*}
	|\nabla_{\bm{\v}_i}f_p(\bm{x})| = \ &  |\nabla_{\bm{\v}_i}f(\bm{x}) + \epsilon|
	\leq  |\nabla_{\bm{\v}_i}f(\bm{x})| + |\epsilon|, &  \epsilon \sim \mathcal{N}(0,\sigma_f^2/\ell^2_{\bm{\v}_i}) 
\end{align*}
It is clear that $|\epsilon|$ follows a half-normal distribution. According to the properties of such a distribution \cite{Leone:1961}, one finds that:
\begin{equation*}  
	p\left(|\epsilon| \leq \frac{3}{\sqrt{2}}\frac{\sigma_f}{\ell_{\bm{\v}_i}} \right) \approx 0.966.
\end{equation*}
So one can expect that (with a probability of 96.6\%): 
\begin{equation*}
	|\nabla_{\bm{\v}_i}f_p(\bm{x})| \leq 
	|\nabla_{\bm{\v}_i}f(\bm{x}) | 
	+ \frac{3}{\sqrt{2}}\frac{\sigma_f}{\ell_{\bm{\v}_i}}.
\end{equation*}
By generalizing this result to all directions of the input space: 
\begin{align} \label{eq:C_13_wjhgvfrcje}
	|\nabla_{\bm{x}}f_p(\bm{x})| \leq \ & 
	|\nabla_{\bm{x}}f(\bm{x}) | 
	+ \frac{3}{\sqrt{2}}\sigma_f\bm{\ell}^{-1}, &
	 \bm{\ell}^{-1}:= \ & [\ell_1^{-1},...,\ell_{n_x}^{-1}]^{\rm T}.
\end{align}

\subsection{The expectation and variance of a gradient of a GP regression (SE-ARD-WN)}
\label{App:Predict_Val_and_Grad}
 
To calculate the expectation and variance of the value and gradient of the regression of $f_p$ at $\bm{x}_*$,  one has to solve the following system:
\begingroup
\allowdisplaybreaks
\begin{align}
	\amsmathbb{E}
	\left[\left(\begin{array}{c}
		f_{p} \\  \nabla_{\bm{x}} f_{p}
	\end{array}
	\right) \right] = \ &
	\left(\begin{array}{c}
		f|_{\bm{x}_*} \\  \nabla_{\bm{x}} f|_{\bm{x}_*}
	\end{array}
	\right) + \bm{\CC}_{*}^{\rm T} \bm{\CC}^{-1}(\bm{y}_p-\bm{y}), \label{eq:C_App_CovMat_F_DF_1} \\
	\amsmathbb{V}
	\left[\left(\begin{array}{c}
		f_{p} \\  \nabla_{\bm{x}} f_{p}
	\end{array}
	\right) \right] = \ & \bm{\CC}_{**} - \bm{\CC}_{*}^{\rm T} \bm{\CC}^{-1}\bm{\CC}_{*}, \label{eq:C_App_CovMat_F_DF}
\end{align}
\endgroup
where $\bm{\CC}$ is defined by \eqref{eq:C_8_verjrhcjdsx} and: 
\begingroup
\allowdisplaybreaks
\begin{align}
	\bm{\CC}_{**} := \ &
	\left(\begin{array}{cc}
		\cc_{**}^{\prime}      & \bm{\cc}_{**}^{\prime \rm T} \nonumber\\
		\bm{\cc}_{**}^{\prime}  & \bm{\CC}_{**}^{\prime} 
	\end{array}\right),\\
	\cc_{**}^{\prime} :=\ &   \k_{\f\f}(\bm{x}_*,\bm{x}_*), \nonumber\\ 
	\bm{\cc}_{**}^{\prime \rm T} := \ & \left[\k_{\f\f^1}(\bm{x}_*,\bm{x}_*,\bm{\v}_1),
	\hdots,
	\k_{\f\f^1}(\bm{x}_*,\bm{x}_*,\bm{\v}_{n_u}) \right],  \nonumber\\
	\bm{\CC}_{**}^{\prime} := \ &
	\left[\begin{array}{ccc}
		\k_{\f^1\f^1}(\bm{x}_*,\bm{x}_*,\bm{\v}_1,\bm{\v}_1)
		& \hdots
		&  \k_{\f^1\f^1}(\bm{x}_*,\bm{x}_*,\bm{\v}_1,\bm{\v}_{n_u}) \\
		\vdots & \ddots & \vdots \\
		\k_{\f^1\f^1}(\bm{x}_*,\bm{x}_*,\bm{\v}_{n_u},\bm{\v}_1) & \hdots &
		\k_{\f^1\f^1}(\bm{x}_*,\bm{x}_*,\bm{\v}_{n_u},\bm{\v}_{n_u})
	\end{array}\right] \nonumber\\
	\bm{\CC}_{*} := \ & 
	\left[\begin{array}{ccc}
		\cc(\bm{x}_*,\bm{x}_1|\bm{\V},\bm{\sigma}) 
		& \hdots
		& \cc(\bm{x}_*,\bm{x}_{1}|\bm{\V},\bm{\sigma}) \\
		\vdots &  \ddots & \vdots \\
		\cc(\bm{x}_*,\bm{x}_{n}|\bm{\V},\bm{\sigma})  & \hdots &
		\cc(\bm{x}_*,\bm{x}_{n}|\bm{\V},\bm{\sigma})
	\end{array}\right]. \label{eq:A_iywbvwe}
\end{align}
\endgroup
Finally, let us name some values of interest of the covariance matrix \eqref{eq:C_App_CovMat_F_DF}:
\begin{equation} \label{eq:C_17_Variances_et_Covariances}
	\bm{\CC}_{**} - \bm{\CC}_{*}^{\rm T} \bm{\CC}^{-1}\bm{\CC}_{*} =  
	\left[\begin{array}{cc}
		\sigma_*^2      & \bm{\varsigma}_*^{\rm T} \\
		\bm{\varsigma}_*  & \bm{\Sigma_*}
	\end{array}\right],
\end{equation}
where $\sigma_*^2\in\mathbb{R}$ is the variance of $f_p$ at $\bm{x}_*$, $\bm{\Sigma}_*\in\amsmathbb{R}^{n_x \times n_x}$ is the covariance matrix of  $\nabla_{\bm{x}} f_p$ at $\bm{x}_*$, and $\bm{\varsigma}_*\in\amsmathbb{R}^{n_x}$ is the covariance vector between the value and gradient of  $f_p$ at $\bm{x}_*$.

\section{Domain on which a PG (SE-ARD) is quasi-affine}
\label{App:AffineLikeDomain}

To determine the size of a domain over which a GP behaves as an affine function, let us take two points and see how much one can separate them until the affine approximation made at one point is not sufficient to predict the location of the other one. Consider the following two points: 
$(\bm{x}_{1}, y_{p,1})$ and $(\bm{x}_{2}, y_{p,2})$, and let's define:
\begin{align*}
	\Delta \bm{x} := & \bm{x}_{2} - \bm{x}_{1}, &
	\Delta \bm{y} := & y_{p,2} - y_{p,1}.
\end{align*}
Then let's express the value of $y_{p,2}$ with a Taylor series of  $f_p$ centered on $\bm{x}_{1}$: 
\begin{align} \label{eq:C_1_Val_reelle}
	y_{p,2} = \  &  f_p(\bm{x}_{1}) 
	+
	\|\Delta\bm{x}\|\nabla_{\Delta \bm{x}}f_p(\bm{x}_{1}) 
	+ 
	\|\Delta\bm{x}\|^2 \frac{\nabla_{\Delta\bm{x}}^2 f_p (\bm{x}_{1})}{2} 
	+ 
	\mathcal{O}(\|\Delta\bm{x}\|^3), 
\end{align}
where $\nabla^{\n}_{\Delta x}f_p(\bm{x}_{1})\in\amsmathbb{R}$ is the $\n$-th directional derivative of $f_p$ at $\bm{x}_{1}$ w.r.t. the direction $\Delta x$. 

If one assumes that $f_p$ has an affine behavior on the segment $[\bm{x}_{1}, \bm{x}_{2}]$, the the expected value $\widetilde{y}_{\{2\}}$ of $y_{p,2}$ would be:
\begin{align} \label{eq:C_2_Val_approx}
	\widetilde{y}_{p,2} = \  &  f_p(\bm{x}_{1}) 
	+
	\|\Delta\bm{x}\| \nabla_{\Delta \bm{x}}f_p(\bm{x}_{1}).
\end{align}

Now, what one would like to know is how well this assumption works w.r.t. $\|\Delta\bm{x}\|$. To do so, let's combine \eqref{eq:C_1_Val_reelle} with \eqref{eq:C_2_Val_approx} and neglect the term $\mathcal{O}(\|\Delta\bm{x}\|^3)$. A series of elementary calculations gives: 
\begin{align} \label{App:C_RelativeLinError}
	\frac{|\widetilde{y}_{p,2} -  y_{p,2}|}{|y_{p,2}|}= \  & 
	\frac{\|\Delta\bm{x}\|^2}{2} \frac{|\nabla_{\Delta\bm{x}}^2 f_p (\bm{x}_{1})|}{|f_p(\bm{x}_{1} + \Delta\bm{x})|} . 
\end{align}
One notices that this relative error depends on the choice of $\bm{x}_1$. However, one is looking for a relative ``global'' error that would work regardless of $\bm{x}_1$. To get what one is looking for let's evaluate separately the distributions of (i) $|\nabla_{\Delta\bm{x}}^2 f_p (\bm{x}_{1})|$ and (ii) $|f_p(\bm{x}_1 + \Delta\bm{x})|$  when $\bm{x}_{1}$ takes all possible values that a $\bm{x}$ can take. 
\begin{itemize}
	\item[(i)] According to \eqref{App:C_ExpectedValuesGrad}:
	\begin{equation} \label{eq:C_App_HessianPlantFunction}
		\nabla_{\Delta\bm{x}}^2 f_p (\bm{x}_{1}) \sim \mathcal{N}(
		\nabla_{\Delta\bm{x}}^2 f (\bm{x}_{1}) ,
		3 \sigma_f^2/\ell_{\Delta\bm{x}}^4
		).
	\end{equation}
	In addition, $f$ can be approximated by a GP (SE-ARD) whose hyperparameters $\sigma_f^{\prime 2}$ and $\bm{\L}$. The matrix $\bm{\L}$ is, according to Hypothesis~\ref{ass:6_2_Meme_Courbusre}, common to the GPs of $f$ and $f_p$. So:
	\begin{equation} \label{eq:C_App_HessianModelFunction}
		\nabla_{\Delta\bm{x}}^2 f(\bm{x}) \sim \mathcal{N}( 0 ,
		3 \sigma_f^{\prime 2}/\ell_{\Delta\bm{x}}^4
		), 
	\end{equation}
	where the parameter $\sigma_f^{\prime 2}$ is the variance of the function $f$, i.e. it is a quantification of the domain of the output space it occupies (note that if this variable plays a role in these ``intermediate'' calculations, it disappears from the final result). Combining \eqref{eq:C_App_HessianPlantFunction} and  \eqref{eq:C_App_HessianModelFunction} obtains an estimate of the Hessian distribution of the plant at a point $\bm{x}\in\amsmathbb{R}^{n_x}$: 
	\begin{equation}
		\nabla_{\Delta\bm{x}}^2 f_p(\bm{x}) \sim \mathcal{N}( 0 ,
		3 (\sigma_f^{\prime 2}+\sigma_f^2) /\ell_{\Delta\bm{x}}^4
		).
	\end{equation}
	From this result, one concludes that $|\nabla_{\Delta\bm{x}}^2 f_p(\bm{x})|$ follows a half-normal distribution whose known properties (\cite{Leone:1961}) can be used to prove that:
	\begin{align*}
		p\left(|\nabla_{\Delta\bm{x}}^2 f_p(\bm{x})| \leq a\right) = \erf\left(a \left(2\cdot 3 (\sigma_f^{\prime 2}+\sigma_f^2) /\ell_{\Delta\bm{x}}^4\right)^{-1/2}\right).
	\end{align*}
	Since $\erf(1.5)\approx 0.966$:
	\begin{equation}  \label{App:C_LimiteCurevature_absFp}
		p\left(|\nabla_{\Delta\bm{x}}^2 f_p(\bm{x}_1)| \leq 1.5 \sqrt{2\cdot 3 (\sigma_f^{\prime 2}+\sigma_f^2) /\ell_{\Delta\bm{x}}^4} \right) \approx 0.966.
	\end{equation}
	So $1.5 \sqrt{2\cdot 3 (\sigma_f^{\prime 2}+\sigma_f^2) /\ell_{\Delta\bm{x}}^4}$ is a probable upper bound (at 96.6\%) of $|\nabla_{\Delta\bm{x}}^2 f_p(\bm{x}_1)|$. 
	\item[(ii)]  With a similar approach: 
	\begin{align*}
		&	\left\{
		\begin{array}{l}
			f_p(\bm{x}_1) \sim \mathcal{N}(
			f(\bm{x}_1),
			\sigma_f^2) \\
			f(\bm{x}) \sim \mathcal{N}(0,
			\sigma_f^{\prime2})
		\end{array}
		\right. ,  & 
		\Rightarrow \	&
		f_p(\bm{x}) \sim \mathcal{N}(
		0,
		\sigma_f^2 + \sigma_f^{\prime2}). 
	\end{align*}
	So, $|f_p(\bm{x})|$ follows a half-normal distribution and according to \cite{Leone:1961}:
	\begin{equation} \label{App:C_Expectancy_Fp}
		\amsmathbb{E}\left[\left| f_p(\bm{x}) \right|\right] = \sqrt{\frac{2}{\pi}} (\sigma_f^2 + \sigma_f^{\prime2}).
	\end{equation} 
\end{itemize}
To estimate \eqref{App:C_RelativeLinError} it is suggested to use \eqref{App:C_LimiteCurevature_absFp} and \eqref{App:C_Expectancy_Fp} as approximations of $|\nabla_{\Delta\bm{x}}^2 f_p (\bm{x}_{1})|$ and $|f_p(\bm{x}_1 + \Delta\bm{x})|$ respectively. Then one finds that:
\begin{align} \label{App:C_RelativeLinError_Global}
	\frac{|\widetilde{y}_ -  y|}{|y|} \leq \  & 
	\frac{3\sqrt{3\pi}}{4}\frac{\|\Delta\bm{x} \|^2}{\ell_{\Delta\bm{x}}^2}, 
\end{align}
is an estimate of the upper bound  of the relative prediction error at a point $\bm{x}_2$ located at a distance $\|\Delta\bm{x}\|$ of a point $\bm{x}_1\in\amsmathbb{R}^{n_x}$ when approximating $y_{p,2} := f_p(\bm{x}_2)$ with \eqref{eq:C_2_Val_approx}.

\section{A derivative}

Let the following Taylor series approximation be used: 
\begin{equation}
	\left\{ 
	\begin{array}{l}
		\phi(\bm{u},\bm{y}+\delta\bm{y}) = 
		\phi(\bm{u},\bm{y}) + \sum_{i=1}^{n_y}\partial_{y_{(i)}} \phi(\bm{u},\bm{y}) \delta y_{(i)}, \\
		\delta y_{(i)} = \epsilon_{i} + \bm{\tau}_i^{\rm T} (\bm{u}-\bm{u}_k), \\
		\bm{y} = \bm{f}(\bm{u}).
	\end{array}
	\right.
\end{equation}
The derivative of $	\phi(\bm{u},\bm{y}+\delta\bm{y})$ is:
\begin{equation*}
	\nabla_{\bm{u}} 
	\phi(\bm{u},\bm{y}+\delta\bm{y}) = 
	\nabla_{\bm{u}}\phi(\bm{u},\bm{y}) + 
	\sum_{i=1}^{n_y} \big(
	\nabla_{\bm{u}} \partial_{y_{(i)}} \phi(\bm{u},\bm{y}) \delta y_{(i)} + 
	\partial_{y_{(i)}} \phi(\bm{u},\bm{y}) \nabla_{\bm{u}} \delta y_{(i)}
	\big).
\end{equation*}
Let's analyze each of the terms: 
\begin{align*}
	\nabla_{\bm{u}} \phi(\bm{u},\bm{y}) = \ & 
	\partial_{\bm{u}} \phi(\bm{u},\bm{y}) +
	\sum_{i=1}^{n_y} \big( 
	\nabla_{\bm{u}}f_{(i)}
	\partial_{y_{(i)}}\phi(\bm{u},\bm{y})
	\big), \\ 
	\nabla_{\bm{u}} \partial_{y_{(i)}} \phi(\bm{u},\bm{y}) = \ & 
	\partial_{\bm{u}} \partial_{y_{(i)}} \phi (\bm{u},\bm{y}) + 
	\nabla_{\bm{u}} \bm{f} \partial_{\bm{y}} \partial_{y_{(i)}} \phi (\bm{u},\bm{y}),  \\
	\nabla_{\bm{u}} \delta y_{(i)} = \ & \bm{\tau}_i.
\end{align*} 
From this one can find that:
\begin{align}
	\nabla_{\bm{u}} 
	\phi(\bm{u},\bm{y}+\delta\bm{y}) = \ & 
	\partial_{\bm{u}} \phi(\bm{u},\bm{y}) +
	\sum_{i=1}^{n_y} \big( 
	\nabla_{\bm{u}}f_{(i)}
	\partial_{y_{(i)}}\phi(\bm{u},\bm{y})
	\big) +  ...  \nonumber \\
	&
	\sum_{i=1}^{n_y}  \Big( 
	\big( 
	\partial_{\bm{u}} \partial_{y_{(i)}} \phi(\bm{u},\bm{y}) + 
	\nabla_{\bm{u}} \bm{f}^{\rm T} \partial_{\bm{y}} \partial_{y_{(i)}}\phi (\bm{u},\bm{y})
	\big)
	\epsilon_{i} + ... \nonumber \\
	& \partial_{y_{(i)}} \phi(\bm{u},\bm{y})
	\bm{\tau}_i
	\Big).
\end{align}

\end{document}